\documentclass[twoside,a4paper,11pt,openright]{report}
\pdfoutput=1 
\usepackage[english]{babel}
\usepackage[utf8]{inputenc}
\usepackage[T1]{fontenc}
\usepackage{titlesec,fancyhdr}
\usepackage{placeins} 
\usepackage{booktabs}
\usepackage{microtype}
\usepackage{subcaption}
\usepackage{lscape,rotating}
\usepackage[twoside,heightrounded=true,inner=4.25cm,outer=3.25cm,top=2.5cm,bottom=3cm]{geometry}
\usepackage{tikz}
\usepackage{tikz-cd}
\usepackage{torbenphd}

\usepackage{mathpazo}
\usepackage{dutchcal} 

\date{March 28, 2024}

\crefname{subsection}{subsection}{subsections}

\let\phi\varphi

\usepackage[backend=bibtex, style=alphabetic]{biblatex}
\renewbibmacro{in:}{%
  \ifentrytype{article}{}{\printtext{\bibstring{in}\intitlepunct}}}
\DeclareFieldFormat{postnote}{#1} 
\DeclareFieldFormat{multipostnote}{#1} 
\DeclareFieldFormat[article]{volume}{\bibstring{jourvol}\addnbspace #1} 

\setlist[enumerate]{labelsep=5pt, itemsep=1pt, parsep=1.5pt, topsep=3pt}

\usetikzlibrary{calc}
\usetikzlibrary{decorations.markings}

\tikzset{
	coxcircle/.style={draw, fill, circle, inner sep=0pt, minimum size=2mm}
}
\tikzset{
	coxdoublearrow/.style={double distance=2pt, decoration={
				markings,
	mark=at position 0.5 with {\node{$ > $};}},
	postaction={decorate}
	}
}
\tikzset{
	nodecircle/.style={draw, fill, solid, circle, inner sep=0pt, minimum size=1.5mm}
}

\renewcommand{\epsilon}{\varepsilon}

\newcommand{\frmon}[1]{\calL(#1)}
\newcommand{\frgrp}[1]{F(#1)}
\newcommand{\mapdot}{\mathord{\:\cdot\:}}
\newcommand{\numint}[2]{[#1, #2]}
\newcommand{\betint}[2]{\langle #1, #2 \rangle}
\newcommand{\Nzero}{\IN_0}
\newcommand{\Npos}{\IN_+}

\newcommand{\listing}[2]{#1_1, \ldots, #1_{#2}}
\newcommand{\tup}[2]{(#1_1, \ldots, #1_{#2})}
\newcommand{\zentrum}{Z}
\DeclareMathOperator{\str}{str}
\newcommand{\inv}[1]{#1^\times}
\newcommand{\delmin}[1]{\delta^-_{#1}}
\newcommand{\delinv}[1]{\delta^{\rinvmap}_{#1}}
\newcommand{\len}{\ell}
\DeclareMathOperator{\sgn}{sgn}
\newcommand{\algclass}[1]{\calC_{#1}}
\newcommand{\ringinvset}[1]{#1^{\times}}
\newcommand{\lexle}{\le_{\text{lex}}}
\newcommand{\lexls}{<_{\text{lex}}}

\DeclareMathOperator{\Weyl}{Weyl}
\DeclareMathOperator{\rootht}{ht}
\newcommand{\refl}[1]{\sigma_{#1}}
\newcommand{\reflbr}[1]{\sigma(#1)}
\newcommand{\cartanint}[2]{\langle #1 \vert #2 \rangle}
\newcommand{\roots}{\Phi}
\newcommand{\rootsub}{\Psi}
\newcommand{\rootbase}{\Delta}
\newcommand{\possys}{\Pi}
\newcommand{\rootint}[2]{\mathord{] #1, #2 [}}
\newcommand{\rootintcox}[2]{\rootint{#1}{#2}_{\text{Cox}}}

\newcommand{\varrootint}[2]{\mathord{\langle #1, #2 \rangle}}
\newcommand{\clrootint}[2]{[#1, #2]}
\newcommand{\clrootintcox}[2]{\clrootint{#1}{#2}_{\text{Cox}}}
\newcommand{\basvec}{e}
\newcommand{\braidword}[1]{P_{#1}}

\newcommand{\switchset}{N}
\newcommand{\indivset}[1]{#1^{\text{indiv}}}
\DeclareMathOperator{\Orb}{Orb}
\newcommand{\IndivOrb}{\Orb^{\text{red}}}
\DeclareMathOperator{\cone}{Cone}
\newcommand{\totorder}{\mathrel{\triangleleft}}
\newcommand{\fixspace}{F}

\newcommand{\weylset}[1]{M_{#1}}
\newcommand{\rootgr}[1]{U_{#1}}
\newcommand{\varrootgr}[1]{\hat{U}_{#1}}
\newcommand{\varvarrootgr}[1]{\bar{U}_{#1}}
\newcommand{\rootgrV}[1]{V_{#1}}
\newcommand{\invset}[1]{\rootgr{#1}^{\smash{\sharp}}}
\newcommand{\invsetV}[1]{\rootgrV{#1}^{\smash{\sharp}}}
\newcommand{\commpart}[4][]{\commutator[#1]{#2}{#3}_{#4}}
\newcommand{\risom}[1]{\theta_{#1}}

\newcommand{\varvarrisom}[1]{\bar{\theta}_{#1}}
\newcommand{\risomnaiv}[1]{\risom{#1}^{\text{n}}}
\newcommand{\whom}[1]{w_{#1}}
\newcommand{\whomnaiv}[1]{\whom{#1}^{\text{n}}}
\newcommand{\weyltrip}{\mathfrak{u}}
\newcommand{\weylmap}[1]{\ifblank{#1}{\mu}{\mu_{#1}}}
\newcommand{\weylleft}[1]{\ifblank{#1}{\kappa}{\kappa_{#1}}}
\newcommand{\weylright}[1]{\ifblank{#1}{\lambda}{\lambda_{#1}}}

\DeclareMathOperator{\St}{St}
\newcommand{\steinhom}[1]{\hat{x}_{#1}}
\newcommand{\steinweyl}[1]{\hat{w}_{#1}}
\newcommand{\steinrootgr}[1]{\hat{U}_{#1}}
\newcommand{\steinproj}{\phi}

\newcommand{\chevstr}{c}
\newcommand{\chevbasis}{\calC}

\newcommand{\chevorder}{\calI}
\DeclareMathOperator{\SL}{SL}
\DeclareMathOperator{\sllie}{\mathfrak{sl}}

\newcommand{\inverparsym}{\eta}
\newcommand{\inverpar}[2]{\inverparsym_{\ifblank{#1}{\mathord{\cdot}}{#1}, \ifblank{#2}{\mathord{\cdot}}{#2}}}
\newcommand{\inverparbr}[2]{\inverparsym(#1, #2)}
\newcommand{\invoparsym}{\mu}
\newcommand{\invopar}[2]{\invoparsym_{#1,#2}}
\newcommand{\invoparbr}[2]{\invoparsym(#1, #2)}
\newcommand{\totalparsym}{\bar{\eta}}
\newcommand{\totalpar}[2]{\totalparsym_{#1, #2}}
\newcommand{\totalparbr}[2]{\totalparsym(#1, #2)}
\newcommand{\inverinvopar}[2]{(\inverparsym \times \invoparsym)_{\ifblank{#1}{\mathord{\cdot}}{#1}, \ifblank{#2}{\mathord{\cdot}}{#2}}}
\newcommand{\inverinvoparprime}[2]{(\inverparsym' \times \invoparsym')_{\ifblank{#1}{\mathord{\cdot}}{#1}, \ifblank{#2}{\mathord{\cdot}}{#2}}}
\newcommand{\twistgroup}{A}
\newcommand{\invogroup}{B}
\newcommand{\word}[1]{\bar{#1}}
\newcommand{\genparmoveset}[4]{#1_{#2, #3 \rightarrow #4}}
\newcommand{\parmoveset}[3]{#1_{#2 \rightarrow #3}}

\newcommand{\blumapsym}{\gamma}
\newcommand{\blumap}[1]{\blumapsym_{#1}}
\newcommand{\blumapGsym}{\tilde{\blumapsym}}
\newcommand{\blumapG}[1]{\blumapGsym_{#1}}
\newcommand{\blutrans}[1]{\psi_{#1}}
\newcommand{\transp}[2]{(#1 \; #2)}

\newcommand{\rootgrmin}[2]{\rootgr{\basvec_{#1} - \basvec_{#2}}}
\newcommand{\rismin}[2]{\risom{\basvec_{#1} - \basvec_{#2}}}
\newcommand{\rootgrplus}[2]{\rootgr{\basvec_{#1} + \basvec_{#2}}}
\newcommand{\risplus}[2]{\risom{\basvec_{#1} + \basvec_{#2}}}

\newcommand{\risminmin}[2]{\risom{-\basvec_{#1} - \basvec_{#2}}}
\newcommand{\rootgrshpos}[1]{\rootgr{\basvec_{#1}}}
\newcommand{\risshpos}[1]{\risom{\basvec_{#1}}}
\newcommand{\rootgrlongpos}[1]{\rootgr{2\basvec_{#1}}}
\newcommand{\rislongpos}[1]{\risom{2\basvec_{#1}}}

\newcommand{\risshneg}[1]{\risom{-\basvec_{#1}}}

\DeclareMathOperator{\E}{E}

\newcommand{\module}{M}
\newcommand{\modinv}[1]{\overline{#1}}
\newcommand{\modinvmap}{\modinv{\mapdot}}
\newcommand{\Bnsub}{\hat{B}_n}
\DeclareMathOperator{\Ortho}{O}
\DeclareMathOperator{\EO}{EO}

\newcommand{\jormod}{J}
\newcommand{\jormodtup}{\calJ}
\newcommand{\jorproj}{\pi}
\newcommand{\jorprojone}{\pi_1}
\DeclareMathOperator{\jorTr}{T}
\DeclareMathOperator{\jorTrone}{T_1}
\newcommand{\jorTroneprime}{\jorTr_1'}
\newcommand{\jorsc}{\phi}
\newcommand{\shortinvo}[1]{#1^*}
\newcommand{\joradd}{\mathbin{\hat{+}}}
\newcommand{\jormin}{\mathbin{\hat{-}}}
\newcommand{\symring}{\Fix(\sigma)}
\newcommand{\ringnorm}{N}
\DeclareMathOperator{\ringTr}{Tr}
\newcommand{\ringzero}{\ring_0}
\DeclareMathOperator{\Unitary}{U}
\DeclareMathOperator{\EU}{EU}
\newcommand{\Cnsub}{\hat{C}_n}
\newcommand{\invsub}{D}
\newcommand{\quotmod}{\module}
\newcommand{\quotepi}{\rho}
\DeclareMathOperator{\Rad}{Rad}

\newcommand{\compalg}{\calA}
\newcommand{\varcompalg}{\calB}
\newcommand{\rootbaseB}{\rootbase^B}
\newcommand{\rootbaseC}{\rootbase^C}
\newcommand{\rootbaseF}{\rootbase^F}
\newcommand{\rootbaseE}{\rootbase^E}
\newcommand{\risomB}[1]{\risom{#1}^B}
\newcommand{\risomC}[1]{\risom{#1}^C}
\newcommand{\inverparsymB}{\inverparsym^B}
\newcommand{\invoparsymB}{\invoparsym^B}
\newcommand{\inverparsymC}{\inverparsym^C}
\newcommand{\varinverparsymC}{\tilde{\inverparsym}^C}
\newcommand{\invoparsymC}{\invoparsym^C}

\newcommand{\inverinvoparB}[2]{(\inverparsymB \times \invoparsymB)_{\ifblank{#1}{\mathord{\cdot}}{#1}, \ifblank{#2}{\mathord{\cdot}}{#2}}}

\newcommand{\inverinvoparC}[2]{(\inverparsymC \times \invoparsymC)_{\ifblank{#1}{\mathord{\cdot}}{#1}, \ifblank{#2}{\mathord{\cdot}}{#2}}}
\newcommand{\compinvmapB}{\compinvmap^B}
\newcommand{\compinvB}[1]{\compinv{#1}^B}
\newcommand{\scmult}{\star}
\newcommand{\zentruminv}{\zentrum^\rinvmap}
\DeclareMathOperator{\compTr}{tr}
\DeclareMathOperator{\compnorm}{n}
\newcommand{\compinv}[1]{\overline{#1}}
\newcommand{\compinvmap}{\compinv{\mapdot}}
\newcommand{\Forder}{\prec}

\newcommand{\ring}{\calR}
\newcommand{\varring}{\calS}
\newcommand{\comring}{\mathcal{k}}
\newcommand{\rmult}{\cdot}
\DeclareMathOperator{\nucleus}{Nucl}
\DeclareMathOperator{\lnucleus}{LNucl}
\DeclareMathOperator{\rnucleus}{RNucl}
\DeclareMathOperator{\mnucleus}{MNucl}
\DeclareMathOperator{\Sym}{Sym}
\newcommand{\rinvmap}{\sigma}
\newcommand{\rinv}[1]{#1^\rinvmap}
\newcommand{\rinvmin}[1]{#1^{-\rinvmap}}

\newcommand{\psgr}{T}

\newcommand{\onpage}[1]{p.~#1}

\makeindex

\index{sufficiently generic|see{generic example group}}
\index{standard commutator relations|see{coordinatisation of a root graded group with standard signs}}

\author{Torben Wiedemann}
\title{Root Graded Groups}

\makeatletter
\newcommand*{\sectionbookmark}{%
  \bookmark[%
    level=section,%
    dest=\@currentHref%
  ]%
}
\makeatother

\begin{document}
	\setlength{\abovedisplayskip}{0.5\abovedisplayskip}
	\setlength{\belowdisplayskip}{0.5\belowdisplayskip}

	\begin{center}
		\thispagestyle{empty}
%
%
%
%
%
		\vspace*{4cm}
		
		{\Huge Root Graded Groups}
		
		\addvspace{1.75cm}
		
%
%
%
		{\Large Torben Wiedemann}
		
		\addvspace{1.75cm}
		
		\makeatletter\@date\makeatother
%
%
%
%
%
%
%
%
%
	\end{center}
	\cleardoublepage

\vspace*{2cm}
\thispagestyle{empty}
\pdfbookmark[1]{Abstract}{Abstract}

\begin{center}
	\Large\bfseries Abstract
\end{center}

\noindent We define and study \emph{root graded groups}, that is, groups graded by finite root systems. This notion generalises several existing concepts in the literature, including in particular Jacques Tits' notion of RGD-systems. The most prominent examples of root graded groups are Chevalley groups over commutative associative rings. Our main result is that every root graded group of rank at least~3 is coordinatised by some algebraic structure satisfying a variation of the Chevalley commutator formula. This result can be regarded as a generalisation of Tits' classification of thick irreducible spherical buildings of rank at least~3 to the case of non-division algebraic structures.

All coordinatisation results in this book are proven in a characteristic-free way. This is made possible by a new computational method that we call the \emph{blueprint technique}.

%
%
%
\cleardoublepage



\pdfbookmark[1]{Acknowledgements}{Acknowledgements}
\chapter*{Acknowledgements}
\markboth{Acknowledgements}{}

First of all, I want to thank my Doktorvater Bernhard Mühlherr for leading me to this beautiful research topic and for supervising me as a PhD student. This book, which is based on my PhD thesis, would not have been possible without his constant support and guidance. He was also the one who initially suggested the idea of the blueprint technique.


Further, I want to thank Richard Weiss for answering and asking questions on this project, for sharing his thoughts on my thesis draft and for showing interest in my work. His initial note on the blueprint computations in $ C_3 $ was a big help when I wrote my Master thesis.

I am also thankful to Holger Petersson for pointing me to the notion of conic algebras and for answering many questions about composition algebras and quadratic forms over rings.

This work was supported by the DFG Grant MU1281/7-1.
\cleardoublepage
	
	\pdfbookmark{\contentsname}{toc1}
	\tableofcontents

	\newpage

\chapter*{Preface}
\addcontentsline{toc}{chapter}{Preface}
\markboth{Preface}{}

Let $ \roots $ be a finite irreducible root system. A \defemph*{$ \roots $-grading} of a group $ G $ is a family of non-trivial subgroups $ (\rootgr{\alpha})_{\alpha \in \roots} $ generating $ G $, called the \defemph*{root groups of $ G $}, such that some commutator relations and a non-degeneracy condition are satisfied and such that so-called Weyl elements exist for each root. If the root system $ \roots $ is not specified, we will also call these objects \defemph*{root gradings}, and the pair $ (G, (\rootgr{\alpha})_{\alpha \in \roots}) $ is called a \defemph*{$ \roots $-graded group} or \defemph*{root graded group}. For a precise definition, see~\ref{rgg-def}. The main result of this work is that every $ \roots $-graded group is coordinatised by some algebraic structure if $ \roots $ is of rank at least~3.

Our definition of $ \roots $-gradings makes sense for any finite root system, crystallographic or not, but all the central results of this book concern mainly the case of crystallographic root systems. We will briefly remark on the situation in the non-crystallographic case in the last part of this preface. Except for this detour, all root systems in the preface are assumed to be crystallographic.

This work is a revised version of the author's PhD thesis \cite{Wiedemann-PhD}. The main differences are in \cref{chap:pre,chap:chev}, which only contain known preliminaries and have been shortened, and \cref{chap:param}, where the long motivational section \cite[4.1]{Wiedemann-PhD} has been replaced by some additional, more concise remarks in the main part of the text.

\pdfbookmark[1]{Motivation and Examples: Chevalley Groups}{pre-chev}
\section*{Motivation and Examples: Chevalley Groups}

The main examples of root graded groups are the Chevalley groups, which were introduced by Chevalley in his famous Tohoku paper \cite{Chev-Tohoku}. By construction, every Chevalley group $ G $ is defined over a commutative associative unital ring~$ \ring $: It is a matrix group with coefficients in $ \ring $, and for each root group $ \rootgr{\alpha} $ there exists a canonical isomorphism $ \map{\risom{\alpha}}{(\ring, +)}{\rootgr{\alpha}}{}{} $. Further, the celebrated \defemph*{Chevalley commutator formula}\index{Chevalley commutator formula} asserts that
\[ \commutator{\risom{\alpha}(a)}{\risom{\beta}(b)} = \prod_{\substack{i,j \ge 1 \\ i\alpha + j\beta \in \roots}} \risom{i\alpha + j\beta}\brackets[\big]{c_{\alpha \beta ij} a^i b^j} \]
for all non-proportional roots $ \alpha, \beta $ and all $ a,b \in \ring $ where $ c_{\alpha \beta ij} $ are integral constants which do not depend on $ \ring $, $ a $ and $ b $. In other words, commutators in $ G $ are described by the ring multiplication. For any root $ \alpha $, an element $ w_\alpha \in \rootgr{-\alpha} \rootgr{\alpha} \rootgr{-\alpha} $ is called an \defemph*{$ \alpha $-Weyl element} if it is a \enquote{lift} of the reflection $ \refl{\alpha} $ in the Weyl group of $ \roots $ to the group $ G $. By this we mean that $ \rootgr{\beta}^{w_\alpha} = \rootgr{\refl{\alpha}(\beta)} $ for all roots $ \beta $. It turns out that an element $ w_\alpha \in G $ is an $ \alpha $-Weyl element if and only if it has the form $ w_\alpha = \risom{-\alpha}(-a^{-1}) \risom{\alpha}(a) \risom{-\alpha}(-a^{-1}) $ where $ a $ is an invertible element of $ \ring $. In particular, $ \alpha $-Weyl elements exist because there exists at least one invertible element in $ \ring $, namely $ 1_\ring $. We conclude that the group-theoretic structure of a Chevalley group $ G $ is intricately connected with the algebraic structure of the ring $ \ring $ over which it is defined.

On the other hand, the definition of root graded groups is purely combinatorial: it does not refer to any underlying algebraic structure. In fact, the axioms of a root grading can be seen as a list of the most important properties of Chevalley groups which can be formulated without reference to the underlying ring $ \ring $.

\pdfbookmark[1]{Goals}{pre:goals}
\section*{Goals}

\pdfbookmark[2]{The Coordinatisation Problem}{pre-coord-prob}
\subsection*{The Coordinatisation Problem}

It is now a natural question whether the combinatorial axioms of a root graded group are enough to construct (or, in the special case of Chevalley groups, reconstruct) a ring \enquote{over which it is defined}. To make this precise, we say that a root graded group $ (G, (\rootgr{\alpha})_{\alpha \in \roots}) $ is \defemph*{coordinatised by the ring $ \ring $} if there exist isomorphisms $ (\map{\risom{\alpha}}{(\ring, +)}{\rootgr{\alpha}}{}{})_{\alpha \in \roots} $ such that the Chevalley commutator formula holds. Thus the question is whether each root graded group is coordinatised by some commutative associative unital ring $ \ring $. While this question is a step in the right direction, it is too naively posed to admit a positive answer. We will explain why in the following.

First of all, the theory of root graded groups seems to be very unruly without additional assumptions in the lower-rank cases. For this reason, we will largely restrict ourselves to the case that $ \roots $ is of rank at least~3.

Secondly, we cannot hope for the coordinatising ring to be commutative in all cases because, for example, (a variation of) the construction of Chevalley groups of type $ A_n $ works for non-commutative rings. There are even constructions of $ A_2 $- and $ C_3 $-graded groups from so-called alternative rings, which form a class of nonassociative rings.

Thirdly, even if we allow a more general class of rings, we make the following observation: If $ G $ is coordinatised by a single ring $ \ring $, then in particular, all root groups must be pairwise isomorphic. However, we will see that the axioms of a root grading merely imply that $ \rootgr{\alpha} $ is isomorphic to $ \rootgr{\beta} $ if $ \alpha $ and $ \beta $ are conjugate under the Weyl group of $ \roots $. Thus in general, we have to allow the distinct orbits of root groups to be coordinatised by distinct algebraic structures. For example, we can construct $ B_n $-graded groups in which the long root groups are coordinatised by a commutative associative unital ring $ \comring $ and the short root groups are coordinatised by a quadratic $ \comring $-module $ (\module, q) $ (that is, a $ \comring $-module $ \module $ with a quadratic form $ \map{q}{\module}{\comring}{}{} $). The commutator formulas in these groups look exactly like the Chevalley commutator formula, except that terms $ a^2 $ for $ a \in \module $ are replaced by $ q(a) $ and terms $ 2ab $ for $ a,b \in \module $ are replaced by $ f(a,b) $ where $ f $ denotes the linearisation of $ q $.

In light of the previous paragraph, we make the following definition where $ \calO $ denotes the set of orbits in $ \roots $: A \defemph*{coordinatisation of $ (G, (\rootgr{\alpha})_{\alpha \in \roots}) $}\index{coordinatisation of a root graded group} consists of a family of (usually abelian) groups $ (M_O)_{O \in \calO} $ called the \defemph*{coordinatising groups}\index{coordinatising group} or \defemph*{coordinatising structures}\index{coordinatising structure}, a family of maps between the coordinatising groups $ (M_O)_{O \in \calO} $ which \enquote{equips $ (M_O)_{O \in \calO} $ with an algebraic structure} and a family of group isomorphisms $ \map{\risom{\alpha}}{M_O}{\rootgr{\alpha}}{}{} $ for $ O \in \calO $ and $ \alpha \in O $. We further require that a generalised version of the Chevalley commutator formula (involving the maps between $ (M_O)_{O \in \calO} $ and the isomorphisms $ (\risom{\alpha})_{\alpha \in \roots} $) is satisfied. For example, in the case of $ B_n $-graded groups, the coordinatising structures are abelian groups $ (\comring, +) $ and $ (\module, +) $ and we have a map $ \map{}{\comring \times \comring}{\comring}{}{} $ which turns $ \comring $ into a commutative associative unital ring, a map $ \map{}{\comring \times \module}{\module}{}{} $ which turns $ \module $ into a $ \comring $-module and a map $ \map{}{\module}{\comring}{}{} $ which is a quadratic form on $ \module $. The connection between Weyl elements in $ G $ and invertible elements in the algebraic structure that we have seen earlier for Chevalley groups remains valid in this more general setting, though of course the precise meaning of \enquote{invertibility} depends on the specific algebraic structure. (For example, an element $ v $ in a quadratic module $ (\module,q) $ is called \enquote{invertible} if $ q(v) $ is invertible, and its \enquote{inverse} is then defined to be $ q(v)^{-1}v $.) In particular, the existence of Weyl elements in $ G $ yields that there exists at least one \enquote{invertible} element in each algebraic structure that we encounter. For this reason, all rings in this book are unital.

With these refinements in mind, we can now formulate the main goal of this work: For each irreducible crystallographic root system $ \roots $ of rank at least~3, show that there exists a class $ \algclass{\roots} $ of (families of) algebraic objects such that each $ \roots $-graded group is coordinatised by an object in $ \algclass{\roots} $. This is the \defemph{coordinatisation problem}. We will discuss in a moment which results in this direction have already been obtained in the literature. In this book, we present a complete solution to the coordinatisation problem which does not rely on any previous work on abstract root gradings. We emphasise that our results are completely characteristic-free, whereas some previous works exclude the case of \enquote{characteristic~2}.

Apart from elementary undergraduate algebra, the list of mathematical prerequisites that we need is short: We will use the basic theory of finite reflection groups and some results on Chevalley groups, though the latter could be replaced by a sequence of long but easy computations in a purely combinatorial setting. Further, we will cite a few elementary results on nonassociative rings without proof. We will also encounter several interesting algebraic structures, such as quadratic modules, composition algebras and involutions on rings, but we will derive all their basic properties in this text.

\pdfbookmark[2]{Main Results}{pre-main-results}
\subsection*{Main Results}

We briefly summarise the specific coordinatisation results that we obtain for each root system. We will elaborate on the history of these results, some which can already be found in the literature in less general (for types $ A $, $ D $ and $ E $) or partial (for type $ C $) versions, in a later section of this preface. For all $ n \in \IN_{\ge 2} $, we show that each $ A_n $-graded group is coordinatised by a (possibly nonassociative noncommutative) unital ring which must be associative if $ n \ge 3 $ (\cref{ADE:thm}). Now let $ n \in \IN_{\ge 3} $. Root graded groups of types $ D_n $ or $ E_n $ are coordinatised by commutative unital rings (again \cref{ADE:thm}). Root graded groups of type $ B_n $ are coordinatised by pairs $ (\module, \comring) $ where $ \comring $ is a commutative associative unital ring and $ \module $ is a quadratic $ \comring $-module (\cref{B:thm}). Root graded groups of type $ F_4 $ are coordinatised by pairs $ (\compalg, \comring) $ where $ \comring $ is a commutative associative unital ring and $ \compalg $ is a multiplicative conic alternative algebra over $ \comring $ in the sense of \cite[Sections~16,~17]{GPR_AlbertRing} (\cref{F4:thm}). The latter objects should be regarded as a generalisation of composition algebras. If $ \comring $ is a field of characteristic not~2 and a certain norm function on $ \compalg $ is non-degenerate, then $ \compalg $ is a composition algebra in the classical sense.

The results for the root systems $ C_n $ and $ BC_n $ are more technical. We will show that every root graded group of type $ C_n $ or $ BC_n $ is coordinatised by a pair $ (\jormod, \ring) $ where $ \ring $ is an alternative ring with nuclear involution and $ \module $ is a Jordan module over $ (\ring, \rinvmap) $ which must be abelian if the root system is $ C_n $ (\cref{BC:thm}). The notion of Jordan modules is new and specifically tailored to be exactly the kind of algebraic structure which coordinatises $ (B)C_n $-graded groups. As important special cases, the class of abelian Jordan modules contains the class of involutory sets which are known from the classification of Moufang quadrangles while the class of non-abelian Jordan modules contains the groups $ T $ which are constructed in \cite[(11.24)]{MoufangPolygons} from pseudo-quadratic modules over $ \ring $. Here an involutory set over an alternative ring $ \ring $ with nuclear involution $ \rinvmap $ is a certain subgroup of the set of fixed points of $ \rinvmap $. If $ 2_\ring $ is invertible, then we can prove the following statement in a purely algebraic manner: Every abelian Jordan module over $ \ring $ is the direct product of an involutory set with an $ \ring $-module and every non-abelian Jordan module can be constructed from a pseudo-quadratic module over $ \ring $ as in \cite[(11.24)]{MoufangPolygons}. 

\pdfbookmark[2]{The Existence Problem}{pre-ex-prob}
\subsection*{The Existence Problem}

As a second goal, it is of course desirable to show that the class $ \algclass{\roots} $ in our solution of the coordinatisation problem is chosen optimally. In other words, for every object $ X $ in $ \algclass{\roots} $ there should exist a $ \roots $-graded group which is coordinatised by $ X $. This is the \defemph{existence problem} or \defemph{construction problem}. In most cases, the objects in $ \algclass{\roots} $ are related to associative rings and we can easily construct $ \roots $-graded groups from these objects by writing down some (generalised) matrices. However, the classes $ \algclass{\roots} $ for $ \roots \in \Set{C_3, BC_3, F_4} $ involve alternative rings. For the subclass of $ \algclass{\roots} $ consisting only of the \enquote{associative objects}, a construction in terms of matrices is usually still possible, but the general construction problem is much more difficult.

For $ F_4 $-graded groups, no general construction is known. However, if the desired coordinatising ring $ \ring $ (which is equipped with the structure of a multiplicative conic algebra) is not only associative but also commutative, then we know that there exists an $ E_6 $-graded group $ G $ which is coordinatised by $ \ring $ (for example, a Chevalley group of type $ E_6 $). Using the general mechanism of \defemph*{foldings of root graded groups}, we obtain an $ F_4 $-graded group which is coordinatised by the conic algebra $ \ring $. This provides a partial solution of the existence problem for this special case.

Now consider $ BC_3 $-graded groups. There are strong indications that, while alternative rings $ \ring $ can technically appear in this setting, only \enquote{the associative part of $ \ring $ is relevant} and thus we can restrict ourselves to the case of associative rings with good conscience. In this situation, we provide a construction which works for all objects in $ \algclass{BC_3} $ for which $ 2 $ is invertible (and also for a large subclass of objects in which $ 2 $ is not invertible).

For $ C_3 $-graded groups, Zhang provides a strategy to solve the existence problem for a certain subclass of $ \algclass{C_3} $ in \cite[Section~4.2]{Zhang}. We are confident that his construction can be adapted to the general case, and we plan to address this in future work.

The class $ \algclass{A_2} $ involves non-associative rings as well. By a construction of Faulkner (see \cite[Section~3]{Faulkner-StableRange} and the appendix of \cite{Faulkner-Barb}), this class is known to contain all alternative rings, but it is unclear whether it contains more general rings. Just like Zhang's construction of $ C_3 $-graded groups, Faulkner's construction realises the desired group as a group of automorphism of the Tits-Kantor-Koecher algebra of a certain Jordan pair.

\pdfbookmark[1]{Root gradings, RGD-systems and Moufang Buildings}{pre-rgd}
\section*{Root gradings, RGD-systems and Moufang Buildings}

We now turn to the history of the subject.
As already said, Chevalley groups were introduced \cite{Chev-Tohoku}. It should be noted, however, that the list of Chevalley groups contains several classical groups which had been studied before, such as the special linear group over a field. These are the earliest examples of root graded groups that have been considered, albeit without the abstract combinatorial framework of root gradings.

The first and most important example of an abstract structure that is similar to (but less general than) root gradings is Tits' notion of \defemph*{RGD-systems}\index{RGD-system}\footnote{The letters \enquote{RGD} stand for \enquote{Root Groups Data}. Tits remarks in \cite[\onpage{258}]{Tits-TwinBuildingsKacMoody} that he chose this name \enquote{for lack of a better idea (or rather, because all appropriate names that I can think of seem to be already taken!)}. It goes without saying that the same could be said about the name \enquote{root grading}, which carries essentially the same information as \enquote{root group data}.}. These objects were introduced under this name in \cite[\onpage{258}]{Tits-TwinBuildingsKacMoody}, but predecessors of the RGD-axioms can be found as early as in \cite[\onpage{140}]{Tits-GrpSemiSimpIsotrop}. A definition which is very similar to the one in \cite{Tits-TwinBuildingsKacMoody} is the one of \defemph*{données radicielles}\index{donnée radicielle} in \cite[(6.1.1)]{BruhatTits1}. Further, RGD-systems are essentially the same thing as the \defemph*{groups with Steinberg relations}\index{Steinberg relations} in \cite{Faulkner-RGD}.

The crucial property that differentiates RGD-systems from root gradings is that in an RGD-system, \emph{every} non-trivial root group element $ b_\alpha \in \rootgr{\alpha} \setminus \compactSet{1_G} $ can be \enquote{extended} to a Weyl element $ w_\alpha = a_{-\alpha} b_\alpha c_{-\alpha} $. In a root grading, we only require that \emph{some} element $ b_\alpha \in \rootgr{\alpha} \setminus \compactSet{1_G} $ has this property. Thus by the correspondence between Weyl elements and invertible elements, any RGD-system that is coordinatised by an algebraic structure must be coordinatised by a \enquote{division structure}. In fact, the coordinatisation problem for RGD-systems is already solved by the work of Tits-Weiss in \cite{MoufangPolygons}, but this fact is usually not stated in the group-theoretic language of RGD-systems. In order to explain this, we first have to make a brief detour through the theory of spherical buildings.

Spherical buildings were introduced by Tits \enquote{in an attempt to give a systematic procedure for the geometric interpretation of the semi-simple Lie groups and, in particular, the exceptional groups} (quotation from \cite[\onpage{V}]{Tits-LectureNotes74}). More generally, they even provide a geometric framework to study simple algebraic groups over arbitrary fields. There are several different ways to define spherical buildings, but the technical details are not relevant for our purposes. Each building has an associated type which is a Coxeter system (or equivalently, a root system). The rank of a spherical building is by definition the rank of its type, and it is called irreducible if its type is irreducible. There also exist buildings which are not spherical (which means that their type is not spherical), but they are not relevant in this context. Further, there exist a notion of thickness of buildings and of the so-called \defemph*{Moufang property}\index{Moufang property}.\footnote{This condition is named after Ruth Moufang, who is known for her work on the class of projective planes which today are called \defemph*{Moufang planes}\index{Moufang plane}. In fact, a projective plane is the same thing as a thick spherical building of type $ A_2 $, and a Moufang plane is a thick spherical building of type $ A_2 $ with the Moufang property.} A spherical building which satisfies the Moufang property is called a Moufang building. It is a highly non-trivial fact, proven in \cite[3.5]{Tits-EndSpiegWeyl}, that every irreducible spherical building of rank at least 3 is automatically Moufang. The notion of \defemph*{generalised polygons}\index{generalised polygon} is equivalent to that of spherical buildings of rank~2, and generalised polygons which satisfy the Moufang property are called \defemph*{Moufang polygons}.

Thick irreducible spherical buildings of rank at least~3 have been classified by Tits in \cite{Tits-LectureNotes74}, and this classification has later been extended to Moufang polygons by Tits-Weiss in \cite{MoufangPolygons}. In fact, the classification of Moufang polygons is completely independent of \cite{Tits-LectureNotes74}, and the classification of higher-rank spherical buildings can be deduced from the classification of Moufang polygons (see \cite[Chapter~40]{MoufangPolygons}). The statement of these classification results is that each thick irreducible spherical Moufang building of rank at least~2 is \enquote{coordinatised} by some algebraic division structure.

The interest in Moufang buildings for the theory of root graded groups stems from the following fact: Modulo technical details, the geometrical notion of thick Moufang buildings of type $ \roots $ is essentially equivalent to the group-theoretic notion of RGD-systems of type $ \roots $ (see \cite[Section~7.8]{AbramenkoBrown-Buildings} for details). More precisely, the Moufang condition ensures the existence of an RGD-system in the automorphism group of such a spherical building, and a spherical building can be constructed abstractly from an RGD-system. Thus the classification results for Moufang buildings that we have cited above provide a complete solution of the coordinatisation problem for RGD-systems of rank at least~3. In this sense, the main results of this book can be seen as a generalisation of the classification of thick irreducible spherical buildings of rank at least~3 to the case of \enquote{non-division structures}. It is noteworthy that the arguments in \cite{MoufangPolygons} are very similar in style to our work: At first Tits-Weiss construct a \defemph*{root group sequence} from any Moufang polygon, which is essentially its RGD-systems, and then they classify root group sequences. The proofs in \cite{Tits-LectureNotes74}, on the other hand, have a more geometric flavour.

\pdfbookmark[1]{Known Results about Root Gradings}{pre-known}
\section*{Known Results about Root Gradings}

We now turn to root gradings which are not necessarily RGD-systems. The terminology of \enquote{$ \roots $-graded groups} appears for the first time in the paper \cite{Shi1993} by Shi. As in our definition, Shi only requires the existence of Weyl elements and not the stronger \enquote{division axiom} that appears in the definition of RGD-systems. In fact, Shi's definition of $ \roots $-graded groups $ (G, (\rootgr{\alpha})_{\alpha \in \roots}) $ is essentially equivalent to our definition except for one additional axiom: that $ G $ contains a homomorphic image of the Steinberg group\index{Steinberg group} of type $ \roots $ over some commutative associative unital ring. Here the Steinberg group is a certain group defined by generators and relations which mimic the commutator relations in Chevalley groups. We will explain the motivation behind this axiom in the following paragraph. One of the main results of \cite{Shi1993} is the solution of the coordinatisation problem for all root graded groups (in Shi's more restrictive sense) of type $ A $, $ D $ or $ E $ and of rank at least~3. That is, Shi coordinatises root graded groups for all simply-laced types of rank at least~3. He arrives at the same conclusion that we do with our more general definition of root gradings. In fact, a non-trivial corollary of our coordinatisation result is that any root graded group in our sense is also a root graded group in Shi's sense.

Shi introduced his notion of root graded groups as a group-theoretic analogue of \defemph*{root graded Lie algebras}\index{root graded Lie algebra}\index{Lie algebra!root graded}. These Lie algebras were introduced in \cite{RGLie-BermanMoody} by Berman and Moody, of whom Shi is a student. A complete classification of root graded Lie algebras is available by the combined work of several authors, see \cite{RGLie-BermanMoody,RGLie-2laced,RGLie-3graded,RGLie-CentExt,RGLie-BCr,RGLie-BC1,RGLie-Prime}. It plays an important role in the classification of semisimple Lie algebras over algebraically closed fields of characteristic at least~5. One of the axioms of a $ \roots $-graded Lie algebra $ L $ is that it contains a subalgebra $ L' $ which is split simple with root system $ \roots $. Clearly this axiom motivated the additional axiom in Shi's definition which we discussed in the previous paragraph. A consequence of this additional axiom is that a certain \defemph*{sign problem}\index{sign problem} which appears in the coordinatisation of root graded groups can be solved very efficiently in Shi's setting. This problem is related to a similar problem about the signs appearing in the Chevalley commutator formula. In order to solve this problem in our situation, we introduce a machinery which is independent of the root system $ \roots $ and which builds upon the solution of the word problem in Coxeter groups. The main result of this machinery is the \defemph{parametrisation theorem} that we will discuss in a moment. See also \cref{param:motiv:parmap-choice} for more details on the sign problem.

The only work known to us which considers root gradings for root systems which are not simply-laced is the PhD thesis \cite{Zhang} of Zhang, a student of Zelmanov. Using essentially the same definition for root gradings as Shi, Zhang proves some partial results concerning $ C_3 $-gradings. Namely, he shows that every $ C_3 $-graded group satisfying some additional conditions is coordinatised by a nonassociative ring $ \ring $ with involution. One of these additional conditions on the group implies that $ 2_\ring $ is invertible in $ \ring $, so it essentially excludes the case of \enquote{characteristic 2}. The general case is significantly more difficult, so this is a serious restriction. Further, Zhang constructs a $ C_3 $-graded group from every alternative ring with involution in which $ 2 $ is invertible. The assumption on the alternativity of the ring is necessary because, as we will show in the main part of this work, the ring which appears in the coordinatisation of $ C_3 $-graded groups must actually be alternative. We prove this using a new method which we call the \defemph{blueprint technique}.

Finally, there are contributions of Faulkner on $ A_2 $-gradings that are important from the historical perspective. In \cite[Section~13.3]{Faulkner-NonAssocProj}, he defines \defemph*{groups of Steinberg type}\index{Steinberg type}, which are essentially the same thing as $ A_2 $-graded groups in our sense. Faulkner proceeds to show that every such group is coordinatised by a nonassociative ring. It is unlikely that arbitrary nonassociative rings can appear here, so this solution of the coordinatisation problem is incomplete. The most general known construction of $ A_2 $-graded groups, given in \cite{Faulkner-StableRange} and in the appendix of \cite{Faulkner-Barb}, starts from an alternative ring.

We can thus summarise the previous state of the literature for root gradings of rank at least 3 as follows. For RGD-systems, the coordinatisation and existence problems are completely solved. For root gradings, all known coordinatisation results use Shi's more restrictive definition of root gradings. With this caveat, the coordinatisation problem for the simply-laced case is completely solved by \cite{Shi1993} and the existence problem is easy in this case. For $ C_n $-gradings, there exist partial results by \cite{Zhang}. For gradings of type $ B $, $ BC $ and $ F_4 $, there are no previous results.

Another noteworthy work in the context of root gradings is \cite{LoosNeherBook}, although its interest does not lie in the coordinatisation and existence problems in our sense. Still, Loos-Neher define \defemph*{groups with $ \roots $-commutator relations} (where $ \roots $ can be a root system, but also any subset of a free abelian group satisfying some axioms), citing \cite{Faulkner-RGD} as their inspiration (see \cite[\onpage{72}]{LoosNeherBook}). Further, they also define the notion of Weyl elements. A $ \roots $-graded group in our sense is essentially the same thing as a group with $ \roots $-commutator relations in which Weyl elements exist and which satisfies a non-degeneracy condition.

\pdfbookmark[1]{A Uniform Approach to Root Graded Groups}{pre-uniform}
\section*{A Uniform Approach to Root Graded Groups}

Since the class $ \algclass{\roots} $ of coordinatising algebraic structures depends on the root system $ \roots $, it is clear that a certain amount of case-by-case analysis is necessary in the study of root graded groups. However, we will encounter two root-system-independent tools which form the cornerstones of our work, and which are both new: the \defemph{parametrisation theorem} and the \defemph{blueprint technique}. The parametrisation theorem is our solution to a very delicate and technical sign problem, and it streamlines the coordinatisation process of root graded groups. The blueprint technique lies at the heart of our work, and it is pivotal for the determination of the class $ \algclass{\roots} $ for each root system $ \roots $. In particular, a characteristic-free solution to the coordinatisation problem would not have been possible without the blueprint technique. Together, these two tools provide a conceptual approach to root gradings which is as uniform as is possible.

\pdfbookmark[2]{The Parametrisation Theorem}{pre-param}
\phantomsection
\subsection*{The Parametrisation Theorem}

\label{subsec:preface:param-coord}
Before we delve into the technicalities of the parametrisation theorem, we should clarify our conventions regarding the notions of \enquote{coordinatisations} and \enquote{para\-me\-tri\-sations}. Let $ G $ be a $ \roots $-graded group and denote the set of all orbits in $ \roots $ under the Weyl group by $ \calO $. Under a \defemph*{parametrisation of $ G $}\index{parametrisation of a root graded group}, we understand a family $ (M_O)_{O \in \calO} $ of groups and for each $ O \in \calO $ a family $ (\map{\risom{\alpha}}{M_O}{\rootgr{\alpha}}{}{})_{\alpha \in O} $ of isomorphisms satisfying a certain consistency condition. Thus a \defemph*{coordinatisation of $ G $}\index{coordinatisation of a root graded group} in the sense that we have defined earlier is simply a parametrisation of $ G $ together with a family of maps between the parametrising groups satisfying two conditions: Firstly, these maps describe the commutator relations between the root groups, and secondly, they equip the parametrising groups with some algebraic structure. Hence essentially, a coordinatisation is a parametrisation by an algebraic structure, and we will use this distinction between \enquote{coordinatisations} and \enquote{parametrisations} throughout this book.

Our approach to the coordinatisation problem consists of first finding a parametrisation of $ G $ and then turning this parametrisation into a coordinatisation. The parametrisation theorem is the one and only criterion that we will use to establish the existence of a parametrisation of $ G $. The main difficulty in the proof of the parametrisation theorem is a certain sign problem which is related to the consistency condition in the definition of parametrisations. This problem could be solved efficiently if the additional axiom in Shi's definition of root gradings were satisfied. In this sense, the parametrisation theorem is the solution of the sign problem in our more general setting.

The requirements of the parametrisation theorem say that, in order to find a parametrisation, we have to do three things.
Firstly, we have to understand the action of $ w_\delta^2 $ on each root group where $ \delta $ is an arbitrary simple root (with respect to some fixed root base $ \rootbase $) and $ w_\delta $ is an arbitrary $ \delta $-Weyl element. Secondly, we have to show that Weyl elements satisfy the \defemph{braid relations}. These first two properties correspond to the homotopy moves which appear in the solution of the word problem in Coxeter groups, and the parametrisation theorem is proven by performing the same kind of moves on the level of Weyl elements. It should be noted that the braid relations for Weyl elements can be verified in a uniform way for all root systems, using the beautiful argument of Tits-Weiss in \cite[(6.9)]{MoufangPolygons}. Thirdly, we have to a find a suitable  \enquote{system of signs} $ \inverparsym $ which is \enquote{consistent}. Technically, $ \inverparsym $ is a map from $ \roots \times \rootbase $ to a finite abelian group $ \twistgroup $ of exponent $ 2 $ which acts on each root group, and \enquote{consistency} refers to a set of combinatorial properties that it should satisfy.

Observe that the system of signs $ \inverparsym $ in the previous paragraph is a map between finite sets. Thus the necessary consistency properties could, in theory, be checked by a long sequence of easy computations. However, there exists a more elegant way of doing this: By verifying that $ \inverparsym $ appears as the system of signs in \emph{some} example of a \enquote{sufficiently generic}\index{generic example group} root graded group. For the simply-laced root systems, the Chevalley groups are sufficient for this purpose. In the general case, a complete solution of the existence problem always produces a sufficiently generic example as a by-product. Even in the cases where no complete solution of the existence problem is available, we are still able to construct a group which is \enquote{sufficiently generic} to produce a suitable map $ \inverparsym $.

After these three problems are tackled, we can apply the parametrisation theorem to obtain a parametrisation of $ G $. In order to turn it into a coordinatisation, we have to equip the parametrising groups with an algebraic structure. This is where the blueprint technique comes into play.

\pdfbookmark[2]{The Blueprint Technique}{pre-blue}
\subsection*{The Blueprint Technique}

The blueprint technique is a powerful tool which reduces the problem of determining the commutator relations in root graded groups with a parametrisation to a straightforward computation involving certain rewriting rules. It is inspired by the work of Ronan-Tits on the construction of buildings in \cite{BuildBuildings} and the key to our char\-ac\-teristic-free approach to root gradings.

The basic idea of the blueprint technique is as follows. We begin with a reduced representation $ f $ of the longest word $ \rho $ in the Weyl group $ W $ of $ \roots $. Then there exists a sequence $ h_1, \ldots, h_r $ of elementary homotopy moves which transforms $ f $ into itself and which, speaking very loosely in terms of algebraic topology, \enquote{moves around the hole in the Cayley graph of $ W $}. For example, the Cayley graph of the Weyl group of $ A_3 $ can be drawn on the 2-sphere with the identity element $ 1_W $ at the north pole and the longest element $ \rho $ at the south pole, and there is a way to move $ f $ (which is a path from $ 1_W $ to $ \rho $) once around the sphere. Further, we take a word $ \tilde{f} $ \enquote{corresponding to $ f $} whose letters are arbitrary elements (in some sense, \enquote{indeterminates}) of the parametrising groups of $ G $.

The blueprint technique associates each elementary homotopy move $ h_i $ on $ f $ to a rewriting rule $ \tilde{h}_i $ on $ \tilde{f} $. This rule $ \tilde{h}_i $ involves the maps appearing in the commutator relations of $ G $. Now the \defemph{blueprint computation} consists of iteratively determining all words which are obtained from $ \tilde{f} $ by applying the rewriting rules $ \tilde{h}_1, \ldots, \tilde{h}_r $. This computation is purely mechanical, albeit lenghty, and perfectly suited to be executed on a computer. For abstract reasons, the end result of the blueprint computation must be the same as the initial word $ \tilde{f} $, so we obtain a sequence of identities (one for each letter in $ \tilde{f} $) which are valid in the parametrising structures. These identities allow us to explicitly compute the commutator relations in $ G $, and the maps involved in these formulas equip the parametrising groups of $ G $ with an algebraic structure.

\pdfbookmark[1]{Organisation of the Book}{pre-orga}
\section*{Organisation of the Book}

In \cref{chap:pre}, we set up some general notation and recall basic preliminaries from the theory of finite root systems. In \cref{chap:rgg}, we introduce the language of root graded groups and prove some general results: the bijectivity of certain product maps and the braid relations for Weyl elements. The definition of root graded groups is given in~\ref{rgg-def}. \Cref{chap:chev} is devoted to a brief summary of the theory of Chevalley groups. We prove the parametrisation theorem in \cref{chap:param}. Afterwards, we begin our investigation of $ \roots $-graded groups for all root systems $ \roots $ of rank at least 3. As explained in the previous section, the general strategy is independent of $ \roots $: We first verify that the requirements of the parametrisation theorem are satisfied, then we apply it and finally, we use the blueprint technique to compute the commutator relations. The only exception to this rule are the simply-laced root systems, which are covered in \cref{chap:simply-laced}, because the commutator relations in this case are so simple that we do not need to use the blueprint technique. In \cref{chap:blue}, we formally introduce the blueprint technique. Afterwards, we study root graded groups of type $ B $ in \cref{chap:B}, of type $ BC $ (which contains type $ C $ as a special case) in \cref{chap:BC} and of type $ F_4 $ in \cref{chap:F}. \Cref{chap:BC-alg} is a prelude to \cref{chap:BC} which introduces Jordan modules, the algebraic structures that coordinatise root gradings of type $ BC $.

\pdfbookmark[1]{Remarks on Related Topics}{pre-related}
\section*{Remarks on Related Topics}

We now turn back to geometrical aspects of the theory.
Recall that Moufang polygons are geometric objects which correspond to RGD-systems of rank 2. In \cite{MW-TitsPolygons}, Mühlherr-Weiss have introduced \defemph*{Tits polygons}, which generalise Moufang polygons. In \cite{MW-RGG2}, they show that Tits polygons correspond bijectively to root gradings of rank 2 which satisfy the so-called \defemph{stability condition}. On the algebraic side of coordinate systems, this condition corresponds to the notion of \defemph*{rings of stable rank 2}\index{stable rank 2}, sometimes simply called \defemph*{stable rings}.

The concept of stability originally stems from $ K $-theory, but it was realised by Veldkamp in \cite{Veldkamp-StabProj} that stable rings are precisely those rings over which projective planes can be defined in a meaningful way. He gives combinatorial axioms for such planes, which today are called \defemph*{Veldkamp planes}\index{Veldkamp plane}, and proceeds to show that a Veldkamp plane is the projective plane of a stable ring if and only if it is Desarguesian. See also \cite{Veldkamp-GeomRings} for additional information.

In \cite[2.8]{MW-VeldPol}, Mühlherr-Weiss define the notion of \defemph*{Veldkamp $ n $-gons}, or simply \defemph*{Veldkamp polygons}\index{Veldkamp polygon}. A Veldkamp plane is exactly the same thing as a Veldkamp triangle (that is, a Veldkamp 3-gon). In this language, a Tits polygon is a Veldkamp polygon which satisfies the Moufang condition. Higher-rank analogues of Veldkamp polygons can be defined, which leads to the notion of \defemph*{Veldkamp buildings (of rank at least~3)}\index{Veldkamp building}. It seems likely that every Veldkamp building of rank at least 3 has an associated stable root graded group (of rank at least~3), and thus the results of this book yield a classification of (or at least strong coordinatisation results for) Veldkamp buildings of rank at least 3. In fact, this relation to Veldkamp buildings is a strong motivation for our work. The absence of stability conditions in this book is simply due to the fact that they are not necessary to build a meaningful theory on the group-theoretic side, at least in the higher-rank situation.

It is possible to define an analogue of Veldkamp buildings without the stability condition, the so-called \defemph*{Faulkner buildings}\index{Faulkner building}. Since the stability condition has proven crucial in the context of projective geometry, it is not at all clear whether Faulkner buildings actually describe a meaningful geometry. However, they are the natural class of spaces which correspond to general root graded groups.

We end with a few words on the non-crystallographic root systems $ H_3 $ and $ H_4 $. It is known from \cite[Hauptsatz]{Tits-EndSpiegWeyl} that there exist no thick buildings (and thus no RGD-systems) of these types, so one might expect the same to be true for root gradings. However, this is wrong: Starting from a $ D_6 $-graded group $ (G, (\rootgr{\alpha})_{\alpha \in D_6}) $, one can fold $ (\rootgr{\alpha})_{\alpha \in D_6} $ to obtain an $ H_3 $-grading $ (V_\beta)_{\beta \in H_3} $ of the same group $ G $ in which each root group $ V_\beta $ equals the product of two commuting root groups in $ (\rootgr{\alpha})_{\alpha \in D_6} $. In particular, if $ (\rootgr{\alpha})_{\alpha \in D_6} $ is coordinatised by a commutative associative ring $ \comring $, then $ (V_\beta)_{\beta \in H_3} $ is coordinatised by the ring $ \comring \times \comring $. A similar construction for $ H_4 $-graded groups is possible. In collaboration with Lennart Berg, we have shown in the joint and yet unpublished work\mywarning{check: H3} \cite{RGG-H3} that every $ \roots $-graded group for $ \roots \in \Set{H_3, H_4} $ is of this form. Incidentally, this provides a new proof of the known fact that RGD-systems of these types cannot exist: The ring $ \comring \times \comring $ is never a field, even when $ \comring $ is a field.

	\chapter{Preliminaries}
	
	\label{chap:pre}
	
	In this chapter, we will recall some standard definitions and facts which will be needed throughout this book. None of the material in this chapter is new. In \cref{sec:commrel}, we set up some general notation. In \cref{sec:rootsys}, we summarise the main results from the theory of finite (non-reduced, non-crystallographic) root systems and their Weyl groups. In \cref{sec:fold}, we briefly discuss foldings of root systems, which will only be needed to construct examples of $ F_4 $-graded groups in \cref{sec:F4:const}. The reader may want to skip \cref{sec:rootsys,sec:fold} and only refer back to them as needed.
	

\section{Elementary Notation and Group-theoretic Facts}

\label{sec:commrel}

\begin{notation}
	We denote the sets of natural numbers (including $ 0 $), positive integers, integers, real numbers and complex numbers by $ \Nzero $, $ \Npos $, $ \IZ $, $ \IR $ and $ \IC $, respectively.
\end{notation}

We will often consider sets of integers of the form $ \Set{n, \ldots, m} $ where $ n \le m $. Since intervals of real numbers are completely absent from this book, the following notation is not ambiguous.

\begin{notation}[Integer intervals]
	For all $ n,m \in \IZ $, we set
	\[ \numint{n}{m} \defl \Set{i \in \IZ \given n \le i \le m} = \begin{cases}
		\emptyset & \text{if } n>m, \\
		\Set{n, \ldots, m} & \text{if } n \le m.
	\end{cases} \]
	Further, we put
	\[ \betint{n}{m} \defl \begin{cases}
		\numint{n}{m} & \text{if } n \le m, \\
		\numint{m}{n} & \text{if } n > m.
	\end{cases} \]
\end{notation}

Note that the statement \enquote{$ i \in \betint{n}{m} $} can be interpreted as \enquote{$ i $ lies between $ n $ and $ m $}, independent of whether $ n<m $ or $ m<n $.

As we have already seen in the preface, we have to adopt a very general definition of rings in this work.

\begin{convention}\label{pre:ring-nonassoc-conv}
	A ring is always understood to be unital, but the multiplication is not assumed to be associative or commutative. All the necessary background on these objects will be introduced in \cref{sec:ring}. However, it suffices to know only the definition of these objects until we reach \cref{sec:ring}.
\end{convention}

\begin{notation}[Matrices]
	For any associative ring $ \ring $ and for any positive integers $ n,m $, we denote by $ M_{nm}(\ring) $ or by $ \ring^{n \times m} $ the set of matrices with $ n $ rows and $ m $ columns. Further, we put $ M_n(\ring) \defl M_{nn}(\ring) $.\index{MnR@$ M_n(\ring) $}\index{MnmR@$ M_{nm}(\ring) $}
\end{notation}

Note the assumption on the ring $ \ring $ to be associative in the definition of matrices. We could, of course, define matrices over arbitrary rings in the same way, but then the matrix multiplication is not necessarily associative. In particular, we cannot construct groups from such matrices, which makes them useless for our purposes.

\begin{definition}[Words]\label{pre:word-def}
	Let $ S $ be a set. A \defemph*{word over $ S $}\index{word} is a tuple $ \tup{s}{k} $ for some $ k \in \Nzero $ such that $ \listing{s}{k} $ lie in $ S $. The \defemph*{free monoid over $ S $}\index{free monoid} is the set $ \frmon{S} $ of words over $ S $ together with the multiplication $ \tup{s}{k} \cdot \tup{t}{l} \defl (\listing{s}{k}, \listing{t}{l}) $. Its neutral element, the empty tuple, is called the \defemph{empty word} and denoted by $ \emptyset $. When there is no danger of confusion, we will also denote a word $ \tup{s}{k} $ by $ s_1 \cdots s_k $. The \defemph*{free group over $ S $}\index{free group} is denoted by $ \frgrp{S} $.
\end{definition}

\begin{notation}\label{braidword-def}
	Let $ M $ be a monoid, let $ m \in \Nzero $ and let $ x,y \in M $. We denote by $ \braidword{m}(x,y) $\index{Pm(x,y)@$ \braidword{m}(x,y) $} the word $ \prod_{i=1}^m z_i $ where $ z_i \defl x $ for all odd $ i \in \numint{1}{m} $ and $ z_i \defl y $ for all even $ i \in \numint{1}{m} $. In particular, $ \braidword{0}(x,y) = 1_M $. If $ x $, $ y $ lie in a set $ S $ which is not equipped with a multiplication, we denote by $ \braidword{m}(x,y) $ the word $ \braidword{m}\brackets[\big]{(x), (y)} = (x, y, x, \ldots) $ in the free monoid over $ S $.
\end{notation}

\begin{notation}[Group-theoretic notions]
	Let $ G $ be a group. For any subset $ U $ of $ G $, we denote by $ \gen{U} $ the subgroup of $ G $ which is \defemph*{generated by~$ U $}\index{generated subgroup}. For all $ g,h \in G $, the \defemph*{commutator of $ g $ and $ h $}\index{commutator} is $ \commutator{g}{h} \defl g^{-1} h^{-1} gh $ and the \defemph*{conjugate of $ g $ by $ h $}\index{conjugation} is $ g^h \defl h^{-1} gh $.
\end{notation}

\begin{definition}\label{act-by}
	Let $ G $ be a group and let $ H $ be a subset of $ G $. For any $ g \in G $, we say that \defemph*{$ g $ acts on $ H $ by inversion} if $ h^g = h^{-1} $ for all $ h \in H $ (where $ h^{-1} $ is not assumed to lie in $ H $) and we say that \defemph*{$ g $ acts trivially on $ H $} if $ h^g = h $ for all $ h \in H $. Further, for any $ a,b \in G $, we say that \defemph*{$ a $ and $ b $ act identically on $ H $} if $ h^a = h^b $ for all $ h \in G $.
\end{definition}

We end this section with some relations which hold in arbitrary groups. The relations in~\ref{group-rel} form the backbone of many computations, and we will often use them without specifically saying so.

\begin{relations}\label{group-rel}
	Let $ G $ be a group and let $ g, g_1, g_2, h, h_1, h_2 \in G $. Then the following relations hold:
	\begin{stenumerate}
		\item \label{group-rel:conj1}$ g^h = g\commutator{g}{h} $.
		\item \label{group-rel:conj2}$ g^h = \commutator{h}{g^{-1}} g $.
		\item $ \commutator{g_1}{g_2}^h = \commutator{g_1^h}{g_2^h} $ and $ (g_1 g_2)^h = g_1^h g_2^h $.
		\item \label{group-rel:conj-comm}$ gh = hg^h $.		
		\item \label{group-rel:comm}$ gh = hg \commutator{g}{h} = \commutator{g^{-1}}{h^{-1}} hg $.
		\item \label{group-rel:inv}$ \commutator{g}{h}^{-1} = \commutator{h}{g} $.
		\item \label{group-rel:add}$ \commutator{g_1 g_2}{h} = \commutator{g_1}{h}^{g_2} \commutator{g_2}{h} $ and $ \commutator{g}{h_1 h_2} = \commutator{g}{h_2} \commutator{g}{h_1}^{h_2} $.
		
		\item \label{group-rel:add-cor}$ \commutator{g_1g_2}{h} =  \commutator{g_2}{h} $ if $ \commutator{g_1}{h} = 1_G $ and $ \commutator{g}{h_1 h_2} = \commutator{g}{h_2} $ if $ \commutator{g}{h_1} = 1_G $.
	\end{stenumerate}
	Further, note that $ \commutator{g}{h} = 1_G $ is equivalent to each of the following conditions: $ gh = hg $, $ g^h = g $, $ \commutator{h}{g} = 1_G $, $ \commutator{g^{-1}}{h} = 1_G $ and $ \commutator{g}{h^{-1}} = 1_G $.
\end{relations}

Another relation which holds in arbitrary groups is the Hall-Witt identity. As a special case, it says that under certain conditions on the group elements, the commutator map is associative. It is used in this specific form in \cite[(2.23)]{Shi1993} to show that for any simply-laced root system $ \roots $ of rank at least~3, any ring which coordinatises a $ \roots $-graded group is associative. Further, the general form of the Hall-Witt identity given in \cref{hall-witt} plays an important role in the computation of the commutator relations in $ C_n $-graded groups (for $ n \ge 3 $) in \cite[3.4.15, 3.4.16, 3.4.18]{Zhang}. Thus the Hall-Witt identity holds an important place in the theory of root graded groups. It is therefore surprising that no important proofs in this book actually rely on the Hall-Witt identity. This is made possible by the blueprint technique, a new computational tool which we will introduce in \cref{chap:blue}.

\begin{lemma}[Hall-Witt identity]\index{Hall-Witt identity}\label{hall-witt}
	Let $ G $ be a group. Then for all $ x,y,z \in G $, we have
	\[ \commutator[\big]{\commutator{x}{y}}{z^x} \commutator[\big]{\commutator{z}{x}}{y^z} \commutator[\big]{\commutator{y}{z}}{x^y} = 1_G. \]
\end{lemma}
\begin{proof}
	This is a straightforward computation, involving only the definition of commutators and conjugates.
\end{proof}

\begin{lemma}[Special case of the Hall-Witt identity]\label{hall-witt-special}
	Let $ G $ be a group and let $ x,z,y \in G $ such that $ \commutator{x}{z} = 1_G $ and $ \commutator[\big]{\commutator{y}{x^{-1}}}{\commutator{y}{z}} = 1_G $. Then
	\[ \commutator[\big]{\commutator{x}{y}}{z} = \commutator[\big]{x}{\commutator{y}{z}}. \]
\end{lemma}
\begin{proof}
	By the Hall-Witt identity, we have $ \commutator[\big]{\commutator{x}{y}}{z} \commutator[\big]{\commutator{y}{z}}{x^y} = 1_G $. Hence
	\begin{align*}
		\commutator[\big]{\commutator{x}{y}}{z} &= \commutator[\big]{\commutator{y}{z}}{x^y}^{-1} = \commutator[\big]{x^y}{\commutator{y}{z}} = \commutator[\big]{\commutator{y}{x^{-1}} x}{\commutator{y}{z}} \\
		&= \commutator[\big]{\commutator{y}{x^{-1}}}{\commutator{y}{z}}^x \commutator[\big]{x}{\commutator{y}{z}} = \commutator[\big]{x}{\commutator{y}{z}},
	\end{align*}
	by \thmitemcref{group-rel}{group-rel:conj2} and~\thmitemref{group-rel}{group-rel:add}.
\end{proof}


\section{Root Systems}

\label{sec:rootsys}

We assume that the reader is familiar with the theory of finite root systems. However, since we want to include both the non-reduced root systems of type $ BC $ as well as the non-crystallographic root systems $ H_3 $ and $ H_4 $ in our theory of root graded groups (the latter only being studied in \cite{RGG-H3}), we need a definition of root systems that is more general than in most standard references. Further, in the setting on non-reduced root systems, we will consider so-called rescaled root bases in place of the more restrictive usual notion of root bases. Most classic results remain valid in this generality, but slight modifications of the statements are necessary in some cases. For this reason, we give a brief survey of all results that we will need throughout this work. In cases in which the cited reference proves only a special case of our assertion, the transfer to the general case is immediate.

A standard reference for the material presented in this section and the following one is \cite{HumphreysCox}, which considers only reduced root systems. Further, we will refer to \cite{BourbakiLie46} and \cite{HumphreysLieAlg} for statements which are specific to crystallographic root systems. It should be noted that root systems in \cite{HumphreysLieAlg} are assumed to be reduced and crystallographic.

\begin{secnotation}
	Starting from \cref{rootsys:conv}, we assume that $ \roots $ is a root system in some Euclidean space $ (V, \cdot) $.
\end{secnotation}

\subsection{Basic Notions}

\begin{convention}\label{rootsys:convention-right}
	In the context of root systems, we use the the convention that automorphisms of a Euclidean space $ (V, \cdot) $ act from the right-hand side. Thus the image of $ v \in V $ under an automorphism $ \phi $ is denoted by $ v^\phi $. See also \cref{rootsys:right-convention}.
\end{convention}

\begin{definition}[Reflection]\label{rootsys:refl-def}
	Let $ (V, \cdot) $ be a Euclidean space. Then for all $ v \in V \setminus \compactSet{0_V} $, the map
	\[ \map{\refl{v}}{V}{V}{x}{x^{\refl{v}} \defl x - 2 \frac{x \cdot v}{v \cdot v} v} \]
	is called the \defemph*{reflection associated to $ v $}\index{reflection} or the \defemph*{reflection along $ v^\perp $}. We will also denote it by $ \reflbr{v} $. Further, for any word $ \word{v} = \tup{v}{m} $ over $ V \setminus \compactSet{0_V} $, we define $ \refl{\word{v}} \defl \refl{v_1 \cdots v_m} \defl \refl{v_1} \cdots \refl{v_m} $ 
	which we also denote by $ \reflbr{\word{v}} $ or $ \reflbr{v_1 \cdots v_m} $.
\end{definition}

\begin{note}[on \cref{rootsys:convention-right}]\label{rootsys:right-convention}
	The reason to let the reflections act from the right side is that this action is connected to the conjugation action by certain \defemph*{Weyl elements} in $ \roots $-graded groups, which is also written from the right. Since the distinction between left and right actions is only relevant when we compose two reflections, we will still sometimes write $ \refl{v}(w) $ instead of $ w^{\refl{v}} $ for $ v,w \in V $. However, we will never write $ \refl{\word{v}}(w) $ in place of $ w^{\refl{\word{v}}} $ when $ \word{v} $ is a word over $ V $.
\end{note}

\begin{definition}[Root system]
	Let $ (V, \cdot) $ be a Euclidean space. A \defemph*{root system in $ (V, \cdot) $}\index{root system} is a finite non-empty subset $ \roots $ of $ V \setminus \compactSet{0_V} $ such that $ \beta^{\reflbr{\alpha}} \in \roots $ for all $ \alpha, \beta \in \roots $.
	The elements of $ \roots $ are called \defemph*{roots (of $ \roots $)}\index{root} and the dimension of the vector space generated by $ \roots $ is called the \defemph*{rank of $ \roots $}\index{root system!rank}. For any root $ \alpha $, its \defemph*{opposite root}\index{root!opposite} is $ -\alpha $. A root system $ \roots $ is called \defemph*{reduced}\index{root system!reduced} if $ \IR \alpha \intersect \roots = \Set{\alpha, -\alpha} $ for all $ \alpha \in \roots $.
\end{definition}

\begin{secnotation}\label{rootsys:conv}
	From now on, we assume that $ \roots $ is a root system in some Euclidean space $ (V, \cdot) $.
\end{secnotation}

There exists no standard terminology for subsystems of root systems. The following definition is the most general one. We will introduce closed and parabolic subsystems in \cref{rootsys:sub-cl-para-def}.

\begin{definition}[Root subsystem]\label{rootsys:sub-def}
	Let $ \roots $ be a root system. A subset $ \roots' $ of $ \roots $ is called a \defemph*{root subsystem of $ \roots $}\index{root system!subsystem} if $ \refl{\alpha}(\roots') = \roots' $ for all $ \alpha \in \roots' $.
\end{definition}

\begin{definition}[Divisible roots]\label{rootsys:redset}
	A root $ \alpha \in \roots $ is called \defemph*{indivisible (in $ \roots $)}\index{root!indivisible} if $ \lambda \alpha $ is not a root (in $ \roots $) for all $ 0<\lambda < 1 $, and it is called \defemph*{divisible (in $ \roots $)}\index{root!divisible} otherwise. For any set $ \rootsub $ of roots, we denote the set of all indivisible roots in $ \rootsub $ by $ \indivset{\rootsub} $. Further, a subset $ \rootsub $ of $ \roots $ is called \defemph*{reduced}\index{root system!reduced} if all roots in $ \rootsub $ are indivisible in $ \roots $.
\end{definition}

\begin{definition}[Positive system]\label{rootsys:possys-def}
	A subset $ \possys $ of $ \roots $ is called a \defemph*{positive system in $ \roots $}\index{root system!positive system} if there exists a total order $ \totorder $ on $ V $ (in the sense of, for example, \cite[\onpage{7}]{HumphreysCox}) such that $ \possys = \Set{\alpha \in \roots \given 0_V \totorder v} $.
\end{definition}

\begin{definition}[Root base, {\cite[Proposition~VI.1.20, Corollaire~3]{BourbakiLie46}}]\label{rootsys:rootbase-def}
	An \defemph*{root base of $ \roots $}\index{root base} is a subset $ \Delta $ of $ \roots $ satisfying the following conditions:
	\begin{stenumerate}
		\item $ \Delta $ is an (ordered) basis of the $ \IR $-vector space spanned by $ \roots $.
		
		\item If $ \listing{\delta}{l} $ are the pairwise distinct elements of $ \Delta $ and $ \listing{\lambda}{l} \in \IR $ are such that $ \sum_{i=1}^l \lambda_i \delta_i $ is a root, then $ \listing{\lambda}{l} $ are all non-negative or all non-positive.
		
		\item \label{rootsys:rootbase-def:indiv}All the roots $ \listing{\delta}{l} $ are indivisible.
	\end{stenumerate}
	A \defemph*{rescaled root base}\index{root base!rescaled} is a set $ \rootbase $ for which there exist positive real numbers $ (\lambda_\delta)_{\delta \in \rootbase} $ such that $ \Set{\lambda_\delta \delta \given \delta \in \rootbase} $ is a root base. A root base $ \rootbase $ will sometimes be called a \defemph*{proper root base}\index{root base!proper} to emphasise that it is not merely a rescaled root base. If a root base $ \rootbase $ is fixed, its elements are called \defemph*{simple roots}\index{root!simple} and the corresponding reflections are called \defemph*{simple reflections}\index{reflection!simple}.
\end{definition}

\begin{note}
	The notion of rescaled root bases is only relevant if $ \roots $ is not reduced. For root systems of type $ BC $, we will in fact mainly consider rescaled root bases which are not proper root bases.
\end{note}

\begin{proposition}[{\cite[Theorem~1.3]{HumphreysCox}}]\label{rootsys:base-pos-bij}
	Every rescaled root base $ \rootbase $ of $ \roots $ lies in a unique positive system $ \possys(\rootbase) $ (namely, $ \possys(\rootbase) \defl \roots \intersect \Set{\sum_{\delta \in \rootbase} \lambda_\delta \delta \given \lambda_\delta \ge 0 \text{ for all } \delta \in \rootbase} $), and every positive system $ \possys $ contains a unique root basis $ \rootbase(\possys) $ (namely, $ \rootbase(\possys) $ is the set of roots in $ \possys $ which cannot be expressed as an $ \IR_{>0} $-linear combinations of two or more roots in $ \possys $ ). In particular, root bases exist.
\end{proposition}

We now turn to crystallographic root systems.

\begin{definition}[Crystallographic root system]\label{rootsys:cartannum-def}
	Let $ \alpha, \beta $ be two elements of a $ V $. The number $ \cartanint{\alpha}{\beta} \defl 2 \frac{\alpha \cdot \beta}{\beta \cdot \beta} \in \IR $ is called the \defemph*{Cartan number for $ (\alpha, \beta) $}.\index{Cartan number} The root system $ \roots $ is called \defemph*{crystallographic}\index{root system!crystallographic} if for all $ \alpha, \beta \in \roots $, the Cartan number $ \cartanint{\alpha}{\beta} $ is an integer, which is then called the \defemph*{Cartan integer for $ (\alpha, \beta) $}. For any root base $ \rootbase $, the matrix $ (\cartanint{\alpha}{\beta})_{\alpha, \beta \in \rootbase} $ (which does not depend on the choice of $ \rootbase $ up to a permutation of the rows and columns) is called the \defemph*{Cartan matrix of $ \roots $ (with respect to $ \rootbase $)}\index{Cartan matrix}.
\end{definition}

\begin{lemma}\label{rootsys:cry-nonred-2}
	Assume that $ \roots $ is crystallographic. Let $ \alpha $ be a root and let $ \lambda > 0 $ such that $ \lambda \alpha $ is also a root. Then $ \lambda \in \Set{1/2, 1, 2} $.
\end{lemma}
\begin{proof}
	This follows from the fact that $ \cartanint{\alpha}{\lambda \alpha} = 2 \lambda^{-1} $ and $ \cartanint{\lambda \alpha}{\alpha} = 2\lambda $.
\end{proof}

\begin{lemma}[{\cite[Théorème VI.1.3, Proposition VI.1.19]{BourbakiLie46}}]\label{rootsys:cry-int-coeff}
	Let $ \rootbase $ be a root base of $ \roots $ and let $ \alpha $ be a positive root with respect to $ \rootbase $. Assume that $ \roots $ is crystallographic. Then there exist $ n \in \Npos $ and $ \listing{\delta}{n} \in \rootbase $ such that $ \alpha = \sum_{i=1}^n \delta_i $ and such that for each $ m \in \numint{1}{n} $, the partial sum $ \sum_{i=1}^m \delta_i $ is a root. In particular, the coefficients of any root with respect to $ \rootbase $ are always integral.
\end{lemma}

\begin{theorem}[Classification of root systems, {\cite[Theorem~11.4]{HumphreysLieAlg}}, {\cite[Theorème VI.4.3, Section~VI.4.14]{BourbakiLie46}}]\label{rootsys:class}
	Let $ \roots $ be an irreducible crystallographic root system. Then $ \roots $ is isomorphic to one of the root systems $ (A_n)_{n \ge 1} $, $ (B_n)_{n \ge 2} $, $ (C_n)_{n \ge 2} $, $ (BC_n)_{n \ge 1} $, $ (D_n)_{n \ge 4} $, $ E_6 $, $ E_7 $, $ E_8 $, $ F_4 $ or $ G_2 $ where the subscript denotes the rank of the root system. It is reduced if and only if it is not of type $ BC $.
\end{theorem}

We will provide an explicit description of the root systems of types $ A $, $ B $, $ BC $, $ C $ and $ F_4 $ in \cref{ADE:simply-laced:An-standard-rep,B:Bn-standard-rep,BC:BCn-standard-rep,BC:Cn-standard-rep,F4:F-standard-rep}, respectively.

\begin{definition}[Simply-laced root system]\label{rootsys:simply-laced-def}
	A crystallographic root system $ \roots $ is called \defemph*{simply-laced}\index{root system!simply-laced} if each irreducible component of $ \roots $ is of type $ A $, $ D $ or $ E $. (We do not assume that all irreducible components have the same type.)
\end{definition}


\subsection{The Weyl Group}

\begin{definition}[Weyl group]\label{weyl:def}
	The \defemph*{Weyl group of $ \roots $}\index{Weyl group} or the \defemph*{Coxeter group of $ \roots $} is the group $ \Weyl(\roots) \defl \gen{\refl{\alpha} \given \alpha \in \roots} \le O(V,\cdot) $.
\end{definition}

\begin{proposition}[{\cite[Theorem 1.5]{HumphreysCox}}]\label{rootsys:simple-gen-weyl}
	Let $ \rootbase $ be a rescaled root base of $ \roots $. Then the Weyl group of $ \roots $ is generated by the reflections $ \Set{\refl{\delta} \given \delta \in \rootbase} $.
\end{proposition}

\begin{lemma}[{\cite[\onpage{7}]{HumphreysCox}}]
	The Weyl group $ \Weyl(\roots) $ acts faithfully on the subset $ \roots $ of $ V $. In particular, it is a finite group.
\end{lemma}

\begin{proposition}[{\cite[Corollary~1.5]{HumphreysCox}}]\label{rootsys:indiv-in-rootbase}
	Let $ \alpha $ be an indivisible root and let $ \rootbase $ be a root base of $ \roots $. Then there exist $ \delta \in \rootbase $ and $ w \in \Weyl(\roots) $ such that $ \alpha = \delta^w $.
\end{proposition}

\begin{proposition}[Corollary of~\ref{rootsys:indiv-in-rootbase} and~\ref{rootsys:simple-gen-weyl}]\label{rootsys:any-in-rootbase}
	Let $ \alpha $ be any root and let $ \rootbase $ be a rescaled root base of $ \roots $. Then there exist $ \delta, \listing{\delta}{n} \in \rootbase $ and $ \lambda \in \IR_{>0} $ and a word $ \word{\delta} $ over $ \rootbase $ such that $ \lambda \alpha = \delta^{\reflbr{\delta_1 \cdots \delta_m}} $. In particular, $ \refl{\alpha} = \refl{\word{\rho}} $ where $ \word{\rho} \defl (-\delta_n, \ldots, -\delta_1, \delta, \delta_1, \ldots, \delta_n) $. If $ \roots $ is reduced, then we have $ \lambda = 1 $.
\end{proposition}

\begin{proposition}[{\cite[Theorem~1.4,~1.8]{HumphreysCox}}]\label{rootsys:weyl-base-trans}
	The Weyl group of $ \roots $ acts simply transitively on the set of root bases and positive systems of~$ \roots $.
\end{proposition}

\begin{lemma}[Corollary of \ref{rootsys:weyl-base-trans}]\label{rootsys:subsys-bas-conj}
	Let $ \roots' $ be a parabolic root subsystem of $ \roots $, let $ \rootbase' $ be a root base of $ \roots' $ and let $ \rootbase $ be a root base of $ \roots $. Then there exists $ w \in \Weyl(\roots) $ such that $ \rootbase' $ is a subset of $ \rootbase^w $.
\end{lemma}

\begin{proposition}[{\cite[Lemma~10.4.C]{HumphreysLieAlg}}]\label{rootsys:cry-orb-length}
	Assume that $ \roots $ is crystallographic and irreducible and let $ \alpha \in \roots $. Then the orbit of $ \alpha $ under the Weyl group is $ \Set{\beta \in \roots \given \norm{\alpha} = \norm{\beta}} $, the set of roots which have the same length as $ \alpha $.
\end{proposition}

\begin{notation}[Orbits]\label{rootsys:orb-def}
	Assume that $ \roots $ is irreducible and either crystallographic or of type $ H $. We put $ \Orb(\roots) \defl (O_1, \ldots, O_k) $ where $ \listing{O}{k} $ are the orbits of $ \roots $ under the Weyl group, ordered by ascending length of roots. This means that for all $ i<j \in \numint{1}{k} $ and all $ \alpha \in O_i $, $ \beta \in O_j $, we have $ \norm{\alpha} < \norm{\beta} $. Further, we denote by $ \IndivOrb(\roots) $ the sub-tuple of $ \Orb(\roots) $ which contains only the indivisible orbits, that is, the orbits consisting of indivisible roots.
\end{notation}

\begin{proposition}[{\cite[Theorem~1.12]{HumphreysCox}}]\label{weyl:stab-ortho}
	Let $ \alpha \in \roots $. Then the stabiliser of $ \alpha $ in the Weyl group $ W $ is the subgroup of $ W $ which is generated by $ \Set{\refl{\beta} \given \beta \in \roots, \beta \cdot \alpha = 0} $.
\end{proposition}

\begin{remark}[Braid relations]
	Let $ W $ be the Weyl group of $ \roots $ and let $ \rootbase $ be any rescaled root base of $ \roots $. Let $ \delta, \delta' \in \rootbase $ be distinct and denote the order of $ \refl{\delta} \refl{\delta'} $ in $ W $ by $ m $. It follows from $ (\refl{\delta} \refl{\delta'})^m = 1_W $ that we have a relation $ \braidword{m}(\refl{\delta}, \refl{\delta'}) = \braidword{m}(\refl{\delta'}, \refl{\delta}) $ which we call a \defemph*{braid relation}\index{braid relations}. Here $ \braidword{m} $ is as in \cref{braidword-def}.
\end{remark}

\begin{definition}[Homotopy]\label{weyl:homotopy-def}
	Let $ \rootbase $ be a rescaled root base and let $ \word{\alpha} $, $ \word{\beta} $ be two words over $ \rootbase \union (-\rootbase) $. We say that $ \word{\alpha}, \word{\beta} $ are \defemph*{square-homotopic}\index{homotopic!square-} if there exist words $ \word{\rho}, \word{\zeta} $ over $ \rootbase \union (-\rootbase) $ and a root $ \delta \in \rootbase $ such that $ \word{\alpha} = \word{\rho} \word{\zeta} $ and $ \word{\beta} = \word{\rho} \delta \delta \word{\zeta} $ or such that $ \word{\beta} = \word{\rho} \word{\zeta} $ and $ \word{\alpha} = \word{\rho} \delta \delta \word{\zeta} $. We say that $ \word{\alpha}, \word{\beta} $ are \defemph*{braid-homotopic}\index{homotopic!braid-} if there exist words $ \word{\rho}, \word{\zeta} $ and simple roots $ \delta, \delta' \in \rootbase $ such that $ \word{\alpha} = \word{\rho} \braidword{m}(\delta, \delta') \word{\zeta} $ and $ \word{\beta} = \word{\rho} \braidword{m}(\delta', \delta) \word{\zeta} $ where $ m $ denotes the order of $ \refl{\delta} \refl{\delta'} $ in the Weyl group and where $ \braidword{m} $ is as in \cref{braidword-def}. We say that $ \word{\alpha} $, $ \word{\beta} $ are \defemph*{elementary homotopic}\index{homotopic!elementary} if they are square-homotopic or braid-homotopic, and we say that they are \defemph{homotopic} if they lie in the symmetric reflexive transitive closure of the relation \enquote{elementary homotopic}.
\end{definition}

\begin{proposition}[{\cite[{Theorem 8.1}]{HumphreysCox}}]\label{word-problem-solution}
	Let $ \rootbase $ be a rescaled root base of $ \roots $ and let $ \word{\alpha}, \word{\beta} $ be two words over $ \rootbase $. Then $ \word{\alpha}, \word{\beta} $ are homotopic if and only if $ \refl{\word{\alpha}} = \refl{\word{\beta}} $.
\end{proposition}

\begin{note}\label{param:homotopy-def-note}
	We define homotopy for words over $ \rootbase \union (-\rootbase) $ and not merely for words over $ \rootbase $. In the context of the word problem, this may seem unnecessary because $ \refl{-\delta} = \refl{\delta} $. However, we will later use words $ \word{\delta} = (\delta_1, \ldots, \delta_k) $ over $ \rootbase \union (-\rootbase) $ to represent group elements $ w_{\word{\delta}} = w_{\delta_1} \cdots w_{\delta_k} $ where $ (w_\delta)_{\delta \in \rootbase} $ is a fixed family of so-called Weyl elements in a root graded group and where $ w_{-\delta} \defl w_\delta^{-1} $ for all $ \delta \in \rootbase $. In this setting, it is not (necessarily) true that $ w_{-\delta} = w_\delta $, and thus it makes sense to consider words over $ \rootbase \union (-\rootbase) $.
	
	However, observe that the words $ \word{\alpha} $ and $ \word{\beta} $ in \cref{word-problem-solution} are not allowed to contain letters from $ -\rootbase $. As a counterexample, the words $ \word{\alpha} \defl (\delta, -\delta) $ (for any $ \delta \in \rootbase $) and the empty word represent the same element in the Weyl group, but they are not homotopic. This fact will be relevant in the proof of \cref{param:inverpar-stab}.
\end{note}


\subsection{Root Intervals and Root Subsystems}

\label{subsec:root-int}

The notion of root intervals will be crucial for the definition of groups with $ \roots $-commutator relations. It is a special case of the notion of commutator sets (of subsets of $ \roots $) that we will introduce in \cref{basic:commset}.

\begin{definition}[Root intervals]\label{rootsys:interval-def}
	Let $ \alpha, \beta $ be non-proportional roots. The \defemph*{open Coxeter root interval for $ (\alpha, \beta) $}\index{root interval} is
	\[ \rootintcox{\alpha}{\beta} \defl \Set{\lambda \alpha + \mu \beta \given \lambda, \mu \in \IR_{>0}} \intersect \roots \]
	and the \defemph*{open crystallographic root interval for $ (\alpha, \beta) $} is
	\[ \rootint{\alpha}{\beta} \defl \Set{\lambda \alpha + \mu \beta \given \lambda, \mu \in \Npos} \intersect \roots \subs \rootintcox{\alpha}{\beta}. \]
	The roots $ \alpha $ and $ \beta $ are called \defemph{adjacent} if $ \rootintcox{\alpha}{\beta} $ is empty, and they are called \defemph*{crystallographically adjacent}\index{adjacent!crystallographically} if $ \rootint{\alpha}{\beta} $ is empty. Further, the sets
	\[ \clrootintcox{\alpha}{\beta} \defl \Set{\alpha, \beta} \union \rootintcox{\alpha}{\beta} \midand \clrootint{\alpha}{\beta} \defl \Set{\alpha, \beta} \union \rootint{\alpha}{\beta} \]
	are called the \defemph*{closed Coxeter root interval for $ (\alpha, \beta) $} and the \defemph*{closed crystallographic root interval for $ (\alpha, \beta) $}, respectively.
\end{definition}

\begin{note}\label{rootsys:loosneher-interval-def}
	In \cite[(1.5.2), (1.6.3)]{LoosNeherBook}, root intervals are defined as follows for all $ \alpha, \beta \in \roots $:
	\begin{align*}
		\clrootint{\alpha}{\beta} &\defl \roots \intersect \Set{m \alpha + n\beta \given \alpha, \beta \in \Nzero, m+n>0}, \\
		\rootint{\alpha}{\beta} &\defl \roots \intersect \Set{m \alpha + n\beta \given \alpha, \beta \in \Npos}.
	\end{align*}
	These definitions differ from ours in two ways. Firstly, there is no assumption on $ \alpha $ and $ \beta $ to be non-proportional. For example, we have $ \clrootint{\alpha}{\alpha} = \roots \intersect \Npos \alpha $ and $ \rootint{\alpha}{\alpha} = \roots \intersect \IN_{\ge 2} \alpha $. Secondly, the closed root interval $ \clrootint{\alpha}{\beta} $ in the sense of Loos-Neher contains, in addition to the closed root interval in our sense, all the roots in $ \rootint{\alpha}{\alpha} \union \rootint{\beta}{\beta} = \roots \intersect (\IN_{\ge 2} \alpha \union \IN_{\ge 2} \beta) $.
\end{note}

\begin{note}\label{rootsys:interval-note}
	Let $ \alpha, \beta $ be non-proportional roots and assume that $ \roots $ is irreducible and crystallographic. Clearly, the crystallographic root interval is a subset of the Coxeter root interval. (In particular, adjacent roots are automatically crystallographically adjacent.) In the general case, it can be smaller. Using the classification of root systems, we can prove the following statements by an inspection of the irreducible root systems $ \roots $:
	\begin{lemenumerate}
		\item If $ \roots $ is simply-laced (that is, of type $ A $, $ D $ or $ E $), then the Coxeter and crystallographic intervals are the same for all pairs of roots.
		
		\item If $ \roots $ is not of type $ G_2 $, then $ \rootintcox{\alpha}{\beta} $ has at most two elements, and if it has two elements, it equals $ \rootint{\alpha}{\beta} $.
	\end{lemenumerate}
	Thus if $ \roots $ is not of type $ G_2 $, then the only difference between Coxeter and crystallographic root intervals is that for some pairs of roots, $ \rootint{\alpha}{\beta} $ is empty but $ \rootintcox{\alpha}{\beta} $ has exactly one element. Specifically, using the standard representations, this happens exactly for the pairs of roots $ \brackets[\big]{\epsilon (\basvec_i - \basvec_j), \epsilon (\basvec_i + \basvec_j)} $ in $ B_n $ and the pairs $ (\epsilon 2\basvec_i, \epsilon' 2\basvec_j) $ in $ (B)C_n $ (where $ i,j \in \numint{1}{n} $ are distinct and $ \epsilon, \epsilon' \in \compactSet{\pm 1} $ are signs).
\end{note}

Since the interval between two roots lies in a subspace of dimension~2, it is geometrically clear that there exist two \enquote{natural} ways to put an order on the roots in this interval. These two natural orders will be called \defemph*{interval orderings}.

\begin{definition}[Interval ordering]\label{rootsys:clockwise-def}
	Let $ \roots $ be root system and let $ S $ be a subset of $ \indivset{\roots} $ which does not contain a pair of opposite roots. Denote by $ k $ the cardinality of $ S $. An \defemph*{interval ordering of $ S $}\index{interval ordering} is a tuple $ (\alpha_1, \ldots, \alpha_k) $ such that $ S = \Set{\alpha_1, \ldots, \alpha_k} $ and $ \indivset{\rootintcox{\alpha_i}{\alpha_j}} = \Set{\alpha_r \given i < r < j} $ for all $ i<j \in \numint{1}{k} $.
\end{definition}

\begin{remark}\label{rootsys:clockwise-rem}
	Let $ S $ be a subset of a root system $ \indivset{\roots} $ and assume that there exists an interval ordering $ (\alpha_1, \ldots, \alpha_k) $ of $ S $. Then we have
	\[ S = \Set{\alpha_1, \ldots, \alpha_k} = \Set{\alpha_1, \alpha_k} \union \indivset{\rootintcox{\alpha_1}{\alpha_k}} = \indivset{\clrootintcox{\alpha_1}{\alpha_k}}. \]
	Thus if there exists an interval ordering of $ S $, then $ S $ is the set of indivisible roots in a (closed) root interval. Conversely, any set of the form $ S = \indivset{\clrootintcox{\alpha}{\beta}} $ for $ \alpha, \beta \in \indivset{\roots} $ has an interval ordering $ (\alpha_1, \ldots, \alpha_k) $, and the only other interval ordering of $ S $ is $ (\alpha_k, \ldots, \alpha_1) $.
\end{remark}

\begin{note}
	We could also define \enquote{crystallographic interval orderings} by requiring that $ \indivset{\rootint{\alpha_i}{\alpha_j}} = \Set{\alpha_r \given i < r < j} $ for all $ i<j \in \numint{1}{k} $. However, there is no need to do this: We will only use interval orderings of root intervals of length 2, and for these there is no difference by \cref{rootsys:interval-note} in our situation.
\end{note}

Using interval orderings, we can introduce a practical way of indexing the roots of $ \roots $ in the rank-2 case, which depends only on the choice of a root base.

\begin{remark}[Rank-2 root systems, {\cite[(4.14), (4.15)]{MoufangPolygons}}]\label{basic:rank2-notation}
	Assume that $ \roots $ is reduced and of rank 2, put $ n \defl \abs{\roots}/2 $ and choose a root base $ (\alpha, \beta) $. Denote by $ \tup{\alpha}{n} $ the unique interval ordering of $ \clrootintcox{\alpha}{\beta} $ with $ \alpha_1 = \alpha $ and $ \alpha_n = \beta $ (see \cref{rootsys:clockwise-rem}). Further, put $ \alpha_{i+n} \defl -\alpha_i $ for all $ i \in \numint{1}{n} $. In this way, we have identified $ \roots $ with the set $ \numint{1}{2n} $. It is practical to define $ \alpha_z \defl \alpha_{z'} $ for all $ z \in \IZ $ where $ z' $ is the unique element of $ \numint{1}{2n} $ which is congruent to $ z $ modulo $ 2n $. Then for all $ i<j \in \IZ $ with $ j-i \le n-1 $, the tuple $ (\alpha_i, \ldots, \alpha_j) $ is the unique interval ordering of $ \rootintcox{\alpha_i}{\alpha_j} $ which starts with $ \alpha_i $. Observe that for all $ z \in \IZ $, the pair $ \rootbase' = (\alpha_z, \alpha_{z+n-1}) $ is a root base, and the labeling which is induced by $ \rootbase' $ is given by $ \alpha_i' = \alpha_{z-1+i} $ for all $ i \in \IZ $. Further, we have
	\[ \alpha_i^{\reflbr{\alpha_j}} = \alpha_{2j+n-i} \]
	for all $ i,j \in \IZ $.
\end{remark}

\begin{definition}[Rank-2 labeling]\label{rootsys:rank2-label-def}
	Assume that $ \roots $ is reduced of rank 2 and let $ (\alpha, \beta) $ be a root base. The \defemph*{rank-2 labeling of $ \roots $ induced by $ (\alpha, \beta) $}\index{rank-2 labeling} is the map $ \map{}{\IZ}{\roots}{i}{\alpha_i} $ which we constructed in \cref{basic:rank2-notation}. A \defemph*{rank-2 labeling of $ \roots $} is the rank-2 labeling induced by some root base.
\end{definition}

The notion of root intervals leads in a natural way to the notion of closed sets of roots. We will study this property more closely in \cref{sec:closed}.

\begin{definition}[Closed set of roots]\label{rootsys:closed-def}
	Let $ \rootsub $ be a subset of $ \roots $. Then $ \rootsub $ is called \defemph*{closed}\index{closed set of roots} if $ \rootintcox{\alpha}{\beta} \subs \rootsub $ for all non-proportional $ \alpha, \beta \in \rootsub $, and it is called \defemph*{crystallographically closed} if $ \rootint{\alpha}{\beta} \subs \rootsub $ for all $ \alpha, \beta \in \rootsub $.
\end{definition}

\begin{remark}
	Every closed set is also crystallographically closed, but the converse is not true.
\end{remark}

We now refine our notion of root subsystems from \cref{rootsys:sub-def}.

\begin{lemma}\label{rootsys:subsys}
	Let $ V' $ be a subspace of the Euclidean space $ V $ surrounding $ \roots $ and put $ \roots' \defl V' \intersect \roots $. Assume that $ \roots' $ is non-empty. Then the following hold:
	\begin{lemenumerate}
		\item \label{rootsys:subsys:sub}$ \roots' $ is a root system in $ V' $ which is closed in $ \roots $. If $ \roots $ is reduced and/or crystallographic, then $ \roots' $ has the same properties.
		
		\item \label{rootsys:subsys:base}Let $ \rootbase' $ be a root base of $ \roots' $. Then there exists a root base $ \rootbase $ of $ \roots $ such that $ \rootbase \intersect \roots' = \rootbase' $.
		
		\item \label{rootsys:subsys:pos}Let $ \possys' $ be a positive system in $ \roots' $. Then there exists a positive system $ \possys $ in $ \roots $ such that $ \possys \intersect \roots' = \possys' $.
	\end{lemenumerate}
\end{lemma}

\begin{definition}[Closed and parabolic root subsystems]\label{rootsys:sub-cl-para-def}
	Let $ \roots' $ be a non-empty subset of $ \roots $. Then $ \roots' $ is called a \defemph*{(crystallographically) closed root subsystem of $ \roots $}\index{root system!subsystem!closed} if it is (crystallographically) closed in the sense of \cref{rootsys:closed-def} and a root subsystem in the sense of \cref{rootsys:sub-def}. It is called a \defemph*{parabolic root subsystem of $ \roots $}\index{root system!subsystem!parabolic} if there exists a subspace $ W $ of $ V $ such that $ \roots' = \roots \intersect W $. For any non-empty subset $ A $ of $ \roots $, the \defemph*{root subsystem spanned by $ A $}\index{root system!subsystem!spanned by a subset} is $ \roots \intersect \gen{A}_\IR $.
\end{definition}

\begin{remark}\label{rootsys:subsys:clos-pos}
	The statement of \thmitemcref{rootsys:subsys}{rootsys:subsys:pos} remains true for closed root subsystems which are not parabolic.
\end{remark}

\begin{lemma}
	Let $ \rootbase $ be a root base of $ \roots $ and let $ \rootbase' $ be a subset of $ \rootbase $. Then $ \rootbase' $ is a root base of the subsystem that it spans.
\end{lemma}


\subsection{The Length Function in Weyl Groups}

\begin{secnotation}
	From now on, we denote by $ \roots $ a root system, by $ \rootbase $ a rescaled root base of $ \roots $, by $ \possys $ the unique positive system in $ \roots $ which contains $ \rootbase $ (as in \cref{rootsys:base-pos-bij}) and by $ W $ the Weyl group of $ \roots $.
\end{secnotation}

In this subsection, we introduce a natural ordering of a certain set $ \switchset(w) \subs \possys $ (for any $ w \in W $), using the notion of root sequences. Standard references are \cite[1.6]{HumphreysCox} and \cite{Howlett}, the latter being closer to our presentation. The main results are \cref{rootsys:switchset-lem,rootsys:switchset-description}.

\begin{definition}[Length function]\label{rootsys:length-def}
	Let $ w \in W $. An \defemph*{expression of $ w $}\index{expression} is a word $ \word{\delta} = (\delta_1, \ldots, \delta_k) $ over $ \rootbase $ such that $ w = \refl{\delta_1} \cdots \refl{\delta_k} $, and the number $ k $ is called the \defemph*{length of $ \word{\delta} $}. If $ \rootbase $ is not clear from the context, then an expression of $ w $ is also called a \defemph*{$ \rootbase $-expression of $ w $}. The \defemph*{length of $ w $}\index{length function} is the minimal length of an expression of $ w $, and it is denoted by $ \len_\rootbase(w) $ (or simply by $ \len(w) $). A \defemph*{reduced expression of $ w $}\index{expression!reduced} is an expression of $ w $ of length $ \len(w) $, and an arbitrary word $ \word{\delta}  $ over $ \rootbase $ is called \defemph*{reduced} if it is a reduced expression of $ \refl{\word{\delta}} $.
\end{definition}

\begin{definition}\label{rootsys:switchset-def}
	For all $ w \in W $, we define
	\[ \switchset(w) \defl \switchset_\possys(w) \defl \Set{\alpha \in \indivset{\possys} \given \alpha^{w^{-1}} \in -\possys}. \]
\end{definition}

\begin{lemma}\label{rootsys:switchset-complete}
	Let $ w \in W $ and let $ \alpha \in \indivset{\roots} \setminus \switchset(w) $. Then there exists a positive system $ \possys' $ which contains $ \switchset(w) $ but not $ \alpha $.
\end{lemma}
\begin{proof}
	Since $ \alpha $ is not contained in $ \switchset(w) $, we have either $ \alpha \in -\possys $ or $ \alpha^{w^{-1}} \in \possys $. In the first case, we can choose $ \possys' \defl \possys $. In the second case, we can choose $ \possys \defl (-\possys)^w $.
\end{proof}

\begin{lemma}[{\cite[Proposition~1.4]{HumphreysCox}}]\label{rootsys:simple-refl-on-pos}
	For all $ \delta \in \rootbase $, we have $ (\possys \setminus \IR_{>0} \delta)^{\reflbr{\delta}} = \possys \setminus \IR_{>0} \delta $. In particular, $ \switchset(\refl{\delta}) = \compactSet{\delta} $.
\end{lemma}

\begin{lemma}[{\cite[Lemma~1.6]{HumphreysCox}}]\label{rootsys:switchset-lem}
	Assume that $ \rootbase $ is a proper root base. Let $ w \in W $ and let $ \delta \in \rootbase $ such that $ \len(w \refl{\delta}) > \len(w) $. Then $ \switchset(w\refl{\delta}) = \switchset(w)^{\reflbr{\delta}} \disjunion \compactSet{\delta} $ (where \enquote{$ \disjunion $} denotes disjoint union).
\end{lemma}

\begin{definition}[Root sequence]\label{rootsys:rootseq-def}
	Let $ \word{\delta} = \tup{\delta}{k} $ be a reduced word over $ \rootbase $. The tuple $ \tup{\beta}{k} $ where $ \beta_i \defl \delta_i^{\reflbr{\alpha_{i+1} \cdots \alpha_k}} $ for all $ i \in \numint{1}{k-1} $ and $ \beta_k \defl \delta_k $ is called the \defemph*{root sequence associated to $ \word{\delta} $}.\index{root sequence} Further, the tuple $ (\beta_k, \ldots, \beta_1) $ is called the \defemph*{inverse root sequence associated to $ \word{\delta} $}.\index{root sequence!inverse}
\end{definition}

\begin{remark}\label{rootsys:rootseq-iter}
	Let $ \word{\delta} = \tup{\delta}{k} $ be a reduced word over $ \rootbase $ such that $ k>0 $ and put $ s \defl \reflbr{\delta_k} $. Denote by $ \tup{\alpha}{k-1} $ the root sequence of $ \tup{\delta}{k-1} $. Then it is clear from \cref{rootsys:rootseq-def} that the root sequence of $ \word{\delta} $ is $ (\alpha_1^s, \ldots, \alpha_{k-1}^s, \delta_k) $.
\end{remark}

\begin{proposition}\label{rootsys:switchset-description}
	Assume that $ \rootbase $ is a proper root base. Let $ w \in W $, let $ \word{\delta} = (\delta_1, \ldots, \delta_k) $ be a reduced expression of $ w $ and let $ \word{\alpha} = \tup{\alpha}{k} $ be the associated root sequence. Then $ \switchset(w) = \Set{\alpha_1, \ldots, \alpha_k} $ and $ \abs{\switchset(w)} = k = \len(w) $.
\end{proposition}
\begin{proof}
	This is clear for $ k=0 $ because $ \switchset(1_W) = \emptyset $. For $ k>0 $, the assertion follows from \cref{rootsys:switchset-lem} by induction on $ k $.
\end{proof}

We will show in \cref{rootorder:switchset-extremal} that the ordering of $ \switchset(w) $ given by an inverse root sequence has the desirable property of being extremal.

\begin{proposition}[{\cite[Theorem~1.8]{HumphreysCox}}]\label{rootsys:longest-el-char}
	There exists a unique element $ \rho $ of $ W $ such that $ \len(\rho) \ge \len(w) $ for all $ w \in W $. It is called the \defemph*{longest element of $ W $ (with respect to $ \rootbase $)}.\index{longest element} It has length $ \abs{\indivset{\possys}} $ and satisfies $ \rho^{-1} = \rho $ and $ \rho(\possys) = -\possys $. In particular, $ \switchset(\rho) = \indivset{\possys} $.
\end{proposition}


\section{Foldings of Root Systems}

\label{sec:fold}

We have seen in \cref{sec:rootsys} that there are many ways to regard root systems as subsystems of other root systems. Foldings, on the other hands, behave more like \enquote{quotients of root systems}. If a root system $ \roots' $ is a folding of $ \roots $, then $ \roots' $-graded groups can be constructed from $ \roots $-graded group in a way which we will describe in \cref{sec:RGG-fold}. We will use this strategy in \cref{sec:F4:const} to construct $ F_4 $-graded groups from $ E_6 $-graded groups.

In this section, we follow the exposition in \cite[13.1, 13.2]{Carter-Chev}.

\begin{secnotation}\label{secnot:rootsys-fold}
	We denote by $ \roots $ a crystallographic root system in the Euclidean space $ (V, \cdot) $ such that $ V $ is generated by $ \roots $, and we choose a root base $ \rootbase $ of $ \roots $. We denote by $ W $ the Weyl group of $ \roots $ and by $ \rho $ an isomorphism of the Coxeter diagram of $ \roots $ with respect to $ \rootbase $. That is, $ \rho $ is a vector space automorphism of $ V $ such that $ \rho(\rootbase) = \rootbase $ and such that the group elements $ \refl{\alpha} \refl{\beta} $ and $ \refl{\rho(\alpha)} \refl{\rho(\beta)} $ have the same order for all $ \alpha, \beta \in \rootbase $. We denote by $ \map{\tau}{V}{V}{}{} $ the unique automorphism of $ V $ which maps any $ v \in V \setminus \compactSet{0} $ to the unique element $ w \in \IR_{>0} \rho(v) $ with $ w \cdot w = v \cdot v $. Further, we put
	\[ \fixspace \defl \Set{v \in V \given \tau(v) = v} \]
	and we denote by $ \map{\pi}{V}{\fixspace}{}{} $ the orthogonal projection on $ \fixspace $. Finally, we set $ \roots' \defl \pi(\roots) $.
\end{secnotation}

\begin{notation}
	Following the conventions in \cite[13.1, 13.2]{Carter-Chev}, we let the Weyl group of $ \roots $ act on $ V $ from the left side in this section.
\end{notation}

\begin{note}
	The map $ \rho $ is not necessarily an isometry of $ V $, and the map $ \tau $ is precisely the unique isometry which maps any $ v \in V \setminus \compactSet{0} $ to a positive scalar multiple of $ \rho(v) $. If all roots in $ \roots $ have the same length, then $ \tau = \rho $. In practice, we will only need the case $ \roots = E_6 $ in which this assumption is satisfied.
\end{note}

\begin{lemma}[{\cite[\onpage{217}]{Carter-Chev}}]\label{fold:proj-average}
	Let $ v \in V $ and denote by $ J $ the orbit of $ v $ under $ \tau $ in $ V $. Then $ \pi(v) $ is the average of the vectors in $ J $. That is, we have
	\[ \pi(v) = \frac{1}{\abs{J}} \sum_{u \in J} u. \]
\end{lemma}
\begin{proof}
	Put $ v' \defl \frac{1}{\abs{J}} \sum_{u \in J} u $. Then $ \pi(v') = v' $, so $ v' $ is contained in $ \fixspace $. Further, since $ \tau $ is an isometry, we have for all $ k \in \Npos $ and $ u \in \fixspace $ that
	\begin{align*}
		u \cdot v = \tau^k(u) \cdot \tau^k(v) = u \cdot \tau^k(v).
	\end{align*}
	Hence $ v' $ is an element of $ \fixspace $ with $ u \cdot v = u \cdot v' $ for all $ u \in \fixspace $. Since $ \pi(v) $ is the unique element with this property, it follows that $ \pi(v) = v' $.
\end{proof}

\begin{definition}[Folded Weyl group]
	We put
	\[ W' \defl \Set{w \in W \given \tau^{-1} w \tau = w} = \Set{w \in W \given w \tau = \tau w}. \]
\end{definition}

\begin{lemma}\label{fold:proj-commute}
	For all $ w \in W' $ and $ v \in V $, we have $ \pi(w(v)) = w(\pi(v)) $.
\end{lemma}
\begin{proof}
	Denote by $ J $ and $ J' $ the orbits of $ v $ and $ w(v) $ under $ \tau $ in $ V $, respectively. Since $ w \tau^k(v) = \tau^k(w(v)) $ for all $ k \in \Npos $ by the definition of $ W' $, we have $ w(J) = J' $. Since $ \pi(w(v)) $ is the average of $ J' $ and $ \pi(v) $ is the average of $ J $ by \cref{fold:proj-average}, the assertion follows.
\end{proof}

\begin{lemma}[{\cite[13.1.1]{Carter-Chev}}]
	For all $ w \in W' $, we have $ w(\fixspace) = \fixspace $, and the action of $ W' $ on $ \fixspace $ is faithful.
\end{lemma}

\begin{notation}
	Let $ J $ be any subset of $ \rootbase $. We denote by $ W_J $ the subgroup of $ W $ which is generated by $ \Set{\refl{\alpha} \given \alpha \in J} $. This is the Weyl group of the root system $ \gen{J}_\IR \intersect \roots $ with respect to the root base $ J $, and we denote its longest element by $ w_0^J $.
\end{notation}

\begin{lemma}[{\cite[13.1.2, 13.1.3]{Carter-Chev}}]\label{fold:W'}
	The following hold:
	\begin{lemenumerate}
		\item \label{fold:W':gen-ident}For any orbit $ J $ of $ \rho $ in $ \rootbase $, the element $ w_0^J $ lies in $ W' $, and its action on $ \fixspace $ agrees with the one of $ \refl{\pi(\alpha)} $ for every $ \alpha \in J $.
		
		\item The set $ \Set{w_0^J \given J \text{ orbit of } \rho \text{ in } \rootbase} $ generates $ W' $.
	\end{lemenumerate}
\end{lemma}

The following result is a simple corollary of \cref{fold:W'}, but its specific form will be useful later.

\begin{lemma}\label{fold:shortest-act}
	Let $ \alpha \in \roots $, $ \beta \in \rootbase $ and denote by $ J $ the orbit of $ \beta $ under $ \rho $. Then
	\[ \pi\brackets[\big]{w_0^J(\alpha)} = \refl{\pi(\beta)}\brackets[\big]{\pi(\alpha)}. \]
\end{lemma}
\begin{proof}
	By \cref{fold:proj-commute}, $ \pi(w_0^J(\alpha)) = w_0^J(\pi(\alpha)) $. By \thmitemcref{fold:W'}{fold:W':gen-ident},
	\[ w_0^J(\pi(\alpha)) = \refl{\pi(\beta)}(\pi(\alpha)). \]
	The assertion follows.
\end{proof}

\begin{proposition}[{\cite[13.2.2]{Carter-Chev}}]\label{fold:rootsys}
	The following hold:
	\begin{proenumerate}
		\item $ \roots' = \pi(\roots) $ is a root system in $ \fixspace $ which spans $ \fixspace $.
		
		\item \label{fold:rootsys:rootbase}The set of indivisible elements in $ \pi(\rootbase) $ is a root base of $ \roots' $.
	\end{proenumerate}
\end{proposition}

\begin{proposition}\label{fold:cor}
	The following hold:
	\begin{proenumerate}
		\item \label{fold:cor:pts-induces}Let $ \possys $ be a positive system in $ \roots $. Then $ \pi(\possys) $ is a positive system in $ \roots' $.
		
		\item \label{fold:cor:ind-surj}Every root base and every positive system in $ \roots' $ is induced by a corresponding object in $ \roots $ via $ \pi $.
		
		\item \label{fold:cor:preim-pts}Let $ \alpha', \beta' \in \roots' $ be non-proportional. Then there exists a positive system in $ \roots $ which contains $ \Set{\gamma \in \roots \given \pi(\gamma) \in \compactSet{\alpha', \beta'}} $.
	\end{proenumerate}
\end{proposition}
\begin{proof}
	Let $ \possys $ be a positive system in $ \roots $ and denote by $ \bar{\rootbase} $ the corresponding root base. By \thmitemcref{fold:rootsys}{fold:rootsys:rootbase}, $ \rootbase' \defl \indivset{\pi(\bar{\rootbase})} $ is a root base of $ \roots' $. Now $ \pi(\possys) $ consists of positive linear combinations of $ \rootbase' $ while $ \pi(-\possys) $ consists of negative linear combinations of $ \rootbase $, and the disjoint union of these two sets is $ \roots' $. Hence $ \pi(\possys) $ is a positive system in $ \roots' $. This proves~\itemref{fold:cor:pts-induces}.
	
	Now let $ \possys' $ be an arbitrary positive system in $ \roots' $ and let $ \possys $ be an arbitrary positive system in $ \roots $. By \cref{rootsys:weyl-base-trans}, there exists an element $ w $ of the Weyl group of $ \roots' $ such that $ \possys' = w(\pi(\possys)) $. By \itemref{fold:W':gen-ident}, we can regard $ w $ as an element of $ W' $. Hence by \cref{fold:proj-commute}, $ w(\pi(\possys)) = \pi(w(\possys)) $. Thus $ \possys' $ is induced by the positive system $ w(\possys) $ in $ \roots $. The assertion for root bases can be proven similarly, so~\itemref{fold:cor:ind-surj} holds.
	
	Now let $ \alpha', \beta' \in \roots' $ be non-proportional. Then there exists a positive system $ \possys' $ in $ \roots' $ which contains $ \alpha' $ and $ \beta' $. Denote by $ \possys $ a positive system in $ \roots $ such that $ \pi(\possys) = \possys' $. Let $ \alpha \in \pi^{-1}(\alpha') \intersect \roots $ be arbitrary, and suppose for a contradiction that $ \alpha \nin \possys $. Then $ \alpha $ lies in $ -\possys $, so
	\[ \alpha' = \pi(\alpha) \in \pi(-\possys) = -\possys'. \]
	This contradicts the choice of $ \possys' $. Hence $ \pi^{-1}(\alpha') \intersect \roots $ is contained in $ \possys $. In a similar way, we can show that $ \pi^{-1}(\beta') \intersect \roots $ is contained in $ \possys $. This finishes the proof.
\end{proof}

\begin{definition}[Folding]
	The root system $ \roots' $ is called a \defemph*{folding of $ \roots $}\index{folding!of root systems}. We will sometimes refer to the map $ \map{\pi}{\roots}{\roots'}{}{} $ as a folding as well.
\end{definition}

	\chapter{Root Graded Groups: Definition and General Observations}
	
	\label{chap:rgg}
	
	In this chapter, we introduce the protagonists of this work: $ \roots $-graded groups where $ \roots $ is an arbitrary root system. These are groups $ G $ together with a family $ (\rootgr{\alpha})_{\alpha \in \roots} $ of subgroups such that, technical details aside, two crucial properties are satisfied: They satisfy some commutator relations and they have Weyl elements. These two conditions are relatively straightforward to define, and we will separately introduce them in \cref{sec:group-commrel,sec:weyl}, respectively.
	
	For technical reasons, the commutator relations and Weyl elements are not enough to build a satisfactory theory of root graded groups. The reason for this is that we want to decompose certain commutators $ \commutator{x_\alpha}{x_\gamma} $ in $ G $ (where $ \alpha, \gamma $ are non-proportional roots and $ x_\alpha \in \rootgr{\alpha} $, $ x_\gamma \in \rootgr{\gamma} $) as products $ x_1 \cdots x_m $ where $ x_i \in \rootgr{\beta_i} $ and $ \listing{\beta}{m} $ are the roots in the root interval $ \rootintcox{\alpha}{\gamma} $. For this decomposition to exist and be unique, we have to assume that the product map (on the set of roots $ \rootintcox{\alpha}{\gamma} $) in $ G $ is bijective, which need not be the case in general. This should be thought of as a non-degeneracy condition which any root graded group has to satisfy.
	
	Unfortunately, it is not at all clear how exactly this non-degeneracy condition should be formulated. There are several possible candidates which seem like a reasonable choice (and which we compare in \cref{rgg:nondeg-comparison}). In order to understand the situation better, we will first introduce and study the purely combinatorial (that is, root-system-theoretic) notion of closed sets of roots in \cref{sec:closed}. This concept is indispensable for a proper study of the product maps in root graded groups. In \cref{sec:prod-bij}, we will investigate several conditions which imply the injectivity or surjectivity of certain product maps. The main result of this section is a simple criterion which guarantees the bijectivity of all product maps on closed sets of roots. This criterion serves as the third axiom for root graded groups, which we can finally introduce in \cref{sec:rgg-def}. In \cref{sec:RGG-fold}, we describe a construction of $ \roots' $-graded groups from $ \roots $-graded groups for every folding $ \map{\pi}{\roots}{\roots'}{}{} $. In \cref{sec:literature}, we discuss various special cases of root gradings which have already been considered in the literature, and we summarise the results which have been obtained.
	

\section{Groups with Commutator Relations}

\label{sec:group-commrel}

\begin{secnotation}
	We denote by $ \roots $ an arbitrary root system.
\end{secnotation}

\subsection{Definition and Comments}

We begin by introducing the notion of $ \roots $-pregradings. Without further assumptions, this is not an interesting concept, but it allows us to avoid writing \enquote{family of subgroups of $ G $} all the time.

\begin{definition}[Pregradings]\label{rgg:pregrade-def}
	Let $ G $ be a group. A \defemph*{$ \roots $-pregrading of $ G $}\index{root pregrading}\index{pregrading|see{root pregrading}} is a family $ (\rootgr{\alpha})_{\alpha \in \roots} $ consisting of subgroups of $ G $. These subgroups are called the \defemph*{root groups of $ G $}\index{root group}.
	For any subset $ S $ of $ \roots $, we denote by $ \rootgr{S} $ the subgroup of $ G $ which is generated by $ (\rootgr{\alpha})_{\alpha \in S} $. For any word $ \word{\alpha} = \tup{\alpha}{k} $ over $ \roots $, we define $ \rootgr{\word{\alpha}} \defl \rootgr{\alpha_1} \times \cdots \times \rootgr{\alpha_k} $. If $ \roots $ is of rank-2 and $ \map{}{\IZ}{\roots}{i}{\alpha_i} $ is a fixed rank-2 labeling of $ \roots $ (in the sense of \cref{rootsys:rank2-label-def}), then we put $ \rootgr{i} \defl \rootgr{\alpha_i} $ for all $ i \in \IZ $.
\end{definition}

\begin{convention}
	Let $ G $ be a group with a $ \roots $-pregrading $ (\rootgr{\alpha})_{\alpha \in \roots} $. If $ X $ is a property that roots may or may not have, we will also refer to $ X $ as a property of the root groups. For example, a root group $ \rootgr{\alpha} $ is called indivisible if $ \alpha $ is an indivisible root. Most of the time, we will use this convention when refering to root lengths: A root group $ \rootgr{\alpha} $ is called short or long if $ \alpha $ is short or long.
\end{convention}

At this point, the reader should recall the notion of root intervals from \cref{subsec:root-int}.

\begin{definition}[Group with commutator relations]\label{rgg:group-commrel-def}
	Let $ G $ be a group with a $ \roots $-pregrading $ (\rootgr{\alpha})_{\alpha \in \roots} $. We say that $ G $ \defemph*{has $ \roots $-commutator relations with root groups $ (\rootgr{\alpha})_{\alpha \in \roots} $}\index{group with commutator relations}\index{commutator relations} if the following conditions are satisfied:
	\begin{description}[font=\normalfont]
		\customitem[(RGG-Com)] \label{rgg-axiom-comm}\index{RGG-Com@(RGG-Com)}
		For all distinct, non-proportional roots $ \alpha, \beta $, we have 
		\[ \commutator{\rootgr{\alpha}}{\rootgr{\beta}} \subs \rootgr{\rootintcox{\alpha}{\beta}}. \]
		
		\customitem[(RGG-Div)] \label{rgg-axiom-div}\index{RGG-Div@(RGG-Div)}For all roots $ \alpha $ for which $ 2\alpha $ is a root, we have
		\[ \rootgr{2\alpha} \subs \rootgr{\alpha}, \qquad \commutator{\rootgr{\alpha}}{\rootgr{2\alpha}} = \compactSet{1_G} \midand \commutator{\rootgr{\alpha}}{\rootgr{\alpha}} \subs \rootgr{2\alpha}. \]
	\end{description}
\end{definition}

\begin{definition}[Group with crystallographic commutator relations]\label{rgg:group-commrel-cry-def}
	Let $ G $ be a group which has $ \roots $-commutator relations with root groups $ (\rootgr{\alpha})_{\alpha \in \roots} $. We say that $ G $ \defemph*{has crystallographic $ \roots $-commutator relations} if the following stronger version of \axiomref{rgg-axiom-comm} is satisfied:
	\begin{description}[font=\normalfont]
		\customitem[(RGG-Com-cry)] \label{rgg-axiom-comm-cry}\index{RGG-Com-cry@(RGG-Com-cry)}
		For all distinct, non-proportional roots $ \alpha, \beta $, we have
		\[ \commutator{\rootgr{\alpha}}{\rootgr{\beta}} \subs \rootgr{\rootint{\alpha}{\beta}}. \]
	\end{description}
\end{definition}

This definition merits several comments.

\begin{note}[Commutators for opposite and equal roots]\label{rgg:rootgr-abelian}
	We have no commutator relations for $ \commutator{\rootgr{\alpha}}{\rootgr{-\alpha}} $ and $ \commutator{\rootgr{\alpha}}{\rootgr{\alpha}} $ unless $ 2\alpha $ is a root, in which case $ \commutator{\rootgr{\alpha}}{\rootgr{\alpha}} \subs \rootgr{2\alpha} $. However, if the additional axioms of a $ \roots $-graded group (see \cref{rgg-def}) are satisfied, the latter problem disappears: We will show (in \cref{ADE:abelian,B:abelian,C:abelian}) that for all irreducible crystallographic root systems $ \roots $ of rank at least~3 as well as for $ A_2 $, any $ \roots $-graded group satisfies $ \commutator{\rootgr{\alpha}}{\rootgr{\alpha}} = \compactSet{1_G} $ (which means that $ \rootgr{\alpha} $ is abelian) for all roots $ \alpha $ for which $ 2\alpha $ is not a root. The same holds for root systems of type $ H $ because every root in $ H_3 $ is contained in a parabolic subsystem of type $ A_2 $.
\end{note}

\begin{note}[The crystallographic assumption]\label{cry-note}
	The notion of crystallographic commutator relations can be defined even when $ \roots $ is not crystallographic. In practice, however, we will not consider crystallographic commutator relations if $ \roots $ is not crystallographic. Conversely, we will not consider non-crystallographic commutator relations if $ \roots $ is crystallographic. The only exception to this rule will be \cref{sec:B:rank2-noncry}, where we study non-crystallographic $ B_2 $-gradings. However, this is not due to a genuine interest in non-crystallographic $ B_n $-gradings, but rather motivated by the fact that certain results about crystallographic $ BC_2 $-gradings can be more easily formulated in the language of non-crystallographic $ B_2 $-gradings.
	
	Recall from \cref{rootsys:interval-note} that if $ \roots $ is crystallographic, then the extra assumption of crystallographic commutator relations says precisely that for some pairs $ (\alpha, \gamma) $ for roots with $ \rootintcox{\alpha}{\gamma} = \compactSet{\beta} $ for some root $ \beta $, we have the trivial commutator relation $ \commutator{\rootgr{\alpha}}{\rootgr{\gamma}} = \compactSet{1} $ instead of the commutator relation $ \commutator{\rootgr{\alpha}}{\rootgr{\gamma}} \subs \rootgr{\beta} $.
\end{note}

\begin{note}[Commutator relations in the sense of Loos-Neher]\label{rgg:commrel:loos-neher-note}
	The specific terminology of \enquote{groups with $ \roots $-commutator relations} appears for the first time in \cite[\onpage{183}]{Faulkner-StableRange}. It also appears in \cite{LoosNeherBook} by Loos-Neher, who cite \cite{Faulkner-RGD} as their inspiration (see \cite[\onpage{72}]{LoosNeherBook}). However, our \cref{rgg:group-commrel-def,rgg:group-commrel-cry-def} differ from the one in \cite[\onpage{23}]{LoosNeherBook} in a few aspects. First of all, while we restrict our attention to root systems $ \roots $, Loos-Neher define arbitrary $ R $-gradings where $ R $ is any subset of a free abelian group $ X $ such that $ R $ contains $ 0 $ and spans $ X $. Many of the facts and concepts that we introduce in this chapter are studied in this broad generality in \cite[Chapter~1]{LoosNeherBook}. In particular, $ R $ is allowed to be infinite.
	
	However, when we restrict our attention to reduced root systems, we see that our definition is actually less restrictive than the one in \cite{LoosNeherBook}. For once, Neher-Loos only consider what we call crystallographic commutator relations. Further, the commutator axiom in \cite{LoosNeherBook} says that $ \commutator{\rootgr{\alpha}}{\rootgr{\beta}} \subs \rootgr{\rootint{\alpha}{\beta}} $ for any so-called nilpotent pair $ (\alpha, \beta) $, which includes the case $ \alpha = \beta $. In our notation, the interval $ \rootint{\alpha}{\alpha} $ is not defined (see \cref{rootsys:interval-def}), but in the notation of \cite{LoosNeherBook}, we have $ \rootint{\alpha}{\alpha} = \IN_{\ge 2} \alpha \intersect \roots $. For reduced root systems, this set is empty. Thus the commutator axioms in \cite{LoosNeherBook} require that each root group is abelian if $ \roots $ reduced, which (for the root systems we are interested in) is redundant by \cref{rgg:rootgr-abelian}.
\end{note}

\begin{note}[Reducedness]
	The axiom \axiomref{rgg-axiom-div} can be omitted if $ \roots $ is reduced. Recall from \cref{rootsys:class} that the only crystallographic root system which is not reduced is $ BC_n $. In fact, \axiomref{rgg-axiom-div} is specifically tailored for this root system. For root systems which are neither reduced nor crystallographic, the axioms have to be modified.\mywarning{mention definition in H3-paper}
\end{note}

\begin{note}[Commutator formulas]\label{rgg:local-global}
	Axioms~\axiomref{rgg-axiom-comm} and~\axiomref{rgg-axiom-comm-cry} can be seen as \enquote{global} conditions: They tell us what happens on the level on root subgroups. However, they provide no \enquote{local} information: They do not give a formula for the commutator of two specific root group elements $ x_\alpha \in \rootgr{\alpha} $ and $ x_\beta \in \rootgr{\beta} $. Our main goal for root graded groups is to find such formulas for all pairs $ (\alpha, \beta) $. See \cref{rgg:coord-gen-def} for more details.
\end{note}

\begin{note}[The rank-1 case]\label{rgg:A1-is-special}
	Assume that $ \roots = A_1 $. In other words, $ \roots $ consists of exactly two roots $ \alpha $ and $ -\alpha $. Then there exist no non-proportional roots in $ \roots $, so that every $ A_1 $-pregrading has $ A_1 $-commutator relations. For this reason, root gradings of type $ A_1 $ (or of type $ A_1 \times \cdots \times A_1 $) are much more difficult to understand than root gradings for other root systems.
	
	As a consequence, we often have to assume that every root in $ \roots $ is contained in a root subsystem which is not of type $ A_1^n $ for some $ n \in \Npos $. Examples in which this assumption is explicit are \cref{weyl:1-not-invertible,basic:weyl-det-by-one-factor}. In many other situations, this assumption is not made explicit because it is always satisfied for irreducible root systems of rank at least~$ 2 $.
\end{note}

Closed root subsystems induce subgroups with commutator relations in a natural way. We will extend this result to root graded groups in \cref{rgg:subgroup}.

\begin{lemma}\label{basic:commrel-subgroup}
	Let $ G $ be a group which has $ \roots $-commutator relations with root groups $ (\rootgr{\alpha})_{\alpha \in \roots} $ and let $ \roots' $ be a subset of $ \roots $. Denote by $ H $ the group which is generated by $ (\rootgr{\alpha})_{\alpha \in \roots'} $. If $ \roots' $ is a closed root subsystem, then $ H $ has $ \roots' $-commutator relations with root groups $ (\rootgr{\alpha})_{\alpha \in \roots'} $. If $ \roots' $ is a crystallographically closed root subsystem and the commutator relations of $ G $ are crystallographic, then $ H $ has crystallographic $ \roots' $-commutator relations with root groups $ (\rootgr{\alpha})_{\alpha \in \roots'} $.
\end{lemma}
\begin{proof}
	This follows from the fact that the interval between $ \alpha, \beta \in \roots' $ does not depend on whether we consider $ \alpha, \beta $ as elements of $ \roots $ or as elements of $ \roots' $.
\end{proof}

\subsection{Commutator and Product Maps}

We will often be faced with situations in which we want to decompose a commutator $ \commutator{x_1 \cdots x_n}{y_1 \cdots y_m} $, where $ \listing{x}{n} \in \rootgr{\alpha} $ and $ \listing{y}{m} \in \rootgr{\gamma} $ for non-proportional roots $ \alpha, \gamma $, into a product of of more \enquote{basic} elements, for example into commutators $ \commutator{x_i}{y_j} $. Such a decomposition is always possible, but its nature depends on the number of roots in $ \rootintcox{\alpha}{\gamma} $. If $ \rootintcox{\alpha}{\gamma} $ is empty, then $ \commutator{\listing{x}{n}}{\listing{y}{m}} = 1_G $. The following \cref{basic:comm-add} describes the situation in which $ \indivset{\rootintcox{\alpha}{\gamma}} $ contains exactly one element. This is always the case in simply-laced root systems, and even in the non-simply-laced case, we can find many pairs of roots $ (\alpha, \gamma) $ which satisfy this condition.

\begin{lemma}\label{basic:comm-add}
	Let $ G $ be a group which has $ \roots $-commutator relations with root groups $ (\rootgr{\alpha})_{\alpha \in \roots} $ and let $ \alpha, \gamma $ be roots such that $ \indivset{\rootintcox{\alpha}{\gamma}} $ has exactly one element. Then we have 
	\[ \commutator{x_\alpha y_\alpha}{x_\gamma} = \commutator{x_\alpha}{x_\gamma} \commutator{y_\alpha}{x_\gamma} \midand \commutator{x_\alpha}{x_\gamma y_\gamma} = \commutator{x_\alpha}{y_\gamma} \commutator{x_\alpha}{x_\gamma} \]
	for all $ x_\alpha, y_\alpha \in \rootgr{\alpha} $ and $ x_\gamma, y_\gamma \in \rootgr{\gamma} $. In particular, 
	\[ \commutator{x_\alpha^{-1}}{x_\gamma} = \commutator{x_\alpha}{x_\gamma}^{-1} = \commutator{x_\alpha}{x_\gamma^{-1}} \]
	for all $ x_\alpha \in \rootgr{\alpha} $ and $ x_\gamma \in \rootgr{\gamma} $.
\end{lemma}
\begin{proof}
	Denote the unique root in $ \indivset{\rootintcox{\alpha}{\gamma}} $ by $ \beta $ and let $ x_\alpha, y_\alpha \in \rootgr{\alpha} $, $ x_\gamma, y_\gamma \in \rootgr{\gamma} $. Then \thmitemcref{group-rel}{group-rel:add} says that
	\begin{align*}
		\commutator{x_\alpha y_\alpha}{x_\gamma} = \commutator{x_\alpha}{x_\gamma}^{y_\alpha} \commutator{y_\alpha}{x_\gamma} \midand \commutator{x_\alpha}{x_\gamma y_\gamma} = \commutator{x_\alpha}{y_\gamma} \commutator{x_\alpha}{x_\gamma}^{y_\alpha}.
	\end{align*}
	Observe that $ \commutator{x_\alpha}{x_\gamma} $ lies in $ \rootgr{\beta} $ by \axiomref{rgg-axiom-comm} and that $ \alpha $ is adjacent to $ \beta $. It follows that $ y_\alpha $ commutes with $ \commutator{x_\alpha}{x_\gamma} $, so the first assertion follows. In particular, we have
	\begin{align*}
		1_G &= \commutator{1_G}{x_\gamma} = \commutator{x_\alpha x_\alpha^{-1}}{x_\gamma} = \commutator{x_\alpha}{x_\gamma} \commutator{x_\alpha^{-1}}{x_\gamma},
	\end{align*}
	so $ \commutator{x_\alpha}{x_\gamma}^{-1} = \commutator{x_\alpha^{-1}}{x_\gamma} $. Similarly, we can see that $ \commutator{x_\alpha}{x_\gamma}^{-1} = \commutator{x_\alpha}{x_\gamma^{-1}} $.
\end{proof}

\begin{remark}\label{basic:comm-add:abel}
	\Cref{basic:comm-add} has the minor flaw that the order of the commutators in the formula $ \commutator{x_\alpha}{x_\gamma y_\gamma} = \commutator{x_\alpha}{y_\gamma} \commutator{x_\alpha}{x_\gamma} $ is reversed. In practice, however, it will always be the case that either $ \rootgr{\beta} $ or $ \rootgr{\gamma} $ is abelian, which implies that $ \commutator{x_\alpha}{x_\gamma y_\gamma} = \commutator{x_\alpha}{x_\gamma} \commutator{x_\alpha}{y_\gamma} $.
\end{remark}

If $ \indivset{\rootint{\alpha}{\gamma}} $ contains more than one root, then the situation is more delicate. Before we can consider this, we have to introduce some notation. 

\begin{notation}[Product map]
	Let $ G $ be a group with a $ \roots $-pregrading $ (\rootgr{\alpha})_{\alpha \in \roots} $, let $ k $ be a positive integer and let $ \alpha_1, \ldots, \alpha_k $ be pairwise distinct roots. The \defemph*{product map on $ \rootgr{\alpha_1} \times \cdots \times \rootgr{\alpha_k} $}\index{product map} (or simply the \defemph*{product map on $ (\alpha_1, \ldots, \alpha_k) $}) is the map
	\[ \map{}{\rootgr{\alpha_1} \times \cdots \times \rootgr{\alpha_k}}{\gen{\rootgr{\alpha_1} \union \cdots \union \rootgr{\alpha_k}}}{(g_1, \ldots, g_k)}{g_1 \cdots g_k}. \]
\end{notation}

\begin{definition}\label{basic:commpart-def}
	Let $ G $ be a group which has $ \roots $-commutator relations with root groups $ (\rootgr{\alpha})_{\alpha \in \roots} $. Let $ \alpha, \beta $ be non-proportional, indivisible roots and denote by $ (\alpha_1, \ldots, \alpha_k) $ the unique interval ordering of $ \indivset{\clrootintcox{\alpha}{\beta}} $ such that $ \alpha_1 = \alpha $ and $ \alpha_k = \beta $ (see \cref{rootsys:clockwise-rem}). Assume that the product map on $ \rootgr{\alpha_2} \times \cdots \times \rootgr{\alpha_{k-1}} $ is bijective. (We will see in \cref{prodmap:commpart-welldef} that this assumption is always satisfied if $ G $ is rank-2-injective in the sense of \cref{basic:rank2inj-def}.) Then we denote by
	\[ \brackets[\big]{\map{\commpart{}{}{\alpha_i}}{\rootgr{\alpha} \times \rootgr{\beta}}{\rootgr{\alpha_i}}{(x_\alpha, x_\beta)}{\commpart{x_\alpha}{x_\beta}{\alpha_i}}}_{i \in \numint{2}{k-1}} \]
	the unique family of maps with the property that
	\[ \commutator{x_\alpha}{x_\beta} = \commpart{x_\alpha}{x_\beta}{\alpha_2} \cdots \commpart{x_\alpha}{x_\beta}{\alpha_{k-1}} \]
	for all $ x_\alpha \in \rootgr{\alpha} $, $ x_\beta \in \rootgr{\beta} $. In other words, for each $ i \in \numint{2}{k-1} $ we denote by $ \commpart{x_\alpha}{x_\beta}{\alpha_i} $ the unique element of $ \rootgr{\alpha_i} $ with the property that
	\begin{equation}\label{eq:basic:commpart-def}
		\commutator{x_\alpha}{x_\beta} \in \rootgr{\alpha_2} \cdots \rootgr{\alpha_{i-1}} \commpart{x_\alpha}{x_\beta}{\alpha_i} \rootgr{\alpha_{i+1}} \cdots \rootgr{\alpha_{k-1}}.
	\end{equation}
\end{definition}

The surjectivity of the product map in \cref{basic:commpart-def} guarantees the existence of an element $ \commpart{x_\alpha}{x_\beta}{\alpha_i} $ with the property~\eqref{eq:basic:commpart-def}, and the injectivity of the product map guarantees that it is uniquely determined. Before we can continue, there is a subtlety in \cref{basic:commpart-def} that we have to point out. 
	
\begin{note}[on \cref{basic:commpart-def}]\label{basic:commpart-def-note}
	We emphasise that $ \commpart{}{}{\alpha_i} $ is a map on the product $ \rootgr{\alpha} \times \rootgr{\beta} $ (on the two arguments $ x_\alpha $ and $ x_\beta $), not a map on $ \rootgr{\rootintcox{\alpha}{\beta}} $ (on the argument $ \commutator{x_\alpha}{x_\beta} $). The reason for this is that we fix an interval ordering of $ \indivset{\clrootintcox{\alpha}{\beta}} $ in \cref{basic:commpart-def}, and this is only possible (without ambiguity of choice) if we declare which of the roots $ \alpha, \beta $ should be the first one and which should be the last one. In other words, we cannot (without ambiguity) define a family of maps
	\[ \brackets[\big]{\map{}{\rootgr{\rootintcox{\alpha}{\beta}}}{\rootgr{\alpha_i}}{g}{g_{\alpha_i}}}_{i \in \numint{2}{k-1}} \]
	because it is not clear whether the elements $ g_{\alpha_2}, \ldots, g_{\alpha_{k-1}} $ are defined by the property $ g = g_{\alpha_2} \cdots g_{\alpha_{k-1}} $ or by the property $ g = g_{\alpha_{k-1}} \cdots g_{\alpha_{2}} $.
\end{note}

\begin{remark}\label{basic:commpart-def-cry}
	The notational subtleties of \cref{basic:commpart-def-note} are insubstantial if $ k = 4 $ because the root groups $ \rootgr{\alpha_2} $ and $ \rootgr{\alpha_3} $ commute by \axiomref{rgg-axiom-comm}. Indeed, it follows from this that for all $ g \in \gen{\rootgr{\alpha_2} \union \rootgr{\alpha_3}} $ there exist unique elements $ g_2 \in \rootgr{\alpha_2} $ and $ g_3 \in \rootgr{\alpha_3} $ such that $ g = g_2 g_3 = g_3 g_2 $. Since closed root intervals in crystallographic root systems of rank at least~$ 3 $ contain at most~$ 4 $ indivisible roots, this means that we can largely ignore the implications of \cref{basic:commpart-def-note} in this work.
\end{remark}

The remarks in \cref{basic:commpart-def-note} are crucial for the validity of the following statement.

\begin{lemma}\label{commpart:inv-switch}
	Let $ G $ be a group which has $ \roots $-commutator relations with root groups $ (\rootgr{\alpha})_{\alpha \in \roots} $. Let $ \alpha, \beta $ be non-proportional, indivisible roots and denote by $ (\alpha_1, \ldots, \alpha_k) $ the unique interval ordering of $ \indivset{\clrootintcox{\alpha}{\beta}} $ such that $ \alpha_1 = \alpha $ and $ \alpha_k = \beta $ (see \cref{rootsys:clockwise-rem}). Assume that the product map on $ \rootgr{\alpha_2} \times \cdots \times \rootgr{\alpha_{k-1}} $ is bijective. Then for all $ x_1 \in \rootgr{\alpha_1} = \rootgr{\alpha} $ and $ x_k \in \rootgr{\alpha_k} = \rootgr{\beta} $, we have $ \commpart{x_1}{x_k}{\alpha_i}^{-1} = \commpart{x_k}{x_1}{\alpha_i} $ for all $ i \in \numint{2}{k-1} $.
\end{lemma}
\begin{proof}
	Let $ x_1 \in \rootgr{\alpha_1} $, $ x_k \in \rootgr{\alpha_k} $. At first, observe that the product map on $ \rootgr{\alpha_{k-1}} \times \cdots \times \rootgr{\alpha_2} $ is bijective as well, so $ \commpart{x_k}{x_1}{\alpha_i} $ is well-defined for all $ i \in \numint{2}{k-1} $. By~\thmitemcref{group-rel}{group-rel:inv},
	\begin{align*}
		\commutator{x_k}{x_1} &= \commutator{x_1}{x_k}^{-1} = \brackets[\big]{\commpart{x_1}{x_k}{\alpha_2} \cdots \commpart{x_1}{x_k}{\alpha_{k-1}}}^{-1} = \commpart{x_1}{x_k}{\alpha_{k-1}}^{-1} \cdots \commpart{x_1}{x_k}{\alpha_2}^{-1}.
	\end{align*}
	Since $ (\alpha_k, \ldots, \alpha_1) $ is the unique interval ordering of $ \indivset{\clrootintcox{\alpha}{\beta}} = \indivset{\clrootintcox{\beta}{\alpha}} $ whose first entry is $ \beta $ and whose last entry is $ \alpha $, we infer that $ \commpart{x_1}{x_k}{\alpha_i}^{-1} = \commpart{x_k}{x_1}{\alpha_i} $ for all $ i \in \numint{2}{k-1} $, as desired.
\end{proof}

Using \cref{basic:commpart-def}, we can prove analogues of \cref{basic:comm-add} for the case that $ \indivset{\rootintcox{\alpha}{\gamma}} $ contains more than one element. These formulas involve not only basic commutators $ \commutator{x_i}{y_j} $ but also nested commutators and the maps from \cref{basic:commpart-def}. We will derive them for the case that $ \indivset{\rootintcox{\alpha}{\gamma}} $ contains exactly~2 elements in \cref{B:comm-add,B:comm-add-cry,BC:comm-add}. Since these formulas will be used a lot, it is essential that the bijectivity assumption on the product map in \cref{basic:commpart-def} is always met. We will see in \cref{prodmap:interval-surj} that we automatically have surjectivity, but injectivity need not hold in general. Thus we introduce the following condition.

\begin{definition}[Rank-2-injective]\label{basic:rank2inj-def}
	Let $ G $ be a group. A $ \roots $-pregrading $ (\rootgr{\alpha})_{\alpha \in \roots} $ of $ G $ is said to be \defemph{rank-2-injective} if for all parabolic root subsystems $ \roots' $ of $ \roots $ of rank $ 2 $, all positive systems $ \possys' $ in $ \roots' $ and all interval orderings $ \tup{\alpha}{k} $ of $ \indivset{(\possys')} $, the product map on $ \tup{\alpha}{k} $ is injective.
\end{definition}

The following result is an easy consequence of rank-2-injectivity.

\begin{lemma}\label{basic:prodmap-triv-intersect}
	Let $ G $ be a group with a $ \roots $-pregrading $ (\rootgr{\alpha})_{\alpha \in \roots} $. If $ \listing{\alpha}{k} $ are pairwise distinct roots such that the product map on $ \tup{\alpha}{k} $ is injective, then for all distinct $ i,j \in \numint{1}{k} $, we have $ \rootgr{\alpha_i} \intersect \rootgr{\alpha_j} = \compactSet{1_G} $. In particular, if $ G $ is rank-2-injective, then $ \rootgr{\alpha} \intersect \rootgr{\beta} = \compactSet{1_G} $ for all non-proportional roots $ \alpha, \beta $.
\end{lemma}
\begin{proof}
	Let $ i,j \in \numint{1}{k} $ be distinct and let $ x \in \rootgr{\alpha_i} \intersect \rootgr{\alpha_j} $. Denote by $ \phi_i(x) $ the tuple $ (1_G, \ldots, 1_G, x, 1_G, \ldots, 1_G) $ with $ x $ at position $ i $, and define $ \phi_j(x) $ in a similar way. Then both $ \phi_i(x) $ and $ \phi_j(x) $ are mapped to $ x $ under the product map on $ \tup{\alpha}{k} $. By the injectivity of this map, it follows that $ x = 1_G $, as desired.
\end{proof}

\subsection{Generalised Commutator Relations}

Groups with commutator relations satisfy a commutator relation of the form
\[ \commutator{\rootgr{\alpha}}{\rootgr{\beta}} \subs \rootgr{\rootintcox{\alpha}{\beta}}. \]
The goal of this subsection is to derive a stronger formula of the form
\[ \commutator{\rootgr{A}}{\rootgr{B}} \subs \rootgr{\rootintcox{A}{B}} \]
where $ A $, $ B $ are subsets of $ \roots $ satisfying some conditions (see \cref{basic:gen-com-rel}). This formula is called the \defemph*{generalised commutator relation}.\index{commutator relations!generalised} We closely follow the arguments in \cite[3.9~(a)]{LoosNeherBook}, but we have to be careful because of the differing conventions concerning root intervals (see \cref{rootsys:loosneher-interval-def}). As a consequence of the generalised commutator relations, we obtain that certain maps between root systems preserve commutator relations (\cref{basic:comm-fold}).

\begin{definition}[Commutator set]\label{basic:commset}
	Let $ A,B $ be two subsets of $ \roots $. Then the sets
	\begin{gather*}
		\rootintcox{A}{B} \defl \roots \intersect \Set*{\sum_{i=1}^n \lambda_i \alpha_i + \sum_{j=1}^m \mu_j \beta_j \given \begin{gathered}
			n,m \in \Npos, \\
			\listing{\lambda}{n}, \listing{\mu}{m} \in \IR_{>0}, \\
			\listing{\alpha}{n} \in A, \listing{\beta}{m} \in B
		\end{gathered}}, \\
		\rootint{A}{B} \defl \roots \intersect \Set*{\sum_{i=1}^n \lambda_i \alpha_i + \sum_{j=1}^m \mu_j \beta_j \given \begin{gathered}
			n,m \in \Npos, \listing{\lambda}{n}, \listing{\mu}{m} \in \Npos, \\
			\listing{\alpha}{n} \in A, \listing{\beta}{m} \in B
		\end{gathered}}
	\end{gather*}
	are called the \defemph*{commutator set of $ (A,B) $}\index{commutator set} and the \defemph*{crystallographic commutator set of $ (A,B) $}\index{commutator set!crystallographic}, respectively.
\end{definition}

\begin{remark}[Root intervals as commutator sets]
	Let $ \alpha, \beta $ be non-proportional roots. Then
	\[ \rootint{\alpha}{\beta} = \rootint{\compactSet{\alpha}}{\compactSet{\beta}} \midand \rootintcox{\alpha}{\beta} = \rootintcox{\compactSet{\alpha}}{\compactSet{\beta}}, \]
	so root intervals are a special case of commutator sets. While root intervals are defined only for pairs of non-proportional roots, there is no such restriction for commutator sets. For example, we have
	\begin{align*}
		\rootint{\compactSet{\alpha}}{\compactSet{\alpha}} = \roots \intersect \IN_{\ge 2} \alpha \midand \rootintcox{\compactSet{\alpha}}{\compactSet{\alpha}} = \roots \intersect \IR_{>0} \alpha.
	\end{align*}
	Note that the set $ \rootint{\compactSet{\alpha}}{\compactSet{\alpha}} $ is precisely the set $ \rootint{\alpha}{\alpha} $ in the notation of \cite{LoosNeherBook}, see \cref{rootsys:loosneher-interval-def}.
\end{remark}

\begin{remark}\label{basic:commset:subset}
	Let $ A,A',B,B' $ be subsets of $ \roots $ such that $ A \subs A' $ and $ B \subs B' $. Then $ \rootint{A}{B} \subs \rootint{A'}{B'} $ and $ \rootintcox{A}{B} \subs \rootintcox{A'}{B'} $. Further, if $ C \defl \rootint{A}{B} $, then $ \rootint{A \union B}{C} \subs C $, and similarly for non-crystallographic commutator sets.
\end{remark}

We begin with an auxiliary lemma from group theory.

\begin{lemma}[{\cite[3.7]{LoosNeherBook}}]\label{basic:normalise-lem}
	Let $ G $ be a group, let $ H $ be a subgroup and let $ X_1, X_2 $ be subsets of $ G $ normalising $ H $ (meaning that $ H^x = H $ for all $ x \in X_1 \union X_2 $) such that $ \commutator{X_1}{X_2} \subs H $. Denote by $ G_1 $ and $ G_2 $ the groups generated by $ X_1 $ and $ X_2 $, respectively. Then $ G_1, G_2 $ normalise $ H $ and they satisfy $ \commutator{G_1}{G_2} \subs H $.
\end{lemma}
\begin{proof}
	It is clear that $ X_1^{-1} $ and $ X_2^{-1} $ normalise $ H $ as well. Hence $ G_1 $ and $ G_2 $ normalise $ H $. Further, it follows from \thmitemcref{group-rel}{group-rel:add} that
	\[ 1_G = \commutator{x_1^{-1} x_1}{x_2} = \commutator{x_1^{-1}}{x_2}^{x_1} \commutator{x_1}{x_2} \]
	for all $ x_1 \in X_1 $, $ x_2 \in X_2 $, which implies that
	\[ \commutator{x_1^{-1}}{x_2} = (\commutator{x_1}{x_2}^{-1})^{x_1^{-1}} \in H^{x_1^{-1}} = H. \]
	Thus we also have $ \commutator{X_1^{-1}}{X_2} \subs H $. In a similar way, it can be shown that $ \commutator{X_1}{X_2^{-1}} \subs H $ and $ \commutator{X_1^{-1}}{X_2^{-1}} \subs H $. Therefore, by replacing $ X_1 $ with $ X_1 \union X_1^{-1} $ and $ X_2 $ with $ X_2 \union X_2^{-1} $, we can assume that $ G_1 $ and $ G_2 $ are generated by $ X_1 $ and $ X_2 $ as monoids, respectively.
	
	Now let $ x_1, y_1 \in X_1 $ and $ x_2 \in X_2 $. Then again by \thmitemcref{group-rel}{group-rel:add}, we have $ \commutator{x_1 y_1}{x_2} = \commutator{x_1}{x_2}^{y_1} \commutator{y_1}{x_2} $ where $ \commutator{x_1}{x_2} $ and $ \commutator{y_1}{x_2} $ lie in $ H $ by assumption. Further, since $ X_1 $ normalises $ H $, we also have that $ \commutator{x_1}{x_2}^{y_1} $ lies in $ H $. Thus $ \commutator{x_1 y_1}{x_2} $ lies in $ H $. By a straightforward induction using a similar argument, we conclude that $ \commutator{g_1}{g_2} $ lies in $ H $ for all $ g_1 \in G_1 $ and $ g_2 \in G_2 $, as desired.
\end{proof}

\begin{lemma}[{\cite[3.9]{LoosNeherBook}}]\label{basic:NL-gencommrel-lem}
	Let $ G $ be a group with a $ \roots $-pregrading $ (\rootgr{\alpha})_{\alpha \in \roots} $ and let $ A,B,C $ be subsets of $ \roots $. Assume that for all $ \alpha \in A $, $ \beta \in B $, we have $ \commutator{\rootgr{\alpha}}{\rootgr{\beta}} \subs \rootgr{C} $ and that for all $ \alpha \in A $, $ \beta \in B $, $ \gamma \in C $, we have $ \commutator{\rootgr{\alpha}}{\rootgr{\gamma}} \subs \rootgr{C} $ and $ \commutator{\rootgr{\beta}}{\rootgr{\gamma}} \subs \rootgr{C} $. Then $ \rootgr{A} $ and $ \rootgr{B} $ normalise $ \rootgr{C} $ and we have $ \commutator{\rootgr{A}}{\rootgr{B}} \subs \rootgr{C} $.
\end{lemma}
\begin{proof}
	Put $ X_1 \defl \bigunion_{\alpha \in A} \rootgr{\alpha} $, $ X_2 \defl \bigunion_{\beta \in B} \rootgr{\beta} $ and $ H \defl \rootgr{C} $. The first assumption on commutation relations in $ G $ yields that $ \commutator{X_1}{X_2} \subs H $ and the second one yields that $ X_1 $ and $ X_2 $ normalise $ H $. Thus the assertion follows from \cref{basic:normalise-lem}.
\end{proof}

\begin{remark}[Cones]
	Let $ W $ be a finite-dimensional real vector space and let $ C $ be a subset of $ W $. Then $ C $ is called a \defemph{cone} if it is closed under addition and scalar multiplication with $ \IR_{>0} $. The \defemph*{cone of $ C $}, denoted by $ \cone(C) $, is the smallest cone which contains $ C $. Equivalently, it is the convex hull of $ \IR_{>0} C $, and it is also the additive semigroup generated by $ \IR_{>0} C $. Note that a cone may, but need not contain $ 0_W $.
\end{remark}

\begin{proposition}[Generalised commutator relation, {\cite[3.9]{LoosNeherBook}}]\label{basic:gen-com-rel}
	Let $ G $ be a group with a $ \roots $-pregrading $ (\rootgr{\zeta})_{\zeta \in \roots} $ and let $ A,B $ be subsets of $ \roots $ such that
	\[ \cone(A) \intersect \cone(-B) = \emptyset. \]
	Assume that $ \roots $ is reduced or crystallographic and that for all roots $ \zeta $ for which $ 2\zeta $ is not a root, the root group $ \rootgr{\zeta} $ is abelian. Then the following hold:
	\begin{proenumerate}
		\item If $ G $ has $ \roots $-commutator relations with root groups $ (\rootgr{\zeta})_{\zeta \in \roots} $, then
		\[ \commutator{\rootgr{A}}{\rootgr{B}} \subs \rootgr{\rootintcox{A}{B}}. \]
		
		\item If $ G $ has crystallographic $ \roots $-commutator relations with root groups $ (\rootgr{\zeta})_{\zeta \in \roots} $, then
		\[ \commutator{\rootgr{A}}{\rootgr{B}} \subs \rootgr{\rootint{A}{B}}. \]
	\end{proenumerate}
\end{proposition}
\begin{proof}
	We only consider the case that $ G $ has $ \roots $-commutator relations with root groups $ (\rootgr{\zeta})_{\zeta \in \roots} $. The other assertion can be proven in the same way. Recall from \cref{rootsys:cry-nonred-2} that the only positive multiples of a root $ \alpha $ which can possibly be a root are $ \frac{1}{2}\alpha $, $ \alpha $ and $ 2\alpha $. Our goal is to verify the conditions in \cref{basic:NL-gencommrel-lem} for $ C \defl \rootintcox{A}{B} $. By the assumption that $ \roots $ is reduced or crystallographic, we have $ \IR_{>0} \alpha \intersect \roots \subs \Set{\frac{1}{2} \alpha, \alpha, 2\alpha} $ for all $ \alpha \in \roots $ (see \cref{rootsys:cry-nonred-2}), and this is the only way in which we will use this assumption.
	
	Let $ \alpha \in A $ and $ \beta \in B $. If $ \alpha, \beta $ are non-proportional, then it follows from the commutator relations of $ G $ and \cref{basic:commset:subset} that
	\[ \commutator{\rootgr{\alpha}}{\rootgr{\beta}} \subs \rootgr{\rootintcox{\alpha}{\beta}} \subs \rootgr{C}. \]
	If $ \beta \in \Set{\frac{1}{2} \alpha, 2 \alpha} $, then $ \commutator{\rootgr{\alpha}}{\rootgr{\beta}} = \compactSet{1_G} $ by Axiom~\axiomref{rgg-axiom-div}. If $ \alpha = \beta $, then $ \commutator{\rootgr{\alpha}}{\rootgr{\beta}} \subs \rootgr{2\alpha} $ if $ 2\alpha $ is a root (by Axiom~\axiomref{rgg-axiom-div}) and $ \commutator{\rootgr{\alpha}}{\rootgr{\beta}} = \compactSet{1_G} $ if $ 2\alpha $ is not a root (by the assumption that $ \rootgr{\alpha} $ is abelian in this case). Further, the case $ \beta \in \Set{-\frac{1}{2} \alpha, -\alpha, -2\alpha} $ cannot occur because $ A \intersect (-B) \subs \cone(A) \intersect \cone(-B) = \emptyset $. Thus we have $ \commutator{\rootgr{\alpha}}{\rootgr{\beta}} \subs \rootgr{C} $ for all $ \alpha \in A $ and $ \beta \in B $.
	
	Now let $ \alpha \in A $ and $ \gamma \in C = \rootintcox{A}{B} $. Again, if $ \alpha $ and $ \gamma $ are non-proportional, then
	\[ \commutator{\rootgr{\alpha}}{\rootgr{\gamma}} \subs \rootgr{\rootintcox{\alpha}{\gamma}} \subs \rootgr{\rootintcox{A}{C}} \subs \rootgr{C} \]
	by \cref{basic:commset:subset}. If $ \alpha = \gamma $, then
	\[ \commutator{\rootgr{\alpha}}{\rootgr{\gamma}} \subs \begin{cases}
			\compactSet{1_G} & \text{if } 2\alpha \nin \roots, \\
			\rootgr{2\alpha} & \text{ if } 2\alpha \in \roots
		\end{cases} \]
	by the same argument as in the previous paragraph. If $ \gamma \in \Set{\frac{1}{2}\alpha, 2\alpha} $, then $ \commutator{\rootgr{\alpha}}{\rootgr{\gamma}} = \compactSet{1_G} $. Hence $ \commutator{\rootgr{\alpha}}{\rootgr{\gamma}} \subs \rootgr{C} $ in any case. Now suppose for a contradiction that $ \gamma = -a \alpha $ for some $ a \in \IR_{>0} $. By the definition of $ C $, there exist $ n,m \in \Npos $, $ \listing{\lambda}{n}, \listing{\mu}{m} \in \IR_{>0} $, $ \listing{\alpha}{n} \in A $ and $ \listing{\beta}{m} \in B $ such that $ \gamma = \sum_{i=1}^n \lambda_i \alpha_i + \sum_{j=1}^m \mu_j \beta_j $. Hence
	\[ a\alpha + \sum_{i=1}^n \lambda_i \alpha_i = -\sum_{j=1}^m \mu_j \beta_j \in \cone(A) \intersect \cone(-B) = \emptyset, \]
	which is impossible. We conclude that $ \commutator{\rootgr{\alpha}}{\rootgr{\gamma}} \subs \rootgr{C} $ for all $ \alpha \in A $ and $ \gamma \in C $. By interchanging the roles of $ A $ and $ B $, we observe that $ \commutator{\rootgr{\beta}}{\rootgr{\gamma}} \subs \rootgr{C} $ for all $ \beta \in B $ and $ \gamma \in C $ as well. Hence the conditions of \cref{basic:NL-gencommrel-lem} are satisfied, and it follows that $ \commutator{\rootgr{A}}{\rootgr{B}} \subs \rootgr{\rootintcox{A}{B}} $.
\end{proof}

\begin{remark}\label{basic:possys-cone-intersect}
	Let $ \possys $ be a positive system in $ \roots $. Then $ \cone(\possys) \intersect \cone(-\possys) = \emptyset $. It follows that for any subsets $ A $, $ B $ of $ \possys $, the assumption in \cref{basic:gen-com-rel} that $ \cone(A) \intersect \cone(-B) = \emptyset $ is satisfied.
\end{remark}

Recall from \cref{rgg:rootgr-abelian} that for $ \roots $-graded groups, the assumption in \cref{basic:gen-com-rel} that certain root groups are abelian is automatically satisfied if $ \roots $ is of rank at least~3. However, this is a non-trivial fact which we still have to prove.

We end this section with a corollary of \cref{basic:gen-com-rel}, which will be used in \cref{sec:RGG-fold} to construct foldings of root graded groups. Before we can state it, we need to introduce some terminology.

\begin{definition}
	Let $ \roots' $ be another root system and let $ \map{\pi}{\roots}{\roots'}{}{} $ be any map. We say that $ \pi $ is \defemph*{compatible with commutator sets} if $ \pi(\rootintcox{A}{B}) \subs \rootintcox{\pi(A)}{\pi(B)} $ for all $ A,B \subs \roots $, and we say that it is \defemph*{compatible with crystallographic commutator sets} if $ \pi(\rootint{A}{B}) \subs \rootint{\pi(A)}{\pi(B)} $ for all $ A,B \subs \roots $.
\end{definition}

\begin{example}
	Let $ V,V' $ be the Euclidean space surrounding root systems $ \roots, \roots' $. If $ \map{\pi}{V}{V'}{}{} $ is linear and $ \pi(\roots) \subs \roots' $, then $ \map{\restrict{\pi}{\roots}}{\roots}{\roots'}{}{} $ is compatible with commutator sets and with crystallographic commutator sets. Similarly, if $ V=V' $ and $ \pi $ is the \defemph*{scaling map} $ \map{\pi}{V}{V'}{v}{\frac{1}{\norm{v}} v} $, then $ \roots'' \defl \pi(\roots) $ is a reduced root system (with the same Weyl group as $ \roots $) and $ \map{\restrict{\pi}{\roots}}{\roots}{\roots'}{}{} $ is compatible with commutator sets and crystallographic commutator sets.
\end{example}

\begin{proposition}\label{basic:comm-fold}
	Let $ \roots' $ be another root system and let $ \map{\pi}{\roots}{\roots'}{}{} $ be a surjective map which is compatible with (crystallographic) commutator sets. Let $ G $ be a group which has (crystallographic) $ \roots $-commutator relations with root groups $ (\rootgr{\alpha})_{\alpha \in \roots} $ and put
	\[ \varrootgr{\alpha'} \defl \gen{\rootgr{\alpha} \given \alpha \in \roots, \pi(\alpha) = \alpha'} \]
	for all $ \alpha' \in \roots' $. Assume that
	for all non-proportional $ \alpha', \beta' \in \roots' $, there exists a positive system $ \possys $ in $ \roots $ which contains $ \roots \intersect (\pi^{-1}(\alpha') \union \pi^{-1}(\beta')) $. Assume further that $ \roots $ is reduced or crystallographic and that for all roots $ \zeta $ for which $ 2\zeta $ is not a root, the root group $ \rootgr{\zeta} $ is abelian. Then $ G $ satisfies \axiomref{rgg-axiom-comm} (respectively, \axiomref{rgg-axiom-comm-cry}) with respect to $ (\varrootgr{\alpha'})_{\alpha' \in \roots'} $. In particular, $ G $ has (crystallographic) $ \roots' $-commutator relations with root groups $ (\varrootgr{\alpha'})_{\alpha' \in \roots'} $ if $ \roots' $ is reduced.
\end{proposition}
\begin{proof}
	We prove only the crystallographic assertion. The proof for the non-crystallographic assertion is similar. Let $ \alpha', \beta' $ be non-proportional roots in $ \roots' $. Put $ A \defl \pi^{-1}(\alpha') $ and $ B \defl \pi^{-1}(\beta') $. By assumption, there exists a positive system containing $ A \union B $, which by \cref{basic:possys-cone-intersect} implies that $ \cone(A) \intersect \cone(-B) = \emptyset $. Thus it follows from \cref{basic:gen-com-rel} that
	\[ \commutator{\varrootgr{\alpha'}}{\varrootgr{\beta'}} = \commutator{\rootgr{A}}{\rootgr{B}} \subs \rootgr{\rootint{A}{B}}. \]
	For any $ \gamma \in \roots $, we have $ \rootgr{\gamma} \subs \varrootgr{\pi(\gamma)} $, so $ \rootgr{\rootint{A}{B}} \subs \varrootgr{\pi(\rootint{A}{B})} $. Since
	\[ \pi(\rootint{A}{B}) \subs \rootint{\pi(A)}{\pi(B)} = \rootint{\compactSet{\alpha'}}{\compactSet{\beta'}} = \rootint{\alpha'}{\beta'}, \]
	we infer that $ \commutator{\varrootgr{\alpha'}}{\varrootgr{\beta'}} \subs \varrootgr{\rootint{\alpha'}{\beta'}} $. Hence $ G $ satisfies \axiomref{rgg-axiom-comm-cry} with respect to $ (\varrootgr{\alpha'})_{\alpha' \in \roots'} $. If $ \roots' $ is reduced, then \axiomref{rgg-axiom-div} is trivially satisfied, so we conclude that $ G $ has crystallographic $ \roots' $-commutator relations with root groups $ (\varrootgr{\alpha'})_{\alpha' \in \roots'} $ in this case.
\end{proof}

\begin{note}
	The assertion of \cref{basic:comm-fold} remains true if $ \roots $ and $ \roots' $ are merely subsets of a Euclidean space and not assumed to be root systems. (Here the notion of \enquote{positive systems in $ \roots $ and $ \roots' $} is defined exactly as in~\ref{rootsys:possys-def}.)
\end{note}


\section{Weyl Elements}

\label{sec:weyl}

\begin{secnotation}
	Unless otherwise specified, we denote by $ \roots $ an arbitrary root system and by $ G $ a group which has $ \roots $-commutator relations with root groups $ (\rootgr{\alpha})_{\alpha \in \roots.} $
\end{secnotation}

In this section, we study Weyl elements. The second axiom of root graded groups will be that $ \alpha $-Weyl elements exist for all roots $ \alpha $.

\subsection{Definition and Basic Observations}

\begin{definition}[Weyl elements, Weyl triples]\label{weyl:weyl-def}
	Let $ G $ be a group with a $ \roots $-pregrading $ (\rootgr{\alpha})_{\alpha \in \roots} $.
	\begin{defenumerate}
		\item Let $ \alpha $ be any root. An \defemph*{$ \alpha $-Weyl element (in $ G $)}\index{Weyl element} is an element $ w_\alpha \in \rootgr{-\alpha} \rootgr{\alpha} \rootgr{-\alpha} $ with the property that $ \rootgr{\beta}^{w_\alpha} = \rootgr{\refl{\alpha}(\beta)} $ for all roots $ \beta $. The set of all $ \alpha $-Weyl elements is denoted by $ \weylset{\alpha} $.
		
		\item A \defemph{Weyl element} is an $ \alpha $-Weyl element for some root $ \alpha $.
		
		\item Let $ \alpha $ be any root. An \defemph*{$ \alpha $-Weyl triple}\index{Weyl triple} is a triple $ (a_{-\alpha}, b_\alpha, c_{-\alpha}) $ consisting of $ a_{-\alpha}, c_{-\alpha} \in \rootgr{-\alpha} $ and $ b_\alpha \in \rootgr{\alpha} $ such that the product $ a_{-\alpha} b_\alpha c_{-\alpha} $ is an $ \alpha $-Weyl element, called the \defemph*{Weyl element corresponding to $ (a_{-\alpha}, b_\alpha, c_{-\alpha}) $}.
		
		\item A \defemph{Weyl triple} is an $ \alpha $-Weyl triple for some root $ \alpha $.
		
		\item \label{weyl:weyl-def:invertible}Let $ \alpha $ be any root. An element $ b_\alpha \in \rootgr{\alpha} $ is called \defemph*{$ \alpha $-invertible}\index{invertible element in a root group} if there exists elements $ a_{-\alpha}, c_{-\alpha} \in \rootgr{-\alpha} $ such that $ (a_{-\alpha}, b_\alpha, c_{-\alpha}) $ is an $ \alpha $-Weyl triple. We denote the set of $ \alpha $-invertible elements of $ \rootgr{\alpha} $ by $ \invset{\alpha} $. We will also write \defemph{Weyl-invertible} in place of $ \alpha $-invertible if the root $ \alpha $ is not specified.
		
		\item \label{weyl:weyl-def:rank2}If $ \roots $ is of rank 2 and $ \map{}{\IZ}{\roots}{i}{\alpha_i} $ is a fixed rank-2 labeling of $ \roots $ (in the sense of \cref{rootsys:rank2-label-def}), then we put $ \invset{i} \defl \invset{\alpha_i} $ and $ \weylset{i} \defl \weylset{\alpha_i} $ for all $ i \in \IZ $.
	\end{defenumerate}
\end{definition}

\begin{note}
	The terminology of \enquote{Weyl elements} appears for the first time in \cite[\onpage{184}]{Faulkner-StableRange}, though of course these elements had already been considered before in the contexts of Chevalley groups, algebraic groups and RGD-systems. Many results about Weyl elements in a more general setting can be found in \cite[Section~5]{LoosNeherBook}.
\end{note}

\begin{note}
	The notion of Weyl elements does not appear to be useful in groups without $ \roots $-commutator relations. However, the way we have phrased \cref{weyl:weyl-def} allows us to discuss Weyl elements in groups which have not yet been proven to have commutator relations.
\end{note}

Before we begin a proper investigation of Weyl elements, we observe that under reasonable assumptions, the trivial element $ 1_G $ is never $ \alpha $-invertible. These assumptions will always be satisfied if $ \roots $ is irreducible of rank at least~2 and the root groups form a root grading in the sense of \cref{rgg-def}.

\begin{lemma}\label{weyl:1-not-invertible}
	Let $ G $ be a group which has $ \roots $-commutator relations with root groups $ (\rootgr{\alpha})_{\alpha \in \roots} $ and let $ \alpha $ be any root. Assume that there exists a root $ \beta $ such that $ (\alpha, \beta) $ is the root base of a closed rank-2 subsystem of $ \roots $ which is not of type $ A_1 \times A_1 $. (In other words, we assume that $ \alpha $ is contained in a closed rank-2 root subsystem $ \roots' $ of $ \roots $ which is not of type $ A_1 \times A_1 $.) Assume further that $ \rootgr{\beta} \intersect \rootgr{\refl{\alpha}(\beta)} = \compactSet{1_G} $ and that $ \rootgr{\beta} \ne \compactSet{1_G} $. Then $ 1_G \nin \invset{\alpha} $.
\end{lemma}
\begin{proof}
	Assume that $ 1_G $ is contained in $ \invset{\alpha} $. Then there exist $ a_{-\alpha}, c_{-\alpha} \in \rootgr{-\alpha} $ such that $ a_{-\alpha} c_{-\alpha} $ is an $ \alpha $-Weyl element. In other words, there exists an $ \alpha $-Weyl element $ w_\alpha $ which is contained in $ \rootgr{-\alpha} $. Choose an element $ x_\beta \in \rootgr{\beta} \setminus \compactSet{1_G} $. Then $ x_\beta^{w_\alpha} \in \rootgr{\refl{\alpha}(\beta)} $. At the same time, since $ \beta $ is adjacent to $ -\alpha $, we have $ x_\beta^{w_\alpha} = x_\beta \in \rootgr{\beta} $. By our assumptions, it follows that $ x_\beta = 1_G $, which is a contradiction.
\end{proof}

The following properties of Weyl elements in \cref{basic:weyl-general} are elementary and will be used all the time. Interestingly, their proof does not rely on the $ \roots $-commutator relations of $ G $ at all. We will prove more specific variations of some of these statements in \cref{basic:weylmap-form1,basic:weylmap-form2} for groups which have \enquote{unique Weyl extensions}.

\begin{proposition}\label{basic:weyl-general}
	Let $ G $ be a group with a $ \roots $-pregrading $ (\rootgr{\alpha})_{\alpha \in \roots} $. Then the following statements hold for all roots $ \alpha $:
	\begin{proenumerate}
		\item \label{basic:weyl-general:inv}If $ (a_{-\alpha}, b_\alpha, c_{-\alpha}) $ is an $ \alpha $-Weyl triple, then $ (c_{-\alpha}^{-1}, b_\alpha^{-1}, a_{-\alpha}^{-1}) $ is an $ \alpha $-Weyl triple as well (with corresponding Weyl element $ w_\alpha^{-1} $). In particular, the sets $ \weylset{\alpha} $ and $ \invset{\alpha} $ are stable under group inversion.
		
		\item \label{basic:weyl-general:conj}Let $ w_\beta $ be a $ \beta $-Weyl element for some root $ \beta $. If $ (a_{-\alpha}, b_\alpha, c_{-\alpha}) $ is an $ \alpha $-Weyl triple with corresponding Weyl element $ w_\alpha $, then $ (a_{-\alpha}^{w_\beta}, b_\alpha^{w_\beta}, c_{-\alpha}^{w_\beta}) $ is a $ \refl{\beta}(\alpha) $-Weyl triple with corresponding Weyl element $ w_\alpha^{w_\beta} $. In particular, $ \weylset{\alpha}^{w_\beta} = \weylset{\refl{\beta}(\alpha)} $ and $ (\invset{\alpha})^{w_\beta} = \invset{\refl{\beta}(\alpha)} $.
		
		\item \label{basic:weyl-general:minus}If $ (a_{-\alpha}, b_\alpha, c_{-\alpha}) $ is an $ \alpha $-Weyl triple with corresponding Weyl element $ w_\alpha \defl a_{-\alpha} b_\alpha c_{-\alpha} $, then $ (c_{-\alpha}^{w_\alpha^{-1}}, a_{-\alpha}, b_\alpha) $ and $ (b_\alpha, c_{-\alpha}, a_{-\alpha}^{w_\alpha}) $ are $ (-\alpha) $-Weyl triples and all three Weyl triples have the same corresponding Weyl element. In particular, $ \weylset{\alpha} = \weylset{-\alpha} $.
		
		\item \label{basic:weyl-general:a-c-inv}If $ (a_{-\alpha}, b_\alpha, c_{-\alpha}) $ is an $ \alpha $-Weyl triple, then $ a_{-\alpha} $ and $ c_{-\alpha} $ are $ (-\alpha) $-invertible.
		
		\item \label{basic:weyl-general:div}If $ \alpha $ is a root such that $ 2\alpha $ is also a root, then every $ 2\alpha $-Weyl triple is also an $ \alpha $-Weyl triple.
	\end{proenumerate}
\end{proposition}
\begin{proof}
	For the whole proof, we fix an arbitrary root $ \alpha $ and an $ \alpha $-Weyl triple $ (a_{-\alpha}, b_\alpha, c_{-\alpha}) $, whose corresponding Weyl element we denote by $ w_\alpha $. For~\itemref{basic:weyl-general:inv}, let $ \beta $ be another arbitrary root. Since $ w_\alpha $ is an $ \alpha $-Weyl element, we have
	\[ \rootgr{\refl{\alpha}(\beta)}^{w_\alpha} = \rootgr{\beta}, \quad \text{so} \quad \rootgr{\beta}^{w_\alpha^{-1}} = \rootgr{\refl{\alpha}(\beta)}. \]
	Further, we have
	\[ w_\alpha^{-1} = c_{-\alpha}^{-1} b_\alpha^{-1} a_{-\alpha}^{-1} \in \rootgr{-\alpha} \rootgr{\alpha} \rootgr{-\alpha}. \]
	Altogether, we conclude that $ (c_{-\alpha}^{-1}, b_\alpha^{-1}, a_{-\alpha}^{-1}) $ is an $ \alpha $-Weyl triple with corresponding Weyl element $ w_\alpha^{-1} $. In particular, $ b_\alpha^{-1} $ lies in $ \invset{\alpha} $. This proves~\itemref{basic:weyl-general:inv}. Using that $ \refl{\beta} \refl{\alpha} \refl{\beta} = \reflbr{\alpha^{\reflbr{\beta}}} $, assertion~\itemref{basic:weyl-general:conj} can be proven in a similar way.
	
	We turn to~\ref{basic:weyl-general:minus}. A simple computation, using only the definition of conjugation in groups, shows that
	\[ b_\alpha c_{-\alpha} a_{-\alpha}^{w_\alpha} = w_\alpha = c_{-\alpha}^{w_\alpha^{-1}} a_{-\alpha} b_{\alpha}. \]
	Further, $ a_{-\alpha}^{w_\alpha} $ and $ c_{-\alpha}^{w_\alpha^{-1}} $ lie in $ \rootgr{\alpha} $ and $ \refl{-\alpha} = \refl{\alpha} $, so it follows that $ (c_{-\alpha}^{w_\alpha^{-1}}, a_{-\alpha}, b_\alpha) $ and $ (b_\alpha, c_{-\alpha}, a_{-\alpha}^{w_\alpha}) $ are $ (-\alpha) $-Weyl triples, both with corresponding Weyl element $ w_\alpha $. This proves~\ref{basic:weyl-general:minus}, and~\ref{basic:weyl-general:a-c-inv} is a consequence of~\ref{basic:weyl-general:minus}. Finally,~\ref{basic:weyl-general:div} holds because $ \rootgr{\pm 2\alpha} $ is contained in $ \rootgr{\pm\alpha} $ by \axiomref{rgg-axiom-div}.
\end{proof}

\begin{lemma}\label{basic:weyl-ex-basis}
	Assume that $ \roots $ is reduced. Let $ \rootbase $ be a root base of $ \roots $ and assume that there exists a $ \delta $-Weyl element $ w_\delta $ for each $ \delta \in \rootbase $. Then there exists an $ \alpha $-Weyl element for each root $ \alpha $.
\end{lemma}
\begin{proof}
	Let $ \alpha $ be an arbitrary root. By \cref{rootsys:indiv-in-rootbase}, there exist $ \delta, \delta_1, \ldots, \delta_n \in \rootbase $ such that $ \alpha = \delta^{\reflbr{\delta_1 \cdots \delta_n}} $. Then it follows from \thmitemcref{basic:weyl-general}{basic:weyl-general:conj} that $ w_\delta^{w_{\delta_1} \cdots w_{\delta_n}} $ is an $ \alpha $-Weyl element, which finishes the proof.
\end{proof}

\begin{remark}\label{basic:weyl-ex-basis:BC}
	\cref{basic:weyl-ex-basis} and its proof remain valid if $ \roots $ is of type $ BC $ and $ \rootbase $ is the standard rescaled root base (see \cref{BC:BCn-standard-rep}). For proper root bases of non-reduced root systems, however, it is not true: The existence of $ \alpha $-Weyl elements does not imply the existence of $ 2\alpha $-Weyl elements.
\end{remark}

\begin{lemma}\label{basic:commpart-conj}
	Assume that $ G $ is rank-2-injective. Let $ \alpha, \beta, \xi $ be roots, let $ \gamma \in \indivset{\rootintcox{\alpha}{\beta}} $ and let $ w_\xi $ be a $ \xi $-Weyl element for some root $ \xi $. Then we have
	\[ \commpart{x_\alpha}{x_\beta}{\gamma}^{w_\xi} = \commpart{x_\alpha^{w_\xi}}{x_\beta^{w_\xi}}{\refl{\xi}(\gamma)} \]
	for all $ x_\alpha \in \rootgr{\alpha} $ and all $ x_\beta \in \rootgr{\beta} $. In particular,
	\[ \commpart{x_\alpha}{x_\beta}{\gamma}^{w_\xi^2} = \commpart{x_\alpha^{w_\xi^2}}{x_\beta^{w_\xi^2}}{\gamma} \]
	for all $ x_\alpha \in \rootgr{\alpha} $ and all $ x_\beta \in \rootgr{\beta} $.
\end{lemma}
\begin{proof}
	Denote by $ \tup{\gamma}{k} $ be the unique interval ordering of $ \indivset{\clrootintcox{\alpha}{\beta}} $ such that $ \gamma_1 = \alpha $ and $ \gamma_k = \beta $, and let $ i \in \numint{1}{k} $ such that $ \gamma = \gamma_i $. Then $ (\gamma_1^{\reflbr{\xi}}, \ldots, \gamma_k^{\reflbr{\xi}}) $ is the unique interval ordering of $ \indivset{\clrootintcox{\alpha^{\reflbr{\xi}}}{\beta^{\reflbr{\xi}}}} $ whose first element is $ \alpha^{\reflbr{\xi}} $ and whose last element in $ \beta^{\reflbr{\xi}} $. By the definition of $ \commpart{x_\alpha}{x_\beta}{\gamma} $, we have
	\begin{align*}
		\commutator{x_\alpha}{x_\beta} \in \rootgr{\gamma_1} \cdots \rootgr{\gamma_{i-1}} \commpart{x_\alpha}{x_\beta}{\gamma} \rootgr{\gamma_{i+1}} \cdots \rootgr{\gamma_{k}}
	\end{align*}
	and $ \commpart{x_\alpha}{x_\beta}{\gamma} \in \rootgr{\gamma} $. It follows that
	\begin{align*}
		\commutator{x_\alpha^{w_\xi}}{x_\beta^{w_\xi}} = \commutator{x_\alpha}{x_\beta}^{w_\xi} &\mathord{}\in \rootgr{\gamma_1}^{w_\xi} \cdots \rootgr{\gamma_{i-1}}^{w_\xi} \commpart{x_\alpha}{x_\beta}{\gamma}^{w_\xi} \rootgr{\gamma_{i+1}}^{w_\xi} \cdots \rootgr{\gamma_{k}}^{w_\xi} \\
		&\hspace{1cm} \mathord{}= \rootgr{\gamma_1^{\reflbr{\xi}}} \cdots \rootgr{\gamma_{i-1}^{\reflbr{\xi}}} \commpart{x_\alpha}{x_\beta}{\gamma}^{w_\xi} \rootgr{\gamma_{i+1}^{\reflbr{\xi}}} \cdots \rootgr{\gamma_{k}^{\reflbr{\xi}}}
	\end{align*}
	and $ \commpart{x_\alpha}{x_\beta}{\gamma}^{w_\xi} \in \rootgr{\refl{\xi}(\gamma)} $. This implies that $ \commpart{x_\alpha^{w_\xi}}{x_\beta^{w_\xi}}{\refl{\xi}(\gamma)} = \commpart{x_\alpha}{x_\beta}{\gamma}^{w_\xi} $.
\end{proof}

The following result is crucial in \cite{MoufangPolygons}, and it will be an important tool in our computations as well.

\begin{proposition}[{\cite[(6.4)]{MoufangPolygons}}]\label{tw:6.4}
	Put $ n \defl \abs{\roots}/2 $. Assume that $ \roots $ has rank 2 with $ n \ge 3 $ and that $ G $ is rank-2-injective. Further, we fix a rank-2 labeling of $ \roots $ and we will use Notation~\thmitemref{weyl:weyl-def}{weyl:weyl-def:rank2}. Let $ z \in \IZ $, let $ b_1 \in \invset{z+1}{} $ and let $ x_n \in \rootgr{z+n}{} $. Let $ a_{1+n}, c_{1+n} \in \invset{z+1+n} $ such that $ (a_{1+n}, b_1, c_{1+n}) $ is a $ (z+1) $-Weyl triple, and denote by $ w_1 $ the associated Weyl element. For all $ i \in \numint{2}{n-1} $, we put $ x_i \defl \commpart{b_1}{x_n^{-1}}{z+i} \in \rootgr{z+i} $, so that $ \commutator{b_1}{x_n^{-1}} = x_2 \cdots x_{n-1} $. Then the following hold:
	\begin{lemenumerate}
		\item \label{tw:6.4:conj}$ x_n^{w_1} = x_2 = \commpart{b_1}{x_n^{-1}}{z+2} $. In particular, $ x_2 $ lies in $ \invset{z+2} $ if $ x_n $ lies in $ \invset{z+n} $.
		
		\item \label{tw:6.4:comm}$ \commutator{x_2}{c_{1+n}^{-1}} = x_3 \cdots x_n $.
	\end{lemenumerate}
\end{proposition}
\begin{proof}
	Put $ w \defl w_1 = a_{1+n} b_1 c_{1+n} $. Since
	\[ \commutator{a_{1+n}^{-1}}{x_n^{-1}} \in \commutator{\rootgr{z+n+1}{}}{\rootgr{z+n}{}} = \compactSet{1_G}, \]
	it follows from \thmitemcref{group-rel}{group-rel:add-cor} that
	\begin{align*}
		\commutator{b_1}{x_{n}^{-1}} &= \commutator{a_{1+n}^{-1} w c_{1+n}^{-1}}{x_n^{-1}} = \commutator{w c_{1+n}^{-1}}{x_n^{-1}}.
	\end{align*}
	Therefore,
	\begin{align*}
		x_2 \cdots x_n &= \commutator{b_1}{x_n^{-1}} x_n = \commutator{w c_{1+n}^{-1}}{x_n^{-1}} x_n = \brackets[\big]{c_{1+n} w^{-1} x_n w c_{1+n}^{-1} x_n^{-1}} x_n \\
		&= c_{1+n} x_n^w c_{1+n}^{-1} = x_n^w (x_n^w)^{-1} c_{1+n} x_n^w c_{1+n}^{-1} = x_n^w \commutator{x_n^w}{c_{1+n}^{-1}}.
	\end{align*}
	Note that $ c_{1+n} $ lies in $ \rootgr{z+n+1} $ and that $ x_n^w $ lies in
	\[ \rootgr{z+n}^w = \rootgr{z+n}^{w_1} = \rootgr{2(z+1) + n - (z+n)} = \rootgr{z+2} \]
	by \cref{basic:rank2-notation}. Thus $ \commutator{x_n^w}{c_{1+n}^{-1}} $ lies in $ \rootgr{\numint{z+3}{z+n}} $. Therefore, it follows from
	\[ x_2 x_3 \cdots x_n = x_n^w \commutator{x_n^w}{c_{1+n}^{-1}} \]
	and the rank-2-injectivity of $ G $ that $ x_2 = x_n^w $ and $ x_3 \cdots x_n = \commutator{x_n^w}{c_{1+n}^{-1}} = \commutator{x_2}{c_{1+n}^{-1}} $. This finishes the proof.
\end{proof}

\begin{note}
	Special cases of \cref{tw:6.4} for specific values of $ n $ can also be proven in a slightly different way. We will see this in \cref{A2Weyl:basecomp-cor,B:basecomp-walpha-cox-cor-delta}. We will also perform a similar computation in \cref{B:weylbeta-on-delta-nonabelian} for a pair of roots which does not form a root base.
\end{note}

\subsection{Balanced Weyl Triples}

	Any Weyl triple $ (a, b, c) $ is associated with four elements of $ G $: $ a $, $ b $, $ c $ and $ w \defl abc $. These four elements are related by the equation $ w=abc $, so in particular, the knowledge of any three of these group elements determines the missing fourth element. We now investigate under which conditions we can find more relations between these four elements.

\begin{definition}[Balanced Weyl triple]\label{weyl:balanced}
	A Weyl triple $ (a, b, c) $ with corresponding Weyl element $ w \defl abc $ is called \defemph*{weakly balanced}\index{Weyl triple!weakly balanced} if $ a=c $, and it is called \defemph*{balanced}\index{Weyl triple!balanced} if, in addition, it satisfies the equivalent conditions of the following \cref{basic:weyl-balanced-char}.
\end{definition}

\begin{note}
	We will also refer to the Weyl element $ w $ in \cref{weyl:balanced} as (weakly) balanced if the corresponding conditions are satisfied, but this is actually abuse of notation: In general, the Weyl triple $ (a,b,c) $ is not necessarily uniquely determined by the Weyl element $ abc $. (For examples, see \cref{B:ex:short-weyl-matrix,BC:ex:short-weyl-matrix}.) However, it will always be clear from the context what we mean.
\end{note}

\begin{lemma}[{\cite[5.15]{LoosNeherBook}}]\label{basic:weyl-balanced-char}
	Let $ \alpha $ be a root and let $ (a_{-\alpha}, b_\alpha, a_{-\alpha}) $ be a weakly balanced $ \alpha $-Weyl triple with corresponding Weyl element $ w_\alpha \defl a_{-\alpha} b_\alpha a_{-\alpha} $. Then the following statements are equivalent:
	\begin{lemenumerate}
		\item \label{basic:weyl-balanced-char:1}$ b_\alpha^{w_\alpha} = a_{-\alpha} $.
		
		\item \label{basic:weyl-balanced-char:2}$ a_{-\alpha}^{w_\alpha} = b_\alpha $.
		
		\item \label{basic:weyl-balanced-char:3}$ w_\alpha = b_\alpha a_{-\alpha} b_\alpha $.
	\end{lemenumerate}
\end{lemma}
\begin{proof}
	We know from \thmitemcref{basic:weyl-general}{basic:weyl-general:minus} that the Weyl elements $ w_\alpha' \defl a_{-\alpha}^{w_\alpha^{-1}} a_{-\alpha} b_\alpha $ and $ w_\alpha'' \defl b_\alpha a_{-\alpha} a_{-\alpha}^{w_\alpha} $ are both equal to $ w_\alpha $. Clearly, $ w_\alpha' $ equals $ b_\alpha a_{-\alpha} b_\alpha $ if and only if $ a_{-\alpha}^{w_\alpha^{-1}} = b_\alpha $, so~\itemref{basic:weyl-balanced-char:1} is equivalent to~\itemref{basic:weyl-balanced-char:3}. Similarly, $ w_\alpha'' $ equals $ b_\alpha a_{-\alpha} b_\alpha $ if and only if $ a_{-\alpha}^{w_\alpha} = b_\alpha $, so~\itemref{basic:weyl-balanced-char:2} is equivalent to~\itemref{basic:weyl-balanced-char:3} as well.
\end{proof}

The notion of balanced Weyl triples is due to \cite[5.15]{LoosNeherBook} while the terminology of weakly balanced Weyl triples is, to our knowledge, new. It is clear that in a weakly balanced Weyl triple, the knowledge of two of the elements $ a $, $ b $ and $ w \defl aba $ is enough to reconstruct the third one.

\begin{remark}
	If $ (a_{-\alpha}, b_\alpha, a_{-\alpha}) $ is a balanced $ \alpha $-Weyl triple with associated Weyl element $ w_\alpha $, then
	\begin{align*}
		b_\alpha^{w_\alpha^2} = a_{-\alpha}^{w_{-\alpha}} = b_\alpha.
	\end{align*}
\end{remark}

A natural question to ask is whether weakly balanced Weyl triples are automatically balanced.
We will see in \thmitemcref{basic:weylmap-form2}{basic:weylmap-form2:balanced} that that this is true if the group has unique Weyl extensions. This extra assumptions will be satisfied in all cases that we are interested in by \cref{basic:unique-weyl-ext-crit}. In more general situations, however, the distinction between balanced and weakly balanced Weyl triples is relevant. Still, we make the following observation.

\begin{proposition}\label{basic:all-balanced}
	If $ \alpha $ is a root such that all $ \alpha $-Weyl triples are weakly balanced, then all $ \alpha $-Weyl triples are balanced. In particular, if all Weyl triples in $ G $ are weakly balanced, then all Weyl triples in $ G $ are balanced.
\end{proposition}
\begin{proof}
	Let $ \alpha $ be a root and assume that all $ \alpha $-Weyl triples are weakly balanced. If there exist no $ \alpha $-Weyl triples, then there is nothing to prove, so we can assume that there exists an $ \alpha $-Weyl triple $ (a_{-\alpha}, b_\alpha, a_{-\alpha}) $ with corresponding Weyl element $ w_\alpha \defl a_{-\alpha} b_\alpha a_{-\alpha} $. If $ (\bar{a}_{\alpha}, \bar{b}_{-\alpha}, \bar{c}_{\alpha}) $ is a $ (-\alpha) $-Weyl triple, then $ (\bar{a}_{\alpha}^{w_\alpha}, \bar{b}_{-\alpha}^{w_\alpha}, \bar{c}_{\alpha}^{w_\alpha}) $ is an $ \alpha $-Weyl triple by \thmitemcref{basic:weyl-general}{basic:weyl-general:conj}, so it follows from our assumptions that $ \bar{a}_{\alpha}^{w_\alpha} = \bar{c}_{\alpha}^{w_\alpha} $. We conclude that all $ (-\alpha) $-Weyl triples are weakly balanced as well. Since $ (b_\alpha, a_{-\alpha}, a_{-\alpha}^{w_\alpha}) $ is a $ (-\alpha) $-Weyl triple by \thmitemcref{basic:weyl-general}{basic:weyl-general:minus}, this implies that $ b_\alpha = a_{-\alpha}^{w_\alpha} $, so $ (a_{-\alpha}, b_\alpha, a_{-\alpha}) $ is balanced. This finishes the proof.
\end{proof}

\begin{note}
	A second natural question is whether every Weyl triple is weakly balanced. We will show in \cref{A2Weyl:weyl} that every Weyl element in $ A_2 $-graded groups is weakly balanced and thus, by \cref{basic:all-balanced} or \thmitemcref{basic:weylmap-form2}{basic:weylmap-form2:balanced}, balanced. This implies that the same is true for all roots $ \alpha $ in an arbitrary root system which are contained in a (crystallographically) closed subsystem of type $ A_2 $. This holds, for example, for the long roots in $ B_n $, the short roots in $ C_n $ (where $ n \ge 3 $) and for all roots in $ H_3 $ and $ H_4 $. However, it does not hold for the short roots in $ B_n $ and the long roots in $ C_n $. We can still prove that all Weyl elements for such roots are balanced (see \cref{B:short-weakly-balanced,BC:long-weakly-balanced}), but this proof relies on computations in the coordinatising structures of these groups. The only remaining case are the short roots in $ BC_n $ for $ n \ge 3 $. For these roots, there do in fact exist examples of Weyl elements which are not weakly balanced. We will encounter them in \cref{BC:ex-weyl-def}.
\end{note}


\subsection{Unique Weyl Extensions}

\begin{definition}[Unique Weyl extensions]\label{basic:mumap-def}
	Let $ \alpha \in \roots $. We say that \defemph*{$ G $ has unique $ \alpha $-Weyl extensions}\index{unique Weyl extensions} if for all $ b_\alpha \in \invset{\alpha} $, there exist \emph{unique} elements $ a_{-\alpha}, c_{-\alpha} \in \rootgr{-\alpha} $ such that $ (a_{-\alpha}, b_\alpha, c_{-\alpha}) $ is an $ \alpha $-Weyl triple. If this is the case, we define $ \map{\weylleft{\alpha}, \weylright{\alpha}}{\invset{\alpha}}{\rootgr{\alpha}}{}{} $ to be the unique maps for which $ \weylleft{\alpha}(b_\alpha) b_\alpha \weylright{\alpha}(b_\alpha) $ is an $ \alpha $-Weyl element for all $ b_\alpha \in \invset{\alpha} $. Further, we define $ \map{\weylmap{\alpha}}{\invset{\alpha}}{\weylset{\alpha}}{b_\alpha}{\weylleft{\alpha}(b_\alpha) b_\alpha \weylright{\alpha}(b_\alpha)} $. We say that \defemph*{$ G $ has unique Weyl extensions} if it has unique $ \alpha $-Weyl extensions for all roots $ \alpha $. We will sometimes leave out the subscripts of the maps $ \weylleft{} $, $ \weylright{} $ and $ \weylmap{} $ if they are clear from the context.
\end{definition}

The property of having unique Weyl extensions says precisely that a Weyl triple is uniquely determined by its \enquote{middle element}. Under mild conditions on the root system and the root groups, it is in fact true that a Weyl triple $ (a,b,c) $ is uniquely determined by any of the elements $ a $, $ b $, $ c $. This is the statement of the following \cref{basic:weyl-det-by-one-factor}, which is a special case of \cite[5.20]{LoosNeherBook}. It is also proven in a less general context in \cite[(6.1)]{MoufangPolygons}. We begin the proof of \cref{basic:weyl-det-by-one-factor} with an auxiliary lemma.

\begin{lemma}\label{basic:centraliser-opp-trivial}
	Let $ \alpha, \beta $ be roots such that $ \rootgr{\alpha} \intersect \rootgr{\refl{\beta}(\alpha)} = \compactSet{1_G} $ and assume that there exists a $ \beta $-Weyl element $ w_\beta $. Then if $ x_\alpha $ is an element of $ \rootgr{\alpha} $ which centralises $ \rootgr{\beta} $ and $ \rootgr{-\beta} $, it follows that $ x_\alpha = 1_G $.
\end{lemma}
\begin{proof}
	By assumption, we have $ x_\alpha = x_\alpha^{w_\beta} $ where the element on the left-hand side lies in $ \rootgr{\alpha} $ while the element on the right-hand side lies in $ \rootgr{\refl{\beta}(\alpha)} $. Since $ \rootgr{\alpha} \intersect \rootgr{\refl{\beta}(\alpha)} = \compactSet{1_G} $, it follows that $ x_\alpha = 1_G $.
\end{proof}

\begin{lemma}[{\cite[5.20]{LoosNeherBook}}]\label{basic:weyl-det-by-one-factor}
	Let $ \alpha $ be a root and let
	\[ \weyltrip = (a_{-\alpha}, b_\alpha, c_{-\alpha}) \midand \tilde{\weyltrip} = (\tilde{a}_{-\alpha}, \tilde{b}_\alpha, \tilde{c}_{-\alpha}) \]
	be two $ \alpha $-Weyl triples. Assume that $ \alpha $ is contained in a parabolic rank-$ 2 $ closed subsystem $ \roots' $ of $ \roots $ which is not of type $ A_1 \times A_1 $ and that there exists a $ \beta $-Weyl element for some root $ \beta \in \roots' $ adjacent to $ \alpha $. Assume further that $ G $ is rank-2-injective. Then if one of the three statements $ a_{-\alpha} = \tilde{a}_{-\alpha} $, $ b_\alpha = \tilde{b}_\alpha $, $ c_{-\alpha} = \tilde{c}_{-\alpha} $ holds, it follows that $ \weyltrip = \tilde{\weyltrip} $.
\end{lemma}
\begin{proof}
	Denote by $ \roots' $ a parabolic rank-2 subsystem of $ \roots $ and by $ \beta $ a root which are chosen as in the assertion. Without loss of generality, we assume that $ \roots' = \roots $, and we fix the unique rank-2 labeling $ \map{}{}{}{i}{\alpha_i} $ of $ \roots $ with $ \alpha_1 = \alpha $ and $ \alpha_2 = \beta $. Further, we put $ m \defl \abs{\roots}/2 $. Since $ \roots $ is not of type $ A_1 \times A_1 $, we have $ m \ge 3 $.
	
	We begin with the case $ b_\alpha = \tilde{b}_\alpha $. Put
	\[ w_\alpha \defl a_{-\alpha} b_\alpha c_{-\alpha} \midand \tilde{w}_\alpha \defl \tilde{a}_{-\alpha} \tilde{b}_\alpha \tilde{c}_{-\alpha} \]
	Then for $ \bar{a}_{-\alpha} \defl \tilde{a}_{-\alpha} a_{-\alpha}^{-1} $, $ \bar{c}_{-\alpha} \defl c_{-\alpha}^{-1} \tilde{c}_{-\alpha} $ and for all $ x_m \in \rootgr{\alpha_m} $, we have
	\begin{align*}
		x_m^{\tilde{w}_\alpha} &= x_m^{\bar{a}_{-\alpha} w_\alpha \bar{c}_{-\alpha}} = (x_m^{w_\alpha})^{\bar{c}_{-\alpha}} = x_m^{w_\alpha} \commutator{x_m^{w_\alpha}}{\bar{c}_{-\alpha}}.
	\end{align*}
	Note that $ x_m^{\tilde{w}_\alpha} $ and $ x_m^{w_\alpha} $ lie in $ \rootgr{\refl{\alpha}(\alpha_m)} = \rootgr{\alpha_2} $ while $ \commutator{x_m^{w_\alpha}}{\bar{c}_{-\alpha}} $ lies in $ \rootgr{\rootintcox{\alpha_2}{\alpha_{m+1}}} = \rootgr{\alpha_3} \cdots \rootgr{\alpha_m} $ for all $ x_m \in \rootgr{\alpha_m} $. Since the product map on $ \word{\alpha} $ is injective, we conclude that $ x_m^{\tilde{w}_\alpha} = x_m^{w_\alpha} $ and $ \commutator{x_m^{w_\alpha}}{\bar{c}_{-\alpha}} = 1_G $ for all $ x_m \in \rootgr{\alpha_m} $. In particular, $ \bar{c}_{-\alpha} $ centralises $ \rootgr{\alpha_m}^{w_\alpha} = \rootgr{\alpha_2} = \rootgr{\beta} $. As $ \alpha $ is adjacent to $ \beta $ (and thus $ -\alpha $ is adjacent to $ -\beta $), we also have that $ \bar{c}_{-\alpha} $ centralises $ \rootgr{-\beta} $. Note further that $ \rootgr{-\alpha} \intersect \rootgr{\refl{\alpha}(\beta)} = \rootgr{\alpha_{m+1}} \intersect \rootgr{\alpha_m} = \compactSet{1_G} $ by \cref{basic:prodmap-triv-intersect}. Thus it follows from \cref{basic:centraliser-opp-trivial} that $ \bar{c}_{-\alpha} = 1_G $. In other words, $ c_{-\alpha} = \tilde{c}_{-\alpha} $. Since $ (c_{-\alpha}^{-1}, b_\alpha^{-1}, a_{-\alpha}^{-1}) $ and $ (\tilde{c}_{-\alpha}^{-1}, \tilde{b}_\alpha^{-1}, \tilde{a}_{-\alpha}^{-1}) $ are $ \alpha $-Weyl elements as well (by \thmitemcref{basic:weyl-general}{basic:weyl-general:inv}), the same chain of arguments yields that $ a_{-\alpha} = \tilde{a}_{-\alpha} $, too. This finishes the proof of the case $ b_\alpha = b_{\alpha}' $.
	
	Now assume that $ c_{-\alpha} = \tilde{c}_{-\alpha} $. By \thmitemcref{basic:weyl-general}{basic:weyl-general:minus}, we also have $ (-\alpha) $-Weyl triples $ (b_\alpha, c_{-\alpha}, a_{-\alpha}^{w_\alpha}) $ and $ (\tilde{b}_\alpha, \tilde{c}_{-\alpha}, \tilde{a}_{-\alpha}^{\tilde{w}_\alpha}) $. Observe that the requirements of the current lemma are also satisfied for $ -\alpha $ in place of $ \alpha $ because the existence of a $ \beta $-Weyl element implies the existence of a $ (-\beta) $-Weyl element, again by \thmitemcref{basic:weyl-general}{basic:weyl-general:minus}. Thus the conclusion of the previous paragraph allows us to infer that $ b_\alpha = \tilde{b}_{\alpha} $ and $ a_{-\alpha}^{w_\alpha} = \tilde{a}_{-\alpha}^{\tilde{w}_\alpha} $. The second equation says precisely that
	\[ c_{-\alpha}^{-1} b_\alpha^{-1} a_{-\alpha}^{-1} a_{-\alpha} a_{-\alpha} b_\alpha c_{-\alpha} = \tilde{c}_{-\alpha}^{-1} \tilde{b}_\alpha^{-1} \tilde{a}_{-\alpha}^{-1} \tilde{a}_{-\alpha} \tilde{a}_{-\alpha} \tilde{b}_\alpha \tilde{c}_{-\alpha}, \]
	which implies that $ a_{-\alpha} =  \tilde{a}_{-\alpha} $ because we already know that $ b_\alpha = \tilde{b}_{\alpha} $ and $ c_{-\alpha} = \tilde{c}_{-\alpha} $. We conclude that $ \weyltrip = \tilde{\weyltrip} $ in this case as well. The case $ a_{-\alpha} =  \tilde{a}_{-\alpha} $ can be proven in the same way by considering the $ (-\alpha) $-Weyl triples $ (c_{-\alpha}^{w_\alpha^{-1}}, a_{-\alpha}, b_\alpha) $ and $ (\tilde{c}_{-\alpha}^{\tilde{w}_\alpha^{-1}}, \tilde{a}_{-\alpha}, \tilde{b}_\alpha) $.
\end{proof}

\begin{note}
	\cref{basic:weyl-det-by-one-factor} does not rule out the possibility of two distinct Weyl triples having the same associated Weyl element. In fact, there exist short Weyl elements in root gradings of types $ B $ and $ BC $ which have two distinct associated Weyl triples: see \cref{B:ex:short-weyl-matrix,BC:ex:short-weyl-matrix}. However, in $ A_2 $-graded groups, every Weyl element has a unique associated Weyl triple (\thmitemcref{A2Weyl:weyl}{A2Weyl:weyl:weylel-gives-triple}).
\end{note}

For most applications, we can reduce \cref{basic:weyl-det-by-one-factor} to the following statement.

\begin{proposition}\label{basic:unique-weyl-ext-crit}
	Assume that $ \roots $ does not have an irreducible component of type $ A_1 $, that $ \invset{\alpha} \ne \emptyset $ for all roots $ \alpha $ and that $ G $ is rank-2-injective. Then $ G $ has unique Weyl extensions.
\end{proposition}
\begin{proof}
	It follows from the assumption on $ \roots $ that every root is contained in a rank-2 subsystem which is not of type $ A_1 \times A_1 $. Thus the assertion is a consequence of \cref{basic:weyl-det-by-one-factor}.
\end{proof}

We now investigate some properties of groups with unique Weyl extensions. The following two statements are essentially reformulations of \cref{basic:weyl-general} using the maps $ \weylmap{} $, $ \weylleft{} $ and $ \weylright{} $. The proofs use the same arguments as in \cite{MoufangPolygons}, but in a more general context.

\begin{lemma}[{\cite[(6.2)]{MoufangPolygons}}]\label{basic:weylmap-form1}
	Let $ \alpha $ be a root such that $ G $ has unique $ \alpha $-Weyl extensions. Then the following hold:
	\begin{lemenumerate}
		\item $ \weylmap{\alpha}(b_\alpha^{-1}) = \weylmap{\alpha}(b_\alpha)^{-1} $, $ \weylleft{\alpha}(b_\alpha^{-1}) = \weylleft{\alpha}(b_\alpha)^{-1} $ and $ \weylright{\alpha}(b_\alpha^{-1}) = \weylright{\alpha}(b_\alpha)^{-1} $ for all $ b_\alpha \in \invset{\alpha} $.
		
		\item \label{basic:weylmap-form1:conj}Let $ \beta $ be a root for which there exists a $ \beta $-Weyl element $ w_\beta $. Then $ G $ has unique $ \refl{\beta}(\alpha) $-Weyl extensions. Further, we have
		\begin{align*}
			\weylmap{\refl{\beta}(\alpha)}(b_\alpha^{w_\beta}) = \weylmap{\alpha}(b_\alpha)^{w_\beta}, \; \weylleft{\refl{\beta}(\alpha)}(b_\alpha^{w_\beta}) = \weylleft{\alpha}(b_\alpha)^{w_\beta}, \; \weylright{\refl{\beta}(\alpha)}(b_\alpha^{w_\beta}) = \weylright{\alpha}(b_\alpha)^{w_\beta}
		\end{align*}
		for all $ b_\alpha \in \invset{\alpha} $.
	\end{lemenumerate}
\end{lemma}
\begin{proof}
	The first assertion follows from \thmitemcref{basic:weyl-general}{basic:weyl-general:inv} and the second one from \thmitemcref{basic:weyl-general}{basic:weyl-general:conj}.
\end{proof}

\begin{lemma}[{\cite[(6.3)]{MoufangPolygons}}]\label{basic:weylmap-form2}
	Let $ \alpha $ be a root such that $ G $ has unique $ \alpha $-Weyl extensions and assume that there exists $ b_\alpha \in \invset{\alpha} $. Then $ G $ has unique $ (-\alpha) $-Weyl extensions, the images of $ \weylleft{\alpha} $ and $ \weylright{\alpha} $ are contained in $ \invset{-\alpha} $ and the following hold:
	\begin{lemenumerate}
		\item \label{basic:weylmap-form2:same-mu}$ \weylmap{\alpha}(b_\alpha) = \weylmap{-\alpha}\brackets[\big]{\weylleft{\alpha}(b_\alpha)} = \weylmap{-\alpha}\brackets[\big]{\weylright{\alpha}(b_\alpha)} $.
		
		\item $ \weylleft{-\alpha}\brackets[\big]{\weylright{\alpha}(b_\alpha)} = b_\alpha = \weylright{-\alpha}\brackets[\big]{\weylleft{\alpha}(b_\alpha)} $.
		
		\item $ \weylright{-\alpha}\brackets[\big]{\weylright{\alpha}(b_\alpha)} = \weylleft{\alpha}(b_\alpha)^{\weylmap{\alpha}(b_\alpha)} $.
		
		\item $ \weylleft{-\alpha}\brackets[\big]{\weylleft{\alpha}(b_\alpha)} = \weylright{\alpha}(b_\alpha)^{\weylmap{\alpha}(b_\alpha)^{-1}} $.
		
		\item \label{basic:weylmap-form2:balanced}If $ \weylright{\alpha}(b_\alpha) = \weylleft{\alpha}(b_\alpha) $, then
		\[ b_\alpha^{\weylmap{\alpha}(b_\alpha)} = \weylleft{\alpha}(b_\alpha) = \weylright{\alpha}(b_\alpha). \]
		That is, every weakly balanced $ \alpha $-Weyl triple is balanced.
	\end{lemenumerate}
\end{lemma}
\begin{proof}
	It follows from \thmitemcref{basic:weylmap-form1}{basic:weylmap-form1:conj} that $ G $ has unique $ (-\alpha) $-Weyl extensions (because $ \refl{\alpha}(\alpha) = -\alpha $) and from \thmitemcref{basic:weyl-general}{basic:weyl-general:minus} that the images of $ \weylleft{\alpha} $ and $ \weylright{\alpha} $ are contained in $ \invset{-\alpha} $. We know from \thmitemcref{basic:weyl-general}{basic:weyl-general:minus} that
	\[ \brackets[\big]{b_\alpha, \weylright{\alpha}(b_\alpha), \weylleft{\alpha}(b_\alpha)^{\weylmap{\alpha}(b_\alpha)}} \]
	is a $ (-\alpha) $-Weyl triple with associated Weyl element $ \weylmap{\alpha}(b_\alpha) $. This implies that
	\begin{align*}
		\weylmap{-\alpha}\brackets[\big]{\weylright{\alpha}(b_\alpha)} = \weylmap{\alpha}(b_\alpha), \quad \weylleft{-\alpha}\brackets[\big]{\weylright{\alpha}(b_\alpha)} = b_\alpha, \quad \weylright{-\alpha}\brackets[\big]{\weylright{\alpha}(b_\alpha)} = \weylleft{\alpha}(b_\alpha)^{\weylmap{\alpha}(b_\alpha)}.
	\end{align*}
	Similarly, we know from \thmitemcref{basic:weyl-general}{basic:weyl-general:minus} that
	\[ \brackets[\big]{\weylright{\alpha}(b_\alpha)^{\weylmap{\alpha}(b_\alpha)^{-1}}, \weylleft{\alpha}(b_\alpha), b_\alpha} \]
	is a $ (-\alpha) $-Weyl triple with associated Weyl element $ \weylmap{\alpha}(b_\alpha) $, which implies that
	\begin{align*}
		\weylmap{-\alpha}\brackets[\big]{\weylleft{\alpha}(b_\alpha)} = \weylmap{\alpha}(b_\alpha), \quad \weylleft{-\alpha}\brackets[\big]{\weylleft{\alpha}(b_\alpha)} = \weylright{\alpha}(b_\alpha)^{\weylmap{\alpha}(b_\alpha)^{-1}}, \quad \weylright{-\alpha}\brackets[\big]{\weylleft{\alpha}(b_\alpha)} = b_\alpha.
	\end{align*}
	This finishes the proof of the first four assertions. If $ \weylright{\alpha}(b_\alpha) = \weylleft{\alpha}(b_\alpha) $, then it follows from the previous assertions that
	\begin{align*}
		b_\alpha &= \weylleft{-\alpha}\brackets[\big]{\weylright{\alpha}(b_\alpha)} = \weylleft{-\alpha}\brackets[\big]{\weylleft{\alpha}(b_\alpha)} = \weylright{\alpha}(b_\alpha)^{\weylmap{\alpha}(b_\alpha)^{-1}}.
	\end{align*}
	Conjugating both sides by $ \weylmap{\alpha}(b_\alpha) $, the final assertion follows.
\end{proof}

\subsection{The Braid Relations for Weyl Elements}

In the proof of the parametrisation theorem in \cref{chap:param}, we will require that a fixed family of Weyl elements satisfies the same braid relations as the Weyl group. We will also use these relations in the computation of certain \enquote{blueprint rewriting rules} which are needed to apply the blueprint technique. To formulate the braid relations for Weyl elements, we first have to define $ \rootbase $-systems of Weyl elements.

\begin{notation}[System of Weyl elements]\label{basic:weyl-system}
	Let $ \rootbase $ be a rescaled root base of $ \roots $. A \defemph*{$ \rootbase $-system of Weyl elements (in $ G $)}\index{Weyl element!system of} is a family $ (w_\delta)_{\delta \in \rootbase} $ such that $ w_\delta $ is a $ \delta $-Weyl element for each $ \delta \in \rootbase $. Given such a $ \rootbase $-system of Weyl elements, we put $ w_{-\delta} \defl w_{\delta}^{-1} $ for any $ \delta \in \rootbase $ and $ w_{\word{\delta}} \defl w_{\delta_1} \cdots w_{\delta_m} $ for any word $ \word{\delta} = \tup{\delta}{m} $ over $ \rootbase \union (-\rootbase) $ (with the convention that $ w_{\word{\delta}} = 1_G $ if $ \word{\delta} $ is the empty word).
\end{notation}

\begin{definition}[Braid relations]\label{basic:braid-def}
	Let $ \rootbase $ be a rescaled root base of $ \roots $ and let $ (w_\delta)_{\delta \in \roots} $ be a $ \rootbase $-system of Weyl elements. We say that \defemph*{$ (w_{\delta})_{\delta \in \rootbase} $ satisfies the braid relations (in $ G $)}\index{braid relations} if for any distinct $ \alpha, \beta \in \rootbase $, we have $ \braidword{o(\refl{\alpha} \refl{\beta})}(w_\alpha, w_\beta) = \braidword{o(\refl{\alpha} \refl{\beta})}(w_\beta, w_\alpha) $ where $ o(\refl{\alpha} \refl{\beta}) $ is the order of $ \refl{\alpha} \refl{\beta} $ in the Weyl group. Similarly, we say that $ (w_{\delta})_{\delta \in \rootbase} $ \defemph*{satisfies the braid relations modulo $ \zentrum(G) $} if the images of $ \braidword{o(\refl{\alpha} \refl{\beta})}(w_\alpha, w_\beta) $ and $ \braidword{o(\refl{\alpha} \refl{\beta})}(w_\beta, w_\alpha) $ in $ G/\zentrum(G) $ are equal for all distinct $ \alpha, \beta \in \rootbase $. Further, we say that \defemph*{$ G $ satisfies the braid relations for Weyl elements (modulo $ \zentrum(G) $)} if for any root base $ \rootbase $, each $ \rootbase $-system of Weyl elements satisfies the braid relations (modulo $ \zentrum(G) $).
\end{definition}

For the purposes of the parametrisation theorem and the blueprint technique, it would be sufficient to find one $ \rootbase $-system of Weyl elements which satisfies the braid relations modulo $ \zentrum(G) $. In practice, we can prove a stronger statement (\cref{braid:all:weyl-ext}): Any $ \rootbase $-system of Weyl elements satisfies the braid relations in $ G $, provided that $ G $ has unique Weyl extensions and that it is rank-2-injective. This is proven in \cite[(6.9)]{MoufangPolygons} under less general assumptions, but the same arguments remain valid in our context. We present this proof in the remaining part of this subsection.

To avoid notational problems, we briefly rule out the trivial case of the root system $ A_1 \times A_1 $.

\begin{proposition}\label{braid:A1xA1}
	Assume that $ \roots = A_1 \times A_1 $. Let $ \rootbase = \Set{\alpha, \beta} $ be a root base of $ A_1 \times A_1 $ and assume that there exist an $ \alpha $-Weyl element $ w_\alpha $ and a $ \beta $-Weyl element $ w_\beta $. Then the family $ (w_\gamma)_{\gamma \in \rootbase} $ satisfies the braid relations.
\end{proposition}
\begin{proof}
	We only have to show that $ w_\alpha w_\beta = w_\beta w_\alpha $, which is trivial because the subgroup $ \gen{\rootgr{\alpha}, \rootgr{-\alpha}} $ commutes with $ \gen{\rootgr{\beta}, \rootgr{-\beta}} $ by the commutator relations.
\end{proof}

\begin{miscthm}[Setup]\label{tw:braid:setup}
	Put $ n \defl \abs{\roots}/2 $. Assume that $ \roots $ is reduced of rank 2 with $ n \ge 3 $ and fix a rank-2 labeling of $ \roots $ (in the sense of \cref{rootsys:rank2-label-def}) as well as elements $ u \in \invset{1} $, $ v \in \invset{n} $. Assume further that $ G $ is rank-2-injective and that it has unique Weyl extensions for all roots. Define sequences $ (e_k)_{k \in \Npos} $ and $ (f_k)_{k \in \Npos} $ in $ G $ by
	\begin{align*}
		e_1 &\defl u \in \invset{1}{}, & f_1 &\defl v \in \invset{n}{} = \invset{1+(n-1)}{} \rightand \\
		e_k &\defl f_{k-1}^{\weylmap{}(e_{k-1})} \in \invset{k}{}, & f_k &\defl \weylright{}(e_{k-1}) \in \invset{k+(n-1)}{} \quad \text{for all } k \in \IN_{\ge 2}.
	\end{align*}
	(All these elements are Weyl-invertible in the sense of \thmitemcref{weyl:weyl-def}{weyl:weyl-def:invertible} by \thmitemcref{basic:weylmap-form1}{basic:weylmap-form1:conj} and \cref{basic:weylmap-form2}.) For all $ k \in \Npos $, we define
	\[ x_{k,1} \defl e_k \in \invset{k}{} \midand x_{k,n} \defl f_k \in \invset{k+(n-1)}{}. \]
	Further, we put
	\[ x_{k,i} \defl \commpart{e_k}{f_k^{-1}}{k+i-1} \in \rootgr{k+i-1} \]
	for all $ k \in \Npos $ and all $ i \in \numint{2}{n-1} $. Thus we have
	\begin{equation}\label{eq:tw:braid:1}
		\commutator{x_{k,1}}{x_{k,n}^{-1}} = x_{k,2} \cdots x_{k,n-1}
	\end{equation}
	for all $ k \in \Npos $. Finally, we put $ w_0 \defl \weylmap{}(v) $ and $ w_k \defl \weylmap{}(e_k) $ for all $ k \in \Npos $. Our goal is to show that $ (w_0, w_1) $ satisfies the braid relations, that is, that $ \braidword{n}(w_0, w_1) = \braidword{n}(w_1, w_0) $.
\end{miscthm}

\begin{lemma}\label{tw:braid:cl1}
	Let everything be as in~\ref{tw:braid:setup}. Then
	\begin{align*}
		f_k &= x_{k,n} = x_{k+1, n-1} = x_{k+2, n-2} = \cdots = x_{k+n-1, 1} = e_{k+n-1}
	\end{align*}
	for all $ k \in \Npos $. In particular, $ f_k = e_{k+n-1} $ for all $ k \in \Npos $.
\end{lemma}
\begin{proof}
	Let $ k \in \Npos $ be arbitrary. By \thmitemcref{tw:6.4}{tw:6.4:conj} and the definition of $ e_{k+1} $, we have
	\begin{align*}
		x_{k+1,1} &= e_{k+1} = f_k^{\weylmap{}(e_k)} = x_{k,n}^{\weylmap{}(x_{k,1})} = x_{k,2}.
	\end{align*}
	Further,
	\begin{align*}
	 x_{k+1,n} = f_{k+1} = \weylright{}(e_k) = \weylright{}(x_{k,1}).
	\end{align*}
	Together with~\eqref{eq:tw:braid:1}, these equations imply that
	\begin{align*}
		x_{k+1,2} \cdots x_{k+1,n-1} &= \commutator{x_{k+1,1}}{x_{k+1,n}^{-1}} = \commutator{x_{k,2}}{\weylright{}(x_{k,1})^{-1}} = x_{k,3} \cdots x_{k,n}.
	\end{align*}
	Since $ G $ is rank-2-injective, we infer that $ x_{k+1,i} = x_{k,i+1} $ for all $ i \in \numint{2}{n-1} $. The assertion follows.
\end{proof}

The following statement is a mere corollary of \cref{tw:braid:cl1}, but it is worth to be pointed out on its own.

\begin{lemma}
	Let $ G $, $ \roots $ and the rank-2 labelling be as in~\ref{tw:braid:setup}. Let $ b_1 \in \invset{1} $ and $ b_n \in \invset{n} $. Then for all $ i \in \numint{2}{n-1} $, the element $ \commpart{b_1}{b_n}{i} $ lies in $ \invset{i} $.
\end{lemma}
\begin{proof}
	Put $ u \defl b_1 \in \invset{1} $ and $ v \defl b_n^{-1} $. By \thmitemcref{basic:weyl-general}{basic:weyl-general:inv}, $ v $ is contained in $ \invset{n} $. Thus we can define sequences $ (e_k)_{k \in \Npos} $, $ (f_k)_{k \in \Npos} $ and $ (x_{k,i})_{k \in \Npos, i \in \numint{2}{n-1}} $ as in~\ref{tw:braid:setup}. Then for all $ i \in \numint{2}{n-1} $, we have
	\begin{align*}
		\commpart{b_1}{b_n}{i} &= x_{1,i} = x_{i, 1} = e_i \in \invset{i},
	\end{align*}
	as desired.
\end{proof}

\begin{lemma}\label{tw:braid:cl2}
	Let everything be as in~\ref{tw:braid:setup}. Then we have $ w_k = w_{k+n} $ for all $ k \in \Nzero $ and $ w_{k-1}^{w_k} = w_{k+1} $ for all $ k \in \Npos $.
\end{lemma}
\begin{proof}
	Let $ k \in \Npos $. It follows from \thmitemcref{basic:weylmap-form2}{basic:weylmap-form2:same-mu}, \cref{tw:braid:cl1} and the definition of $ f_{k+1} $ that
	\begin{align*}
		w_k &= \weylmap{}(e_k) = \weylmap{}\brackets[\big]{\weylright{}(e_k)} = \weylmap{}(f_{k+1}) = \weylmap{}(e_{k+n}) = w_{k+n}.
	\end{align*}
	Similarly,
	\[ w_0 = \weylmap{}(v) = \weylmap{}(f_1) = \weylmap{}(e_n) = w_n. \]
	This proves the first claim. Further, we have
	\begin{align*}
		e_{k+1} &= f_k^{\weylmap{}(e_k)} = e_{k+n-1}^{w_k}.
	\end{align*}
	By \thmitemcref{basic:weylmap-form1}{basic:weylmap-form1:conj} and the previous claim, this implies that
	\[ w_{k+1} = \weylmap{}(e_{k+1}) = \weylmap{}(e_{k+n+1}^{w_k}) = \weylmap{}(e_{k+n+1})^{w_k} = w_{k+n-1}^{w_k} = w_{k-1}^{w_k}, \]
	as desired.
\end{proof}

\begin{lemma}\label{tw:braid:cl3}
	Let everything be as in~\ref{tw:braid:setup}. Then for all $ k \in \Npos $, we have
	\[ w_{2k} = w_0^{w_1 (w_0 w_1)^{k-1}} \midand w_{2k+1} = w_1^{(w_0 w_1)^k}. \]
\end{lemma}	
\begin{proof}
	We prove that statement by induction on $ k $. For $ k=1 $, it is a direct consequence of \cref{tw:braid:cl2} that
	\begin{align*}
		w_{2k} &= w_2 = w_0^{w_1} = w_0^{w_1 (w_0 w_1)^{k-1}} \rightand \\
		w_{2k+1} &= w_3 = w_1^{w_2} = w_1^{w_1^{-1} w_0 w_1} = w_1^{w_0 w_1} = w_1^{(w_0 w_1)^k}.
	\end{align*}
	Now assume that the desired statement is true for some $ k \in \Npos $. We proceed to show that it is true for $ k+1 $ as well. Again by \cref{tw:braid:cl2}, we have
	\begin{align*}
		w_{2k+2} = w_{2k}^{w_{2k+1}} = w_{2k+1}^{-1} w_{2k} w_{2k+1}.
	\end{align*}
	By the induction hypothesis, it follows that
	\begin{align*}
		w_{2k+2} &= \brackets[\big]{w_1^{(w_0 w_1)^k}}^{-1} w_0^{w_1 (w_0 w_1)^{k-1}} w_1^{(w_0 w_1)^k} = \brackets[\big]{(w_1^{w_0})^{-1} w_0 (w_1^{w_0})}^{w_1 (w_0 w_1)^{k-1}} \\
		&= \brackets[\big]{w_0^{-1} w_1^{-1} w_0 w_0 w_0^{-1} w_1 w_0}^{{w_1 (w_0 w_1)^{k-1}}} = \brackets[\big]{w_0^{-1} w_1^{-1} w_0 w_1 w_0}^{{w_1 (w_0 w_1)^{k-1}}} \\
		&= (w_0^{w_1 w_0})^{{w_1 (w_0 w_1)^{k-1}}} = w_0^{w_1(w_0 w_1)^k}.
	\end{align*}
	In the same way, we have
	\begin{align*}
		w_{2k+3} &= w_{2k+1}^{w_{2k+2}} = w_{2k+2}^{-1} w_{2k+1} w_{2k+2} = \brackets[\big]{w_0^{w_1 (w_0 w_1)^k}}^{-1} w_1^{(w_0 w_1)^k} w_0^{w_1 (w_0 w_1)^k} \\
		&= \brackets[\big]{(w_0^{w_1})^{-1} w_1 w_0^{w_1}}^{(w_0 w_1)^k} = \brackets[\big]{w_1^{-1} w_0^{-1} w_1 w_1 w_1^{-1} w_0 w_1}^{(w_0 w_1)^k} \\
		&= \brackets[\big]{w_1^{-1} w_0^{-1} w_1 w_0 w_1}^{(w_0 w_1)^k} = (w_1^{w_0 w_1})^{(w_0 w_1)^k} = w_1^{(w_0 w_1)^{k+1}}.
	\end{align*}
	This finishes the proof.
\end{proof}

\begin{remark}\label{braid:nonred}
	Assume that $ \bar{\roots} $ is a root system which is crystallographic and not reduced and let $ \bar{G} $ be a group which has $ \bar{\roots} $-commutator relations with root groups $ (\bar{U}_\alpha)_{\alpha \in \bar{\roots}} $. We do not assume that these commutator relations are crystallographic. For any root $ \alpha $, we denote by $ \alpha' $ the unique indivisible root in $ \IR_{>0} \alpha $. Then $ \alpha' \in \Set{\alpha, \alpha/2} $ for all roots $ \alpha $ by \cref{rootsys:cry-nonred-2}. Thus it follows from Axiom~\axiomref{rgg-axiom-div} that $ \bar{U}_{\alpha} \subs \bar{U}_{\alpha'} $ for all roots $ \alpha $. Hence every $ \alpha $-Weyl element is also an $ \alpha' $-Weyl element for all roots $ \alpha $ and the groups $ (\bar{U}_\alpha)_{\alpha \in \indivset{\bar{\roots}}} $ satisfy $ \indivset{\bar{\roots}} $-commutator relations. We conclude that, in order to verify the braid relations for Weyl elements in $ (\bar{G}, (\bar{U}_\alpha)_{\alpha \in \bar{\roots}}) $, it suffices to verify the braid relations in $ (\bar{G}, (\bar{U}_\alpha)_{\alpha \in \indivset{\bar{\roots}}}) $.
\end{remark}

\begin{theorem}[{\cite[(6.9)]{MoufangPolygons}}]\label{braid:all:weyl-ext}
	Assume that $ \roots $ is crystallographic or reduced, that $ G $ is rank-2-injective and that $ \invset{\alpha} \ne \emptyset $ for all roots $ \alpha $. Then $ G $ satisfies the braid relations for Weyl elements.
\end{theorem}
\begin{proof}
	If $ \roots $ is of rank 1, there is nothing to show, so we assume that it has rank at least~2. Let $ \rootbase $ be any root base and let $ \alpha, \beta \in \roots $ be distinct. It suffices to consider the root subsystem spanned by $ \Set{\alpha, \beta} $, so we can assume that $ \roots $ is of rank~2. Further, we can assume by \cref{braid:nonred} that $ \roots $ is reduced.
	
	Put $ n \defl \abs{\roots}/2 $. If $ n=2 $, then the assertion holds by \cref{braid:A1xA1}, so assume that $ n \ge 3 $. In this situation, we know from \cref{basic:unique-weyl-ext-crit} that $ G $ has unique Weyl extensions. Fix the rank-2 labeling of $ \roots $ which is induced by $ (\alpha, \beta) $, so that $ \rootgr{\alpha} = \rootgr{1} $ and $ \rootgr{\beta} = \rootgr{n} $. Let $ u \in \invset{1} $, $ v \in \invset{n} $ and put $ w_1 \defl \weylmap{}(u) $, $ w_0 \defl \weylmap{}(v) $. Then we are in the situation of~\ref{tw:braid:setup}, and we have to show that $ \braidword{n}(w_1, w_0) = \braidword{n}(w_0, w_1) $.
	
	At first, assume that $ n $ is even and put $ l \defl n/2 \in \Npos $. Then it follows from \cref{tw:braid:cl3} that
	\[ w_0 = w_n = w_{2l} = w_0^{w_1 (w_0 w_1)^{l-1}} = \brackets[\big]{w_1 (w_0 w_1)^{l-1}}^{-1} \cdot w_0 \cdot w_1 (w_0 w_1)^{l-1}. \]
	Multiplying this term from the left side with $ w_1 (w_0 w_1)^{l-1} $, we infer that
	\begin{align*}
		w_1 (w_0 w_1)^{l-1} w_0 =  w_0 w_1 (w_0 w_1)^{l-1}.
	\end{align*}
	In other words, $ \braidword{n}(w_1, w_0) = \braidword{n}(w_0, w_1) $, as desired.

	Now assume that $ n $ is odd, so that $ n=2l+1 $ for some $ l \in \IN_{\ge 1} $. Then
	\[ w_0 = w_n = w_{2l+1} = w_1^{(w_0 w_1)^l} = \brackets[\big]{(w_0 w_1)^l}^{-1} \cdot w_1 \cdot (w_0 w_1)^l \]
	by \cref{tw:braid:cl3}. Multiplying from the left side with $ (w_0 w_1)^l $, we infer that
	\[ (w_0 w_1)^l w_0 = w_1 (w_0 w_1)^l. \]
	Again, this says precisely that $ \braidword{n}(w_1, w_0) = \braidword{n}(w_0, w_1) $, which finishes the proof of the braid relations.
\end{proof}


\section{Closed and Ordered Sets of Roots}

\label{sec:closed}

\begin{secnotation}
	We denote by $ \roots $ a root system and we will frequently use \cref{rootorder:cry-brackets-note}.
\end{secnotation}

Our next goal is to study under which conditions the product map on a positive system in $ \roots $ is bijective. The notion of closed sets of roots, which we have already introduced in \cref{rootsys:closed-def}, will be indispensable in this context. The current section is dedicated to a proper investigation of these objects. This can be done on a purely combinatorial, root-system-theoretic level without reference to any $ \roots $-graded group. We could have already covered this subject in \cref{chap:pre}, but we decided not to because it is less standard than the other topics of \cref{chap:pre}. Further, the motivation is more obvious after having introduced groups with $ \roots $-commutator relations and product maps.

Apart from some basic properties, our goal is to show that every positive system of $ \roots $ has an extremal ordering (except when $ \roots $ is not reduced, in which case we have to remove roots until only one root from each ray remains). This is \cref{rootorder:extremal-exists}. In \cref{prodmap:extremal-bij}, we will show that under a suitable condition on a group $ G $ with $ \roots $-commutator relations, the product map for any extremal ordering is bijective. 

None of the material in this section is new or surprising, but the literature on this topic can be confusing because no single reference provides all the facts that we will list here. Further, there are subtle differences in the used definitions. This concerns, in particular, the case of non-reduced root systems. Some references on this topic are \cite[Section~VI.1.7]{BourbakiLie46} and \cite[Lemmas~16--18]{Steinberg-ChevGroups}.

\begin{definition}[Ideals]\label{rootorder:ideal-def}
	Let $ \rootsub $ be a subset of $ \roots $. A subset $ I $ of $ \rootsub $ is called an \defemph*{ideal of $ \rootsub $}\index{ideal (in root systems)} if $ \rootintcox{\alpha}{\beta} \subs I $ for all non-proportional $ \alpha \in I $ and $ \beta \in \rootsub $, and it is called a \defemph*{crystallographic ideal of $ \rootsub $}\index{ideal (in root systems)!crystallographic} if $ \rootint{\alpha}{\beta} \subs I $ for all non-proportional $ \alpha \in I $ and $ \beta \in \rootsub $.
\end{definition}

\begin{warning}
	The terminology \enquote{crystallographic ideal} makes it sound like being crystallographic is a property that an ideal may or may not have, but it is actually the other way around: Every ideal is also a crystallographic ideal, but not every crystallographic ideal is an ideal.
\end{warning}

\begin{example}
	For any element $ w $ of the Weyl group and any positive system $ \possys $ in $ \roots $, the set $ \switchset(w) \defl \switchset_\possys(w) \defl \Set{\alpha \in \indivset{\possys} \given \alpha^{w^{-1}} \in -\possys} $ from \cref{rootsys:switchset-def} is closed.
\end{example}

The same remarks as in \cref{cry-note} apply to (crystallographically) closed set and to (crystallographic) ideals: The crystallographic properties are well-defined even in the non-crystallographic setting, but they are only useful in the crystallographic setting.

\begin{convention}\label{rootorder:cry-brackets-note}
	In the following, we will often put the words \enquote{crystallographic} and \enquote{crystallographically} in brackets. Any assertion of this form should be interpreted as two statements, once with and once without \emph{all} words in brackets. For example, the following \cref{rootorder:ideal-basic} asserts that any ideal is closed and that any crystallographic ideal is crystallographically closed. On the contrary, it does not assert that any crystallographic ideal is closed. By the end of this section, it should be clear that the underlying principle is always the same. In \cref{rootsys-with-int}, we will indicate how the crystallographic and the non-crystallographic theory could be treated in a more uniform way.
\end{convention}

\begin{remark}\label{rootorder:ideal-basic}
	Let $ \rootsub $ be a subset of $ \roots $ and let $ I $ be a subset of $ \rootsub $. If $ I $ is a (crystallographic) ideal of $ \rootsub $, then clearly $ I $ is (crystallographically) closed. Conversely, if $ I = \rootsub \setminus \compactSet{\alpha} $ for some $ \alpha \in \rootsub $ and $ \rootsub $ is (crystallographically) closed, then $ I $ is automatically a (crystallographic) ideal of $ \rootsub $ if $ I $ is (crystallographically) closed. Further, $ \rootsub $ is a (crystallographic) ideal of $ \rootsub $ if and only if $ \rootsub $ is (crystallographically) closed.
\end{remark}

\begin{lemma}[{\cite[Lemma~9.4]{HumphreysLieAlg}}]\label{rootorder:pair-cry-lem}
	Assume that $ \roots $ is crystallographic and let $ \alpha, \beta $ be non-proportional roots. If $ \alpha \cdot \beta > 0 $, then $ \alpha-\beta $ is a root. If $ \alpha \cdot \beta < 0 $, then $ \alpha+\beta $ is a root.
\end{lemma}

\begin{lemma}\label{rootorder:root-minus-lem}
	Assume that $ \roots $ is crystallographic. Let $ \alpha, \beta $ be non-proportional roots and let $ i,j \in \Npos $ such that $ i\alpha + j\beta $ is a root. Then either $ (i-1)\alpha + j\beta $ or $ i\alpha + (j-1)\beta $ is a root.
\end{lemma}
\begin{proof}
	Set $ \gamma \defl i\alpha + j\beta \in \roots $. If $ \gamma \cdot \alpha \le 0 $ and $ \gamma \cdot \beta \le 0 $, then
	\[ \gamma \cdot \gamma = \gamma \cdot (i\alpha + j\beta) = i(\gamma \cdot \alpha) + j(\gamma \cdot \beta) \le 0 \]
	which contradicts the fact that $ \gamma $ is a root. Hence $ \gamma \cdot \alpha > 0 $ or $ \gamma \cdot \beta > 0 $, so the assertion follows from \cref{rootorder:pair-cry-lem}.
\end{proof}

The following result shows that the condition in the definition of crystallographic ideals can be slightly weakened.

\begin{lemma}\label{rootorder:cry-ideal-char}
	Assume that $ \roots $ is crystallographic. Let $ \rootsub $ be a crystallographically closed subset of $ \roots $ and let $ I $ be a subset of $ \rootsub $ with the property that for all non-proportional $ \alpha \in I $ and $ \beta \in \rootsub $ for which $ \alpha+\beta $ is a root, $ \alpha + \beta $ lies in $ I $. Then $ I $ is a crystallographic ideal of $ \rootsub $.
\end{lemma}
\begin{proof}
	Let $ \alpha \in I $, $ \beta \in \rootsub $ and $ \gamma \in \rootint{\alpha}{\beta} $. Then there exist $ i,j \in \Npos $ such that $ \gamma = i\alpha + j\beta $. By induction on $ i+j $, we show that $ \gamma $ lies in $ I $. If $ i+j=2 $, then $ \gamma = \alpha + \beta $ and so $ \gamma $ lies in $ I $ by the assumption on $ I $. Now assume that $ i+j \ge 3 $. By \cref{rootorder:root-minus-lem}, we have $ (i-1)\alpha + j\beta \in \roots $ or $ i\alpha + (j-1)\beta \in \roots $. The induction hypothesis now implies that $ (i-1)\alpha + j\beta \in I $ or $ i\alpha + (j-1)\beta \in I $. In the first case, it follows that
	\[ \gamma = \brackets[\big]{(i-1)\alpha + j\beta} + \alpha \in I \]
	because $ \gamma \in \roots $, $ (i-1)\alpha + j\beta \in I $ and $ \alpha \in \rootsub $. Similarly, in the second case it follows that
	\[ \gamma = \brackets[\big]{i\alpha + (j-1)\beta} + \beta \in I. \]
	This finishes the proof.
\end{proof}

We make a brief detour into a slightly different notion of closedness, which will be needed in the proof of \cref{rootorder:clos-in-pos} (and only there). Alternatively, we could prove \cref{rootorder:clos-in-pos} by citing \cite[1.14]{LoosNeherBook}, but this would require a similar amount of notational setup.

\begin{definition}[Bourbaki-closed]
	Assume that $ \roots $ is crystallographic. We say that a subset $ \rootsub $ of $ \roots $ is \defemph*{Bourbaki-closed}\index{closed set of roots!Bourbaki-} if it is closed in the sense of \cite[Définition~VI.1.4, \onpage{160}]{BourbakiLie46}, which means that for all $ \alpha, \beta \in \rootsub $ with $ \alpha + \beta \in \roots $, we have $ \alpha + \beta \in \rootsub $.
\end{definition}

\begin{lemma}\label{rootorder:bourbaki-clos-lem}
	Assume that $ \roots $ is crystallographic and let $ \rootsub $ be a subset of $ \roots $. Then the following hold:
	\begin{lemenumerate}
		\item \label{rootorder:bourbaki-clos-lem:is-cry}If $ \rootsub $ is Bourbaki-closed, then $ \rootsub $ is crystallographically closed.
		
		\item \label{rootorder:bourbaki-clos-lem:iff}Assume that $ \rootsub $ is crystallographically closed. Then $ \rootsub $ is Bourbaki-closed if and only if it satisfies the following conditions for all roots $ \alpha $ for which $ 2\alpha $ is also a root: If $ \alpha \in \rootsub $ then $ 2\alpha \in \rootsub $; if $ 2\alpha, -\alpha \in \rootsub $ then $ \alpha \in \rootsub $.
		
		\item \label{rootorder:bourbaki-clos-lem:red}If $ \roots $ is reduced, then $ \rootsub $ is crystallographically closed if and only if it is Bourbaki-closed.
		
		\item \label{rootorder:bourbaki-clos-lem:subset}Assume that $ \rootsub $ is crystallographically closed and that $ -\lambda \alpha \nin \rootsub $ for all $ \alpha \in \rootsub $ and $ \lambda > 0 $. Then $ \rootsub $ is contained in a Bourbaki-closed set $ \rootsub' $ with $ \rootsub' \intersect (-\rootsub') = \emptyset $.
	\end{lemenumerate}
\end{lemma}
\begin{proof}
	By \cref{rootorder:cry-ideal-char}, any Bourbaki-closed set $ \rootsub $ is a crystallographic ideal of itself, which by \cref{rootorder:ideal-basic} means precisely that $ \rootsub $ is crystallographically closed. This proves~\itemref{rootorder:bourbaki-clos-lem:is-cry}. However, a crystallographically closed set $ \rootsub $ is not necessarily Bourbaki-closed because the roots $ \alpha $, $ \beta $ in \cref{rootsys:closed-def} are required to be non-proportional. Thus a crystallographically closed set $ \rootsub $ is Bourbaki-closed if and only if for all roots $ \alpha \in \rootsub $ and for all $ \lambda \in \IR $ such that $ \lambda \alpha \in \rootsub $ and $ (1+\lambda)\alpha \in \roots $, we have $ (1+\lambda) \alpha \in \rootsub $. By \cref{rootsys:cry-nonred-2}, we can restrict to $ \lambda \in \Set{1/2, 1, 2} $, which yields precisely the statement of~\itemref{rootorder:bourbaki-clos-lem:iff}. Assertion~\itemref{rootorder:bourbaki-clos-lem:red} follows from~\itemref{rootorder:bourbaki-clos-lem:iff}. In~\itemref{rootorder:bourbaki-clos-lem:subset}, we can take $ \rootsub' \defl (\rootsub \union 2\rootsub) \intersect \roots $.
\end{proof}

\begin{proposition}\label{rootorder:clos-in-pos}
	Assume that $ \roots $ is crystallographic and let $ \rootsub $ be a crystallographically closed subset of $ \roots $ such that for all $ \alpha \in \rootsub $, we have $ -\lambda \alpha \nin \rootsub $ for all $ \lambda > 0 $. (If $ \roots $ is reduced, this simply means that $ \rootsub \intersect (-\rootsub) = \emptyset $.) Then there exists a positive system $ \possys $ in $ \roots $ which contains $ \rootsub $.
\end{proposition}
\begin{proof}
	This is proven in \cite[Proposition~VI.1.22, \onpage{163}]{BourbakiLie46} under the slightly modified assumptions that $ \rootsub $ is Bourbaki-closed and $ \rootsub \intersect (-\rootsub) = \emptyset $. We have proved in \thmitemcref{rootorder:bourbaki-clos-lem}{rootorder:bourbaki-clos-lem:subset} that $ \rootsub $ is contained in a set $ \rootsub' $ with these properties, so the assertion follows.
\end{proof}

\begin{definition}[Extremal root]
	Let $ \rootsub $ be a subset of $ \roots $ and let $ \alpha \in \rootsub $. Then $ \alpha $ is called \defemph*{extremal in $ \rootsub $}\index{root!extremal} if there exists a point $ p $ in the Euclidean space $ (V, \cdot) $ surrounding $ \roots $ which is not orthogonal to any root in $ \roots $ such that $ \alpha \cdot p > 0 $ and $ \beta \cdot p <0 $ for all $ \beta \in \rootsub \setminus \compactSet{\alpha} $.
\end{definition}

\begin{remark}
	Let $ \rootsub $ be a subset of $ \roots $. Then a root $ \alpha \in \rootsub $ is extremal in $ \rootsub $ if and only if there exists a positive system $ \possys $ in $ \roots $ which contains $ \rootsub \setminus \compactSet{\alpha} $ but not~$ \alpha $.
\end{remark}

\begin{definition}[Root ordering]
	Let $ \rootsub $ be a subset of $ \roots $. An \defemph*{ordering of $ \rootsub $}\index{root ordering} is a tuple $ (\alpha_1, \ldots, \alpha_m) $ such that $ \listing{\alpha}{m} $ are pairwise distinct and $ \rootsub = \Set{\listing{\alpha}{m}} $.
\end{definition}

\begin{definition}[Properties of orderings]\label{rootorder:prop-def}
	Let $ \rootsub $ be a subset of $ \roots $ and let $ \word{\alpha} = (\listing{\alpha}{m}) $ be an ordering of $ \rootsub $.
	\begin{defenumerate}
		\item $ \word{\alpha} $ is called \defemph*{extremal}\index{root ordering!extremal} if $ \alpha_i $ is extremal in $ \Set{\alpha_i, \ldots, \alpha_m} $ for all $ i \in \numint{1}{m} $.
		
		\item $ \word{\alpha} $ is called \defemph*{(crystallographically) normal}\index{root ordering!normal}\index{root ordering!crystallographically normal} if $ \Set{\alpha_i, \ldots, \alpha_m} $ is a (crystallographic) ideal of $ \rootsub $ for all $ i \in \numint{1}{m} $.
		
		\item $ \word{\alpha} $ is called \defemph*{(crystallographically) subnormal}\index{root ordering!subnormal}\index{root ordering!crystallographically subnormal} if $ \rootsub $ is (crystallographically) closed and $ \Set{\alpha_{i+1}, \ldots, \alpha_m} $ is a (crystallographic) ideal of $ \Set{\alpha_i, \ldots, \alpha_m} $ for all $ i \in \numint{1}{m-1} $.
		
		\item $ \word{\alpha} $ is called a \defemph*{height ordering}\index{root ordering!height ordering} if $ \rootsub $ is contained in some positive system $ \possys $ and $ \rootht(\alpha_1) \le \cdots \le \rootht(\alpha_m) $ where $ \rootht $ denotes the height function with respect to~$ \possys $.
	\end{defenumerate}
\end{definition}

We will now investigate various basic properties of the notions defined in \cref{rootorder:prop-def}.

\begin{remark}[Basic properties of extremal orderings]\label{rootorder:ext-prop}
	Let $ \word{\alpha} = \tup{\alpha}{m} $ be an ordering of a subset $ \rootsub $ of $ \roots $.
	\begin{remenumerate}
		\item If $ \word{\alpha} $ is extremal, then there cannot exist distinct $ i,j \in \numint{1}{m} $ and $ \lambda > 0 $ such that $ \alpha_i = \lambda \alpha_j $.
		
		\item \label{rootorder:ext-prop:scale}Let $ i \in \numint{1}{m} $ and $ \lambda \in \IR_{>0} $ such that $ \lambda \alpha_i $ is a root. Then $ \tup{\alpha}{m} $ is extremal if and only if $ (\alpha_1, \ldots, \lambda \alpha_i, \ldots, \alpha_m) $ is extremal.
		
		\item \label{rootorder:ext-prop:ortho}Let $ \psi $ be an orthogonal automorphism of the Euclidean space $ (V, \cdot) $ surrounding $ \roots $. Then $ \tup{\alpha}{m} $ is an extremal ordering of $ \rootsub $ if and only if $ (\alpha_1^\psi, \ldots, \alpha_m^\psi) $ is an extremal ordering of $ \rootsub^\psi $.
		
		\item \label{rootorder:ext-prop:delete}Any tuple which is obtained from $ \word{\alpha} $ by deleting an arbitrary number of entries is also extremal.
	\end{remenumerate}
\end{remark}

\begin{remark}[Basic properties of (sub-) normal orderings]\label{rootorder:subnorm-prop}
	Let $ \rootsub $ be a subset of $ \roots $ and let $ \word{\alpha} = (\listing{\alpha}{m}) $ be an ordering of $ \rootsub $.
	\begin{remenumerate}
		\item If $ \word{\alpha} $ is (crystallographically) normal, then it is also (crystallographically) subnormal.
		
		\item If $ \word{\alpha} $ is (sub-) normal, then it is also crystallographically (sub-) normal. 
		
		\item If $ \word{\alpha} $ is (sub-) normal, then for all for all $ i \in \numint{1}{m} $, the tuple $ (\alpha_i, \ldots, \alpha_m) $ is a (sub-) normal ordering of $ \Set{\alpha_i, \ldots, \alpha_m} $. The same assertion holds for the crystallgraphic properties.
		
		\item If $ \word{\alpha} $ is (crystallographically) normal, then $ \rootsub $ is a (crystallographic) ideal of itself, and so it is (crystallographically) closed by \cref{rootorder:ideal-basic}.
		
		\item \label{rootorder:subnorm-prop:subnorm-crit}It follows from \cref{rootorder:ideal-basic} that $ \word{\alpha} $ is (crystallographically) subnormal if and only if all sets $ \Set{\alpha_i, \ldots, \alpha_m} $ for $ i \in \numint{1}{m} $ are (crystallographically) closed.
	\end{remenumerate}
	Observe that we had to require in the definition of (crystallographically) subnormal orderings that $ \rootsub $ is (crystallographically) closed, which was not necessary for (crystallographically) normal orderings.
\end{remark}

\begin{example}[of a subnormal extremal ordering]\label{rootorder:interval-subnormal}
	Let $ S = \indivset{\clrootintcox{\alpha}{\beta}} $ for some indivisible roots $ \alpha $, $ \beta $ and let $ \word{\alpha} = \tup{\alpha}{k} $ be an interval ordering of $ S $ in the sense of \cref{rootsys:clockwise-def}. Then it follows from \thmitemcref{rootorder:subnorm-prop}{rootorder:subnorm-prop:subnorm-crit} that $ \word{\alpha} $ is subnormal. Further, it is clear that $ \word{\alpha} $ is extremal.
\end{example}

\begin{lemma}\label{rootorder:height-is-normal}
	Assume that $ \roots $ is crystallographic, let $ \rootsub $ be a crystallographically closed subset of $ \roots $ and let $ \word{\alpha} = \tup{\alpha}{m} $ be a height ordering of $ \rootsub $. Then $ \word{\alpha} $ is a crystallographically normal ordering of $ \rootsub $.
\end{lemma}
\begin{proof}
	Let $ i \in \numint{1}{m} $. We have to show that $ I \defl \Set{\alpha_i, \ldots, \alpha_m} $ is a crystallographic ideal of $ \rootsub $. By the definition of height orderings, there exists a positive subsystem $ \possys $ of $ \roots $ which contains $ \rootsub $ such that $ \rootht(\alpha_1) \le \cdots \le \rootht(\alpha_m) $. By \cref{rootorder:cry-ideal-char}, it suffices to show that for non-proportional $ \alpha \in \rootsub $ and $ \beta \in I $ for which $ \alpha+\beta $ is a root, $ \alpha+\beta $ lies in $ I $. Since $ \rootsub $ is crystallographically closed, $ \alpha+\beta $ lies in $ \rootsub $, so there exists $ j \in \numint{1}{m} $ such that $ \alpha+\beta = \alpha_j $. Now $ \rootht(\alpha_j) = \rootht(\alpha) + \rootht(\beta) > \rootht(\beta) \ge \rootht(\alpha_i) $, so $ j>i $. This implies that $ \alpha_j $ lies in $ I $, which finishes the proof.
\end{proof}

\begin{lemma}\label{rootorder:extremal-subnormal}
	Let $ \rootsub $ be a (crystallographically) closed subset of $ \roots $. Then every extremal ordering $ \word{\alpha} = \tup{\alpha}{m} $ of $ \rootsub $ is (crystallographically) subnormal, and for all $ i \in \numint{1}{m} $, the set $ \Set{\alpha_i, \ldots, \alpha_m} $ is (crystallographically) closed.
\end{lemma}
\begin{proof}
	Since $ (\alpha_2, \ldots, \alpha_m) $ is also extremal by \thmitemcref{rootorder:ext-prop}{rootorder:ext-prop:delete}, it suffices by induction to show that $ \Set{\alpha_2, \ldots, \alpha_m} $ is (crystallographically) closed and a (crystallographic) ideal of $ \rootsub $. By \cref{rootorder:ideal-basic}, we only have to prove the first statement because it implies the second one. Let $ i,j \in \numint{2}{m} $ such that $ \alpha_i, \alpha_j $ are non-proportional. Write $ I \defl \rootint{\alpha_i}{\alpha_j} $ for the proof of the crystallographic assertion and $ I \defl \rootintcox{\alpha_i}{\alpha_j} $ for the proof of the non-crystallographic assertion. Since $ \rootsub $ is (crystallographically) closed, we have $ I \subs \rootsub $, so it suffices to show that $ \alpha_1 $ is not contained in $ I $. Assume that $ \alpha_1 = \lambda \alpha_i + \mu \alpha_j $ for some $ \lambda, \mu \in \IR_{>0} $ (where $ \lambda, \mu \in \Npos $ for the crystallographic assertion). Since $ \alpha_1 $ is extremal in $ \rootsub $, there exists a point $ p $ in the surrounding Euclidean space of $ \roots $ such that $ \alpha_1 \cdot p > 0 $ and $ \alpha_i \cdot p, \alpha_j \cdot p < 0 $. Then
	\[ 0 < \alpha_1 \cdot p = \lambda (\alpha_i \cdot p) + \mu(\alpha_j \cdot p) < 0, \]
	which is a contradiction.
\end{proof}

We now prove some existence results concerning extremal orderings. The following result is essentially a stronger version of \cref{rootsys:switchset-description}: It provides not only a description of the set $ \switchset(w) $, but even an extremal ordering of $ \switchset(w) $.

\begin{proposition}\label{rootorder:switchset-extremal}
	Choose a root base $ \rootbase $ of $ \roots $ and denote by $ \possys $ the corresponding positive system. Let $ w $ be an element of the Weyl group of $ \roots $, let $ \word{\delta} = \tup{\delta}{k} $ be a $ \rootbase $-expression of $ w $ and let $ \word{\alpha} = (\alpha_k, \ldots, \alpha_1) $ be the inverse associated root sequence from \cref{rootsys:rootseq-def}. Then $ \word{\alpha} $ is an extremal ordering of the set $ \switchset(w) = \Set{\beta \in \indivset{\possys} \given \beta^{w^{-1}} \in -\possys} $ from \cref{rootsys:switchset-def}.
\end{proposition}
\begin{proof}
	If $ k=0 $, then $ w = 1_W $ and there is nothing to prove. We proceed by induction and assume that $ k>0 $. Put $ \delta \defl \delta_k $, $ \word{\delta}' \defl \tup{\delta}{k-1} $ and $ w' \defl \refl{\word{\delta}'} $. Then $ \word{\delta}' $ is a reduced expression of $ w' $ and $ \len(w' \refl{\delta}) = \len(w) > \len(w') $. Denote the root sequence of $ \word{\delta}' $ by $ \tup{\beta}{k-1} $. By the induction hypothesis, $ (\beta_{k-1}, \ldots, \beta_1) $ is an extremal ordering of $ \switchset(w') $. Thus by \thmitemcref{rootorder:ext-prop}{rootorder:ext-prop:ortho}, $ (\beta_{k-1}^{\reflbr{\delta}}, \ldots, \beta_1^{\reflbr{\delta}}) $ is an extremal ordering of $ \switchset(w')^{\reflbr{\delta}} $. Recall from \cref{rootsys:rootseq-iter} that the root sequence of $ \word{\delta} $ is $ (\beta_1^{\reflbr{\delta}}, \ldots, \beta_{k-1}^{\reflbr{\delta}}, \delta) $ and from \cref{rootsys:switchset-lem} that $ \switchset(w) = \switchset(w')^{\reflbr{\delta}} \disjunion \compactSet{\delta} $. Hence it remains to show that $ \delta $ is extremal in $ \compactSet{\delta} \union \switchset(w')^{\reflbr{\delta}} $. Equivalently, we have to show that $ \delta^{\reflbr{\delta}} $ is extremal in $ \compactSet{\delta^{\reflbr{\delta}}} \union \switchset(w') $, which holds by \cref{rootsys:switchset-complete}. This finishes the proof.
\end{proof}

\begin{proposition}\label{rootorder:switchset-extremal-possys}
	Choose a root base $ \rootbase $ of $ \roots $ and denote by $ \possys $ the corresponding positive system. Put $ m \defl \abs{\indivset{\possys}} $, let $ \word{\delta} = \tup{\delta}{m} $ be any reduced expression of the longest element $ \rho $ of $ W $ with respect to $ \rootbase $ and denote by $ \word{\alpha} $ the inverse root sequence associated to $ \word{\delta} $. Let $ \rootsub $ be any subset of $ \indivset{\possys} $ and let $ \word{\beta} $ be the tuple which is obtained from $ \word{\alpha} $ by deleting all entries in $ \possys \setminus \rootsub $. Then $ \word{\beta} $ is an extremal ordering of $ \rootsub $.
\end{proposition}
\begin{proof}
	Since $ \switchset(\rho) = \indivset{\possys} $ by \cref{rootsys:longest-el-char}, it is a consequence of \cref{rootorder:switchset-extremal} that $ \word{\alpha} $ is an extremal ordering of $ \indivset{\possys} $. Thus by \thmitemcref{rootorder:ext-prop}{rootorder:ext-prop:delete}, $ \word{\beta} $ is an extremal ordering of $ \rootsub $.
\end{proof}

\begin{proposition}\label{rootorder:extremal-exists}
	Let $ \possys $ be a positive system in $ \roots $ and let $ \rootsub $ be a subset of $ \possys $ consisting of pairwise non-proportional roots. Then there exists an extremal ordering of $ \rootsub $.
\end{proposition}
\begin{proof}
	For any root $ \alpha \in \rootsub $, there exists a unique indivisible root $ \alpha' $ (which might be equal to $ \alpha $) such that $ \IR\alpha = \IR\alpha' $. Then $ \rootsub' \defl \Set{\alpha' \given \alpha \in \rootsub} $ is a subset of $ \indivset{\possys} $. Thus \cref{rootorder:switchset-extremal-possys} yields an extremal ordering $ \word{\alpha}' $ of $ \rootsub $. Scaling the roots in $ \word{\alpha}' $ appropriately and applying \thmitemcref{rootorder:ext-prop}{rootorder:ext-prop:scale}, we obtain an extremal ordering of~$ \rootsub $.
\end{proof}

\begin{note}\label{rootsys-with-int}
	Continuing the remarks in \cref{rootorder:cry-brackets-note}, we briefly describe how the crystallographic and the non-crystallographic considerations in this section (and in the rest of this book) could be treated in a more uniform way. Define a \defemph*{root system with intervals}\index{root system!with intervals} to be a tuple $ \Psi = (\roots, (\varrootint{\alpha}{\beta})_{\alpha,\beta}) $ where $ \roots $ is a root system in the regular sense and $ (\varrootint{\alpha}{\beta})_{\alpha, \beta} $ is a family of subsets of $ \roots $ where $ \alpha, \beta $ runs over all pairs of non-proportional roots in $ \roots $. We require that $ \varrootint{\alpha}{\beta} \subs \rootintcox{\alpha}{\beta} $ for all non-proportional roots $ \alpha, \beta $. Now a \defemph*{group with $ \Psi $-commutator relations} is defined as in \cref{rgg:group-commrel-def}, except that we must have $ \commutator{\rootgr{\alpha}}{\rootgr{\beta}} \subs \rootgr{\varrootint{\alpha}{\beta}} $. Similarly, the notions of ideals and closed sets of roots are defined with respect to the sets $ \varrootint{\alpha}{\beta} $.
	
	Now let $ \roots $ be a root system. There are only two root systems with intervals corresponding to $ \roots $ that are interesting:
	\[ \tilde{\roots} \defl \brackets[\big]{\roots, (\rootintcox{\alpha}{\beta})_{\alpha, \beta}} \midand \tilde{\roots}^{\text{cry}} \defl \brackets[\big]{\roots, (\rootint{\alpha}{\beta})_{\alpha, \beta}}. \]
	Note that in this setup, the information \enquote{whether we use the crystallographic terminology or not} is part of the datum of a root system with intervals. For example, \enquote{$ C_3 $ regarded as a root system in the Coxeter sense} and \enquote{$ C_3 $ regarded as a crystallographic root systems} are two distinct objects, namely $ \tilde{C}_3 $ and $ \tilde{C}_3^{\text{cry}} $. Thus $ C_3 $-graded groups are precisely $ \tilde{C}_3 $-graded groups and crystallographic $ C_3 $-graded groups are precisely $ \tilde{C}_3^{\text{cry}} $-graded groups.
	
	Many (though not all) results in this book which concern $ \roots $-graded groups for arbitrary root systems $ \roots $ can be stated and proven for $ \Psi $-graded groups where $ \Psi = (\roots, (\varrootint{\alpha}{\beta})_{\alpha,\beta}) $ is an arbitrary root system with intervals. In this setup, root gradings, root intervals, ideals and closed sets of roots are understood to be defined with respect to the sets $ \varrootint{\alpha}{\beta} $. Thus the crystallographic and the non-crystallographic case would not have to be considered separately most of the time. There are only a few results for which this strategy would not work, such as \cref{prodmap:ideal-normal}.
	
	Despite the advantages of this approach, we will not use the language of root systems with intervals in the sequel. The reason for this is that the regular notion of root systems is well-established and that the distinction between the crystallographic and the non-crystallographic case is only a minor nuisance.
\end{note}


\section{Bijectivity of the Product Map}

\label{sec:prod-bij}

\begin{secnotation}
	We denote by $ \roots $ a root system and by $ G $ a group with $ \roots $-commutator relations with root groups $ (\rootgr{\alpha})_{\alpha \in \roots} $. Further, we will frequently use \cref{rootorder:cry-brackets-note}.
\end{secnotation}

\begin{definition}
	Let $ \rootsub $ be a subset of $ \roots $, let $ \word{\alpha} = (\listing{\alpha}{m}) $ be an ordering of $ \rootsub $ and let $ G $ be a group with $ \roots $-commutator relations with root groups $ (\rootgr{\alpha})_{\alpha \in \roots} $. Then $ \word{\alpha} $ is called \defemph*{$ G $-injective}\index{root ordering!injective}, \defemph*{$ G $-surjective}\index{root ordering!surjective} or \defemph*{$ G $-bijective}\index{root ordering!bijective} if the product map on $ \word{\alpha} $ is injective, surjective or bijective, respectively.
\end{definition}

The goal of this section is to find a criterion which guarantees the existence of a $ G $-bijective ordering on a (crystallographically) closed subset $ \rootsub $ of $ \roots $. In the crystallographic setting, we will even find that every ordering of $ \rootsub $ is $ G $-bijective.

\begin{note}[The product map in Chevalley groups]
	Assume that $ \roots $ is crystallographic and reduced and that $ G $ is a Chevalley group of type $ \roots $ (which implies that $ G $ has crystallographic $ \roots $-commutator relations). Let $ \rootsub $ be a crystallographically closed subset of $ \roots $ such that $ \rootsub \intersect (-\rootsub) = \emptyset $. It is shown in \cite[Lemma~17]{Steinberg-ChevGroups} that every height ordering of $ \rootsub $ is $ G $-bijective. Further, it is shown in \cite[Lemma~18]{Steinberg-ChevGroups} that the existence of a crystallographically normal $ G $-bijective ordering of $ \rootsub $ implies that every ordering of $ \rootsub $ is $ G $-bijective. Since every height ordering is crystallographically normal by \cref{rootorder:height-is-normal}, it follows that every ordering of $ G $ is bijective.
\end{note}

Our approach in this section is motivated by the strategy for Chevalley groups outlined above, but some of the arguments do not carry over to our more general setting. We begin this section with an investigation of product maps in arbitrary groups $ H $. This includes Steinberg's result that the existence of a crystallographically normal $ G $-bijective ordering of $ \rootsub $ implies that every ordering of $ \rootsub $ is $ G $-bijective (\cref{prodmap:abstract-normal-bijective}). Unfortunately, this result is not applicable in our situation, and we state it purely for completeness. In its place, we will use \cref{prodmap:tits}, which is due to Tits and which says that a similar assertion holds if $ G $ has a central series with certain properties. \Cref{prodmap:tits} came to our attention in the form of \cite[3.10]{LoosNeherBook}.

In the second part of this section, we prove the existence of a $ G $-bijective ordering of $ \rootsub $ (where $ \rootsub $ is any subset of a positive system of $ \roots $ whose roots are pairwise non-proportional). By the same arguments as in \cite[Lemma~17]{Steinberg-ChevGroups}, every subnormal ordering of $ \rootsub $ is $ G $-surjective (\cref{prodmap:subnormal-surj}). However, Steinberg's proof of $ G $-injectivity uses the natural action of a Chevalley group on a certain module, so it cannot be generalised to $ \roots $-graded groups. Instead, we will show in \cref{prodmap:extremal-inj} that every extremal ordering of $ \rootsub $ is $ G $-injective. Since extremal orderings exist by \cref{rootorder:extremal-exists} and are subnormal by \cref{rootorder:extremal-subnormal}, we infer that there exists a $ G $-bijective ordering of $ \rootsub $ (\cref{prodmap:bij-ex}). In general, this existence result is the best we can do because Steinberg's \cref{prodmap:abstract-normal-bijective} does not hold for subnormal orderings and the conditions of Tits' \cref{prodmap:tits} are not satisfied for general $ \roots $-gradings. However, Tits' lemma does apply in the crystallographic setting, so we infer that every ordering of $ \rootsub $ is $ G $-bijective in this case.

\subsection{Purely Group-theoretic Observations}

\begin{lemma}\label{prodmap:switch-order}
	Let $ H $ be a group, let $ m \in \IN_{\ge 1} $ and let $ \listing{U}{m} $ be subgroups of $ H $. Then the product map $ \map{\mu_1}{U_1 \times \cdots \times U_m}{H}{}{} $ is injective, surjective or bijective if and only if the product map $ \map{\mu_2}{U_m \times \cdots \times U_1}{H}{}{} $ is injective, surjective or bijective, respectively.
\end{lemma}
\begin{proof}
	Let $ g_i \in U_i $ for all $ i \in \numint{1}{m} $. Then we have
	\[ \mu_1(g_1, \ldots, g_m) = \mu_2(g_m^{-1}, \ldots, g_1^{-})^{-1}. \]
	The assertion follows.
\end{proof}

\begin{lemma}\label{prodmap:abstract-subnormal-surjective}
	Let $ H $ be a group, let $ m \in \IN_{\ge 1} $ and let $ \listing{U}{m} $ be subgroups of $ H $ such that $ U_1 \union \cdots \union U_m $ generates $ H $. Put $ V_i \defl \gen{U_j \given j \in \numint{i}{m}} $ for all $ i \in \numint{1}{m} $. If $ V_{i+1} $ is normal in $ V_i $ for all $ i \in \numint{1}{m-1} $, then the product map $ \map{}{U_1 \times \cdots \times U_m}{H}{}{} $ is surjective.
\end{lemma}
\begin{proof}
	For $ m=1 $, there is nothing to prove. Now assume that $ m=2 $. Since $ V_2 $ is normal in $ V_1 $, we have that $ U_1 $ normalises $ U_2 $. Hence by \thmitemcref{group-rel}{group-rel:conj-comm}, $ U_1 U_2 $ is a subgroup of $ H $. In fact, $ U_1 U_2 = H $ because $ U_1 \union U_2 $ generates $ H $. Thus the assertion holds for $ m=2 $. The assertion for arbitrary $ m $ follows by induction.
\end{proof}

For completeness, here is the lemma used in \cite{Steinberg-ChevGroups} to deduce $ G $-bijectivity of arbitrary orderings.

\begin{lemma}[{\cite[Lemma~18]{Steinberg-ChevGroups}}]\label{prodmap:abstract-normal-bijective}
	Let $ H $ be a group, let $ m \in \IN_{\ge 1} $ and let $ \listing{U}{m} $ be subgroups of $ H $ such that $ U_1 \union \cdots \union U_m $ generates $ H $. Put $ V_i \defl \gen{U_j \given j \in \numint{i}{m}} $ for all $ i \in \numint{1}{m} $. Assume that $ V_i $ is normal in $ H $ for all $ i \in \numint{1}{m} $ and that the product map $ \map{}{U_1 \times \cdots \times U_m}{H}{}{} $ is injective. Then for all permutations $ \map{\sigma}{\numint{1}{m}}{\numint{1}{m}}{}{} $, the product map $ \map{}{U_{\sigma(1)} \times \cdots \times U_{\sigma(m)}}{H}{}{} $ is bijective.
\end{lemma}

To state Tits' lemma, we briefly recall the notion of central series.

\begin{definition}[Central series]
	Let $ H $ be a group. A \defemph*{central series of $ H $}\index{central series} is a finite subgroup sequence
	\[ \compactSet{1_G} = Z_n \normsub Z_{n-1} \normsub \cdots \normsub Z_1 = H \]
	such that $ \commutator{H}{Z_i} \subs Z_{i+1} $ for all $ i \in \numint{1}{n-1} $. In particular, $ \listing{Z}{n} $ are normal subgroups of $ H $.
\end{definition}

\begin{lemma}
	Let $ H, K $ be groups and let $ \map{\pi}{H}{K}{}{} $ be a surjective homomorphism. If
	\[ \compactSet{1_H} = Z_n \normsub Z_{n-1} \normsub \cdots \normsub Z_1 = H \]
	is a central series in $ H $, then
	\[ \compactSet{1_K} = \pi(Z_n) \normsub \pi(Z_{n-1}) \normsub \cdots \normsub \pi(Z_1) = K \]
	is a central series in $ K $.
\end{lemma}
\begin{proof}
	For all $ i \in \numint{1}{n-1} $, we have
	\[ \commutator{K}{\pi(Z_i)} = \commutator{\pi(H)}{\pi(Z_i)} = \pi(\commutator{H}{Z_i}) \subs \pi(Z_{i+1}), \]
	as desired.
\end{proof}

\begin{lemma}[{\cite[4.7, Lemma~2, \onpage{559}]{Tits-KacMoodyUniquePres}}]\label{prodmap:tits}
	Let $ H $ be a group and let $ \listing{H}{n} $ be subgroups such that $ H $ is generated by $ H_1 \union \cdots \union H_n $. Assume that $ H $ possesses a central series
	\[ \compactSet{1} = Z_{h+1} \normsub Z_h \normsub \cdots \normsub Z_1 = H \]
	such that for all $ j \in \Set{1, \ldots, h} $, we have $ Z_j \subs \gen{H_{i(j)}, Z_{j+1}} $ for some $ i(j) \in \Set{1, \ldots, n} $. Then the following hold:
	\begin{lemenumerate}
		\item For every permutation $ \sigma $ of the set $ \Set{1, \ldots, n} $, the product map
		\[ \map{\mu_\sigma}{H_{\sigma(1)} \times \cdots \times H_{\sigma(n)}}{X}{}{} \]
		is surjective.
		
		\item If the product map $ \mu_{\sigma} $ is injective for some permutation $ \sigma $, then it is injective for all permutations $ \sigma $.
	\end{lemenumerate}
\end{lemma}

\subsection{Existence of a Bijective Ordering}

We now turn away from arbitrary groups and back to the group $ G $ with $ \roots $-commutator relations. We begin with our only criterion for $ G $-injectivity.

\begin{lemma}\label{prodmap:extremal-inj}
	Let $ \word{\alpha} = (\alpha_1, \ldots, \alpha_m) $ be an extremal ordering of $ \rootsub $. Assume that for any positive system $ \possys $ in $ \roots $ and any $ \beta \in -\possys $, we have $ \rootgr{\beta} \intersect \rootgr{\possys} = \compactSet{1_G} $. Then $ \word{\alpha} $ is $ G $-injective.
\end{lemma}
\begin{proof}
	We prove by induction on $ m $ that the product map
	\[ \map{\mu_{\word{\alpha}}}{\rootgr{\alpha_1} \times \cdots \times \rootgr{\alpha_m}}{\rootgr{\rootsub}}{}{} \]
	is injective. The case $ m=1 $ is trivial, so assume that $ m>1 $. Let $ \tup{g}{m} $ and $ \tup{h}{m} $ be two elements of $ \rootgr{\alpha_1} \times \cdots \times \rootgr{\alpha_m} $ such that $ g_1 \cdots g_m = h_1 \cdots h_m $. Then
	\begin{equation}\label{eq:prodmap:extremal-inj:1}
		h_1^{-1} g_1 = h_2 \cdots h_m (g_2 \cdots g_m)^{-1} \in \rootgr{\alpha_1} \intersect \gen{\rootgr{\alpha_j} \given j \in \numint{2}{m}}.
	\end{equation}
	Since $ \alpha_1 $ is extremal in $ \rootsub $, there exists a positive system $ \possys $ which contains $ \Set{\alpha_2, \ldots, \alpha_m} $ but not $ \alpha_1 $, so
	\begin{equation}\label{eq:prodmap:extremal-inj:2}
		\rootgr{\alpha_1} \intersect \gen{\rootgr{\alpha_j} \given j \in \numint{2}{m}} \subs \rootgr{\alpha_1} \intersect \rootgr{\possys} = \compactSet{1_G}.
	\end{equation}
	Putting~\eqref{eq:prodmap:extremal-inj:1} and~\eqref{eq:prodmap:extremal-inj:2} together, we infer that $ g_1 = h_1 $ and $ g_2 \cdots g_m = h_2 \cdots h_m $. By the induction hypothesis, it follows that $ g_j = h_j $ for all $ j \in \numint{1}{m} $. Thus $ \mu_{\word{\alpha}} $ is injective.
\end{proof}

As a next step, we make the connection between ideals of $ \roots $ and normal subgroups of $ G $. This is done in \cref{prodmap:ideal-reduced} for reduced root systems and in \cref{prodmap:ideal-nonreduced} for non-reduced crystallographic root systems. Both lemmas make essentially the same assertion, and they will be summarised in \cref{prodmap:ideal-normal}.

\begin{lemma}\label{prodmap:ideal-reduced}
	Assume that $ \roots $ is reduced. Let $ \rootsub $ be a subset of $ \roots $ with $ \rootsub \intersect (-\rootsub) = \emptyset $ and let $ I $ be a subset of $ \rootsub $. Assume that $ I $ is a (crystallographic) ideal of $ \rootsub $ (and that $ G $ has crystallographic $ \roots $-commutator relations). Then $ \rootgr{I} $ is normal in $ \rootgr{\rootsub} $.
\end{lemma}
\begin{proof}
	Let $ \alpha \in I $ and let $ \beta \in \rootsub $. Let $ x_\alpha \in \rootgr{\alpha} $, $ x_\beta \in \rootgr{\beta} $. We have to show that $ \commutator{x_\alpha}{x_\beta} $ lies in $ \rootgr{I} $. If $ \beta $ lies in $ I $, then this is clear, so we can assume that $ \beta $ lies in $ \rootsub \setminus I $. In particular, $ \alpha \ne \beta $. Since $ \rootsub \intersect (-\rootsub) = \emptyset $, we also have $ \alpha \ne -\beta $. As $ \roots $ is reduced, this implies that $ \alpha $ and $ \beta $ are non-proportional. If $ G $ has crystallographic $ \roots $-commutator relations, then $ \commutator{x_\alpha}{x_\beta} $ lies in $ \rootgr{\rootint{\alpha}{\beta}} $, so it lies in $ \rootgr{I} $ if $ I $ is a crystallographic ideal of $ \rootsub $. In general, $ \commutator{x_\alpha}{x_\beta} $ lies in $ \rootgr{\rootintcox{\alpha}{\beta}} $, so it lies in $ \rootgr{I} $ if $ I $ is an ideal. This finishes the proof.
\end{proof}

\begin{lemma}\label{prodmap:ideal-nonreduced}
	Assume that $ \roots $ is crystallographic and not reduced. Let $ \rootsub $ be a subset of $ \roots $ such that $ -\lambda \alpha \nin \rootsub $ for all $ \alpha \in \rootsub $ and all $ \lambda > 0 $. Let $ I $ be a subset of $ \rootsub $. Assume that $ I $ is a (crystallographic) ideal of $ \rootsub $ (and that $ G $ has crystallographic $ \roots $-commutator relations). Then $ \rootgr{I} $ is normal in $ \rootgr{\rootsub} $.
\end{lemma}
\begin{proof}
	Just like in the proof of \cref{prodmap:ideal-reduced}, we only have to show that $ \commutator{x_\alpha}{x_\beta} $ lies in $ \rootgr{I} $ where $ \alpha \in I $, $ \beta \in \rootsub \setminus I $, $ x_\alpha \in \rootgr{\alpha} $ and $ x_\beta \in \rootgr{\beta} $. If $ \alpha $ and $ \beta $ are non-proportional, we can proceed as in the proof of \cref{prodmap:ideal-reduced}. Otherwise we have $ \beta = \lambda \alpha $ for some $ \lambda \in \IR_{>0} \setminus \compactSet{1} $. Since $ \roots $ is crystallographic, we must have $ \lambda = 2 $ or $ \lambda = 1/2 $ by \cref{rootsys:cry-nonred-2}. If $ \lambda = 2 $, then $ x_\alpha $ and $ x_\beta $ both lie in $ \rootgr{\alpha} $ by \axiomref{rgg-axiom-div}, and thus $ \commutator{x_\alpha}{x_\beta} \in \rootgr{\alpha} \subs \rootgr{I} $. If $ \lambda = 1/2 $, then
	\[ \commutator{x_\alpha}{x_\beta} \in \commutator{\rootgr{2\beta}}{\rootgr{\beta}} \subs \commutator{\rootgr{\beta}}{\rootgr{\beta}} \subs \rootgr{2\beta} = \rootgr{\alpha} \subs \rootgr{I}, \]
	again by \axiomref{rgg-axiom-div}. This finishes the proof.
\end{proof}

\begin{proposition}\label{prodmap:ideal-normal}
	Assume that $ \roots $ is reduced or crystallographic. Let $ \rootsub $ be a subset of $ \roots $ such that $ -\lambda \alpha \nin \rootsub $ for all $ \alpha \in \rootsub $ and all $ \lambda > 0 $. Let $ I $ be a subset of $ \rootsub $. Assume that $ I $ is a (crystallographic) ideal of $ \rootsub $ (and that $ G $ has crystallographic $ \roots $-commutator relations). Then $ \rootgr{I} $ is normal in $ \rootgr{\rootsub} $.
\end{proposition}
\begin{proof}
	This follows from \cref{prodmap:ideal-reduced,prodmap:ideal-nonreduced}.
\end{proof}

As a consequence of \cref{prodmap:ideal-normal}, we obtain a criterion for $ G $-surjectivity.

\begin{lemma}\label{prodmap:subnormal-surj}
	 Assume that $ \roots $ is reduced or crystallographic. Let $ \rootsub $ be a subset of $ \roots $ such that $ -\lambda \alpha \nin \rootsub $ for all $ \alpha \in \rootsub $ and all $ \lambda > 0 $. Let $ \word{\alpha} = \tup{\alpha}{m} $ be an ordering of $ \rootsub $. Assume that $ \word{\alpha} $ is (crystallographically) subnormal (and that $ G $ has crystallographic $ \roots $-commutator relations). Then $ \word{\alpha} $ is $ G $-surjective.
\end{lemma}
\begin{proof}
	Put $ U_i \defl \rootgr{\alpha_i} $ for all $ i \in \numint{1}{m} $ and $ V_i \defl \gen{U_j \given j \in \numint{i}{m}} $ for all $ i \in \numint{1}{m} $. Since $ \word{\alpha} $ is (crystallographically) subnormal, it follows from \cref{prodmap:ideal-normal} that $ V_{i+1} $ is normal in $ V_i $ for all $ i \in \numint{1}{m-1} $. Hence by \cref{prodmap:abstract-subnormal-surjective}, the product map $ \map{}{\rootgr{\alpha_1} \times \cdots \times \rootgr{\alpha_m}}{G}{}{} $ is surjective, as desired.
\end{proof}

A consequence of \cref{prodmap:subnormal-surj} is that the product map on every closed root interval is surjective. In particular, if $ G $ is rank-2-injective (in the sense of \cref{basic:rank2inj-def}), then it is actually \enquote{rank-2-bijective}.

\begin{lemma}\label{prodmap:interval-surj}
	Let $ S = \indivset{\clrootintcox{\alpha}{\beta}} $ for non-proportional $ \alpha, \beta \in \indivset{\roots} $ and let $ \word{\alpha} = \tup{\alpha}{k} $ be an interval ordering of $ S $ in the sense of \cref{rootsys:clockwise-def}. Then $ \word{\alpha} $ is $ G $-surjective.
\end{lemma}
\begin{proof}
	We know from \cref{rootorder:interval-subnormal} that $ \word{\alpha} $ is a subnormal ordering of $ S $. Thus the assertion follows from \cref{prodmap:subnormal-surj}.
\end{proof}

\begin{remark}\label{prodmap:commpart-welldef}
	It follows from \cref{prodmap:interval-surj} that the maps in \cref{basic:commpart-def} are always well-defined if $ G $ is rank-2-injective.
\end{remark}

Putting everything together, we obtain a criterion for $ G $-bijectivity.

\begin{lemma}\label{prodmap:extremal-bij}
	Let $ \rootsub $ be a subset of $ \possys $. Assume that $ \rootsub $ is (crystallographically) closed (and that $ G $ has crystallographic $ \roots $-commutator relations). Assume further that for any positive system $ \possys' $ in $ \roots $ and any $ \gamma \in -\possys' $, we have $ \rootgr{\gamma} \intersect \rootgr{\possys'} = \compactSet{1_G} $. Then any extremal ordering of $ \rootsub $ is $ G $-bijective.
\end{lemma}
\begin{proof}
	Any such ordering is $ G $-injective by \cref{prodmap:extremal-inj} and $ G $-surjective by \cref{rootorder:extremal-subnormal,prodmap:subnormal-surj}.
\end{proof}

In particular, $ G $-bijective orderings exist under natural conditions.

\begin{proposition}\label{prodmap:bij-ex}
	Let $ \possys $ be a positive system in $ \roots $ and let $ \rootsub $ be a subset of $ \possys $ consisting of pairwise non-proportional roots. Assume that $ \rootsub $ is (crystallographically) closed (and that $ G $ has crystallographic $ \roots $-commutator relations). Assume further that for any positive subsystem $ \possys' $ of $ \roots $ and any $ \gamma \in -\possys $, we have $ \rootgr{\gamma} \intersect \rootgr{\possys'} = \compactSet{1_G} $. Then there exists a $ G $-bijective ordering of $ \rootsub $.
\end{proposition}
\begin{proof}
	By \cref{rootorder:extremal-exists}, there exists an extremal ordering of $ \rootsub $ and this ordering is $ G $-bijective by \cref{prodmap:extremal-bij}.
\end{proof}

\begin{note}
	The assumption in \cref{prodmap:bij-ex} that $ \rootsub $ is (crystallographically) closed is essential for the $ G $-surjectivity of $ \rootsub $. Even if $ \rootsub $ is not (crystallographically) closed, there exists an extremal ordering of $ \rootsub $ by \cref{rootorder:extremal-exists}, but we cannot apply \cref{prodmap:extremal-bij}. However, such an extremal ordering is still $ G $-injective by \cref{prodmap:extremal-inj}.
\end{note}

\subsection{All Orderings are Bijective for Crystallographic Root Systems}

In the crystallographic setting, we can apply \cref{prodmap:tits} to infer that all orderings are $ G $-bijective. Here we have to assume that certain root groups are abelian. For $ \roots $-graded groups of rank at least~3, this will be proven to automatically be true (see \cref{rgg:rootgr-abelian}), and the proof of this fact does not rely on \cref{prodmap:cry-all-bij}. The first two assertions of the following result are essentially a special case of \cite[3.12]{LoosNeherBook}.

\begin{proposition}\label{prodmap:cry-all-bij}
	Assume that $ \roots $ is crystallographic, that $ G $ has crystallographic $ \roots $-commutator relations and that for all roots $ \alpha $ for which $ 2\alpha $ is not a root, the root group $ \rootgr{\alpha} $ is abelian. Let $ \rootsub $ be a crystallographically closed subset of $ \roots $ such that $ -\lambda \alpha \nin \rootsub $ for all $ \alpha \in \rootsub $ and all $ \lambda > 0 $. Assume further that the roots in $ \rootsub $ are pairwise non-proportional. Then the following hold:
	\begin{proenumerate}
		\item \label{prodmap:cry-all-bij:sur}Any ordering of $ \rootsub $ is $ G $-surjective.
		
		\item \label{prodmap:cry-all-bij:ex-all}If some ordering of $ \rootsub $ is $ G $-injective, then every ordering of $ \rootsub $ is $ G $-bijective.
		
		\item \label{prodmap:cry-all-bij:all}If $ \rootgr{\gamma} \intersect \rootgr{\possys'} = \compactSet{1_G} $ for any positive system $ \possys' $ in $ \roots $ and any $ \gamma \in -\possys $, then every ordering of $ \rootsub $ is $ G $-bijective.
	\end{proenumerate}
\end{proposition}
\begin{proof}
	By \cref{rootorder:clos-in-pos}, there exists a positive system $ \possys $ in $ \roots $ which contains $ \rootsub $. Denote by $ \map{\rootht}{\roots}{\IZ}{}{} $ the height function with respect to $ \possys $, choose a height ordering $ \word{\alpha} = \tup{\alpha}{m} $ of $ \rootsub $ and define $ Z_i \defl \gen{\rootgr{\alpha_j} \given j \in \numint{i}{m}} $ for all $ i \in \numint{1}{m} $. Then by the generalised commutator relation (\cref{basic:gen-com-rel}), which we can apply by the assumption on certain root groups to be abelian, we have
	\begin{align*}
		\commutator{\rootgr{\rootsub}}{Z_i} &\subs \gen{\rootgr{\alpha} \given \alpha \in \roots, \rootht(\alpha) > \rootht(\alpha_i)} \subs Z_{i+1}
	\end{align*}
	for all $ i \in \numint{1}{m-1} $. Hence
	\[ \compactSet{1_G} \normsub Z_m \normsub Z_{m-1} \normsub \cdots \normsub Z_1 = \rootgr{\rootsub} \]
	is a central series of $ \rootgr{\rootsub} $. Further, it is clear that $ Z_i = \gen{\rootgr{\alpha_i}, Z_{i+1}} $ for all $ i \in \numint{1}{m-1} $. Thus the conditions of \cref{prodmap:tits} are satisfied. The assertions of~\itemref{prodmap:cry-all-bij:sur} and~\itemref{prodmap:cry-all-bij:ex-all} follow. Using \cref{prodmap:bij-ex}, assertion~\itemref{prodmap:cry-all-bij:all} follows from~\itemref{prodmap:cry-all-bij:ex-all}.
\end{proof}


\section{Root Graded Groups}

\label{sec:rgg-def}

\begin{secnotation}
	We denote by $ \roots $ an arbitrary root system.
\end{secnotation}

In \cref{prodmap:bij-ex,prodmap:cry-all-bij}, we have found a simple criterion for the existence of $ G $-bijective orderings in positive systems of $ \roots $. This criterion is the final axiom for root graded groups, which should be thought of as a \enquote{non-degeneracy condition}. In \cref{rgg:nondeg-comparison}, we will discuss several variations of this condition and how they relate to each other.

\begin{definition}[Root graded group]\label{rgg-def}
	Let $ G $ be a group. A \defemph*{$ \roots $-grading of $ G $}\index{root grading|see{root graded group}}\index{root graded group} is a $ \roots $-pregrading $ (\rootgr{\alpha})_{\alpha \in \roots} $ with the following properties:
	\begin{stenumerate}
		\item $ G $ is generated by $ (\rootgr{\alpha})_{\alpha \in \roots} $ and $ \rootgr{\alpha} \ne \compactSet{1_G} $ for all roots $ \alpha $.
		
		\item $ G $ has $ \roots $-commutator relations with root groups $ (\rootgr{\alpha})_{\alpha \in \roots} $ (in the sense of \cref{rgg:group-commrel-def}).
		
		\item For all $ \alpha \in \roots $, there exists an $ \alpha $-Weyl element in $ G $ (in the sense of \cref{weyl:weyl-def}).
		
		\item \label{rgg-def:nondeg}For all positive systems $ \possys $ and all $ \alpha \in \roots \setminus \possys $, we have $ \rootgr{\possys} \intersect \rootgr{\alpha} = \compactSet{1_G} $.
	\end{stenumerate}
	If, in addition, $ G $ has crystallographic $ \roots $-commutator relations with root groups $ (\rootgr{\alpha})_{\alpha \in \roots} $ (in the sense of \cref{rgg:group-commrel-cry-def}), then $ (\rootgr{\alpha})_{\alpha \in \roots} $ is called a \defemph*{crystallographic $ \roots $-grading of $ G $}. In this context, the pair $ (G, (\rootgr{\alpha})_{\alpha \in \roots}) $ is called a (crystallographic) \defemph*{$ \roots $-graded group}.\index{root graded group}
\end{definition}

\begin{notation}
	If $ (G, (\rootgr{\alpha})_{\alpha \in \roots}) $ is a $ \roots $-graded group, then the root system $ \roots $ is called the \defemph*{type of $ (G, (\rootgr{\alpha})_{\alpha \in \roots}) $}\index{type of a root graded group} or simply the \defemph*{type of $ G $}. We adopt the convention that if $ X $ is a property of $ \roots $, then $ G $ is also said to have property $ X $. For example, a root graded group is said to be irreducible, simply-laced or of rank~$ r $ if its type has the corresponding property.
\end{notation}

	Note that the definition of root graded groups is purely combinatorial: it only involves a root system and its Weyl group. There is no algebraic structure (like a ring or a module) hidden in the axioms, and there is at first glance no reason to believe that such a thing should exist. However, in \cref{chap:chev} we will see that Chevalley groups provide a large and well-known class of examples of crystallographic root graded groups which are \enquote{defined over a commutative associative ring}. The goal of this book is to show that every crystallographic root graded group of rank at least~3 is \enquote{of algebraic origin} in a similar way. The following remark makes this statement more precise.

\begin{goal}[Coordinatisation problem]\label{rgg:coord-gen-def}
	Let $ G $ be a group with a $ \roots $-grading $ (\rootgr{\alpha})_{\alpha \in \roots} $ and denote by $ \calO \defl \IndivOrb(\roots) = (O_1, \ldots, O_k) $ the set of indivisible orbits of $ \roots $ (see \cref{rootsys:orb-def}). Let $ \calX = (X_1, \ldots, X_k) $ be a family of algebraic structures (each of which has an underlying group structure) and put $ X_\alpha \defl X_{i(\alpha)} $ for all $ \alpha \in \indivset{\roots} $ where $ i(\alpha) $ is the unique index such that $ \alpha \in O_{i(\alpha)} $. A \defemph*{coordinatisation of $ G $ by $ \calX $}\index{coordinatisation of a root graded group} is a family of maps $ (\risom{\alpha})_{\alpha \in \roots} $ with the following properties:
	\begin{stenumerate}
		\item For all $ \alpha \in \indivset{\roots} $, the map $ \map{\risom{\alpha}}{X_\alpha}{\rootgr{\alpha}}{}{} $ is an isomorphism between the underlying group of the algebraic structure $ X_\alpha $ and the root group $ \rootgr{\alpha} $.
		
		\item Let $ (\alpha, \beta) $ be a pair of non-proportional indivisible roots and denote by $ \tup{\gamma}{p} $ an interval ordering of $ \rootintcox{\alpha}{\beta} $. Then there exist maps
		\[ (\map{f_i}{X_\alpha \times X_\beta}{X_{\gamma_i}}{}{})_{i \in \numint{1}{p}} \]
		such that the following generalisation of the Chevalley commutator formula (see \cref{chev:comm-formula}) holds:
		\[ \commutator{\risom{\alpha}(a)}{\risom{\beta}(b)} = \prod_{i=1}^p \risom{\gamma_i}\brackets[\big]{f_i(a,b)} \]
		for all $ a \in X_\alpha $ and $ b \in X_\beta $. Further, the maps $ f_1, \ldots, f_p $ can be described explicitly using only the structural maps of the algebraic structures in $ \calX $.
	\end{stenumerate}
	The preceding definition should be regarded more as a template than as a precise definition. For each type of $ \roots $ in the classification of root systems, we will introduce a more concise definition of coordinatisations of $ \roots $-graded groups: see \cref{ADE:param-def} for types $ A $, $ D $ and $ E $, \cref{B:standard-param-def} for type $ B $, and \cref{BC:standard-param-def} for type $ BC $ (which includes type $ C $ as a special case). For type $ F_4 $, technicalities prevent us from using a definition in the exact same spirit, see \cref{F4:ex:coord:stsigns,F4:ex:coord:stsigns:note}.
	
	For all crystallographic root systems $ \roots $ of rank at least $ 3 $, the goal of this book is to find a class $ \algclass{\roots} $ of algebraic structures such that every crystallographic $ \roots $-graded groups is coordinatised by an object $ \calX \in \algclass{\roots} $. This is precisely the content of \cref{ADE:thm,B:thm,BC:thm,F4:thm}. We refer to this problem as the \defemph{coordinatisation problem}.
\end{goal}

\begin{goal}[Existence problem]\label{rgg:existence-problem}
	A secondary goal of this work is a solution of the \defemph{existence problem}: Having found a suitable class $ \algclass{\roots} $ of algebraic structures as in \cref{rgg:coord-gen-def}, we want to show that $ \algclass{\roots} $ is chosen optimally. This means that for all algebraic structures $ \calX \in \algclass{\roots} $, we want to show that there exists a $ \roots $-graded group $ G_\calX $ which is coordinatised by $ \calX $. If $ \calX $ is, in a suitable sense, \enquote{associative}, we can always find a matrix group $ G_\calX $ with this property. In the case that $ \algclass{\roots} $ allows non-associative structures, more work is needed, which is beyond the scope of this book. In some cases, no complete solution of the existence problem is available.
\end{goal}

\begin{note}
	The axiom that $ G $ is generated by $ (\rootgr{\alpha})_{\alpha \in \roots} $ is harmless: If $ G $ is a group with a $ \roots $-pregrading $ (\rootgr{\alpha})_{\alpha \in \roots} $ which satisfies all axioms except for this one, then $ (\rootgr{\alpha})_{\alpha \in \roots} $ is a $ \roots $-grading of the group $ G' $ which is generated by $ (\rootgr{\alpha})_{\alpha \in \roots} $. Since our interest lies in the coordinatisation of the root groups, everything outside of $ G' $ is irrelevant for our purposes, so we simply require that $ G=G' $.
\end{note}

We will never need the notion of homomorphisms of root graded groups, but there is no harm in stating the obvious definition.

\begin{definition}[Homomorphisms of root graded groups]
	Let $ (G, (\rootgr{\alpha})_{\alpha \in \roots}) $ and $ (G', (\rootgr{\alpha}')_{\alpha \in \roots}) $ be $ \roots $-graded groups. A group homomorphism $ \map{f}{G}{G'}{}{} $ is called a \defemph*{homomorphism of root graded groups}\index{root graded group!homomorphism} if $ f(\rootgr{\alpha}) \subs \rootgr{\alpha}' $ for all roots $ \alpha $. A homomorphism of root graded groups is called an \defemph*{isomorphism of root graded groups}\index{root graded group!isomorphism} if it is an isomorphism of groups whose inverse is also a homomorphism of root graded groups.
\end{definition}

We now record some special cases of previous results from this chapter.

\begin{proposition}\label{rgg:prodmap-Nw}
	Let $ G $ be a group with a $ \roots $-grading $ (\rootgr{\alpha})_{\alpha \in \roots} $. Choose a root base $ \rootbase $ of $ \roots $ and denote by $ \possys $ the corresponding positive system. Let $ w $ be an element of the Weyl group of $ \roots $, let $ \word{\delta} = \tup{\delta}{k} $ be a reduced $ \rootbase $-expression of $ w $ and let $ \word{\alpha} = (\alpha_1, \ldots, \alpha_k) $ be the associated root sequence from \cref{rootsys:rootseq-def}. Then the product map on $ \word{\alpha} $ in $ G $ is bijective.
\end{proposition}
\begin{proof}
	Denote by $ \word{\alpha}' \defl (\alpha_k, \ldots, \alpha_1) $ the inverse root sequence. Then $ \switchset(w) = \Set{\listing{\alpha}{k}} $ by \cref{rootsys:switchset-description}, where $ \switchset $ is defined as in \cref{rootsys:switchset-def}. Further, $ \word{\alpha} $ is an extremal ordering by \cref{rootorder:switchset-extremal} and $ \switchset(w) $ is clearly closed, so the product map on $ \word{\alpha}' $ in $ G $ is bijective by \cref{prodmap:extremal-bij}. Using \cref{prodmap:switch-order}, we infer that the product map on $ \word{\alpha} $ is bijective as well.
\end{proof}

\begin{proposition}\label{rgg:prodmap-cry}
	Assume that $ \roots $ is crystallographic and let $ G $ be a group with a crystallographic $ \roots $-grading $ (\rootgr{\alpha})_{\alpha \in \roots} $. Assume further that for all roots $ \alpha $ for which $ 2\alpha $ is not a root, the root group $ \rootgr{\alpha} $ is abelian. Let $ \possys $ be a positive system in $ \roots $ and let $ \rootsub $ be a crystallographically closed subset of $ \possys $ whose elements are pairwise non-proportional. Then the product map on any ordering of $ \rootsub $ is bijective.
\end{proposition}
\begin{proof}
	This follows from \thmitemcref{prodmap:cry-all-bij}{prodmap:cry-all-bij:all}.
\end{proof}

By \cref{rgg:rootgr-abelian}, the assumption in \cref{rgg:prodmap-cry} on certain root groups to be abelian is always satisfied if $ \roots $ is of rank at least~3. However, we still have to prove this fact.

\begin{lemma}\label{rgg:rank2inj}
	Let $ G $ be a group with a $ \roots $-grading $ (\rootgr{\alpha})_{\alpha \in \roots} $. Then $ G $ is rank-2-injective.
\end{lemma}
\begin{proof}
	Let $ \possys' $ be a positive system in a parabolic rank-2 subsystem of $ \roots $ and let $ \word{\alpha} = \tup{\alpha}{k} $ be an interval ordering of $ \indivset{(\possys')} $. Then $ \rootsub \defl \Set{\listing{\alpha}{k}} $ is closed and $ \word{\alpha} $ is an extremal ordering of $ \rootsub $. Thus it follows from \cref{prodmap:extremal-bij} that the product map on $ \word{\alpha} $ is bijective, which finishes the proof.
\end{proof}

\begin{theorem}[Braid relations]\label{braid:all}
	Let $ G $ be a group with a $ \roots $-grading $ (\rootgr{\alpha})_{\alpha \in \roots} $ and assume that $ \roots $ is crystallographic or reduced. Then $ G $ satisfies the braid relations for Weyl elements.
\end{theorem}
\begin{proof}
	By \cref{rgg:rank2inj}, $ G $ is rank-2-injective. Hence the assertion is a special case of \cref{braid:all:weyl-ext}.
\end{proof}

The following result, the analogue of \cref{basic:commrel-subgroup} for root gradings, guarantees the existence of many root graded subgroups in a root graded group. It allows us to prove many properties of root graded groups by restricting to the cases of rank 2 or rank 3.

\begin{proposition}\label{rgg:subgroup}
	Let $ G $ be a group with a $ \roots $-grading $ (\rootgr{\alpha})_{\alpha \in \roots} $ and let $ \roots' $ be a subset of $ \roots $. Denote by $ H $ the group which is generated by $ (\rootgr{\alpha})_{\alpha \in \roots'} $. If $ \roots' $ is a closed root subsystem, then$ (\rootgr{\alpha})_{\alpha \in \roots'} $ is a $ \roots' $-grading of $ H $. If $ \roots' $ is a crystallographically closed root subsystem and the commutator relations of $ G $ are crystallographic, then $ (\rootgr{\alpha})_{\alpha \in \roots'} $ is a crystallographic $ \roots' $-grading of $ H $.
	For any root $ \alpha $ in $ \roots' $, every $ \alpha $-Weyl element in $ G $ is also an $ \alpha $-Weyl element in $ H $.
\end{proposition}
\begin{proof}
	The assertion about Weyl elements is clear. The subgroup $ H $ has the desired (crystallographic) $ \roots' $-commutator relations by \cref{basic:commrel-subgroup}. Further, Axiom~\thmitemref{rgg-def}{rgg-def:nondeg} is satisfied by \cref{rootsys:subsys:clos-pos}.
\end{proof}

The existence of Weyl elements implies that many pairs of root groups are isomorphic. However, it does not provide canonical isomorphisms between these root groups. The purpose of the parametrisation theorem, which we will study in \cref{chap:param}, is precisely to rectify this problem.

\begin{lemma}\label{rgg:orbit-isom}
	Let $ G $ be a group with a $ \roots $-grading $ (\rootgr{\alpha})_{\alpha \in \roots} $ and let $ \alpha $, $ \beta $ be roots which lies in the some orbit under the Weyl group. Then $ \rootgr{\alpha} $ is isomorphic to $ \rootgr{\beta} $.
\end{lemma}
\begin{proof}
	Choose an element $ u $ of the Weyl group such that $ \alpha^u = \beta $. By the definition of the Weyl group, there exist roots $ \listing{\gamma}{k} $ such that $ u = \reflbr{\gamma_1} \cdots \reflbr{\gamma_k} $. For each $ i \in \numint{1}{k} $, we choose a $ \gamma_i $-Weyl element $ w_i $, and we put $ w \defl w_1 \cdots w_k $. Then the map $ \map{}{}{}{x}{x^w} $ is an isomorphism from $ \rootgr{\alpha} $ to $ \rootgr{\beta} $.
\end{proof}

We end this section with a discussion of Axiom~\thmitemref{rgg-def}{rgg-def:nondeg}.

\begin{note}[Possible choices of the non-degeneracy condition]\label{rgg:nondeg-comparison}
	Let $ G $ be a group which has $ \roots $-commutator relations with root groups $ (\rootgr{\alpha})_{\alpha \in \roots} $. We consider the following conditions on $ G $:
	\begin{stenumerate}
		\item For all positive systems $ \possys $ in $ \roots $, we have $ \rootgr{\possys} \intersect \rootgr{-\possys} = \compactSet{1_G} $.
		
		\item Axiom~\thmitemref{rgg-def}{rgg-def:nondeg}: For all positive systems $ \possys $ and all $ \alpha \in \roots \setminus \possys $, we have $ \rootgr{\possys} \intersect \rootgr{\alpha} = \compactSet{1_G} $.
		
		\item For any positive system $ \possys $, every ordering of $ \possys $ is $ G $-bijective.
		
		\item For any positive system $ \possys $, there exists a $ G $-bijective ordering of $ \possys $.
		
		\item $ G $ is rank-2-injective.

		\item For all non-proportional roots $ \alpha, \beta $, we have $ \rootgr{\alpha} \intersect \rootgr{\beta} = \compactSet{1_G} $.
	\end{stenumerate}
	Then we have implications \enquote{(i) $ \implies $ (ii) $ \implies $ (iv)} and \enquote{(iii) $ \implies $ (v) $ \implies $ (vi)} where the non-trivial implications hold by \cref{prodmap:bij-ex,basic:prodmap-triv-intersect}. Further, if $ \roots $ is crystallographic and $ G $ has crystallographic $ \roots $-commutator relations, then we also have \enquote{(iii) $ \Longleftrightarrow $ (iv)} by \thmitemcref{prodmap:cry-all-bij}{prodmap:cry-all-bij:ex-all}. Thus under crystallographic assumptions, we have a clean chain
	\[ \text{(i)} \implies \text{(ii)} \implies \text{(iii)} \Longleftrightarrow \text{(iv)} \implies \text{(v)} \implies \text{(vi)}. \]
	
	Most of our proofs work under the relatively weak assumption (v) on rank-2-injectivity. However, the blueprint technique itself requires the stronger condition (ii), which is why we have chosen it as an axiom for root graded groups. It is interesting to note that the coordinatisation of simply-laced root graded groups in \cref{chap:simply-laced} (which does not rely on the blueprint technique) even works under the weaker condition (vi). This suggests that it might be possible to prove the coordinatisation results on non-simply-laced root systems under slightly weaker assumptions than (ii) if one found a way to work without the blueprint technique. However, even Zhang's computations on $ C_n $-graded groups in \cite{Zhang} require condition~(iii), so it seems that this is the minimal condition which has to be used in any case.
	
	As a final remark, it is worth mentioning that all specific examples of root graded groups that we construct in this book (including all Chevalley groups) even satisfy the strongest condition~(i). In these examples, $ \rootgr{\possys} $ can be regarded as a group of upper triangular (generalised) matrices while $ \rootgr{-\possys} $ consists of lower triangular matrices.
\end{note}


\section{Foldings of Root Graded Groups}

\label{sec:RGG-fold}

\begin{secnotation}
	We use the same notation as in~\ref{secnot:rootsys-fold} for $ \roots $, $ V $, $ \rootbase $, $ \rho $, $ \tau $, $ \fixspace \defl \Set{v \in V \given \tau(v) = v} $, $ \map{\pi}{V}{\fixspace}{}{} $ and $ \roots' \defl \pi(\roots) $. In addition, we denote by $ G $ a group with a $ \roots $-pregrading $ (\rootgr{\alpha})_{\alpha \in \roots} $.
\end{secnotation}

In this section, we describe how the folding $ \map{\pi}{\roots}{\roots'}{}{} $ gives rise to a construction of $ \roots' $-graded groups from $ \roots $-graded groups. This construction is closely related to the construction of the twisted groups from Chevalley groups (see \cite[13.4]{Carter-Chev}), but more general. We will use it in \cref{sec:F4:const} to construct examples of $ F_4 $-graded groups as foldings of $ E_6 $-graded groups.

We use a similar method in \cite{RGG-H3} to construct $ H_3 $-graded groups as foldings of $ D_6 $-graded groups. However, this requires a more general approach because $ H_3 $ cannot be realised as a folding of $ D_6 $ which is induced by a diagram automorphism, but only by a folding in the more general sense of \cite{Muehlherr_Foldings}.

By construction, roots in $ \roots' $ can be identified with their preimages under $ \pi $ in $ \roots $, that is, with sets of roots in $ \roots $. Thus the following definition of a $ \roots' $-pregrading $ (\varrootgr{\alpha'})_{\alpha' \in \roots'} $ is not surprising.

\begin{definition}[Folded pregrading]
	For all $ \alpha' \in \roots' $, we put
	\[ \varrootgr{\alpha'} \defl \gen{\rootgr{\alpha} \given \alpha \in \roots, \pi(\alpha) = \alpha'}. \]
	The $ \roots' $-pregrading $ (\varrootgr{\alpha'})_{\alpha' \in \roots'} $ is called a \defemph*{folding of $ (\rootgr{\alpha})_{\alpha \in \roots} $}.\index{folding!of root graded groups}
\end{definition}

The goal of this section is to find conditions on $ (G, (\rootgr{\alpha})_{\alpha \in \roots}) $ which guarantee that $ (G, (\varrootgr{\alpha'})_{\alpha' \in \roots}) $ is a $ \roots' $-graded group. (Note that group $ G $ remains the same in both gradings, only the family of root subgroups changes.) Clearly, we should require that $ (G, (\rootgr{\alpha})_{\alpha \in \roots}) $ itself is a $ \roots $-graded group, but it turns out that we need a bit more.

\begin{lemma}[Commutator relations]\label{rgg-fold:comm}
	Assume that $ \roots' $ is reduced and that for all roots $ \zeta $ for which $ 2\zeta $ is not a root, the root group $ \rootgr{\zeta} $ is abelian. If $ G $ has (crystallographic) $ \roots $-commutator relations with root groups $ (\rootgr{\alpha})_{\alpha \in \roots} $, then it also has (crystallographic) $ \roots' $-commutator relations with root groups $ (\varrootgr{\alpha'})_{\alpha' \in \roots'} $.
\end{lemma}
\begin{proof}
	This follows from \cref{basic:comm-fold}, whose assumptions are satisfied by \thmitemcref{fold:cor}{fold:cor:preim-pts}.
\end{proof}

\begin{lemma}[Non-degeneracy]\label{rgg-fold:nondeg}
	If $ \rootgr{\possys} \intersect \rootgr{-\possys} = \compactSet{1_G} $ for all positive systems $ \possys $ in $ \roots $, then also $ \varrootgr{\possys'} \intersect \varrootgr{-\possys'} = \compactSet{1_G} $ for all positive systems $ \possys' $ in $ \roots' $.
\end{lemma}
\begin{proof}
	Let $ \possys' $ be a positive system in $ \roots' $. By \thmitemcref{fold:cor}{fold:cor:ind-surj}, there exists a positive system $ \possys $ in $ \roots $ such that $ \pi(\possys) = \possys' $. Then for all $ \beta \in \roots $, we have $ \beta \in \possys $ if and only if $ \pi(\beta) \in \possys' $ because $ \roots $ is the disjoint union of $ \possys $ and $ -\possys $. It follows that
	\begin{align*}
		\varrootgr{\possys'} \intersect \varrootgr{-\possys'} \subs \rootgr{\possys} \intersect \rootgr{-\possys} = \compactSet{1_G},
	\end{align*}
	as desired.
\end{proof}

\begin{note}
	We cannot show that the validity of the more general Axiom~\thmitemref{rgg-def}{rgg-def:nondeg} for $ (\rootgr{\alpha})_{\alpha \in \roots} $ (that is, $ \rootgr{\possys} \intersect \rootgr{\alpha} = \compactSet{1_G} $ for all $ \alpha \in \roots \setminus \possys $) implies the validity of Axiom~\thmitemref{rgg-def}{rgg-def:nondeg} for $ (\varrootgr{\alpha'})_{\alpha' \in \roots'} $. The reason for this is that root groups in $ (\varrootgr{\alpha'})_{\alpha' \in \roots'} $ are not necessarily contained in root groups of $ (\rootgr{\alpha})_{\alpha \in \roots} $.
\end{note}

\begin{lemma}[Weyl elements]\label{rgg-fold:weyl}
	Let $ J $ be an orbit of $ \rootbase $ under $ \rho $ and let $ u_\alpha $ be an $ \alpha $-Weyl element (with respect to $ (\rootgr{\gamma})_{\gamma \in \roots} $) for all $ \alpha \in J $. Choose $ k \in \Npos $ and $ \listing{\alpha}{k} \in J $ such that $ \word{\alpha} = \tup{\alpha}{k} $ is an expression of the longest element $ w_0 \defl w_0^J $ of the subgroup $ W_J \defl \gen{\refl{\alpha} \given \alpha \in J} $ of $ \Weyl(\roots) $. Put $ u \defl \prod_{i=1}^k u_{\alpha_i} $ and denote by $ \alpha' $ the unique element in $ \pi(J) $. Then $ \varrootgr{\beta'}^u = \varrootgr{(\beta')^{\reflbr{\alpha'}}} $ for all $ \beta' \in \roots' $. If, in addition, the roots in $ J $ are pairwise orthogonal, then $ u $ is an $ \alpha' $-Weyl element (with respect to $ (\varrootgr{\gamma'})_{\gamma' \in \roots'} $).
\end{lemma}
\begin{proof}
	Let $ \beta \in \roots $ and put $ \beta' \defl \pi(\beta) $. As a first step, we prove that $ \rootgr{\beta}^u \subs \varrootgr{(\beta')^{\reflbr{\alpha'}}} $. Since $ (u_{\gamma})_{\gamma \in J} $ are Weyl elements, $ \rootgr{\beta}^u $ lies in the root group associated to
	\[ \beta^{\reflbr{\alpha_1} \cdots \reflbr{\alpha_k}} = \beta^{w_0} \in \roots. \]
	Now it follows from \cref{fold:shortest-act} that
	\begin{align*}
		\pi\brackets[\big]{\beta^{w_0}} = \pi(\beta)^{\reflbr{\alpha'}} = (\beta')^{\reflbr{\alpha'}}.
	\end{align*}
	We conclude that $ \rootgr{\beta}^u $ is contained in $ \varrootgr{(\beta')^{\reflbr{\alpha'}}} $.
	
	It follows from the conclusion of the previous paragraph that $ \varrootgr{\beta'}^u \subs \varrootgr{(\beta')^{\reflbr{\alpha'}}} $ for all $ \beta' \in \roots' $. Note that $ (\alpha_k, \ldots, \alpha_1) $ is an expression of $ (w_0)^{-1} = w_0 $, and that $ u_\alpha^{-1} $ is also an $ \alpha $-Weyl element for all $ \alpha \in J $. Hence we can also apply the conclusion of the previous paragraph to $ \bar{u} \defl u_{\alpha_k}^{-1} \cdots u_{\alpha_1}^{-1} = u^{-1} $, which yields that $ \varrootgr{\beta'}^{\bar{u}} \subs \varrootgr{(\beta')^{\reflbr{\alpha'}}} $ for all $ \beta' \in \roots' $. Replacing $ \beta' $ by $ (\beta')^{\reflbr{\alpha'}} $ and using that $ \bar{u} = u^{-1} $, we infer that $ \varrootgr{(\beta')^{\reflbr{\alpha'}}} \subs \varrootgr{\beta'}^u $ for all $ \beta' \in \roots' $. We conclude that $ \varrootgr{\beta'}^u = \varrootgr{(\beta')^{\reflbr{\alpha'}}} $ for all $ \beta' \in \roots' $, which proves the first part of the claim.
	
	To prove that $ u $ is an $ \alpha' $-Weyl element, it remains to show that it is contained in $ \varrootgr{-\alpha'} \varrootgr{\alpha'} \varrootgr{-\alpha'} $. By construction, it is contained in $ \prod_{i=1}^k \rootgr{-\alpha_i} \rootgr{\alpha_i} \rootgr{-\alpha_i} $. The assumption on $ J $ implies that it is a root base of type $ A_1^k $. In particular, $ \listing{\alpha}{k} $ must be pairwise distinct and $ J = \Set{\listing{\alpha}{k}} $. We conclude that for all distinct $ i,j \in \numint{1}{k} $, the groups $ \rootgr{\alpha_i} $ and $ \gen{\rootgr{\alpha_j}, \rootgr{-\alpha_j}} $ commute. Hence the factors of $ u $ can be reordered in a way which yields
	\[ u \in \prod_{i=1}^k \rootgr{-\alpha_i} \prod_{i=1}^k \rootgr{\alpha_i} \prod_{i=1}^k \rootgr{-\alpha_i} \subs \varrootgr{-\alpha'} \varrootgr{\alpha'} \varrootgr{-\alpha'}. \]
	The assertion follows.
\end{proof}

We can summarise the previous results as follows.

\begin{proposition}\label{rgg-fold:thm}
	Let $ (G, (\rootgr{\alpha})_{\alpha \in \roots}) $ be a (crystallographic) $ \roots $-graded group with $ \rootgr{\possys} \intersect \rootgr{-\possys} = \compactSet{1_G} $ for all positive systems $ \possys $ in $ \roots $. Assume that $ \roots' $ is reduced and that for all roots $ \zeta $ for which $ 2\zeta $ is not a root, the root group $ \rootgr{\zeta} $ is abelian. Assume further that for each orbit $ J $ of $ \rootbase $ under $ \rho $, all roots in $ J $ are pairwise orthogonal. Then $ (G, (\varrootgr{\alpha'})_{\alpha' \in \roots'}) $ is a (crystallographic) $ \roots' $-graded group with $ \varrootgr{\possys'} \intersect \varrootgr{-\possys'} = \compactSet{1_G} $ for all positive systems $ \possys' $ in $ \roots' $.
\end{proposition}
\begin{proof}
	This is the culmination of \cref{rgg-fold:comm,rgg-fold:nondeg,rgg-fold:weyl}. We also use \cref{basic:weyl-ex-basis} to ensure that Weyl elements exist for all roots and not merely for the simple roots.
\end{proof}

Recall from \cref{rgg:rootgr-abelian} that the assumption on certain root groups to be abelian is automatically satisfied if $ \roots $ is irreducible of rank at least~3.


\section{Root Graded Groups in the Literature}

\label{sec:literature}

In this section, we briefly review several different notions of root gradings which have been discussed in the literature. All of them are special cases of our definition. We have already given an overview of this topic in the preface of this book, so we restrict ourselves to filling in the technical details that are missing. Some references to later parts of this book will be necessary at certain points.

We begin with Tits' notion of RGD-systems.

\begin{definition}[RGD-system, {\cite[\onpage{258}]{Tits-TwinBuildingsKacMoody}}]\label{rgd-def}
	Let $ \roots $ be a root system. An \defemph*{RGD-system of type $ \roots $}\index{RGD-system} is a $ \roots $-graded group $ (G, (\rootgr{\alpha})_{\alpha \in \roots}) $ such that $ \invset{\alpha} = \rootgr{\alpha} \setminus \compactSet{1_G} $ for all roots~$ \alpha $.
\end{definition}

As discussed in the preface, predecessors of the RGD-axioms can be found as early as in \cite[\onpage{140}]{Tits-GrpSemiSimpIsotrop}, and essentially the same definition is given in \cite{Faulkner-RGD}.

\begin{note}
	\Cref{rgd-def} does not agree precisely with the one given in \cite[\onpage{258}]{Tits-TwinBuildingsKacMoody}. First of all, the root system $ \roots $ is not assumed to be finite in \cite{Tits-TwinBuildingsKacMoody}, which makes the formulation of the commutator relation axiom slightly more technical. Secondly, an RGD-system in the sense of Tits is not assumed to be generated by $ (\rootgr{\alpha})_{\alpha \in \roots} $, but by $ (\rootgr{\alpha}) \union H $ where $ H $ denotes the intersection of the normalisers of all root groups. Since our interest lies in coordinatising the root groups of a root graded group, we are free to replace an RGD-system in Tits' sense by the subgroup generated by all root groups. Thirdly, Axiom~\thmitemref{rgg-def}{rgg-def:nondeg} is replaced by the slightly weaker assumption that for any root base $ \rootbase $ and for all $ \delta \in \rootbase $, the root group $ \rootgr{-\delta} $ is not contained in $ \rootgr{\possys} $ where $ \possys $ denotes the positive system corresponding to $ \rootbase $. Fourthly, the definition in \cite{Tits-TwinBuildingsKacMoody} only requires that $ \rootgr{\alpha} \setminus \compactSet{1_G} $ is a subset of $ \invset{\alpha} $, but by \cref{weyl:1-not-invertible}, this is equivalent to our requirement that $ \invset{\alpha} = \rootgr{\alpha} \setminus \compactSet{1_G} $ (except possibly for root systems with irreducible components of type $ A_1 $).
\end{note}

\begin{note}[Invertible elements in root gradings]
	Let $ G $ be a root graded group which is coordinatised by some algebraic structure $ \calA $. We will show in \cref{ADE:stsign:weyl-char,B:stsign:weyl-char,BC:stsign:weyl-char,F4:stsign:weyl-char} that for any root $ \alpha $, the elements of $ \invset{\alpha} $ are in bijective correspondence with the elements of $ \calA $ which are invertible (in a suitable algebraic sense). Thus RGD-systems are precisely the root graded groups which are coordinatised by \enquote{division structures}.
\end{note}

We now turn to notions of root gradings which do not have \enquote{division assumptions}.

\begin{note}[Faulkner's $ A_2 $-graded groups]
	In \cite[13.3]{Faulkner-NonAssocProj}, Faulkner defines \defemph*{groups of Steinberg type}\index{group of Steinberg type}. These are the same as $ A_2 $-graded groups, except for a few minor differences. Firstly, the commutator axiom~\axiomref{rgg-axiom-comm} in \cref{rgg:group-commrel-def} is required to hold for $ \alpha = \beta $ as well, which says precisely that all root groups are abelian. (That is, Faulkner uses the same definition of commutator relations as Loos-Neher, see \cref{rgg:commrel:loos-neher-note}.) We will show in \cref{ADE:abelian} that this axiom is superfluous. Secondly, $ \alpha $-Weyl elements in our sense are $ (-\alpha) $-Weyl elements in Faulkner's sense, but this is practically irrelevant by \thmitemcref{basic:weyl-general}{basic:weyl-general:minus}. Thirdly and finally, Axiom~\thmitemref{rgg-def}{rgg-def:nondeg} is not part of the definition of a group of Steinberg type. Instead, the weaker axiom that $ \rootgr{\alpha} \intersect \rootgr{\beta} = \compactSet{1_G} $ for all distinct roots $ \alpha, \beta $ is an additional assumption in Faulkner's coordinatisation theorem. We will see that our proof of the coordinatisation theorem for $ A_2 $-gradings (\cref{ADE:thm}) works under this more general assumption as well. We conclude that Faulkner's definition is essentially the same as our definition of $ A_2 $-graded groups.
\end{note}

We now turn to Shi's definition of root graded groups. It is formulated in \cite{Shi1993} only for the simply-laced root systems, but the generalisation to arbitrary (finite) reduced root systems is immediate. Essentially the same definition is used in \cite[3.1.5]{Zhang}. Before we can phrase it, we have to introduce Steinberg groups. Their definition relies on some structure constants in Chevalley groups which will be properly introduced in \cref{chap:chev}.

\begin{definition}[Steinberg group]\label{rgg-lit:steinberg-def}
	Let $ \comring $ be a commutative associative ring, let $ \roots $ be a reduced crystallographic root system and let $ \chevstr = (\chevstr_{\alpha, \beta})_{\alpha, \beta \in \roots} $ be a family of Chevalley structure constants of type $ \roots $ (in the sense of \cref{chev:struc-family-def}). The \defemph*{Steinberg group of type $ \roots $ over $ \comring $ with respect to $ \chevstr $}\index{Steinberg group} is the abstract group defined by the following presentation: The generators are symbols $ \steinhom{\alpha}(r) $ for $ \alpha \in \roots $ and $ r \in \comring $ and the following relations hold.
	\begin{stenumerate}
		\item $ \steinhom{\alpha}(r+s) = \steinhom{\alpha}(r) \steinhom{\alpha}(s) $ for all $ \alpha \in \roots $ and all $ r,s \in \comring $.
		
		\item Let $ \alpha $, $ \beta $ be non-proportional roots and define
		\[ \chevorder(\alpha, \beta) \defl \Set{(i,j) \in \Npos^2 \given i\alpha + j\beta \in \roots} \]
		as in \cref{chevorder:def}. Then for all $ r,s \in \comring $, we have a relation
		\[ \commutator{\steinhom{\alpha}(r)}{\steinhom{\beta}(s)} = \prod_{(i,j) \in \chevorder} \steinhom{i\alpha + j\beta}\brackets[\big]{\chevstr_{\alpha, \beta, i, j} r^i s^j} \]
		which has exactly the same form a the Chevalley commutator formula (\cref{chev:comm-formula}). Here the integers $ \chevstr_{\alpha, \beta, i, j} $ can be computed from $ \chevstr $ as in \cref{chev:struc-const-formulas}.
	\end{stenumerate}
	We denote this group by $ \St_{\roots, \chevstr}(\comring) $ or simply by $ \St_{\roots}(\comring) $. For each root $ \alpha $, the subgroup $ \steinrootgr{\alpha} \defl \Set{\steinhom{\alpha}(r) \given r \in \comring} $ is called the \defemph*{root group of $ \St_{\roots, \chevstr}(\comring) $ corresponding to~$ \alpha $}. For any root $ \alpha $ and any invertible element $ r $ in $ \comring $, we define
	\[ \steinweyl{\alpha}(r) \defl \steinhom{-\alpha}(-r^{-1}) \steinhom{\alpha}(r) \steinhom{-\alpha}(-r^{-1}). \]
\end{definition}

\begin{definition}[Shi's notion of root gradings, {\cite[(2.1)]{Shi1993}}]\label{shi-def}
	Let $ \roots $ be a reduced crystallographic root system and let $ G $ be a group. A \defemph*{$ \roots $-grading of $ G $ in Shi's sense} is a $ \roots $-pregrading $ (\rootgr{\alpha})_{\alpha \in \roots} $ for which there exist a family $ \chevstr = (\chevstr_{\alpha, \beta})_{\alpha, \beta \in \roots} $ of Chevalley structure constants, a commutative associative non-zero ring $ \comring $ and a homomorphism $ \map{\steinproj}{\St_\roots(\comring)}{G}{}{} $ with the following properties:
	\begin{stenumerate}
		\item There exists a positive system $ \possys $ in $ \roots $ such that the restriction of $ \steinproj $ to $ \gen{\steinrootgr{\alpha} \given \alpha \in \possys} $ is injective.
		
		\item $ \steinproj(\steinrootgr{\alpha}) $ is contained in $ \rootgr{\alpha} $ for all roots $ \alpha $.
		
		\item $ G $ is generated by $ (\rootgr{\alpha})_{\alpha \in \rootbase} $.
		
		\item $ G $ has crystallographic $ \roots $-commutator relations, and all root groups are abelian.
		
		\item For all distinct roots $ \alpha $, $ \beta $, we have $ \rootgr{\alpha} \intersect \rootgr{\beta} = \compactSet{1_G} $.
		
		\item \label{shi-def:weyl}For all roots $ \alpha $, the element $ \steinproj(\steinweyl{\alpha}(1_\comring)) $ is an $ \alpha $-Weyl element in $ G $.
	\end{stenumerate}
\end{definition}

Observe that all root groups in \cref{shi-def} are not equal to $ \compactSet{1} $ because $ \comring $ is non-zero.

\begin{note}
	Let $ \roots $ be a reduced crystallographic root system, let $ \alpha $ be an arbitrary root and put $ w_\alpha \defl \steinproj(\steinweyl{\alpha}(1_\IZ)) $. In \cite{Shi1993}, it is not assumed that $ w_\alpha $ is an $ \alpha $-Weyl element, but only that for all roots $ \beta $ such that $ (\alpha, \beta) $ is a root base of the subsystem that it spans, we have $ \rootgr{\beta}^{w_\alpha} = \rootgr{\refl{\alpha}(\beta)} $. It is then shown in \cite[(2.4)~(ii)]{Shi1993} that the same equation holds for all roots $ \beta $, which implies that $ w_\alpha $ is indeed an $ \alpha $-Weyl element. Thus \cref{shi-def} is equivalent to Shi's definition in \cite[(2.1)]{Shi1993} for simply-laced root systems. However, it is not at all clear that \cite[(2.1)]{Shi1993} also holds for non-simply-laced root systems. Hence \cref{shi-def} seems to be the correct way to formulate Shi's definition for arbitrary (reduced crystallographic) root systems.
\end{note}

\begin{remark}[Comparison]
	Let $ \roots $ be a reduced crystallographic root system. It is immediate that every $ \roots $-graded group in Shi's sense is a crystallographic $ \roots $-graded group in our sense. Conversely, assume that $ \roots $ is simply-laced of rank at least~2 and that $ (G, (\rootgr{\alpha})_{\alpha \in \roots}) $ is a $ \roots $-graded group in our sense. Our coordinatisation result for such groups (\cref{ADE:thm}) says precisely that there exists a ring $ \ring $ such that $ G $ satisfies the same relations as in \cref{rgg-lit:steinberg-def}. In other words, $ G $ is a quotient of \enquote{$ \St_{\roots, \chevstr}(\ring) $} for an appropriate $ \chevstr $, though $ \ring $ need not be commutative associative. Denote by $ \comring $ the commutative associative subring of $ \ring $ which is generated by $ 1_\ring $, so that $ \comring $ is isomorphic to $ \IZ $ or to $ \IZ / n\IZ $ for some $ n \in \IN_{\ge 2} $. Then there exists a homomorphism $ \map{\steinproj}{\St_{\roots, \chevstr}(\comring)}{G}{}{} $ with the properties in \cref{shi-def}. Here the injectivity of $ \steinproj $ on $ \gen{\steinrootgr{\alpha} \given \alpha \in \possys} $ follows from the fact that the commutator relations on $ \gen{\steinrootgr{\alpha} \given \alpha \in \possys} $ (and thus on $ \gen{\rootgr{\alpha} \given \alpha \in \possys} $) already determine the group multiplication on this subset. This can be proven as in \cite[(8.13)]{MoufangPolygons}, and it is also a consequence of the arguments in \cref{sec:prod-bij}.
	
	We conclude that Shi's definition of root gradings is in fact equivalent to ours for simply-laced root systems of rank at least~2. For non-simply-laced root systems, the situation is more delicate because the involved coordinatising structures are no longer mere rings. For example, assume that $ (G, (\rootgr{\alpha})_{\alpha \in B_3}) $ is a crystallographic $ B_3 $-graded group. \Cref{B:thm} yields a pointable quadratic module $ (\module, q) $ over a commutative associative ring $ \comring $ which coordinatises $ G $ via root isomorphisms $ (\risom{\alpha})_{\alpha \in B_3} $. Choose $ v_0 \in \module $ with $ q(v_0) = 1_\comring $. Then we have a homomorphism $ \map{\steinproj}{\St_{\roots, \chevstr}(\comring)}{G}{}{} $ (for an appropriate family $ \chevstr $) whose image is
	\[ G' \defl \gen{\Set{\risom{\alpha}(r) \given \alpha \text{ long}, r \in \comring} \union \Set{\risom{\beta}(r v_0) \given \beta \text{ short}, r \in \comring}}.  \]
	However, it is not clear that the restriction of $ \steinproj $ to $ \gen{\steinrootgr{\alpha} \given \alpha \in \possys} $ is injective for some positive system $ \possys $. For example, if $ \comring = \IZ / n\IZ $, it is not obvious that $ mv_0 \ne 0 $ for every proper divisor $ m $ of $ n $. This illustrates that \cref{shi-def} is not adequate to capture the more general phenomena which may occur in non-simply-laced root graded groups.
\end{remark}

While the involvement of the Steinberg group makes Shi's definition more complicated, it has practical benefits: It simplifies the (sign problem in the) parametrisation of these groups. See \cref{param:motiv:simply-laced-shi} for a few more details.

\begin{miscthm}[Coordinatisation results in the literature]\label{rgg-lit:coord}
	Let $ \roots $ be a root system. Recall from the preface that RGD-systems of type $ \roots $ are essentially the same thing as Moufang buildings of type $ \roots $ (see \cite[Section~7.8]{AbramenkoBrown-Buildings} for details). Thus Tits' classification of irreducible Moufang buildings for root systems of rank at least~3 in \cite{Tits-LectureNotes74} yields coordinatisation results for irreducible RGD-systems of rank at least~3. It should be noted, however, that the language of \cite{Tits-LectureNotes74} is rather different from the one in this book. In contrast, the classification of irreducible Moufang buildings of rank~2 (that is, of Moufang polygons) in \cite{MoufangPolygons} uses a group-theoretic language which is very similar to ours and which provides explicit commutator formulas. In \cite[(40.22)]{MoufangPolygons}, the classification of irreducible Moufang buildings of rank at least~3 is deduced from the classification of Moufang polygons, and explicit commutator relations are provided for some (but not all) pairs of roots. By Tits' famous Theorem 4.1.2 (see \cite[4.1.2]{Tits-LectureNotes74}), the commutator relations on all pairs of (non-proportional) roots are uniquely determined by the commutator relations in \cite[(40.22)]{MoufangPolygons}.
	
	We now turn to $ \roots $-gradings which are not assumed to be RGD-systems. Our coordinatisation result for simply-laced root gradings of rank at least~2 (\cref{ADE:thm}) is given in essentially the same way in \cite[(2.3)]{Shi1993} for root gradings in Shi's sense. For root gradings of types $ B $, $ BC $ and $ F_4 $, there exist no prior results.
	
	For $ C_n $-graded groups in Shi's sense with $ n \ge 3 $, a weaker version of our coordinatisation result (\cref{BC:thm}) can be found in 3.4.19 of Zhang's PhD thesis \cite{Zhang}. The result in \cite{Zhang} is obtained under the additional assumption that the root groups are 2-divisible, and it remains a conjecture that the coordinatising ring is alternative. Further, Zhang overlooks the examples of $ C_n $-graded groups which are coordinatised not only by a ring $ \ring $ with involution $ \rinvmap $, but by a tuple $ (\ring, \module) $ where $ \module $ is an additional module over $ \ring $. See \cref{sec:BC-example} and, in particular, \cref{BC:ex:C} for more details on these additional examples. We can solve all these problems in the setting of $ BC_n $-gradings (which properly generalise $ C_n $-gradings) using the blueprint technique.
\end{miscthm}

	\chapter{Chevalley Groups}
	
	\label{chap:chev}
	The most well-known examples of root graded groups are the Chevalley groups, which were introduced by Chevalley in his famous paper \cite{Chev-Tohoku} in 1955. By construction, every Chevalley group is coordinatised by a commutative associative ring $ \ring $ in the sense of \cref{rgg:coord-gen-def}: Its root groups are isomorphic to the additive group $ (\ring, +) $ and commutators of root group elements are described by the celebrated Chevalley commutator formula.
	
	In this chapter, we briefly review the foundations of the theory of Chevalley groups. While these are not needed in a technical sense to develop our theory of root graded groups, they provide a welcome motivation and illustration of several important concepts and notions. In \cref{sec:chev:basic}, we set up our notation concerning Chevalley groups and state the Chevalley commutator formula. In \cref{sec:chev:weyl}, we take a closer look at Weyl elements in Chevalley groups.
	
	Standard references on the subject are \cite{HumphreysLieAlg} for the basic theory of Lie algebras, which we assume the reader is familiar with, and \cite{Steinberg-ChevGroups,Carter-Chev} for the material on Chevalley groups. It should be noted that \cite[Chapter~VII]{HumphreysLieAlg} also contains the construction of Chevalley groups (but not many properties of these groups beyond that) and that \cite{Carter-Chev} only covers adjoint Chevalley groups, not general Chevalley groups. Further, all three references have the disadvantage that they only introduce Chevalley groups over fields and not over commutative associative rings. However, all the basic results in this chapter remain valid over commutative associative rings, and only minor modifications of the proofs are necessary. A survey of Chevalley groups in this general setting is \cite{VavilovPlotkin}. It contains most of the results which we introduce in this chapter, albeit usually without proofs.
	

\section{Basic Properties of Chevalley Groups}

\label{sec:chev:basic}

Throughout this section, we use the notation of the first two chapters of \cite{HumphreysLieAlg} on the basic theory of semisimple Lie algebras. Further, we assume that the reader is familiar with this theory.

\begin{secnotation}
	We denote by $ \comring $ an algebraically closed field of characteristic~$ 0 $, by $ L $ a semisimple Lie algebra over $ \comring $, by $ H $ a maximal toral subalgebra of $ L $ and by $ \roots \defl \roots(L,H) \subs \Hom_\comring(H, \comring) $ the (crystallographic reduced) root system of $ L $ with respect to $ H $. We denote the rank of $ \roots $ by $ \ell $.
\end{secnotation}

\begin{theorem}[Root space decomposition, {\cite[8.1, 8.2, 8.4]{HumphreysLieAlg}}]
	As a $ \comring $-vector space, we have the following \defemph*{root space decomposition of $ L $}\index{Lie algebra!root space decomposition}:
	\[ L = H \dirsum \bigoplus_{\alpha \in \roots(L,H)} L_\alpha \]
	where $ L_\alpha \defl \Set{x \in L \given \lie{h}{x} = \alpha(h) x \text{ for all } h \in H} $ for all $ \alpha \in \Hom_\comring(H, \comring) $. Further, $ \dim L_\alpha = 1 $ for all roots $ \alpha $ and $ \commutator{L_\alpha}{L_\beta} = L_{\alpha+\beta} $ for all roots $ \alpha, \beta $.
\end{theorem}

\begin{definition}[Chevalley basis, {\cite[Proposition 25.2]{HumphreysLieAlg}}]\label{chev:chevbasis-def}
	Choose an ordered root base $ \rootbase = \tup{\alpha}{\ell} $ of $ \roots $. A \defemph*{Chevalley basis of $ L $ of type $ \rootbase $ (with respect to $ H $)}\index{Chevalley basis} is a family $ \chevbasis = (x_\alpha)_{\alpha \in \roots} \union (h_i)_{i \in \numint{1}{\ell}} $ with the following properties:
	\begin{stenumerate}
		\item $ x_\alpha \in L_\alpha $ for all $ \alpha \in \roots $.
		
		\item $ h_i = h_{\alpha_i} $ for all $ i \in \numint{1}{\ell} $, where $ h_{\alpha_i} $ is as in \cite[Proposition~8.3]{HumphreysLieAlg}.
		
		\item $ \lie{x_\alpha}{x_{-\alpha}} = h_\alpha $ for all $ \alpha \in \roots $.
		
		\item \label{chev:chevbasis-def:const}Let $ \alpha, \beta $ be roots such that $ \alpha+\beta $ is a root. If the scalars $ \chevstr_{\alpha, \beta} , \chevstr_{-\alpha, -\beta} \in \comring $ are defined by $ \lie{x_\alpha}{x_\beta} = \chevstr_{\alpha, \beta} x_{\alpha+\beta} $ and $ \lie{x_{-\alpha}}{x_{-\beta}} = \chevstr_{-\alpha, -\beta} x_{-\alpha-\beta} $, then $ \chevstr_{-\alpha, -\beta} = -\chevstr_{\alpha, \beta} $.
	\end{stenumerate}
	A \defemph*{Chevalley basis of $ L $ (with respect to $ H $)} is a Chevalley basis of some type. The scalars $ (\chevstr_{\alpha, \beta})_{\alpha, \beta \in \roots} $ defined as in~\ref{chev:chevbasis-def:const} (and by $ \chevstr_{\alpha, \beta} \defl 0 $ if $ \alpha+ \beta $ is not a root) are called the \defemph*{(Chevalley) structure constants of $ L $ (with respect to $ \chevbasis $)}.\index{Chevalley structure constants}
\end{definition}

\begin{example}\label{chev:An-chevbas}
	The \enquote{standard basis} of the Lie algebra $ \sllie_{\ell+1}(\IC) $ (see \cite[1.2]{HumphreysLieAlg}) is a Chevalley basis, and the corresponding \enquote{standard} family of structure constants (of type $ A_\ell $) is given by
	\[ \chevstr_{\basvec_i - \basvec_j, \basvec_j - \basvec_k} = 1, \qquad \chevstr_{\basvec_j - \basvec_k, \basvec_i - \basvec_j} = -1 \]
	and $ \chevstr_{\alpha, \beta} = 0 $ for all other pairs $ \alpha, \beta $ of roots.
\end{example}

We now cite some essential properties of Chevalley bases.

\begin{lemma}[{\cite[{Proposition 25.2}]{HumphreysLieAlg}}]
	There exists a Chevalley basis of $ L $.
\end{lemma}

\begin{lemma}\label{chev:strconst-prop}
	Let $ \chevbasis $ be a Chevalley basis of $ L $. Then the Chevalley structure constants of $ L $ with respect to $ \chevbasis $ have the following properties:
	\begin{lemenumerate}
		\item \label{chev:strconst-prop:switch}For all roots $ \alpha, \beta $, we have $ \chevstr_{\alpha, \beta} = -\chevstr_{\beta, \alpha} $ and $ \chevstr_{-\alpha, -\beta} = -\chevstr_{\alpha, \beta} $.
		
		\item \label{chev:strconst-prop:string}If $ \alpha, \beta $ are roots such that $ \alpha+\beta $ is a root, then $ \chevstr_{\alpha,\beta} \in \Set{\pm (r+1)} $ where $ r $ is the maximal non-negative integer such that $ \beta-r\alpha $ is a root.
		
		\item \label{chev:strconst-prop:int}All structure constants of $ L $ with respect to $ \chevbasis $ are integers, and their absolute values are uniquely determined by the root system $ \roots $.
		
		\item \label{chev:strconst-prop:sign}If $ \chevbasis' $ is another Chevalley basis of $ L $ with structure constants $ (\chevstr'_{\alpha, \beta})_{\alpha, \beta \in \roots} $, then $ \chevstr_{\alpha, \beta}' \in \Set{\pm \chevstr_{\alpha, \beta}} $ for all roots $ \alpha, \beta $.
		
		\item \label{chev:strconst-prop:pos}If $ \chevbasis' $ is another Chevalley basis of $ L $ with structure constants $ (\chevstr'_{\alpha, \beta})_{\alpha, \beta \in \roots} $ such that $ \chevstr_{\alpha, \beta} = \chevstr_{\alpha, \beta}' $ for all roots $ \alpha, \beta $ lying in some positive system $ \possys $ in $ \roots $, then $ \chevstr_{\alpha, \beta} = \chevstr_{\alpha, \beta}' $ for all roots $ \alpha, \beta $.
	\end{lemenumerate}
\end{lemma}
\begin{proof}
	The property $ \chevstr_{\alpha, \beta} = -\chevstr_{\beta, \alpha} $ follows from the anti-commutativity of the Lie bracket while the second part of~\itemref{chev:strconst-prop:switch} is part of the definition of a Chevalley basis. A proof of~\itemref{chev:strconst-prop:string} can be found in \cite[{Theorem 25.2}]{HumphreysLieAlg}. Parts~\itemref{chev:strconst-prop:int} and~\itemref{chev:strconst-prop:sign} are consequences of~\itemref{chev:strconst-prop:string}. Property~\itemref{chev:strconst-prop:pos} holds by \cite[(14.8)]{VavilovPlotkin}.
\end{proof}

\begin{remark}\label{chev:simply-laced-struct-1}
	It follows from \thmitemcref{chev:strconst-prop}{chev:strconst-prop:string} that in a simply-laced root system, we have $ \chevstr_{\alpha, \beta} \in \compactSet{\pm 1} $ for all roots $ \alpha,\beta $ for which $ \alpha+\beta $ is a root.
\end{remark}

\begin{definition}\label{chev:struc-family-def}
	Let $ \roots' $ be a reduced crystallographic root system and let $ \chevstr = (\chevstr_{\alpha, \beta})_{\alpha, \beta \in \roots} $ be a family of integers. We say that $ \chevstr $ is a \defemph*{family of Chevalley structure constants of type $ \roots' $}\index{Chevalley structure constants} if there exist a semisimple complex Lie algebra $ L' $ with root system $ \roots' $ and a Chevalley basis $ \chevbasis $ of $ L' $ such that $ (\chevstr_{\alpha, \beta})_{\alpha, \beta \in \roots} $ are the Chevalley structure constants of $ L' $ with respect to $ \chevbasis $.
\end{definition}

\begin{note}[The sign problem for Chevalley bases]\label{chev:signprob-chevbasis}
	We emphasise the statement of \thmitemcref{chev:strconst-prop}{chev:strconst-prop:sign}: The structure constants of a Chevalley basis are unique only up to a sign. However, not every choice of signs corresponds to a Chevalley basis. This raises the following natural question: How much freedom do we have in the choice of signs, and can we describe an algorithm which produces an (or any) explicit choice of signs? This problem is well-known in the theory of Chevalley groups, and we will refer to it as the \defemph*{sign problem (for Chevalley groups)}\index{sign problem!for Chevalley bases}. In \cite[Section~4.2]{Carter-Chev}, it is shown that the signs can be chosen arbitrarily for so-called extraspecial pairs of roots, and then all other signs are determined by these arbitrary choices. More information can also be found in \cite[Sections~14,~15]{VavilovPlotkin}.
	
	A direct consequence of the sign problem for Chevalley bases is that the signs in the Chevalley commutator formula are not unique either. Similarly, the signs which appear in certain conjugation formulas for Weyl elements depend on the signs of the Chevalley basis. See \cref{chev:signprob-commformula,chev:signprob-weyl} for more details. Further, a variation of the sign problem is the main difficulty in the parametrisation of root graded groups. We will elaborate on this in \cref{param:motiv:parmap-choice}.
\end{note}

\begin{definition}\label{chev:def-chevgrp}
	Let $ \ring $ be an associative commutative ring and let $ \roots' $ be a crystallographic reduced root system. A \defemph*{Chevalley group of type $ \roots' $ over $ \ring $}\index{Chevalley group} is a group $ G $ defined as in \cite[27.4]{HumphreysLieAlg}, \cite[Section~6]{VavilovPlotkin} or \cite[Chapter~3]{Steinberg-ChevGroups} (with respect to some semisimple Lie algebra $ L' $ over $ \comring $ whose root system is $ \roots' $, some Chevalley basis $ \chevbasis $ of $ L' $, some finite-dimensional faithful $ L' $-module $ V $ and some admissible lattice in $ V $ with respect to $ \chevbasis $). By construction, a Chevalley group $ G $ over $ \ring $ is equipped with a family of subgroups $ (\rootgr{\alpha})_{\alpha \in \roots'} $, called the \defemph*{root subgroups}, and a family $ (\map{\risom{\alpha}}{(\ring, +)}{\rootgr{\alpha}}{}{})_{\alpha \in \roots'} $ of group isomorphisms which we call \defemph*{root isomorphisms}\index{root homomorphism}. If $ \chevstr = (\chevstr_{\alpha, \beta})_{\alpha, \beta \in \roots'} $ is the family of Chevalley structure constants of a Chevalley basis with respect to which $ G $ is defined, we say that $ G $ is defined \defemph*{with respect to $ \chevstr $}. If $ G $ is defined with respect to the adjoint module $ V=L $, then $ G $ is called an \defemph*{adjoint Chevalley group}\index{Chevalley group!adjoint}.
\end{definition}

\begin{remark}\label{chev:data-remark}
	As is evident from \cref{chev:def-chevgrp}, it is incorrect to refer to \emph{the} Chevalley group of type $ \roots $ over $ \ring $ because many more input data are required to construct a Chevalley group. However, all properties that we will be interested only depend on $ \roots $ and $ \ring $, so we do not need to distinguish any further. In fact, many properties are even independent of the base ring in a suitable sense.
\end{remark}

\begin{remark}
	Any Chevalley group over a commutative associative ring $ \ring $ is a group of automorphisms of some free $ \ring $-module of finite rank, and hence a matrix group over $ \ring $.
\end{remark}

\begin{note}[Chevalley groups over rings]
	In many references, for example in \cite{HumphreysLieAlg,Carter-Chev,Steinberg-ChevGroups}, Chevalley groups are only defined over fields, and many results do not extend to rings. For example, adjoint Chevalley groups are nearly always simple (see \cite[Chapter~4]{Steinberg-ChevGroups}), but Chevalley groups over rings containing a proper ideal are not simple. However, the basic construction of Chevalley groups remains the same over rings, and our definition agrees with the one in \cite[Section~6]{VavilovPlotkin}.
\end{note}

\begin{example}[Chevalley groups of type $ A_\ell $]\label{chev:ex-SL}
	Let $ \ring $ be an associative commutative ring. We use the standard representation $ A_\ell = \Set{\basvec_i - \basvec_j \given i \ne j \in \numint{1}{\ell+1}} $ of the root system $ A_\ell $ where $ (\basvec_i)_{i \in \numint{1}{\ell+1}} $ denotes the standard basis of $ \IR^{\ell+1} $. For all roots $ \alpha = \basvec_i - \basvec_j \in A_\ell $, we put
	\[ \map{\risom{\alpha}}{(\ring, +)}{M_{\ell+1}(\ring)}{r}{\id_{\ell+1} + re_{ij}} \]
	where $ e_{ij} $ denotes the matrix with $ 1_\ring $ at position $ (i,j) $ and $ 0_\ring $ at every other position. This defines a family $ (\risom{\alpha})_{\alpha \in A_{\ell}} $ of injective group homomorphisms, whose images we denote by $ (\rootgr{\alpha})_{\alpha \in A_{\ell}} $. Matrices in $ \bigunion_{\alpha \in A_\ell} \rootgr{\alpha} $ are called \defemph*{elementary matrices}\index{elementary matrix} and the group $ \E_{\ell+1}(\ring) $ generated by all elementary matrices is called the \defemph*{elementary group (over $ \ring $)}\index{elementary group}. It is contained in the special linear group $ \SL_{\ell+1}(\ring) $, and we even have $ \SL_{\ell+1}(\ring) = \E_{\ell+1}(\ring) $ if $ \ring $ is a field. In any case, $ \E_{\ell+1}(\ring) $ is a Chevalley group of type $ A_\ell $ over $ \ring $ with root groups $ (\rootgr{\alpha})_{\alpha \in A_{\ell}} $ and root isomorphisms $ (\risom{\alpha})_{\alpha \in A_{\ell}} $. It can be constructed from the Lie algebra $ \sllie_{\ell+1}(\IC) $ and its natural module $ \IC^{\ell +1} $.
	
	The elementary group satisfies the following commutator relations:
	\[ \commutator{\risom{\basvec_i - \basvec_j}(r)}{\risom{\basvec_j - \basvec_k}(s)} = \risom{\basvec_i - \basvec_k}(rs) \midand \commutator{\risom{\basvec_j - \basvec_k}(r)}{\risom{\basvec_i - \basvec_j}(s)} = \risom{\basvec_i - \basvec_k}(-rs) \]
	for all pairwise distinct $ i,j,k \in \numint{1}{\ell+1} $ and all $ r,s \in \ring $. Note that each of these two formulas is a simple consequence of the other one by the relation $ \commutator{x}{y}^{-1} = \commutator{y}{x} $.
	
	Observe that elementary matrices and the elementary group can also be defined for an arbitrary associative ring $ \ring $ which is not assumed to be commutative. The resulting group is no longer a Chevalley group, but it is still an $ A_\ell $-graded group. In this case, the commutator relation $ \commutator{\risom{\basvec_i - \basvec_j}(r)}{\risom{\basvec_j - \basvec_k}(s)} = \risom{\basvec_i - \basvec_k}(rs) $ remains valid, but we have $ \commutator{\risom{\basvec_j - \basvec_k}(r)}{\risom{\basvec_i - \basvec_j}(s)} = \risom{\basvec_i - \basvec_k}(-sr) $.
\end{example}

Finally, we can state the Chevalley commutator formula, which is the most important result in this chapter. Not only does it say that a Chevalley group has crystallographic $ \roots $-commutator relations with respect to its root groups, but it also provides explicit formulas for the commutator of two root group elements.

\begin{notation}\label{chevorder:def}
	For any pair of non-proportional roots $ \alpha $, $ \beta $, we put\index{I(alpha, beta)@$ \chevorder(\alpha, \beta) $}
	\[ \chevorder(\alpha, \beta) \defl \Set{(i,j) \in \Npos^2 \given i\alpha + j\beta \in \roots} = \Set{(i,j) \in \Npos^2 \given i\alpha + j\beta \in \rootint{\alpha}{\beta}}. \]
\end{notation}

\begin{theorem}[Chevalley commutator formula, {\cite[Lemma 15]{Steinberg-ChevGroups}}]\label{chev:comm-formula}
	Let $ \alpha, \beta $ be non-proportional roots and put $ \chevorder \defl \chevorder(\alpha, \beta) $. Endow the set $ \chevorder $ with an arbitrary but fixed order. Then there exists a unique family $ (\chevstr_{\alpha, \beta, i, j})_{(i,j) \in \chevorder} $ of integers (which may depend on the chosen order on $ \chevorder $) such that for all commutative associative rings $ \ring $ and for all Chevalley groups $ G $ of type $ \roots $ over $ \ring $, the corresponding root homomorphisms $ (\risom{\gamma})_{\gamma \in \roots} $ satisfy
	\[ \commutator{\risom{\alpha}(r)}{\risom{\beta}(s)} = \prod_{(i,j) \in \chevorder} \risom{i\alpha + j\beta}\brackets[\big]{\chevstr_{\alpha, \beta, i, j} r^i s^j} \]
	for all $ r,s \in \ring $. In particular, $ G $ has crystallographic $ \roots $-commutator relations with root groups $ (\rootgr{\alpha})_{\alpha \in \roots} $.
\end{theorem}

\begin{definition}[Chevalley structure constants]\label{chev:grp-const-def}
	The integers $ \chevstr_{\alpha, \beta, i, j} $ (where $ \alpha, \beta $ are non-proportional roots and $ (i,j) \in \chevorder(\alpha, \beta) $) in \cref{chev:comm-formula} are called the \defemph*{(Chevalley) structure constants of (the Chevalley group) $ G $}.
\end{definition}

\begin{remark}\label{chev:struc-const-formulas}
	The following slightly different version of the commutator formula is given in \cite[5.2.2]{Carter-Chev}:
	\[ \risom{\alpha}(u)^{-1} \risom{\beta}(t)^{-1} \risom{\alpha}(u) \risom{\beta}(t) = \commutator{\risom{\alpha}(u)}{\risom{\beta}(t)} = \prod_{(j, i) \in \chevorder} \risom{i\beta + j\alpha}\brackets[\big]{\chevstr_{i,j,\beta, \alpha}' (-t)^i u^j}. \]
	Clearly, the structure constants in this expression are related to our structure constants by $ \chevstr_{\alpha,\beta,i,j} = (-1)^i \chevstr_{i,j,\beta,\alpha}' $. By the formulas in \cite[5.2.2]{Carter-Chev}, the structure constants $ \chevstr_{\alpha, \beta, i, j} $ can be computed from the structure constants $ (\chevstr_{\alpha, \beta})_{\alpha, \beta \in \roots} $, provided that the order $ \prec $ on $ \chevorder(\alpha, \beta) $ has the property that $ i+j < i'+j' $ implies $ (i,j) \prec (i', j') $. In this case, we have $ \chevstr_{\alpha, \beta, 1, 1} = \chevstr_{\alpha,\beta} $. Further, the values of the structure constants $ \chevstr_{\alpha,\beta,i,j} $ always lie in $ \Set{\pm 1, \pm 2, \pm 3} $.
\end{remark}

\begin{note}[The sign problem in the Chevalley commutator formula]\label{chev:signprob-commformula}
	We have seen in \cref{chev:struc-const-formulas} that the structure constants in the Chevalley commutator formula can be computed from the structure constants of the Chevalley bases. Since the latter are only unique up to a sign by \cref{chev:signprob-chevbasis}, it follows that the former are only unique up to a sign as well. This is the second occurrence of the sign problem in the theory of Chevalley groups.\index{sign problem!for the Chevalley commutator formula}
\end{note}

\begin{remark}[Injectivity of the root homomorphisms]\label{chev:roothom-inj}
	Let $ \alpha $ be any root. The homomorphism $ \map{\risom{\alpha}}{(\ring,+)}{\rootgr{\alpha}}{}{} $ is surjective by the definition of $ \rootgr{\alpha} $, and it is in fact also injective. This is proven in \cite[Corollary~1 of Lemma~17]{Steinberg-ChevGroups} for fields, and it remains true over commutative associative rings. For example, this is stated in \cite[below (9.1)]{VavilovPlotkin}.
\end{remark}


\section{Weyl Elements in Chevalley Groups}

\label{sec:chev:weyl}

\begin{secnotation}\label{chev:conv:weyl}
	We denote by $ \roots $ a crystallographic reduced root system, by $ (\chevstr_{\alpha, \beta})_{\alpha, \beta \in \roots} $ a family of Chevalley structure constants of type $ \roots $, by $ \ring $ a commutative associative ring and by $ G $ a Chevalley group of type $ \roots $ over $ \ring $ with respect to $ (\chevstr_{\alpha, \beta})_{\alpha, \beta \in \roots} $. Further, we denote the root groups of $ G $ by $ (\rootgr{\alpha})_{\alpha \in \roots} $ and its root homomorphisms by $ (\risom{\alpha})_{\alpha \in \roots} $.
\end{secnotation}

Since we already know from the Chevalley commutator formula that $ G $ has crystallographic $ \roots $-commutator relations with root groups $ (\rootgr{\alpha})_{\alpha \in \roots} $, it remains to establish the existence of Weyl elements in $ G $. However, we will not settle for a mere existence result. Instead, we will give explicit formulas for the conjugation action of Weyl elements on the root groups, and we will derive properties of the signs which appear in these formulas. In particular, we will encounter the square formula for Weyl elements (see \cref{chev:square-formula-note}) which plays an important role in the theory of root graded groups (see Definitions~\thmitemref{parmap:prop-def}{parmap:prop-def:square-formula} and~\ref{param:square-formula-rgg-def}).

\begin{note}[Validity of the results over rings]\label{chev:steinberg-ring-note}
	Essentially all results in this sections are consequences of \cite[Lemma~19~(a),~(c)]{Steinberg-ChevGroups}. It is slightly problematic that Lemma~19 is only formulated for fields. However, the statements of Lemma~19~(a) and~(c) remain valid over commutative associative rings: For each equation in Lemma~19~(a) and~(c), we can find a certain integral Laurent polynomial
	\[ f \in \IZ[t_1, \ldots, t_n, u_1, \ldots, u_m, u_1^{-1}, \ldots, u_m^{-1}] \]
	such that the validity of the desired equation (over $ \ring $) is equivalent to the assertion that $ f $ vanishes at all points in $ \ring^n \times (\inv{\ring})^m $ (where $ \inv{\ring} $ is the set of invertible elements in $ \ring $). If the desired results holds over $ \IC $, then $ f $ must be the zero polynomial, and hence it vanishes at all points in any commutative associative ring. Thus if Lemma~19~(a) and~(c) hold over $ \IC $, then they hold over all commutative associative rings. In fact, Steinberg only proves \cite[Lemma~19]{Steinberg-ChevGroups} over $ \IC $ and then invokes the same argument that we have just described to infer that it holds over arbitrary fields.
	
	Similarly, it can be shown that statement of Lemma~19~(b) in \cite{Steinberg-ChevGroups} remains true over commutative associative rings. However, the existence assertion in Lemma~19~(b) makes the required argument more subtle. In any case, we will not need this result.
	
	The validity of \cite[Lemma~20~(b),~(c)]{Steinberg-ChevGroups} over commutative associative rings can be proven by similar arguments. Alternatively, since these results are direct consequences of Lemma~19~(a) and~(c), one can also apply the same proof as in \cite{Steinberg-ChevGroups}, but with references to Lemma~19 replaced by references to the more general assertions for commutative associative rings.
\end{note}

We will show that the elements in the following definition are Weyl elements, and investigate their properties. One can even show that all Weyl elements in $ G $ are of this form.

\begin{definition}\label{chev:weyl-def}
	Let $ \alpha $ be a root and let $ r $ be an invertible element of $ \ring $. We define the following two elements of $ G $:
	\[ w_\alpha(r) \defl \risom{-\alpha}(-r^{-1}) \risom{\alpha}(r) \risom{-\alpha}(-r^{-1}), \qquad w_\alpha'(r) \defl \risom{\alpha}(r) \risom{-\alpha}(-r^{-1}) \risom{\alpha}(r). \]
	Further, we put $ w_\alpha \defl w_\alpha(1) $ and $ w_\alpha' \defl w'_\alpha(1) $. The elements $ (w_\alpha)_{\alpha \in \roots} $ are called the \defemph*{standard Weyl elements in $ G $}\index{Weyl element!standard}.
\end{definition}

\begin{note}
	Clearly, we have the following relations:
	\begin{align*}
		w_\alpha(r)^{-1} &= w_\alpha(-r), & w'_\alpha(r)^{-1} &= w'_\alpha(-r), \\
		 w_\alpha(r) &= w_{-\alpha}'(-r^{-1}), & w_\alpha &= w_{-\alpha}'(-1) = (w_{-\alpha}')^{-1}.
	\end{align*}
	The reason to introduce both $ w_\alpha(r) $ and $ w_\alpha'(r) $ is that \cite{Steinberg-ChevGroups} uses the elements $ w'_\alpha(r) $ while we prefer to work with the elements $ w_\alpha(r) $. For this reason, we will first cite the results in \cite{Steinberg-ChevGroups} involving $ w_\alpha'(r) $ and then state the analogous results involving $ w_\alpha(r) $ as a corollary.
\end{note}

The essential result for our purposes is the following one.

\begin{lemma}\label{chev:squarerel-steinberg}
	For any pair of roots $ \alpha, \beta $, there exists a unique $ \inverparsym_\ring'(\alpha, \beta) \in \compactSet{\pm 1_\ring} $ (independent of all input data over which the Chevalley group $ G $ is defined, except for the choice of $ \chevstr $ and $ \ring $) such that
	\[ w_\alpha' \risom{\beta}(\lambda) (w_\alpha')^{-1} = \risom{\refl{\alpha}(\beta)}\brackets[\big]{\inverparsym_\ring'(\alpha, \beta) \lambda} \]
	for all $ \lambda \in \ring $. Further, these constants satisfy
	\[ \inverparsym_\ring'(\alpha, \beta) = \inverparsym_\ring'(\alpha, -\beta) \midand \inverparsym_\ring'(\alpha, \beta) \inverparsym_\ring'\brackets[\big]{\alpha, \refl{\alpha}( \beta)} = (-1)^{\cartanint{\beta}{\alpha}}. \]
\end{lemma}
\begin{proof}
	The first statement is \cite[Lemma 20 (b)]{Steinberg-ChevGroups} and the relation $ \inverparsym_\ring'(\alpha, \beta) = \inverparsym_\ring'(\alpha, -\beta) $ is part of \cite[Lemma 19 (a)]{Steinberg-ChevGroups}. For any invertible element $ r $ of $ \ring $, we define $ h_\alpha'(r) \defl w_\alpha'(r) w_\alpha'(-1) $ (as in \cite[Lemma 19]{Steinberg-ChevGroups}). Then \cite[Lemma~20~(c)]{Steinberg-ChevGroups} says that
	\[ h_\alpha'(r) \risom{\beta}(\lambda) h_\alpha'(r)^{-1} = \risom{\beta}(r^{\cartanint{\beta}{\alpha}} \lambda) \]
	for all invertible $ r \in \ring $ and all $ \lambda \in \ring $. Applying this statement with $ r \defl -1_\ring $ (so that $ h_\alpha'(r) = w_\alpha'(-1)^2 = (w_\alpha')^{-2} $), we infer that
	\[ (w_\alpha')^{-2} \risom{\beta}(\lambda) (w_\alpha')^2 = \risom{\beta}((-1)^{\cartanint{\beta}{\alpha}} \lambda) \]
	for all $ \lambda \in \ring $. Thus
	\begin{align*}
		\risom{\beta}(\lambda) &= (w_\alpha')^2 \risom{\beta}\brackets[\big]{(-1)^{\cartanint{\beta}{\alpha}} \lambda} (w_\alpha')^{-2} = w_\alpha' \risom{\refl{\alpha}(\beta)}\brackets[\big]{\inverparsym_\ring'(\alpha, \beta)(-1)^{\cartanint{\beta}{\alpha}} \lambda} (w_\alpha')^{-1} \\
		&=\risom{\beta}\brackets[\big]{\inverparsym_\ring'(\alpha, \refl{\alpha}(\beta))\inverparsym_\ring'(\alpha, \beta)(-1)^{\cartanint{\beta}{\alpha}} \lambda}
	\end{align*}
	for all $ \lambda \in \ring $. Since $ \risom{\beta} $ is injective by \cref{chev:roothom-inj}, the last assertion follows by putting $ \lambda \defl 1_\ring $.
\end{proof}

\begin{note}
	Lemma~19 in \cite{Steinberg-ChevGroups} is actually formulated for all elements $ w_\alpha'(\lambda) $ where $ \lambda \in \ring $ is invertible, and not just for $ \lambda = 1_\ring $. Lemma~20, on the other hand, is not stated in this generality. An adaption of the proof of Lemma~20 makes it possible to show that
	\[ w_\alpha'(r) \risom{\beta}(\lambda) (w_\alpha'(r))^{-1} = \risom{\refl{\alpha}(\beta)}\brackets[\big]{\inverparsym_\ring'(\alpha, \beta) r^{-\cartanint{\beta}{\alpha}} \lambda} \]
	for all invertible $ r \in \ring $. However, the weaker statement of \cref{chev:squarerel-steinberg} will be sufficient for our purposes.
\end{note}

In order to suit the conventions of this work, we will now write down the assertion of \cref{chev:squarerel-steinberg} for the elements $ w_\alpha $ in place of $ w_\alpha' $ and for $ \inverparsym_\ring(\beta, \alpha) $ in place of $ \inverparsym_\ring'(\alpha, \beta) $.

\begin{lemma}\label{chev:squarerel}
	For any pair of roots $ \alpha, \beta $, there exists a unique $ \inverparsym_\ring(\beta, \alpha) \in \compactSet{\pm 1_\ring} $ (independent of all input data over which the Chevalley group $ G $ is defined, except for the choices of $ \roots $, $ \chevstr $ and $ \ring $) such that
	\[ w_\alpha^{-1} \risom{\beta}(\lambda) w_\alpha = \risom{\refl{\alpha}(\beta)} \brackets[\big]{\inverparsym_\ring(\beta, \alpha) \lambda} \]
	for all $ \lambda \in \ring $. Further, these constants satisfy
	\[ \inverparsym_\ring(\beta, \alpha) = \inverparsym_\ring(-\beta, \alpha) \midand \inverparsym_\ring(\beta, \alpha) \inverparsym_\ring(\refl{\alpha}(\beta), \alpha) = (-1)^{\cartanint{\beta}{\alpha}}. \]
\end{lemma}
\begin{proof}
	Using that $ w_\alpha^{-1} = w_{-\alpha}' $, the first assertion follows from \cref{chev:squarerel-steinberg}, and we see that $ \inverparsym_\ring(\beta, \alpha) = \inverparsym_\ring'(-\alpha, \beta) $. Now
	\[ \inverparsym_\ring(\beta, \alpha) = \inverparsym_\ring'(-\alpha, \beta) = \inverparsym_\ring'(-\alpha, -\beta) = \inverparsym_\ring(-\beta, \alpha) \]
	and
	\begin{align*}
		\inverparsym_\ring(\beta, -\alpha) \inverparsym_\ring\brackets[\big]{\refl{\alpha}(\beta), -\alpha} &= \inverparsym_\ring'(\alpha, \beta) \inverparsym_\ring'\brackets[\big]{\alpha, \refl{\alpha} (\beta)} = (-1)^{\cartanint{\beta}{\alpha}}.
	\end{align*}
	Since $ \refl{\alpha} = \refl{-\alpha} $, the latter equality implies that
	\begin{align*}
		\inverparsym_\ring(\beta, \alpha) \inverparsym_\ring\brackets[\big]{\refl{\alpha}(\beta), \alpha} = (-1)^{\cartanint{\beta}{-\alpha}} = (-1)^{-\cartanint{\beta}{\alpha}} = (-1)^{\cartanint{\beta}{\alpha}},
	\end{align*}
	which finishes the proof.
\end{proof}

\begin{note}[The sign problem for Weyl elements in Chevalley groups]\label{chev:signprob-weyl}
	Recall from \cref{chev:signprob-chevbasis} that the signs of the structure constants $ (\chevstr_{\alpha, \beta})_{\alpha, \beta \in \roots} $ are not uniquely determined by $ \roots $. It follows that the same holds for the \enquote{signs} $ \inverparsym_\ring'(\alpha, \beta) $ and $ \inverparsym_\ring(\beta, \alpha) $ in \cref{chev:squarerel-steinberg,chev:squarerel}. This is the third and final occurrence of the sign problem in the theory of Chevalley groups (see also \cref{chev:signprob-chevbasis,chev:signprob-commformula}).\index{sign problem!for Weyl elements in Chevalley groups} It is essentially a special case of the sign problem that we will be faced with in the parametrisation of root graded groups. See \cref{param:motiv:parmap-choice} for more details.
\end{note}

\begin{remark}\label{chev:conj-const-rem}
	Note that the numbers $ \inverparsym_\ring(\beta, \alpha) $ and $ \inverparsym_\ring'(\alpha, \beta) $ in \cref{chev:squarerel-steinberg,chev:squarerel} lie in $ \compactSet{\pm 1_\ring} $, not in $ \compactSet{\pm 1_\IZ} $. We can of course lift them to numbers $ \tilde{\inverparsym}_\ring(\alpha, \beta) \in \compactSet{\pm 1_\IZ} $ and $ \tilde{\inverparsym}'_\ring(\beta, \alpha) \in \compactSet{\pm 1_\IZ} $, but this lift is only uniquely determined if $ 1_\ring \ne -1_\ring $ (that is, if $ 2_\ring \ne 0_\ring $). Nonetheless, the uniqueness assertions in \cref{chev:squarerel-steinberg,chev:squarerel} remain true if $ 2_\ring = 0_\ring $ because in this case, $ \compactSet{\pm 1_\ring} $ is a singleton and thus, of course, contains a unique element.
	
	In \cite[Lemma~19~(a)]{Steinberg-ChevGroups}, it is shown that if $ \ring $ is a field, then the numbers $ \inverparsym_\ring'(\alpha, \beta) $ \enquote{do not depend on $ \ring $} in the sense that there exists $ c(\alpha, \beta) \in \compactSet{\pm 1_\IZ} $ such that $ \inverparsym_\ring'(\alpha, \beta) = \str_\ring(c(\alpha, \beta)) $ for all fields $ \ring $ (where $ \str_\ring $ denotes the unique ring homomorphism from $ \IZ $ to $ \ring $). The same remains true over commutative associative rings (see \cref{chev:steinberg-ring-note}), and we then of course have $ c(\alpha, \beta) = \inverparsym_\IZ'(\alpha, \beta) $. Similar assertions hold for the elements $ \inverparsym_\ring(\beta, \alpha) $ in \cref{chev:squarerel}.
	
	We conclude that we have constants $ \inverparsym(\beta, \alpha) \defl \inverparsym_\IZ(\beta, \alpha) $ which depend only on $ \roots $ and the choice of a family $ \chevstr $ of Chevalley structure constants. These constants describe (an important part of) the structure of all Chevalley groups of type $ \roots $ which are defined with respect to $ \chevstr $. We have seen that these constants satisfy the relations
	\[ \inverparsym(\beta, \alpha) = \inverparsym(-\beta, \alpha) \midand \inverparsym(\beta, \alpha) \inverparsym(\refl{\alpha}(\beta), -\alpha) = (-1)^{\cartanint{\beta}{\alpha}}. \]
\end{remark}

\begin{definition}[Chevalley parity map]\label{chev:parmap}
	The map $ \map{\inverparsym}{\roots \times \roots}{\compactSet{\pm 1_\IZ}}{}{} $ from \cref{chev:conj-const-rem} is called the \defemph*{Chevalley parity map for $ \chevstr = (\chevstr_{\alpha, \beta})_{\alpha, \beta \in \roots} $}.\index{parity map!Chevalley}\index{Chevalley parity map}
\end{definition}

We will see in \cref{chev:inverpar-comp} that $ \inverparsym $ can be computed from $ \chevstr $, which justifies the name.

\begin{remark}\label{chev:weyl-canon-extend}
	It follows from \cref{chev:squarerel} that for all roots $ \alpha $, $ \beta $, we have
	\begin{align*}
		w_\alpha^{w_\beta} &= \risom{-\alpha}(-1)^{w_\beta} \risom{\alpha}(1)^{w_\beta} \risom{-\alpha}(-1)^{w_\beta} \\
		&= \risom{-\refl{\beta}(\alpha)}\brackets[\big]{-\inverparbr{-\alpha}{\beta}} \risom{\refl{\beta}(\alpha)}\brackets[\big]{\inverparbr{\alpha}{\beta}} \risom{-\refl{\beta}(\alpha)}\brackets[\big]{-\inverparbr{-\alpha}{\beta}} \\
		&= \risom{-\refl{\beta}(\alpha)}\brackets[\big]{-\inverparbr{\alpha}{\beta}} \risom{\refl{\beta}(\alpha)}\brackets[\big]{\inverparbr{\alpha}{\beta}} \risom{-\refl{\beta}(\alpha)}\brackets[\big]{-\inverparbr{\alpha}{\beta}} = w_{\refl{\beta}(\alpha)}^{\inverparbr{\alpha}{\beta}}.
	\end{align*}
	Recall from \cref{rootsys:any-in-rootbase} that every root can be written as $ \delta^{\reflbr{\word{\delta}}} $ for some $ \delta \in \rootbase $ and some word $ \word{\delta} $ over $ \rootbase $. It follows that the whole family $ (w_\alpha)_{\alpha \in \roots} $ of standard Weyl elements can be reconstructed from $ (w_\delta)_{\delta \in \rootbase} $ using only conjugation and group inversion.
	
	We will later see that for any $ \rootbase $-system $ (w_\delta)_{\delta \in \rootbase} $ of Weyl elements in a root graded group $ H $ of rank at least~3, we can find a coordinatisation of $ H $ such that elements $ (w_\delta)_{\delta \in \rootbase} $ behave in a similar way as the standard Weyl elements in Chevalley groups. More precisely, the coordinatisation is constructed in a way which ensures that the elements  $ (w_\delta)_{\delta \in \rootbase} $ satisfy a similar formula as in \cref{chev:squarerel}.
\end{remark}

\begin{note}\label{chev:square-formula-note}
	The formula $ \inverparsym(\beta, \alpha) \inverparsym(\refl{\alpha}(\beta), -\alpha) = (-1)^{\cartanint{\beta}{\alpha}} $ implies that
	\[ \risom{\beta}(\lambda)^{w_\alpha^2} = \risom{\beta}\brackets[\big]{(-1)^{\cartanint{\beta}{\alpha}} \lambda} \]
	for all $ \lambda \in \ring $, so it describes the action of squares of (certain) Weyl elements on the root groups. For this reason, we will also call it the \defemph*{square formula for Weyl elements}.\index{square formula!for Chevalley groups} It will turn out that, except for a few special cases, the action of squares of arbitrary Weyl elements in root graded groups of rank at least~$ 3 $ adheres to the same formula. See \cref{A2Weyl:cartan-comp,B:square-act:summary,BC:square-act:summary} In fact, these are the main intermediate results which have to be obtained in order to parametrise root graded groups.
	
	Observe that the square formula can be formulated purely in terms of the group $ G $ and its root groups $ (\rootgr{\alpha})_{\alpha \in \roots} $, without reference to the root isomorphisms: We have $ x_\beta^{w} = x_\beta^{\epsilon} $
	for all roots $ \alpha, \beta $ and all $ x_\beta \in \rootgr{\beta} $, where $ w \defl w_\alpha^2 $ and $ \epsilon \defl (-1)^{\cartanint{\beta}{\alpha}} $.
\end{note}

\begin{lemma}\label{chev:const-on-itself}
	For all $ \alpha \in \roots $, we have $ \inverparbr{\alpha}{\alpha} = \inverparbr{-\alpha}{\alpha} = -1 $ where $ \inverparbr{\alpha}{\alpha} $ is the constant from \cref{chev:parmap}.
\end{lemma}
\begin{proof}
	Since the statement involves only a single root $ \alpha $, we can reduce to the case that $ \roots = A_1 $. Further, since there do not exist roots $ \alpha, \beta \in A_1 $ such that $ \alpha+\beta $ is a root, there is only one possible way to choose the family $ (\chevstr_{\alpha, \beta})_{\alpha, \beta \in A_1} $. Hence it suffices to verify that in some Chevalley group of type $ A_1 $ over $ \ring $, we have $ \risom{\alpha}(1_\ring)^{w_\alpha} = \risom{-\alpha}(-1_\ring) $ and $ \risom{-\alpha}(1_\ring)^{w_\alpha} = \risom{\alpha}(-1_\ring) $. It is a straightforward matrix computation to do this for the elementary group $ \E_2(\ring) $ from \cref{chev:ex-SL}.
\end{proof}

\begin{remark}[Computation of the Chevalley parity map]\label{chev:inverpar-comp}
	Let $ \alpha $, $ \beta $ be roots. If $ \alpha \in \compactSet{\pm \beta} $, then we know from \cref{chev:const-on-itself} that $ \inverparbr{\beta}{\alpha} = -1 $. If $ \alpha $, $ \beta $ are linearly independent, then we can perform computations similar to the ones in \cref{A2Weyl:basecomp-cor,B:basecomp-walpha-cox-cor-delta} to determine $ \inverparbr{\beta}{\alpha} $ in terms of the structure constants which appear in the Chevalley commutator formula. Since the structure constants in the Chevalley commutator formula can be computed from the structure constants $ \chevstr = (\chevstr_{\gamma, \delta})_{\gamma, \delta \in \roots} $ of the Chevalley basis by \cref{chev:struc-const-formulas}, we conclude that there exists an algorithm to compute $ \inverparbr{\beta}{\alpha} $ from $ \chevstr $. For the simply-laced root systems, the resulting formulas can be found in \cref{ADE:parity-comp}.
\end{remark}

In order to conclude that $ G $ is a $ \roots $-graded group, we have to make a final observation.

\begin{proposition}\label{chev:nondeg}
	There exist a positive system $ \possys $ in $ \roots $ such that the matrices of elements of $ \rootgr{\possys} $ are unipotent upper-triangular while matrices of elements of $ \rootgr{-\possys} $ are unipotent lower-triangular. In particular, $ \rootgr{\possys} \intersect \rootgr{-\possys} = \compactSet{1} $.
\end{proposition}
\begin{proof}
	This is proven in \cite[Corollary~3 of Lemma~17]{Steinberg-ChevGroups} for fields, but the same proof is valid for arbitrary commutative associative rings.
\end{proof}

\begin{theorem}\label{chev:rgg}
	$ (\rootgr{\alpha})_{\alpha \in \roots} $ is a crystallographic $ \roots $-grading of $ G $.
\end{theorem}
\begin{proof}
	By definition, $ G $ is generated by $ (\rootgr{\alpha})_{\alpha \in \roots} $. The Chevalley commutator formula (\cref{chev:comm-formula}) yields that $ G $ has crystallographic $ \roots $-commutator relations with root groups $ (\rootgr{\alpha})_{\alpha \in \roots} $. By \cref{chev:squarerel}, $ (w_\alpha)_{\alpha \in \roots} $ is a family of Weyl elements in $ G $. Finally, Axiom~\thmitemref{rgg-def}{rgg-def:nondeg} is satisfied by \cref{chev:nondeg}.
\end{proof}

\begin{note}[see also \cref{chev:steinberg-ring-note}]
	In its generality, \cref{chev:rgg} relies on the fact that everything we cite from \cite{Steinberg-ChevGroups} remains true over the arbitrary commutative associative ring $ \ring $. However, this fact is only needed to ensure the existence of Weyl elements. Alternatively, the existence of Weyl elements can be deduced from the concrete form of the Chevalley commutator formula. We will see a sketch of the necessary computation in \cref{ADE:stsign:weyl-char}, and we will also use this fact in \cref{B:stsign:weyl-char,BC:stsign:weyl-char}.
\end{note}

	\chapter{The Parametrisation Theorem}
	\label{chap:param}
	
	As explained in the preface and in~\ref{rgg:coord-gen-def}, the goal of this work is to construct a coordinatisation for any crystallographic root graded group $ (G, (\rootgr{\alpha})_{\alpha \in \roots}) $ of rank at least 3. As a first step, we construct a \defemph*{parametrisation of $ (G, (\rootgr{\alpha})_{\alpha \in \roots}) $:} A family of groups $ (M_\alpha)_{\alpha \in \roots} $ (where $ M_\alpha = M_\beta $ if $ \alpha, \beta $ are conjugate under the Weyl group) and a family $ (\map{\risom{\alpha}}{M_\alpha}{\rootgr{\alpha}}{}{})_{\alpha \in \roots} $ of isomorphisms such that
	\begin{equation}\label{eq:param:intro}
		\risom{\alpha}(x)^{w_\delta} = \risom{\refl{\delta}(\alpha)}(\inverpar{\alpha}{\delta}.x)
	\end{equation}
	for all $ \alpha \in \roots $, $ \delta \in \rootbase $ and $ x \in M_\alpha $. Here $ \rootbase $ is a fixed (rescaled) root base of $ \roots $, $ (w_\delta)_{\delta \in \rootbase} $ is a fixed $ \rootbase $-system of Weyl elements and $ (\inverpar{\alpha}{\delta})_{\alpha \in \roots, \delta \in \rootbase} $ is a fixed family of elements of some finite abelian group $ \twistgroup $. The group $ \twistgroup $ is equipped with an action on each of the underlying sets $ (M_\alpha)_{\alpha \in \roots} $, and the pair $ (\twistgroup, (M_\alpha)_{\alpha \in \roots}) $ is called a \defemph*{parameter system for $ G $ with twisting group $ \twistgroup $}. The map $ \map{\inverparsym}{\roots \times \rootbase}{\twistgroup}{}{} $ is called a \defemph*{parity map with values in $ \twistgroup $}.
	
	In the special case of Chevalley groups, each of the groups $ (M_\alpha)_{\alpha \in \roots} $ is the (additive group of the) commutative associative ring $ \ring $ over which $ G $ is defined. The group $ \twistgroup $ is $ \compactSet{\pm 1} $, the action of $ -1 $ on $ (\ring, +) $ is group inversion, the parity map $ \inverparsym $ is the Chevalley parity map from \cref{chev:parmap} (or rather, its restriction to $ \roots \times \rootbase $) and \eqref{eq:param:intro} is precisely the formula in \cref{chev:squarerel}. The essential property that distinguishes general root graded groups from Chevalley groups is that distinct orbits $ O $ in $ \roots $ may be equipped with distinct parametrising groups $ (M_\alpha)_{\alpha \in O} $ and that more general forms of twistings (that is, a twisting group $ \twistgroup $ which is larger than $ \compactSet{\pm 1} $) may occur.
	
	The goal of this chapter is to state and prove the \defemph*{parametrisation theorem}, which provides a practical criterion for the existence of a parametrisation of $ (G, (\rootgr{\alpha})_{\alpha \in \roots}) $. It is formulated and proven for arbitrary root systems $ \roots $ and without reference to the classification of root systems. In this way, the parametrisation theorem allows us to approach the parametrisation of $ \roots $-graded groups in a uniform way, independent of $ \roots $. Similarly, the blueprint technique will later provide a uniform way to equip the parametrising groups $ (M_\alpha)_{\alpha \in \roots} $ that the parametrisation theorem yields with an algebraic structure, thereby turning the parametrisation into a coordinatisation (in the sense of \cref{rgg:coord-gen-def}).
	
	The statement and proof of the parametrisation theorem need some technical machinery. The general idea goes as follows. For each indivisible orbit $ O $ of $ \roots $, choose a root $ \alpha_O \in O $, a group $ M_O $ which is isomorphic to $ \rootgr{\alpha_O} $ and an isomorphism $ \map{\risom{\alpha_O}}{M_O}{\rootgr{\alpha_O}}{}{} $. Put $ M_\beta \defl M_O $ for all indivisible orbits $ O $ and all $ \beta \in O $. Assume that we are given a group $ \twistgroup $ acting on all root groups and a parity map $ \inverparsym $ with values in $ \twistgroup $ such that certain consistency conditions are satisfied. The key observation in the proof of the parametrisation theorem is that for all $ \alpha \in \roots $ and all words $ \word{\delta}, \word{\rho} $ over $ \rootbase \union (-\rootbase) $, we have $ x_\alpha^{w_{\word{\delta}}} = a.(x_\beta^{w_{\word{\rho}}}) $ for all $ x_\alpha \in \rootgr{\alpha} $ where $ a $ is an element of $ \twistgroup $ which is uniquely determined by $ \alpha $, $ \word{\delta} $, $ \word{\rho} $ and $ \inverparsym $ (\cref{param:abstract-thm}). In other words, the action of $ w_{\word{\delta}} $ on $ \rootgr{\alpha} $ depends only on $ \alpha^{\reflbr{\word{\delta}}} $ and $ \inverpar{\alpha}{\word{\delta}} $. Once this is proven, we can use conjugation by the Weyl elements $ (w_\delta)_{\delta \in \rootbase} $ and the parity map $ \inverparsym $ to extend the family $ (\risom{\alpha_O})_{O \in \IndivOrb(\roots)} $ to a family $ (\risom{\alpha})_{\alpha \in \roots} $ which satisfies~\eqref{eq:param:intro}.
	
	The first three sections of this chapter are foundational, setting up the necessary language and machinery: We study parity maps and their combinatorial properties in \cref{sec:parmaps}, parametrisations of $ \roots $-gradings in \cref{sec:param:twist} and certain compatibility conditions between parity maps and parametrisations in \cref{sec:param:comp}. In \cref{sec:param:proof}, we state and prove the parametrisation theorem. The assumptions of the parametrisation theorem involve several of the compatibility conditions from \cref{sec:param:comp}. We provide some practical criteria to prove these conditions for specific twisting structures in \cref{sec:param:crit}. In \cref{sec:param:outline}, we outline our uniform approach to the coordinatisation problem for root graded groups, which builds heavily on the parametrisation theorem.
	
	An elaborate explanation of the motivation behind the precise formulation of several definitions and results in this chapter can be found in \cite[4.1]{Wiedemann-PhD}.
	

\section{Parity Maps}

\label{sec:parmaps}

\begin{secnotation}
	In this section, $ \roots $ is a root system and $ \rootbase $ is a rescaled root base of $ \roots $.
\end{secnotation}

In this section, we study the purely combinatorial notion of parity maps, without reference to any $ \roots $-graded group.

\begin{note}[Rescaled root bases]\label{param:motiv:rescaled}
	The reason to allow $ \rootbase $ to be a rescaled root base (and not a proper root base) is that if $ \roots = BC_n $, we will use the rescaled root base $ \rootbase = (\basvec_1 - \basvec_2, \ldots, \basvec_{n-1} - \basvec_n, 2\basvec_n) $ and not the proper root base $ \rootbase' = (\basvec_1 - \basvec_2, \ldots, \basvec_{n-1} - \basvec_n, \basvec_n) $. Hence when we choose a $ \rootbase $-system $ (w_\delta)_{\delta \in \rootbase} $ of Weyl elements, the element $ w_{2\basvec_n} $ is a $ 2\basvec_n $-Weyl element, which is a stronger property than being an $ \basvec_n $-Weyl element (see \thmitemcref{basic:weyl-general}{basic:weyl-general:div}).
\end{note}

\begin{definition}[$ \rootbase $-expressions]\label{param:Delta-expr}
	For any word $ \word{\alpha} = \tup{\alpha}{m} $ over $ \rootbase \union (-\rootbase) $, we define the corresponding \defemph{inverse word} by $ \word{\alpha}^{-1} \defl (-\alpha_m, \ldots, -\alpha_1) $. For any root $ \alpha $, a \defemph*{$ \rootbase $-expression of $ \alpha $}\index{Delta-expression@$ \rootbase $-expression} is a word over $ \rootbase \union (-\rootbase) $ of the form $ (\word{\delta}^{-1}, \delta, \word{\delta}) $ where $ \delta \in \rootbase $ and $ \word{\delta} $ is a word over $ \rootbase \union (-\rootbase) $ such that $ \delta^{\reflbr{\word{\delta}}} = \alpha $.
\end{definition}

Note that a $ \rootbase $-expression of the root $ \alpha $ is automatically a $ \rootbase $-expression of the element $ \reflbr{\alpha} $ in the sense of \cref{rootsys:length-def}.

\begin{remark}\label{param:Delta-exp-weyl}
	Let $ G $ be a group with a $ \roots $-pregrading $ (\rootgr{\alpha})_{\alpha \in \roots} $, let $ (w_\delta)_{\delta \in \rootbase} $ be a $ \rootbase $-system of Weyl elements and let $ \alpha $ be a root. Then by \thmitemcref{basic:weyl-general}{basic:weyl-general:conj}, $ w_{\word{\rho}} $ is an $ \alpha $-Weyl element for any $ \rootbase $-expression $ \word{\rho} $ of $ \alpha $. Further, $ w_{\word{\alpha}^{-1}} = w_{\word{\alpha}}^{-1} $ for any word $ \word{\alpha} $ over $ \rootbase \union (-\rootbase) $. In fact, these two observation motivate \cref{param:Delta-expr}. Further, note that we have already seen $ \rootbase $-expressions in \cref{rootsys:any-in-rootbase}.
\end{remark}

\begin{definition}[Parity map]\label{param:parmap-def}
	Let $ A $ be an abelian group. A \defemph*{$ \rootbase $-parity map with values in $ A $}\index{parity map} is a map $ \map{\inverparsym}{\roots \times \rootbase}{\twistgroup}{(\alpha, \delta)}{\inverpar{\alpha}{\delta} = \inverparsym(\alpha, \delta)} $ with the property that $ \inverpar{\alpha}{\delta} = \inverpar{\lambda \alpha}{\delta} $ whenever $ \lambda $ is a positive real number such that $ \alpha $ and $ \lambda \alpha $ are roots. Given any such map, we define $ \inverparbr{\alpha}{-\delta} \defl \inverparbr{\alpha^{\reflbr{\delta}}}{\delta}^{-1} $ for all $ \alpha \in \roots $, $ \delta \in \rootbase $ and
	\[ \inverpar{\alpha}{\word{\delta}} \defl \prod_{i=1}^m \inverparbr{\alpha^{\reflbr{\delta_1 \cdots \delta_{i-1}}}}{\delta_i} = \inverparbr{\alpha}{\delta_1} \cdot \inverparbr{\alpha^{\reflbr{\delta_1}}}{\delta_2} \cdots \inverparbr{\alpha^{\reflbr{\delta_1 \cdots \delta_{m-1}}}}{\delta_m} \]
	for any word $ \word{\delta} = \tup{\delta}{m} $ over $ \rootbase \union (-\rootbase) $. We use the convention that $ \inverpar{\alpha}{\emptyset} = 1_A $ for the empty word $ \emptyset $. This defines a map
	\[ \map{\inverparsym}{\roots \times \frmon{\rootbase \union (-\rootbase)}}{\twistgroup}{}{} \]
	which we will sometimes call the \defemph*{extended $ \rootbase $-parity map}\index{parity map!extended}. Here $ \frmon{\rootbase \union (-\rootbase)} $ denotes the free monoid over the set $ \rootbase \union (-\rootbase) $, see \cref{pre:word-def}. When $ \rootbase $ and $ A $ are clear from the context, we will simply say that $ \inverparsym $ is a \defemph{parity map}.
\end{definition}

\begin{example}
	Assume that $ \roots $ is crystallographic and reduced.	For any family $ \chevstr $ of Chevalley structure constants (see \cref{chev:struc-family-def}), we have a map $ \map{\inverparsym}{\roots \times \roots}{\compactSet{\pm 1}}{}{} $ from \cref{chev:parmap} which is called the Chevalley parity map for $ \chevstr $. Its restriction to $ \roots \times \rootbase $ is a $ \rootbase $-parity map with values in $ \compactSet{\pm 1} $.
\end{example}

\begin{note}
	For each root system $ \roots $, we will later define a suitable parity map $ \inverparsym_\roots $ by \enquote{reading it off} from a given $ \roots $-graded group $ G $ with a coordinatisation. See \cref{param:motiv:parmap-choice} for more details.
\end{note}

We begin with some straightforward consequences of \cref{param:parmap-def}.

\begin{remark}\label{param:inverpar-minus-cancel}
	Let $ \inverparsym $ be a parity map with values in some abelian group $ \twistgroup $ and let $ \alpha \in \roots $, $ \delta \in \rootbase $. Then
	\[ \inverparbr{\alpha^{\reflbr{\delta}}}{-\delta} = \inverparbr{\alpha^{\reflbr{\delta} \reflbr{\delta}}}{\delta}^{-1} = \inverparbr{\alpha}{\delta}^{-1}. \]
	This implies that
	\begin{align*}
		\inverpar{\alpha}{(\delta, -\delta)} = \inverparbr{\alpha}{\delta} \inverparbr{\alpha^{\reflbr{\delta}}}{-\delta} = \inverparbr{\alpha}{\delta} \inverparbr{\alpha^{\reflbr{\delta} \reflbr{\delta}}}{\delta}^{-1} = 1_\twistgroup
	\end{align*}
	and that, similarly, $ \inverpar{\alpha}{(-\delta, \delta)} = 1_A $. Thus for any word $ \word{\delta} $ over $ \rootbase \union (-\rootbase) $, we can delete any occurence of $ (\delta, -\delta) $ for any $ \delta \in \rootbase \union (-\rootbase) $ without changing the value of $ \inverpar{\beta}{\word{\delta}} $ for any $ \beta \in \roots $. In other words, the extended parity map $ \map{}{\roots \times \frmon{\rootbase \union (-\rootbase)}}{\twistgroup}{}{} $ defined by $ \inverparsym $ factors through a map $ \map{}{\roots \times \frgrp{\rootbase}}{\twistgroup}{}{} $ where $ \frgrp{\rootbase} $ denotes the free group on $ \rootbase $. Further, since $ \inverpar{\alpha}{-(-\delta)} = \inverpar{\alpha}{\delta} = \inverparbr{\alpha^{\reflbr{\delta}}}{-\delta}^{-1} $, we see that $ \inverparbr{\alpha}{-\rho} = \inverparbr{\alpha^{\reflbr{\rho}}}{\rho}^{-1} $ holds for all $ \rho \in \rootbase \union (-\rootbase) $ (and not merely for all $ \rho \in \rootbase $).
\end{remark}

\begin{lemma}\label{param:inverpar-backwards}
	Let $ \inverparsym $ be a parity map with values in some abelian group $ A $, let $ \alpha \in \roots $ and let $ \word{\beta} = (\beta_1, \ldots, \beta_m) $ be a word over $ \rootbase \union (-\rootbase) $ where $ m \ge 1 $. Then $ \inverparbr{\alpha^{\reflbr{\word{\beta}}}}{\word{\beta}^{-1}} = \inverparbr{\alpha}{\word{\beta}}^{-1} $. In particular, $ \inverpar{\alpha}{\word{\beta}^{-1}} = \inverpar{\alpha}{\word{\beta}}^{-1} $ if $ \alpha^{\reflbr{\word{\beta}}} = \alpha $.
\end{lemma}
\begin{proof}
	Write $ \word{\beta}^{-1} = (\beta_1', \ldots, \beta_m') = (-\beta_m, \ldots, -\beta_1) $. On the one hand,
	\[ \inverpar{\alpha}{\word{\beta}}^{-1} = \brackets*{\prod_{i=1}^m \inverparbr{\alpha^{\reflbr{\beta_1 \cdots \beta_{i-1}}}}{\beta_i}}^{-1} = \prod_{i=1}^m \inverparbr{\alpha^{\reflbr{\beta_1 \cdots \beta_{i-1}}}}{\beta_i}^{-1} \]
	(with the convention that $ \reflbr{\beta_1 \cdots \beta_{i-1}} = \id_\roots $ for $ i=1 $). On the other hand,
	\begin{align*}
		\inverparbr{\alpha^{\reflbr{\word{\beta}}}}{\word{\beta}^{-1}} &= \prod_{i=1}^m \inverparbr{\alpha^{\reflbr{\word{\beta}} \reflbr{\beta_1' \cdots \beta_{i-1}'}}}{\beta_i'} = \prod_{i=1}^m \inverparbr{\alpha^{\reflbr{\beta_1 \cdots \beta_{m+1-i}}}}{-\beta_{m+1-i}}.
	\end{align*}
	Performing the index shift $ j \defl m+1-i $, we infer that
	\begin{align*}
		\inverparbr{\alpha^{\reflbr{\word{\beta}}}}{\word{\beta}^{-1}} &= \prod_{j=1}^m \inverparbr{\alpha^{\reflbr{\beta_1 \cdots \beta_j}}}{-\beta_j}.
	\end{align*}
	By \cref{param:inverpar-minus-cancel}, we have $ \inverparbr{\alpha^{\reflbr{\beta_1 \cdots \beta_j}}}{-\beta_j} = \inverparbr{\alpha^{\reflbr{\beta_1 \cdots \beta_{j-1}}}}{\beta_j}^{-1} $ for all $ j \in \Set{1, \ldots, m} $, so it follows that $ \inverparbr{\alpha^{\reflbr{\word{\beta}}}}{\word{\beta}^{-1}} = \inverparbr{\alpha}{\word{\beta}}^{-1} $.
\end{proof}

\begin{lemma}\label{parmap:stab-conj}
	Let $ \inverparsym $ be a parity map with values in some abelian group $ \twistgroup $, let $ \alpha $ be a root and let $ \word{\delta}, \word{\rho} $ be words over $ \rootbase \union (-\rootbase) $ such that $ \alpha^{\reflbr{\word{\delta}}} = \alpha $. Then $ \inverpar{\alpha^{\reflbr{\word{\rho}}}}{\word{\rho}^{-1} \word{\delta} \word{\rho}} = \inverpar{\alpha}{\word{\delta}} $.
\end{lemma}
\begin{proof}
	Write $ \word{\rho} = \tup{\rho}{m} $. If $ m=0 $, the assertion is trivial, so assume that $ m>0 $. Put $ \word{\rho}' \defl (\rho_2, \ldots, \rho_m) $, $ \word{\delta}' \defl (-\rho_1, \word{\delta}, \rho_1) $ and $ \alpha' \defl \alpha^{\reflbr{\rho_1}} $. Then
	\begin{align*}
		\inverparbr{\alpha^{\reflbr{\word{\rho}}}}{\word{\rho}^{-1} \word{\delta} \word{\rho}} &= \inverparsym\brackets[\big]{{(\alpha')^{\reflbr{\word{\rho}'}}},{(\word{\rho}')^{-1} \word{\delta}' \word{\rho}'}}
	\end{align*}
	where $ (\alpha')^{\reflbr{\word{\delta}'}} = \alpha^{\reflbr{\rho_1}} = \alpha' $. Hence induction on $ m $ yields $ \inverparbr{\alpha^{\reflbr{\word{\rho}}}}{\word{\rho}^{-1} \word{\delta} \word{\rho}} = \inverparbr{\alpha'}{\word{\delta}'} $. Using the definition of the extended parity map, we infer that
	\begin{align*}
		\inverparbr{\alpha^{\reflbr{\word{\rho}}}}{\word{\rho}^{-1} \word{\delta} \word{\rho}} &= \inverparsym\brackets[\big]{{\alpha^{\reflbr{\rho_1}}},{(-\rho_1, \word{\delta}, \rho_1)}} = \inverparbr{\alpha^{\reflbr{\rho_1}}}{-\rho_1} \inverparbr{\alpha}{\word{\delta}} \inverparbr{\alpha^{\reflbr{\word{\delta}}}}{\rho_1}.
	\end{align*}
	Since $ \inverparbr{\alpha^{\reflbr{\rho_1}}}{-\rho_1} = \inverparbr{\alpha}{\rho_1}^{-1} $ by \cref{param:inverpar-minus-cancel} and since $ \alpha^{\reflbr{\word{\delta}}} = \alpha $, the assertion follows.
\end{proof}

We now define some properties of parity maps that will be needed throughout this chapter.

\begin{note}
	Every property $ X $ that we define in this chapter will be of the form \enquote{$ \alpha $-$ X $ holds for all $ \alpha \in \roots $} where $ \alpha $-$ X $ is another property which is defined for every root $ \alpha $. However, we will only define the properties $ \alpha $-$ X $ for those $ X $ for which this terminology will be needed later, and simply define $ X $ directly in all other cases.
\end{note}

\begin{definition}\label{parmap:prop-def}
	Let $ A $ be an abelian group and let $ \inverparsym $ be a $ \rootbase $-parity map with values in $ A $.
	\begin{defenumerate}
		\item We say that $ \inverparsym $ is \defemph*{trivial}\index{parity map!trivial} if $ \inverpar{\alpha}{\delta} = 1 $ for all $ \alpha \in \roots $ and $ \delta \in \rootbase $.
		\item For any root $ \gamma $, we say that $ \inverparsym $ is \defemph*{$ \gamma $-braid-invariant}\index{parity map!braid-invariant} if for all $ \alpha, \beta \in \rootbase $, we have $ \inverparbr{\gamma}{\braidword{m}(\alpha, \beta)} = \inverparbr{\gamma}{\braidword{m}(\beta, \alpha)} $ where $ m $ denotes the order of $ \refl{\alpha} \refl{\beta} $ in the Weyl group and $ \braidword{m} $ is as in \cref{braidword-def}. We say that $ \inverparsym $ is \defemph*{braid-invariant} if it is $ \gamma $-braid-invariant for all $ \gamma \in \roots $.
		
		\item We say that $ \inverparsym $ is \defemph*{square-invariant}\index{parity map!square-invariant} if for all $ \alpha \in \roots, \delta \in \rootbase $, we have $ \inverpar{\alpha}{\delta \delta} = 1_A $.
		
		\item We say that $ \inverparsym $ is \defemph*{Weyl-invariant}\index{parity map!Weyl-invariant} if it is braid- and square-invariant.
		
		\item \label{parmap:prop-def:square-formula}For any root $ \alpha $, we say that \defemph*{$ \inverparsym $ satisfies the square formula for $ \alpha $}\index{square formula!for parity maps} if $ \roots $ is crystallographic, $ \twistgroup = \compactSet{\pm 1}^p $ for some $ p \in \Npos $ and $ \inverpar{\alpha}{\delta \delta} = (-1, 1, \ldots, 1)^{\cartanint{\alpha}{\delta}} $ for all $ \delta \in \rootbase $ where $ \cartanint{\alpha}{\delta} $ is the Cartan integer from \cref{rootsys:cartannum-def}. Further, we say that \defemph*{$ \inverparsym $ satisfies the square formula} if it satisfies the square formula for all $ \alpha \in \roots $.
		
		\item For any root $ \alpha \in \roots $, we say that $ \inverparsym $ is \defemph*{$ \alpha $-adjacency-trivial}\index{parity map!adjacency-trivial} if for any $ \beta \in \roots $ such that $ \alpha $ is adjacent to $ \beta $ and $ -\beta $ and for any $ \rootbase $-expression $ \word{\rho} $ of $ \beta $, we have $ \inverpar{\alpha}{\word{\rho}} = 1_\twistgroup $. Further, we say that $ \inverparsym $ is \defemph*{adjacency-trivial} if it is $ \alpha $-adjacency-trivial for all roots $ \alpha $.
	\end{defenumerate}
\end{definition}

\begin{remark}
	Let $ \alpha, \beta $ be non-proportional roots. Then $ \alpha $ is adjacent to $ \beta $ and $ -\beta $ if and only if the parabolic subsystem spanned by $ \alpha $ and $ \beta $ is of type $ A_1 \times A_1 $. Now assume that $ \roots $ is crystallographic. If $ \alpha $ is crystallographically adjacent to $ \beta $ and $ -\beta $, then the subsystem spanned by these roots need not be of type $ A_1 \times A_1 $, but an inspection of the possible rank-2 subsystems $ A_1 \times A_1 $, $ A_2 $, $ B_2 = C_2 $, $ BC_2 $ and $ G_2 $ shows that we still have $ \alpha \cdot \beta = 0 $.
\end{remark}

\begin{note}
	Braid-invariance and square-invariance are the most essential properties in \cref{parmap:prop-def}. The square formula (which we have already seen, in a different form, for Chevalley groups in \cref{chev:square-formula-note}) is merely a useful tool for proving square compatibility for groups which satisfy the square formula for Weyl elements, see \cref{param:square-comp-def,param:square-formula-rgg-def}. In a similar way, adjacency-triviality is merely a tool for proving stabiliser compatibility in certain cases, see \cref{param:stabcomp-def,param:stabcomp-crit-ortho,param:adj-implies-stab}. Both the square formula and adjacency-triviality have the useful property that they can be phrased purely in terms of parity maps, without reference to a specific root graded group. We will elaborate on this in \cref{param:combi-grp-separation-note}.
\end{note}

\begin{remark}[Weyl-invariant parity maps]\label{param:Weyl-invar}
	Let $ \inverparsym $ be a square-invariant parity map with values in some abelian group $ \twistgroup $. Then for all $ \alpha \in \roots $ and $ \delta \in \rootbase $, we have
	\begin{align*}
		\inverpar{\alpha}{-\delta} &= \inverpar{\alpha}{-\delta} \inverpar{\alpha^{\reflbr{-\delta}}}{\delta \delta} = \inverpar{\alpha}{(-\delta, \delta, \delta)}.
	\end{align*}
	Applying \cref{param:inverpar-minus-cancel}, we infer that $ \inverpar{\alpha}{-\delta} = \inverpar{\alpha}{\delta} $ for all $ \alpha \in \roots $ and $ \delta \in \rootbase $. Now assume that $ \inverparsym $ is, in addition, braid-invariant (so that it is Weyl-invariant). By the solution of the word problem in Coxeter groups (\cref{word-problem-solution}), we then have that $ \inverpar{\alpha}{\word{\beta}} $ depends only on $ \alpha $ and $ \reflbr{\word{\beta}} $ but not on the representation $ \word{\beta} $ of $ \reflbr{\word{\beta}} $, which justifies the name \enquote{Weyl-invariant}. Thus in this case, we can and will regard $ \inverparsym $ as a map defined on $ \roots \times \Weyl(\roots) $.
\end{remark}

We now turn to the somewhat technical notion of transporter sets.

\begin{definition}[Transporter sets]\label{param:comp-def}
	Let $ \twistgroup $ be an abelian group and let $ \inverparsym $ be a $ \rootbase $-parity map with values in $ A $. For any roots $ \alpha, \beta $ which lie in the same orbit under the Weyl group, the set
	\[ \twistgroup_{\inverparsym, \alpha \rightarrow \beta} \defl \Set{\inverpar{\alpha}{\word{\delta}} \given \word{\delta} \text{ word over } \rootbase \union (-\rootbase) \text{ with } \alpha^{\reflbr{\word{\delta}}} = \beta} \subs \twistgroup \]
	is called the \defemph*{transporter set for $ (\alpha, \beta) $ with respect to $ \inverparsym $}\index{transporter set}. If $ \inverparsym $ is clear from the context, we will usually write $ \parmoveset{\twistgroup}{\alpha}{\beta} $ for $ \twistgroup_{\inverparsym, \alpha \rightarrow \beta} $.
\end{definition}

\begin{remark}\label{param:transport-mult}
	Let $ \inverparsym $ be a parity map with values in some abelian group $ \twistgroup $. Then for all roots $ \alpha $, $ \beta $, $ \gamma $ which lie in the same orbit as $ \alpha $ under the Weyl group, we have $ \parmoveset{\twistgroup}{\alpha}{\beta} \parmoveset{\twistgroup}{\beta}{\gamma} \subs \parmoveset{\twistgroup}{\alpha}{\gamma} $, $ \parmoveset{\twistgroup}{\alpha}{\beta}^{-1} = \parmoveset{\twistgroup}{\beta}{\alpha} $ and $ 1_\twistgroup \in \parmoveset{\twistgroup}{\alpha}{\alpha} $. In particular, $ \parmoveset{\twistgroup}{\alpha}{\beta} $ is not necessarily a subgroup of $ \twistgroup $, but $ \parmoveset{\twistgroup}{\alpha}{\alpha} $ always is.
\end{remark}

In practice, we are only interested in parity maps for which all transporter sets within any orbit are the same (though different orbits may have different transporter sets). The following \cref{parmap:transport-invar-char} provides some (seemingly weaker) criteria from which this property follows. Here we always denote by $ \hat{\alpha} $, $ \hat{\beta} $ roots which are prefixed with an existential quantifier and by $ \alpha $, $ \beta $ roots which are prefixed with a universal quantifier.

\begin{lemma}\label{parmap:transport-invar-char}
	Let $ \inverparsym $ be a parity map with values in some abelian group $ \twistgroup $ and let $ O $ be an orbit of $ \roots $ under the Weyl group. Then the following properties are equivalent:
	\begin{stenumerate}
		\item \label{parmap:transport-invar-char:id}There exists $ \hat{\alpha} \in O $ such that for all $ \beta \in O $, we have $ 1_\twistgroup \in \parmoveset{\twistgroup}{\hat{\alpha}}{\beta} $.
		
		\item \label{parmap:transport-invar-char:id2}There exists $ \hat{\alpha} \in O $ such that for all $ \beta \in O $, we have $ 1_\twistgroup \in \parmoveset{\twistgroup}{\beta}{\hat{\alpha}} $.
		
		\item \label{parmap:transport-invar-char:weak-bound}There exist $ \hat{\alpha}, \hat{\beta} \in O $ such that $ \parmoveset{\twistgroup}{\hat{\alpha}}{\hat{\beta}} $ is a subgroup of $ \twistgroup $ which contains $ \parmoveset{\twistgroup}{\hat{\beta}}{\beta} $ for all $ \beta \in O $.
		
		\item \label{parmap:transport-invar-char:bound}There exist $ \hat{\alpha}, \hat{\beta} \in O $ such that $ \parmoveset{\twistgroup}{\hat{\alpha}}{\hat{\beta}} $ is a subgroup of $ \twistgroup $ which contains $ \parmoveset{\twistgroup}{\alpha}{\beta} $ for all $ \alpha, \beta \in O $.
		
		\item \label{parmap:transport-invar-char:same}For all $ \alpha, \beta, \gamma, \delta \in O $, we have $ \parmoveset{\twistgroup}{\alpha}{\beta} = \parmoveset{\twistgroup}{\gamma}{\delta} $.
	\end{stenumerate}
\end{lemma}
\begin{proof}
	The implications \enquote{\itemref{parmap:transport-invar-char:same} $ \implies $ \itemref{parmap:transport-invar-char:bound} $ \implies $ \itemref{parmap:transport-invar-char:weak-bound}} and \enquote{\itemref{parmap:transport-invar-char:id} $ \Longleftrightarrow $ \itemref{parmap:transport-invar-char:id2}} are clear. Now assume that~\itemref{parmap:transport-invar-char:weak-bound} holds, and choose $ \hat{\alpha} $, $ \hat{\beta} $ as in the assertion. We want to prove~\itemref{parmap:transport-invar-char:id}. For all $ \beta \in O $, we have
	\[ \parmoveset{\twistgroup}{\hat{\alpha}}{\hat{\beta}} \parmoveset{\twistgroup}{\hat{\beta}}{\beta} \subs \parmoveset{\twistgroup}{\hat{\alpha}}{\beta} \]
	by \cref{param:transport-mult}. By assumption, $ \parmoveset{\twistgroup}{\hat{\alpha}}{\hat{\beta}} $ contains $ \parmoveset{\twistgroup}{\hat{\beta}}{\beta} $. Since it is a subgroup of $ \twistgroup $, it also contains $ \parmoveset{\twistgroup}{\hat{\beta}}{\beta}^{-1} $. Thus the inclusion above implies that $ 1_\twistgroup $ is contained in $ \parmoveset{\twistgroup}{\hat{\alpha}}{\beta} $ for all $ \beta \in O $, which proves~\itemref{parmap:transport-invar-char:id}.
	
	Now assume that~\itemref{parmap:transport-invar-char:id} holds. We want to prove~\itemref{parmap:transport-invar-char:same}. Let $ \hat{\alpha} $ be as in~\itemref{parmap:transport-invar-char:id}. Note that property~\itemref{parmap:transport-invar-char:id2} is also satisfied for $ \hat{\alpha} $, too, because $ \parmoveset{\twistgroup}{\beta}{\hat{\alpha}} = \parmoveset{\twistgroup}{\hat{\alpha}}{\beta}^{-1} $. Further, for all $ \beta \in O $, it follows from \cref{param:transport-mult} that
	\[ \parmoveset{\twistgroup}{\hat{\alpha}}{\hat{\alpha}} \parmoveset{\twistgroup}{\hat{\alpha}}{\beta} \subs \parmoveset{\twistgroup}{\hat{\alpha}}{\beta} \midand \parmoveset{\twistgroup}{\hat{\alpha}}{\beta} \parmoveset{\twistgroup}{\beta}{\hat{\alpha}} \subs \parmoveset{\twistgroup}{\hat{\alpha}}{\hat{\alpha}}. \]
	Since $ \parmoveset{\twistgroup}{\hat{\alpha}}{\beta} $ and $ \parmoveset{\twistgroup}{\beta}{\hat{\alpha}} $ both contain $ 1_\twistgroup $, we infer that $ \parmoveset{\twistgroup}{\hat{\alpha}}{\hat{\alpha}} \subs \parmoveset{\twistgroup}{\hat{\alpha}}{\beta} \subs \parmoveset{\twistgroup}{\hat{\alpha}}{\hat{\alpha}} $, so $ \parmoveset{\twistgroup}{\hat{\alpha}}{\hat{\alpha}} = \parmoveset{\twistgroup}{\hat{\alpha}}{\beta} $ for all $ \beta \in O $. In particular, $ \parmoveset{\twistgroup}{\hat{\alpha}}{\beta} $ is a group, so $ \parmoveset{\twistgroup}{\beta}{\hat{\alpha}} = \parmoveset{\twistgroup}{\hat{\alpha}}{\beta}^{-1} = \parmoveset{\twistgroup}{\hat{\alpha}}{\beta} = \parmoveset{\twistgroup}{\hat{\alpha}}{\hat{\alpha}} $ for all $ \beta \in O $ as well. Now for arbitrary $ \beta, \gamma $, we have
	\[ \parmoveset{\twistgroup}{\beta}{\hat{\alpha}} \parmoveset{\twistgroup}{\hat{\alpha}}{\gamma} \subs \parmoveset{\twistgroup}{\beta}{\gamma} \midand \parmoveset{\twistgroup}{\hat{\alpha}}{\beta} \parmoveset{\twistgroup}{\beta}{\gamma} \subs \parmoveset{\twistgroup}{\hat{\alpha}}{\gamma} \]
	where $ 1_\twistgroup $ lies in $ \parmoveset{\twistgroup}{\beta}{\hat{\alpha}} $ and in $ \parmoveset{\twistgroup}{\hat{\alpha}}{\beta} $. By similar arguments as above, it follows that $ \parmoveset{\twistgroup}{\hat{\alpha}}{\gamma} = \parmoveset{\twistgroup}{\beta}{\gamma} $ for all $ \beta, \gamma \in O $. Since we already know that $ \parmoveset{\twistgroup}{\hat{\alpha}}{\gamma} = \parmoveset{\twistgroup}{\hat{\alpha}}{\hat{\alpha}} $, we conclude that all the groups $ (\parmoveset{\twistgroup}{\beta}{\gamma})_{\beta, \gamma \in O} $ are equal to $ \parmoveset{\twistgroup}{\hat{\alpha}}{\hat{\alpha}} $. In particular, they are pairwise equal. This finishes the proof of~\itemref{parmap:transport-invar-char:same}.
\end{proof}

\begin{remark}\label{parmap:transport-invar-crit}
	Note that property~\thmitemref{parmap:transport-invar-char}{parmap:transport-invar-char:bound} is automatically satisfied if $ \parmoveset{\twistgroup}{\hat{\alpha}}{\hat{\beta}} = \twistgroup $ for some $ \hat{\alpha}, \hat{\beta} \in O $. In this case, we have $ \parmoveset{\twistgroup}{\alpha}{\beta} = \twistgroup $ for all $ \alpha, \beta \in O $.
\end{remark}

\begin{definition}[Transporter properties]\label{param:transporter-prop-def}
	Let $ \twistgroup $ be an abelian group and let $ \inverparsym $ be a $ \rootbase $-parity map with values in $ A $.
	\begin{defenumerate}
		\item We say that $ \inverparsym $ is \defemph*{complete}\index{parity map!complete} if $ \parmoveset{\twistgroup}{\alpha}{\alpha} = \twistgroup $ for all roots $ \alpha $.
		
		\item We say that $ \inverparsym $ is \defemph*{semi-complete}\index{parity map!semi-complete} if for all roots $ \alpha $, the group $ \parmoveset{\twistgroup}{\alpha}{\alpha} $ has a group-theoretic complement in $ \twistgroup $. That is, there exists a subgroup $ C_\alpha $ of $ \twistgroup $ such that $ \parmoveset{\twistgroup}{\alpha}{\alpha} \intersect C_\alpha = \compactSet{1} $ and $ \parmoveset{\twistgroup}{\alpha}{\alpha} C_\alpha = \twistgroup $.
		
		\item For any orbit $ O $ of $ \roots $ under the Weyl group, we say that $ \inverparsym $ is \defemph*{transporter-invariant on $ O $}\index{transporter-invariance} if it satisfies the equivalent conditions in \cref{parmap:transport-invar-char} for $ O $. Further, we say that $ \inverparsym $ is \defemph*{transporter-invariant} if it is transporter-invariant on all orbits of $ \roots $.
	\end{defenumerate}
\end{definition}

\begin{remark}\label{param:transp-invar-1}
	Let $ \inverparsym $ be a transporter-invariant $ \rootbase $-parity map with values in some abelian group $ A $ and let $ \alpha, \beta $ be roots which lie in the same orbit. Since $ \parmoveset{\twistgroup}{\alpha}{\beta} = \parmoveset{\twistgroup}{\alpha}{\alpha} $, it follows from \cref{param:transport-mult} that $ 1_\twistgroup $ is contained in $ \parmoveset{\twistgroup}{\alpha}{\beta} $.
\end{remark}

\begin{note}\label{param:transporter-prop-notes}
	Clearly, any complete parity map is also transporter-invariant and semi-complete. In practice, we would like our parity maps to be complete, but there will be situations in which this is not the case. We delay an explanation of this fact until \cref{param:twist-for-orbit-note}. There we will also sketch an alternative approach which would allow us to only consider complete parity maps, but which introduces other difficulties.
	
	The seemingly strange notion of semi-completeness will be needed in \cref{param:invo-def}. In all practical cases, the group $ \twistgroup $ will be of the form $ \twistgroup = \compactSet{\pm 1}^p $ for some $ p \in \Set{1,2,3} $ and every transporter set in $ \twistgroup $ will simply be the product of some of the $ p $ components of $ \twistgroup $. In this situation, we can simply choose the complement of any transporter set to be the product of the remaining components.
\end{note}

We now consider products of parity maps.

\begin{notation}
	Let $ \inverparsym, \invoparsym $ be two parity maps with values in some abelian groups $ \twistgroup, \invogroup $, respectively. Then we denote by $ \inverparsym \times \invoparsym $ the parity map
	\[ \map{}{\roots \times \rootbase}{\twistgroup \times \invogroup}{(\alpha, \beta)}{\brackets{\inverpar{\alpha}{\beta}, \invopar{\alpha}{\beta}}}. \]
\end{notation}

\begin{remark}\label{param:parmap-split}
	If $ \twistgroup, \invogroup $ are two abelian groups and $ \totalparsym $ is a parity map with values in $ \twistgroup \times \invogroup $, then there exist unique parity maps $ \inverparsym, \invoparsym $ with values in $ \twistgroup, \invogroup $, respectively, such that $ \totalparsym = \inverparsym \times \invoparsym $.
\end{remark}

\begin{definition}[Independent parity maps]\label{param:parmap-indep-def}
	Let $ \twistgroup, \invogroup $ be two abelian groups and let $ \inverparsym, \invoparsym $ be parity maps with values in $ \twistgroup, \invogroup $, respectively. The parity maps $ \inverparsym $ and $ \invoparsym $ are called \defemph*{independent}\index{parity map!independent} if $ \genparmoveset{(\twistgroup \times \invogroup)}{\inverparsym \times \invoparsym}{\alpha}{\alpha} = \genparmoveset{\twistgroup}{\inverparsym}{\alpha}{\alpha} \times \genparmoveset{\invogroup}{\invoparsym}{\alpha}{\alpha} $ for all $ \alpha \in \roots $.
\end{definition}

\begin{remark}\label{param:transporter-proj}
	In the notation of \cref{param:parmap-indep-def}, we always have
	\[ \genparmoveset{(\twistgroup \times \invogroup)}{\inverparsym \times \invoparsym}{\alpha}{\beta} \subs \genparmoveset{\twistgroup}{\inverparsym}{\alpha}{\beta} \times \genparmoveset{\invogroup}{\invoparsym}{\alpha}{\beta}. \]
	Hence $ \inverparsym $ and $ \invoparsym $ are independent if and only if the reverse inclusion holds.
\end{remark}

\begin{lemma}\label{param:prod-eq-lem}
	Let $ \inverparsym, \invoparsym $ be two parity maps with values in respective abelian groups $ \twistgroup $ and $ \invogroup $, let $ \alpha $ be any root and let $ O $ be any orbit in $ \roots $.
	\begin{lemenumerate}
		\item \label{param:prod-eq-lem:eq}For all the following properties, it is true that $ \inverparsym \times \invoparsym $ has this property if and only if  $ \inverparsym $ and $ \invoparsym $ both have this property: square-invariance, ($ \alpha $-) braid-invariance, Weyl-invariance, (crystallographic) ($ \alpha $-) adjacency-triviality.
		
		\item \label{param:prod-eq-lem:comp}For all the following properties, it is true that $ \inverparsym $ and $ \invoparsym $ both have this property if $ \inverparsym \times \invoparsym $ has this property: completeness, transporter-invariance (on $ O $).
	\end{lemenumerate}
\end{lemma}
\begin{proof}
	The assertions of~\itemref{param:prod-eq-lem:eq} are easy to verify. The assertions of \itemref{param:prod-eq-lem:comp} follow from \cref{param:transporter-proj}, using criterion~\thmitemref{parmap:transport-invar-char}{parmap:transport-invar-char:id} to prove transporter-invariance.
\end{proof}


\section{Parametrisations and Twisting Structures}

\label{sec:param:twist}

\begin{secnotation}\label{param:secnot:param-def}
	In this section, $ \roots $ is a root system and $ G $ is a group with a $ \roots $-pregrading $ (\rootgr{\alpha})_{\alpha \in \roots} $ such that $ \rootgr{\lambda \alpha} \subs \rootgr{\alpha} $ for all $ \alpha \in \roots $ and all $ \lambda \in \IR_{>1} $ for which $ \lambda \alpha $ is a root. Further, we fix a rescaled root base $ \rootbase $ of $ \roots $ and we assume that there exists a $ \rootbase $-system $ (w_\delta)_{\delta \in \rootbase} $ of Weyl elements in $ G $, which we also fix.
\end{secnotation}

From now on, we study the interplay of the pair $ (G , (\rootgr{\alpha})_{\alpha \in \roots}) $ with parity maps. We begin with the notion of parametrisations.

\begin{definition}[Parameter system]\label{param:parsys-def}
	Write $ \IndivOrb(\roots) = \tup{O}{k} $ (see \cref{rootsys:orb-def}). A \defemph*{parameter system (of type $ \roots $)}\index{parameter system} is a tuple $ \calP = (\twistgroup, \listing{M}{k}) $ consisting of an abelian group $ \twistgroup $ and groups $ \listing{M}{k} $ such that $ \twistgroup $ acts on each of the underlying sets $ \listing{M}{k} $. The group $ \twistgroup $ is called the \defemph*{twisting group of $ \calP $}\index{twisting group!of a parameter system} and the actions of $ \twistgroup $ on $ \listing{M}{k} $ are called the \defemph*{twisting actions}\index{twisting action}. Given a parameter system $ \calP = (\twistgroup, \listing{M}{k}) $ and a root $ \alpha \in \roots $, we denote by $ M_\alpha $ the group $ M_{i(\alpha)} $ where $ i(\alpha) \in \numint{1}{k} $ is the unique index such that $ \alpha \in O_{i(\alpha)} $.
\end{definition}

\begin{note}
	The root system $ \roots $ appears in the definition of parameter systems of type $ \roots $ only insofar that the number $ k $ of groups $ \listing{M}{k} $ is the number of indivisible orbits in $ \roots $. However, the definition of the groups $ (M_\alpha)_{\alpha \in \roots} $ clearly depends on $ \roots $. For this reason, we use the terminology of \enquote{parameter systems of type $ \roots $}.
	
	Further, observe that we do not require $ \twistgroup $ to act on $ \listing{M}{k} $ by group automorphisms. Thus we need not have $ a.(xy) = (a.x)(a.y) $ for all $ a \in \twistgroup $ and $ x,y $ in one of the parametrising groups. However, this additional condition will often be satisfied.
\end{note}

\begin{definition}[Parametrisation]\label{param:param-def}
	Write $ \IndivOrb(\roots) = \tup{O}{k} $, let $ \calP = (\twistgroup, \listing{M}{k}) $ be a parameter system of type $ \roots $ and let $ \inverparsym $ a $ \rootbase $-parity map with values in $ \twistgroup $. A \defemph*{parametrisation of $ G $ by $ \calP $ with respect to $ \inverparsym $ and $ (w_{\delta})_{\delta \in \rootbase} $}\index{parametrisation of a root graded group} is a family
	\[ \brackets[\big]{\map{\risom{\alpha}}{M_\alpha}{\rootgr{\alpha}}{}{}}_{\alpha \in \indivset{\roots}} \]
	such that for all roots $ \alpha \in \indivset{\roots} $, $ \delta \in \rootbase $ and all $ x \in M_\alpha $, we have
	\[ \risom{\alpha}(x)^{w_{\delta}} = \risom{\refl{\delta}(\alpha)}(\inverpar{\alpha}{\delta}.x). \]
	If $ \beta $ is a root which is not indivisible and $ \beta' $ is the unique indivisible root in $ \IR_{>0} \beta $, then $ \rootgr{\beta} \subs \rootgr{\beta'} $ by \cref{param:secnot:param-def}, so we can define $ M_\beta \defl \risom{\beta'}^{-1}(\rootgr{\beta'}) $ and $ \map{\risom{\beta} \defl \restrict{\risom{\beta'}}{M_\beta}}{M_\beta}{\rootgr{\beta}}{}{} $. The family $ (\risom{\alpha})_{\alpha \in \roots} $ will also be called a parametrisation of $ G $, and the maps in this family are called the \defemph*{root isomorphisms of $ G $}\index{root isomorphism}. We will also say that \defemph*{$ G $ is parametrised by $ \calP $ with respect to $ \inverparsym $ and $ (w_{\delta})_{\delta \in \rootbase} $ and root isomorphisms $ (\risom{\alpha})_{\alpha \in \roots} $}.
\end{definition}

\begin{remark}\label{param:conj-formula-word}
	Let everything be as in \cref{param:param-def}. Let $ \alpha \in \roots $ and $ x \in M_\alpha $. It follows from the equation
	\[ \risom{\alpha}(x)^{w_{\delta}} = \risom{\refl{\delta}(\alpha)}(\inverpar{\alpha}{\delta}.x) \]
	that
	\[ \risom{\alpha}(x) = \risom{\refl{\delta}(\alpha)}(\inverpar{\alpha}{\delta}.x)^{w_\delta^{-1}}. \]
	Replacing $ \alpha $ by $ \refl{\delta}(\alpha) $ and $ x $ by $ \inverpar{\refl{\delta}(\alpha)}{\delta}^{-1}.x $, we infer that
	\[ \risom{\refl{\delta}(\alpha)}(\inverpar{\refl{\delta}(\alpha)}{\delta}^{-1}.x) = \risom{\alpha}(x)^{w_\delta^{-1}}. \]
	Thus
	\[ \risom{\alpha}(x)^{w_{-\delta}} = \risom{\alpha}(x)^{w_\delta^{-1}} = \risom{\refl{\delta}(\alpha)}\brackets[\big]{\inverpar{\refl{\alpha}(\alpha)}{\delta}^{-1}.x} = \risom{\refl{\delta}(\alpha)}\brackets[\big]{\inverpar{\alpha}{-\delta}.x}. \]
	This implies that for all words $ \word{\delta} $ over $ \rootbase \union (-\rootbase) $, we have
	\[ \risom{\alpha}(x)^{w_{\word{\delta}}} = \risom{\alpha^{\reflbr{\word{\delta}}}}(\inverpar{\alpha}{\word{\delta}}.x). \]
\end{remark}

We now turn to the notion of faithfulness for parameter systems.

\begin{definition}[Faithfulness with respect to parity maps]\label{param:faithful-act-def}
	Let $ \twistgroup $ be any group acting on a set $ M $ and let $ \inverparsym $ be a $ \rootbase $-parity map with values in $ \twistgroup $. For any root $ \alpha $, we say that the action of $ \twistgroup $ on $ M $ is \defemph*{$ \alpha $-faithful with respect to $ \inverparsym $}, or simply \defemph*{$ \inverpar{\alpha}{} $-faithful}\index{eta-faithful@$ \inverpar{\alpha}{} $-faithful} if for all words $ \word{\delta} $, $ \word{\rho} $ over $ \rootbase \union (-\rootbase) $ such that $ \inverpar{\alpha}{\word{\delta}} $ and $ \inverpar{\alpha}{\word{\rho}} $ act identically on $ M $, we have $ \inverpar{\alpha}{\word{\delta}} = \inverpar{\alpha}{\word{\rho}} $.
\end{definition}

\begin{remark}
	Equivalently, $ \twistgroup $ acts $ \alpha $-faithfully on $ M $ with respect to $ \inverparsym $ if and only if the subset $ \bigunion_{\beta \in O} \genparmoveset{\twistgroup}{\inverparsym}{\alpha}{\beta} $ acts faithfully on $ M $ where $ O $ denotes the orbit of $ \alpha $ in $ \roots $. If $ \inverparsym $ is transporter-invariant (which will always be the case in practice), then $ \bigunion_{\beta \in O} \genparmoveset{\twistgroup}{\inverparsym}{\alpha}{\beta} = \genparmoveset{\twistgroup}{\inverparsym}{\alpha}{\alpha} $.
\end{remark}

\begin{definition}[Faithfulness for parameter systems]\label{param:parsys-faithful}
	Let $ \calP = (\twistgroup, \listing{M}{k}) $ be a parameter system of type $ \roots $. We say that $ \calP $ is \defemph*{faithful}\index{parameter system!faithful} if $ \twistgroup $ acts faithfully on all groups $ \listing{M}{k} $. For any $ \rootbase $-parity map $ \inverparsym $ with values in $ \twistgroup $, we say that $ \calP $ is \defemph*{$ \inverparsym $-faithful}\index{parameter system!faithful!eta-@$ \inverparsym $-} if for all roots $ \alpha $, the group $ \twistgroup $ acts $ \alpha $-faithfully with respect to $ \inverparsym $ on all groups $ \listing{M}{k} $.
\end{definition}

\begin{note}\label{param:faithful-note}
	Let $ \calP = (\twistgroup, \listing{M}{k}) $ be a parameter system of type $ \roots $ and let $ \inverparsym $ a $ \rootbase $-parity map with values in $ \twistgroup $. In all practical situations, the map
	\[ \map{\inverparsym}{\roots \times \frmon{\rootbase \union (-\rootbase)}}{\twistgroup}{}{} \]
	will be surjective because otherwise, we would simply shrink the co\-do\-main to the image of $ \inverparsym $. However, this does not imply that the map
	\begin{equation}\label{eq:param:faithful-note:inv-alpha}
		\map{\inverpar{\alpha}{}}{\frmon{\rootbase \union (-\rootbase)}}{\twistgroup}{}{}
	\end{equation}
	is surjective for all $ \alpha \in \roots $. This is illustrated by \cref{param:B-twist-ex}. As a consequence, faithfulness of all twisting actions of $ \twistgroup $ is a stronger property than $ \inverparsym $-faithfulness of $ \calP $. However, $ \inverparsym $-faithfulness of $ \calP $ is precisely the property that we need to show that a parity map which is read off from some $ \roots $-graded group automatically has several useful properties (see \cref{param:motiv:parmap-choice,param:param-parmap-has-properties}).
\end{note}

\begin{example}[see also \cref{quadmod:standard-param}]\label{param:B-twist-ex}
	Let $ \comring $ be a commutative associative ring and let $ \module $ be a $ \comring $-module. Assume that there exist a $ \comring $-quadratic form $ \map{q}{\module}{\comring}{}{} $ and an element $ v_0 \in \module $ with $ q(v_0) = 1_\comring $, and denote the reflection corresponding to $ v_0 $ by $ \map{\refl{v_0}}{\module}{\module}{}{} $. (We will give precise definitions of quadratic forms and reflections in \cref{quadmod:quadform-def,quadmod:refl-def}, but they are actually not relevant for the moment.) Put $ \twistgroup \defl \compactSet{\pm 1}^2 $ and define actions of $ \twistgroup $ on $ \comring $ and $ \module $ by
	\begin{align*}
		(\epsilon_1, \epsilon_2).a \defl \epsilon_1 a \midand (\epsilon_1, \epsilon_2).v \defl \begin{cases}
			\epsilon_1 v & \text{if } \epsilon_2 = 1, \\
			\epsilon_1 \refl{v_0}(v) & \text{if } \epsilon_2 = -1
		\end{cases}
	\end{align*}
	for all $ (\epsilon_1, \epsilon_2) \in \twistgroup $, all $ a \in \comring $ and all $ v \in \module $. In other words, the first component of $ \twistgroup $ acts on $ \comring $ and $ \module $ by additive inversion while the second component acts trivially on $ \comring $ and by $ \refl{v_0} $ on $ \module $. Thus the action of $ \twistgroup $ on $ \module $ is faithful in the generic case (that is, unless we specifically choose $ \module $ and $ v_0 $ to satisfy, for example, $ \refl{v_0} = -\id_\module $), but the action of $ \twistgroup $ on $ \comring $ is never faithful. There is no way to avoid this in our setup because we need two twisting actions on $ \module $ but only one twisting action on $ \comring $, so there must be a part of $ \twistgroup $ which acts trivially on $ \comring $.
	
	In \cref{sec:B-example}, we will construct a $ B_3 $-graded group $ H $ which is parametrised by the parameter system $ \calP \defl (\twistgroup, \module, \comring) $. Since the action of $ \twistgroup $ on $ \comring $ is not faithful, there is no unique parity map $ \inverparsym $ with respect to which $ H $ is parametrised: We can freely choose the second component of $ \inverpar{\alpha}{\delta} $ for all long roots $ \alpha $ and all $ \delta \in \rootbase $. However, we can choose a parity map $ \inverparsym $ by declaring that the second component should always be $ 1 $ when it is not uniquely determined. Then $ \calP $ is $ \inverparsym $-faithful for this choice of $ \inverparsym $.
\end{example}

\begin{note}\label{param:twist-for-orbit-note}
	The problems in \cref{param:B-twist-ex} could be avoided by introducing a separate twisting group $ \twistgroup_O $ and \enquote{partial parity map} $ \map{\inverparsym_O}{O \times \rootbase}{\twistgroup_O}{}{} $ for each orbit $ O $ of roots, that is, for each parametrising group in the parameter system. In \cref{param:B-twist-ex}, we would have $ A_L = \compactSet{\pm 1} $ for the orbit $ L $ of long roots and $ A_S = \compactSet{\pm 1}^2 $ for the orbit $ S $ of short roots. In this setup, all (partial) parity maps we deal with would be complete, so we would not have to bother with the weaker properties of semi-completeness and transporter-invariance. However, this approach has the disadvantage of making our notation more unwieldy by introducing more twisting groups and parity maps, which is why we will not follow it.
\end{note}

\begin{remark}[Twisted parametrisations]\label{param:parmap-twist}
	Let $ \calP = (\twistgroup, \listing{M}{k}) $ be a parameter system of type $ \roots $, let $ \inverparsym $ be a $ \rootbase $-parity map with values in $ \twistgroup $ and let $ (\risom{\gamma})_{\gamma \in \roots} $ be a parametrisation of $ G $ by $ \calP $ with respect to $ \inverparsym $ and $ (w_\delta)_{\delta \in \rootbase} $. Fix an arbitrary root $ \alpha $ and an arbitrary element $ a \in \twistgroup $. Assume that the action of $ a $ on $ M_\alpha $ is a group endomorphism (and thus an isomorphism). We define a twisted root isomorphism for $ \rootgr{\alpha} $ by
	\[ \map{\risom{\alpha}'}{M_\alpha}{\rootgr{\alpha}}{x}{\risom{\alpha}(a.x)}. \]
	The assumption on $ a $ assures that $ \risom{\alpha}' $ is indeed a homomorphism. Further, define $ \risom{\beta}' \defl \risom{\beta} $ for all roots $ \beta $ distinct from $ \alpha $. Then $ (\risom{\gamma}')_{\gamma \in \roots} $ is also a parametrisation of $ G $ by $ \calP $ with respect to $ (w_\delta)_{\delta \in \rootbase} $, but the corresponding parity map $ \inverparsym' $ is distinct from $ \inverparsym $. Namely, we have $ \inverpar{\alpha}{\delta}' = a \inverpar{\alpha}{\delta} $
	for all $ \delta \in \rootbase $ and
	$ \inverpar{\beta}{\delta}' = a^{-1} \inverpar{\beta}{\delta} $
	for all $ \beta \in \roots $ and $ \delta \in \rootbase $ such that $ \refl{\delta}(\beta) = \alpha $. We say that $ (\risom{\gamma}')_{\gamma \in \roots} $ is \defemph*{obtained from $ (\risom{\gamma})_{\gamma \in \roots} $ by twisting}\index{parametrisation!twisted}, and similarly for $ \inverparsym' $ and $ \inverparsym $.\index{parity map!twisted}
	
	More generally, we say that a parametrisation or a parity map is obtained from another one by twisting if it can be obtained by a finite number of twistings in the previous sense. It is clear that, once we have found one parametrisation and a corresponding parity map, we can construct a large number of different parametrisations and parity maps by twisting. In \cref{param:gpt-twist}, we will show that all parity maps which satisfy the necessary compatibility conditions for some root graded group $ G $ are the same up to twisting if $ G $ is sufficiently generic (which, in more technical terms, translates to some faithfulness assumptions).
\end{remark}

We now turn the the structures which are needed as a stepping stone in the construction of a parametrisation: twisting groups and (partial) twisting systems.

\begin{definition}[Twisting group]\label{param:twist-grp-def}
	A \defemph*{twisting group for $ (G, (w_\delta)_{\delta \in \rootbase}) $ with respect to $ \rootbase $}\index{twisting group!for a root graded group} is a tuple $ (\twistgroup, (\omega_\alpha)_{\alpha \in \roots}) $ with the following properties:
	\begin{stenumerate}
		\item $ \twistgroup $ is an abelian group.
		
		\item For each root $ \alpha $, the map $ \map{\omega_\alpha}{\twistgroup \times \rootgr{\alpha}}{\rootgr{\alpha}}{(a,g)}{a.g} $ is a group action of $ \twistgroup $ on the underlying set of $ \rootgr{\alpha} $, called the \defemph*{twisting action of $ \twistgroup $ on $ \rootgr{\alpha} $}\index{twisting action}.
		
		\item \label{param:twist-grp-def:comm}The twisting action on $ G $ commutes with conjugation by the fixed set of Weyl elements and their inverses: For all $ \alpha \in \roots $, $ \delta \in \rootbase \union (-\rootbase) $ and $ a \in \twistgroup $, we have $ (a.x_\alpha)^{w_\delta} = a.(x_\alpha^{w_\delta}) $ for all $ x_\alpha \in \rootgr{\alpha} $.
	\end{stenumerate}
	We will usually simply say that $ \twistgroup $ a twisting group for $ G $, leaving the twisting actions, the root base and the Weyl elements implicit.
\end{definition}

\begin{note}
	Given a parametrisation $ (\risom{\alpha})_{\alpha \in \roots} $ of $ G $ by a parameter system $ \calP = (\twistgroup, \listing{M}{K}) $ with respect to $ (w_\delta)_{\delta \in \rootbase} $, we can easily turn the twisting group $ \twistgroup $ of $ \calP $ into a twisting group for $ (G, (w_\delta)_{\delta \in \rootbase}) $ by putting $ a.\risom{\alpha}(x) \defl \risom{\alpha}(a.x) $ for all $ a \in \twistgroup $, all $ \alpha \in \roots $ and all $ x \in M_\alpha $. (We will give a proof of this fact in \cref{param:param-parmap-has-properties}.) The content of the parametrisation theorem is precisely that under certain conditions, a twisting group $ \twistgroup $ for $ G $ can be turned into the twisting group of a parameter system $ (\twistgroup, \listing{M}{k}) $ by which $ G $ is parametrised. Hence there is a close connection between \cref{param:parsys-def,param:twist-grp-def}, which justifies that we use the terminology of \enquote{twisting groups} in both situations.
\end{note}

The following twisting group will play an important role for every root system that we consider.

\begin{example}\label{param:pargroup-inv-example}
	Consider the multiplicative group $ \twistgroup \defl \compactSet{\pm 1} $ of order $ 2 $ and define that the non-trivial element of $ \twistgroup $ acts on all root groups by group inversion. That is, $ a.x_\alpha \defl x_\alpha^{-1} $ for all $ \alpha \in \roots $, $ x_\alpha \in \rootgr{\alpha} $ and the non-trivial element $ a \in \twistgroup \setminus \compactSet{1_\twistgroup} $. Then $ \twistgroup $ is a twisting group for $ G $. Note that for all $ \alpha \in \roots $, the twisting action on $ \rootgr{\alpha} $ is compatible with the multiplication on $ \rootgr{\alpha} $ (that is, $ \twistgroup $ acts by group automorphisms) if and only if $ \rootgr{\alpha} $ is abelian.
\end{example}

\begin{note}
	In every example, property~\thmitemref{param:twist-grp-def}{param:twist-grp-def:comm} will follow from the fact that the twisting action can be expressed \enquote{naturally} without making any choices. For example, the twisting group in \cref{param:pargroup-inv-example} acts by inversion, and it is clear that $ (x_\alpha^{w_\delta})^{-1} = (x_\alpha^{-1})^{w_\delta} $. It is somewhat surprising that this always works because in some situations, the twisting action will be defined as conjugation by a certain group element. A priori, this involves a choice, but it can be shown that the twisting action is independent of this choice.
	In fact, we will always see that property~\thmitemref{param:twist-grp-def}{param:twist-grp-def:comm} is not only satisfied for the fixed $ \rootbase $-system $ (w_\delta)_{\delta \in \rootbase} $ of Weyl elements but in fact for all Weyl elements. However, the weaker condition~\thmitemref{param:twist-grp-def}{param:twist-grp-def:comm} will be sufficient to prove the parametrisation theorem.
\end{note}

\begin{definition}
	Let $ \twistgroup $ be a twisting group for $ G $, let $ g \in G $ and let $ a \in \twistgroup $. For any root $ \alpha $, we say that \emph{$ g $ acts on $ \rootgr{\alpha} $ by $ a $} if $ g $ normalises $ \rootgr{\alpha} $ and $ x_\alpha^g = a.x_\alpha $ for all $ x_\alpha \in \rootgr{\alpha} $.
\end{definition}

The notion of (partial) twisting systems captures some of the additional assumptions that are required in the parametrisation theorem.

\begin{definition}[Partial twisting system]\label{param:partwist-def}
	A tuple $ (\twistgroup, \inverparsym, \invogroup, \invoparsym) $ is called a \defemph*{partial twisting system for $ (G, (w_\delta)_{\delta \in \rootbase}) $ with respect to $ \rootbase $}\index{twisting system!partial} if it satisfies the following properties:
	\begin{stenumerate}
		\item $ \twistgroup $ is a twisting group for $ (G, (w_\delta)_{\delta \in \rootbase}) $ and $ \inverparsym $ is a parity map with values in~$ \twistgroup $.
		
		\item $ \invogroup $ is an abelian group and $ \invoparsym $ is a parity map with values in $ \invogroup $.
		
		\item $ \inverparsym $ is braid-invariant.
		
		\item $ \invoparsym $ is Weyl-invariant and semi-complete.
		
		\item \label{param:partwist-def:stab-comp}$ \inverparsym \times \invoparsym $ is transporter-invariant and $ \inverparsym $, $ \invoparsym $ are independent in the sense of \cref{param:transporter-prop-def,param:parmap-indep-def}.
	\end{stenumerate}
	We will usually leave the root base $ \rootbase $ implicit.
\end{definition}

Recall from \thmitemcref{param:prod-eq-lem}{param:prod-eq-lem:comp} that, as a consequence of Axiom~\thmitemref{param:partwist-def}{param:partwist-def:stab-comp}, the parity maps $ \inverparsym $ and $ \invoparsym $ in a partial twisting system are transporter-invariant as well.

\begin{definition}[Twisting system]
	A \defemph*{twisting system for $ (G, (w_\delta)_{\delta \in \rootbase}) $ with respect to $ \rootbase $}\index{twisting system} is a tuple $ (\twistgroup, \inverparsym) $ where $ \twistgroup $ is a twisting group for $ (G, (w_\delta)_{\delta \in \rootbase}) $ and $ \inverparsym $ is a transporter-invariant braid-invariant parity map with values in $ \twistgroup $. In other words, it is a tuple $ (\twistgroup, \inverparsym) $ such that $ (\twistgroup, \inverparsym, \compactSet{1}, \map{}{}{}{(\alpha, \delta)}{1}) $ is a partial twisting system, which is then called the \defemph*{partial twisting system associated to $ (\twistgroup, \inverparsym) $}.
\end{definition}

\begin{note}
	The appearance of a second group $ \invogroup $ and of a second parity map $ \invoparsym $ in \cref{param:partwist-def} can be explained as follows. In our examination of $ \roots $-graded groups for specific choices of $ \roots $ in the later chapters, we will encounter two kinds of twistings of the root groups: Some which are relatively easy to define and some whose mere definition is more challenging. Our strategy is to simply ignore twistings of the second kind in the beginning, which is why the group $ \invogroup $ is, unlike $ \twistgroup $, not equipped with actions on the root groups. The first step in the proof of the parametrisation theorem is to construct an action of $ \invogroup $ on all root groups, thereby turning the partial twisting system $ (\twistgroup, \inverparsym, \invogroup, \invoparsym) $ into a twisting system $ (\twistgroup \times \invogroup, \inverparsym \times \invoparsym) $. In this way, the construction of the action of $ \invogroup $ on the root groups is done in a uniform way as part of the parametrisation theorem, so we do not have to do it \enquote{by hand} for each specific choice of $ \roots $.
	
	In practice, we will always have $ \twistgroup = \compactSet{\pm 1}^p $ for some $ p \in \Set{1,2} $ and $ \invogroup = \compactSet{1}^q $ for some $ \Set{0,1} $. The case $ q=0 $ appears precisely for those root systems in which pairs of orthogonal roots are automatically adjacent: in other words, for the simply-laced types $ A $, $ D $, $ E $ and for the types $ H_3 $ and $ H_4 $ (which are covered in \cite{RGG-H3}). The case $ q=1 $ appears for the types $ B $, $ (B)C $ and $ F_4 $.
\end{note}

\begin{example}[Parametrisations of Chevalley groups]\label{param:chev-ex}
	We briefly investigate what the notions defined in this section mean in the context of Chevalley groups. Assume that $ \roots $ is crystallographic and reduced. Let $ \ring $ be a commutative associative ring, let $ \chevstr = (\chevstr_{\alpha, \beta})_{\alpha, \beta \in \roots} $ be a family of Chevalley structure constants of type $ \roots $ (in the sense of \cref{chev:struc-family-def}) and let $ G $ be a Chevalley group of type $ \roots $ over $ \ring $ with respect to $ \chevstr $. Denote by $ (\rootgr{\alpha})_{\alpha \in \roots} $ the root groups of $ G $, by $ (\risom{\alpha})_{\alpha \in \roots} $ the root isomorphisms of $ G $ and by $ (w_\alpha)_{\alpha \in \roots} $ the standard Weyl elements from \cref{chev:weyl-def}. Then the Chevalley parity map
	\[ \map{\inverparsym}{\roots \times \roots}{\twistgroup}{}{} \]
	from \cref{chev:parmap} induces a $ \rootbase $-parity map $ \inverparsym' \defl \restrict{\inverparsym}{\roots \times \rootbase} $ with values in $ \twistgroup \defl \compactSet{\pm 1_\IZ} $. Letting $ \twistgroup $ act on all root groups $ (\rootgr{\alpha})_{\alpha \in \roots} $ by inversion, we turn it into a twisting group for $ (G, (w_\delta)_{\delta \in \rootbase}) $. Even more, $ (\twistgroup, \inverparsym') $ is a twisting system for $ (G, (w_\delta)_{\delta \in \rootbase}) $. Further, the tuple $ \calP \defl (\twistgroup, \ring, \ldots, \ring) $ is a parameter system of type $ \roots $, where $ \ring $ appears $ \abs{\Orb(\roots)} $ times in the tuple and $ \twistgroup $ acts on all copies of $ \ring $ by additive inversion. This parameter system is faithful if and only if $ 1_\ring \ne -1_\ring $, or in other words, $ 2_\ring \ne 0_\ring $. Finally, by \cref{chev:squarerel}, $ (\risom{\alpha})_{\alpha \in \roots} $ is a parametrisation of $ G $ by $ \calP $ with respect to $ (w_\delta)_{\delta \in \rootbase} $ and $ \inverparsym' $.
\end{example}

\begin{remark}[Choice of a parity map: The sign problem]\label{param:motiv:parmap-choice}
	Later on, we want to find a parity map $ \inverparsym_\roots $ such that every $ \roots $-graded group $ G $ is parametrised by some parameter system with respect to $ \inverparsym_\roots $. We refer to the problem of finding a suitable $ \inverparsym_\roots $ as the \defemph*{sign problem (for root graded groups)}\index{sign problem!for root graded groups}. When we specialise to Chevalley groups, it is precisely the sign problem that we have already seen in this context (see \cref{chev:signprob-weyl} and also \cref{chev:signprob-chevbasis,chev:signprob-commformula}).
	
	The best way to find a suitable parity map $ \inverparsym_\roots $ is to \enquote{read it off} from an appropriate \enquote{model} $ H $, that is, from a $ \roots $-graded group $ H $ which is known to be coordinatised by some algebraic structure. In \cref{param:param-parmap-has-properties,param:param-parmap-has-properties2}, we will show that a parity map defined in this way automatically has many desirable properties. For this strategy to work, we have to require that the parity map $ \inverparsym_\roots $ can actually be read off from the group $ H $ in a unique way. As a counterexample, consider a Chevalley group $ H $ defined over a ring $ \ring $ with $ 1_\ring = -1_\ring $. In this group, we have
	\[ \risom{\alpha}(\lambda)^{w_\delta} = \risom{\refl{\delta}(\alpha)}(\lambda) = \risom{\refl{\delta}(\alpha)}(-\lambda) \]
	for all roots $ \alpha \in \roots $, $ \delta \in \rootbase $ and all $ \lambda \in \ring $. Thus we cannot read off the sign $ \inverparbr{\alpha}{\delta} $ from this specific group $ H $. This example illustrates that the group $ H $ has to be \defemph*{sufficiently generic}\index{generic example group}, which roughly means that \enquote{all twisting phenomena which may arise in general $ \roots $-graded groups can be observed in $ H $}. Hence we have to require that $ H $ satisfies the faithfulness assumptions from \cref{param:faithful-act-def,param:parsys-faithful}.
	
	Recall \cref{rgg:coord-gen-def,rgg:existence-problem} for the precise formulation of the coordinatisation and existence problems. We conclude form the previous paragraph that our solution of the sign problem (and thus of the coordinatisation problem) is closely tied to the existence problem: Before we can coordinatise \emph{all} $ \roots $-graded groups, we have to construct \emph{one} sufficiently generic $ \roots $-graded group $ H $. In other words, we have to solve the existence problem for one sufficiently generic algebraic structure $ \calX $ in the class $ \algclass{\roots} $ which appears in the solution of the coordinatisation problem. However, we do not have to solve the existence problem for all algebraic structures in $ \algclass{\roots} $. In particular, in cases where $ \algclass{\roots} $ allows non-associative structures, it will be sufficient to solve the existence problem for a sufficiently generic associative $ \calX $ because the non-associativity is irrelevant for the twisting actions.
	
	For simply-laced root systems $ \roots $, the class $ \algclass{\roots} $ consists of rings and the only type of twisting which occurs is additive inversion. Hence we can use a Chevalley group $ H $ over a (commutative associatve) ring $ \ring $ with $ 2_\ring \ne 0_\ring $ as a model to construct $ \inverparsym_\roots $. For the other root systems, more general types of twistings occur and thus Chevalley groups are not sufficiently generic in these cases.
\end{remark}


\section{Compatibility Conditions}

\label{sec:param:comp}

\begin{secnotation}
	We continue to use \cref{param:secnot:param-def}.
\end{secnotation}

\begin{note}\label{param:combi-grp-separation-note}
	Until now, all conditions that we imposed on parity maps were independent of the twisting actions on $ G $, and vice versa. This clean separation is desirable: We will prove in \cref{param:param-parmap-has-properties,param:param-parmap-has-properties2} that every parity map $ \inverparsym $ which is \enquote{read off} from a group $ H $ automatically has all the properties that we want, but only with respect to $ H $. For properties which are phrased purely in terms of parity maps, such as braid-invariance, this is no restriction. However, the compatibility conditions that we now introduce relate $ \inverparsym $ only to the specific $ \roots $-graded group $ H $, so we still have to prove that arbitrary $ \roots $-graded groups (equipped with suitably defined twisting actions) are compatible with $ \inverparsym $. For square compatibility, this is usually easy because in most cases, the square formula provides a clean separation between the combinatoric side and the group-theoretic side. Stabiliser compatibility, on the other hand, will cause more difficulties.
\end{note}

\begin{definition}[Square compatibility]\label{param:square-comp-def}
	Let $ \inverparsym $ be a $ \rootbase $-parity map with values in some twisting group $ \twistgroup $ for $ G $. For any root $ \alpha $, we say that $ G $ is \defemph*{$ \alpha $-square-compatible with respect to $ \inverparsym $ (and $ (w_{\delta})_{\delta \in \rootbase} $)}\index{square-compatibility} if for all $ \delta \in \rootbase $, the element $ w_\delta^2 $ acts on $ \rootgr{\alpha} $ by $ \inverpar{\alpha}{\delta \delta} $. Further, we say that $ G $ is \defemph*{square-compatible with respect to $ \inverparsym $ (and $ (w_{\delta})_{\delta \in \rootbase} $)} if it is $ \alpha $-square-compatible with respect to $ \inverparsym $ (and $ (w_{\delta})_{\delta \in \rootbase} $) for all $ \alpha \in \roots $.
\end{definition}

\begin{definition}[Square formula]\label{param:square-formula-rgg-def}
	Let $ \alpha $ be a root. We say that $ G $ \defemph*{satisfies the square formula (for Weyl elements) for $ \alpha $}\index{square formula!for root graded groups} if for all roots $ \beta $, we have $ x_\alpha^w = x_\alpha^\epsilon $ for all $ x_\alpha \in \rootgr{\alpha} $ where $ w $ is the square of any $ \beta $-Weyl element and $ \epsilon \defl (-1)^{\cartanint{\alpha}{\beta}} $. Further, we say that $ G $ \defemph*{satisfies the square formula (for Weyl elements)} if it satisfies the square formula for all roots.
\end{definition}

\begin{definition}[Stabiliser compatibility]\label{param:stabcomp-def}
	Let $ \inverparsym, \invoparsym $ be two $ \rootbase $-parity maps with values in some abelian groups $ \twistgroup, \invogroup $ where $ \twistgroup $ is a twisting group for $ G $ and $ \invoparsym $ is Weyl-invariant. For any root $ \alpha \in \roots $, we say that $ G $ is \defemph*{$ \alpha $-stabiliser-compatible with respect to $ (\inverparsym, \invoparsym) $ (and $ (w_{\delta})_{\delta \in \rootbase} $)}\index{stabiliser-compatibility} if for all $ u \in \Weyl(\roots) $ such that $ \invopar{\alpha}{u} = 1_\invogroup $ and $ \alpha^u=\alpha $, there exists a word $ \word{\delta} $ over $ \rootbase \union (-\rootbase) $ such that $ \refl{\word{\delta}} = u $ and $ w_{\word{\delta}} $ acts on $ \rootgr{\alpha} $ by $ \inverpar{\alpha}{\word{\delta}} $. We say that $ G $ is \defemph*{stabiliser-compatible with respect to $ (\inverparsym, \invoparsym) $ (and $ (w_{\delta})_{\delta \in \rootbase} $)} if it is $ \alpha $-stabiliser-compatible with respect to $ (\inverparsym, \invoparsym) $ (and $ (w_{\delta})_{\delta \in \rootbase} $) for all $ \alpha \in \roots $. If these conditions are satisfied for the trivial parity map $ \invoparsym $, we will simply say that  $ G $ is \defemph*{($ \alpha $-) stabiliser-compatible with respect to $ \inverparsym $ (and $ (w_{\delta})_{\delta \in \rootbase} $)}
\end{definition}

\begin{remark}\label{param:stab-only-inver-rem}
	If $ \inverparsym $ is a parity map with values in some abelian group $ \twistgroup $ and $ G $ is $ \alpha $-stabiliser-compatible with respect to $ \inverparsym $ and $ (w_{\delta})_{\delta \in \roots} $, then for any abelian group $ \invogroup $ and any parity map $ \invoparsym $ with values in $ \invogroup $, the group $ G $ is $ \alpha $-stabiliser-compatible with respect to $ (\inverparsym, \invoparsym) $ and $ (w_{\delta})_{\delta \in \rootbase} $.
\end{remark}

\begin{note}
	For root systems with the property that any pair of orthogonal roots is (crystallographically) adjacent, we will see in \cref{param:adj-implies-stab} that $ G $ is stabiliser-compatible with respect to every (crystallographically) adjacency-trivial parity map.
\end{note}

As announced in \cref{param:motiv:parmap-choice}, we want to define the parity maps $ \inverparsym $ and $ \invoparsym $ by reading them off from an example group. More precisely, we will define a parity map $ \totalparsym $ in this way and then use \cref{param:parmap-split} to split it into two parity maps $ \inverparsym $, $ \invoparsym $ such that $ \totalparsym = \inverparsym \times \invoparsym $. The following result says that a parity map which is defined in this way automatically has many desirable properties. However, the restrictions of \cref{param:combi-grp-separation-note} apply. For this reason, only \cref{param:param-parmap-has-properties2} is truly useful, whereas \cref{param:param-parmap-has-properties} merely serves to illustrate the connection between some of the concepts we have introduced, and to show that the imposed conditions are natural.

\begin{lemma}\label{param:param-parmap-has-properties}
	Write $ \IndivOrb(\roots) = \tup{O}{k} $ (see \cref{rootsys:orb-def}). Let $ \calP = (\twistgroup, \listing{M}{k}) $ be a parameter system of type $ \roots $ and let $ \inverparsym $ be a $ \rootbase $-parity map with values in $ \rootbase $. Assume that there exists a parametrisation $ (\risom{\alpha})_{\alpha \in \roots} $ of $ G $ by $ \calP $ with respect to $ \inverparsym $ and $ (w_\delta)_{\delta \in \rootbase} $. For each root $ \alpha $, we define an action $ \omega_\alpha $ of $ \twistgroup $ on $ \rootgr{\alpha} $ by
	\[ a.\risom{\alpha}(m) \defl \risom{\alpha}(a.m) \quad \text{for all } m \in M_\alpha. \]
	Then the following hold:
	\begin{lemenumerate}
		\item \label{param:param-parmap-has-properties:twist}$ (\twistgroup, (\omega_\alpha)_{\alpha \in \roots}) $ is a twisting group for $ G $.
		
		\item \label{param:param-parmap-has-properties:square}$ G $ is square-compatible and stabiliser-compatible with respect to $ \inverparsym $ and $ (w_\delta)_{\delta \in \rootbase} $.
	\end{lemenumerate}
\end{lemma}
\begin{proof}
	For~\itemref{param:param-parmap-has-properties:twist}, the only non-trivial thing to show is that the twisting actions commute with conjugation by the fixed Weyl elements. Let $ \alpha \in \roots $, $ \delta \in \rootbase \union (-\rootbase) $ and let $ x_\alpha \in \rootgr{\alpha} $. Let $ a \in \twistgroup $ and put $ m \defl \risom{\alpha}^{-1}(x_\alpha) $. Then
	\begin{align*}
		(a.x_\alpha)^{w_\delta} &= (a.\risom{\alpha}(m))^{w_\delta} =  \risom{\alpha}(a.m)^{w_\delta} = \risom{\refl{\delta}(\alpha)}(\inverpar{\alpha}{\delta} a.m) \rightand \\
		a.(x_\alpha^{w_\alpha}) &= a.\brackets[\big]{\risom{\alpha}(m)^{w_\delta}} = a.\risom{\refl{\delta}(\alpha)}(\inverpar{\alpha}{\delta}.m) = \risom{\refl{\delta}(\alpha)}(a\inverpar{\alpha}{\delta}.m).
	\end{align*}
	Since $ \twistgroup $ is abelian, we conclude that $ (a.x_\alpha)^{w_\delta} = a.(x_\alpha^{w_\alpha}) $, as desired.

	For \itemref{param:param-parmap-has-properties:square}, let $ \alpha \in \roots $, $ \delta \in \rootbase $ and $ x_\alpha \in \rootgr{\alpha} $. Then we have
	\begin{align*}
		x_\alpha^{w_\delta^2} &= \risom{\alpha}(m)^{w_\delta^2} = \risom{\refl{\delta}(\alpha)}(\inverpar{\alpha}{\delta}.m)^{w_\delta} =  \risom{\alpha}(\inverpar{\refl{\delta}(\alpha)}{\delta} \inverpar{\alpha}{\delta}.m) = \risom{\alpha}(\inverpar{\alpha}{\delta \delta}.m) = \inverpar{\alpha}{\delta \delta}.x_\alpha
	\end{align*}
	where $ m \defl \risom{\alpha}^{-1}(x_\alpha) $. Thus $ G $ is square-compatible. For stabiliser compatibility, we have to show that for all $ \alpha \in \roots $ and all $ u \in \Weyl(\roots) $ with $ \alpha^u = \alpha $, we can find a representation $ \word{\delta} $ of $ u $ over $ \rootbase $ such that $ w_{\word{\delta}} $ acts on $ \rootgr{\alpha} $ by $ \inverpar{\alpha}{\word{\delta}} $. Since $ G $ is parametrised with respect to $ \inverparsym $ and $ (w_\delta)_{\delta \in \rootbase} $, this actually holds for all representations $ \word{\delta} $ of $ u $.
\end{proof}

\begin{proposition}\label{param:param-parmap-has-properties2}
	Let everything be as in \cref{param:param-parmap-has-properties} and let $ \alpha $ be a root. Assume that the action of $ \twistgroup $ on $ M_\alpha $ is $ \alpha $-faithful with respect to $ \inverparsym $. Then the following hold:
	\begin{proenumerate}
		\item \label{param:param-parmap-has-properties2:adj-triv}If $ G $ has (crystallographic) $ \roots $-commutator relations with root groups $ (\rootgr{\alpha})_{\alpha \in \roots} $, then$ \inverparsym $ is (crystallographically) $ \alpha $-adjacency-trivial.

		\item \label{param:param-parmap-has-properties2:braid}If the family $ (w_\delta)_{\delta \in \rootbase} $ satisfies the braid relations modulo $ \zentrum(G) $, then $ \inverparsym $ is $ \alpha $-braid-invariant.
		
		\item \label{param:param-parmap-has-properties2:square}Assume that $ G $ satisfies the square formula for Weyl elements and that $ \twistgroup = \compactSet{\pm 1}^p $ for some $ p \in \Npos $. Assume further that $ (-1, 1, \ldots, 1) $ lies in $ \genparmoveset{\twistgroup}{\inverparsym}{\alpha}{\alpha} $ and that it acts on all parametrising groups $ \listing{M}{k} $ by inversion. Then $ \inverparsym $ satisfies the square formula for $ \alpha $.
	\end{proenumerate}
\end{proposition}
\begin{proof}
	For \itemref{param:param-parmap-has-properties2:adj-triv}, let $ \beta $ be a root such that $ \alpha $ is (crystallographically) adjacent to $ \beta $ and $ -\beta $ and let $ \word{\rho} $ be a $ \rootbase $-expression of $ \beta $. Then $ \reflbr{\word{\rho}} = \reflbr{\beta} $ and $ w \defl w_{\word{\rho}} $ is a $ \beta $-Weyl element. By adjacency assumption, it follows that $ w $ acts trivially on $ \rootgr{\alpha} $, so that
	\[ \risom{\alpha}(m) = \risom{\alpha}(m)^w = \risom{\alpha^{\reflbr{\word{\rho}}}}(\inverpar{\alpha}{\word{\rho}}.m) = \risom{\alpha^{\reflbr{\beta}}}(\inverpar{\alpha}{\word{\rho}}.m) = \risom{\alpha}(\inverpar{\alpha}{\word{\rho}}.m) \]
	for all $ m \in M_\alpha $. Since the action of $ \twistgroup $ on $ M_\alpha $ is $ \alpha $-faithful with respect to $ \inverparsym $, this implies that $ \inverpar{\alpha}{\word{\rho}} = 1_\twistgroup $. Thus $ \inverparsym $ is (crystallographically) $ \alpha $-adjacency-trivial.
	
	For~\itemref{param:param-parmap-has-properties2:braid}, let $ \beta, \gamma \in \rootbase $ be distinct, denote the order of $ \refl{\beta} \refl{\gamma} $ by $ o $ and put $ \word{\rho} \defl \braidword{o}(\beta, \gamma) $, $ \word{\zeta} \defl \braidword{o}(\gamma, \beta) $. By the braid relation in the Weyl group, we have $ \refl{\word{\rho}} = \refl{\word{\zeta}} $ and we denote this element by $ \phi $. Now for all $ m \in M_\alpha $, we see that $ \risom{\alpha}(m)^{w_{\word{\rho}}} = \risom{\alpha^\phi}(\inverpar{\alpha}{\word{\rho}}.m) $ and $ \risom{\alpha}(m)^{w_{\word{\zeta}}} = \risom{\alpha^\phi}(\inverpar{\alpha}{\word{\zeta}}.m) $. Since $ (w_\delta)_{\delta \in \rootbase} $ satisfies the braid relations modulo $ \zentrum(G) $, these two group elements are equal. We infer that $ \inverpar{\alpha}{\word{\rho}} $ and $ \inverpar{\alpha}{\word{\zeta}} $ act identically on $ M_\alpha $. Since the action of $ \twistgroup $ on $ M_\alpha $ is $ \alpha $-faithful with respect to $ \inverparsym $, it follows that $ \inverpar{\alpha}{\word{\rho}} = \inverpar{\alpha}{\word{\zeta}} $, so $ \inverparsym $ is $ \alpha $-braid-invariant.
	
	For~\itemref{param:param-parmap-has-properties2:square}, let $ \delta \in \rootbase $ and put $ \epsilon \defl (-1)^{\cartanint{\alpha}{\delta}} $ and $ w \defl w_\delta^2 $. Let $ m \in M_\alpha $ be arbitrary. Since $ G $ satisfies the square formula for Weyl elements, we have
	\begin{align*}
		\risom{\alpha}(\inverpar{\alpha}{\delta \delta}.m) = \risom{\alpha}(m)^w = \risom{\alpha}(\epsilon m) = \risom{\alpha}\brackets[\big]{(\epsilon, 1, \ldots, 1).m}.
	\end{align*}
	Hence $ \inverpar{\alpha}{\delta \delta} $ and $ (\epsilon, 1, \ldots, 1) $ act identically on $ M_\alpha $. Since both elements lie in $ \parmoveset{\twistgroup}{\alpha}{\alpha} $ (by our additional assumption if $ \epsilon = -1 $), we infer that they are equal. Thus $ \inverparsym $ satisfies the square formula for $ \alpha $.
\end{proof}


\section{Proof of the Parametrisation Theorem}

\label{sec:param:proof}

\begin{secnotation}\label{param:partwist-conv}
	We denote by $ \roots $ a root system and by $ G $ is a group with a $ \roots $-pregrading $ (\rootgr{\alpha})_{\alpha \in \roots} $ such that $ \rootgr{\lambda \alpha} \subs \rootgr{\alpha} $ for all $ \alpha \in \roots $ and all $ \lambda \in \IR_{>1} $ for which $ \lambda \alpha $ is a root. Further, we denote by $ \rootbase $ a rescaled root base of $ \roots $, by $ (w_\delta)_{\delta \in \rootbase} $ a $ \rootbase $-system of Weyl elements in $ G $ and by $ (\twistgroup, \inverparsym, \invogroup, \invoparsym) $ a partial twisting system for $ (G, (w_\delta)_{\delta \in \rootbase}) $. We assume that $ G $ is square-compatible with respect to $ \inverparsym $ and stabiliser-compatible with respect to $ (\inverparsym, \invoparsym) $ and that $ (w_\delta)_{\delta \in \rootbase} $ satisfies the braid relations modulo $ \zentrum(G) $. Further, \cref{param:thm-notation} holds from the point where it is introduced.
\end{secnotation}

Finally, we have gathered all the necessary tools to prove the para\-me\-tri\-sation theorem (\cref{param:thm}). Our first intermediate goal is to prove \cref{param:abstract-thm}.

\begin{note}[Braid relations]
	If $ G $ is a $ \roots $-graded group and $ \roots $ is reduced or crystallographic, then it automatically satisfies the braid relations for Weyl elements by \cref{braid:all}. However, since the proof of the parametrisation theorem works without any reference to the commutator relations and Axiom~\thmitemref{rgg-def}{rgg-def:nondeg}, we formulate it for more general groups with a $ \roots $-pregrading. Hence we need to make the assumption that $ (w_\delta)_{\delta \in \rootbase} $ satisfies the braid relations modulo $ \zentrum(G) $. Note that the braid relations are only required for this specific system of Weyl elements, not for any system, and that they are only required to hold modulo $ \zentrum(G) $.
\end{note}

\begin{lemma}\label{param:homotopic-eq}
	Let $ \alpha $ be a root and let $ \word{\beta} = (\beta_1, \ldots, \beta_m) $, $ \word{\gamma} = (\gamma_1, \ldots, \gamma_k) $ be two homotopic words over $ \rootbase \union (-\rootbase) $ (in the sense of \cref{weyl:homotopy-def}) such that $ \alpha^{\reflbr{\word{\beta}}} = \alpha = \alpha^{\reflbr{\word{\gamma}}} $. Then the following assertions are equivalent:
	\begin{stenumerate}
		\item $ w_{\word{\beta}} $ acts on $ \rootgr{\alpha} $ by $ \inverpar{\alpha}{\word{\beta}} $.
		
		\item $ w_{\word{\gamma}} $ acts on $ \rootgr{\alpha} $ by $ \inverpar{\alpha}{\word{\gamma}} $.
	\end{stenumerate}
\end{lemma}
\begin{proof}
	By induction, we can assume that $ \word{\gamma} $ and $ \word{\beta} $ are either braid-homotopic or square-homotopic. In the case of braid-homotopy, we have $ \inverpar{\alpha}{\word{\beta}} = \inverpar{\alpha}{\word{\gamma}} $ (because $ \inverparsym $ is braid-invariant) and the conjugation actions of $ w_{\word{\beta}} $ and $ w_{\word{\gamma}} $ on $ G $ are identical (because $ (w_\delta)_{\delta \in \rootbase} $ satisfies the braid relations modulo $ \zentrum(G) $). Hence the assertion is clear for braid-homotopy. In the case of square-homotopy, we can assume that
	\[ \word{\gamma} = (\beta_1, \ldots, \beta_{p}, \delta, \delta, \beta_{p+1}, \ldots, \beta_m) \]
	for some $ p \in \numint{0}{m} $ and some $ \delta \in \rootbase $. Put
	\[ \word{\omega} \defl (\beta_1, \ldots, \beta_p), \qquad \word{\zeta} \defl (\beta_{p+1}, \ldots, \beta_m) \quad \text{and} \quad \rho \defl \alpha^{\reflbr{\word{\omega}}}. \]
	Observe that
	\begin{align*}
		\inverparbr{\alpha}{\word{\gamma}} &= \inverparbr{\alpha}{\word{\omega}} \inverparbr{\rho}{\delta \delta} \inverparbr{\rho}{\word{\zeta}} = \inverparbr{\alpha}{\word{\omega}} \inverparbr{\rho}{\word{\zeta}} \inverparbr{\rho}{\delta \delta} = \inverparbr{\alpha}{\word{\beta}} \inverparbr{\rho}{\delta \delta}
	\end{align*}
	and $ w_{\word{\gamma}} = w_{\word{\omega}} w_{\delta}^2 w_{\word{\zeta}} $. Since $ G $ is square-compatible with respect to $ \inverparsym $ and the action of $ \twistgroup $ commutes with conjugation by $ w_{\word{\zeta}} $ (by Axiom~\thmitemref{param:twist-grp-def}{param:twist-grp-def:comm}), we have
	\begin{align*}
		x_\alpha^{w_{\word{\gamma}}} &= \brackets[\big]{\inverpar{\rho}{\delta\delta}.x_\alpha^{w_{\word{\omega}}}}^{w_{\word{\zeta}}} = \inverpar{\rho}{\delta \delta}.x_\alpha^{w_{\word{\omega}} w_{\word{\zeta}}} = \inverpar{\rho}{\delta \delta}.x_\alpha^{w_{\word{\beta}}}.
	\end{align*}
	Since $ \inverpar{\alpha}{\word{\gamma}} = \inverpar{\alpha}{\word{\beta}} \inverpar{\rho}{\delta \delta} $, the assertion follows.
\end{proof}

The following statement allows us to restrict our attention to words over $ \rootbase $ (and not over $ \rootbase \union (-\rootbase) $).

\begin{lemma}\label{param:homotopic-minus}
	Let $ \alpha $ be a root and let $ \word{\beta} = (\beta_1, \ldots, \beta_m) $ be a word over $ \rootbase \union (-\rootbase) $ such that $ \alpha^{\reflbr{\word{\beta}}} = \alpha $. Choose an arbitrary $ i \in \numint{1}{m} $ and put
	\[ \word{\beta}' \defl (\beta_1, \ldots, \beta_{i-1}, -\beta_i, \beta_{i+1}, \ldots, \beta_m) \]
	Then the following assertions are equivalent:
	\begin{stenumerate}
		\item $ w_{\word{\beta}} $ acts on $ \rootgr{\alpha} $ by $ \inverpar{\alpha}{\word{\beta}} $.
		
		\item \label{param:homotopic-minus:minus}$ w_{\word{\beta}'} $ acts on $ \rootgr{\alpha} $ by $ \inverpar{\alpha}{\word{\beta}'} $.
	\end{stenumerate}
\end{lemma}
\begin{proof}
	Note that the word $ \word{\beta}'' \defl (\beta_1, \ldots, \beta_{i-1}, -\beta_i, \beta_i, \beta_i, \beta_{i+1}, \ldots, \beta_m) $ is square-homotopic to $ \word{\beta}' $. Thus it follows from \cref{param:homotopic-eq} that~\itemref{param:homotopic-minus:minus} is equivalent to the statement \enquote{$ w_{\word{\beta}''} $ acts on $ \rootgr{\alpha} $ by $ \inverpar{\alpha}{\word{\beta}''} $}. However, $ w_{\word{\beta}''} = w_{\word{\beta}} $ (because $ w_{-\beta_i} = w_{\beta_i}^{-1} $) and $ \inverpar{\alpha}{\word{\beta}''} = \inverpar{\alpha}{\word{\beta}} $ (because of \cref{param:inverpar-minus-cancel}). The assertion follows.
\end{proof}

\begin{lemma}\label{param:inverpar-stab}
	Let $ \alpha \in \roots $ and let $ \word{\delta} = \tup{\delta}{m} $ be a word over $ \rootbase \union (-\rootbase) $ such that $ \alpha^{\reflbr{\word{\delta}}} = \alpha $ and $ \invopar{\alpha}{\word{\delta}} = 1_\invogroup $. Then $ w_{\word{\delta}} $ acts on $ \rootgr{\alpha} $ by $ \inverpar{\alpha}{\word{\delta}} $.
\end{lemma}
\begin{proof}
	Since $ G $ is stabiliser-compatible with respect to $ (\inverparsym, \invoparsym) $, there exists a word $ \word{\rho} $ over $ \rootbase \union (-\rootbase) $ such that $ \refl{\word{\rho}} = \refl{\word{\delta}} $ and such that $ w_{\word{\rho}} $ acts on $ \rootgr{\alpha} $ by $ \inverpar{\alpha}{\word{\rho}} $. Denote by $ \word{\rho}' $ and $ \word{\delta}' $ the words obtained from $ \word{\rho} $ and $ \word{\delta} $ by replacing each letter in $ -\rootbase $ by its negative. Then $ \refl{\word{\rho}'} = \refl{\word{\rho}} = \refl{\word{\delta}} = \refl{\word{\delta}'} $. Hence it follows from \cref{word-problem-solution} that $ \word{\rho}' $ and $ \word{\delta}' $ are homotopic. (Here we use that $ \word{\rho}' $ and $ \word{\delta}' $ are words over $ \rootbase $ and not over $ \rootbase \union (-\rootbase) $. See also \cref{param:homotopy-def-note}.) Since $ w_{\word{\rho}} $ acts on $ \rootgr{\alpha} $ by $ \inverpar{\alpha}{\word{\rho}} $, it follows from \cref{param:homotopic-minus} that $ w_{\word{\rho}'} $ acts on $ \rootgr{\alpha} $ by $ \inverpar{\alpha}{\word{\rho}'} $. By \cref{param:homotopic-eq}, this implies that $ w_{\word{\delta}'} $ acts on $ \rootgr{\alpha} $ by $ w_{\word{\delta}'} $. Again by \cref{param:homotopic-minus}, it follows that $ w_{\word{\delta}} $ acts on $ \rootgr{\alpha} $ by $ w_{\word{\delta}} $. This finishes the proof.
\end{proof}

\begin{proposition}\label{param:abstract-thm}
	Let $ \alpha \in \roots $ and let $ \word{\beta} = \tup{\beta}{m} $, $ \word{\gamma} = (\gamma_1, \ldots, \gamma_k) $ be two words over $ \rootbase \union (-\rootbase) $ such that $ \alpha^{\reflbr{\word{\beta}}} = \alpha^{\reflbr{\word{\gamma}}} $ and $ \invopar{\alpha}{\word{\beta}} = \invopar{\alpha}{\word{\gamma}} $. Put $ a \defl \inverpar{\alpha}{\word{\beta}} \cdot \inverpar{\alpha}{\word{\gamma}}^{-1} $. Then $ x_\alpha^{w_{\word{\beta}}} = (a.x_\alpha)^{w_{\word{\gamma}}} $ for all $ x_\alpha \in \rootgr{\alpha} $.
\end{proposition}
\begin{proof}
	Consider $ \word{\gamma}^{-1} = (-\gamma_k, \ldots, -\gamma_1) $ and $ \word{\delta} \defl (\word{\beta}, \word{\gamma}^{-1}) $. Since $ \rootgr{\alpha}^{\reflbr{\word{\beta}}} = \rootgr{\alpha}^{\reflbr{\word{\gamma}}} $, we have $ \alpha^{w_{\word{\delta}}} = \alpha $. Further, the following hold by \cref{param:inverpar-backwards}:
	\begin{align*}
		\inverparbr{\alpha}{\word{\delta}} &= \inverparbr{\alpha}{\word{\beta}} \inverparbr{\alpha^{\reflbr{\word{\beta}}}}{\word{\gamma}^{-1}} = \inverparbr{\alpha}{\word{\beta}} \inverparbr{\alpha^{\reflbr{\word{\gamma}}}}{\word{\gamma}^{-1}} = \inverparbr{\alpha}{\word{\beta}} \inverparbr{\alpha}{\word{\gamma}}^{-1} = a, \\
		\invoparbr{\alpha}{\word{\delta}} &= \invoparbr{\alpha}{\word{\beta}} \invoparbr{\alpha^{\reflbr{\word{\beta}}}}{\word{\gamma}^{-1}} = \invoparbr{\alpha}{\word{\beta}} \invoparbr{\alpha^{\reflbr{\word{\gamma}}}}{\word{\gamma}^{-1}} = \invoparbr{\alpha}{\word{\beta}} \invoparbr{\alpha}{\word{\gamma}}^{-1} = 1_\invogroup.
	\end{align*}
	It follows from the second equation and \cref{param:inverpar-stab} that $ w_{\word{\delta}} $ acts on $ \rootgr{\alpha} $ by $ a $, that is, by $ \inverparbr{\alpha}{\word{\gamma}}^{-1} $. Therefore, we have shown that
	\[ x_\alpha^{w_{\word{\beta}} w_{\word{\gamma}^{-1}}} = x_\alpha^{w_{\word{\delta}}} = a.x_\alpha \]
	for all $ x_\alpha \in \rootgr{\alpha} $. This finishes the proof.
\end{proof}

The following special case of \cref{param:abstract-thm} will often be useful.

\begin{proposition}\label{param:inver-invo-eq}
	Let $ \alpha \in \roots $ and let $ \word{\beta} = \tup{\beta}{m} $, $ \word{\gamma} = (\gamma_1, \ldots, \gamma_k) $ be two words over $ \rootbase \union (-\rootbase) $ such that $ \alpha^{\reflbr{\word{\beta}}} = \alpha^{\reflbr{\word{\gamma}}} $, $ \inverpar{\alpha}{\word{\beta}} = \inverpar{\alpha}{\word{\gamma}} $ and $ \invopar{\alpha}{\word{\beta}} = \invopar{\alpha}{\word{\gamma}} $. Then $ x_\alpha^{w_{\word{\beta}}} = x_\alpha^{w_{\word{\gamma}}} $ for all $ x_\alpha \in \rootgr{\alpha} $.
\end{proposition}

\begin{secnotation}\label{param:thm-notation}
	From now on, we write $ \IndivOrb(\roots) = \tup{O}{k} $ and for each root $ \alpha $, we fix a complement $ C_\alpha $ of $ \genparmoveset{\invogroup}{\invoparsym}{\alpha}{\alpha} $ in $ \invogroup $ (which exists because $ \invoparsym $ is semi-complete).
\end{secnotation}

Thanks to \cref{param:inver-invo-eq}, we can now define an action of $ \invogroup $ on the root groups. At this point, the semi-completeness of $ \invoparsym $ finally becomes relevant.

\begin{definition}\label{param:invo-def}
	Let $ \alpha \in \roots $. Denote by $ \map{\pi_\alpha}{\invogroup}{\parmoveset{\invogroup}{\alpha}{\alpha}}{}{} $ the canonical projection with respect to the decomposition $ \invogroup = \parmoveset{\invogroup}{\alpha}{\alpha} \times C_\alpha $.	We define an action of $ \invogroup $ on the set $ \rootgr{\alpha} $ by $ b.x_\alpha \defl x_\alpha^{w_{\word{\delta}}} $ for all $ b \in \invogroup $ and $ x_\alpha \in \rootgr{\alpha} $ where $ \word{\delta} $ is any word over $ \rootbase \union (-\rootbase) $ such that $ \alpha^{\reflbr{\word{\delta}}} = \alpha $, $ \inverpar{\alpha}{\word{\delta}} = 1_\twistgroup $ and $ \invopar{\alpha}{\word{\delta}} = \pi_\alpha(b) $. (This action is well-defined by \cref{param:inver-invo-eq} and because $ \inverparsym $ and $ \invoparsym $ are independent.) If $ \word{\delta} $ is any word with these properties, we say that \defemph*{$ \word{\delta} $ induces the action of $ b $ on $ \rootgr{\alpha} $}.
\end{definition}

\begin{remark}
	Equivalently, we could define the action in \cref{param:invo-def} by first defining an action of $ \parmoveset{\invogroup}{\alpha}{\alpha} $ on $ \rootgr{\alpha} $ in the same way, but without using $ \pi_\alpha $, and then declaring that $ C_\alpha $ acts trivially on $ \rootgr{\alpha} $.
\end{remark}

\begin{note}\label{param:semi-comp-note}
	Without the assumption that $ \invoparsym $ is semi-complete, it is not clear that the action of $ \parmoveset{\invogroup}{\alpha}{\alpha} $ in \cref{param:invo-def} can be extended to all of $ \invogroup $. Since a different choice of the complement $ C_\alpha $ would, in general, yield a different action, we have fixed the choice of $ (C_\alpha)_{\alpha \in \roots} $ in \cref{param:thm-notation}. In practice, we will only ever be interested in the action of elements from $ \parmoveset{\invogroup}{\alpha}{\alpha} $, so this choice of complements is harmless in practice.
	
	This observation has two interesting consequences: Firstly, instead of declaring that $ C_\alpha $ acts trivially on $ \rootgr{\alpha} $, we could choose an arbitrary action and the results of this section would remain true (because these results do not make any assertions about this action). Secondly, it implies that suitable modifications of the statements in this section remain true even if $ \invoparsym $ is not semi-complete if we take care to replace $ \invogroup $ by $ \parmoveset{\invogroup}{\alpha}{\alpha} $ in appropriate places. However, the assumption of semi-completeness is satisfied in all practical examples and it simplifies notation, which is we why we use it.
\end{note}

As a next goal, we want to show that $ \twistgroup \times \invogroup $ is a twisting group for $ (G, (w_\delta)_{\delta \in \rootbase}) $ (\cref{param:product-twistgrp}).

\begin{lemma}
	For each $ \alpha \in \roots $, the action $ \omega_\alpha $ of $ \invogroup $ on $ \rootgr{\alpha} $ in \cref{param:invo-def} is a well-defined group action. With these actions, $ (\invogroup, (\omega_\alpha)_{\alpha \in \roots}) $ is a twisting group for $ (G, (w_\delta)_{\delta \in \rootbase}) $.
\end{lemma}
\begin{proof}
	Let $ \alpha \in \roots $. For every $ b \in \invogroup $ there exists a word $ \word{\delta} $ which induces the action of $ b $ on $ \rootgr{\alpha} $. To show that this defines a group action, we have to show that $ 1_\invogroup.x_\alpha = x_\alpha $ and $ b_1 b_2.x_\alpha = b_1.(b_2.x_\alpha) $ for all $ x_\alpha \in \rootgr{\alpha} $ and all $ b_1, b_2 \in \invogroup $. The first assertion is clear because the empty word $ \emptyset $ satisfies $ \invopar{\alpha}{\emptyset} = 1_\invogroup $, so that
	\[ 1_\invogroup.x_\alpha = x_\alpha^{w_{\emptyset}} = x_\alpha^{1_G} = x_\alpha \quad \text{for all } x_\alpha \in \rootgr{\alpha}. \]
	For the second assertion, let $ \word{\delta}_1, \word{\delta}_2 $ be two words over $ \rootbase \union (-\rootbase) $ such that
	\[ \alpha^{\reflbr{\word{\delta}_i}} = \alpha, \qquad \inverpar{\alpha}{\word{\delta}_i} = 1_\twistgroup \quad \text{and} \quad  \invopar{\alpha}{\word{\delta}_i} = \pi_\alpha(b_i) \]
	for all $ i \in \Set{1,2} $. Then the word $ \word{\delta} \defl (\word{\delta}_1, \word{\delta}_2) $ satisfies
	\[ \alpha^{\reflbr{\word{\delta}}} = \alpha, \qquad \inverpar{\alpha}{\word{\delta}} = 1_\twistgroup \quad \text{and} \quad \invopar{\alpha}{\word{\delta}} = \pi_\alpha(b_1 b_2), \]
	so that
	\[ b_1 b_2.x_\alpha = x_\alpha^{w_{\word{\delta}}} = \brackets[\big]{x_\alpha^{w_{\word{\delta}_1}}}^{w_{\word{\delta}_2}} = b_2.(b_1.x_\alpha) \quad \text{for all } x_\alpha \in \rootgr{\alpha}. \]
	Since $ \invogroup $ is abelian, we conclude that $ \omega_\alpha $ is indeed an action.
	
	To see that $ (\invogroup, (\omega_\alpha)_{\alpha \in \roots}) $ is a twisting group for $ G $, it remains to show that the twisting action is compatible with conjugation by the fixed set of Weyl elements and their inverses. Let $ \alpha \in \roots $, $ \delta \in \rootbase \union (-\rootbase) $, $ x_\alpha \in \rootgr{\alpha} $ and $ b \in \invogroup $. Choose words $ \word{\omega}, \word{\zeta} $ over $ \rootbase \union (-\rootbase) $ such that $ \word{\omega} $ induces the action of $ b $ on $ \rootgr{\alpha} $ and $ \word{\zeta} $ induces the action of $ b $ on $ \rootgr{\alpha^{\reflbr{\delta}}} $. Then
	\[ (b.x_\alpha)^{w_{\delta}} = x_\alpha^{w_{\word{\omega}} w_{\delta}} \midand b.(x_\alpha^{w_\delta}) = x_\alpha^{w_{\delta} w_{\word{\zeta}}}, \]
	and also
	\[ \inverpar{\alpha}{\word{\omega}} = 1_\twistgroup = \inverpar{\alpha^{\reflbr{\delta}}}{\word{\zeta}} \midand \alpha^{\reflbr{\word{\omega}}} = \alpha = \alpha^{\reflbr{\word{\zeta}}}. \]
	We now compute that
	\begin{align*}
		\inverpar{\alpha}{\word{\omega} \delta} &= \inverpar{\alpha}{\word{\omega}} \inverpar{\alpha}{\delta} = \inverpar{\alpha}{\delta}, & \inverpar{\alpha}{\delta \word{\zeta}} &= \inverpar{\alpha}{\delta} \inverpar{\alpha^{\reflbr{\delta}}}{\word{\zeta}} = \inverpar{\alpha}{\delta}, \\
		\invopar{\alpha}{\word{\omega} \delta} &= \invopar{\alpha}{\word{\omega}} \invopar{\alpha}{\delta} = b \invopar{\alpha}{\delta}, &	\invopar{\alpha}{\delta \word{\zeta}} &= \invopar{\alpha}{\delta} \invopar{\alpha^{\reflbr{\delta}}}{\word{\zeta}} = \invopar{\alpha}{\delta} b = b \invopar{\alpha}{\delta}.
	\end{align*}
	Thus $ \inverpar{\alpha}{\word{\omega} \delta} = \inverpar{\alpha}{\delta \word{\omega}} $ and $ \invopar{\alpha}{\word{\omega} \delta} = \invopar{\alpha}{\delta \word{\omega}} $. By an application of \cref{param:inver-invo-eq}, we conclude that the actions of $ w_{\word{\omega}} w_{\delta} $ and $ w_{\delta} w_{\word{\zeta}} $ on $ \rootgr{\alpha} $ are identical, as desired.
\end{proof}

\begin{lemma}\label{param:twist-invo-commute}
	Let $ \alpha \in \roots $. Then the actions of $ \twistgroup $ and $ \invogroup $ on $ \rootgr{\alpha} $ commute. That is, $ a.(b.x_\alpha) = b.(a.x_\alpha) $ for all $ a \in \twistgroup $, $ b \in \invogroup $ and $ x_\alpha \in \rootgr{\alpha} $.
\end{lemma}
\begin{proof}
	Let $ a \in \twistgroup $, $ b \in \invogroup $ and $ x_\alpha \in \rootgr{\alpha} $. Choose a word $ \word{\zeta} $ over $ \rootbase \union (-\rootbase) $ which induces the action of $ b $ on $ \rootgr{\alpha} $. Using Axiom~\thmitemref{param:twist-grp-def}{param:twist-grp-def:comm}, we see that
	\begin{align*}
		a.(b.x_\alpha) &= a.x_\alpha^{w_{\word{\zeta}}} = (a.x_\alpha)^{w_{\word{\zeta}}} = b.(a.x_\alpha),
	\end{align*}
	as desired.
\end{proof}

\begin{proposition}\label{param:product-twistgrp}
	$ \twistgroup \times \invogroup $ is a twisting group for $ G $ (with the natural action induced by the actions of $ \twistgroup $ and $ \invogroup $) and $ (\twistgroup \times \invogroup, \inverparsym \times \invoparsym) $ is a twisting system for $ G $.
\end{proposition}
\begin{proof}
	By \cref{param:twist-invo-commute}, the actions of $ \twistgroup $ and $ \invogroup $ induce a natural action of $ \twistgroup \times \invogroup $ on each root group. Since the actions of $ \twistgroup $ and $ \invogroup $ commute with conjugation by Weyl elements, so does the action of $ \twistgroup \times \invogroup $. Thus $ \twistgroup \times \invogroup $ is indeed a twisting group for $ G $. Further, it is clear that $ \inverparsym \times \invoparsym $ is braid-invariant (because both $ \inverparsym $ and $ \invoparsym $ are, see \thmitemcref{param:prod-eq-lem}{param:prod-eq-lem:eq}) and transporter-invariant (by Axiom~\thmitemref{param:partwist-def}{param:partwist-def:stab-comp}), so $ (\twistgroup \times \invogroup, \inverparsym \times \invoparsym) $ is a twisting system for $ G $.
\end{proof}

\begin{notation}
	We define $ \totalparsym \defl \inverparsym \times \invoparsym $.
\end{notation}

\begin{theorem}[Parametrisation theorem]\label{param:thm}
	There exist a parameter system $ \calP = (\twistgroup \times \invogroup, \listing{M}{k}) $ of type $ \roots $ (in the sense of \cref{param:parsys-def}) whose twisting group is $ \twistgroup \times \invogroup $ and a parametrisation $ (\map{\risom{\alpha}}{M_\alpha}{\rootgr{\alpha}}{}{})_{\alpha \in \roots} $ of $ G $ by $ \calP $ with respect to $ \totalparsym $ and $ (w_\delta)_{\delta \in \rootbase} $ (in the sense of \cref{param:param-def}) such that for all roots $ \alpha $, the twisting action of $ \twistgroup \times \invogroup $ on $ M_\alpha $ is compatible with the twisting action of $ \twistgroup \times \invogroup $ on $ \rootgr{\alpha} $, meaning that
	\[ c.\risom{\alpha}(m) = \risom{\alpha}(c.m) \]
	for all $ c \in \twistgroup \times \invogroup $ and all $ m \in M_\alpha $. 
\end{theorem}
\begin{proof}
	For the whole proof, we fix an arbitrary $ i \in \numint{1}{k} $. Choose an arbitrary root $ \alpha_i \in O_i $, a group $ M_i $ which is isomorphic to $ \rootgr{\alpha_i} $ and an isomorphism $ \map{\risom{\alpha_i}}{M_i}{\rootgr{\alpha_i}}{}{} $. For each $ \beta \in O_i $, we define an isomorphism
	\[ \map{\risom{\beta}}{M_i}{\rootgr{\beta}}{x}{\risom{\alpha_i}(x)^{w_{\word{\omega}}}} \]
	where $ \word{\omega} $ is any word over $ \rootbase \union (-\rootbase) $ such that
	\[ \alpha_i^{\reflbr{\word{\omega}}} = \beta \midand \totalpar{\alpha}{\word{\omega}} = (1_\twistgroup, 1_\invogroup). \]
	Since $ \totalparsym $ is transporter-invariant, such a word $ \word{\omega} $ exists by \cref{param:transp-invar-1}. Further, we know from \cref{param:inver-invo-eq} that $ \risom{\beta} $ does not depend on the choice of $ \word{\omega} $. Thus $ \risom{\beta} $ is well-defined.
	
	Now we define an action of $ \twistgroup \times \invogroup $ on $ M_i $ by $ c.x \defl \risom{\alpha_i}^{-1}(c.\risom{\alpha_i}(x)) $ for all $ c \in \twistgroup \times \invogroup $ and $ x \in M_i $, so that $ \risom{\alpha_i}(c.x) = c.\risom{\alpha_i}(x) $. By definition, this action is compatible with the twisting action of $ \twistgroup $. Since $ \twistgroup \times \invogroup $ is a twisting group for $ G $, it follows for all $ \beta \in O_i $ that
	\begin{align*}
		\risom{\beta}(c.x) &= \risom{\alpha_i}(c.x)^{w_{\word{\delta}}} = \brackets[\big]{c.\risom{\alpha_i}(x)}^{w_{\word{\delta}}} = c.\risom{\alpha_i}(x)^{w_{\word{\delta}}} = c.\risom{\beta}(x)
	\end{align*}
	for all $ x \in M_i $ and $ c \in \twistgroup \times \invogroup $, where $ \word{\delta} $ is an arbitrary word over $ \rootbase \union (-\rootbase) $ such that $ \risom{\beta}(y) = \risom{\alpha_i}(y)^{w_{\word{\delta}}} $ for all $ y \in M_i $. Thus the action of $ \twistgroup \times \invogroup $ on $ M_i $ is compatible with all root isomorphisms $ (\risom{\beta})_{\beta \in O_i} $.
	
	Now let $ \beta \in O_i $, $ \delta \in \rootbase $ and $ x \in M_i $. We want to show that $ \risom{\beta}(x)^{w_{\delta}} = \risom{\gamma}(c.x) $ where $ \gamma \defl \beta^{\reflbr{\delta}} $ and $ c \defl \totalpar{\beta}{\delta} $. By its definition, $ c $ lies in $ \parmoveset{(\twistgroup \times \invogroup)}{\beta}{\gamma} $, so it also lies in $ \parmoveset{(\twistgroup \times \invogroup)}{\gamma}{\gamma} $ because $ \totalparsym $ is transporter-invariant. Thus there exists a word $ \word{\omega} $ over $ \rootbase \union (-\rootbase) $ such that
	\[ \totalparbr{\gamma}{\word{\omega}} = c \midand \gamma^{\reflbr{\word{\omega}}} = \gamma. \]
	Further, by \cref{param:transp-invar-1}, there also exists a word $ \word{\zeta} $ over $ \rootbase \union (-\rootbase) $ such that
	\[ \totalparbr{\alpha_i}{\word{\zeta}} = (1_\twistgroup, 1_\invogroup) \midand \alpha_i^{\reflbr{\word{\zeta}}} = \beta. \]
	By their respective definitions, these words satisfy
	\begin{equation}\label{eq:param:thm}
		c.\risom{\gamma}(x) = \risom{\gamma}(c.x) = \risom{\gamma}(x)^{w_{\word{\omega}}} \midand \risom{\alpha_i}(x)^{w_{\word{\zeta}}} = \risom{\beta}(x).
	\end{equation}
	Further
	\begin{align*}
		\totalparbr{\alpha_i}{\word{\zeta} \delta \word{\omega}^{-1}} &= \totalparbr{\alpha_i}{\word{\zeta}} \totalparbr{\beta}{\delta} \totalparbr{\gamma}{\word{\omega}^{-1}} = c \totalparbr{\gamma^{\reflbr{\word{\omega}}^{-1}}}{\word{\omega}}^{-1} = c \totalparbr{\gamma}{\word{\omega}}^{-1} \\
		&= cc^{-1} = (1_\twistgroup, 1_\invogroup).
	\end{align*}
	By the definition of $ \risom{\gamma} $, it follows that $ w \defl w_{\word{\zeta}} w_\delta w_{\word{\omega}}^{-1} $ satisfies $ \risom{\alpha_i}(x)^w = \risom{\gamma}(x) $. Therefore,
	\begin{align*}
		\risom{\gamma}(x)^{w_{\word{\omega}}} &= \risom{\alpha_i}(x)^{w_{\word{\zeta}} w_\delta}.
	\end{align*}
	Applying both equations in~\eqref{eq:param:thm}, we infer that
	\[ \risom{\gamma}(c.x) = \risom{\beta}(x)^{w_\delta}, \]
	which finishes the proof.
\end{proof}

\begin{remark}\label{param:thm-strengthen-rem}
	It is clear from the first lines of the proof that the statement of \cref{param:thm} can be slightly strengthened: For all $ i \in \numint{1}{k} $ choose a root $ \alpha_i \in O_i $, a group $ M_i $ which is isomorphic to $ \rootgr{\alpha_i} $ and an isomorphism $ \map{\risom{\alpha_i}}{M_i}{\rootgr{\alpha_i}}{}{} $. Then there exist an action of $ \twistgroup \times \invogroup $ on each of the underlying sets $ \listing{M}{k} $ and a family of isomorphisms $ (\risom{\alpha})_{\alpha \in \roots} $ extending $ (\risom{\alpha_i})_{i \in \numint{1}{k}} $ such that the assertions in \cref{param:thm} hold. In other words, whenever we can construct multiple parametrisations of a group (for different parity maps, for example), we can always arrange it so that they are parametrised by the same groups $ \listing{M}{k} $. This technical observation will be useful in the parametrisation of $ F_4 $-graded groups, in which we will parametrise certain $ B_3 $- and $ C_3 $-graded subgroups and then consider their \enquote{overlap}.
\end{remark}

The essential content of the following remark is that different choices of parity maps yield the same parametrisations up to twisting.

\begin{remark}\label{param:gpt-twist}
	We briefly consider what changes the effect of choosing another pair $ \inverparsym' $, $ \invoparsym' $ of $ \rootbase $-parity maps with values in $ \twistgroup $ and $ \invogroup $, respectively. Assume that $ (\twistgroup, \inverparsym', \invogroup, \invoparsym') $ satisfies the same assumptions as $ (\twistgroup, \inverparsym, \invogroup, \invoparsym) $ in \cref{param:partwist-conv}. Specifically, this means that $ (\twistgroup, \inverparsym', \invogroup, \invoparsym') $ is a partial twisting system for $ (G, (w_\delta)_{\delta \in \rootbase}) $ and that $ G $ is square-compatible with respect to $ \inverparsym' $ and stabiliser-compatible with respect to $ (\inverparsym', \invoparsym') $. Assume further that $ \genparmoveset{\invogroup}{\invoparsym}{\alpha}{\alpha} = \genparmoveset{\invogroup}{\invoparsym'}{\alpha}{\alpha} $ for all roots $ \alpha $, and choose the same complement $ C_\alpha $ of these groups in $ \invogroup $.
	
	The parametrisation theorem yields two parameter systems
	\[ \calP = (\twistgroup \times \invogroup, \listing{M}{k}) \midand \calP' = (\twistgroup \times \invogroup, M_1', \ldots, M_k') \]
	of type $ \roots $ and corresponding parametrisations $ (\risom{\alpha})_{\alpha \in \roots} $ and $ (\risom{\alpha}')_{\alpha \in \roots} $ of $ G $ with respect to $ (w_\delta)_{\delta \in \rootbase} $ and $ \inverparsym \times \invoparsym $ or $ \inverparsym' \times \invoparsym' $, respectively. For each $ i \in \numint{1}{k} $, we fix a root $ \alpha_i \in O_i $. Using \cref{param:thm-strengthen-rem}, we can assume that $ M_i = M_i' $ and $ \risom{\alpha_i} = \risom{\alpha_i}' $ for all $ i \in \numint{1}{k} $. However, observe that the parameter systems $ \calP $ and $ \calP' $ are not the same because the actions of $ \twistgroup \times \invogroup $ on $ \listing{M}{k} $ may be different. We solve this notational issue by writing $ \omega_i(c, x) $ for the action of $ \calP $ and $ \omega_i'(c,x) $ for the action of $ \calP' $, where $ i \in \numint{1}{k} $, $ x \in M_\alpha $ and $ c \in \twistgroup \times \invogroup $.
	
	We begin by comparing the different twisting actions. Fix $ i \in \numint{1}{k} $ for the scope of this paragraph. For all $ a \in \twistgroup $, we have
	\[ \risom{\alpha_i}'\brackets[\big]{\omega_i'(a,x)} = a.\risom{\alpha_i}(x) = \risom{\alpha_i}\brackets[\big]{\omega_i(a,x)} \]
	for all $ x \in M_i $ because the twisting action of $ \twistgroup $ on $ \rootgr{\alpha} $ is fixed. It follows that the two actions of $ \twistgroup $ on $ M_i $ are identical. Now let $ b \in \genparmoveset{\invogroup}{\invoparsym}{\alpha_i}{\alpha_i} = \genparmoveset{\invogroup}{\invoparsym'}{\alpha_i}{\alpha_i} $ and choose a word $ \word{\delta} $ over $ \rootbase \union (-\rootbase) $ such that $ \reflbr{\word{\delta}} $ stabilises $ \alpha_i $ and such that $ \inverinvopar{\alpha_i}{\word{\delta}}' = b $. Then
	\begin{align*}
		\risom{\alpha_i}'\brackets[\big]{\omega_i'(b, x)} &= \risom{\alpha_i}'(x)^{w_{\word{\delta}'}} = \risom{\alpha_i}(x)^{w_{\word{\delta}'}} = \risom{\alpha_i}\brackets[\big]{\omega_i(\inverinvopar{\alpha_i}{\word{\delta}'},x)} \\
		&= \risom{\alpha_i}'\brackets[\big]{\omega_i(\inverinvopar{\alpha_i}{\word{\delta}'},x)}
	\end{align*}
	for all $ x \in M_i $. We conclude that, while the actions $ \omega_i(b, \mapdot) $ and $ \omega_i'(b, \mapdot) $ for individual elements $ b $ might be distinct, each action $ \omega_i(b, \mapdot) $ can be expressed as $ \omega_i'(b', \mapdot) $ for some $ b' \in \invogroup $. In fact, one can show that the map $ \map{}{}{}{b}{b'} $ is an automorphism of the group $ \genparmoveset{\invogroup}{\invoparsym}{\alpha_i}{\alpha_i} $. By choosing the same complement $ C_\alpha $ for $ \genparmoveset{\invogroup}{\invoparsym}{\alpha_i}{\alpha_i} $ and $ \genparmoveset{\invogroup}{\invoparsym'}{\alpha_i}{\alpha_i} $ in \cref{param:invo-def}, we can extend $ \map{}{}{}{b}{b'} $ to an automorphism of $ \invogroup $ which still has the property that $ \omega_i(b, \mapdot) = \omega_i'(b', \mapdot) $ for all $ b \in \invogroup $.
	
	We now compare the different root isomorphisms. Again, let $ i \in \numint{1}{k} $ be arbitrary. Let $ \beta $ be any root in $ O_i $ and let $ \word{\delta} $ be a word over $ \rootbase \union (-\rootbase) $ such that $ \reflbr{\word{\delta}} $ maps $ \alpha_i $ to $ \beta $. Since $ \risom{\alpha_i} = \risom{\alpha_i}' $, we then have
	\begin{align*}
		\risom{\beta}'\brackets[\big]{\omega_i(\inverinvoparprime{\alpha}{\word{\delta}}, x)} &= \risom{\alpha_i}'(x)^{w_{\word{\delta}}} = \risom{\alpha_i}(x)^{w_{\word{\delta}}} = \risom{\beta}\brackets[\big]{\omega_i'(\inverinvopar{\alpha}{\word{\delta}}, x)}
	\end{align*}
	for all $ x \in M_i $. Using the conclusion of the previous paragraph, it follows that there exists $ c \in \twistgroup \times \invogroup $ such that
	\[ \risom{\beta}'(x) = \risom{\beta}\brackets[\big]{\omega_i(c,x)} \]
	for all $ x \in M_i $. We conclude that $ (\risom{\gamma}')_{\gamma \in \roots} $ can be obtained from $ (\risom{\gamma})_{\gamma \in \roots} $ by twisting (in the sense of \cref{param:parmap-twist}). Observe that $ c $ automatically has the property that $ \omega_i(c, \mapdot) $ is a group automorphism of $ M_i $ because both $ \risom{\beta} $ and $ \risom{\beta}' $ are homomorphisms.
	
	Finally, assume in addition that $ \calP $ is $ (\inverparsym \times \invoparsym) $-faithful and $ \calP' $ is $ (\inverparsym' \times \invoparsym') $-faithful. (By the previous observations, this is already true if it holds for one of the parameter systems.) Since then $ \inverparsym \times \invoparsym $ and $ \inverparsym' \times \invoparsym' $ are determined by their corresponding parametrisations $ (\risom{\gamma})_{\gamma \in \roots} $ and $ (\risom{\gamma}')_{\gamma \in \roots} $, it follows that $ \inverparsym' \times \invoparsym' $ can be obtained from $ \inverparsym \times \invoparsym $ by twisting as well.
\end{remark}

\begin{note}
	In all examples of root graded groups of rank at least $ 3 $, the group $ \invogroup $ is $ \compactSet{\pm 1} $ or trivial. It follows that the automorphism $ \map{}{\invogroup}{\invogroup}{b}{b'} $ in \cref{param:gpt-twist} must be the identity map. In an imaginary example where $ \invogroup $ is of the form $ \compactSet{\pm 1}^p $, it would however be possible that $ \map{}{}{}{b}{b'} $ interchanges some components of $ \invogroup $.
\end{note}

\begin{note}\label{param:parmap-choice-note}
	A consequence of \cref{param:gpt-twist} is that the specific choice of the parity maps $ \inverparsym $ and $ \invoparsym $ is ultimately not relevant. Thus there is no harm in choosing explicit parity maps $ \inverparsym_\roots $ and $ \invoparsym_\roots $ during the parametrisation process of $ \roots $-graded groups. The reader who prefers a different choice of parity maps $ \inverparsym_\roots' $ and $ \invoparsym_\roots' $ does not have to repeat all our proofs: Instead, they can simply take our parametrisation $ (\risom{\alpha})_{\alpha \in \roots} $ (which respects $ \inverparsym_\roots $ and $ \invoparsym_\roots $) and twist it appropriately to obtain a parametrisation $ (\risom{\alpha}')_{\alpha \in \roots} $ with respect to $ \inverparsym_\roots' $ and $ \invoparsym_\roots' $.
\end{note}

We end this section with some remarks on how the parametrisation procedure has been carried out in some special cases in the literature.

\begin{note}[The rank-2 case]
	Assume that $ \roots $ is irreducible crystallographic of rank~2 and choose a positive system $ \possys $ in $ \roots $. Recall that irreducible RGD-systems of type $ \roots $ were classified by Tits-Weiss in \cite{MoufangPolygons}. More precisely, they classify \defemph*{root group sequences}\index{root group sequence}. Here root group sequences are pairs $ (\rootgr{\possys}, (\rootgr{\alpha})_{\alpha \in \possys}) $ consisting of a group $ \rootgr{\possys} $ and a family $ (\rootgr{\alpha})_{\alpha \in \possys} $ satisfying some axioms. For any RGD-system $ (G, (\rootgr{\alpha})_{\alpha \in \roots}) $, the pair $ (\rootgr{\possys}, (\rootgr{\alpha})_{\alpha \in \possys}) $ is a root group sequence, and any RGD-system is essentially determined by its root group sequence. For more details, see (8.7), (8.10) and~(8.11) in \cite{MoufangPolygons}.
	
	By inspecting the classification of irreducible crystallographic root systems of rank~2, we see that $ \roots $ has at most 16 elements, so $ \possys $ has at most 8 elements. This low number of roots allows Tits-Weiss to parametrise the root groups \enquote{by hand}, without using a machinery similar to the one in the parametrisation theorem. The same is true for \cite[13.7]{Faulkner-NonAssocProj} in which Faulkner parametrises $ A_2 $-graded groups.
\end{note}

\begin{remark}[Parametrisations for simply-laced root systems]\label{param:motiv:simply-laced-shi}
	Assume that $ \roots $ is simply-laced. In this setting, the property of stabiliser compatibility can be shown to be trivial (see \cref{param:adj-implies-stab}). Using this simplification and his more restrictive definition of root gradings (\ref{shi-def}), Shi proves in \cite[(2.11)]{Shi1993} that the assertion of \cref{param:abstract-thm} holds for the specific Weyl elements $ \phi(\steinweyl{\alpha}(1_\comring)) $ in Axiom~\thmitemref{shi-def}{shi-def:weyl} and for $ \invoparsym $ being trivial. This allows him to parametrise the root groups by an abelian group $ (\ring, +) $. His proof uses in a significant way that these Weyl elements come from the homomorphic image of a Steinberg group. This makes his argument more concise but less general than our proof of the parametrisation theorem. In particular, the role of the parity map $ \inverparsym $ (and thus of the sign problem) is not prominent in Shi's work because any $ \roots $-graded group in Shi's sense comes, by definition, with a family $ \chevstr = (\chevstr_{\alpha, \beta})_{\alpha, \beta \in \roots} $ of Chevalley structure constants which determine the parity map.
\end{remark}


\section{Criteria for the Compatibility Conditions}

\label{sec:param:crit}

\begin{secnotation}
	We denote by $ \roots $ an arbitrary root system, by $ \rootbase $ a rescaled root base of $ \roots $, by $ G $ a group with a $ \roots $-pregrading $ (\rootgr{\alpha})_{\alpha \in \roots} $ and by $ (w_\delta)_{\delta \in \rootbase} $ a $ \rootbase $-system of Weyl elements in $ G $.
\end{secnotation}

In this section, we collect some criteria which will be used to prove the compatibility conditions in assumptions of the parametrisation theorem. We begin with a criterion for square compatibility. Surprisingly, it will be applicable in many situations, and similar ideas work in all remaining situations.

\begin{lemma}\label{param:square-formula-comp}
	Let $ \twistgroup \defl \compactSet{\pm 1}^n $ be a twisting group for $ (G, (w_\delta)_{\delta \in \rootbase}) $ for some $ n \in \Npos $ and assume that the first component of $ \twistgroup $ acts on all root groups by inversion. Let $ \alpha $ be any root and let $ \inverparsym $ be a $ \rootbase $-parity map with values in $ \twistgroup $. Assume that both $ \inverparsym $ and $ G $ satisfy the square formula for $ \alpha $. Then $ G $ is $ \alpha $-square compatible with respect to $ \inverparsym $ (and any $ \rootbase $-system of Weyl elements).
\end{lemma}
\begin{proof}
	Let $ \delta \in \rootbase $, set $ w \defl w_\delta^2 $ for some $ \delta $-Weyl element $ w_\delta $ and put $ \epsilon \defl (-1)^{\cartanint{\alpha}{\delta}} $. Then $ x_\alpha^w = x_\alpha^\epsilon $ for all $ x_\alpha \in \rootgr{\alpha} $ because $ G $ satisfies the square formula for $ \alpha $ and $ \inverpar{\alpha}{\delta \delta} = (\epsilon, 1, \ldots, 1) $ because $ \inverparsym $ satisfies the square formula for $ \alpha $. It follows that $ w $ acts on $ \rootgr{\alpha} $ by $ \inverpar{\alpha}{\delta \delta} $, as desired.
\end{proof}

The following result is the main criterion for stabiliser compatibility. It is very general and can be applied for all root systems. However, it is unnecessarily complicated in many situations, which is why we state a simpler, less general version in \cref{param:adj-implies-stab}. We use \cref{rootorder:cry-brackets-note} to formulate both assertions.

\begin{proposition}\label{param:stabcomp-crit-ortho}
	Assume that $ G $ has (crystallographic) $ \roots $-commutator relations with root groups $ (\rootgr{\alpha})_{\alpha \in \roots} $. Let $ \twistgroup $ be a twisting group for $ G $, let $ \inverparsym $ be a (crystallographically) adjacency-trivial parity map with values in $ \twistgroup $ such that $ G $ is square-compatible with respect to $ \inverparsym $ and $ (w_\delta)_{\delta \in \rootbase} $ and let $ \invoparsym $ be a Weyl-invariant parity map with values in the abelian group $ \invogroup \defl \compactSet{\pm 1_\invogroup} $. Fix a root $ \alpha $ and define the following sets:
	\begin{align*}
		\calO \defl{}& \Set{\beta \in \roots \given \alpha \cdot \beta = 0}, \\
		\calA \defl{}& \Set{\beta \in \calO \given \alpha \text{ is (crystallographically) adjacent to } \beta \text{ and } -\beta}, \\
		={}& \Set{\beta \in \calO \given \alpha \text{ is (crystallographically) adjacent to } \beta} \\
		\bar{\calA} \defl{}& \calO \setminus \calA.
	\end{align*}
	Assume that we have $ \invopar{\alpha}{\reflbr{\beta}} = 1_\invogroup $ for all $ \beta \in \calA $ and $ \invopar{\alpha}{\reflbr{\beta}} = -1_\invogroup $ for all $ \beta \in \bar{\calA} $. Assume further that for all $ \rootbase $-positive roots $ \beta, \beta' \in \bar{\calA} $, there exist words $ \word{\delta}, \word{\delta}' $ over $ \rootbase \union (-\rootbase) $ such that $ \reflbr{\word{\delta}} = \reflbr{\beta} $, $ \reflbr{\word{\delta}'} = \reflbr{\beta'} $, $ \inverpar{\alpha}{\word{\delta}} = \inverpar{\alpha}{\word{\delta}'} $ and such that $ w_{\word{\delta}} $ and $ w_{\word{\delta}'} $ act identically on $ \rootgr{\alpha} $. Then $ G $ is $ \alpha $-stabiliser-compatible with respect to $ (\inverparsym, \invoparsym) $ and $ (w_\delta)_{\delta \in \rootbase} $.
\end{proposition}
\begin{proof}
	Let $ u \in \Weyl(\roots) $ such that $ \alpha^u = \alpha $ and $ \invopar{\alpha}{u} = 1_\invogroup $. By \cref{weyl:stab-ortho}, there exist $ \listing{\rho}{k} \in \calO $ such that $ u = \refl{\rho_1} \cdots \refl{\rho_k} $. Since $ \refl{\rho_i} = \refl{-\rho_i} $ for all $ i \in \numint{1}{n} $, we can choose $ \listing{\rho}{k} $ so that all of them are $ \rootbase $-positive. At first, we only consider the special case that all roots $ \listing{\rho}{k} $ lie in $ \bar{\calA} $. We will show that for each $ i \in \numint{1}{k} $, there exists a word $ \word{\delta}^i $ over $ \rootbase \union (-\rootbase) $ such that $ \reflbr{\word{\delta}^i} = \reflbr{\rho_i} $ and such that $ w_{\word{\delta}^1 \cdots \word{\delta}^k} $ acts on $ \rootgr{\alpha} $ by $ \inverpar{\alpha}{\word{\delta}^1 \cdots \word{\delta}^k} $. By the assumption on the map $ \invoparsym $, we have
	\begin{align*}
		1_\invogroup &= \invopar{\alpha}{u} = \invopar{\alpha}{\reflbr{\rho_1} \cdots \reflbr{\rho_k}} = \prod_{i=1}^k \invopar{\alpha}{\reflbr{\rho_i}} = (-1_\invogroup)^k,
	\end{align*}
	which implies that $ k $ is even. For all $ i \in \numint{1}{k/2} $, we can now by assumption choose words $ \word{\delta}^{2i-1} $, $ \word{\delta}^{2i} $ over $ \rootbase \union (-\rootbase) $ such that
	\[ \reflbr{\word{\delta}^{2i-1}} = \reflbr{\rho_{2i-1}}, \qquad \reflbr{\word{\delta}^{2i}} = \reflbr{\rho_{2i}}, \qquad \inverpar{\alpha}{\word{\delta}^{2i-1}} = \inverpar{\alpha}{\word{\delta}^{2i}} \]
	and such that $ w_{\word{\delta}^{2i-1}} $ and $ w_{\word{\delta}^{2i}} $ act identically on $ \rootgr{\alpha} $. Then $ w_{\word{\delta}^1 \cdots \word{\delta}^k} $ and $ \prod_{i=1}^{k/2} w_{\word{\delta}^{2i}}^2 $ act identically on $ \rootgr{\alpha} $. Since $ G $ is square-compatible with respect to $ \inverparsym $, it follows that $ w_{\word{\delta}^1 \cdots \word{\delta}^k} $ acts on $ \rootgr{\alpha} $ by $ \prod_{i=1}^{k/2} \inverpar{\alpha}{\word{\delta}^{2i} \word{\delta}^{2i}} $. Since $ \inverpar{\alpha}{\word{\delta}^{2i-1}} = \inverpar{\alpha}{\word{\delta}^{2i}} $ and $ \refl{\word{\delta}^{2i}}(\alpha) = \alpha = \refl{\word{\delta}^{2i-1}}(\alpha) $ for all $ i \in \numint{1}{k/2} $, we also have that
	\begin{align*}
		\prod_{i=1}^{k/2} \inverpar{\alpha}{\word{\delta}^{2i} \word{\delta}^{2i}} &= \prod_{i=1}^{k/2} \inverpar{\alpha}{\word{\delta}^{2i}} \inverpar{\alpha}{\word{\delta}^{2i}} = \prod_{i=1}^{k/2} \inverpar{\alpha}{\word{\delta}^{2i-1}} \inverpar{\alpha}{\word{\delta}^{2i}} = \prod_{j=1}^k \inverpar{\alpha}{\word{\delta}^j} = \inverpar{\alpha}{\word{\delta}^1 \cdots \word{\delta}^k}.
	\end{align*}
	It follows that $ w_{\word{\delta}^1 \cdots \word{\delta}^k} $ acts on $ \rootgr{\alpha} $ by $ \inverpar{\alpha}{\word{\delta}^1 \cdots \word{\delta}^k} $, as desired.
	
	Now we consider the general case. Let $ i(1), \ldots, i(n), j(1), \ldots, j(m) $ be pairwise distinct indices from $ \numint{1}{k} $ such that
	\begin{align*}
		\Set{p \in \numint{1}{k} \given \rho_p \in \calA} &= \Set{i(1), \ldots, i(n)}, && i(1) < \cdots < i(n), \\
		\Set{p \in \numint{1}{k} \given \rho_p \in \bar{\calA}} &= \Set{j(1), \ldots, j(m)}, && j(1) < \cdots < j(m).
	\end{align*}
	Put $ u' \defl \reflbr{\rho_{j(1)}} \cdots \reflbr{\rho_{j(m)}} $. By the assumption on $ \invoparsym $, we have $ \invopar{\alpha}{\reflbr{\rho_{i(p)}}} = 1_\invogroup $ for all $ p \in \numint{1}{n} $, and thus $ \invopar{\alpha}{u'} = \invopar{\alpha}{u} = 1_\invogroup $. Hence by the conclusion of the previous paragraph, we can for each $ p \in \numint{1}{m} $ find a word $ \word{\delta}^p $ over $ \rootbase \union (-\rootbase) $ such that $ \reflbr{\word{\delta}^p} = \reflbr{\rho_{j(p)}} $ and such that $ w_{\word{\delta}^1 \cdots \word{\delta}^m} $ acts on $ \rootgr{\alpha} $ by $ \inverpar{\alpha}{\word{\delta}^1 \cdots \word{\delta}^m} $. For every $ p \in \numint{1}{n} $, we choose an arbitrary $ \rootbase $-expression $ \word{\zeta}^p $ of $ \rho_{i(p)} $. For each $ q \in \numint{1}{k} $, we can now define $ \word{\xi}^q \defl \word{\delta}^{p} $ if $ q = j(p) $ for some $ p \in \numint{1}{m} $ and $ \word{\xi}^q \defl \word{\zeta}^p $ if $ q = i(p) $ for some $ p \in \numint{1}{n} $. Note that for each $ p \in \numint{1}{n} $, we have $ \inverpar{\alpha}{\word{\zeta}^p} = 1_\twistgroup $ because $ \inverparsym $ is (crystallographically) adjacency-trivial and we also have that $ w_{\word{\zeta}^p} $ acts trivially on $ \rootgr{\alpha} $ because it is a $ \rho_{i(p)} $-Weyl element and $ \rho_{i(p)} $ lies in $ \calA $. This implies that $ w_{\word{\xi}^1 \cdots \word{\xi}^k} $ and $ w_{\word{\delta}^1 \cdots \word{\delta}^m} $ act identically on $ \rootgr{\alpha} $ and that $ \inverpar{\alpha}{\word{\xi}^1 \cdots \word{\xi}^k} = \inverpar{\alpha}{\word{\delta}^1 \cdots \word{\delta}^m} $. We conclude that $ w_{\word{\xi}^1 \cdots \word{\xi}^k} $ acts on $ \rootgr{\alpha} $ by $ \inverpar{\alpha}{\word{\xi}^1 \cdots \word{\xi}^k} $, as desired.
\end{proof}

In the following result, we consider the special case that $ \bar{\calA} $ is empty. In this situation, we can easily drop the assumptions on the parity map $ \invoparsym $. Again, we use \cref{rootorder:cry-brackets-note}.

\begin{lemma}\label{param:adj-implies-stab}
	Assume that $ G $ has (crystallographic) $ \roots $-commutator relations with root groups $ (\rootgr{\alpha})_{\alpha \in \roots} $. Let $ \inverparsym $ be a $ \rootbase $-parity map with values in some twisting group $ \twistgroup $ for $ (G, (w_\delta)_{\delta \in \rootbase}) $ and let $ \alpha \in \roots $. Put
	\[ \rootbase^W \defl \Set{\delta^u \given \delta \in \rootbase, u \in \Weyl(\roots)}. \]
	Assume that $ \inverparsym $ is (crystallographically) $ \alpha $-adjacency-trivial and that any root in $ \rootbase^W $ which is orthogonal to $ \alpha $ is also (crystallographically) adjacent to $ \alpha $. Then $ G $ is $ \alpha $-stabiliser-compatible with respect to~$ \inverparsym $.
\end{lemma}
\begin{proof}
	Let $ u $ be an element of the Weyl group of $ \roots $ such that $ \alpha^u = \alpha $. Then by \cref{weyl:stab-ortho}, there exist roots $ \beta_1, \ldots, \beta_m \in \roots $ which are all orthogonal to $ \alpha $ such that $ u = \refl{\beta_1} \cdots \refl{\beta_m} $. For the rest of this paragraph, we fix an arbitrary $ i \in \numint{1}{m} $. Applying \cref{rootsys:any-in-rootbase}, we can find a $ \rootbase $-expression $ \word{\rho}^i $ of $ \lambda \beta^i $ for some $ \lambda \in \IR_{>0} $. By \cref{param:Delta-exp-weyl}, $ w_{\word{\rho}^i} $ is a $ \lambda \beta_i $-Weyl element. Since $ \lambda \beta^i $ has a $ \rootbase $-expression, it lies in $ \rootbase^W $, and it is orthogonal to $ \alpha $ because $ \beta^i $ is. The same assertions hold for $ -\lambda \beta^i $. Thus it follows from our assumptions that $ \alpha $ is (crystallographically) adjacent to $ \lambda\beta_i $ and $ -\lambda\beta_i $. Hence $ w_{\word{\rho}^i} $ centralises $ \rootgr{\alpha} $. Further, $ \inverpar{\alpha}{\word{\rho}^i} = 1_\twistgroup $ because $ \inverparsym $ is $ \alpha $-adjacency-trivial. We conclude that $ w_{\word{\rho}^i} $ acts on $ \rootgr{\alpha} $ by $ \inverpar{\alpha}{\word{\rho}^i} $.
	
	It follows from the conclusion of the previous paragraph that $ w_{\word{\rho}} $ acts on $ \rootgr{\alpha} $ by $ \inverpar{\alpha}{\word{\rho}} $ where $ \word{\rho} \defl (\word{\rho}^1, \ldots, \word{\rho}^m) $. Since $ u = \refl{\beta_1} \cdots \refl{\beta_m} = \refl{\word{\rho}^1} \cdots \refl{\word{\rho}^m} = \refl{\word{\rho}} $, $ \word{\rho} $ is a representation of $ u $. The assertion follows.
\end{proof}

\begin{remark}\label{param:adj-implies-stab:length-rem}
	If $ \roots $ is reduced, then the set $ \rootbase^W $ in \cref{param:adj-implies-stab} is the whole root system by \cref{rootsys:any-in-rootbase}. In contrast, assume now that $ \roots $ is the non-reduced root system $ BC_n $ in standard representation for some $ n \ge 2 $ and that $ \rootbase $ is the standard rescaled root base of $ BC_n $, as in \cref{BC:BCn-standard-rep}. Then $ \rootbase^W $ consists of all roots which are long or of medium length. Let $ \alpha $ be an short root. Then $ \alpha $ does indeed have the property that any root $ \beta $ in $ \rootbase^W $ which is orthogonal to $ \alpha $ must be crystallographically adjacent to $ \alpha $. This conclusion does not hold if $ \beta $ is short: For example, the roots $ \alpha \defl \basvec_1 $ and $ \beta \defl \basvec_2 $ are orthogonal but not crystallographically adjacent.
\end{remark}


\section{Remarks on the General Strategy}

\label{sec:param:outline}

As emphasised before, the parametrisation theorem allows us to approach $ \roots $-graded groups for different root systems $ \roots $ in a uniform way. For this reason, \cref{chap:simply-laced,chap:B,chap:BC-alg,chap:BC,chap:F} all follow a common outline. In this section, we describe this outline. We always denote by $ \roots $ the specific root system which is studied in the corresponding chapter (or one of these root systems if there are multiple).

\begin{miscthm}[The coordinatising algebraic structures]\label{param:strategy:alg}
	Each of the aforementioned chapters begins with an investigation of the algebraic structures that coordinatise $ \roots $-graded groups. Most of the time, these are well-known structures with a highly developed underlying theory. We confine ourselves to covering the basics of these theories, though we will sometimes go beyond what is strictly necessary for our purposes. In the case of root gradings of type $ (B)C $, the algebraic structures in question are so involved that they are covered in a separate chapter.
	
	For each algebraic structure, we define a \defemph*{standard parameter system}\index{parameter system!standard}, which is a parameter system in the sense of \cref{param:parsys-def}. We will show that every $ \roots $-graded group can be parametrised by a parameter system of this form.
\end{miscthm}

\begin{miscthm}[The root system]\label{param:strategy:rootsys}
	As a next step, we study some (either known or basic) combinatorial properties of $ \roots $. Most importantly, we compute the Cartan integers, which will later allow us to prove the square formula for Weyl elements in $ \roots $-graded groups (in most cases). Further, we will prove some minor lemmas concerning the relationship between arbitrary roots $ \alpha $ and $ \beta $.
\end{miscthm}

\begin{miscthm}[Construction of a generic example]\label{param:strategy:ex}
	Before we turn to the investigation of arbitrary $ \roots $-graded groups, we construct examples of such groups. In these examples, we will see the commutator relations that we ultimately want to establish in every $ \roots $-graded group. Further, we will also see the standard parameter system from~\ref{param:strategy:alg} in action.
	
	Ideally, we would construct a $ \roots $-graded group for any object in the class of algebraic structures which were studied in~\ref{param:strategy:alg}. However, we will restrict ourselves to certain special cases whenever a general construction would be much more complicated. In some cases, it is even true that no general construction is known. Yet, the groups that we construct will be sufficiently generic (in the sense of \cref{param:motiv:parmap-choice}) to yield adequate parity maps. These parity maps will later be used in the parametrisation of $ \roots $-graded groups.
\end{miscthm}

\begin{miscthm}[Computations in $ \roots $-graded groups]\label{param:strategy:comp}
	As a next step, we study the action of Weyl elements in $ \roots $-graded groups on the root groups, independently of any choice of a parity map. These results form the mathematical core of the study of $ \roots $-graded groups for any specific root system $ \roots $. Our goal is to obtain some fundamental results which will later be used to prove the compatibility conditions which appear in the parametrisation theorem. This means that we have to investigate the actions of squares of Weyl elements on all root groups and the action of $ \beta $-Weyl elements on any root group which is orthogonal to $ \beta $.
	
	In order to prove the aforementioned statements, we need two computational tools. At first, we investigate how the commutator maps $ \commpart{\mapdot}{\mapdot}{\alpha} $ behave under products in the first and second argument. We will see that these maps are always additive or, in a loose sense, \enquote{quadratic} in each argument, which is precisely the behaviour that we have seen in the commutator relations in~\ref{param:strategy:ex}. Secondly, we need the formulas from \cref{tw:6.4}. Using the identities from the previous step, we will give alternative, more explicit proofs of these formulas for each root system.
	
	With these tools at hand, we can approach the problems from the first paragraph. While all computations in the previous paragraph take place in root graded subgroups of rank~2, we will now have to use the rank-3 assumptions in a few places. However, some partial results still hold in the rank-2 situation, and we will prove as much as we can in this generality.
\end{miscthm}

\begin{miscthm}[Standard signs in $ \roots $-graded groups]\label{param:strategy:stsigns}
	We say that a $ \roots $-pregrading is \defemph*{coordinatised with standard signs} if there exist an algebraic structure as in~\ref{param:strategy:alg} and root isomorphisms which satisfy the precise commutator relations of the example in~\ref{param:strategy:ex} (or a suitable generalisation of these relations). The terminology is due to the fact that the signs in these relations are not canonically determined (see also \cref{param:parmap-twist}), so we (somewhat arbitrarily) declare that the signs which occur in~\ref{param:strategy:ex} are the \enquote{standard ones}. Using the parity maps from~\ref{param:strategy:ex}, our goal is to show that every $ \roots $-graded group is coordinatised with standard signs.
	
	Note that the notion of coordinatisations with standard signs is defined for arbitrary pregradings. It is clear that any pregrading which is coordinatised in this way automatically has (crystallographic) $ \roots $-commutator relations, but the existence of Weyl elements is not obvious. However, we can use our specific knowledge of the commutator relations to show that any element of the coordinatising structure which is invertible (in a suitable sense) induces a Weyl element, and that every Weyl element is of this form. Thus any $ \roots $-pregrading which has a coordinatisation with standard signs is automatically a $ \roots $-grading, except that the validity of Axiom~\thmitemref{rgg-def}{rgg-def:nondeg} is not clear.
\end{miscthm}

\begin{note}[On the choice of a parity map]\label{param:strategy:parmap-choice}
	After working with parity maps in a purely abstract way for the whole chapter, it might seem disappointing that we ultimately prove our coordinatisation theorems only for a specific choice of standard signs and not in a more abstract way. However, observe that as soon as we have established the existence of one coordinatisation $ \theta = (\risom{\alpha})_{\alpha \in \roots} $ (with standard signs), we can easily twist $ \theta $ (as in \cref{param:parmap-twist}) to obtain coordinatisations of the same group with non-standard signs. By \cref{param:gpt-twist}, essentially every coordinatisation can be obtained in this way. Thus it is not a serious restriction to only consider standard signs. See also \cref{param:parmap-choice-note}.
\end{note}

\begin{miscthm}[Admissible partial twisting systems and the parametrisation of $ \roots $-graded groups]\label{param:strategy:admit-twist}
	Recall that the parametrisation theorem starts from a partial twisting system $ (\twistgroup, \inverparsym, \invogroup, \invoparsym) $ and a $ \roots $-graded group $ G $ which satisfies some compatibility conditions. It is clear that these compatibility conditions impose some restrictions on the partial twisting system, so not every imaginable partial twisting system can be used to parametrise $ \roots $-graded groups. A \defemph*{$ \roots $-admissible partial twisting system}\index{twisting system!partial!admissible} is, by definition, an admissible partial twisting which satisfies some additional conditions which depend on $ \roots $. These conditions are, unlike the compatibility conditions, easy to verify if the parity maps $ \inverparsym $ and $ \invoparsym $ are explicitly given. 
	
	The \defemph*{standard partial twisting system (of type $ \roots $)}\index{twisting system!partial!standard} is defined to be the partial twisting system $ (\twistgroup, \inverparsym, \invogroup, \invoparsym) $ that we see in~\ref{param:strategy:ex}, except that we \enquote{forget} the actions of $ \invogroup $ on the root groups. We will show that the standard partial twisting system of type $ \roots $ is $ \roots $-admissible. Further, we will prove that every $ \roots $-graded group $ G $ satisfies the compatibility conditions with respect to any $ \roots $-admissible partial twisting system. As a consequence, we can apply the parametrisation theorem to $ G $ and its standard partial twisting system, which yields a parametrisation of $ G $.
\end{miscthm}

\begin{note}[Standard signs and twisting systems for simply-laced groups]
	Our approach to standard signs and standard (partial) twisting systems for the simply-laced root graded groups will be slightly different than described in~\ref{param:strategy:stsigns} and~\ref{param:strategy:admit-twist} because we want to cover the types $ A $, $ D $ and $ E $ at the same time. We will elaborate on this in \cref{sec:ADE:sttwist,sec:ADE:stsigns}.
\end{note}

\begin{miscthm}[Blueprint computations]
	As a final step, we apply the blueprint technique to compute the commutator relations with respect to the chosen parametrisation. We will elaborate on this in \cref{chap:blue}. For these computations, it will be crucial that the signs in our parametrisation are explicitly given. We emphasise that there is nothing specific about the standard signs that makes the blueprint computation work: The only relevant point is that we make \emph{some} explicit choice of signs.
\end{miscthm}

	\chapter{Root Gradings of Simply-laced Type}
	
	\label{chap:simply-laced}
	
	In this chapter, we begin our investigation of $ \roots $-graded groups for specific root systems $ \roots $ with the case that $ \roots $ is simply-laced (that is, of type $ A $, $ D $ or $ E $) and of rank at least 2. The main result of this section is \cref{ADE:thm}: Every $ \roots $-graded group is coordinatised by a ring which must be associative if $ \roots $ is of rank at least~3 and commutative if $ \roots $ is of type $ D $ or $ E $. It turns out that, after our hard work to establish the parametrisation theorem, there is actually not much left to do to prove this. We can even carry out the complete coordinatisation of these groups without using the blueprint technique.
	
	The coordinatisation result of this section has already been proven by Shi in \cite{Shi1993}, but for a more restrictive definition of root gradings. (See also \cref{sec:literature}.) Hence our results are more general than Shi's, though many of the arguments are the same. The main additional difficulty in our generality is the sign problem in the sense of \cref{param:motiv:parmap-choice}, which we solve using the parametrisation theorem. A corollary of our coordinatisation result is that every simply laced root graded group in our sense is also a root graded group in Shi's sense. Thus the two notions of root gradings are actually equivalent, but this is a non-trivial fact.
	
	This chapter is organised in the way described in \cref{sec:param:outline}. We begin with a brief introduction to nonassociative rings. In the following section, we study some purely combinatorial properties of the simply-laced root systems and compute the Cartan integers for these root systems. In \cref{sec:ADE:ex}, we discuss the existence problem for simply-laced root graded groups. The main mathematical work of this chapter happens in \cref{sec:ADE:weyl}: We prove that $ \roots $-graded groups satisfy the square formula for Weyl elements. In \cref{sec:ADE:sttwist,sec:ADE:stsigns}, we introduce standard twisting systems and standard signs for root graded groups, respectively. Further, we show that every ring which coordinatises a $ \roots $-graded group (with standard signs) must be associative if $ \roots $ is of rank at least~3 and commutative if $ \roots $ is, in addition, of type $ D $ or $ E $. In \cref{sec:ADE:param}, we carry out the coordinatisation of simply-laced root graded groups.

\newcommand{\leftmul}[1]{l_{#1}}
\newcommand{\rightmul}[1]{r_{#1}}

\section{Nonassociative Rings}

\label{sec:ring}

One of our first announcements in this book was the frightening \cref{pre:ring-nonassoc-conv} that rings are not assumed to be associative. We will sometimes call these objects \enquote{nonassociative rings} to emphasise this convention. In this section, we introduce the most basic notions which are relevant in the theory of nonassociative rings. Technically, only the definition of these objects is needed in this chapter, but it seems prudent to go a bit beyond that. Some of the results in this section will only be needed in later chapters.

Of particular interest to us are rings which satisfy the so-called alternative laws, and which are called alternative rings. While we can merely show that every $ A_2 $-graded group is coordinatised by a nonassociative ring with no further properties, all known examples of rings which appear as coordinatising rings of $ A_2 $-graded groups are in fact alternative. See \cref{sec:ADE:ex} for a few more details. However, since we do not need the notion of alternative rings until we investigate root graded groups of type $ (B)C $ in \cref{chap:BC}, we delay their introduction until \cref{sec:BC:altring}.

A standard reference for most of the material in this section is \cite{Schafer}.

\begin{definition}[Ring]\label{ring-def}
	A \defemph{ring} is a triple $ (\ring, +, \rmult) $ consisting of a set $ \ring $ and two binary operations
	\[ \map{+, \rmult}{\ring \times \ring}{\ring}{}{} \]
	such that the following conditions are satisfied:
	\begin{stenumerate}
		\item $ (\ring, +) $ is an abelian group.
		
		\item \label{ring-def:id}There exists an element $ 1_\ring \in \ring $ such that $ 1_\ring \rmult x = x = x \rmult 1_\ring $ for all $ x \in \ring $.
		
		\item The distributive laws are satisfied. That is,
		\begin{align*}
			(x+y) \rmult z = (x \rmult z) + (y \rmult z) \midand x \rmult (y+z) = (x \rmult y) + (x \rmult z)
		\end{align*}
		for all $ x,y,z \in \ring $.
	\end{stenumerate}
	A \defemph*{subring of $ \ring $}\index{ring!sub-} is a subset of $ \ring $ which is closed under addition, additive inversion and multiplication and which contains $ 1_\ring $. An \defemph*{ideal of $ \ring $}\index{ideal (in a ring)} is a subgroup $ I $ of $ (\ring, +) $ such that $ ax \in I $ and $ xa \in I $ for all $ a \in \ring $ and $ x \in I $.
\end{definition}

\begin{definition}[Homomorphism of rings]
	Let $ \ring $, $ \varring $ be two rings. A map $ \map{f}{\ring}{\varring}{}{} $ is called a \defemph*{homomorphism of rings} if it preserves addition and multiplication and $ f(1_\ring) = 1_\varring $.
\end{definition}

\begin{secnotation}
	For the rest of this section, $ (\ring, +, \rmult) $ is a ring.
\end{secnotation}

\begin{notation}
	The neutral element of the group $ (\ring, +) $ will always be denoted by $ 0_\ring $ and the identity element from Axiom~\thmitemref{ring-def}{ring-def:id} will always be denoted by $ 1_\ring $. Further, we will often leave out the multiplication dot.
\end{notation}

\begin{remark}
	For any ideal $ I $ of $ \ring $, the quotient $ \ring / I $ has a canonical ring structure, as in the associative setting.
\end{remark}

The standard parameter system for a ring can be defined immediately.

\begin{definition}[Standard parameter system]
	Put $ \twistgroup \defl \compactSet{\pm 1} $ and declare that $ \twistgroup $ acts on $ (\ring, +) $ by inversion. That is, we put $ 1_\twistgroup.r \defl r $ and $ -1_\twistgroup.r \defl -r $. Then the pair $ (\twistgroup, \ring) $ is called the \defemph*{standard parameter system for $ \ring $}.\index{parameter system!standard!simply-laced type}
\end{definition}

It is a well-established fact that the commutator is a useful tool for studying noncommutative structures. Naturally, we can apply the same idea to nonassociative structures.

\begin{definition}[Associator]
	The \defemph*{associator map} or simply the \defemph{associator} is
	\[ \map{}{\ring \times \ring \times \ring}{\ring}{(a,b,c)}{\assoc{a}{b}{c} \defl (ab)c - a(bc).} \]
	A ring element of the form $ \assoc{x}{y}{z} $ for some $ x,y,z \in \ring $ will also be called an \defemph*{associator}. For any subsets $ A,B,C $ of $ \ring $, we denote by $ \assoc{A}{B}{C} $ the ideal which is generated by $ \Set{\assoc{a}{b}{c} \given a \in A, b \in B, c \in C} $.
\end{definition}

By the distributive laws, the associator map is additive in each component. Many of the additional properties which can be imposed on rings, such as the alternative laws, can be phrased in terms of vanishing of certain associators. We will get back to this in \cref{ring:alternative-def}.

\begin{remark}
	Clearly, $ \ring / \assoc{\ring}{\ring}{\ring} $ is an associative ring, and any ideal $ I $ for which $ \ring/I $ is associative contains $ \assoc{\ring}{\ring}{\ring} $.
\end{remark}

It is well-known that any associative ring has a center, which is itself an associative commutative ring. We can define a similar object for nonassociative rings, called the nucleus, but here the definition is a bit more delicate because we have to distinguish between left, middle and right nucleus. We will see in \cref{ring:alternative-nucleus} that this distinction is superfluous for alternative rings.

\begin{definition}[Nucleus and center]\label{ring:nucleus-def}
	We define the following subsets of $ \ring $, called the \defemph*{left nucleus}, \defemph*{middle nucleus} and \defemph*{right nucleus:}\index{nucleus!left}\index{nucleus!middle}\index{nucleus!right}
	\begin{align*}
		\lnucleus(\ring) &\defl \Set{x \in \ring \given \assoc{x}{y}{z} = 0 \text{ for all } y,z \in \ring}, \\
		\mnucleus(\ring) &\defl \Set{y \in \ring \given \assoc{x}{y}{z} = 0 \text{ for all } x,z \in \ring}, \\
		\rnucleus(\ring) &\defl \Set{z \in \ring \given \assoc{x}{y}{z} = 0 \text{ for all } x,y \in \ring}.
	\end{align*}
	Further, we define the \defemph*{nucleus of $ \ring $}\index{nucleus} by
	\[ \nucleus(\ring) \defl \lnucleus(\ring) \intersect \mnucleus(\ring) \intersect \rnucleus(\ring). \]
	Finally, the \defemph*{center of $ \ring $}\index{center} is
	\[ \zentrum(\ring) \defl \Set{x \in \nucleus(\ring) \given xy = yx \text{ for all } x,y \in \ring} \subs \nucleus(\ring). \]
\end{definition}

\begin{remark}[{\cite[{\onpage{13}}]{Schafer}}]\label{ring:nucleus-closure}
	By the distributive laws, all the nuclei defined in~\ref{ring:nucleus-def} are closed under addition and additive inverse. Less obviously, they are also closed under multiplication: Let $ x,x', y, z \in \ring $. If $ x,x' \in \lnucleus{\ring} $, then
	\[ \brackets[\big]{(xx')y} z = \brackets[\big]{x(x'y)}z = x \brackets[\big]{(x'y)z} = x \brackets[\big]{x' (yz)} = (xx') (yz), \]
	so $ xx' \in \lnucleus(\ring) $. If $ x,x' \in \mnucleus(\ring) $, then
	\[ \brackets[\big]{y(xx')}z = \brackets[\big]{(yx)x'}z = (yx) (x'z) = y \brackets[\big]{x(x'z)} = y \brackets[\big]{(xx')z}, \]
	so $ xx' \in \mnucleus(\ring) $. If $ x,x' \in \rnucleus(\ring) $, then
	\[ (yz) (xx') = \brackets[\big]{(yz)x}x' = \brackets[\big]{y(zx)}x' = y\brackets[\big]{(zx)x'} \]
	so $ xx' \in \rnucleus(\ring) $. Further, the center of $ \ring $ has the same properties because
	\[ (xx')y = x(x'y) = x(yx') = (xy)x' = (yx)x' = y(xx') \]
	for all $ x,x' \in \zentrum(\ring) $ and $ y \in \ring $. Since $ 1_\ring $ is clearly contained in the center of $ \ring $, it follows that all the sets defined in~\ref{ring:nucleus-def} are subrings of $ \ring $.
\end{remark}

\begin{remark}[{\cite[21.2.1]{McCrimmon_TasteJordan}}]\label{ring:nucl-slip}
	Since nuclear elements, in the words of McCrimmon, \enquote{slip in and out of parentheses}, we have
	\begin{align*}
		n \assoc{x}{y}{z} &= \assoc{nx}{y}{z}, & \assoc{xn}{y}{z} &= \assoc{x}{ny}{z}, \\
		\assoc{x}{yn}{z} &= \assoc{x}{y}{nz}, & \assoc{x}{y}{zn} &= \assoc{x}{y}{z}n
	\end{align*}
	for all $ x,y,z \in \ring $ and all $ n \in \nucleus(\ring) $. For example,
	\begin{align*}
		n \assoc{x}{y}{z} &= n \brackets[\big]{(xy)z} - n \brackets[\big]{x(yz)} = \brackets[\big]{n(xy)}z - (nx)(yz) \\
		&= \brackets[\big]{(nx)y}z - (nx)(yz) = \assoc{nx}{y}{z}.
	\end{align*}
	In \cref{altring:nucl-slip}, we will see stronger versions of these formulas for alternative rings .
\end{remark}

We now turn to invertibility in nonassociative rings, which is a slightly subtle concept.

\begin{definition}[Invertibility, {\cite[\onpage{190}]{Faulkner-NonAssocProj}}]\label{ring:def-invertible}
	For all $ x \in \ring $, we define \defemph*{left and right multiplication maps} $ \map{\leftmul{x}, \rightmul{x}}{\ring}{\ring}{}{} $ by $ \leftmul{x}(y) \defl x \rmult y $ and $ \rightmul{x}(y) \defl y \rmult x $ for all $ y \in \ring $. For any $ x,y \in \ring $, we say that $ y $ is an \defemph*{inverse of $ x $}\index{inverse} if $ \leftmul{x} $ and $ \rightmul{x} $ are invertible maps such that $ \leftmul{x}^{-1} = \leftmul{y} $ and $ \rightmul{x}^{-1} = \rightmul{y} $. Further, we say that $ x $ is \defemph*{invertible} if there exists $ z \in \ring $ such that $ z $ is an inverse of $ x $. 
\end{definition}

\begin{note}\label{ring:inv-note}
	There exists no standard notion of invertibility in nonassociative rings. For example, elements $ a,b \in \ring $ are called \defemph*{inverses of each other} in \cite[\onpage{38}]{Schafer} if $ ab = 1_\ring = ba $, which is a weaker notion of invertibility than ours. However, we will see in \cref{altring:invert-char} that the two notions are equivalent for alternative rings.
	
	Our notion of invertibility is precisely the one which relates to Weyl-invertible elements in root gradings (in the sense of \thmitemcref{weyl:weyl-def}{weyl:weyl-def:invertible}). This fact was first observed by Faulkner in \cite[Theorem~13.8]{Faulkner-NonAssocProj}, and we will prove it in \thmitemcref{ADE:stsign:weyl-char}{ADE:stsign:weyl-char:bij}.
\end{note}

\begin{remark}
	If $ a $ is invertible with inverse $ b $, then it follows from an equation of the form $ ax=y $ that $ by = b(ax) = \leftmul{b} \leftmul{a}(x) = x $. This is not necessarily the case for the weaker notion of invertibility in \cref{ring:inv-note}.
\end{remark}

\begin{remark}
	The product of two invertible elements is not necessarily invertible. This is, of course, a major drawback of our notion of invertibility.
\end{remark}

\begin{definition}[Division ring]\label{ring:division-ring-def}
	A \defemph*{division ring}\index{ring!division} is a ring in which every non-zero element is invertible.
\end{definition}

\begin{remark}\label{ring:inv-def-rem}
	Let $ x,y \in \ring $. Then $ y $ is an inverse of $ x $ if and only if all the maps $ \leftmul{x} \circ \leftmul{y} $, $ \leftmul{y} \circ \leftmul{x} $, $ \rightmul{x} \circ \rightmul{y} $ and $ \rightmul{y} \circ \rightmul{x} $ are trivial, which is to say that
	\begin{equation*}
		x(yz) = y(xz) = z = (zx)y = (zy)x \quad \text{for all } z \in \ring.
	\end{equation*}
	Thus clearly, if $ y $ is an inverse of $ x $, then $ x $ is an inverse of $ y $.
\end{remark}

\begin{lemma}\label{ring:inv-char}
	For all $ x,y \in \ring $, the following assertions are equivalent:
	\begin{stenumerate}
		\item $ x $ is an inverse of $ y $.
		
		\item \label{ring:inv-char:assoc}$ xy = yx = 1_\ring $ and $ \assoc{x}{y}{z} = \assoc{y}{x}{z} = \assoc{z}{x}{y} = \assoc{z}{y}{x} = 0_\ring $ for all $ z \in \ring $.
	\end{stenumerate}
\end{lemma}
\begin{proof}
	At first, assume that $ x $ is an inverse of $ y $. Putting $ z \defl 1_\ring $ in \cref{ring:inv-def-rem} yields that $ xy = 1_\ring = yx $. Further, it follows from \cref{ring:inv-def-rem} that $ x(yz) = z = 1_\ring z = (xy)z $ for all $ z \in \ring $, so $ \assoc{x}{y}{z} = 0_\ring $. The remaining statements in~\itemref{ring:inv-char:assoc} can be proven similarly. In essentially the same way, it follows from~\itemref{ring:inv-char:assoc} that $ x $ is an inverse of $ y $. 
\end{proof}

\begin{lemma}\label{ring:inverse-eq}
	Let $ x \in \ring $ be invertible. Then the following statements are equivalent for all $ y \in \ring $:
	\begin{stenumerate}
		\item $ y $ is an inverse of $ x $.
		
		\item $ xy=1_\ring $.
		
		\item $ yx = 1_\ring $.
	\end{stenumerate}
	Further, there exists a unique $ y \in \ring $ satisfying these properties.
\end{lemma}
\begin{proof}
	If $ y $ is an inverse of $ x $, we have already shown in \cref{ring:inv-char} that $ xy = 1_\ring = yx $. Conversely, assume that $ y \in \ring $ satisfies $ xy = 1_\ring $. Since $ x $ is invertible, there exists an element $ x^{-1} $ which is an inverse of $ x $. Multiplying the equation $ xy = 1_\ring = xx^{-1} $ by $ x^{-1} $ from the left, we obtain $ y = x^{-1} $. Thus any $ y \in \ring $ satisfying $ xy=1_\ring $ equals $ x^{-1} $. Similarly, we can prove that any $ y \in \ring $ satisfying $ yx = 1_\ring $ equals $ x^{-1} $. In particular, every inverse of $ x $ equals $ x^{-1} $, so all inverses of $ x $ are identical. This finishes the proof. 
\end{proof}

\begin{notation}
	If $ x $ is an invertible element of $ \ring $, we denote its unique inverse by $ x^{-1} $.
\end{notation}

\begin{note}\label{ring:inv-warning}
	The assumption that $ x $ is invertible in \cref{ring:inverse-eq} is crucial. In general, if $ x $ and $ y $ are arbitrary elements of $ \ring $ such that $ xy=1_\ring = yx $, then it is not clear that $ x $ and $ y $ are invertible (see \cref{ring:inv-char}). Further, the product of two invertible ring elements is not necessarily invertible.
\end{note}

For nuclear elements, the warning in \cref{ring:inv-warning} does not apply.

\begin{lemma}
	Let $ x \in \nucleus{\ring} $ and $ y \in \ring $ such that $ xy = 1_\ring = yx $. Then $ x $ is invertible and $ x^{-1} = y $.
\end{lemma}
\begin{proof}
	This follows from \cref{ring:inv-char}.
\end{proof}

We end this section with a brief introduction of modules over nonassociative rings which will be needed in \cref{chap:BC}. The definition is word by word the same as for associative rings. However, for reasons that we will see in \cref{ring:module-obs}, it is not standard in nonassociative algebra. In fact, we are not aware of a single reference which uses this definition.

\begin{definition}[Module]\label{ring:module}
	A \defemph*{right $ \ring $-module}\index{module} is an abelian group $ (\module, +) $ together with a scalar multiplication $ \map{\scmult}{\module \times \ring}{\module}{}{} $ satisfying the following conditions for all $ r,s \in \ring $ and all $ v,w \in \module $:
	\begin{stenumerate}
		\item $ (v+w) \scmult r = (v \scmult r) + (w \scmult r) $.
		
		\item $ v \scmult (r+s) = (v \scmult r) + (v \scmult s) $.
		
		\item \label{ring:module:assoc}$ v \scmult (r \rmult s) = (v \scmult r) \scmult s $.
		
		\item $ v \scmult 1_\ring = v $.
	\end{stenumerate}
	It is called \defemph*{faithful}\index{module!faithful} if for any $ r \in \ring $ with $ v \scmult r = 0 $ for all $ v \in \module $, we have $ r=0_\ring $. Further, \defemph*{left $ \ring $-modules} are defined similarly.
\end{definition}

\begin{observation}\label{ring:module-obs}
	Let $ \module $ be a module over $ \ring $. Let $ a,b,c \in \ring $ and let $ v \in \module $. Using only Axiom~\thmitemref{ring:module}{ring:module:assoc}, we see that
	\begin{align*}
		v \scmult \brackets[\big]{(ab)c} &= \brackets[\big]{v \scmult (ab)} \scmult c = \brackets[\big]{(v \scmult a) \scmult b} \scmult c = (v \scmult a) \scmult (bc) = v \scmult \brackets[\big]{a(bc)}.
	\end{align*}
	Thus $ \module $ is actually a module over the associative ring $ \ring / \assoc{\ring}{\ring}{\ring} $. This explains why \cref{ring:module} is not standard in the theory of nonassociative rings. However, it is of use in the theory of root graded groups: We will see root graded groups of type $ (B)C $ in which some root groups are coordinatised by a ring $ \ring $ while other root groups are coordinatised by an algebraic structure which involves a module over $ \ring $. In this setting, we do not have the \enquote{freedom} to replace $ \ring $ by $ \ring/ \assoc{\ring}{\ring}{\ring} $ because the full ring $ \ring $ is needed as a coordinatising structure for some root groups.
	
	On a final note, recall that a module over an associative ring $ \ring $ can also be seen as an abelian group $ (\module, +) $ together with a homomorphism $ \map{}{\ring}{\End(\module, +)}{}{} $. Since the ring $ \End(\module, +) $ is associative, it is natural from this viewpoint that any module structure over $ \ring $ should factor through $ \ring / \assoc{\ring}{\ring}{\ring} $.
\end{observation}


\section{The Simply-laced Root Systems}

\label{sec:simply-laced}

\begin{secnotation}
	Unless otherwise specified, we denote by $ \roots $ a simply-laced root system in the sense of \cref{rootsys:simply-laced-def}.
\end{secnotation}

In this section, we collect some properties of the simply-laced root systems. Most importantly, we compute their Cartan integers. All results in this section are straightforward to verify, usually by an inspection of the root system in its standard representation, so we usually leave out proofs.

We begin with the standard representation of $ A_n $. Similar representations of the root systems of types $ D $ and $ E $ can be found in \cite[Section~12.1]{HumphreysLieAlg}, but we will not need them.

\begin{remark}[Standard representation of $ A_n $]\label{ADE:simply-laced:An-standard-rep}
	Let $ n \in \Npos $ and let $ V $ be a Euclidean space of dimension $ n+1 $ with orthonormal basis $ \tup{\basvec}{n+1} $. The \defemph*{standard representation of $ A_n $}\index{standard representation!of An@of $ A_n $} is
	\[ A_n \defl \Set{\basvec_i - \basvec_j \given i \ne j \in \numint{1}{n+1}}. \]
	Note that a root $ \basvec_i - \basvec_j $ is uniquely determined by the tuple $ (i,j) $ of indices. For this reason, we will often write $ ij $ in place of $ \basvec_i - \basvec_j $. For example, we could write $ \rootgr{ij} $ for the root group $ \rootgr{\basvec_i - \basvec_j} $. The \defemph*{standard (ordered) root base} is
	\[ \rootbase \defl \Set{\basvec_i - \basvec_{i+1} \given i \in \numint{1}{n}} \qquad \brackets[\big]{\text{or } \rootbase_{\text{ord}} \defl (\basvec_1 - \basvec_2, \ldots, \basvec_n - \basvec_{n+1})} \]
	and the corresponding positive system is
	\[ \possys \defl \Set{\basvec_i - \basvec_{j} \given i<j \in \numint{1}{n+1}}. \]
	Observe that $ V $ is not spanned by $ A_n $.
\end{remark}

\begin{definition}[$ A_2 $-pairs and $ A_2 $-triples]
	Let $ \roots $ be any root system. An \defemph*{$ A_2 $-pair (in $ \roots $)}\index{A2-pair@$ A_2 $-pair} is a tuple $ (\alpha, \gamma) $ of roots such that $ \Set{\alpha, \gamma} $ is a root base of the parabolic root subsystem that it spans. An \defemph*{$ A_2 $-triple (in $ \roots $)}\index{A2-triple@$ A_2 $-triple} is a triple $ (\alpha, \beta, \gamma) $ of roots such that $ (\alpha, \gamma) $ is an $ A_2 $-pair and $ \beta = \alpha + \gamma $.
\end{definition}
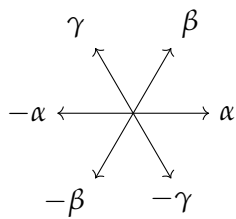
\begin{figure}[htb]
	\centering\begin{tikzpicture}
		\draw[->] (0,0) -- (0:1);
		\draw[->] (0,0) -- (60:1);
		\draw[->] (0,0) -- (120:1);
		\draw[->] (0,0) -- (180:1);
		\draw[->] (0,0) -- (240:1);
		\draw[->] (0,0) -- (300:1);
		
		\node[right] at (0:1){$ \alpha $};
		\node[above right] at (60:1){$ \beta $};
		\node[above left] at (120:1){$ \gamma $};
		\node[left] at (180:1){$ -\alpha $};
		\node[below left] at (240:1){$ -\beta $};
		\node[below] at (300:1){$ -\gamma $};
	\end{tikzpicture}
	\caption{An $ A_2 $-triple $ (\alpha, \beta, \gamma) $.}
	\label{fig:A2-trip}
\end{figure}

\begin{remark}
	In an $ A_2 $-triple $ (\alpha, \beta, \gamma) $, the roots $ \alpha, \beta, \gamma $ are exactly the positive roots of the corresponding $ A_2 $-subsystem (with respect to the root base $ (\alpha, \gamma) $). Further, the reflection $ \refl{\alpha} $ interchanges $ \beta $ with $ \gamma $, $ -\beta $ with $ -\gamma $ and $ \alpha $ with $ -\alpha $. See \cref{fig:A2-trip} for an illustration.
\end{remark}

\begin{note}
	The terminology of \enquote{$ A_2 $-pairs} is borrowed from \cite[(1.4)]{Shi1993}. We will introduce obvious variants of this notion for the root systems $ B_2 $ and $ BC_2 $ as well (\cref{B:B2-pair-def,BC:BC2-pair-def}).
\end{note}

We now collect some basic results.

\begin{remark}\label{rootsys:simply-laced-int}
	For all non-proportional roots $ \alpha, \beta $, we have $ \rootintcox{\alpha}{\beta} \subs \compactSet{\alpha+\beta} $. In fact, this property characterises simply-laced root systems. In combination with \cref{basic:comm-add}, this is the main reason why root graded groups of simply-laced type are the easiest class to study. Note that, in particular, we always have $ \rootint{\alpha}{\beta} = \rootintcox{\alpha}{\beta} $, so that any $ \roots $-graded group is a crystallographic $ \roots $-graded group.
\end{remark}

\begin{remark}[The Weyl group of $ A_n $]\label{ADE:A-weylgrp}
	Let $ n \in \Npos $ and let $ \rootbase $ be the standard root base of $ A_n $. An easy computation shows that
	\[ \basvec_i^{\reflbr{\basvec_j - \basvec_k}} = \basvec_{i^{\transp{j}{k}}} \]
	for all $ i,j,k \in \numint{1}{n+1} $ with $ j \ne k $ where $ \transp{j}{k} $ denotes the transposition which interchanges $ j $ and $ k $. It follows that there exists an isomorphism $ \phi $ from the Weyl group $ W $ of $ A_n $ to the symmetric group on $ \numint{1}{n+1} $ which maps $ \reflbr{\basvec_j - \basvec_k} $ to $ \transp{j}{k} $ for all distinct $ j,k \in \numint{1}{n+1} $. This isomorphism satisfies $ \basvec_i^w = \basvec_{i^{\phi(w)}} $ for all $ w \in W $ and $ i \in \numint{1}{n+1} $. Observe that $ \phi $ depends on the choice of the ordered orthonormal base $ (\basvec_1, \ldots, \basvec_{n+1}) $.
\end{remark}

\begin{remark}\label{ADE:simply-laced:all-in-A2}
	If $ \roots $ is irreducible and of rank at least~2 (or more generally, if $ \roots $ does not contain an irreducible component of type $ A_1 $), then every root lies in an $ A_2 $-triple.
\end{remark}

\begin{remark}\label{ADE:rank2-is-parabolic}
	Every subsystem of $ \roots $ of rank at most 2 is parabolic. In particular, we do not have to distinguish between $ A_2 $-subsystems, closed $ A_2 $-subsystems and parabolic $ A_2 $-subsystems of $ \roots $.
\end{remark}

\begin{lemma}\label{ADE:ortho-adjacent}
	If $ \alpha, \beta \in \roots $ are orthogonal, then $ \alpha $ is adjacent to $ \beta $ and $ -\beta $.
\end{lemma}

An inspection of the root system $ A_2 $ yields the following characterisation of Cartan integers in $ \roots $. It will play an important role in the proof of the square formula for Weyl elements (\cref{A2Weyl:cartan-comp}).

\begin{proposition}[Cartan integers]\label{ADE:cartan-char}
	Let $ \alpha, \beta $ be roots. The Cartan integer $ \cartanint{\beta}{\alpha} = 2 \frac{\beta \cdot \alpha}{\alpha \cdot \alpha} $ is determined as follows:
	\begin{lemenumerate}
		\item $ \cartanint{\beta}{\alpha} = 0 $ if and only if $ \alpha $ and $ \beta $ are orthogonal (or equivalently, if and only if they lie in no common $ A_2 $-subsystem).
		
		\item $ \cartanint{\beta}{\alpha} = -1 $ if and only if $ (\alpha, \beta) $ is an $ A_2 $-pair.
		
		\item $ \cartanint{\beta}{\alpha} = 1 $ if and only if $ (\alpha, -\beta) $ is an $ A_2 $-pair.
		
		\item $ \cartanint{\beta}{\alpha} = 2 $ if and only if $ \alpha =\beta $.
		
		\item $ \cartanint{\beta}{\alpha} = -2 $ if and only if $ \alpha = -\beta $.
	\end{lemenumerate}
	Further, $ \cartanint{\alpha}{\beta} = \cartanint{\beta}{\alpha} $.
\end{proposition}

We end this section with the notion of positive $ A_2 $-pairs. They will be needed to properly define coordinatisations of $ A_n $-graded groups (\cref{ADE:param-def}).

\begin{definition}[Positive $ A_2 $-pair]\label{ADE:pos-pair}
	Let $ (\alpha, \beta) $, $ (\alpha', \beta') $ be two $ A_2 $-pairs in $ \roots $. We say that $ (\alpha', \beta') $ is \defemph*{$ (\alpha, \beta) $-positive}\index{A2-pair!positive@positive} if they lie in the same orbit under the Weyl group $ \Weyl(\roots) $, that is, if there exists $ u \in \Weyl(\roots) $ such that $ \alpha^u = \alpha' $ and $ \beta^u = \beta' $.
\end{definition}

The following result is an easy computation. It is the underlying reason why any ring which coordinatises a root graded group of type $ D $ or $ E $ is commutative while the same is not true for root graded groups of type $ A $.

\begin{lemma}[{\cite[1.4]{RGLie-BermanMoody}}]\label{ADE:A2-pair-orbits}
	Let $ (\alpha_0, \beta_0) $ be an $ A_2 $-pair in $ \roots $, and assume that $ \roots $ is irreducible.
	\begin{lemenumerate}
		\item If $ \roots $ is of type $ A $, then the set of $ A_2 $-pairs has exactly two orbits under the Weyl group of $ \roots $: The set of $ (\alpha_0, \beta_0) $-positive $ A_2 $-pairs and the set of $ (\beta_0, \alpha_0) $-positive $ A_2 $-pairs.
		
		\item \label{ADE:A2-pair-orbits:DE}If $ \roots $ is of type $ D $ or $ E $, then every $ A_2 $-pair is $ (\alpha_0, \beta_0) $-positive.
	\end{lemenumerate}
\end{lemma}

The following result on positive $ A_2 $-pairs will be needed to prove that the coordinatising ring of an $ A_3 $-graded group is associative (\cref{ADE:assoc}).

\begin{lemma}\label{ADE:rank-3-has-A3}
	Assume that $ \roots $ is of rank at least~3 and choose some $ A_2 $-pair $ (\alpha_0, \gamma_0) $. Then there exist roots $ \alpha, \beta, \gamma $ with the following properties:
	\begin{stenumerate}
		\item $ \gen{\alpha, \beta, \gamma}_\IR \intersect \roots $ is a root subsystem of type $ A_3 $ with root base $ \Set{\alpha, \beta, \gamma} $.
		
		\item $ (\alpha, \beta) $ and $ (\beta, \gamma) $ are $ (\alpha_0, \gamma_0) $-positive $ A_2 $-pairs.
	\end{stenumerate}
	Further, any such choice of roots has the property that $ (\alpha+\beta, \gamma) $ and $ (\alpha, \beta+\gamma) $ are $ (\alpha_0, \gamma_0) $-positive $ A_2 $-pairs.
\end{lemma}
\begin{proof}
	Since $ \roots $ is of type $ A $, $ D $ or $ E $, it is easy to see that $ \roots $ has a parabolic subsystem of type $ A_3 $. Thus we can choose roots $ \alpha, \beta, \gamma $ such that $ \Set{\alpha, \beta, \gamma} $ is a root base of a parabolic $ A_3 $-subsystem of $ \roots $ and such that $ (\alpha, \beta) $, $ (\beta, \gamma) $ are $ A_2 $-pairs. We will show that either both $ (\alpha, \beta) $ and $ (\beta, \gamma) $ or both $ (\gamma, \beta) $ and $ (\beta, \alpha) $ are $ (\alpha_0, \gamma_0) $-positive. Note that $ (\alpha, \beta) $ and $ (\beta, \gamma) $ lie in the same orbit under the Weyl group because $ \refl{\gamma} \refl{\beta} \refl{\alpha} $ maps $ (\alpha, \beta) $ to $ (\beta, \gamma) $:
	\[ \begin{tikzcd}[column sep=normal, row sep=tiny]
		\alpha \arrow[r, "\reflbr{\gamma}"] & \alpha \arrow[r, "\reflbr{\beta}"] & \alpha + \beta \arrow[r, "\reflbr{\alpha}"] & -\alpha + (\alpha+\beta) = \beta, \\
		\beta \arrow[r, "\reflbr{\gamma}"] & \beta+\gamma \arrow[r, "\reflbr{\beta}"] & -\beta + (\beta+\gamma) = \gamma \arrow[r, "\reflbr{\alpha}"] & \gamma.
	\end{tikzcd} \]
	By the same argument, $ (\beta, \alpha) $ and $ (\gamma, \beta) $ lie in the same orbit. Thus for the existence assertion, it only remains to show that $ (\alpha, \beta) $ or $ (\beta, \alpha) $ is $ (\alpha_0, \gamma_0) $-positive. This holds by \cref{ADE:A2-pair-orbits}.
	
	Now let $ \alpha, \beta, \gamma $ be any roots with the desired properties. Then $ \refl{\alpha} $ maps $ (\beta, \gamma) $ to $ (\alpha+\beta, \gamma) $, so $ (\alpha+\beta, \gamma) $ is an $ (\alpha_0, \gamma_0) $-positive $ A_2 $-pair (because $ (\beta, \gamma) $ has the same properties). Similarly, $ \refl{\gamma} $ maps $ (\alpha, \beta) $ to $ (\alpha, \beta+\gamma) $, so $ (\alpha, \beta+\gamma) $ is also an $ (\alpha_0, \gamma_0) $-positive $ A_2 $-pair. This finishes the proof.
\end{proof}


\section{Construction of Simply-laced Root Graded Groups}

\label{sec:ADE:ex}

For simply-laced root graded groups, the existence problem is (nearly) completely solved by the theory of Chevalley groups. In this sections, we collect some remarks on this topic. Recall that the existence problem was defined in \cref{rgg:existence-problem}.

\begin{miscthm}[The existence problem for types $ D $ and $ E $]
	For any root system $ \roots $, the Chevalley groups of type $ \roots $ constitute examples of $ \roots $-graded groups for any commutative associative ring. Since we will show that every root graded group of type $ D $ or $ E $ is coordinatised by such a ring, the Chevalley groups thus provide a complete solution of the existence problem for these types. We emphasise that this does not mean that every root graded group of type $ D $ or $ E $ is isomorphic to a Chevalley group. Instead, it means that every such group $ G $ is coordinatised by a commutative associative ring $ \ring $ which also coordinatises some Chevalley group $ H $ such that $ G $ and $ H $ satisfy the same commutator relations.
\end{miscthm}

\begin{miscthm}[The existence problem for type $ A $ in rank at least~3]
	We will show that $ A_n $-graded groups for $ n \ge 3 $ are coordinatised by associative rings which need not be commutative. Thus the Chevalley groups only provide a partial solution of the existence problem. However, a complete solution of this problem is not far away: Recall from \cref{chev:ex-SL} that for any commutative ring $ \ring $, the group $ \E_{n+1}(\ring) $ is a Chevalley group of type $ A_n $ and that the definition of this group makes sense even when $ \ring $ is not assumed to be commutative. It is not difficult to see that $ \E_{n+1}(\ring) $ is $ A_n $-graded for any associative ring $ \ring $. Thus this group solves the existence problem for type $ A_n $.
	
	It is noteworthy that, even though Chevalley groups do not solve the existence problem for type $ A_n $, they are still sufficiently generic in the sense of \cref{param:motiv:parmap-choice} to read off a parity map for the parametrisation of arbitrary $ A_n $-graded groups.
\end{miscthm}

\begin{note}[Existence problem for $ A_2 $-graded groups]\label{ADE:existence-A2}
	The rank-2 case (which, we recall, is only of secondary concern in this book) is more difficult. For any alternative ring $ \ring $, Faulkner constructs in \cite[Section~3]{Faulkner-StableRange} a Jordan pair $ V = V(\ring) $, a Tits-Kantor-Koecher algebra $ T = \operatorname{TKK}(V) $ and a group $ G(\ring) = G(T) $ such that $ G(\ring) $ is $ A_2 $-graded and coordinatised by $ \ring $. Essentially the same construction, but phrased without the language of Jordan pairs, can also be found in the appendix of \cite{Faulkner-Barb}. This is the most general known construction of $ A_2 $-graded groups. There is no known example of an $ A_2 $-graded group which is coordinatised by a ring that is not alternative. However, the problem whether every such ring is alternative remains open.
\end{note}


\section{Rank-2 Computations}

\label{sec:ADE:weyl}

\begin{secnotation}
	We denote by $ \roots $ any irreducible simply-laced root system of rank at least~$ 2 $ and by $ G $ a group which has $ \roots $-commutator relations with root groups $ (\rootgr{\alpha})_{\alpha \in \roots} $ (in the sense of \cref{rgg:group-commrel-def}). Further, we assume that all for all non-proportional roots $ \alpha $, $ \beta $, we have $ \rootgr{\alpha} \intersect \rootgr{\beta} = \compactSet{1_G} $.
\end{secnotation}

In this section, we study the action of Weyl elements in simply-laced root graded groups. Our goal is to show that $ G $ satisfies the square formula for Weyl elements (\cref{A2Weyl:cartan-comp}). This is also essentially the content of \cite[(2.4)(iii)]{Shi1993} and \cite[(19.3)]{MoufangPolygons}, but we prove it in a more general setting (see \cref{sec:literature}). The core arguments, however, are mostly the same.

Most of the time, we will fix some $ A_2 $-triple and perform all computations in the corresponding $ A_2 $-graded subgroup (see \cref{rgg:subgroup}). Since every root lies in a parabolic subsystem of type $ A_2 $ by \cref{ADE:simply-laced:all-in-A2}, this is no restriction of generality.

Observe that we do not assume that $ \invset{\alpha} $ is non-empty for all roots $ \alpha $. However, we will most of the time make statements about arbitrary Weyl elements, and these statements are of course empty if no such elements exist. All statements which do not involve Weyl elements are independent of their existence, however. Further, note that we do not require Axiom~\thmitemref{rgg-def}{rgg-def:nondeg}, simply because all proofs work without difficulty under the weaker assumption that $ \rootgr{\alpha} \intersect \rootgr{\beta} = \compactSet{1_G} $ for all non-proportional roots $ \alpha $, $ \beta $. See also \cref{rgg:nondeg-comparison}.

We begin with a computational result which will be our main tool in this section. It is essentially the special case of $ n=3 $ in \cref{tw:6.4} (which is the equivalent of \cite[(6.4)]{MoufangPolygons}), but we prove it again in a slightly different way. The same arguments can be found in the proof of \cite[(2.4)]{Shi1993}. The same result for RGD-systems is proven in \cite[(19.1)]{MoufangPolygons}.

\begin{proposition}\label{A2Weyl:basecomp-cor}
	Let $ (\alpha, \beta, \gamma) $ be an $ A_2 $-triple. Assume that $ (a_{-\alpha}, b_\alpha, c_{-\alpha}) $ is an $ \alpha $-Weyl triple and denote its corresponding Weyl element by $ w_\alpha $. Then the following statements hold:
	\begin{proenumerate}
		\item \label{A2Weyl:basecomp-cor:1}For any $ x_\gamma \in \rootgr{\gamma} $, we have $ x_\gamma = \commutator{\commutator{x_\gamma}{b_\alpha}}{c_{-\alpha}}^{-1} $ and $ x_\gamma^{w_\alpha} = \commutator{x_\gamma}{b_\alpha} $.
		
		\item \label{A2Weyl:basecomp-cor:2}For any $ x_\beta \in \rootgr{\beta} $, we have $ x_\beta = \commutator{\commutator{x_\beta}{a_{-\alpha}}}{b_\alpha}^{-1} $ and $ x_\beta^{w_\alpha} = \commutator{x_\beta}{a_{-\alpha}} $.
	\end{proenumerate}
\end{proposition}
\begin{proof}
	Let $ x_\gamma \in \rootgr{\gamma} $. Using the commutator axiom of root graded groups and \thmitemcref{group-rel}{group-rel:conj1}, we perform the following computation:
	\begin{align*}
		x_\gamma^{w_\alpha} &= x_\gamma^{b_\alpha c_{-\alpha}} = \brackets[\big]{x_\gamma \commutator{x_\gamma}{b_\alpha}}^{c_{-\alpha}} = x_\gamma^{c_{-\alpha}} \commutator{x_\gamma}{b_\alpha}^{c_{-\alpha}} = x_\gamma \commutator{x_\gamma}{b_\alpha} \commutator[\big]{\commutator{x_\gamma}{b_\alpha}}{c_{-\alpha}}.
	\end{align*}
	Note that $ \commutator{\commutator{x_\gamma}{b_\alpha}}{c_{-\alpha}} $ lies in $ \rootgr{\gamma} $ while $ \commutator{x_\gamma}{b_\alpha} $ lies in $ \rootgr{\beta} $, so these two terms commute. Thus we can write
	\[ x_\gamma^{w_\alpha} = x_\gamma \commutator[\big]{\commutator{x_\gamma}{b_\alpha}}{c_{-\alpha}} \commutator{x_\gamma}{b_\alpha} \]
	with $ x_\gamma \commutator{\commutator{x_\gamma}{b_\alpha}}{c_{-\alpha}} \in \rootgr{\gamma} $ and $ \commutator{x_\gamma}{b_\alpha} \in \rootgr{\beta} $. Since $ w_\alpha $ is an $ \alpha $-Weyl element, $ x_\gamma^{w_\alpha} $ lies in $ \rootgr{\beta} $. Hence
	\[ x_\gamma^{w_\alpha} \commutator{x_\gamma}{b_\alpha}^{-1} = x_\gamma \commutator[\big]{\commutator{x_\gamma}{b_\alpha}}{c_{-\alpha}} \]
	where the left-hand side lies in $ \rootgr{\beta} $ and the right-hand side lies in $ \rootgr{\gamma} $. Since $ \rootgr{\beta} \intersect \rootgr{\gamma} = \compactSet{1_G} $ by assumption, we infer that $ x_\gamma \commutator{\commutator{x_\gamma}{b_\alpha}}{c_{-\alpha}} = 1_G $. In other words, $ x_\gamma = \commutator{\commutator{x_\gamma}{b_\alpha}}{c_{-\alpha}}^{-1} $. This finishes the proof of~\itemref{A2Weyl:basecomp-cor:1}.
	
	Now let $ x_\beta \in \rootgr{\beta} $. We perform the same kind of computation as in~\itemref{A2Weyl:basecomp-cor:1}, but the result is slightly more complicated because $ x_\beta $ and $ a_{-\alpha} $ do not commute:
	\begin{align*}
		x_\beta^{w_\alpha} &= \brackets[\big]{x_\beta \commutator{x_\beta}{a_{-\alpha}}}^{b_\alpha c_{-\alpha}} = \brackets[\big]{x_\beta \commutator{x_\beta}{a_{-\alpha}} \commutator[\big]{\commutator{x_\beta}{a_{-\alpha}}}{b_\alpha}}^{c_{-\alpha}} \\
		&= x_\beta \commutator{x_\beta}{c_{-\alpha}} \commutator{x_\beta}{a_{-\alpha}} \commutator[\big]{\commutator{x_\beta}{a_{-\alpha}}}{b_\alpha} \commutator[\bigg]{\commutator[\big]{\commutator{x_\beta}{a_{-\alpha}}}{b_\alpha}}{c_{-\alpha}}.
	\end{align*}
	Note that $ x_\beta $ and $ \commutator{\commutator{x_\beta}{a_{-\alpha}}}{b_\alpha} $ lie in $ \rootgr{\beta} $ while all the other factors on the right-hand side lie in $ \rootgr{\gamma} $. Since $ x_\beta^{w_\alpha} $ lies in $ \rootgr{\gamma} $, we infer that
	\[ x_\beta \commutator[\big]{\commutator{x_\beta}{a_{-\alpha}}}{b_\alpha} = x_\beta^{w_\alpha} \brackets[\bigg]{\commutator{x_\beta}{c_{-\alpha}} \commutator{x_\beta}{a_{-\alpha}} \commutator[\bigg]{\commutator[\big]{\commutator{x_\beta}{a_{-\alpha}}}{b_\alpha}}{c_{-\alpha}}}^{-1} \]
	where the left-hand side lies in $ \rootgr{\beta} $ and the right-hand side lies in $ \rootgr{\gamma} $. Since $ \rootgr{\beta} \intersect \rootgr{\gamma} = \compactSet{1_G} $, it follows that $ x_\beta = \commutator{\commutator{x_\beta}{a_{-\alpha}}}{b_\alpha}^{-1} $ and
	\[ x_\beta^{w_\alpha} = \commutator{x_\beta}{c_{-\alpha}} \commutator{x_\beta}{a_{-\alpha}} \commutator[\bigg]{\commutator[\big]{\commutator{x_\beta}{a_{-\alpha}}}{b_\alpha}}{c_{-\alpha}} = \commutator{x_\beta}{c_{-\alpha}} \commutator{x_\beta}{a_{-\alpha}} \commutator{x_\beta^{-1}}{c_{-\alpha}}. \]
	Note that $ \commutator{x_\beta}{c_{-\alpha}} $ lies in $ \rootgr{\gamma} $, so it commutes with $ x_\beta $ and $ a_{-\alpha} $ and thus also with $ \commutator{x_\beta}{a_{-\alpha}} $. Hence we have $ x_\beta^{w_\alpha} = \commutator{x_\beta}{a_{-\alpha}} \commutator{x_\beta}{c_{-\alpha}} \commutator{x_\beta^{-1}}{c_{-\alpha}} $ where $ \commutator{x_\beta^{-1}}{c_{-\alpha}} = \commutator{x_\beta}{c_{-\alpha}}^{-1} $ by \cref{basic:comm-add}, so $ x_\beta^{w_\alpha} = \commutator{x_\beta}{a_{-\alpha}} $. This finishes the proof of~\itemref{A2Weyl:basecomp-cor:2}.
\end{proof}

\begin{note}\label{A2Weyl:spirit}
	In the following, we will often use that for any $ A_2 $-triple $ (\alpha, \beta, \gamma) $, $ (-\alpha, \gamma, \beta) $ is an $ A_2 $-triple as well. This is easily seen from \cref{fig:A2-trip}. Thus \thmitemcref{A2Weyl:basecomp-cor}{A2Weyl:basecomp-cor:1} not only tells us how $ \alpha $-Weyl elements act on $ \rootgr{\gamma} $ but also  how $ (-\alpha) $-Weyl elements act on $ \rootgr{\beta} $. Since any $ \alpha $-Weyl element is a $ (-\alpha) $-Weyl element as well (by \thmitemcref{basic:weyl-general}{basic:weyl-general:minus}), this is a powerful tool. In fact, this provides us with an alternative way to prove \thmitemcref{A2Weyl:basecomp-cor}{A2Weyl:basecomp-cor:2}. See \cref{A2Weyl:left-right-comm-altproof} for more details.
\end{note}

We will show in \thmitemcref{A2Weyl:weyl}{A2Weyl:weyl:leftright} that every Weyl triple is weakly balanced (and thus balanced by \cref{basic:all-balanced} or by \thmitemcref{basic:weylmap-form2}{basic:weylmap-form2:balanced}). The following result is a first step into this direction.

\begin{lemma}\label{A2Weyl:left-right-comm}
	Let $ (\alpha, \beta, \gamma) $ be an $ A_2 $-triple and assume that $ (a_{-\alpha}, b_\alpha, c_{-\alpha}) $ is an $ \alpha $-Weyl triple. Then $ \commutator{x_\beta}{a_{-\alpha}} = \commutator{x_\beta}{c_{-\alpha}} $ for all $ x_\beta \in \rootgr{\beta} $.
\end{lemma}
\begin{proof}
	Let $ x_\beta \in \rootgr{\beta} $ and denote by $ w_\alpha \defl a_{-\alpha} b_\alpha c_{-\alpha} $ the Weyl element corresponding to $ (a_{-\alpha}, b_\alpha, c_{-\alpha}) $. We already know from \thmitemcref{A2Weyl:basecomp-cor}{A2Weyl:basecomp-cor:2} that $ x_\beta^{w_\alpha} = \commutator{x_\beta}{a_{-\alpha}} $. Further, we know from \thmitemcref{basic:weyl-general}{basic:weyl-general:minus} that $ (b_\alpha, c_{-\alpha}, a_{-\alpha}^{w_\alpha}) $ is a $ (-\alpha) $-Weyl triple with corresponding Weyl element $ w_\alpha $. Applying \thmitemcref{A2Weyl:basecomp-cor}{A2Weyl:basecomp-cor:1} to the $ A_2 $-triple $ (-\alpha, \gamma, \beta) $ and the $ (-\alpha) $-Weyl triple $ (b_\alpha, c_{-\alpha}, a_{-\alpha}^{w_\alpha}) $, we obtain that $ x_\beta^{w_\alpha} = \commutator{x_\beta}{c_{-\alpha}} $. This finishes the proof.
\end{proof}

\begin{note}\label{A2Weyl:left-right-comm-altproof}
	Let everything be as in the proof of \cref{A2Weyl:left-right-comm}. By choosing the Weyl triple $ (c_{-\alpha}^{w_\alpha^{-1}}, a_{-\alpha}, b_\alpha) $ instead of $ (b_\alpha, c_{-\alpha}, a_{-\alpha}^{w_\alpha}) $, we obtain that $ x_\beta^{w_\alpha} = \commutator{x_\beta}{a_{-\alpha}} $. This provides a different proof of (the main result of) \thmitemcref{A2Weyl:basecomp-cor}{A2Weyl:basecomp-cor:2}.
\end{note}

Recall that the axioms of a group with $ \roots $-commutator relations merely require that $ \commutator{\rootgr{\alpha}}{\rootgr{\gamma}} $ is a subset of $ \rootgr{\alpha+\gamma} $ if $ (\alpha, \gamma) $ is an $ A_2 $-triple. By the following result, we actually have equality, provided that there exist enough Weyl elements.

\begin{proposition}[{\cite[(2.4)(i)]{Shi1993}}]\label{A2Weyl:rootisom}
	Let $ (\alpha, \beta, \gamma) $ be an $ A_2 $-triple. Then for all $ b_\alpha \in \invset{\alpha} $, the map
	\[ \map{}{\rootgr{\gamma}}{\rootgr{\beta}}{x_\gamma}{\commutator{x_\gamma}{b_\alpha}} \]
	is an isomorphism of groups. In particular, $ \commutator{\rootgr{\alpha}}{\rootgr{\gamma}} = \rootgr{\beta} $ if $ \invset{\alpha} $ is non-empty.
\end{proposition}
\begin{proof}
	Choose $ a_{-\alpha}, c_{-\alpha} \in \rootgr{-\alpha} $ such that $ w_\alpha \defl a_{-\alpha} b_\alpha c_{-\alpha} $ is an $ \alpha $-Weyl element. By \thmitemcref{A2Weyl:basecomp-cor}{A2Weyl:basecomp-cor:1}, the map in the assertion is simply the map
	\[ \map{}{\rootgr{\gamma}}{\rootgr{\beta}}{x_\gamma}{x_\gamma^{w_\alpha}} \]
	which is clearly a homomorphism. Further, it is bijective because $ w_\alpha $ is an $ \alpha $-Weyl element.
\end{proof}

As a consequence of the previous result, the root groups must be abelian if we have sufficiently many Weyl elements.

\begin{lemma}\label{ADE:abelian-lem}
	Let $ (\alpha, \beta, \gamma) $ be an $ A_2 $-triple and assume that $ \invset{\alpha} $ is non-empty. Then $ \rootgr{\beta} $ is abelian.
\end{lemma}
\begin{proof}
	We have $ \rootgr{\beta} = \commutator{\rootgr{\alpha}}{\rootgr{\gamma}} $ by \cref{A2Weyl:rootisom}. Since $ \rootgr{\beta} $ commutes with $ \rootgr{\alpha} $ and $ \rootgr{\gamma} $ by the commutator axiom, we infer that $ \rootgr{\beta} $ commutes with itself. In other words, $ \rootgr{\beta} $ is abelian.
\end{proof}

\begin{proposition}\label{ADE:abelian}
	Assume that $ \invset{\alpha} $ is non-empty for all roots $ \alpha $. Then all root groups of $ G $ are abelian.
\end{proposition}
\begin{proof}
	Every root lies in an $ A_2 $-triple by \cref{ADE:simply-laced:all-in-A2}, so the assertion follows from \cref{ADE:abelian-lem}.
\end{proof}

We are now in a position to prove many desirable properties of Weyl elements in $ G $. While we already know from \cref{basic:unique-weyl-ext-crit} that $ G $ has unique Weyl extensions, we will give a new proof of this fact to show how it can be derived from the previous results in this section.

\begin{proposition}\label{A2Weyl:weyl}
	Assume that $ \invset{\beta} $ is non-empty for all roots $ \beta $. Then for any root $ \alpha $, the following statements hold:
	\begin{proenumerate}
		\item \label{A2Weyl:weyl:leftright}Every $ \alpha $-Weyl triple is balanced.
		
		\item \label{A2Weyl:weyl:weylel-gives-triple}If $ (a_{-\alpha}, b_\alpha, c_{-\alpha}) $ and $ (a_{-\alpha}', b_\alpha', c_{-\alpha}') $ are two $ \alpha $-Weyl triples whose corresponding Weyl elements are identical, then $ (a_{-\alpha}, b_\alpha, c_{-\alpha}) = (a_{-\alpha}', b_\alpha', c_{-\alpha}') $.
		
		\item \label{A2Weyl:weyl:invel-gives-triple}$ G $ has unique $ \alpha $-Weyl extensions.
	\end{proenumerate}
\end{proposition}
\begin{proof}
	Choose roots $ \beta, \gamma $ such that $ (\alpha, \beta, \gamma) $ is an $ A_2 $-triple and choose invertible elements $ x_\beta \in \invset{\beta} $, $ x_\gamma \in \invset{\gamma} $. We will keep these choices throughout the whole proof.
	
	For~\itemref{A2Weyl:weyl:leftright}, let $ (a_{-\alpha}, b_\alpha, c_{-\alpha}) $ be an $ \alpha $-Weyl triple. By \cref{A2Weyl:left-right-comm}, we have $ \commutator{x_\beta}{c_{-\alpha}} = \commutator{x_\beta}{a_{-\alpha}} $, which implies that
	\[ \commutator{c_{-\alpha}}{x_\beta} = \commutator{x_\beta}{c_{-\alpha}}^{-1} = \commutator{x_\beta}{a_{-\alpha}}^{-1} = \commutator{a_{-\alpha}}{x_\beta}. \]
	Since $ (-\alpha, \gamma, \beta) $ is an $ A_2 $-triple and $ x_\beta $ is invertible, we know from \cref{A2Weyl:rootisom} that the map $ \map{}{\rootgr{-\alpha}}{\rootgr{\gamma}}{x_{-\alpha}}{\commutator{x_{-\alpha}}{x_\beta}} $ is an isomorphism. Thus it follows that $ a_{-\alpha} = c_{-\alpha} $, proving that every $ \alpha $-Weyl triple is weakly balanced. By \cref{basic:all-balanced}, it follows that every $ \alpha $-Weyl triple is balanced, so~\itemref{A2Weyl:weyl:leftright} holds.
	
	Now let $ (a_{-\alpha}, b_\alpha, c_{-\alpha}) $ and $ (a_{-\alpha}', b_\alpha', c_{-\alpha}') $ be two $ \alpha $-Weyl triples whose corresponding Weyl elements are identical. Denote their common Weyl element by $ w_\alpha $. Then by \cref{A2Weyl:basecomp-cor}, we have $ \commutator{x_\beta}{a_{-\alpha}} = x_\beta^{w_\alpha} = \commutator{x_\beta}{a_{-\alpha}'} $ and $ \commutator{x_\gamma}{b_\alpha} = x_\gamma^{w_\alpha} = \commutator{x_\gamma}{b_\alpha'} $. As in the proof of~\itemref{A2Weyl:weyl:leftright}, this implies that $ a_{-\alpha} = a_{-\alpha}' $ and $ b_\alpha = b_\alpha' $. It follows that $ c_{-\alpha} = (a_{-\alpha} b_\alpha)^{-1} w_\alpha = (a_{-\alpha}' b_\alpha')^{-1} w_\alpha = c_{-\alpha}' $. This finishes the proof of~\itemref{A2Weyl:weyl:weylel-gives-triple}.
	
	Now let $ b_\alpha \in \invset{\alpha} $. By the definition of $ \invset{\alpha} $ and by~\itemref{A2Weyl:weyl:leftright}, there exists $ a_{-\alpha} \in \rootgr{-\alpha} $ such that $ (a_{-\alpha}, b_\alpha, a_{-\alpha}) $ is an $ \alpha $-Weyl triple. To show the uniqueness of $ a_{-\alpha} $, we take another element $ a_{-\alpha}' \in \rootgr{-\alpha} $ such that $ (a_{-\alpha}', b_\alpha, a_{-\alpha}') $. (Note that $ a_{-\alpha} b_\alpha a_{-\alpha} $ and $ a_{-\alpha}' b_\alpha a_{-\alpha}' $ are allowed to be distinct.) By \thmitemcref{A2Weyl:basecomp-cor}{A2Weyl:basecomp-cor:2}, we have
	\[ \commutator[\big]{\commutator{x_\beta}{a_{-\alpha}}}{b_\alpha} = x_\beta^{-1} = \commutator[\big]{\commutator{x_\beta}{a_{-\alpha}'}}{b_{\alpha}}. \]
	Using the same arguments as in~\itemref{A2Weyl:weyl:leftright} and~\itemref{A2Weyl:weyl:weylel-gives-triple} and the fact that $ b_\alpha $ is invertible, we infer that $ \commutator{x_\beta}{a_{-\alpha}} = \commutator{x_\beta}{a_{-\alpha}'} $. Another application of the same argument, using that $ x_\beta $ is invertible, yields $ a_{-\alpha} = a_{-\alpha}' $. This finishes the proof of~\itemref{A2Weyl:weyl:invel-gives-triple}.
\end{proof}

\begin{remark}[Braid relations]\label{ADE:braid-rem}
	We could also use \cref{A2Weyl:basecomp-cor} to give an explicit proof of the braid relations for Weyl elements in $ G $, which are already known to hold by the general argument in \cref{braid:all}.
\end{remark}

Our next goal is to prove the square formula for Weyl elements (\cref{A2Weyl:cartan-comp}). To do this, we consider arbitrary roots $ \alpha $, $ \beta $ and compute for each possible \enquote{configuration} of $ (\alpha, \beta) $ how the squares of $ \alpha $-Weyl elements act on $ \rootgr{\beta} $. Here by \enquote{configurations} we essentially mean the different cases in \cref{ADE:cartan-char}.

\begin{reminder}[see \cref{act-by}]
	Let $ w $ be an element of $ G $ and let $ \alpha $ be a root. We say that \defemph*{$ w $ acts trivially on $ \rootgr{\alpha} $} if $ x_\alpha^w = x_\alpha $ for all $ x_\alpha \in \rootgr{\alpha} $ and we say that \defemph*{$ w $ acts on $ \rootgr{\alpha} $ by inversion} if $ x_\alpha^w = x_\alpha^{-1} $ for all $ x_\alpha \in \rootgr{\alpha} $.
\end{reminder}

The following proof uses the ideas from \cref{A2Weyl:spirit}.

\begin{lemma}\label{A2Weyl:square-act-lemma-A2}
	Let $ (\alpha, \beta, \gamma) $ be an $ A_2 $-triple and assume that $ w_\alpha $ is an $ \alpha $-Weyl element. Then $ w_\alpha^2 $ acts on $ \rootgr{\beta} $ and on $ \rootgr{\gamma} $ by inversion.
\end{lemma}
\begin{proof}
	Choose an $ \alpha $-Weyl triple $ (a_{-\alpha}, b_\alpha, a_{-\alpha}) $ such that $ w_\alpha = a_{-\alpha} b_\alpha a_{-\alpha} $ and let $ x_\beta \in \rootgr{\beta} $ be arbitrary. Put $ w \defl w_\alpha^2 $. Using the formulas from \cref{A2Weyl:basecomp-cor} two times, we see that
	\begin{align*}
		x_\beta^w &= \commutator{x_\beta}{a_{-\alpha}}^{w_\alpha} = \commutator{\commutator{x_\beta}{a_{-\alpha}}}{b_\alpha}.
	\end{align*}
	By \thmitemcref{A2Weyl:basecomp-cor}{A2Weyl:basecomp-cor:2}, this means that $ x_\beta^w = x_\beta^{-1} $, so $ w $ acts on $ \rootgr{\beta} $ by inversion. Since $ w_\alpha $ is also a $ (-\alpha) $-Weyl element and $ (-\alpha, \gamma, \beta) $ is an $ A_2 $-triple, it follows that $ w $ acts on $ \rootgr{\gamma} $ by inversion as well.
\end{proof}

\cref{A2Weyl:square-act-lemma-A2} was relatively straightforward to prove because we only had to apply the formulas from \cref{A2Weyl:basecomp-cor}. For the action of $ w_\beta $ on its \enquote{own} root group $ \rootgr{\beta} $, however, we have no such formula. Instead, we have use the fact that $ \rootgr{\beta} $ can be written as the commutator of two adjacent root groups, and then apply \cref{A2Weyl:square-act-lemma-A2}.

\begin{lemma}\label{A2Weyl:square-act-lemma-self}
	Assume that $ \invset{\alpha} $ is non-empty for all roots $ \alpha $ and let $ \beta $ be any root. Then for all $ \beta $-Weyl elements $ w_\beta $, $ w_\beta^2 $ acts trivially on $ \rootgr{\beta} $ and $ \rootgr{-\beta} $.
\end{lemma}
\begin{proof}
	By \cref{ADE:simply-laced:all-in-A2}, there exist roots $ \alpha, \gamma $ such that $ (\alpha, \beta, \gamma) $ is an $ A_2 $-pair. Let $ x_\beta \in \rootgr{\beta} $ be arbitrary and put $ w \defl w_\beta^2 $. By \cref{A2Weyl:rootisom}, there exist $ x_\alpha \in \rootgr{\alpha} $ and $ x_\gamma \in \rootgr{\gamma} $ such that $ x_\beta = \commutator{x_\alpha}{x_\gamma} $. Thus it follows from \cref{A2Weyl:square-act-lemma-A2,basic:comm-add} that
	\[ x_\beta^w = \commutator{x_\alpha^w}{x_\gamma^w} = \commutator{x_\alpha^{-1}}{x_\gamma^{-1}} = \commutator{x_\alpha}{x_\gamma} = x_\gamma. \]
	That is, $ w $ acts trivially on $ \rootgr{\beta} $. Since $ w_\beta $ is also a $ (-\beta) $-Weyl element by \thmitemcref{basic:weyl-general}{basic:weyl-general:minus}, this implies that $ w $ acts trivially on $ \rootgr{-\beta} $, too.
\end{proof}

\begin{lemma}\label{A2Weyl:ortho}
	Let $ \alpha, \beta $ be orthogonal roots and assume that $ w_\alpha $ is an $ \alpha $-Weyl element. Then $ w_\alpha^2 $ acts trivially on $ \rootgr{\beta} $.
\end{lemma}
\begin{proof}
	It follows from the assumptions that $ \beta $ is adjacent to $ \alpha $ and $ -\alpha $. This implies that $ \rootgr{\beta} $ commutes with $ \gen{\rootgr{\alpha}, \rootgr{-\alpha}} $ and thus with $ w_\alpha $. The assertion follows.
\end{proof}

The previous results cover all possible \enquote{configurations} of pairs of roots $ (\alpha, \beta) $. It is now a straightforward computation to check that the square formula is satisfied in each case.

\begin{proposition}\label{A2Weyl:cartan-comp}
	Assume that $ \invset{\gamma} $ is non-empty for all roots $ \gamma $ and let $ \alpha, \beta $ be two roots. Let $ w_\alpha $ be an $ \alpha $-Weyl element, put $ w \defl w_\alpha^2 $ and set $ \epsilon \defl (-1)^{\cartanint{\beta}{\alpha}} $ where $ \cartanint{\beta}{\alpha} = 2 \frac{\beta \cdot \alpha}{\alpha \cdot \alpha} $ is the Cartan integer for $ (\beta, \alpha) $. Then $ x_\beta^w = x_\beta^\epsilon $ for all $ x_\beta \in \rootgr{\beta} $.  In other words, $ G $ satisfies the square formula for Weyl elements (see \cref{param:square-formula-rgg-def}).
\end{proposition}
\begin{proof}
	We only have to put \cref{A2Weyl:square-act-lemma-A2,A2Weyl:square-act-lemma-self,A2Weyl:ortho} together and apply \cref{ADE:cartan-char}.
\end{proof}


\section{Standard Twisting Systems}

\label{sec:ADE:sttwist}

\begin{secnotation}
	We denote by $ \roots $ an irreducible simply-laced root system and by $ \chevstr = (\chevstr_{\alpha, \beta})_{\alpha, \beta \in \roots} $ a family of Chevalley structure constants of type $ \roots $ (in the sense of \cref{chev:struc-family-def}).
\end{secnotation}

In this section, we collect some properties of the Chevalley structure constants of simply-laced type and of the standard twisting systems which they define.

\begin{note}\label{ADE:stand-twist-note}
	In \ref{param:strategy:stsigns} and~\ref{param:strategy:admit-twist}, we declared that the \enquote{standard (partial) twisting system} and the \enquote{standard signs} for $ \roots $-graded groups are the ones which occur in some kind of \enquote{standard example}. For the root systems of types $ B $, $ BC $, $ C $ and $ F $, we will state explicit formulas for the involved signs and parity maps. In the setting of simply-laced root graded groups, however, this approach is not possibe because we want to cover three types of root systems at the same time. For instance, a parity map $ \inverparsym $ defined on $ A_n $ cannot be used to parametrise root graded groups of type $ D $ or $ E $.
	
	Instead, we take a more general approach in this chapter: We define not a single standard twisting system, but rather a twisting system for every family $ \chevstr $ of Chevalley structure constants of simply-laced type (\cref{ADE:standard-twist-sys-def}). Similarly, we will define coordinatisations of $ \roots $-graded groups with signs $ \chevstr $ in \cref{ADE:param-def}. In other words, we can cover all simply-laced root systems at the same time by working with the family $ \chevstr $ (which is not explicitly given). An additional advantage of this approach is that it makes the connection to Chevalley groups more evident.
	
	It should be noted that the strategy described above is only possible because simply-laced root graded groups have the exact same twisting systems and commutator relations as Chevalley groups. The more complicated root systems will require a different approach.
\end{note}

\begin{remark}\label{ADE:chev-C-choice}
	We will, throughout this section, often consider a Chevalley group $ H $ of type $ \roots $ over $ \IC $ with respect to $ \chevstr $, and we will denote its root groups, its root isomorphisms and the standard Weyl elements from \cref{chev:weyl-def} by $ (\rootgr{\alpha})_{\alpha \in \roots} $, $ (\risom{\alpha})_{\alpha \in \roots} $ and $ (w_\alpha)_{\alpha \in \roots} $, respectively. By choosing the complex numbers as a base ring, we ensure that the set $ \compactSet{\pm 1_\IZ} $ embeds into the base ring. Thus we can write $ \risom{\alpha}(a) $ for any $ a \in \compactSet{\pm 1_\IZ} $, and $ \risom{\alpha}(a) = \risom{\alpha}(b) $ implies $ a=b $ for all $ a,b \in \compactSet{\pm 1_\IZ} $ and all $ \alpha \in \roots $. This is the only property of the base ring that we need, so we could as well replace $ \IC $ by any commutative associative ring $ \ring $ with $ 2_\ring \ne 0_\ring $. However, the choice of the complex numbers has the additional benefit that the construction of Chevalley groups over $ \IC $ is more straightforward than over arbitrary rings.
\end{remark}

At first, we show how the parity map in simply-laced Chevalley groups can be computed from the family $ \chevstr $. Similar computations could be done for all other crystallographic root systems as well.

\begin{lemma}\label{ADE:parity-comp}
	Let $ \map{\inverparsym}{\roots \times \roots}{\compactSet{\pm 1_\IZ}}{}{} $ the Chevalley parity map for $ \chevstr $ from \cref{chev:parmap}. Then for all roots $ \alpha, \delta $, the value of $ \inverpar{\alpha}{\delta} $ is determined by $ \chevstr $ in the following way:
	\begin{lemenumerate}
		\item \label{ADE:parity-comp:self}$ \inverparbr{\alpha}{\alpha} = \inverparbr{\alpha}{-\alpha} = -1 $.
		
		\item $ \inverparbr{\alpha}{\delta} = 1 $ if $ \alpha $, $ \delta $ are orthogonal.
		
		\item $ \inverparbr{\alpha}{\delta} = \chevstr_{\alpha, \delta} $ if $ (\alpha, \delta) $ is an $ A_2 $-pair.
		
		\item $ \inverparbr{\alpha}{\delta} = -\chevstr_{\alpha, -\delta} $ if $ (\alpha, -\delta) $ is an $ A_2 $-pair.
	\end{lemenumerate}
\end{lemma}
\begin{proof}
	Let $ H $, $ (\rootgr{\alpha})_{\alpha \in \roots} $, $ (\risom{\alpha})_{\alpha \in \roots} $ and $ (w_\alpha)_{\alpha \in \roots} $ be as in \cref{ADE:chev-C-choice}. The first assertion holds by \cref{chev:const-on-itself}. If $ \alpha $ and $ \delta $ are orthogonal, then $ w_\delta $ commutes with $ \rootgr{\alpha} $ and so $ \inverparbr{\alpha}{\delta} = 1 $, which implies the second assertion. Now assume that $ (\alpha, \delta) $ is an $ A_2 $-pair, and put $ \beta \defl \alpha+\delta $. Then it follows from \thmitemcref{A2Weyl:basecomp-cor}{A2Weyl:basecomp-cor:1} that
	\[ \risom{\beta}(\inverpar{\alpha}{\delta}) = \risom{\alpha}(1)^{w_\delta} = \commutator{\risom{\alpha}(1)}{\risom{\delta}(1)} = \risom{\beta}(\chevstr_{\alpha, \delta}), \]
	so $ \inverpar{\alpha}{\delta} = \chevstr_{\alpha, \delta} $. Now assume that $ (\alpha, -\delta) $ is an $ A_2 $-pair, and put $ \gamma \defl \alpha-\delta $. Then $ (\delta, \alpha, \gamma) $ is an $ A_2 $-triple, so it follows from \thmitemcref{A2Weyl:basecomp-cor}{A2Weyl:basecomp-cor:2} that
	\[ \risom{\gamma}(\inverpar{\alpha}{\delta}) = \risom{\alpha}(1)^{w_\delta} = \commutator{\risom{\alpha}(1)}{\risom{-\delta}(-1)} = \risom{\gamma}(-\chevstr_{\alpha, -\delta}). \]
	Thus $ \inverpar{\alpha}{\delta} = -\chevstr_{\alpha, -\delta} $. This finishes the computation of $ \inverparsym $.
\end{proof}

\begin{example}\label{ADE:An-parmap}
	We consider the \enquote{standard} family of Chevalley structure constants for $ A_n $ from \cref{chev:An-chevbas}, which is given by
	\[ \chevstr_{\basvec_i - \basvec_j, \basvec_j - \basvec_k} = 1, \qquad \chevstr_{\basvec_j - \basvec_k, \basvec_i - \basvec_j} = -1 \]
	and $ \chevstr_{\alpha, \beta} = 0 $ for all other pairs $ \alpha, \beta $ of roots. Let $ i,j,k,l \in \numint{1}{n+1} $ such that $ i \ne j $ and $ k \ne l $. Write $ \chevstr_{ij,kl} $ for $ \chevstr_{\basvec_i - \basvec_j, \basvec_k - \basvec_l} $. We want to compute $ a \defl \inverpar{\basvec_i - \basvec_j}{\basvec_k - \basvec_l} $. By the first two assertions in \cref{ADE:parity-comp}, we have $ a=-1 $ if $ \abs{\compactSet{i,j} \intersect \compactSet{k,l}} = 2 $ and $ a = 1 $ if $ \compactSet{i,j} \intersect \compactSet{k,l} = \emptyset $. From now on, we assume that $ \abs{\compactSet{i,j} \intersect \compactSet{k,l}} = 1 $. Then
	\[ a = \begin{cases}
		\chevstr_{ij,jl} = 1 & \text{if } j=k, \\
		-\chevstr_{ij,jk} = -1 & \text{if } j=l, \\
		-\chevstr_{ij,li} = 1 & \text{if } i=k, \\
		\chevstr_{ij,ki} = -1 & \text{if } i=l.
	\end{cases} \]
	In summary, we have
	\[ \inverpar{\basvec_i - \basvec_j}{\basvec_k - \basvec_l} = \begin{cases}
		-1 & \text{if } l \in \compactSet{i,j}, \\
		1 & \text{otherwise.}
	\end{cases} \]
\end{example}

\begin{definition}[Standard twisting system]\label{ADE:standard-twist-sys-def}
	Let $ \rootbase $ be a root base of $ \roots $. The \defemph*{standard twisting system of type $ \roots $ defined by $ \chevstr $ with respect to $ \rootbase $}\index{twisting system!standard} is the pair $ (\twistgroup, \restrict{\inverparsym}{\roots \times \rootbase}) $ where $ \twistgroup \defl \compactSet{\pm 1_\IZ} $ and $ \inverparsym $ is the Chevalley parity map for $ \chevstr $ (which is given by the formulas in \cref{ADE:parity-comp}). If $ G $ is a group with a $ \roots $-pregrading, then the \defemph*{standard twisting system for $ G $ defined by $ \chevstr $ with respect to $ \rootbase $} is the same pair together with the additional information that $ \twistgroup $ acts on all root groups of $ G $ by inversion.
\end{definition}

\begin{lemma}\label{ADE:standard-is-twist}
	Let $ G $ be a group with a $ \roots $-pregrading, let $ \rootbase $ be a root base of $ \roots $ and let $ (w_\delta)_{\delta \in \rootbase} $ be a $ \rootbase $-system of Weyl elements in $ G $. Denote the standard twisting system for $ G $ defined by $ \chevstr $ with respect to $ \rootbase $ by $ (\twistgroup, \inverparsym) $. Then $ (\twistgroup, \inverparsym) $ is a twisting system for $ (G, (w_\delta)_{\delta \in \rootbase}) $. Further, $ \inverparsym $ is complete and adjacency-trivial and it satisfies the square formula.
\end{lemma}
\begin{proof}
	By \cref{param:pargroup-inv-example}, $ \twistgroup $ is a twisting group for $ (G, (w_\delta)_{\delta \in \rootbase}) $. We know from \thmitemcref{ADE:parity-comp}{ADE:parity-comp:self} that $ \parmoveset{\twistgroup}{\alpha}{\alpha} = \twistgroup $ for all roots $ \alpha $, so $ \inverparsym $ is complete. In particular, it is transporter-invariant by \cref{parmap:transport-invar-crit}. Using that the map $ \inverparsym $ is read off from the Chevalley group $ H $ in \cref{ADE:chev-C-choice} and that Weyl elements in this group satisfy the braid relations by \cref{braid:all}, it follows from  \thmitemcref{param:param-parmap-has-properties2}{param:param-parmap-has-properties2:braid} that $ \inverparsym $ is braid-invariant as well. Similarly, since $ H $ satisfies the square formula for Weyl elements by \cref{A2Weyl:cartan-comp}, it follows from \thmitemcref{param:param-parmap-has-properties2}{param:param-parmap-has-properties2:square} that $ \inverparsym $ satisfies the square formula. Finally, $ \inverparsym $ is adjacency-trivial by \thmitemcref{param:param-parmap-has-properties2}{param:param-parmap-has-properties2:adj-triv}. This finishes the proof.
\end{proof}

We end this section with two properties of $ \chevstr $ and $ \inverparsym $ which, at the moment, are completely unmotivated. \Cref{ADE:chev-assoc} will be used in the proof of the associativity of the coordinatising ring (\cref{ADE:assoc}). \Cref{ADE:chev-conj-formula} will be applied when we transport the commutator formula on one fixed $ A_2 $-pair $ (\alpha_0, \gamma_0) $ to all other $ A_2 $-pairs (\cref{ADE:comm-formula}).

\begin{lemma}\label{ADE:chev-assoc}
	Assume that $ \roots $ is of rank at least~3 and let $ \alpha, \beta, \gamma $ be roots such that $ \Set{\alpha, \beta, \gamma} $ is a root base of a parabolic $ A_3 $-subsystem of $ \roots $ and such that $ \alpha \cdot \gamma = 0 $. Then $ \chevstr_{\alpha, \beta} \chevstr_{\alpha+\beta, \gamma} = \chevstr_{\beta, \gamma} \chevstr_{\alpha, \beta+\gamma} $.
\end{lemma}
\begin{proof}
	By the Chevalley commutator formula (\cref{chev:comm-formula}), we have
	\begin{align*}
		\commutator{\risom{\alpha}(1)}{\risom{\gamma}(1)} &= 1_H \rightand \\
		\commutator[\big]{\commutator{\risom{\beta}(1)}{\risom{\alpha}(1)^{-1}}}{\commutator{\risom{\beta}(1)}{\risom{\gamma}(1)}} &\in \commutator{\rootgr{\alpha+\beta}}{\rootgr{\beta+\gamma}} = \compactSet{1_H}
	\end{align*}
	Thus it follows from the variant of the Hall-Witt identity given in \cref{hall-witt-special} that
	\[ \commutator[\big]{\commutator{\risom{\alpha}(1)}{\risom{\beta}(1)}}{\risom{\gamma}(1)} = \commutator[\big]{\risom{\alpha}(1)}{\commutator{\risom{\beta}(1)}{\risom{\gamma}(1)}}. \]
	By the Chevalley commutator formula, this means precisely that
	\[ \risom{\alpha+\beta+\gamma}(\chevstr_{\alpha, \beta} \chevstr_{\alpha+\beta, \gamma}) = \risom{\alpha+\beta+\gamma}(\chevstr_{\beta, \gamma} \chevstr_{\alpha, \beta+\gamma}). \]
	The assertion follows.
\end{proof}

\begin{lemma}\label{ADE:chev-conj-formula}
	Let $ \inverparsym $ be the Chevalley parity map for $ \chevstr $, let $ (\alpha, \gamma) $ be any $ A_2 $-pair and let $ \word{\delta} $ be any word over $ \roots $. Then $ \inverpar{\alpha}{\word{\delta}} \inverpar{\beta}{\word{\delta}} \chevstr_{\alpha^{\reflbr{\word{\delta}}}, \gamma^{\reflbr{\word{\delta}}}} = \inverpar{\alpha+\gamma}{\word{\delta}} \chevstr_{\alpha,\gamma} $.
\end{lemma}
\begin{proof}
	Let $ H $, $ (\rootgr{\alpha})_{\alpha \in \roots} $, $ (\risom{\alpha})_{\alpha \in \roots} $ and $ (w_\alpha)_{\alpha \in \roots} $ be as in \cref{ADE:chev-C-choice}. By the Chevalley commutator formula (\cref{chev:comm-formula}), we have 
	\[  \commutator{\risom{\alpha}(1)}{\risom{\gamma}(1)} = \risom{\alpha + \gamma}(\chevstr_{\alpha, \gamma}) \]
	Conjugating this equation by $ w_{\word{\delta}} $, we obtain
	\[  \commutator{\risom{\alpha^{\reflbr{\word{\delta}}}}(\inverpar{\alpha}{\word{\delta}})}{\risom{\gamma^{\reflbr{\word{\delta}}}}(\inverpar{\gamma}{\word{\delta}})} = \risom{\alpha^{\reflbr{\word{\delta}}} + \gamma^{\reflbr{\word{\delta}}}}(\chevstr_{\alpha, \gamma} \inverpar{\alpha+\gamma}{\word{\delta}}). \]
	However, we also have
	\[  \commutator{\risom{\alpha^{\reflbr{\word{\delta}}}}(\inverpar{\alpha}{\word{\delta}})}{\risom{\gamma^{\reflbr{\word{\delta}}}}(\inverpar{\gamma}{\word{\delta}})} = \risom{\alpha^{\reflbr{\word{\delta}}} + \gamma^{\reflbr{\word{\delta}}}}(\chevstr_{\alpha^{\reflbr{\word{\delta}}}, \gamma^{\reflbr{\word{\delta}}}} \inverpar{\alpha}{\word{\delta}} \inverpar{\gamma}{\word{\delta}}). \]
	Since $ \risom{\alpha^{\reflbr{\word{\delta}}} + \gamma^{\reflbr{\word{\delta}}}} $ is injective, the assertion follows.
\end{proof}


\section{Standard Signs}

\label{sec:ADE:stsigns}

\begin{secnotation}
	We denote by $ \roots $ an irreducible simply-laced root system of rank at least~2, by $ (\alpha_0, \gamma_0) $ an arbitrary $ A_2 $-pair in $ \roots $, by $ \chevstr = (\chevstr_{\alpha, \beta})_{\alpha, \beta \in \roots} $ a family of Chevalley structure constants of type $ \roots $ (in the sense of \cref{chev:struc-family-def}) and by $ G $ a group with a $ \roots $-pregrading $ (\rootgr{\alpha})_{\alpha \in \roots} $. We assume that there exists a coordinatisation of $ G $ by a ring $ \ring $ based at $ (\alpha_0, \gamma_0) $ with signs $ \chevstr $ (in the sense of the following \cref{ADE:param-def}), and we fix this coordinatisation.
\end{secnotation}

In this section, we define the notion of coordinatisations of $ G $ with signs $ \chevstr $, as explained in \cref{ADE:stand-twist-note}. Further, we will show that any ring which coordinatises $ G $ must be associative if $ \operatorname{rank}(\roots) \ge 3 $ and, if $ \roots $ is of type $ D $ or $ E $, even commutative. In particular, these properties are independent of the specific construction the we will undertake in \cref{sec:ADE:param}.

\begin{definition}[Coordinatisations]\label{ADE:param-def}
	Choose an $ A_2 $-pair $ (\alpha_0, \gamma_0) $ and let $ \ring $ be a ring. A \defemph*{coordinatisation of $ G $ by $ \ring $ based at $ (\alpha_0, \gamma_0) $ with signs $ \chevstr $}\index{coordinatisation of a root graded group!with standard signs!simply-laced type} is a family
	\[ (\map{\risom{\alpha}}{(\ring,+)}{\rootgr{\alpha}}{}{})_{\alpha \in \roots} \]
	of isomorphisms such that for all $ r,s \in \ring $ and all $ (\alpha_0, \gamma_0) $-positive $ A_2 $-pairs $ (\alpha,\gamma) $ (in the sense of \cref{ADE:pos-pair}), we have the commutator formula
	\[ \commutator{\risom{\alpha}(r)}{\risom{\gamma}(s)} = \risom{\alpha+\gamma}(\chevstr_{\alpha, \gamma} rs). \]
\end{definition}

\begin{remark}\label{ADE:DE-basepoint}
	If $ \roots $ is of type $ D $ or $ E $, then any $ A_2 $-pair is $ (\alpha_0, \gamma_0) $-positive for any other $ A_2 $-pair $ (\alpha_0, \gamma_0) $ by \thmitemcref{ADE:A2-pair-orbits}{ADE:A2-pair-orbits:DE}. Hence we can drop the \enquote{based at $ (\alpha_0, \gamma_0) $} part in \cref{ADE:param-def} for these root systems.
\end{remark}

\begin{remark}\label{ADE:param-comm-inv}
	Let $ (\alpha,\gamma) $ be an $ A_2 $-pair which is not $ (\alpha_0, \gamma_0) $-positive. Then $ (\gamma, \alpha) $ is $ (\alpha_0, \gamma_0) $-positive by \cref{ADE:A2-pair-orbits}. Since $ \chevstr_{\gamma, \alpha} = -\chevstr_{\alpha, \gamma} $ by \thmitemcref{chev:strconst-prop}{chev:strconst-prop:switch}, it follows that
	\[ \commutator{\risom{\alpha}(r)}{\risom{\gamma}(s)} = \commutator{\risom{\gamma}(s)}{\risom{\alpha}(r)}^{-1} = \risom{\alpha+\gamma}(-\chevstr_{\gamma, \alpha} sr) = \risom{\alpha+\gamma}(\chevstr_{\alpha, \gamma} sr). \]
	for all $ r,s \in \ring $. Thus non-positive $ A_2 $-pairs satisfy the same commutator relations as positive $ A_2 $-pairs except that $ rs $ is replaced by $ sr $. This is the key observation that we will use in \cref{ADE:comm} to show that $ \ring $ must be commutative if $ \roots $ is of type $ D $ or $ E $.
\end{remark}

\begin{example}\label{ADE:A-param-standard}
	Assume that $ \roots $ is the root system $ A_n $ in standard representation for some $ n \ge 2 $ and let $ \chevstr $ be the \enquote{standard} family of Chevalley structure constants from \cref{chev:An-chevbas}. Further, let $ (\alpha_0, \gamma_0) \defl (\basvec_i - \basvec_j, \basvec_j - \basvec_k) $ for some pairwise distinct $ i,j,k \in \numint{1}{n+1} $. Then a coordinatisation $ (\risom{\alpha})_{\alpha \in A_n} $ by a ring $ \ring $ based at $ (\alpha_0, \gamma_0) $ with signs $ \chevstr $ satisfies the same commutator relations as the group $ \E_{n+1}(\ring) $ in \cref{chev:ex-SL}:
	\[ \commutator{\risom{\basvec_i - \basvec_j}(r)}{\risom{\basvec_j - \basvec_k}(s)} = \risom{\basvec_i - \basvec_k}(rs) \midand \commutator{\risom{\basvec_j - \basvec_k}(r)}{\risom{\basvec_i - \basvec_j}(s)} = \risom{\basvec_i - \basvec_k}(-sr) \]
	for all pairwise distinct $ i,j,k \in \numint{1}{\ell+1} $ and all $ r,s \in \ring $. In this case, the set of $ (\alpha_0, \gamma_0) $-positive $ A_2 $-pairs (that is, the set of $ A_2 $-pairs which have $ rs $ in the commutator formula and not $ sr $) coincides with the set of $ A_2 $-pairs which have a positive sign in the commutator formula. However, for other choices of $ (\alpha_0, \gamma_0) $ and $ \chevstr $, these sets may be distinct.
\end{example}

Note that we have defined coordinatisations for arbitrary $ \roots $-pregradings. However, it is clear from the definition that any pregrading with a coordinatisation has $ \roots $-commutator relations. The following result says that pregradings with standard coordinatisations also have Weyl elements.

\begin{proposition}\label{ADE:stsign:weyl-char}
	Assume that $ G $ is rank-2-injective. For all roots $ \alpha $ and all invertible $ r \in \ring $, define
	\[ w_\alpha(r) \defl \risom{-\alpha}(-r^{-1}) \risom{\alpha}(r) \risom{-\alpha}(-r^{-1}). \]
	Then the following hold:
	\begin{lemenumerate}
		\item \label{ADE:stsign:weyl-char:bij}The maps
		\[ \map{}{\ringinvset{\ring}}{\invset{\alpha}}{r}{\risom{\alpha}(r)} \midand \map{}{\ringinvset{\ring}}{\weylset{\alpha}}{r}{w_\alpha(r)} \]
		are well-defined bijections for all roots $ \alpha $. Here $ \ringinvset{\ring} $, $ \invset{\alpha} $ and $ \weylset{\alpha} $ denote the sets of invertible elements in $ \ring $ (in the sense of \cref{ring:def-invertible}), the set of $ \alpha $-invertible elements in $ \rootgr{\alpha} $ and the set of $ \alpha $-Weyl elements, respectively.
		
		\item \label{ADE:stsign:weyl-char:weyl-ex}There exists an $ \alpha $-Weyl element for each root $ \alpha $.
		
		\item \label{ADE:stsign:weyl-char:param}Let $ \rootbase $ be any root base of $ \roots $ and denote by $ (\twistgroup, \inverparsym) $ the standard twisting system for $ G $ defined by $ \chevstr $ with respect to $ \rootbase $ (as in \cref{ADE:standard-twist-sys-def}). Then $ G $ is parametrised by $ (\twistgroup, \ring) $ with respect to $ \inverparsym $ and $ (w_\delta(1_\ring))_{\delta \in \rootbase} $.
	\end{lemenumerate}
\end{proposition}
\begin{proof}
	Let $ \alpha $ be any root and let $ w_\alpha = \risom{-\alpha}(r) \risom{\alpha}(s) \risom{-\alpha}(t) $ be an $ \alpha $-Weyl element for some $ r,s,t \in \ring $. Choose a root $ \gamma $ such that $ (\alpha, \gamma) $ is a $ (\alpha_0, \gamma_0) $-positive $ A_2 $-pair, and put $ \beta \defl \alpha + \gamma $. Then it follows from \thmitemcref{A2Weyl:basecomp-cor}{A2Weyl:basecomp-cor:1} that for all $ a \in \ring $, we have
	\begin{align*}
		\risom{\gamma}(a) &= \commutator[\big]{\commutator{\risom{\gamma}(a)}{\risom{\alpha}(s)}}{\risom{-\alpha}(t)} = \commutator{\risom{\beta}(\chevstr_{\gamma, \alpha} sa)}{\risom{-\alpha}(t)} = \risom{\gamma}\brackets[\big]{\chevstr_{\beta, -\alpha} \chevstr_{\gamma, \alpha} \brackets[\big]{t(sa)}}.
	\end{align*}
	By performing the same computation in a Chevalley group of type $ \roots $, one can show that $ \chevstr_{\beta, -\alpha} \chevstr_{\gamma, \alpha} = -1 $. Thus we infer that $ s $ is invertible with inverse $ -t $. In a similar way, one can show that the same assertion holds for $ t $ replaced by $ r $. Hence every Weyl element in $ G $ has the desired form. This says precisely that the second map in~\itemref{ADE:stsign:weyl-char:bij} is surjective (if it is well-defined).
	
	Now let $ \beta $ be an arbitrary root, let $ r \in \ring $ be invertible and consider $ w_\beta \defl w_\beta(r) $. By the same computation as in \cref{A2Weyl:basecomp-cor} and by similar ideas as in the previous paragraph, we see that $ \rootgr{\zeta}^{w_\beta} = \rootgr{\refl{\beta}(\zeta)} $ for all $ \zeta \in \roots \setminus \Set{\beta, -\beta} $. It remains to show the same statement for $ \zeta \in \Set{\beta, -\beta} $. Choose an $ (\alpha_0, \gamma_0) $-positive $ A_2 $-pair $ (\alpha, \gamma) $ such that $ \beta = \alpha+\gamma $. Then $ \rootgr{\beta} = \commutator{\rootgr{\alpha}}{\rootgr{\gamma}} $ because $ \risom{\beta}(a) = \commutator{\risom{\alpha}(\chevstr_{\alpha, \beta}^{-1} 1_\ring)}{\risom{\gamma}(a)} $ for all $ a \in \ring $. Similarly, $ \rootgr{-\beta} = \commutator{\rootgr{-\alpha}}{\rootgr{-\gamma}} $. It follows that
	\begin{align*}
		\rootgr{\pm \beta}^{w_\beta} &= \commutator{\rootgr{\pm \alpha}}{\rootgr{\pm \gamma}}^{w_\beta} = \commutator{\rootgr{\pm \alpha}^{w_\beta}}{\rootgr{\pm \gamma}^{w_\beta}} = \commutator{\rootgr{\mp \gamma}}{\rootgr{\mp \alpha}} = \rootgr{\mp \beta}.
	\end{align*}
	Hence $ w_\beta(r) $ is indeed a $ \beta $-Weyl element. This shows that the second map in~\itemref{ADE:stsign:weyl-char:bij} is well-defined. Further, it is injective by \thmitemcref{A2Weyl:weyl}{A2Weyl:weyl:weylel-gives-triple}. It follows that the first map in~\itemref{ADE:stsign:weyl-char:bij} is a well-defined bijection as well.
	
	Assertion~\itemref{ADE:stsign:weyl-char:weyl-ex} follows from~\itemref{ADE:stsign:weyl-char:bij} because $ \ringinvset{\ring} $ contains $ 1_\ring $. Finally, assertion~\itemref{ADE:stsign:weyl-char:param} can be proven by similar computations as above, using \cref{A2Weyl:basecomp-cor}.
\end{proof}

Observe that \thmitemcref{ADE:stsign:weyl-char}{ADE:stsign:weyl-char:param} says precisely that the Weyl elements $ (w_\alpha(1_\ring))_{\alpha \in \roots} $ satisfy the same conjugation formulas as in Chevalley groups (\cref{chev:squarerel}).

\begin{remark}\label{ADE:stsign:weyl-char:rem}
	If $ G $ is not rank-2-injective, then we cannot apply \cref{A2Weyl:basecomp-cor} in the proof of \cref{ADE:stsign:weyl-char}. However, it is still possible to perform the same computations as in \cref{A2Weyl:basecomp-cor,ADE:stsign:weyl-char} for the elements $ w_\alpha(r) $. Using the standard commutator relations, we obtain that $ w_\alpha(r) $ is an $ \alpha $-Weyl element for all $ \alpha \in \roots $ and all $ r \in \ringinvset{\ring} $ even if $ G $ is not rank-2-injective. In other words, the maps in \thmitemcref{ADE:stsign:weyl-char}{ADE:stsign:weyl-char:bij} are still well-defined, but it is no longer clear that they are bijective.
\end{remark}

\begin{remark}
	Assume that $ G $ is rank-2-injective. It follows from \cref{ADE:stsign:weyl-char} that $ G $ satisfies $ \invset{\alpha} = \rootgr{\alpha} \setminus \compactSet{1_G} $ (the additional condition of being an RGD-system) if and only if $ \ring $ is a division ring.
\end{remark}

\begin{example}\label{ADE:stsign:A-standard-conjformula}
	Assume that $ \roots $ is the root system $ A_n $ in standard representation for some $ n \ge 2 $ and let $ (\risom{\alpha})_{\alpha \in \roots} $ be a coordinatisation of $ G $ by a ring $ \ring $ with signs and base point chosen as in \cref{ADE:A-param-standard}. Then by \cref{ADE:An-parmap}, the conjugation formula in \thmitemcref{ADE:stsign:weyl-char}{ADE:stsign:weyl-char:param} says that
	\[ \risom{kl}(r)^{w_{ij}(1)} = \begin{cases}
		\risom{\tau(k), \tau(l)}(r) & \text{if } j \nin \compactSet{k,l}, \\
		\risom{\tau(k), \tau(l)}(-r) & \text{if } j \in \compactSet{k,l}
	\end{cases} \]
	for all $ i,j,k,l \in \numint{1}{n+1} $ with $ i \ne j $ and $ k \ne l $ and all $ r \in \ring $, where $ \tau $ is the transposition which interchanges $ i $ and $ j $.
\end{example}

\begin{proposition}\label{ADE:assoc}
	Let $ \ring $ be a ring which coordinatises $ G $ (with signs $ \chevstr $ and based at some $ A_2 $-pair $ (\alpha_0, \gamma_0) $) and assume that $ \roots $ is of rank at least~3. Then the ring $ \ring $ is associative.
\end{proposition}
\begin{proof}
	Choose roots $ \alpha, \beta, \gamma $ as in \cref{ADE:rank-3-has-A3}, and let $ r,s,t \in \ring $. By the same arguments as in the proof of~\cref{ADE:chev-assoc}, the Hall-Witt identity in the form of \cref{hall-witt-special} yields that
	\[ \commutator[\big]{\commutator{\risom{\alpha}(r)}{\risom{\beta}(s)}}{\risom{\gamma}(t)} = \commutator[\big]{\risom{\alpha}(r)}{\commutator{\risom{\beta}(s)}{\risom{\gamma}(t)}}. \]
	By the properties of the roots $ \alpha, \beta, \gamma $ from \cref{ADE:rank-3-has-A3}, this says precisely that
	\[ \risom{\alpha+\beta+\gamma}\brackets[\big]{\chevstr_{\alpha, \beta} \chevstr_{\alpha+\beta, \gamma} (r \rmult s) \rmult t} = \risom{\alpha+\beta+\gamma}\brackets[\big]{\chevstr_{\beta, \gamma} \chevstr_{\alpha, \beta+\gamma} r \rmult (s \rmult t)}. \]
	Since $ \chevstr_{\alpha, \beta} \chevstr_{\alpha+\beta, \gamma} = \chevstr_{\beta, \gamma} \chevstr_{\alpha, \beta+\gamma} $ by \cref{ADE:chev-assoc}, it follows that $ (r \rmult s) \rmult t = r \rmult (s \rmult t) $. Thus $ \ring $ is associative.
\end{proof}

\begin{remark}[Invertible elements are Moufang elements]
	In \cite[Theorem~13.8]{Faulkner-NonAssocProj}, it is shown that, in addition to the statement of \thmitemcref{ADE:stsign:weyl-char}{ADE:stsign:weyl-char:bij}, every invertible element $ a \in \ring $ is a \defemph{Moufang element}. This means that for all $ y,z \in \ring $, the Moufang identities hold:
	\begin{align*}
		a \brackets[\big]{y(az)} &= \brackets[\big]{a(ya)}z, \\
		\brackets[\big]{(za)y}a &= z \brackets[\big]{(ay)a}, \\
		(ay)(za) &= \brackets[\big]{a(yz)}a.
	\end{align*}
	By \cref{ADE:assoc}, this statement is only of interest if $ \roots $ is of rank 2, which means that $ \roots = A_2 $. We will see in \cref{ring:alternative-char} that a ring satisfies the Moufang identities for all $ a,y,z \in \ring $ if and only if it is an alternative ring. That is, $ \ring $ is alternative if and only if every element is a Moufang element. Since $ 0_\ring $ is always a Moufang element, it follows that $ \ring $ must be alternative if it is a division ring. In other words, any ring which coordinatises an RGD-system of type $ A_2 $ must be alternative. For general $ A_2 $-gradings, this is not known to be true. See also \cref{ADE:existence-A2}.
\end{remark}

\begin{proposition}\label{ADE:comm}
	Assume that $ \roots $ is of type $ D $ or $ E $. Then any ring $ \ring $ which coor\-dina\-tises $ G $ (with signs $ \chevstr $) is associative and commutative.
\end{proposition}
\begin{proof}
	Since root systems of type $ D $ are $ E $ have rank at least~4, the associativity of $ \ring $ follows from \cref{ADE:assoc}. For the commutativity, let $ r,s \in \ring $ and let $ (\alpha, \gamma) $ be any $ A_2 $-pair. By \cref{ADE:param-comm-inv}, the reversed $ A_2 $-pair $ (\gamma, \alpha) $ satisfies the commutator formula
	\[ \commutator[\big]{\risom{\gamma}(r)}{\risom{\alpha}(s)} = \risom{\alpha+\gamma}\brackets[\big]{\chevstr_{\gamma, \alpha} sr}. \]
	On the other hand, since $ (\gamma, \alpha) $ is itself positive with respect to any $ A_2 $-pair $ (\alpha_0, \gamma_0) $ by \thmitemcref{ADE:A2-pair-orbits}{ADE:A2-pair-orbits:DE}, we also have the usual commutator formula
	\[ \commutator[\big]{\risom{\gamma}(r)}{\risom{\alpha}(s)} = \risom{\alpha+\gamma}\brackets[\big]{\chevstr_{\gamma, \alpha} rs}. \]
	It follows that $ rs = sr $, as desired.
\end{proof}

\begin{note}
	\cref{ADE:assoc,ADE:comm} are not specific to standard signs: By \cref{param:strategy:parmap-choice}, every coordinatisation of $ G $ with non-standard signs can be obtained from a standard coordinatisation via twisting, and this process does not change the coordinatising ring. Thus \cref{ADE:assoc,ADE:comm} hold for non-standard coordinatisations as well.
\end{note}

\begin{remark}[A universal root graded group]\label{universal:ADE}
	Let $ \varring $ be an associative ring and assume that $ \roots = A_n $ for some $ n \ge 3 $. Denote by $ \hat{G}(\varring) $ the abstract group which is generated by symbols $ \steinhom{\alpha}(r) $ for $ \alpha \in \roots $ and $ r \in \varring $ with respect to the two kinds of relations: Firstly, $ \steinhom{\alpha}(r+s) = \steinhom{\alpha}(r) \steinhom{\alpha}(s) $ for all $ \alpha \in \roots $, $ r,s \in \varring $, and secondly, precisely the commutator relations that are required in \cref{ADE:param-def}. (Note that this is the non-commutative analogue of the Steinberg group in \cref{rgg-lit:steinberg-def}.) Put $ \steinrootgr{\alpha} \defl \Set{\steinhom{\alpha}(r) \given r \in \varring} $ for all $ \alpha \in \roots $.
	
	We want to show that $ (\hat{G}(\varring), (\steinrootgr{\alpha})_{\alpha \in \roots}) $ is a $ \roots $-graded group and that $ (\steinhom{\alpha})_{\alpha \in \roots} $ is a coordinatisation of $ \hat{G}(\varring) $ by $ \varring $ with signs $ \chevstr $ based at $ (\alpha_0, \gamma_0) $. By \cref{ADE:stsign:weyl-char:rem}, $ \hat{G}(\varring) $ has Weyl elements, and it has the desired commutator relations by definition. It remains to verify Axiom~\thmitemref{rgg-def}{rgg-def:nondeg} and that the maps $ (\steinhom{\alpha})_{\alpha \in \roots} $ are injective. To do this, we choose any $ \roots $-graded group $ (G, (\rootgr{\alpha})_{\alpha \in \roots}) $ which has a coordinatisation $ (\risom{\alpha})_{\alpha \in \roots} $ by $ \varring $ with signs $ \chevstr $ based at $ (\alpha_0, \gamma_0) $. We know that such a group exists: We can take the elementary group from \cref{chev:ex-SL} for the \enquote{standard} choice of $ \chevstr $, and an appropriate twisted coordinatisation of this group (in the sense of \cref{param:parmap-twist}) in the general case. Since $ G $ satisfies the defining relations of $ \hat{G}(\varring) $, there exists a unique group epimorphism $ \map{\pi}{\hat{G}(\varring)}{G}{}{} $ such that $ \pi(\steinhom{\alpha}(r)) = \risom{\alpha}(r) $ for all $ r \in \varring $ and $ \alpha \in \roots $. This implies that the maps $ (\steinhom{\alpha})_{\alpha \in \roots} $ are injective and that $ \hat{G}(\varring) $ satisfies Axiom~\thmitemref{rgg-def}{rgg-def:nondeg} (because $ G $ has these properties). Hence $ (\hat{G}(\varring), (\steinrootgr{\alpha})_{\alpha \in \roots}) $ is a $ \roots $-graded group and $ (\steinhom{\alpha})_{\alpha \in \roots} $ is a coordinatisation of $ \hat{G}(\varring) $ by $ \varring $ with signs $ \chevstr $ based at $ (\alpha_0, \gamma_0) $.
	
	Our main result on $ \roots $-graded groups (\cref{ADE:thm}) says that every $ \roots $-graded group is a quotient of $ \hat{G}(\varring) $ for some associative ring $ \varring $. For this reason, we call $ (\hat{G}(\varring), (\steinrootgr{\alpha})_{\alpha \in \roots}) $ the \defemph*{universal $ A_n $-graded group over $ \varring $ (with signs $ \chevstr $ based at $ (\alpha_0, \gamma_0) $)}\index{root graded group!universal}. Note that $ \hat{G}(\varring) $ is functorial in $ \varring $: For every homomorphism $ \map{\phi}{\varring}{\varring'}{}{} $ of rings, we have an induced homomorphism $ \map{\hat{G}(\phi)}{\hat{G}(\varring)}{\hat{G}(\varring')}{}{} $ of $ A_n $-graded groups (with $ \hat{G}(\phi)(\steinhom{\alpha, \varring}(r)) = \steinhom{\alpha, \varring'}(\phi(r)) $ for all $ \alpha \in A_n $ and all $ r \in \varring $).
	
	The same arguments as above apply to root systems of type $ D $ or $ E $, except that the ring $ \varring $ has to be commutative and that we have to choose a different group $ (G, (\rootgr{\alpha})_{\alpha \in \roots}) $ in our proof. (We can, for example, choose a Chevalley group over $ \varring $.) In fact, similar arguments and constructions work for all root systems for which the existence problem is completely solved. See also \cref{universal:B,universal:BC}.
\end{remark}


\section{The Coordinatisation}

\label{sec:ADE:param}

\begin{secnotation}\label{secnot:ADE:param}
	We denote by $ \roots $ any irreducible simply-laced root system of rank at least~$ 2 $, by $ \rootbase $ a root base of $ \roots $ and by $ G $ a group which has $ \roots $-commutator relations with root groups $ (\rootgr{\alpha})_{\alpha \in \roots} $ and which satisfies $ \invset{\alpha} \ne \emptyset $ for all roots $ \alpha $. We assume that all for all non-proportional roots $ \alpha $, $ \beta $, we have $ \rootgr{\alpha} \intersect \rootgr{\beta} = \compactSet{1_G} $. Further, we choose a family $ \chevstr = (\chevstr_{\alpha, \beta})_{\alpha, \beta \in \roots} $ of Chevalley structure constants of type $ \roots $ and we denote by $ (\twistgroup, \inverparsym) $ the standard twisting system for $ G $ defined by $ \chevstr $ (in the sense of \cref{ADE:standard-twist-sys-def}). From \cref{ADE:ring-mult-const} on, we will also fix an arbitrary $ A_2 $-pair $ (\alpha_0, \gamma_0) $.
\end{secnotation}

\begin{reminder}
	We know from \cref{ADE:standard-is-twist} that $ \inverparsym $ is complete and adjacency-trivial.
\end{reminder}

In this section, we can finally construct the coordinatising ring for $ G $. We begin by verifying that the assumptions of the parametrisation theorem are satisfied.

\begin{lemma}
	$ G $ is stabiliser-compatible with respect to $ \inverparsym $ and $ (w_\delta)_{\delta \in \roots} $.
\end{lemma}
\begin{proof}
	By \cref{ADE:ortho-adjacent}, any pair of orthogonal roots in a simply-laced root system is adjacent. Further, we know that $ \inverparsym $ is adjacency-trivial. Thus the assertion follows from \cref{param:adj-implies-stab}.
\end{proof}

\begin{lemma}
	$ G $ is square-compatible with respect to $ \inverparsym $ and $ (w_\delta)_{\delta \in \roots} $.
\end{lemma}
\begin{proof}
	Since $ G $ satisfies the square formula for Weyl elements by \cref{A2Weyl:cartan-comp} and $ \inverparsym $ satisfies the square formula, this follows from \cref{param:square-formula-comp}.
\end{proof}

\begin{construction}[of the coordinatising ring]\label{ADE:ring-const}
	Having verified that all the conditions in \cref{param:partwist-conv} are satisfied, we can now apply the parametrisation theorem (\cref{param:thm}). We conclude that there exist a group $ (\ring, +) $ and a parametrisation $ (\risom{\alpha})_{\alpha \in \rootbase} $ of $ G $ by the parameter system $ \calP = (\twistgroup, \ring) $ with respect to $ \inverparsym $ and $ (w_\delta)_{\delta \in \roots} $. By \cref{ADE:abelian}, the root groups of $ G $ are abelian, so $ (\ring, +) $ is abelian as well. Further, the action of $ \twistgroup $ on the set $ \ring $ is given by the equation
	\[ a.\risom{\alpha}(r) = \risom{\alpha}(a.r) \]
	for all $ a \in \twistgroup $, all roots $ \alpha $ and all $ r \in \ring $. It follows that $ 1_\twistgroup.r = r $ and $ (-1_\twistgroup).r = -r $ for all $ r \in \ring $. Since $ (\ring,+) $ is abelian, this implies that $ \twistgroup $ acts on $ \ring $ by group automorphisms, so that $ a.(r+s) = (a.r) + (a.s) $ for all $ a \in \twistgroup $ and all $ r,s \in \ring $.
\end{construction}

\begin{construction}[of the ring multiplication]\label{ADE:ring-mult-const}
	For the remainder of this section, we fix an arbitrary $ A_2 $-pair $ (\alpha_0, \gamma_0) $. Since $ \chevstr_{\alpha_0, \gamma_0} $ lies in $ \compactSet{\pm 1_\IZ} $ by \cref{chev:simply-laced-struct-1}, we can regard it as an element of $ \twistgroup $ and thus let it act on $ \ring $. With this convention, we now define a multiplication on $ \ring $ by
	\[ \map{\rmult}{\ring \times \ring}{\ring}{(a,b)}{\chevstr_{\alpha_0, \gamma_0}. \risom{\alpha_0+\gamma_0}^{-1}\brackets[\big]{\commutator{\risom{\alpha_0}(a)}{\risom{\gamma_0}(b)}}}. \]
	Put differently, we define $ \rmult $ to be the unique multiplication on $ \ring $ which satisfies
	\[ \commutator{\risom{\alpha_0}(a)}{\risom{\gamma_0}(b)} = \risom{\alpha_0 + \gamma_0}(\chevstr_{\alpha_0, \gamma_0}. (a \rmult b)) \]
	for all $ a,b \in \ring $. (Here we have used that $ \chevstr_{\alpha_0, \gamma_0}^{-1} = \chevstr_{\alpha_0, \gamma_0} $.)
	Since $ \chevstr_{\alpha_0, \gamma_0, 1, 1} = \chevstr_{\alpha_0,\gamma_0} $ by \cref{chev:struc-const-formulas}, this is the same formula that is satisfied in Chevalley groups, see \cref{chev:ex-SL}. We will later show that this multiplication turns $ \ring $ into a ring.
\end{construction}

\begin{note}
	A priori, the action of $ \twistgroup $ on $ \ring $ is defined independently of the multiplication on $ \ring $, and it does not rely on the existence of a \enquote{unit element $ 1_\ring $}. However, we will soon show that $ \ring $ does indeed have a unit element $ 1_\ring $, and then it is clear that $ a.r = i(a) \rmult r $ for all $ a \in \twistgroup $ and $ r \in \ring $ where $ i(1_\twistgroup) \defl 1_\ring $ and $ i(-1_\twistgroup) \defl -1_\ring $.
\end{note}

It remains to verify two properties: Firstly, that $ \ring $ is a ring. Secondly, that the commutator formula is satisfied for all $ (\alpha_0, \gamma_0) $-positive $ A_2 $-pairs and not merely for $ (\alpha_0, \gamma_0) $. We begin with the distributive laws. They are proven in essentially the same way in \cite[(2.19)]{Shi1993} and \cite[(19.7)]{MoufangPolygons}.

\begin{lemma}[Distributive laws]\label{ADE:distributivity}
	The multiplication on $ \ring $ satisfies
	\[ a \rmult (b+c) = a \rmult b + a \rmult c \midand (a+b) \rmult c = a \rmult c + b \rmult c \]
	for all $ a,b,c \in \ring $. In particular, we have $ (-a) \rmult b = -(a \rmult b) = a \rmult (-b) $ for all $ a,b \in \ring $.
\end{lemma}
\begin{proof}
	Let $ a,b,c \in \ring $ and denote by $ (\alpha,\gamma) \defl (\alpha_0, \gamma_0) $ the fixed $ A_2 $-pair from \cref{ADE:ring-mult-const}. By \cref{basic:comm-add}, we have
	\[ \commutator{\risom{\alpha}(a+b)}{\risom{\gamma}(c)} = \commutator{\risom{\alpha}(a) \risom{\alpha}(b)}{\risom{\gamma}(c)} = \commutator{\risom{\alpha}(a)}{\risom{\gamma}(c)} \commutator{\risom{\alpha}(b)}{\risom{\gamma}(c)}. \]
	We conclude that
	\begin{align*}
		\risom{\alpha+\gamma}\brackets[\big]{\chevstr_{\alpha,\gamma}.((a+b) \rmult c)} &= \commutator{\risom{\alpha}(a+b)}{\risom{\gamma}(c)} = \commutator{\risom{\alpha}(a)}{\risom{\gamma}(c)} \commutator{\risom{\alpha}(b)}{\risom{\gamma}(c)} \\
		&= \risom{\alpha+\gamma}(\chevstr_{\alpha, \gamma}.(a \rmult c)) \risom{\alpha+\gamma}(\chevstr_{\alpha, \gamma}.(b \rmult c)) \\
		&= \risom{\alpha+\gamma}\brackets[\big]{\chevstr_{\alpha, \gamma}.(a \rmult c) + \chevstr_{\alpha, \gamma}.(b \rmult c)}.
	\end{align*}
	Thus
	\[ \chevstr_{\alpha,\gamma}.((a+b) \rmult c) = \chevstr_{\alpha, \gamma}.(a \rmult c) + \chevstr_{\alpha, \gamma}.(b \rmult c). \]
	Since the right-hand side equals $ \chevstr_{\alpha, \gamma}. ((a \rmult c) + (b \rmult c)) $ by \cref{ADE:ring-const}, it follows that $ (a+b) \rmult c = (a \rmult c) + (b \rmult c) $. The right distributive law can be proven in a similar way.
\end{proof}

\begin{note}
	From now on, we will simply write $ ar $ in place of $ a.r $ for $ a \in \twistgroup $ and $ r \in \ring $. Since $ a.(r \rmult s) = (a.r) \rmult s = r \rmult (a.s) $ for all $ r,s \in \ring $ and all $ a \in \twistgroup $ by \cref{ADE:distributivity}, this will not cause any confusion.
\end{note}

\begin{lemma}\label{ADE:comm-formula}
	Let $ (\alpha, \gamma) $ be an $ A_2 $-pair which lies in the orbit of $ (\alpha_0, \gamma_0) $ under the Weyl group of $ \roots $ (where $ (\alpha_0, \gamma_0) $ is the fixed $ A_2 $-pair from \cref{ADE:ring-mult-const}). Then we have
	\[ \commutator{\risom{\alpha}(r)}{\risom{\gamma}(s)} = \risom{\alpha + \gamma}(\chevstr_{\alpha, \gamma} r \rmult s) \midand \commutator{\risom{\gamma}(s)}{\risom{\alpha}(r)} = \risom{\alpha+\gamma}(\chevstr_{\gamma, \alpha} r \rmult s) \]
	for all $ r,s \in \ring $. In particular, if $ (\alpha_0, \gamma_0) $ and $ (\alpha_0', \gamma_0') $ are two $ A_2 $-pairs which are conjugate under the Weyl group, then replacing $ (\alpha_0, \gamma_0) $ by $ (\alpha_0', \gamma_0') $ in \cref{ADE:ring-mult-const} does not change the multiplication on $ \ring $.
\end{lemma}
\begin{proof}
	By \cref{ADE:ring-mult-const}, the first formula is satisfied for $ (\alpha, \gamma) = (\alpha_0, \gamma_0) $. Now let $ u $ be an arbitrary element of the Weyl group of $ \roots $ and put $ \alpha \defl \alpha_0^u $, $ \beta \defl \beta_0^u $. Since the simple reflections generate the Weyl group (\cref{rootsys:simple-gen-weyl}), we can choose a word $ \word{\delta} $ over $ \rootbase $ such that $ u = \refl{\word{\delta}} $, and we put $ w \defl w_{\word{\delta}} $. Let $ r,s \in \ring $ be arbitrary. Conjugating the equation
	\[ \commutator{\risom{\alpha_0}(r)}{\risom{\gamma_0}(s)} = \risom{\alpha_0 + \gamma_0}(\chevstr_{\alpha_0, \gamma_0} r \rmult s) \]
	by $ w $, we obtain that
	\[ \commutator{\risom{\alpha}(\inverpar{\alpha_0}{\word{\delta}}r)}{\risom{\gamma}(\inverpar{\gamma_0}{\word{\delta}}s)} = \risom{\alpha + \gamma}(\inverpar{\alpha_0 + \gamma_0}{\word{\delta}}\chevstr_{\alpha_0, \gamma_0} r \rmult s). \]
	By \cref{ADE:chev-conj-formula}, this says that
	\[ \commutator{\risom{\alpha}(\inverpar{\alpha_0}{\word{\delta}}r)}{\risom{\gamma}(\inverpar{\gamma_0}{\word{\delta}}s)} = \risom{\alpha + \gamma}(\inverpar{\alpha_0}{\word{\delta}} \inverpar{\gamma_0}{\word{\delta}} \chevstr_{\alpha, \gamma} r \rmult s). \]
	Replacing $ r $ by $ \inverpar{\alpha_0}{\word{\delta}}^{-1}r $ and $ s $ by $ \inverpar{\gamma_0}{\word{\delta}}^{-1}s $, this gives us the desired first formula.
	
	For the second formula, observe that by \thmitemcref{group-rel}{group-rel:inv}, the first formula yields that
	\[ \commutator{\risom{\gamma}(s)}{\risom{\alpha}(r)} = \commutator{\risom{\alpha}(r)}{\risom{\gamma}(s)}^{-1} = \risom{\alpha+\gamma}(-\chevstr_{\alpha, \gamma}. r \rmult s). \]
	Since $ \chevstr_{\gamma, \alpha} = -\chevstr_{\alpha, \gamma} $ by \thmitemcref{chev:strconst-prop}{chev:strconst-prop:switch}, the second formula follows.
\end{proof}

\begin{construction}[of the identity element]\label{ADE:1-const}
	Fix an arbitrary simple root $ \delta_0 $. By \cref{A2Weyl:weyl}, there exist unique $ a_{-\delta_0} \in \rootgr{-\delta_0} $ and $ b_{\delta_0} \in \rootgr{\delta_0} $ such that $ w_{\delta_0} = a_{-\delta_0} b_{\delta_0} a_{-\delta_0} $. We can thus define an element $ 1_\ring \defl \risom{\delta_0}^{-1}(b_{\delta_0}) $ which, a priori, depends on the choice of $ \delta_0 $.
\end{construction}

\begin{lemma}[{\cite[(19.9)]{MoufangPolygons}}]\label{ADE:ring-identity}
	We have $ r \rmult 1_\ring = r $ and $ 1_\ring \rmult r = r $ for all $ r \in \ring $. In particular, $ 1_\ring $ does not depend on the choice of $ \delta_0 $ in \cref{ADE:1-const}.
\end{lemma}
\begin{proof}
	Let $ r \in \ring $ and write $ w_{\delta_0} = a_{-\delta_0} b_{\delta_0} a_{-\delta_0} $ as in \cref{ADE:1-const}. Since the Weyl group acts transitively on $ \roots $, there exists a root $ \rho $ such that $ (\rho, \delta_0) $ is an $ A_2 $-pair which lies in the orbit of $ (\alpha_0, \gamma_0) $. Recall that $ \refl{\delta_0}(\rho) = \rho + \delta_0 $. Now it follows from \cref{ADE:comm-formula,A2Weyl:basecomp-cor} that
	\begin{align*}
		\risom{\rho + \delta_0}(\chevstr_{\rho, \delta_0} r \rmult 1_\ring) &= \commutator{\risom{\rho}(r)}{\risom{\delta_0}(1_\ring)} = \commutator{\risom{\rho}(r)}{b_{\delta_0}} \\
		&= \risom{\rho}(r)^{w_{\delta_0}} = \risom{\rho + \delta_0}(\inverpar{\rho}{\delta_0} r).
	\end{align*}
	Since $ (\rho, \delta_0) $ is an $ A_2 $-pair, we know from \cref{ADE:parity-comp} that $ \inverpar{\rho}{\delta_0} = \chevstr_{\rho, \delta_0} $. Thus $ r \rmult 1_\ring = r $.
	
	The proof of the identity $ 1_\ring \rmult r = r $ is similar, but with a twist. This time we choose a root $ \xi $ such that $ (\delta_0, \xi) $ is an $ A_2 $-pair in the orbit of $ (\alpha_0, \gamma_0) $. Then again by \cref{ADE:comm-formula,A2Weyl:basecomp-cor}, we have
	\begin{align*}
		\risom{\delta_0 + \xi}\brackets{\chevstr_{\delta_0, \xi} 1_\ring \rmult r} &= \commutator{\risom{\delta_0}(1_\ring)}{\risom{\xi}(r)} = \commutator{\risom{\xi}(r)}{b_{\delta_0}}^{-1} \\
		&= (\risom{\xi}(r)^{w_{\delta_0}})^{-1} = \risom{\xi+\delta_0}(-\inverpar{\xi}{\delta_0} r).
	\end{align*}
	As above, we know from \cref{ADE:parity-comp} that $ \inverpar{\xi}{\delta_0} = \chevstr_{\xi, \delta_0} $. Since $ \chevstr_{\xi, \delta_0} = -\chevstr_{\delta_0, \xi} $ by \thmitemcref{chev:strconst-prop}{chev:strconst-prop:switch}, we infer that $ 1_\ring \rmult r = r $.
	
	Now assume that $ 1_\ring' $ is the element which is defined by the choice of another simple root $ \delta_0' $. Then it follows from the previous results that $ 1_\ring' = 1_\ring' \rmult 1_\ring = 1_\ring $, which finishes the proof.
\end{proof}

\begin{proposition}\label{ADE:ring}
	$ (\ring, +, \rmult) $ is a ring.
\end{proposition}
\begin{proof}
	This is a consequence of \cref{ADE:distributivity,ADE:ring-identity}.
\end{proof}

We can now state our main result for simply-laced root graded groups. For simplicity, we phrase it for $ \roots $-graded groups (which have to satisfy Axiom~\thmitemref{rgg-def}{rgg-def:nondeg}), but we have seen that it is actually true for any group $ G $ as in \cref{secnot:ADE:param}. Recall further from \cref{ADE:DE-basepoint} that the \enquote{base point} of a coordinatisation is only relevant if $ \roots $ is of type $ A $ and that in this case, it is possible to choose $ (\alpha_0, \gamma_0) $ and $ \chevstr $ in a way which yields the simple commutator formulas in \cref{ADE:A-param-standard}.

\begin{theorem}[Coordinatisation theorem for $ A_n $, $ D_n $, $ E_n $]\label{ADE:thm}
	Let $ \roots $ be an irreducible simply-laced root system of rank at least~$ 2 $ and let $ (G, (\rootgr{\alpha})_{\alpha \in \roots}) $ be a $ \roots $-graded group. Choose any $ A_2 $-pair $ (\alpha_0, \gamma_0) $ and any family $ \chevstr = (\chevstr_{\alpha, \beta})_{\alpha, \beta \in A_n} $ of Chevalley structure constants (in the sense of \cref{chev:struc-family-def}). Then there exist a ring $ \ring $ and a coordinatisation of $ G $ by $ \ring $ based at $ (\alpha_0, \gamma_0) $ with signs $ \chevstr $. If $ \roots $ is of rank at least 3, then $ \ring $ is necessarily associative. If $ \roots $ is of type $ D $ or $ E $, then $ \ring $ is necessarily associative and commutative. Further, if we fix a $ \rootbase $-system of Weyl elements in $ G $, then we can choose the root isomorphisms $ (\risom{\alpha})_{\alpha \in \roots} $ so that $ w_\delta = \risom{-\delta}(-1_\ring) \risom{\delta}(1_\ring) \risom{-\delta}(-1_\ring) $ for all $ \delta \in \rootbase $.
\end{theorem}
\begin{proof}
	We have constructed $ \ring $ in \cref{ADE:ring-const,ADE:ring-mult-const}, shown that it is a ring in \cref{ADE:ring} and proven the commutator formulas in \cref{ADE:comm-formula}. By \cref{ADE:assoc,ADE:comm}, the assertions about the associativity and commutativity of $ \ring $ hold. By \cref{ADE:ring-identity}, the Weyl elements have the desired form.
\end{proof}

	\chapter{The Blueprint Technique}
	
	\label{chap:blue}
	
	In \cite{BuildBuildings}, Ronan-Tits proved the existence of a large class of thick buildings by constructing them from so-called \enquote{blueprints}. Here a blueprint is, essentially, a local parameter system for a building which encompasses certain commutator relations in the rank-2 residues. A blueprint is said to be realisable if there exists a building which conforms to it, and Ronan-Tits give a simple realisability criterion for blueprints. This criterion involves computations with certain rewriting rules which are induced by a self-homotopy of the longest element in the Coxeter group. For example, they show that a blueprint of type $ A_3 $ is realisable if and only if it is associative (in a suitable sense). The corresponding computation is the same one that we will perform in \cref{sec:A3-blue}.
	
	The blueprint technique is a novel method which consists of performing the same kind of computations as in \cite{BuildBuildings} in the context of root graded groups. It makes sense for arbitrary root systems, but it will become apparent in this chapter that it only produces meaningful results if the root system is of rank at least 3. While the computations in the blueprint technique follow the same algorithm as in \cite{BuildBuildings}, the underlying logic is, in a certain sense, inverted: Ronan-Tits start from explicit commutator relations and then perform computations to derive conditions on these relations. These extra conditions guarantee that a building with the given commutator relations exist. On the other hand, we start with an existing (arbitrary) root graded group $ G $ (which one should think of as a \enquote{generalised building}) with a parametrisation $ (\risom{\alpha})_{\alpha \in \roots} $ (in the sense of \cref{param:param-def}) and then perform the blueprint computations to obtain information about the commutator relations in $ G $.
	
	The initial idea to \enquote{turn around} the computations of Ronan-Tits in this way was suggested by Bernhard Mühlherr. Originally, the only objective of this approach was to prove that coordinatising rings for $ C_3 $-graded groups are alternative. It later turned out that the blueprint technique can be used in a more general way to not only prove identities in the coordinatising algebraic structure (such as the alternative law), but to actually equip the parametrising groups with an algebraic structure. See \cref{blue:original} for more details.
	
	We begin this chapter with a section which describes the general idea of the blueprint technique in a slightly informal way. In the following section, we state and prove the results which make this technique work. In the third section, we show how the blueprint technique can be used to obtain a new proof of the fact that the coordinatising rings for root graded groups of type $ A $ in rank at least~3 are associative. We end this chapter with some concluding remarks which, in \cref{blue:summary}, also contain a summary of the blueprint technique. Later in \cref{chap:B,chap:BC}, we will apply the blueprint technique to root graded groups of types $ B $ and $ (B)C $, respectively. It should be noted that the specific computations for types $ B $ and $ (B)C $ have not been performed in \cite{BuildBuildings}.
	

\section{The Idea of the Blueprint Technique}

\label{sec:blue:idea}

\begin{secnotation}
	We denote by $ \roots $ a root system of rank $ n $, by $ \rootbase =(\delta_1, \ldots, \delta_n) $ a rescaled ordered root base of $ \roots $ and by $ W \defl \Weyl(\roots) $ the Weyl group. Further, we denote by $ G $ a $ \roots $-graded group with root groups $ (\rootgr{\alpha})_{\alpha \in \roots} $.
\end{secnotation}

In this section, we describe the idea of the blueprint technique. At some points, we will restrict ourselves to the root system $ \roots = A_3 $ and the standard root base $ \rootbase $ to simplify notation, but it will always be clear how our arguments should be generalised to arbitrary root systems. Recall from \cref{ADE:A-weylgrp} that the Weyl group of $ A_3 $ can be identified with the symmetric group of $ \numint{1}{4} $ in a canonical way (as soon as a suitable root basis has been fixed).

\begin{notation}\label{blue:index-notation}
	In the context of the blueprint technique, we will often write $ i $ for the simple root $ \delta_i $. For example, we write $ \rootgr{i} $ for the root group $ \rootgr{\delta_i} $ or $ 121 $ for the word $ (\delta_1, \delta_2, \delta_1) $.
\end{notation}

At first, we have to introduce homotopy cycles.

\begin{definition}[Homotopy cycles]
	Let $ w $ be an arbitrary element of $ W $ and let $ \word{\alpha} $ be any expression of $ w $ (in the sense of \cref{rootsys:length-def}). A \defemph*{homotopy cycle of $ \word{\alpha} $}\index{homotopy cycle} is a sequence $ \tup{\word{\alpha}}{k} $ of expressions of $ w $ such that $ \word{\alpha}_1 = \word{\alpha}_k $ and such that $ \word{\alpha}_i $ is braid-homotopic to $ \word{\alpha}_{i-1} $ (in the sense of \cref{weyl:homotopy-def}) for all $ i \in \numint{1}{k-1} $. A homotopy cycle is said to be \defemph*{trivial}\index{homotopy cycle!trivial} if every homotopy step is \enquote{reversed} by a later step, and \defemph*{non-trivial} otherwise.
\end{definition}

\begin{note}
	Our definition of triviality for homotopy cycles is rather vague. However, it will only be used in informal discussions, so this poses no problem.
\end{note}

\begin{example}[Homotopy cycles in $ A_3 $]\label{blue:homotopy-ex-A3}
	Consider the root system $ \roots = A_3 $. Here we have the braid moves
	\[ 121 \rightarrow 212, \qquad 232 \rightarrow 323, \qquad 13 \rightarrow 31 \]
	(and their inverses) which correspond to the braid relations
	\begin{gather*}
		\reflbr{\delta_1} \reflbr{\delta_2} \reflbr{\delta_1} = \reflbr{\delta_2} \reflbr{\delta_1} \reflbr{\delta_2}, \qquad \reflbr{\delta_2} \reflbr{\delta_3} \reflbr{\delta_2} = \reflbr{\delta_3} \reflbr{\delta_2} \reflbr{\delta_3}, \\
		\reflbr{\delta_1} \reflbr{\delta_3} = \reflbr{\delta_3} \reflbr{\delta_1}.
	\end{gather*}
	Thus the sequences $ 121 \rightarrow 212 \rightarrow 121 $ and $ 12113 \rightarrow 21213 \rightarrow 21231 \rightarrow 12131 \rightarrow 12113 $ are examples of homotopy cycles. Clearly, they are trivial. A reduced expression of the longest element $ \rho $ in $ W $ (in the sense of \cref{rootsys:longest-el-char}) is given by $ \word{\alpha}_1 = 123121 $. A non-trivial homotopy cycle of $ \word{\alpha}_1 $ is given in \cref{fig:blue:A3-cycle}.
\end{example}

\premidfigure
\begin{figure}[htb]
	\centering\begin{tabular}{rc@{\hspace{2cm}}rc}
		(1) & 12\underline{31}21 & (8) & 321\underline{323} \\
		(2) & \underline{121}321 & (9) & 3\underline{212}32 \\
		(3) & 21\underline{232}1 & (10) & \underline{31}2\underline{13}2 \\
		(4) & 2\underline{13}2\underline{31} & (11) & 1\underline{323}12 \\
		(5) & 23\underline{121}3 & (12) & 123\underline{212} \\
		(6) & \underline{232}123 & (13) & 123121 \\
		(7) & 32\underline{31}23 & &
	\end{tabular}
	\caption{A non-trivial homotopy cycle for the longest word in $ \Weyl(A_3) $ (taken from \cite[Fig.~8, \onpage{472}]{MoufangPolygons}).}
	\label{fig:blue:A3-cycle}
\end{figure}
\postmidfigure

\begin{figure}
	\begin{subfigure}{\linewidth}
				\centering\begin{tikzpicture}[auto, font=\scriptsize, densely dashed, scale=0.939]
			\foreach \angle [count=\ani] in {0, 90, 180, 270}{
				\begin{scope}[rotate=\angle]
					\node[nodecircle](\ani-1) at (1, 1){};
					\node[nodecircle](\ani-2) at (2, 2){};
					\node[nodecircle](\ani-3) at (2, 4){};
					\node[nodecircle](\ani-4) at (4, 4){};
					\node[nodecircle](\ani-5) at (4, 2){};
					\node[nodecircle](\ani-6) at (5.5, 5.5){};
					\draw (\ani-1) edge node{2} (\ani-2)
						(\ani-4) edge node{2} (\ani-6);
				\end{scope}
			}

			\foreach \angle/\ani in {0/1, 180/3}{
				\draw (\ani-2) edge node{3} (\ani-3)
					(\ani-3) edge node{1} (\ani-4)
					(\ani-4) edge node{3} (\ani-5)
					(\ani-5) edge node{1} (\ani-2);
			}
			
			\foreach \angle/\ani in {90/2, 270/4}{
				\draw (\ani-2) edge node{1} (\ani-3)
					(\ani-3) edge node{3} (\ani-4)
					(\ani-4) edge node{1} (\ani-5)
					(\ani-5) edge node{3} (\ani-2);
			}
			
			\draw (1-6) edge node{1} (2-6)
				(2-6) edge node{3} (3-6)
				(3-6) edge node{1} (4-6)
				(4-6) edge node{3} (1-6);
			\draw (1-1) edge node{3} (2-1)
					(2-1) edge node{1} (3-1)
					(3-1) edge node{3} (4-1)
					(4-1) edge node{1} (1-1);
			\draw (1-3) edge node{2} (2-5)
					(2-3) edge node{2} (3-5)
					(3-3) edge node{2} (4-5)
					(4-3) edge node{2} (1-5);
					
			\node[below left, font=\normalfont] at (2-6){$ \id $};
			\node[below left, font=\normalfont] at (4-1){$ \rho $};
					
			\draw[thick, black, solid] (2-6) -- (1-6) -- (1-4) -- (1-5) -- (1-2) -- (1-1) -- (4-1);
		\end{tikzpicture}
		\caption{The path $ 123121 $ (rows 1 and 13 in \cref{fig:blue:A3-cycle} on page~\pageref{fig:blue:A3-cycle}).}
		\label{fig:blue:An:cycle:1}
	\end{subfigure}
	
	\addvspace{0.4cm}
	\begin{subfigure}{\linewidth}
				\centering\begin{tikzpicture}[auto, font=\scriptsize, densely dashed, scale=0.939]
			\foreach \angle [count=\ani] in {0, 90, 180, 270}{
				\begin{scope}[rotate=\angle]
					\node[nodecircle](\ani-1) at (1, 1){};
					\node[nodecircle](\ani-2) at (2, 2){};
					\node[nodecircle](\ani-3) at (2, 4){};
					\node[nodecircle](\ani-4) at (4, 4){};
					\node[nodecircle](\ani-5) at (4, 2){};
					\node[nodecircle](\ani-6) at (5.5, 5.5){};
					\draw (\ani-1) edge node{2} (\ani-2)
						(\ani-4) edge node{2} (\ani-6);
				\end{scope}
			}

			\foreach \angle/\ani in {0/1, 180/3}{
				\draw (\ani-2) edge node{3} (\ani-3)
					(\ani-3) edge node{1} (\ani-4)
					(\ani-4) edge node{3} (\ani-5)
					(\ani-5) edge node{1} (\ani-2);
			}
			
			\foreach \angle/\ani in {90/2, 270/4}{
				\draw (\ani-2) edge node{1} (\ani-3)
					(\ani-3) edge node{3} (\ani-4)
					(\ani-4) edge node{1} (\ani-5)
					(\ani-5) edge node{3} (\ani-2);
			}
			
			\draw (1-6) edge node{1} (2-6)
				(2-6) edge node{3} (3-6)
				(3-6) edge node{1} (4-6)
				(4-6) edge node{3} (1-6);
			\draw (1-1) edge node{3} (2-1)
					(2-1) edge node{1} (3-1)
					(3-1) edge node{3} (4-1)
					(4-1) edge node{1} (1-1);
			\draw (1-3) edge node{2} (2-5)
					(2-3) edge node{2} (3-5)
					(3-3) edge node{2} (4-5)
					(4-3) edge node{2} (1-5);
					
			\node[below left, font=\normalfont] at (2-6){$ \id $};
			\node[below left, font=\normalfont] at (4-1){$ \rho $};
					
					323123
			\draw[thick, black, solid] (2-6) -- (3-6) -- (3-4) -- (3-5) -- (3-2) -- (3-1) -- (4-1);
		\end{tikzpicture}
		\caption{The path $ 323123 $ (row 7 in \cref{fig:blue:A3-cycle} on page~\pageref{fig:blue:A3-cycle}).}
		\label{fig:blue:An:cycle:2}
	\end{subfigure}
	\caption{The Cayley graph of $ \Sym(\numint{1}{4}) $ for the standard generators $ s_1  $, $ s_2 $, $ s_3 $ and two reduced paths from $ \id $ to the longest element $ \rho $. See \cref{blue:homo-cycle-cayley}.}
	\label{fig:cayley-A3}
\end{figure}
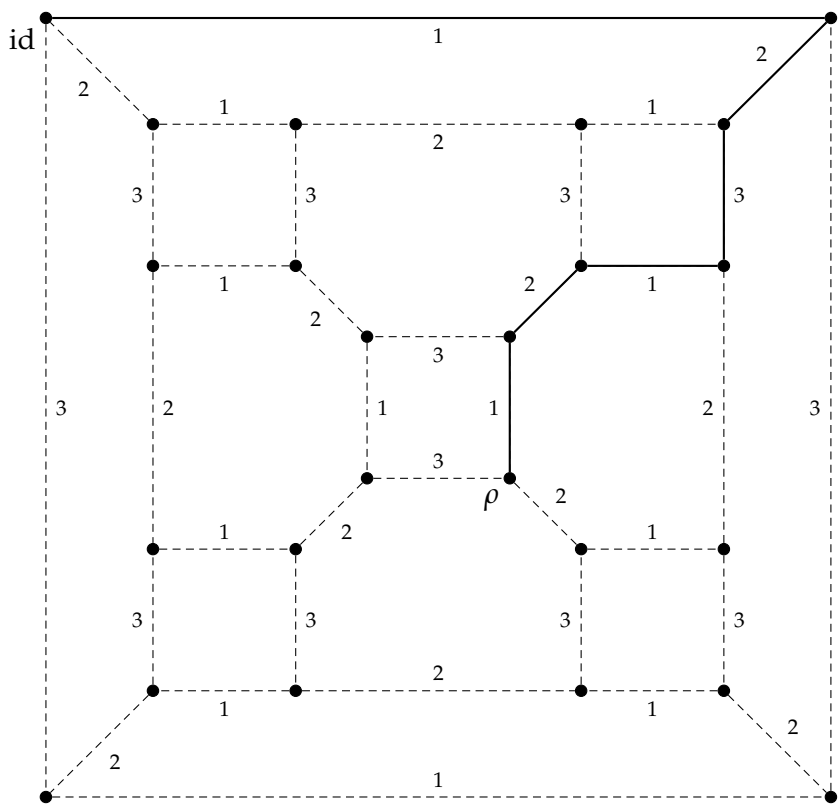
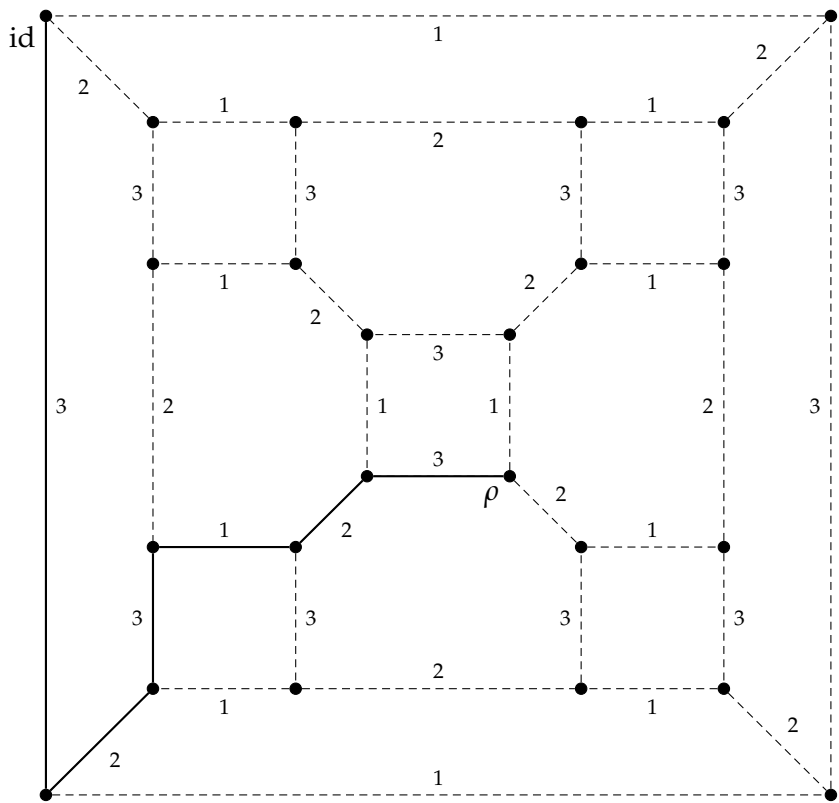

\begin{remark}[Homotopy cycles in Cayley graphs]\label{blue:homo-cycle-cayley}
	The \defemph*{Cayley graph of $ W $ with respect to the generators $ \reflbr{\delta_1}, \ldots, \reflbr{\delta_n} $} is, by definition, the graph whose vertex set is $ W $ and whose edge set is
	\[ \Set{\compactSet{w, w \reflbr{\delta_i}} \given w \in W, i \in \numint{1}{n}}. \]
	Further, any edge of the form $ \compactSet{w, w \reflbr{\delta_i}} $ is labelled with $ i $. As an example, we consider the Weyl group of $ A_3 $. Recall from \cref{ADE:A-weylgrp} that this group can be identified with $ \Sym(\numint{1}{4}) $ such that $ \refl{\basvec_1 - \basvec_2} = \transp{1}{2} $, $ \refl{\basvec_2 - \basvec_3} = \transp{2}{3} $, $ \refl{\basvec_3 - \basvec_4} = \transp{3}{4} $. The Cayley graph of $ \Sym(\numint{1}{4}) $ with respect to these standard generators is depicted in \cref{fig:cayley-A3} on page~\pageref{fig:cayley-A3}. Now the (reduced) expressions of any element $ w $ correspond to the paths from $ \id $ to $ w $ (of minimal length) and the braid homotopy moves correspond to certain transformation rules for these paths.
	
	The homotopy cycle of the longest element $ \rho $ in \cref{fig:blue:A3-cycle} can be interpreted as follows. The Cayley graph in \cref{fig:cayley-A3} can also be drawn on the sphere $ S^2 $ with the large outer square (whose top-left corner is $ \id $) on the top and the small inner square (whose bottom-right corner is $ \rho $) on the bottom. Then $ \id $ and $ \rho $ lie on opposite points of $ S^2 $. Now any minimal path from $ \id $ to $ w $ goes in an approximate semicircle from the top to the bottom, and the homotopy moves in \cref{fig:blue:A3-cycle} \enquote{push this semicircle once around (the hole in) $ S^2 $}. The first and the seventh path of this sequence are depicted in \cref{fig:blue:An:cycle:1,fig:blue:An:cycle:2}, respectively.
	
	As a corollary of the observations in the previous paragraph, we observe that $ \rho $ is the only element of $ W $ for which non-trivial homotopy cycles exist: reduced expression of any other word are \enquote{too short to be pushed around the sphere}. This statement remains true if we replace $ \roots $ by any other root system, and it is the reason why we will only consider reduced expressions of the longest element in the blueprint technique.
\end{remark}

\begin{note}
	We have observed in \cref{blue:homo-cycle-cayley} that every element of $ W $ which is not the longest element has only trivial homotopy cycles. If $ \roots $ has rank 2, then even the longest element has only trivial homotopy cycles: If $ \Set{\delta, \delta'} $ is the root base of such a root system, then the only reduced expressions of the longest element $ \rho $ are $ \braidword{m}(\delta, \delta') $ and $ \braidword{m}(\delta', \delta) $ for $ m \defl \abs{\roots}/2 $. Hence the only homotopy cycles of $ \rho $ are $ \braidword{m}(\delta, \delta') \rightarrow \braidword{m}(\delta', \delta) \rightarrow \braidword{m}(\delta, \delta') $ and $ \braidword{m}(\delta', \delta) \rightarrow \braidword{m}(\delta, \delta') \rightarrow \braidword{m}(\delta, \delta') $ (and iterations of these cycles). Clearly, these cycles are trivial.
\end{note}

Using the language of homotopy cycles, we can now describe the idea of the blueprint technique.

\begin{miscthm}[The blueprint technique]\label{blue-idea}
	For ease of presentation, we consider only the root system $ \roots = A_3 $ in this introduction, but the underlying principle applies to all root systems. We use the standard representation of $ A_3 $ from \cref{ADE:simply-laced:An-standard-rep} and we denote by $ \rootbase = \Set{\delta_1, \delta_2, \delta_3} $ the standard root base of $ A_3 $, so that $ \delta_i = \basvec_i - \basvec_{i+1} $ for all $ i \in \numint{1}{3} $. By \cref{ADE:thm}, we have a ring $ (\ring, +, \rmult) $ and isomorphisms $ (\map{\risom{\alpha}}{(\ring, +)}{\rootgr{\alpha}}{}{})_{\alpha \in \roots} $. This ring must be associative by \cref{ADE:thm}, but we proceed as if we did not know this, and we try to derive a new proof of this fact.

	For each $ i \in \numint{1}{13} $, we denote by $ \word{\alpha}_i $ the reduced expression of the longest element $ \rho $ in the $ i $-th column of \cref{fig:blue:A3-cycle}. For example, $ \word{\alpha}_1 = (\delta_1, \delta_2, \delta_3, \delta_1, \delta_2, \delta_1) = \word{\alpha}_{13} $ or, in simplified notation, $ \word{\alpha}_1 = 123121 = \word{\alpha}_{13} $.	Further, we choose arbitrary ring elements $ a,b,c,d,e,f \in \ring $ and consider the tuple
	\[ y_1 \defl \brackets[\big]{\risom{1}(a), \risom{2}(b), \risom{3}(c), \risom{1}(d), \risom{2}(e), \risom{1}(f)}. \]
	In other words, $ y_1 $ denotes an arbitrary element of
	\[ \rootgr{\word{\alpha}_1} = \rootgr{1} \times \rootgr{2} \times \rootgr{3} \times \rootgr{1} \times \rootgr{2} \times \rootgr{1}. \]
	
	Recall that \cref{fig:blue:A3-cycle} depicts a non-trivial homotopy cycle of $ \word{\alpha}_1 $ which uses the braid moves from \cref{blue:homotopy-ex-A3}. Using the theoretical framework of the blueprint technique, we can compute a \enquote{sensible} rewriting rule on the tuple $ y_1 $ for each of these homotopy moves. For examples, the braid moves $ 13 \rightarrow 31 $ and $ 121 \rightarrow 212 $ correspond to rewriting rules
	\begin{align*}
		\map{}{\rootgr{1} \times \rootgr{3}&}{\rootgr{3} \times \rootgr{1}}{\\\brackets[\big]{\risom{1}(x), \risom{3}(y)}&}{\brackets[\big]{\risom{3}(y), \risom{1}(x)}} \rightand \\
		\map{}{\rootgr{1} \times \rootgr{2} \times \rootgr{1}&}{\rootgr{2} \times \rootgr{1} \times \rootgr{2}}{\\\brackets[\big]{\risom{1}(x), \risom{2}(y), \risom{1}(z)}&}{\brackets[\big]{\risom{2}(z), \risom{1}(-y-zx), \risom{2}(x)}}.
	\end{align*}
	Now the main idea of the blueprint technique is to \enquote{work down} the column in \cref{fig:blue:A3-cycle} and to apply the corresponding rewriting rule to the tuple $ y_1 $ in each step. We refer to this procedure as the \defemph{blueprint computation}. In the first step, we apply the elementary homotopy $ 31 \rightarrow 13 $ to $ \word{\alpha}_1 $, which produces $ \word{\alpha}_2 = 121321 $. Simultaneously, we apply the corresponding rewriting rule $ \map{}{\rootgr{3} \times \rootgr{1}}{\rootgr{1} \times \rootgr{3}}{}{} $ to $ y_1 $. This yields a tuple
	\[ y_2 = \brackets[\big]{\risom{1}(a), \risom{2}(b), \risom{1}(c'), \risom{3}(d'), \risom{2}(e), \risom{1}(f)} \in \rootgr{\word{\alpha}_2} \]
	where $ c',d' \in \ring $ are new ring elements which depend on $ c $ and $ d $. In the next step, we apply the elementary homotopy $ 121 \rightarrow 212 $ to $ \word{\alpha}_2 $ to obtain the word $ \word{\alpha}_3 = 212321 $. Again we use the associated rewriting rule to compute a new tuple
	\[ y_3 = \brackets[\big]{\risom{2}(a'), \risom{1}(b'), \risom{2}(c''), \risom{3}(d'), \risom{2}(e), \risom{1}(f)} \in \rootgr{\word{\alpha}_2} \]
	where $ a',b',c'' $ depend on $ a $, $ b $ and $ c' $. After applying twelve transformations in total, we end up with $ \word{\alpha}_{13} = 123121 = \word{\alpha}_1 $ and a tuple
	\[ y_{13} = \brackets[\big]{\risom{1}(\tilde{a}), \risom{2}(\tilde{b}), \risom{3}(\tilde{c}), \risom{1}(\tilde{d}), \risom{2}(\tilde{e}), \risom{1}(\tilde{f})} \in \rootgr{\word{\alpha}_{13}} = \rootgr{\word{\alpha}_1}. \]
	For reasons which will be explained later, we know that we must have $ y_1 = y_{13} $, so $ \tilde{a} = a $, $ \tilde{b} = b, \ldots, \tilde{f} = f $. Since $ \tilde{a}, \ldots, \tilde{f} $ can be calculated explicitly (using our rewriting rules), we obtain six relations which hold for arbitrary elements $ a, \ldots, f \in \ring $. When we actually perform these calculations, five of these relations turn out to be trivial, but we will see that the last relation is equivalent to the associativity law. Therefore, the blueprint technique yields that any ring which coordinatises an $ A_3 $-graded group must be associative.
\end{miscthm}

\begin{remark}[Formal justification]\label{blue-form-just}
	In order to actually perform the calculations in \cref{blue-idea}, we need to state the rewriting rules on our tuples. These rules should be chosen in a way which guarantees that, at the end of the computation, we have $ y_{13} = y_1 $. The idea is as follows: For every reduced expression $ \word{\beta} = \tup{\beta}{m} $ of the longest element, we define a certain map $ \map{\blumap{\word{\beta}}}{\rootgr{\beta_1} \times \cdots \times \rootgr{\beta_m}}{G}{}{} $ and we choose our rewriting rules so that $ \blumap{\word{\alpha}_{i+1}}(y_{i+1}) = \blumap{\word{\alpha}_i}(y_i) $ for all $ i $. In other words, we can think of $ \blumapsym $ as an invariant of the involved tuples, called the \defemph{blueprint invariant}, and the rewriting rules are chosen precisely to leave the invariant unchanged. Thus we have
	\begin{align*}
		\blumap{\word{\alpha}_1}(y_1) = \blumap{\word{\alpha}_2}(y_2) = \cdots = \blumap{\word{\alpha}_{k-1}}(y_{k-1}) = \blumap{\word{\alpha}_{k}}(y_{k}) = \blumap{\word{\alpha}_1}(y_k)
	\end{align*}
	where $ k $ is the length of the homotopy cycle in question. Now it remains to show that $ \blumap{\word{\alpha}_1} $ is injective, which is done in \cref{blue:inj}.
\end{remark}

\begin{remark}[Simplification of the blueprint computation]\label{blue-half}
	Instead of working down an entire homotopy cycle, we will in practice perform two separate computations: One working \enquote{halfway down}, which transforms the initial tuple $ y_1 $ into some tuple $ x $, and another one working \enquote{halfway up}, which transforms $ y_1 $ into some other tuple $ x' $. By the same arguments as in \cref{blue-form-just}, we have $ \blumap{}(x') = \blumap{}(y_1) = \blumap{}(x'') $, so we infer that $ x' = x'' $. We will do this because the tuple entries in the blueprint computation become more complicated with each step, so it is more efficient to perform two short computations than to perform one long computation.
\end{remark}


\section{The Formal Framework of the Blueprint Technique}

\label{sec:rewrite}
\begin{secnotation}
	We denote by $ \roots $ a root system of rank $ n $, by $ \rootbase =(\delta_1, \ldots, \delta_n) $ a rescaled ordered root base of $ \roots $, by $ \possys $ the corresponding positve system in $ \roots $ and by $ W \defl \Weyl(\roots) $ the Weyl group. Further, we denote by $ G $ a $ \roots $-graded group with root groups $ (\rootgr{\alpha})_{\alpha \in \roots} $ and we fix a $ \rootbase $-system $ (w_\delta)_{\delta \in \rootbase} $ of Weyl elements in $ G $.
\end{secnotation}

The first part of this section is devoted to the introduction of the invariant~$ \blumapsym $ that we announced in \cref{blue-form-just} and the proof that it is an injective map (\cref{blue:inj}). In the second part, we state some basic results concerning rewriting rules which will be used later on.

\begin{reminder}
	Recall from \cref{rgg:pregrade-def} that for any word $ \word{\alpha} = \tup{\alpha}{k} $ over $ \roots $, we put
	$ \rootgr{\word{\alpha}} \defl \rootgr{\alpha_1} \times \cdots \times \rootgr{\alpha_k} $.
	This should not be confused with the notation $ \rootgr{S} \defl \gen{\rootgr{\alpha} \given \alpha \in S} $ for any subset $ S $ of $ \roots $.
\end{reminder}

\begin{definition}[Blueprint invariant]\label{blumap-def}
	Let $ \word{\delta} $ be a word over $ \rootbase $. The \defemph*{blueprint invariant of type $ \word{\delta} $ (with respect to $ (w_\delta)_{\delta \in \rootbase} $)}\index{blueprint invariant} is the map
	\[ \map{\blumapG{\word{\delta}}}{\rootgr{\word{\delta}}}{G}{\tup{g}{m}}{\prod_{i=1}^m w_{\delta_i} g_i = w_{\delta_1} g_1 \cdots w_{\delta_m} g_m}. \]
	Further, the \defemph*{reduced blueprint invariant of type $ \word{\delta} $ (with respect to $ (w_\delta)_{\delta \in \rootbase} $)}\index{blueprint invariant!reduced} is the composition of $ \blumapG{\word{\delta}} $ with the canonical projection to $ G/\zentrum(G) $:
	\[ \map{\blumap{\word{\delta}}}{\rootgr{\word{\alpha}}}{G/\zentrum(G)}{\tup{g}{m}}{\blumapG{\listing{g}{m}} \zentrum(G)}. \]
\end{definition}

\begin{definition}[Blueprint rewriting rules]
	Let $ \word{\alpha} $, $ \word{\beta} $ be words over $ \rootbase $. A \defemph*{rewriting rule of type $ (\word{\alpha}, \word{\beta}) $ (with respect to $ (w_\delta)_{\delta \in \rootbase} $)}\index{rewriting rule} is a map $ \map{}{\rootgr{\word{\alpha}}}{\rootgr{\word{\beta}}}{}{} $. A rewriting rule $ \phi $ of type $ (\word{\alpha}, \word{\beta}) $ is called a \defemph*{blueprint rewriting rule}\index{rewriting rule!blueprint|see{blueprint rewriting rule}}\index{blueprint rewriting rule} if it respects the reduced blueprint invariant: $ \blumap{\word{\beta}} \circ \phi = \blumap{\word{\alpha}} $.
\end{definition}

\premidfigure
\begin{figure}[htb]
	\centering$ \begin{tikzcd}
		\rootgr{\word{\alpha}} \arrow[rr, "\phi"] \arrow[dr, "\blumap{{\word{\alpha}}}"'] && \rootgr{\word{\beta}} \arrow[dl, "\blumap{{\word{\beta}}}"] \\
		& G/\zentrum(G) &
	\end{tikzcd} $
	\caption{The condition of being a blueprint rewriting rule.}
	\label{fig:blu-comp}
\end{figure}
\postmidfigure

The main result in this section is the fact that the reduced blueprint invariant is an injective map. For the non-reduced blueprint invariant, this is a direct consequence of the injectivity of the product maps (\cref{rgg:prodmap-Nw}), but we have to work a bit more to see that the same holds modulo the center. 

\begin{lemma}\label{blue:Nw-center}
	Let $ \word{\alpha} = \tup{\alpha}{m} $ be a reduced word over $ \rootbase $ (where $ m \ge 1 $) and let $ v \defl \reflbr{\word{\alpha}} $. Then $ \rootgr{\switchset(v)} \intersect \zentrum(G) = \compactSet{1_G} $ where $ \zentrum(G) $ denotes the center of $ G $.
\end{lemma}
\begin{proof}
	It suffices to prove the statement for the root base corresponding to $ \rootbase $, so we can assume that $ \rootbase $ is a proper root base. This allows us to apply \cref{rootsys:switchset-lem,rootsys:switchset-description}.
	The statement is trivial if $ m=1 $, so we assume that $ m>1 $. Let $ \word{\beta} = \tup{\beta}{m} $ be the root sequence associated to $ \word{\alpha} $ as in \cref{rootsys:rootseq-def}, so that $ \switchset(w) = \Set{\beta_1, \cdots, \beta_m} $ by \cref{rootsys:switchset-description}.
	Let $ x \in \rootgr{\switchset(w)} \intersect \zentrum(G) $. We want to show that $ x = 1_G $. By \cref{rgg:prodmap-Nw}, the product map on $ \word{\beta} $ is bijective, so there exist $ x_1 \in \rootgr{\beta_1}, \ldots, x_m \in \rootgr{\beta_m} $ such that $ x = x_1 \cdots x_m $. Choose an arbitrary $ \beta_m $-Weyl element $ w_m $, and observe that $ \beta_m = \alpha_m $ by the definition of root sequences. In particular, $ \beta_m $ is a simple root. Now $ x_1 \cdots x_m = x_1^{w_m} \cdots x_m^{w_m} $ because $ x_1 \cdots x_m $ is central, which implies that
	\[ (x_1^{w_m} \cdots x_{m-1}^{w_m})^{-1} x_1 \cdots x_m = x_m^{w_m}. \]
	Note that the element on the left-hand side lies in $ \rootgr{\possys} $ (because $ (\possys \setminus \IR_{>0} \gamma)^{\reflbr{\gamma}} = \possys \setminus \IR_{>0}\gamma \subs \possys $ for any $ \gamma \in \rootbase $ by \cref{rootsys:simple-refl-on-pos}) while the element on the right-hand side lies in $ \rootgr{-\beta_m} $. Since $ \rootgr{\possys} \intersect \rootgr{-\beta_m} = \compactSet{1_G} $ by Axiom~\thmitemref{rgg-def}{rgg-def:nondeg}, it follows that $ x_m^{w_m} = 1_G $, so $ x_m = 1_G $.
	
	We conclude from the previous paragraph that $ x = x_1 \cdots x_{m-1} $. Put $ v' \defl v \reflbr{\alpha_m} $. Note that $ (\alpha_1, \ldots, \alpha_{m-1}) $ is a reduced expression of $ v' $. By \cref{rootsys:switchset-lem},
	\[ \Set{\beta_1, \cdots, \beta_{m-1}} = \switchset(v) \setminus \compactSet{\beta_m} = \switchset(v) \setminus \compactSet{\alpha_m} = \switchset(w')^{\reflbr{\alpha_m}}. \]
	Thus $ x_1 \cdots x_{m-1} $ lies in $ \rootgr{\switchset(v')^{\reflbr{\alpha_m}}} = \rootgr{\switchset(v')}^{w_m} $, so $ x_1^{w_m} \cdots x_{m-1}^{w_m} $ lies in $ \rootgr{\switchset(v')} $. It follows that $ x = x_1^{w_m} \cdots x_{m-1}^{w_m} \in \rootgr{\switchset(w')} \intersect \zentrum(G) $. Using induction, we infer that $ x=1_G $, as desired.
\end{proof}

\begin{proposition}\label{blue:pos-center}
	$ \rootgr{\possys} \intersect \zentrum(G) = \compactSet{1_G} $.
\end{proposition}
\begin{proof}
	The longest word $ \rho $ satisfies $ \switchset(\rho) = \indivset{\possys} $ by \cref{rootsys:longest-el-char}. Since $ \rootgr{\indivset{\possys}} = \rootgr{\possys} $, the assertion follows from \cref{blue:Nw-center}.
\end{proof}

\begin{proposition}\label{blue:inj}
	Let $ \word{\alpha} = \tup{\alpha}{m} $ be a reduced word over $ \rootbase $. Then the reduced blueprint invariant $ \blumap{\word{\alpha}} $ of type $ \word{\alpha} $ with respect to $ (w_\delta)_{\delta \in \rootbase} $ is injective.
\end{proposition}
\begin{proof}
	As in \cref{blue:index-notation}, we put $ \rootgr{i} \defl \rootgr{\delta_i} $ and $ w_i \defl w_{\delta_i} $ for all $ i \in \numint{1}{n} $.
	Let $ g =\tup{g}{m} $ and $ h=\tup{h}{m} $ be elements of $ U_1 \times \cdots \times U_m $ such that $ \blumap{\word{\alpha}}(g) = \blumap{\word{\alpha}}(h) $. Then there exists $ z \in \zentrum(G) $ such that
	\[ w_1 g_1 w_2 g_2 \cdots w_{m-1} g_{m-1} w_m g_m = w_1 h_1 w_2 h_2 \cdots w_{m-1} h_{m-1} w_m h_mz. \]
	We want to show that $ g=h $. Using that $ ab = b a^b $ for all $ a,b \in G $, we see that the product on the left side equals
	\[ w_1 w_2 \cdots w_m g_1^{w_2 \cdots w_m} g_2^{w_3 \cdots w_m} \cdots g_{m-1}^{w_m} g_m. \]
	A similar transformation can be done on the right-hand side. We infer that
	\begin{equation}\label{eq:blue:inj}
		g_1^{w_2 \cdots w_m} g_2^{w_3 \cdots w_m} \cdots g_{m-1}^{w_m} g_m = h_1^{w_2 \cdots w_m} h_2^{w_3 \cdots w_m} \cdots h_{m-1}^{w_m} h_m z.
	\end{equation}
	Put
	\[ \word{g} \defl \brackets{g_1^{w_2 \cdots w_m}, g_2^{w_3 \cdots w_m}, \ldots, g_{m-1}^{w_m}, g_m}, \qquad \word{h} \defl \brackets{h_1^{w_2 \cdots w_m}, h_2^{w_3 \cdots w_m}, \ldots, h_{m-1}^{w_m}, h_m}. \]
	Note that for each $ i \in \numint{1}{m} $, $ g_i^{w_{i+1} \cdots w_m} $ and $ h_i^{w_{i+1} \cdots w_m} $ lie in the root group associated to $ \alpha_i^{\reflbr{\alpha_{i+1} \cdots \alpha_{m}}} $. It follows that $ \word{g} $ and $ \word{h} $ lie in $ \rootgr{\word{\beta}} $ where $ \word{\beta} = \tup{\beta}{m} $ denotes the root sequence associated to $ \word{\alpha} $. Thus~\eqref{eq:blue:inj} says that $ \mu_{\word{\beta}}(\word{g}) = \mu_{\word{\beta}}(\word{h}) z $ where $ \mu_{\word{\beta}} $ denotes the product map on $ \word{\beta} $. Hence
	\[ z = \mu_{\word{\beta}}(\word{h})^{-1} \mu_{\word{\beta}}(\word{g}) \in \rootgr{\possys} \intersect \zentrum(G). \]
	By \cref{blue:pos-center}, it follows that $ z=1_G $, so we actually have $ \mu_{\word{\beta}}(\word{g}) = \mu_{\word{\beta}}(\word{h}) $. Using \cref{rgg:prodmap-Nw}, we conclude that $ \word{g} = \word{h} $. Since conjugation by Weyl elements is an automorphism, this yields $ g = h $, so $ \blumap{\word{\alpha}} $ is indeed injective.
\end{proof}

We can now state the essence of the blueprint technique as follows.

\begin{theorem}\label{blue:thm}\index{blueprint technique}
	Let $ \word{\alpha}, \word{\beta} $ be two reduced words over $ \rootbase $ such that $ \refl{\word{\alpha}} = \refl{\word{\beta}} $, let $ x \in \rootgr{\word{\alpha}} $ and let $ \map{\psi_1, \psi_2}{\rootgr{\word{\alpha}}}{\rootgr{\word{\beta}}}{}{} $ be two blueprint rewriting rules. Then $ \psi_1(x) = \psi_2(x) $.
\end{theorem}
\begin{proof}
	Since $ \psi_1 $ and $ \psi_2 $ are $ \blumap{} $-compatible, we have
	\[ \blumap{\word{\beta}}\brackets[\big]{\psi_1(x)} = \blumap{\word{\alpha}}(x) = \blumap{\word{\beta}}\brackets[\big]{\psi_2(x)}. \]
	Since $ \blumap{\word{\beta}} $ is injective by \cref{blue:inj}, the assertion follows.
\end{proof}

\begin{note}
	In practice, all blueprint rewriting rules in this book will be compatible with the non-reduced blueprint invariant, and not merely with the reduced blueprint invariant. 
\end{note}

The rest of this section is a collection of elementary facts about blueprint rewriting rules.

\begin{remark}\label{blumap-conc}
	Let $ \word{\alpha} = \tup{\alpha}{n} $ and $ \word{\beta} = \tup{\beta}{m} $ be two words over $ \rootbase $, let $ g_i \in \rootgr{\alpha_i} $ for all $ i \in \numint{1}{n} $ and let $ h_j \in \rootgr{\beta_j} $ for all $ j \in \numint{1}{m} $. Then
	\begin{align*}
		\blumap{\word{\alpha}\word{\beta}}\brackets{\listing{g}{n}, \listing{h}{m}} &= \blumap{\word{\alpha}}(\listing{g}{n}) \blumap{\word{\beta}}(\listing{h}{m}).
	\end{align*}
\end{remark}

\begin{remark}\label{blumap-comp-inv}
	Let $ \word{\alpha} $, $ \word{\beta} $ be words over $ \roots $ and let $ \map{\phi}{\rootgr{\word{\alpha}}}{\rootgr{\word{\beta}}}{}{} $ be a bijective blueprint rewriting rule. Then $ \blumap{\word{\beta}} \circ \phi = \blumap{\word{\alpha}} $, so
	\[ \blumap{\word{\beta}} = \blumap{\word{\beta}} \circ \phi \circ \phi^{-1} = \blumap{\word{\alpha}} \circ \phi^{-1}. \]
	That is, $ \phi^{-1} $ is also a blueprint rewriting rule.
\end{remark}

\begin{remark}
	Clearly, the concatenation of two blueprint rewriting rules is again a blueprint rewriting rule.
\end{remark}

In the following \cref{blumap-glue} and its proof, we will identify $ (\word{g}, \word{h}) $ (for any tuples $ \word{g} = \tup{g}{n} $, $ \word{h} = \tup{h}{m} $) with the tuple $ (\listing{g}{n}, \listing{h}{m}) $.

\begin{lemma}\label{blumap-glue}
	Let $ \word{\alpha} $, $ \word{\beta} $, $ \word{\alpha}' $, $ \word{\beta}' $ be words over $ \roots $ and let
	\[ \map{\phi}{\rootgr{\word{\alpha}}}{\rootgr{\word{\alpha}'}}{\word{g}}{\phi(\word{g})} \midand \map{\psi}{\rootgr{\word{\beta}}}{\rootgr{\word{\beta}'}}{\word{h}}{\psi(\word{h})} \]
	be blueprint rewriting rules. Further, define
	\begin{align*}
		\map{\phi \times \psi}{\rootgr{\word{\alpha} \word{\beta}} = \rootgr{\word{\alpha}} \times \rootgr{\word{\beta}}}{\rootgr{\word{\alpha}' \word{\beta}'} = \rootgr{\word{\alpha}'} \times \rootgr{\word{\beta}'}}{(\word{g}, \word{h})}{\brackets[\big]{\phi(\word{g}), \psi(\word{h})}}.
	\end{align*}
	Then $ \phi \times \psi $ is a blueprint rewriting rule.
\end{lemma}
\begin{proof}
	Let $ \word{g} \in \rootgr{\word{\alpha}} $ and let $ \word{h} \in \rootgr{\word{\beta}} $. By \cref{blumap-conc}, we have
	\begin{equation*}
		\blumap{\word{\alpha} \word{\beta}}(\word{g}, \word{h}) = \blumap{\word{\alpha}}(\word{g}) \blumap{\word{\beta}}(\word{h}).
	\end{equation*}
	On the other hand,
	\begin{align*}
		\brackets[\big]{\blumap{\word{\alpha}' \word{\beta}'} \circ (\phi \times \psi)}(\word{g}, \word{h}) &= \blumap{\word{\alpha}' \word{\beta}'}\brackets[\big]{\phi(\word{g}), \psi(\word{h})} = \blumap{\word{\alpha}'}\brackets[\big]{\phi(\word{g})} \blumap{\word{\beta}'}\brackets[\big]{\psi(\word{h})}.
	\end{align*}
	Since $ \phi $ and $ \psi $ are blueprint rewriting rules, the assertion follows.
\end{proof}

The following result says that it is sufficient to compute blueprint rewriting rules for the basic braid moves, for example for $ 121 \rightarrow 212 $ in $ A_3 $. These basic rewriting rules can then be extended to rewriting rules which are defined on reduced expressions of the longest word.

\begin{proposition}\label{blue-kr-id}
	Let $ \word{\alpha}_1, \word{\alpha}_2, \word{\alpha}_3, \word{\alpha}_2' $ be words over $ \Delta $ and let $ \map{\phi}{\rootgr{\word{\alpha}_2}}{\rootgr{\word{\alpha}_2'}}{}{} $ be a blueprint rewriting rule. Then the map
	\[ \map{\id \times \phi \times \id}{\rootgr{\word{\alpha}_1} \times \rootgr{\word{\alpha}_2} \times \rootgr{\word{\alpha}_3}}{\rootgr{\word{\alpha}_1} \times \rootgr{\word{\alpha}_2'} \times \rootgr{\word{\alpha}_3}}{(\word{g}_1, \word{g}_2, \word{g}_3)}{(\word{g}_1, \phi(\word{g}_2), \word{g}_3)} \]
	is a blueprint rewriting rule.
\end{proposition}
\begin{proof}
	The identity map is clearly a blueprint rewriting rule, so this is a consequence of \cref{blumap-glue}.
\end{proof}

The following blueprint rewriting rule is easy to compute and will be used for all root systems.

\begin{lemma}\label{blue:rewriting-switch}
	Let $ \alpha, \beta \in \roots $ be such that $ \beta $ is adjacent to $ \alpha $ and to $ -\alpha $. Then the map
	\[ \map{\phi}{\rootgr{\alpha} \times \rootgr{\beta}}{\rootgr{\beta} \times \rootgr{\alpha}}{(x_\alpha, x_\beta)}{(x_\beta, x_\alpha)}. \]
	is a blueprint rewriting rule.
\end{lemma}
\begin{proof}
	Note that $ -\beta $ is adjacent to $ -\alpha $ because $ \beta $ is adjacent to $ \alpha $. Similarly, $ -\beta $ is adjacent to $ \alpha $ because $ \beta $ is adjacent to $ -\alpha $. Thus $ \gen{\rootgr{\alpha}, \rootgr{-\alpha}} $ commutes with $ \gen{\rootgr{\beta}, \rootgr{-\beta}} $. Hence
	\begin{align*}
		\blumap{\alpha \beta}(x_\alpha, x_\beta) &= \whom{\alpha} x_\alpha \whom{\beta} x_\beta = \whom{\beta} x_\beta \whom{\alpha} x_\alpha \rightand \\
		(\blumap{\beta \alpha} \circ \phi)(x_\alpha, x_\beta) &= \blumap{\beta \alpha}(x_\beta, x_\alpha) = \whom{\beta} x_\beta \whom{\alpha} x_\alpha
	\end{align*}
	for all $ x_\alpha \in \rootgr{\alpha} $ and all $ x_\beta \in \rootgr{\beta} $. In other words, $ \phi $ is a blueprint rewriting rule.
\end{proof}


\newcommand{\basroot}[1]{\delta_{#1}}

\section{Blueprint Computations for \texorpdfstring{$ A_3 $}{A\_3}}

\label{sec:A3-blue}
\begin{secnotation}
	We consider the root system $ A_n $ for $ n \ge 3 $ in standard representation and with its standard ordered root base $ \rootbase = \tup{\delta}{n} $. Further, we denote by $ G $ a group with an $ A_n $-grading $ (\rootgr{\alpha})_{\alpha \in A_n} $ and by $ (\risom{\alpha})_{\alpha \in A_n} $ a coordinatisation of $ G $ by some ring $ \ring $ with signs and base point chosen as in \cref{ADE:A-param-standard}.
	and by $ (w_\delta)_{\delta \in \rootbase} $ a fixed $ \rootbase $-system of Weyl elements in $ G $.

\end{secnotation}

\begin{notation}\label{blue:A3:root-notation}
	We will usually (and in contrast to \cref{blue:index-notation}) identify any root $ \basvec_i - \basvec_j $ with the pair $ (i,j) $ and we will always write $ ij $ instead of $ (i,j) $. In other words, we set $ \rootgr{ij} \defl \rootgrmin{i}{j} $ and $ \risom{ij} \defl \rismin{i}{j} $ for all distinct $ i,j \in \numint{1}{n+1} $. Further, we denote by
	\[ w_{ij} \defl w_{ij}(1_\ring) = \risom{ji}(-1_\ring) \risom{ij}(1_\ring) \risom{ij}(1_\ring) \]
	the Weyl element from \cref{ADE:stsign:weyl-char}.
\end{notation}

In this section, we will demonstrate how the blueprint technique can be used to show that the coordinatising ring $ \ring $ must be associative. In other words, we give a new proof of the $ A_n $-case in \cref{ADE:assoc}. All computations in this section will take place in the parabolic root subsystem
\[ A_3 \defl \Set{\basvec_i - \basvec_j \given i \ne j \in \numint{1}{4}} \]
of $ A_n $ and in the corresponding $ A_3 $-graded subgroup of $ G $. No further information can be obtained from the blueprint technique by considering larger subsystems.

At first, we have to compute the blueprint rewriting rules for $ A_3 $. We begin by explicitly stating them. In the proof of \cref{blue:A3-rewriting-comp}, we will see how to arrive at these specific formulas.

\begin{definition}[Rewriting rules for $ A_3 $]\label{blue-A3-mapdef}
	We define the following rewriting rules:
	\begin{align*}
		\map{\blutrans{123}}{\rootgr{12} \times \rootgr{23} \times \rootgr{12}&}{\rootgr{23} \times \rootgr{12} \times \rootgr{12}}{\\\brackets[\big]{\risom{12}(a), \risom{23}(b), \risom{12}(c)}&}{\brackets[\big]{\risom{23}(c), \risom{12}(-b-ca), \risom{23}(a)}}; \\
		\map{\blutrans{234}}{\rootgr{23} \times \rootgr{34} \times \rootgr{23}&}{\rootgr{34} \times \rootgr{23} \times \rootgr{34}}{\\\brackets[\big]{\risom{23}(a), \risom{34}(b), \risom{23}(c)}&}{\brackets[\big]{\risom{34}(c), \risom{23}(-b-ca), \risom{34}(a)}}; \rightand \\
		\map{\phi}{\rootgr{12} \times \rootgr{34}&}{\rootgr{34} \times \rootgr{12}}{\\\brackets[\big]{\risom{12}(a), \risom{34}(b)}&}{\brackets[\big]{\risom{34}(b), \risom{12}(a)}}.
	\end{align*}
	It is easy to see that these maps are bijections whose inverses are given by
	\begin{align*}
		\map{\blutrans{123}^{-1}}{\rootgr{23} \times \rootgr{12} \times \rootgr{12}&}{\rootgr{12} \times \rootgr{23} \times \rootgr{12}}{\\\brackets[\big]{\risom{23}(a), \risom{12}(b), \risom{23}(c)}&}{\brackets[\big]{\risom{12}(c), \risom{23}(-b-ac), \risom{12}(a)}}; \\
		\map{\blutrans{234}^{-1}}{\rootgr{34} \times \rootgr{23} \times \rootgr{34}&}{\rootgr{23} \times \rootgr{34} \times \rootgr{23}}{\\\brackets[\big]{\risom{34}(a), \risom{23}(b), \risom{34}(c)}&}{\brackets[\big]{\risom{23}(c), \risom{34}(-b-ac), \risom{23}(a)}}; \rightand \\
		\map{\phi^{-1}}{\rootgr{34} \times \rootgr{12}&}{\rootgr{12} \times \rootgr{34}}{\\\brackets[\big]{\risom{34}(b), \risom{12}(a)}&}{\brackets[\big]{\risom{12}(a), \risom{34}(b)}}.
	\end{align*}
\end{definition}

\begin{lemma}\label{blue:A3-rewriting-comp}
	$ \blutrans{123} $, $ \blutrans{123}^{-1} $, $ \blutrans{234} $, $ \blutrans{234}^{-1} $, $ \phi $ and $ \phi^{-1} $ are blueprint rewriting rules (with respect to the Weyl elements $ \whom{12} $, $ \whom{23} $, $ \whom{34} $).
\end{lemma}
\begin{proof}
	By \cref{blumap-comp-inv}, we only need to prove that $ \blutrans{123} $, $ \blutrans{234} $ and $ \phi $ are blueprint rewriting rules. Further, $ \phi $ is a blueprint rewriting rule by \cref{blue:rewriting-switch}. We proceed to prove that $ \blutrans{123} $ is a blueprint rewriting rule (and the proof for $ \blutrans{234} $ is identical). Let $ a,b,c \in \ring $ and set
	\[ \word{\alpha} \defl (\basvec_1 - \basvec_2, \basvec_2 - \basvec_3, \basvec_1 - \basvec_2) \midand \word{\beta} \defl (\basvec_2 - \basvec_3, \basvec_1 - \basvec_2, \basvec_2 - \basvec_3). \]
	On the one hand, we have
	\begin{align*}
		x \defl{}& (\blumap{\word{\beta}} \circ \blutrans{123})\brackets[\big]{\risom{12}(a), \risom{23}(b), \risom{12}(c)} = \blumap{\word{\beta}}\brackets[\big]{\risom{23}(c), \risom{12}(-b-ca), \risom{23}(a)} \\
		={}& \whom{23} \risom{23}(c) \whom{12} \risom{12}(-b-ca) \whom{23} \risom{23}(a) \\
		={}& \whom{23} \whom{12} \whom{23} \risom{23}(c)^{\whom{12} \whom{23}} \risom{12}(-b-ca)^{\whom{23}} \risom{23}(a).
	\end{align*}
	Using the conjugation formula from \cref{ADE:stsign:A-standard-conjformula}, we can compute that
	\begin{align*}
		\risom{23}(c)^{\whom{12} \whom{23}} \risom{12}(-b-ca)^{\whom{23}} \risom{23}(a) &= \risom{13}(-c)^{\whom{23}} \risom{13}(-b-ca) \risom{23}(a) \\
		&= \risom{12}(c) \risom{13}(-b-ca) \risom{23}(a)
	\end{align*}
	On the other hand, we have
	\begin{align*}
		y \defl{}& \blumap{\word{\alpha}}\brackets[\big]{\risom{12}(a), \risom{23}(b), \risom{12}(c)} = \whom{12} \risom{12}(a) \whom{23} \risom{23}(b) \whom{12} \risom{12}(c) \\
		={}& \whom{12} \whom{23} \whom{12} \risom{12}(a)^{\whom{23} \whom{12}} \risom{23}(b)^{\whom{12}} \risom{12}(c)
	\end{align*}
	where
	\begin{align*}
		\risom{12}(a)^{\whom{23} \whom{12}} \risom{23}(b)^{\whom{12}} \risom{12}(c) = \risom{13}(a)^{\whom{12}} \risom{13}(-b) \risom{12}(c) = \risom{23}(a) \risom{13}(-b) \risom{12}(c).
	\end{align*}
	We have to prove that $ x=y $. Since $ \whom{23} \whom{12} \whom{23} = \whom{12} \whom{23} \whom{12} $ by \cref{braid:all}, it only remains to show that
	\[ \risom{12}(c) \risom{13}(-b-ca) \risom{23}(a) = \risom{23}(a) \risom{13}(-b) \risom{12}(c). \]
	In order to prove this, we apply the commutator relations in $ A_3 $-graded groups:
	\begin{align*}
		\risom{23}(a) \risom{13}(-b) \risom{12}(c) &= \risom{23}(a) \risom{12}(c) \risom{13}(-b) = \risom{12}(c) \risom{23}(a) \commutator{\risom{23}(a)}{\risom{12}(c)} \risom{13}(-b) \\
		&= \risom{12}(c) \risom{23}(a) \commutator{\risom{12}(c)}{\risom{23}(a)}^{-1} \risom{13}(-b) \\
		&= \risom{12}(c) \risom{23}(a) \risom{13}(-ca) \risom{13}(-b) = \risom{12}(c) \risom{13}(-b-ca) \risom{23}(a).
	\end{align*}
	Thus $ x=y $, so $ \blutrans{123} $ is a blueprint rewriting rule. This finishes the proof.
\end{proof}

Now that we have determined the blueprint rewriting rules for $ A_3 $, we can apply the ideas from~\ref{blue-idea} and~\ref{blue-half} in practice.

\begin{theorem}\label{A3-assoc}
	The multiplication $ \rmult $ on $ \ring $ is associative.
\end{theorem}
\begin{proof}
	In this proof, we drop \cref{blue:A3:root-notation} and go back to \cref{blue:index-notation} instead. That is, we write $ \risom{i} $ in place of $ \risom{i,i+1} $. Let $ a,b,c,d,e,f \in \ring $ be arbitrary and set
	\[ x \defl \brackets[\big]{\risom{1}(a), \risom{2}(b), \risom{3}(c), \risom{1}(d), \risom{2}(e), \risom{1}(f)}. \]
	A homotopy cycle for the longest word in $ A_3 $ is given in \cref{fig:blue:An:cycle} on page~\pageref{fig:blue:An:cycle}.
	Working down rows~1 to~7 and applying the corresponding blueprint rewriting rules in the process (as in \cref{fig:blue:An:down}), we obtain a tuple $ y = \tup{y}{6} $ where
	\begin{align*}
		y_1 &= \risom{3}(f), & y_2 &= \risom{2}(-e-fd), & y_3 &= \risom{3}(d), \\
		y_4 &= \risom{1}\brackets[\big]{c+ea+f(b+da)}, & y_5 &= \risom{2}(-b-da), & y_6 &= \risom{3}(a).
	\end{align*}
	Further, working up rows 13 to 7 yields a tuple $ z = \tup{z}{6} $ where
	\begin{align*}
		z_1 &= \risom{3}(f), & z_2 &= \risom{2}(-e-fd), & z_3 &= \risom{3}(d), \\
		z_4 &= \risom{1}\brackets[\big]{c+fb+(e+fd)a}, & z_5 &= \risom{2}(-b-da), & z_6 &= \risom{3}(a).
	\end{align*}
	(For the intermediate steps of these two calculations, see \cref{fig:blue:An:down,fig:blue:An:up}.) Since each used rewriting rule is a blueprint rewriting rule, an application of \cref{blue:thm} yields $ y_i = z_i $ for all $ i \in \numint{1}{6} $. In particular, $ y_4 = z_4 $, so
	\[ c+ea+f(b+da) = c+fb+(e+fd)a \qquad \text{for all } a,b,c,d,e,f \in \ring. \]
	Setting $ b = c = e = 0 $, we conclude that $ f(da) = (fd)a $ for all $ a,f,d \in \ring $. That is, $ \ring $ is associative.
\end{proof}

\premidfigure
\begin{sidewaysfigure}
	\begin{minipage}{0.84\linewidth}
		\begin{subfigure}{\linewidth}
			\hspace{0.2cm}\begin{tabular}{rc|cccccc}
				(1) & 12\underline{31}21 & $ \risom{1}(a) $ & $ \risom{2}(b) $ & $ \risom{3}(c) $ & $ \risom{1}(d) $ & $ \risom{2}(e) $ & $ \risom{1}(f) $ \\
				(2) & \underline{121}321 & $ \risom{1}(a) $ & $ \risom{2}(b) $ & $ \risom{1}(d) $ & $ \risom{3}(c) $ & $ \risom{2}(e) $ & $ \risom{1}(f) $ \\
				(3) & 21\underline{232}1 & $ \risom{2}(d) $ & $ \risom{1}(-b-da) $ & $ \risom{2}(a) $ & $ \risom{3}(c) $ & $ \risom{2}(e) $ & $ \risom{1}(f) $ \\
				(4) & 2\underline{13}2\underline{31} &$ \risom{2}(d) $ & $ \risom{1}(-b-da) $ & $ \risom{3}(e) $ & $ \risom{2}(-c-ea) $ & $ \risom{3}(a) $ & $ \risom{1}(f) $ \\
				(5) & 23\underline{121}3 & $ \risom{2}(d) $ & $ \risom{3}(e) $ & $ \risom{1}(-b-da) $ & $ \risom{2}(-c-ea) $ & $ \risom{1}(f) $ & $ \risom{3}(a) $ \\
				(6) & \underline{232}123 & $ \risom{2}(d) $ & $ \risom{3}(e) $ & $ \risom{2}(f) $ & $ \risom{1}(c+ea+f(b+da)) $ & $ \risom{2}(-b-da) $ & $ \risom{3}(a) $ \\
				(7) & 323123 & $ \risom{3}(f) $ & $ \risom{2}(-e-fd) $ & $ \risom{3}(d) $ & $ \risom{1}\brackets[\big]{c+ea+f(b+da)} $ & $ \risom{2}(-b-da) $ & $ \risom{3}(a) $
			\end{tabular}
			\caption{The first part of the blueprint computation: rows 1 to 7.}
			\label{fig:blue:An:down}
		\end{subfigure}
		
		\addvspace{0.75cm}
		\begin{subfigure}{\linewidth}
			\begin{tabular}{rc|cccccc}
				(13) & 123\underline{121} & $ \risom{1}(a) $ & $ \risom{2}(b) $ & $ \risom{3}(c) $ & $ \risom{1}(d) $ & $ \risom{2}(e) $ & $ \risom{1}(f) $ \\
				(12) & 1\underline{232}12 & $ \risom{1}(a) $ & $ \risom{2}(b) $ & $ \risom{3}(c) $ & $ \risom{2}(f) $ & $ \risom{1}(-e-fd) $ & $ \risom{2}(d) $ \\
				(11) & \underline{13}2\underline{31}2 & $ \risom{1}(a) $ & $ \risom{3}(f) $ & $ \risom{2}(-c-fb) $ & $ \risom{3}(b) $ & $ \risom{1}(-e-fd) $ & $ \risom{2}(d) $ \\
				(10) & 3\underline{121}32 & $ \risom{3}(f) $ & $ \risom{1}(a) $ & $ \risom{2}(-c-fb) $ & $ \risom{1}(-e-fd) $ & $ \risom{3}(b) $ & $ \risom{2}(d) $ \\
				(9) & 321\underline{232} & $ \risom{3}(f) $ & $ \risom{2}(-e-fd) $ & $ \risom{1}\brackets[\big]{c+fb+(e+fd)a} $ & $ \risom{2}(a) $ & $ \risom{3}(b) $ & $ \risom{2}(d) $ \\
				(8) & 32\underline{13}23 & $ \risom{3}(f) $ & $ \risom{2}(-e-fd) $ & $ \risom{1}\brackets[\big]{c+fb+(e+fd)a} $ & $ \risom{3}(d) $ & $ \risom{2}(-b-da) $ & $ \risom{3}(a) $ \\
				(7) & 323123 & $ \risom{3}(f) $ & $ \risom{2}(-e-fd) $ & $ \risom{3}(d) $ & $ \risom{1}\brackets[\big]{c+fb+(e+fd)a} $ & $ \risom{2}(-b-da) $ & $ \risom{3}(a) $
			\end{tabular}
			\caption{The first part of the blueprint computation: rows 13 to 7.}
			\label{fig:blue:An:up}
		\end{subfigure}
	\end{minipage}%
	\begin{subfigure}{0.16\linewidth}
		\centering\begin{tabular}{rc}
			(1) & 12\underline{31}21 \\
			(2) & \underline{121}321 \\
			(3) & 21\underline{232}1 \\
			(4) & 2\underline{13}2\underline{31} \\
			(5) & 23\underline{121}3 \\
			(6) & \underline{232}123 \\
			(7) & 32\underline{31}23 \\
			(8) & 321\underline{323} \\
			(9) & 3\underline{212}32 \\
			(10) & \underline{31}2\underline{13}2 \\
			(11) & 1\underline{323}12 \\
			(12) & 123\underline{212} \\
			(13) & 123121
		\end{tabular}
		\caption{A homotopy cycle.}
		\label{fig:blue:An:cycle}
	\end{subfigure}
	\caption{The blueprint technique applied to $ A_3 $.}
	\label{fig:blue:An}
\end{sidewaysfigure}
\postmidfigure

\begin{note}
	We could set $ b,c $ and $ e $ to zero right at the beginning of the computation in \cref{fig:blue:An} to obtain an easier (and yet flawless) proof of \cref{A3-assoc}. Of course, this is only possible with the knowledge which variables turn out to be immaterial, and we can only obtain this knowledge by performing the full computation.
\end{note}

\begin{note}\label{blue:A3:distr}
	Observe that we have never applied the distributive laws in the computation in \cref{fig:blue:An}. In fact, these laws follow from the blueprint computations: Putting $ c=b=0 $ and $ f=1 $ in the equation $ y_4 = z_4 $ in the proof of \cref{A3-assoc}, we obtain that
	\[ ea + da = (e+d)a. \]
	Conversely, putting $ c=d=0 $ and $ d=1 $, we infer that
	\[ f(b+a) = fb+fa. \]
	This proof is independent of the one in \cref{ADE:distributivity}. We will get back to this observation in \cref{blue:identity-classes}.
\end{note}

\begin{note}
	In the same manner as for $ A_3 $, we can use the blueprint technique to show that any ring which coordinatises a $ D_4 $-graded group is commutative. Since all the root systems $ E_6 $, $ E_7 $, $ E_8 $ and $ D_n $ for $ n \ge 4 $ contain $ D_4 $ as a parabolic subsystem, it follows that the same assertion holds for these root systems as well.
\end{note}

\FloatBarrier


\section{Concluding Remarks}

We end this chapter with some remarks.

\begin{remark}[Commutation maps]\label{blue:commmap-rem}
	In \cref{sec:A3-blue}, we have shown that any ring $ (\ring, +, \rmult) $ which coordinatises an $ A_3 $-graded group (with standard signs) is associative. We can also phrase this as follows: We have shown in \cref{ADE:ring-identity,A3-assoc,blue:A3:distr} that for any group $ (\ring, +) $ which coordinatises an $ A_3 $-graded group with standard signs, the commutation map $ \map{f}{\ring \times \ring}{\ring}{}{} $ satisfies the identities which turn $ (\ring, +) $ into an associative ring $ (\ring, +, f) $. This viewpoint -- that the blueprint computations yield identities for the commutation maps which are precisely the axioms of some algebraic structure -- will generalise without problems to the more complicated root systems.
\end{remark}

\begin{remark}[Classes of identities]\label{blue:identity-classes}
	It will be beneficial to informally distinguish two classes of identities for commutation maps. We say that an identity is a \defemph*{rank-2 identity}\index{rank-2 identity} if it follows from computations \enquote{within a single rank-2 subgroup}, and we say that it is a \defemph*{higher-rank identity}\index{higher-rank identity} if it requires computations in a subgroup of rank 3 or higher. For example, the \enquote{bi-additivity} of commutators in \cref{basic:comm-add} is proven in the rank-2 setting, and so the distributive laws for the ring multiplication in \cref{ADE:distributivity} are rank-2 identities. We will see other examples of rank-2 identities in \cref{B:comm-add,B:comm-add-cry,BC:comm-add}. On the other hand, the associativity and commutativity laws in \cref{ADE:assoc,ADE:comm} are higher-rank identities.
	
	Both kinds of identities can in principle be derived with the blueprint technique, as we have seen in \cref{A3-assoc,blue:A3:distr}. However, since the rank-2 identities can be proven in an easy and straightforward way, it is usually more efficient to derive them directly and then to use them during the blueprint computation to simplify terms.
	
	Another important type of identities are \defemph*{Weyl identities}\index{Weyl identity}. These are, by definition, identities which hold for the specific elements that parametrise the fixed $ \rootbase $-system $ (w_\delta)_{\delta \in \rootbase} $ of Weyl elements. In the $ A_3 $-case, the identity $ 1_\ring \rmult r = r = r \rmult 1_\ring $ (where $ 1_\ring $ is defined as in \cref{ADE:1-const}) is an (and the only) example of a Weyl identity. These identities always follow from formulas as in \cref{A2Weyl:basecomp-cor}. Technically, they are rank-2 identities, but it is beneficial to consider them separately because they do not follow from the blueprint computations. The reason for this is that the blueprint computations only produce identities which hold for arbitrary parameters whereas the Weyl identities involve the specific elements which parametrise $ (w_\delta)_{\delta \in \rootbase} $.
\end{remark}

\begin{strategy}\label{blue:summary}
	Let $ \roots $ be a root system of rank at least~3 and let $ G $ be a $ \roots $-graded group. We can summarise the steps in the blueprint technique for $ \roots $-graded groups as follows:
	\begin{remenumerate}
		\item Compute a reduced expression of the longest element in $ \Weyl(\roots) $ and a non-trivial homotopy cycle of this expression.
		
		\item Define commutation maps in $ G $, as in \cref{ADE:ring-mult-const} (and later in \cref{B:commmap-def,BC:commmap-def}).
		
		\item Compute the rank-2 identities of the commutation maps, including the Weyl identities.
		
		\item For each braid homotopy move in $ \Weyl(\roots) $, compute a blueprint rewriting rule. These rewriting rules involve the commutation maps in $ G $.
		
		\item Perform the blueprint computations, as in \cref{fig:blue:An}.
		
		\item Show that the set of identities which result from the blueprint computations is equivalent to a set of \enquote{smaller, easier-to-understand} identities. This is done by putting several of the involved variables to zero (or to $ 1 $ or to similar \enquote{canonical} elements) and then using the resulting identities to simplify the original identities.
	\end{remenumerate}
\end{strategy}

\begin{note}[The original conception of the blueprint technique]\label{blue:original}
	Using the Hall-Witt identity (\cref{hall-witt}), Zhang shows in \cite[3.4.9 to 3.4.18]{Zhang} that the commutation maps in certain $ C_3 $-graded groups can be described explicitly in terms of the algebraic structure of a ring with involution. The original purpose of the blueprint technique was to show that this ring must be alternative. Using Zhang's computations, the blueprint rewriting rules for these $ C_3 $-graded groups can be expressed in terms of only the multiplication and the involution on the ring, which makes them less complicated. As a consequence, the blueprint computation is less cumbersome, essentially because many of the identities which would \enquote{normally} result from the blueprint computation have already been verified \enquote{by hand} with the Hall-Witt identity. The result of the blueprint computation in this setting is that the ring in question is indeed alternative.
	
	It is a crucial observation that the blueprint computations work just as well if the rewriting rules are given only in terms of abstract commutation maps (see also \cref{blue:commmap-rem}). In other words, the structure of a parametrisation of $ (G, (\rootgr{\alpha})_{\alpha \in \roots}) $ (in the sense of \cref{param:param-def}) is sufficient to apply the blueprint computations. This makes it possible to circumvent any Hall-Witt computations, which significantly streamlines our approach to root graded groups.
\end{note}

\begin{remark}[Rank higher than~3]\label{blue:rank-3}
	The length of the longest word in $ \Weyl(\roots) $ and thus of the blueprint computation increases with the rank of $ \roots $, so there is a high incentive to keep this rank as low as possible. It turns out that, in all types of root systems for which we perform the blueprint computation, there is no need to go higher than rank~3. For example, if we were to perform the blueprint computation for $ A_4 $, we would obtain no additional identities. For this reason, we will always restrict the blueprint computations to the rank-3 case.
	
	We emphasise that this restriction does not lead to a restriction of the generality of our coordinatising results. To see why, observe that every root system of type $ X_n $ for $ X \in \Set{A, B, C, BC} $ and $ n \ge 3 $ contains a subsystem of type $ X_3 $. Thus every $ X_n $-graded group $ G $ contains an $ X_3 $-graded subgroup $ G' $. The blueprint computations in rank~3 show that the commutation maps in $ G' $ satisfy certain identities. Using similar ideas as in the proof of \cref{ADE:comm-formula}, we can describe all commutation maps in $ G $ using only the commutation maps in $ G' $. Thus the blueprint computation yields information about the computation maps in all of $ G $, not only in $ G' $.
\end{remark}

	\chapter{Root Gradings of Type \texorpdfstring{$ B $}{B}}
	
	In this chapter, we investigate $ B_n $-graded groups for $ n \ge 2 $, though we will restrict ourselves to the case $ n \ge 3 $ for the main results. Except for the case of RGD-systems, we are not aware of any prior literature on this topic. The standard reference for RGD-systems of type $ B_2 $ is Chapter~23 in \cite{MoufangPolygons}. Many of the arguments in \cite{MoufangPolygons} no longer work in the setting of root graded groups. Instead, we have to deploy several \enquote{rank-3 arguments}.
	
	Unlike the simply-laced root systems, $ B_n $ (and every other root systems that we study from this point on) has multiple orbits of roots, so we will for the first time see root graded groups with non-isomorphic pairs of root groups. Since every long root in $ B_n $ is contained in an $ A_2 $-subsystem (\cref{B:rootsys:long-in-A2}), we should expect the long root groups to be coordinatised by some ring $ \comring $. In the end, we will show that this ring must be commutative associative and that the short root groups are coordinatised by a quadratic module $ \module $ over $ \comring $.
	
	More precisely, we prove the following two statements. Firstly, there exist an abelian group $ (\comring, +) $ which parametrises the long root groups and another abelian group $ (\module, +) $ which parametrises the short root groups. Secondly, the commutation maps will be described by a map $ \map{}{\comring \times \comring}{\comring}{}{} $ which turns $ \comring $ into a commutative associative ring, a map $ \map{}{\comring \times \module}{\module}{}{} $ which turns $ \module $ into a left $ \comring $-module, a $ \comring $-quadratic form $ \map{q}{\module}{\comring}{}{} $ and a $ \comring $-bilinear form $ \map{f}{\module \times \module}{\comring}{}{} $ which is the linearisation of $ q $. The commutator formula has the exact same form as the Chevalley commutator formula, except that expressions of the form $ v^2 $ and $ 2vw $ for $ v,w \in \module $ (which do not make sense for module elements) are replaced by $ q(v) $ and $ f(v,w) $, respectively.
	
	This chapter is organised in the way described in \cref{sec:param:outline}. We begin with a brief introduction to quadratic modules over commutative associative rings. In \cref{sec:B:rootsys}, we study some purely combinatorial properties of the root system $ B_n $ and compute its Cartan integers. In \cref{sec:B-example}, we show how a so-called elementary orthogonal group can be constructed from any quadratic module, which provides a complete solution of the existence problem. The rank-2 computations are split into two parts: We begin with the general case of (non-crystallographic) $ B_2 $-gradings in \cref{sec:B:rank2-noncry} and continue with crystallographic $ B_2 $-gradings in \cref{sec:B:rank2-cry}. Results from both sections will be used in \cref{chap:BC} as well.
	
	From the second half of the chapter on, we consider $ B_n $-gradings for $ n \ge 3 $. In \cref{sec:B-rank3}, we prove that Weyl elements in these groups satisfy the square formula. In \cref{sec:B:stsigns,sec:B:sttwist}, we define the notions of standard partial twisting systems and standard signs for $ B_n $-graded groups, respectively. In \cref{sec:B-param}, we apply the results from the previous sections and the parametrisation theorem to construct parametrising groups $ (\comring, +) $ and $ (\module, +) $. The commutation maps, their rank-2 identities and the blueprint rewriting rules are defined, derived and computed in \cref{sec:B-bluerules}. Finally, we perform the blueprint computations for $ B_3 $ in \cref{sec:B-bluecomp}, and we state our final result in \cref{B:thm}.
	
	\label{chap:B}
	

\section{Quadratic Maps and Modules}

\label{sec:quadmod}

\begin{secnotation}
	We denote by $ \comring $ an arbitrary commutative associative ring. Unless otherwise specified, all modules are understood to be left modules over $ \comring $ (in the standard sense of \cref{ring:module}).
\end{secnotation}

In this section, we introduce the basic terminology of quadratic maps, forms and modules. These objects are often only considered over base fields, but the basic notions (which are all that we need) translate to the case of commutative associative base rings without a change. With this caveat, all the material in this section is standard and can be found, for example, in \cite[Section~7]{ElmanKarpenkoMerkurjev-QuadForms}. A reference which treats quadratic modules over arbitrary commutative associative rings is \cite{Knus_QuadFormRing}, though it only considers quadratic modules which are finitely generated projective.

Some of the results in this section will be generalised in \cref{sec:weak-quad} to the more general weakly quadratic maps.

We begin with the notion of quadratic maps between arbitrary $ \comring $-modules.

\begin{definition}[Polarisation, {\cite[7.1]{ElmanKarpenkoMerkurjev-QuadForms}}]\label{quadmod:polarisation-def}
	Let $ (M,+), (N,+) $ be two abelian groups and let $ \map{f}{M}{N}{}{} $ be any map. The \defemph*{polarisation of $ f $}\index{polarisation} is the map
	\[ \map{}{M \times M}{N}{(v,w)}{f(v+w) - f(v) - f(w)}. \]
	More generally, let $ n $ be a positive integer, let $ (M_1, +), \ldots, (M_n, +), (N,+) $ be abelian groups and let $ \map{f}{M_1 \times \cdots \times M_n}{N}{}{} $ be any map. Then for all $ i \in \numint{1}{n} $, the \defemph*{polarisation of $ f $ at position $ i $}\index{polarisation!at position i@at position $ i $} is the map
	\begin{gather*}
		\map{}{M_1 \times \cdots \times M_{i-1} \times M_i \times M_i \times M_{i+1} \times \cdots \times M_n}{N}{\\ (v_1, \ldots, v_{i-1}, v_i, w_i, v_{i+1}, \ldots, v_n)}{\tilde{f}(v_i + w_i) - \tilde{f}(v_i) - \tilde{f}(w_i)}
	\end{gather*}
	where $ \tilde{f}(u_i) \defl f(v_1, \ldots, v_{i-1}, u_i, v_{i+1}, \ldots, v_n) $.
\end{definition}

\begin{definition}[Quadratic map, {\cite[7.1]{ElmanKarpenkoMerkurjev-QuadForms}}, {\cite[(5.3.5)]{Knus_QuadFormRing}}]\label{quadmod:quadmap-def}
	Let $ M,N $ be $ \comring $-modules. A map $ \map{q}{M}{N}{}{} $ is called \defemph*{$ \comring $-quadratic} (or simply \defemph*{quadratic}\index{quadratic map}) if it satisfies the following two conditions:
	\begin{stenumerate}
		\item \label{quadmod:quadmap-def:scalar}$ q(\lambda v) = \lambda^2 q(v) $ for all $ \lambda \in \comring $, $ v \in M $.
		
		\item The polarisation $ \map{f}{M \times M}{N}{(v,w)}{q(v+w) - q(v) - q(w)} $ of $ f $ is $ \comring $-bilinear. This map is also called the \defemph*{linearisation of $ q $}\index{linearisation}.
	\end{stenumerate}
\end{definition}

In practice, we will only be concerned with the case that $ N = \comring $, which we study from \cref{quadmod:quadform-def} on.

\begin{note}[Polarisations and linearisations]\label{quadmod:terminology-polar-linear}
	The terms \enquote{polarisation} and \enquote{linearisation} are often used interchangeably in the literature. In this book, we will use the word \enquote{polarisation} as a general descriptor for the maps defined in \cref{quadmod:polarisation-def} while only the polarisations of quadratic maps will be called \enquote{linearisations}. In other words, a polarisation will only be called a linearisation if it is actually linear (or rather, bilinear).
\end{note}

\begin{lemma}\label{quadmod:basiclem}
	Let $ M,N $ be $ \comring $-modules, let $ \map{q}{M}{N}{}{} $ be a quadratic map and let $ f $ denote its linearisation. Then the following statements hold:
	\begin{lemenumerate}
		\item \label{quadmod:basiclem:lin}$ q(v+w) = q(v) + q(w) + f(v,w) $ for all $ v,w \in M $.
		
		\item \label{quadmod:basiclem:sum-formula}For all $ n \in \Npos $ and all $ v_1, \ldots, v_n \in M $, we have
		\[ q\brackets*{\sum_{i=1}^n v_i} = \sum_{i=1}^n q(v_i) + \sum_{i=1}^n \sum_{j=i+1}^n f(v_i, v_j). \]
		
		\item $ f $ is symmetric, that is, $ f(v,w) = f(w,v) $ for all $ v,w \in M $.
		
		\item $ q(0_M) = 0_N $.
		
		\item \label{quadmod:basiclem:fvv}$ f(v,v) = 2_\comring q(v) $ for all $ v \in M $.
	\end{lemenumerate}
\end{lemma}
\begin{proof}
	The first assertion is simply a reformulation of the definition of $ f $ and the second assertion follows from the first one by induction. The third assertion is clear. The fourth assertion follows from the equation
	\[ q(0_M) = q(0_M + 0_M) = q(0_M) + q(0_M) + f(0_M, 0_M) = 2_\comring q(0_M). \]
	For the last assertion, let $ v \in M $. On the one hand,
	\[ q(2_\comring v) = 2_\comring^2 q(v) = 4_\comring q(v). \]
	On the other hand,
	\[ q(2_\comring v) = q(v+v) = q(v) + q(v) + f(v,v) = 2_\comring q(v) + f(v,v). \]
	Thus $ 4_\comring q(v) = 2_\comring q(v) + f(v,v) $, or in other words, $ 2_\comring q(v) = f(v,v) $.
\end{proof}

The following \cref{quadmod:symbil-bij} says that for commutative associative rings in which $ 2 $ is invertible, a quadratic map from $ M $ to $ N $ is essentially the same thing as a symmetric bilinear map from $ M \times M $ to $ N $. Thus the study of quadratic maps reduces to the study of symmetric bilinear maps in this case. In the general situation, however, we will see in \cref{quadmod:square-ex} that a quadratic map \enquote{contains more information} than its linearisation. For this reason, quadratic maps are better suited to the study of general commutative associative rings than symmetric bilinear forms.

\begin{lemma}\label{quadmod:symbil-bij}
	Assume that $ 2_\comring $ is invertible and let $ \map{f}{M \times M}{N}{}{} $ be a symmetric $ \comring $-bilinear map between $ \comring $-modules $ M,N $. Then
	\[ \map{q}{M}{N}{v}{\frac{f(v,v)}{2}} \]
	is a $ \comring $-quadratic map with linearisation $ f $, and it is the unique map with these properties. Thus we have a bijective correspondence between symmetric $ \comring $-bilinear maps $ \map{}{M \times M}{N}{}{} $ and $ \comring $-quadratic maps $ \map{}{M}{N}{}{} $ if $ 2_\comring $ is invertible.
\end{lemma}
\begin{proof}
	It is easy to check that $ \map{}{}{}{v}{\frac{f(v,v)}{2}} $ defines a quadratic map, and it follows from \thmitemcref{quadmod:basiclem}{quadmod:basiclem:fvv} that a quadratic map is uniquely determined by its linearisation.
\end{proof}

\begin{lemma}
	Let $ M,N $ be $ \comring $-modules. Then the set of quadratic maps from $ M $ to $ N $ is closed under addition and $ \comring $-scalar multiplication.
\end{lemma}
\begin{proof}
	This is clear from the definition of quadratic maps.
\end{proof}

From now on, we specialise to the case that the codomain of the quadratic map in question is $ \comring $. This is the only case of relevance for our purposes.

\begin{definition}[Quadratic forms and modules]\label{quadmod:quadform-def}
	Let $ M $ be a $ \comring $-module. A \defemph*{quadratic form on $ M $}\index{quadratic form} is a $ \comring $-quadratic map $ \map{q}{M}{\comring}{}{} $ where $ \comring $ is regarded as a $ \comring $-module in the natural way. It is called \defemph*{anisotropic}\index{quadratic form!anisotropic} if $ q(v) = 0 $ implies $ v=0 $ for all $ v \in \module $. The pair $ (M,q) $ is called a \defemph*{quadratic module (over $ \comring $)}\index{quadratic module}.
\end{definition}

\begin{example}\label{quadmod:square-ex}
	The map $ \map{q}{\comring}{\comring}{r}{r^2} $ is a quadratic form with linearisation $ \map{f}{\comring \times \comring}{\comring}{(r,s)}{2rs} $. In fact, the group that we will construct in \cref{sec:B-example} for an arbitrary quadratic module $ (\module, q) $ will be a Chevalley group if we take $ \module = \comring $ and $ \map{q}{}{}{r}{r^2} $. 
	
	More generally, the map
	\[ \map{q}{\comring^n}{\comring}{v}{\sum_{i=1}^n v_i^2} \]
	is a quadratic form with linearisation
	\[ \map{f}{\comring^n \times \comring^n}{\comring}{(v,w)}{2\sum_{i=1}^n v_i w_i}. \]
	Note that $ f=0 $ if $ 2_\comring = 0_\comring $. In particular, if $ 2_\comring = 0_\comring $, then $ q $ and the zero map are two quadratic forms which are distinct (unless $ \comring $ is the zero ring) but have the same linearisation. This shows that in general, a quadratic form \enquote{contains more information} than its linearisation.
\end{example}

\begin{lemma}[Direct sum of quadratic modules]\label{quadmod:dirsum}
	Let $ (\module, q) $, $ (\module', q') $ be two quadratic modules and let $ f,f' $ denote the linearisations of $ q $ and $ q' $, respectively. Then the direct sum $ \module \dirsum \module' $ together with the map
	\[ \map{q \dirsum q'}{\module \dirsum \module'}{\comring}{(v,v')}{q(v) + q'(v')} \]
	is a quadratic module, called the \defemph*{direct sum of $ (\module, q) $, $ (\module', q') $}\index{quadratic module!direct sum}. The linearisation of $ q \dirsum q' $ is
	\[ \map{f \dirsum f'}{(\module \dirsum \module') \times (\module \dirsum \module')}{\comring}{\brackets[\big]{(u,u'), (v, v')}}{f(u,v) + f(u',v')}. \]
\end{lemma}
\begin{proof}
	This follows from a simple computation.
\end{proof}

\begin{remark}[Inner direct sum]\label{quadmod:inner-dirsum}
	Let $ (\module, q) $ be a quadratic module with subspaces $ A,B $ and with linearisation $ f $. We say that $ \module $ is the \defemph*{(inner) direct sum of $ A $ and $ B $}\index{quadratic module!direct sum!inner}, and we write $ \module = A \dirsum B $, if $ \module $ is the (inner) direct sum of $ A $ and $ B $ as $ \comring $-modules and $ f(a,b) = 0 $ for all $ a \in A $, $ b \in B $. If this is the case, then $ \module $ is isomorphic to the outer direct sum of $ A $ and $ B $ in the sense of \cref{quadmod:dirsum}. A subspace $ C $ of $ \module $ is called a \defemph*{direct summand of $ \module $}\index{direct summand} if there exists a subspace $ D $ such that $ \module = C \dirsum D $ (as quadratic modules).
\end{remark}

\begin{example}[Hyperbolic space, {\cite[\onpage{40}]{ElmanKarpenkoMerkurjev-QuadForms}}, {\cite[(5.6.2)]{Knus_QuadFormRing}}]\label{quadmod:hyperbolic-ex}
	Let $ n \in \IN_{\ge 1} $ and let $ V \defl \comring^n $. Denote by $ V^* $ the dual space of $ V $. The \defemph*{hyperbolic space of dimension $ 2n $ over $ \comring $}\index{hyperbolic space} is the $ \comring $-module $ \module \defl V \dirsum V^* $ together with the quadratic form
	\[ \map{q}{\module}{\comring}{(v, f)}{f(v)}. \]
	If we denote the canonical basis of $ V $ by $ \tup{e}{n} $ and the corresponding dual basis by $ (e_{-1}, \ldots, e_{-n}) $, then $ q $ is given on coordinates by
	\[ \map{q}{\comring^n \dirsum \comring^n}{\comring}{\brackets[\big]{\tup{\lambda}{n}, (\lambda_{-1}, \ldots, \lambda_{-n})}}{\sum_{i=1}^n \lambda_i \lambda_{-i}}. \]
	The pairs $ (e_i, e_{-i}) $ are called \defemph*{hyperbolic pairs}\index{hyperbolic pair}. More abstractly, a pair $ (v,w) $ of elements in a quadratic module $ (\module', q') $ with linearisation $ f' $ is called a hyperbolic pair if $ q'(v) = 0 = q'(w) $, $ f'(v,w) = 1 $ and $ \gen{v,w}_\comring $ is a direct summand of $ \module' $ (in the sense of \cref{quadmod:inner-dirsum}). Such a pair is automatically $ \comring $-linearly independent: If $ \lambda, \mu \in \comring $ satisfy $ \lambda v + \mu w = 0 $, then
	\begin{align*}
		0_\comring &= f'(v,0) = f'\brackets[\big]{v, \lambda v + \mu w} = \lambda f'(v,v) + \mu f'(v,w) = 2\lambda q'(v) + \mu = \mu,
	\end{align*}
	and similarly $ \lambda = 0_\comring $. Thus $ v,w $ span a direct summand of $ \module $ which is isomorphic to the hyperbolic plane $ \comring^2 $.
\end{example}

\begin{note}
	The condition in \cref{quadmod:hyperbolic-ex} that the module spanned by a hyperbolic pair is a direct summand will not be relevant in our context. It is absent from the definition of hyperbolic pairs in \cite[11.17]{GPR_AlbertRing}.
\end{note}

\begin{lemma}\label{quadmod:2-faithful}
	Let $ (\module, q) $ be a quadratic $ \comring $-module which contains an element $ v_0 \in \module $ such that $ q(v_0) $ is invertible. Then for all $ a \in \comring $ with $ av_0 = 0_\module $, we have $ 2a=0 $. In particular, $ \module $ is a faithful $ \comring $-module if $ 2_\comring $ is not a zero divisor.
\end{lemma}
\begin{proof}
	Let $ a \in \comring $ such that $ av_0 = 0 $ and denote the linearisation of $ q $ by $ f $. Then $ a f(v_0, v_0) = f(av_0, v_0) = 0 $. At the same time, $ af(v_0, v_0) = 2aq(v_0) $ by \thmitemcref{quadmod:basiclem}{quadmod:basiclem:fvv}. Since $ q(v_0) $ is invertible, it follows that $ 2a=0 $, as desired. In particular, any $ a \in \comring $ with $ a \module = \compactSet{0} $ satisfies $ a v_0 = 0 $ and thus $ 2a=0 $. Hence $ \module $ is faithful if $ 2_\comring $ is not a zero divisor.
\end{proof}

The notions of orthogonality and orthogonal groups which are known for inner products transfer to the setting of quadratic forms without surprises.

\begin{definition}[Orthogonal group]\label{quadmod:ortho-grp}
	For any quadratic module $ (\module, q) $, we call
	\[ \Ortho(q) \defl \Set{\phi \in \Aut_\comring(\module) \given q(\phi(v)) = q(v) \text{ for all } v \in \module} \]
	the \defemph*{orthogonal group of $ q $}\index{orthogonal group}. We use the convention that it acts on $ \module $ from the left side, so that the composition $ \phi \circ \psi $ of $ \phi, \psi \in \Ortho(q) $ is the map $ \map{}{}{}{x}{\phi(\psi(x))} $.
\end{definition}

\begin{remark}
	Let $ (\module, q) $ be a quadratic module and let $ \phi \in \Ortho(q) $. Then we also have $ f(\phi(u), \phi(v)) = f(u,v) $ for all $ u,v \in \module $ where $ f $ is the linearisation of $ q $.
\end{remark}

\begin{note}
	Recall from \cref{rootsys:convention-right,rootsys:right-convention} that in the context of root systems in a Euclidean space $ (V, \cdot) $, we always let automorphisms of $ V $ act from the right side. The reason for this is that we often need reflections in the context of Weyl elements, whose conjugation action is also written on the right side. Since this specific reasoning does not apply to quadratic forms, we have the freedom to let $ \Ortho(q) $ act on $ \module $ from the left side, \enquote{as usual}.
\end{note}

\begin{definition}[Orthogonality]
	For any quadratic module $ (\module, q) $, we say that $ v,u \in \module $ are \defemph*{orthogonal}\index{orthogonality} if $ f(u,v) = 0 $ where $ f $ is the linearisation of $ q $. For any $ \comring $-subspace $ U $ of $ \module $, we define $ U^\perp \defl \Set{v \in \module \given f(u,v) = 0 \text{ for all } u \in U} $ and we put $ v^\perp \defl \gen{v}_\comring^\perp $ for all $ v \in \module $.
\end{definition}

Even the notion of reflections, which we defined in \cref{rootsys:refl-def} for inner products, can be generalised to the setting of quadratic forms. However, we have to be more careful at this point. Recall from \cref{rootsys:refl-def} that the reflection corresponding to a non-zero vector $ v $ in a Euclidean space $ (V, \cdot) $ is the map
\[ \map{\refl{v}}{V}{V}{x}{x^{\refl{v}} \defl x - 2 \frac{x \cdot v}{v \cdot v} v}. \]
We can write down the same definition for arbitrary symmetric bilinear forms $ f $ as long as we assume that $ f(v,v) $ is invertible. However, since $ f(v,v) = 2q(v) $ by \thmitemcref{quadmod:basiclem}{quadmod:basiclem:fvv}, this assumption is never satisfied if $ 2_\comring $ is not invertible. The solution to this problems consists of removing the coefficient $ 2 $ from the definition of $ \refl{v} $. If $ 2_\comring $ is invertible, then the resulting definition is equivalent to the naive definition above.

\begin{definition}[Reflection, {\cite[7.2]{ElmanKarpenkoMerkurjev-QuadForms}}, {\cite[(6.1.2)]{Knus_QuadFormRing}}]\label{quadmod:refl-def}
	Let $ (\module, q) $ be a quadratic module over $ \comring $, denote by $ f $ the linearisation of $ q $ and let $ v \in \module $ such that $ q(v) $ is invertible in $ \comring $. The \defemph*{reflection corresponding to $ v $}\index{reflection} is the map
	\[ \map{\refl{v}}{\module}{\module}{u}{u - q(v)^{-1} f(u,v) v}. \]
\end{definition}

Reflection in quadratic modules satisfy some of the main properties that we are used to from reflections in Euclidean spaces. The most important caveat is that in general, we cannot expect the module $ \module $ to decompose as $ \gen{v}_\comring \dirsum v^\perp $.

\begin{lemma}\label{quadmod:refl-O}
	Let $ (\module, q) $ be a quadratic module over $ \comring $, denote by $ f $ the linearisation of $ q $ and let $ v \in \module $ such that $ q(v) $ is invertible in $ \comring $. Then $ \refl{v} $ lies in $ \Ortho(q) $ and satisfies $ \refl{v}(u) = u $ for all $ u \in v^\perp $, $ \refl{v}(v) = -v $ and $ \refl{v}^2 = \id_\module $.
\end{lemma}
\begin{proof}
	It is clear that $ \refl{v} $ is $ \comring $-linear and that it is the identity on $ v^\perp $. Further, since $ f(v,v) = 2 q(v) $ by \thmitemcref{quadmod:basiclem}{quadmod:basiclem:fvv}, we have
	\[ \refl{v}(v) = v - q(v)^{-1} f(v,v)v = v - 2v = -v. \]
	This implies that for all $ u \in \module $,
	\begin{align*}
		\refl{v}^2(u) &= \refl{v}\brackets[\big]{u - q(v)^{-1} f(u,v) v} = \refl{v}(u) - q(v)^{-1} f(u,v) \refl{v}(v) \\
		&= u - q(v)^{-1} f(u,v) v + q(v)^{-1} f(u,v) v = u.
	\end{align*}
	Finally, by \thmitemcref{quadmod:basiclem}{quadmod:basiclem:lin},
	\begin{align*}
		q\brackets[\big]{\refl{v}(u)} &= q\brackets[\big]{u - q(v)^{-1} f(u,v) v} \\
		&= q(u) + q\brackets[\big]{- q(v)^{-1} f(u,v) v} + f\brackets[\big]{u, -q(v)^{-1} f(u,v)v} \\
		&= q(u) + q(v)^{-2} f(u,v)^2 q(v) - q(v)^{-1} f(u,v) f(u,v) = q(u).
	\end{align*}
	This finishes the proof.
\end{proof}

We now turn to pointed quadratic modules.

\begin{definition}[Pointed quadratic module]
	A \defemph*{pointed quadratic module over $ \comring $}\index{quadratic module!pointed} is a triple $ (\module, q, e) $ where $ (\module, q) $ is a quadratic module over $ \comring $ and $ e $ is an element of $ \module $ with $ q(e) = 1_\comring $. The element $ e $ is called the \defemph*{base point}\index{quadratic module!base point}. A quadratic module $ (\module, q) $ is called \defemph*{pointable}\index{quadratic module!pointable} if there exists $ e \in \module $ such that $ (\module, q, e) $ is a pointed quadratic module.
\end{definition}

\begin{note}
	Recall from \cref{ADE:ring-identity} that in simply-laced root graded groups, the ring element $ 1_\comring $ which appears in the decomposition of the fixed Weyl elements turns out to be the unit element in the coordinatising ring. We will observe a similar phenomenon in $ B_n $-graded groups: The module element $ e $ which appears in the decomposition of the fixed short Weyl element will ultimately satisfy $ q(e) = 1_\comring $. Thus the quadratic module that we construct from any $ B_n $-graded group is pointable, and the choice of a $ \rootbase $-system of Weyl elements determines a base point. However, note that the unit element of a ring is uniquely determined by the ring structure whereas a pointable quadratic module may have many possible base points.
\end{note}

\begin{note}[see {\cite[11.14]{GPR_AlbertRing}}]\label{quadmod:pointed-tr-conj}
	Let $ (\module, q, e) $ be a pointed quadratic module and denote by $ f $ the linearisation of $ q $. One can define a \defemph*{trace map}\index{trace!on a pointed quadratic module}
	\[ \map{\compTr}{\module}{\comring}{u}{f(e,u) = f(u,e)} \]
	and a \defemph*{conjugation map}\index{conjugation!on a pointed quadratic module}
	\[ \map{\compinvmap}{\module}{\module}{u}{\compinv{u} \defl \compTr(u)e - u}. \]
	Observe that, since $ q(e) = 1_\comring $, we have $ \compinvmap = -\refl{e} $ where $ \refl{e} $ denotes the reflection corresponding to $ e $.
	
	We will not need the trace and conjugation maps in this generality. However, we will encounter them in \cref{conic:conj-tr-def} in the special case that $ (\module, q) $ is a conic algebra. Further, we will see similar concepts in \cref{sec:rinv}.
\end{note}

Finally, we can define the standard parameter system for pointed quadratic modules.

\begin{definition}[Standard parameter system]\label{quadmod:standard-param}
	Let $ (\module, q, e) $ be a pointed qua\-dra\-tic module. Put $ \twistgroup \defl \invogroup \defl \compactSet{\pm 1} $. Declare that $ \twistgroup $ acts on $ \comring $ and $ \module $ by inversion and that $ \invogroup $ acts trivially on $ \comring $ and by $ \refl{e} $ on $ \module $. More precisely, this means that
	\[ -1_\twistgroup.a = -a, \quad -1_\twistgroup.m=-m, \quad -1_\invogroup.a=a, \quad -1_\invogroup.m = \refl{e}.m \]
	for all $ a \in \comring $ and all $ m \in \module $. Then the triple $ (\twistgroup \times \invogroup, \module, \comring) $ is called the \defemph*{standard parameter system for $ (\module, q, e) $}.\index{parameter system!standard!type Bn@type $ B_n $}
\end{definition}

\begin{remark}
	To obtain an induced action of $ \twistgroup \times \invogroup $ on $ \module $ and $ \comring $, we have to observe that the actions of $ \twistgroup $ and $ \invogroup $ commute. This holds because the reflection map $ \refl{v_0} $ is linear.
\end{remark}


\section{Root Systems of Type \texorpdfstring{ $ B $}{B}}

\label{sec:B:rootsys}

\begin{secnotation}
	We denote by $ n $ an integer at least $ 2 $. When we say \enquote{root}, we always mean \enquote{root in $ B_n $} unless otherwise specified.
\end{secnotation}

In this section, we collect some basic facts about the root system $ B_n $ which will be needed later on. Note that we exclude the case $ n=1 $ because the resulting root system $ B_1 $ would be isomorphic to $ A_1 $. As usual, most results follow from a straightforward inspection of the standard representation of $ B_n $.

\begin{remark}[Standard representation of $ B_n $]\label{B:Bn-standard-rep}
	Let $ V $ be a Euclidean space of dimension $ n $ with orthonormal basis $ \tup{\basvec}{n} $. The \defemph*{standard representation of $ B_n $}\index{standard representation!of Bn@of $ B_n $} is
	\[ B_n \defl \Set{\epsilon_1 \basvec_i +\epsilon_2 \basvec_j \given i \ne j \in \numint{1}{n}, \epsilon_1, \epsilon_2 \in \compactSet{\pm 1}} \union \Set{\epsilon \basvec_i \given i \in \numint{1}{n}, \epsilon \in \compactSet{\pm 1}}. \]
	The long roots are exactly those which lie in the first set and the short roots are exactly those in the second one. The \defemph*{standard root base} is
	\[ \rootbase \defl \Set{\basvec_i - \basvec_{i+1} \given i \in \numint{1}{n-1}} \union \compactSet{\basvec_n} \]
	and the corresponding positive system is
	\[ \possys \defl \Set{\basvec_i - \basvec_{j} \given i<j \in \numint{1}{n}} \union \Set{\basvec_i + \basvec_j \given i \ne j \in \numint{1}{n}} \union \Set{\basvec_i \given i \in \numint{1}{n}}. \]
\end{remark}

\begin{definition}[$ B_2 $-pairs and $ B_2 $-quadruples]\label{B:B2-pair-def}
	Let $ \roots $ be any root system. A \defemph*{$ B_2 $-pair (in $ \roots $)}\index{B2-pair@$ B_2 $-pair} is a pair $ (\alpha, \delta) $ of roots such that $ \gen{\alpha, \delta}_\IR \intersect \roots $ is a root subsystem of $ \roots $ of type $ B_2 $ with root base $ (\alpha, \delta) $ with $ \alpha $ being the long simple root in this subsystem and $ \delta $ being the short simple root. A \defemph*{$ B_2 $-quadruple (in $ \roots $)}\index{B2-quadruple@$ B_2 $-quadruple} is a quadruple $ (\alpha, \beta, \gamma, \delta) $ of roots in $ \roots $ such that $ (\alpha, \delta) $ is a $ B_2 $-pair, $ \beta = \alpha + \delta $ and $ \gamma = \alpha + 2\delta $.
\end{definition}

\premidfigure
\begin{figure}[htb]
	\centering\begin{tikzpicture}
		\draw[->] (0,0) -- (1,-1);
		\draw[->] (0,0) -- (1,0);
		\draw[->] (0,0) -- (1,1);
		\draw[->] (0,0) -- (0,1);
		\draw[->] (0,0) -- (-1,1);
		\draw[->] (0,0) -- (-1,0);
		\draw[->] (0,0) -- (-1,-1);
		\draw[->] (0,0) -- (0,-1);
		
		\node[below right] at (1,-1){$ \alpha $};
		\node[right] at (1,0){$ \beta $};
		\node[above right] at (1,1){$ \gamma $};
		\node[above] at (0,1){$ \delta $};
		\node[above left] at (-1,1){$ -\alpha $};
		\node[left] at (-1,0){$ -\beta $};
		\node[below left] at (-1,-1){$ -\gamma $};
		\node[below] at (0,-1){$ -\delta $};
	\end{tikzpicture}
	\caption{A $ B_2 $-quadruple $ (\alpha, \beta, \gamma, \delta) $.}
	\label{fig:B2-quad}
\end{figure}
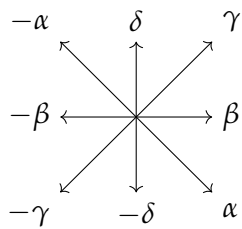
\postmidfigure

\begin{remark}
	In a $ B_2 $-quadruple $ (\alpha, \beta, \gamma, \delta) $, the roots $ \alpha, \beta, \gamma, \delta $ are exactly the positive roots of the corresponding $ B_2 $-subsystem (with respect to the root base $ (\alpha, \delta) $). All reflection maps on this subsystem can be read off from the diagram in \cref{fig:B2-quad}.
\end{remark}

\begin{remark}[compare \cref{BC:C2-cry-crit}]\label{B:B2-cry-crit}
	Let $ G $ be a group with a $ B_2 $-grading $ (\rootgr{\alpha})_{\alpha \in B_2} $ and let $ (\alpha, \beta, \gamma, \delta) $ be a $ B_2 $-quadruple. Then this grading is crystallographic if and only if the commutators $ \commutator{\rootgr{\epsilon\alpha}}{\rootgr{\sigma\gamma}} $ are trivial for all $ \epsilon,\sigma \in \compactSet{\pm 1} $.
\end{remark}

We now state some simple results, mostly on the possible ways in which to roots can \enquote{relate} to each other. Since every pair of roots lies in a parabolic rank-2 subsystem, there is only a handful of cases which can occur.

\begin{remark}\label{B:parabolic-sub}
	Every $ A_2 $-subsystem and every $ B_2 $-subsystem of $ B_n $ is parabolic. The same is not true for subsystems of type $ A_1 \times A_1 $, however: The short roots in $ B_2 $ form such a subsystem which is neither parabolic nor crystallographically closed.
\end{remark}

\begin{lemma}\label{B:rootsys:short-ortho}
	Let $ \beta $, $ \delta $ be two short roots such that $ \beta \nin \compactSet{\pm \delta} $. Then there exist long roots $ \alpha $, $ \gamma $ such that $ (\alpha, \beta, \gamma, \delta) $ is a $ B_2 $-quadruple. In fact, we must have $ \alpha = \beta - \delta $ and $ \gamma = \beta + \delta $.
\end{lemma}

\begin{lemma}\label{B:rootsys:ortho}
	Let $ \rho, \zeta $ be two orthogonal roots in $ B_n $. Then one of the following two conditions is satisfied:
	\begin{stenumerate}
		\item \label{B:rootsys:ortho:cry-adj}$ \rho $ is crystallographically adjacent to $ \zeta $ and $ -\zeta $.
		
		\item \label{B:rootsys:ortho:short}$ \rho $ and $ \zeta $ are both short and there exist long roots $ \alpha, \gamma $ such that $ (\alpha, \rho, \gamma, \zeta) $ is a $ B_2 $-quadruple.
	\end{stenumerate}
\end{lemma}
\begin{proof}
	Denote by $ \roots' $ the root subsystem which is spanned by $ \rho $, $ \zeta $. Then $ \roots' $ must be of type $ A_1 \times A_1 $, $ A_2 $ or $ B_2 $. In the first case,~\itemref{B:rootsys:ortho:cry-adj} holds. The second case is not possible because there are no orthogonal roots in $ A_2 $. Now assume that we are in the third case. Then $ \rho $ and $ \zeta $ must have the same length. If they are both short, then~\itemref{B:rootsys:ortho:short} holds by \cref{B:rootsys:short-ortho}. Otherwise~\itemref{B:rootsys:ortho:cry-adj} holds. This finishes the proof.
\end{proof}

\begin{lemma}\label{B:rootsys:long-in-A2}
	If $ n \ge 3 $, then every long root in $ B_n $ lies in an $ A_2 $-subsystem.
\end{lemma}

\begin{lemma}\label{B:rootsys:commonA2}
	If $ n \ge 3 $, then two long roots $ \alpha, \gamma $ in $ B_n $ lie in a common $ A_2 $-subsystem if and only if they are not orthogonal.
\end{lemma}

\begin{note}
	If $ n=2 $, then long roots $ \alpha, \gamma $ with $ \gamma \in \Set{\pm \alpha} $ do not lie in a common $ A_2 $-subsystem (because $ B_2 $ has no $ A_2 $-subsystem) but they are not orthogonal. Hence we have to assume that $ n \ge 3 $ in \cref{B:rootsys:commonA2}.
\end{note}

\begin{lemma}\label{B:rootsys:long-short-ortho}
	Let $ \alpha $ be a long root and let $ \beta $ be a short root. Then $ \alpha, \beta $ lie in a common $ B_2 $-subsystem if and only if they are not orthogonal.
\end{lemma}

In the computation of the Cartan integers, it is practical to consider each orbit of roots separately.

\begin{lemma}\label{B:rootsys:cartan-int-longlong}
	Let $ \alpha, \gamma $ be two long roots. Then the Cartan integer $ \cartanint{\alpha}{\gamma} = 2 \frac{\alpha \cdot \gamma}{\gamma \cdot \gamma} $ is determined as follows:
	\begin{lemenumerate}
		\item \label{B:rootsys:cartan-int-longlong:0}$ \cartanint{\alpha}{\gamma} = 0 $ if and only if $ \alpha $ and $ \gamma $ are orthogonal. (If $ n  \ge 3 $, this is equivalent by \cref{B:rootsys:commonA2} to $ \alpha $ and $ \gamma $ lying in no common $ A_2 $-subsystem).
		
		\item $ \cartanint{\alpha}{\gamma} = -1 $ if and only if $ (\alpha, \gamma) $ is an $ A_2 $-pair.
		
		\item $ \cartanint{\alpha}{\gamma} = 1 $ if and only if $ (\alpha, -\gamma) $ is an $ A_2 $-pair.
		
		\item $ \cartanint{\alpha}{\gamma} = 2 $ if and only if $ \alpha = \gamma $.
		
		\item $ \cartanint{\alpha}{\gamma} = -2 $ if and only if $ \alpha=-\gamma $.
	\end{lemenumerate}
\end{lemma}

\begin{lemma}\label{B:rootsys:cartan-int-shortshort}
	Let $ \beta, \delta $ be two short roots. Then the Cartan integer $ \cartanint{\beta}{\delta} $ is determined as follows:
	\begin{lemenumerate}
		\item $ \cartanint{\beta}{\delta} = 0 $ if and only if $ \beta $ and $ \delta $ are orthogonal (or equivalently, if and only if $ \beta \nin \Set{\pm \delta} $).
		
		\item $ \cartanint{\beta}{\delta} = 2 $ if and only if $ \beta = \delta $.
		
		\item $ \cartanint{\beta}{\delta} = -2 $ if and only if $ \beta = -\delta $.
	\end{lemenumerate}
\end{lemma}

\begin{lemma}\label{B:rootsys:cartan-int-longshort}
	Let $ \alpha $ be a long root and let $ \beta $ be a short root. Then the Cartan integer $ \cartanint{\alpha}{\beta} $ is determined as follows:
	\begin{lemenumerate}
		\item $ \cartanint{\alpha}{\beta} = 0 $ if and only if $ \alpha $ and $ \beta $ are orthogonal (or, equivalently by \cref{B:rootsys:long-short-ortho}, if and only if they do not lie in a common $ B_2 $-subsystem).
		
		\item $ \cartanint{\alpha}{\beta} = 2 $ if and only if $ (\alpha, \beta) $ is a $ B_2 $-pair.
		
		\item $ \cartanint{\alpha}{\beta} = -2 $ if and only if $ (\alpha, -\beta) $ is a $ B_2 $-pair.
	\end{lemenumerate}
\end{lemma}

\begin{lemma}
	Let $ \alpha $ be a long root and let $ \beta $ be a short root. Then the Cartan integer $ \cartanint{\beta}{\alpha} $ is determined as follows:
	\begin{lemenumerate}
		\item $ \cartanint{\beta}{\alpha} = 0 $ if and only if $ \alpha $ and $ \beta $ are orthogonal (or equivalently, if and only if they do not lie in a common $ B_2 $-subsystem).
		
		\item $ \cartanint{\beta}{\alpha} = 1 $ if and only if $ (\alpha, \beta) $ is a $ B_2 $-pair.
		
		\item $ \cartanint{\beta}{\alpha} = -1 $ if and only if $ (\alpha, -\beta) $ is a $ B_2 $-pair.
	\end{lemenumerate}
\end{lemma}

Since we will usually only be interested in the number $ (-1)^{\cartanint{\alpha}{\beta}} $ (or, in other words, the parity of the Cartan integer), the following summary will be useful.

\begin{proposition}\label{B:rootsys:cartan-int-parity}
	Let $ \rho, \zeta $ be two roots in $ B_n $. Then $ \cartanint{\rho}{\zeta} $ is an even number if and only if one of the following conditions is satisfied:
	\begin{stenumerate}
		\item $ \zeta $ is short.
	
		\item $ \rho \in \Set{\pm \zeta} $.
		
		\item $ \rho $ and $ \zeta $ are orthogonal.
	\end{stenumerate}
	Equivalently, we can replace the last statement by the following one:
	\begin{lemenumerate}
		\item[(iii')] $ \rho $ is long and and $ \rho, \zeta $ do not lie in a common subsystem of type $ A_2 $ or $ B_2 $.
	\end{lemenumerate}
\end{proposition}

We end this section with the definition of a certain subset $ \Bnsub $ of $ B_n $ which will prove useful later on.

\begin{remark}\label{B:Bnsub-def}
	Let $ (w_\delta)_{\delta \in \rootbase} $ be a $ \rootbase $-system of Weyl elements in a $ B_n $-graded group and let $ \alpha, \beta $ be roots in the same orbit. Then there exists a word $ \word{\delta} $ over $ \rootbase $ such that $ \rootgr{\alpha}^{w_{\word{\delta}}} = \rootgr{\beta} $, but the word $ \word{\delta} $ might be rather long. For this reason, it will in some situations be useful to consider a $ \Bnsub $-extension of $ (w_\delta)_{\delta \in \rootbase} $ (in the sense of the following \cref{B:Bnsub-ext-def}) where $ \Bnsub $ is the subset
	\[ \Bnsub \defl \Set{\basvec_i - \basvec_j \given i \ne j \in \numint{1}{n}} \union \Set{\basvec_i \given i \in \numint{1}{n}} \]
	of $ B_n $.\index{Bn@$ \Bnsub $} In other words, $ \Bnsub $ is the union of the canonical subsystem of type $ A_{n-1} $ and the set of positive short roots. Observe that this set is neither closed nor a root system.
	
	Now the set $ \bar{S} \defl \Set{\refl{\alpha} \given \alpha \in \Bnsub} $ is a practical set of generators for the Weyl group for the following reason: For all roots $ \alpha, \beta \in B_n $ which lie in the same orbit, there exist $ m \in \numint{0}{4} $ and a word $ \word{\delta} $ over $ \Bnsub $ of length $ m $ such that $ \alpha^{\reflbr{\word{\delta}}} = \beta $. The idea here is to use the subgroup $ \Weyl(A_{n-1}) $ of $ \Weyl(B_n) $ to permute the basis vectors $ (\basvec_i)_{i \in \numint{1}{n}} $ (which requires at most two reflections) and the reflections $ (\refl{\basvec_i})_{i \in \numint{1}{n}} $ to change the signs of the basis vectors (which also requires at most two reflections).
\end{remark}

\begin{definition}[Standard $ \rootbase $-expression]\label{B:Bnsub-ext-word}
	Let $ \rootbase $ denote the standard root base of $ B_n $. For any root $ \alpha \in \Bnsub $, we define a $ \rootbase $-expression $ \word{\rho}^{\alpha} $ of $ \alpha $ (in the sense of \cref{param:Delta-expr}) as follows:
	\begin{defenumerate}
		\item If $ \alpha \in \rootbase $, we put $ \word{\rho}^\alpha \defl (\alpha) $.
		
		\item If $ \alpha = \basvec_i $ for some $ i \in \numint{1}{n-1} $, we define $ \word{\rho}^{\basvec_i} \defl (\basvec_i - \basvec_{i+1}, \word{\rho}^{\basvec_{i+1}}, \basvec_{i+1}-\basvec_i) $.
		
		\item If $ \alpha = \basvec_i - \basvec_j $ for some $ i<j \in \numint{1}{n-1} $ with $ i+1<j $, we define $ \word{\rho}^{\basvec_i - \basvec_j} \defl (\basvec_j - \basvec_{j-1}, \word{\rho}^{\basvec_i - \basvec_{j-1}}, \basvec_{j-1} - \basvec_j) $.
		
		\item If $ \alpha = \basvec_j - \basvec_i $ for some $ i<j \in \numint{1}{n-1} $ with $ i+1<j $, we define $ \word{\rho}^{\basvec_j - \basvec_i} \defl (\word{\rho}^{\basvec_i - \basvec_j})^{-1} $.
	\end{defenumerate}
	The word $ \word{\rho}^{\alpha} $ will also be called the \defemph*{standard $ \rootbase $-expression of $ \alpha $}\index{standard Delta-expression@standard $ \rootbase $-expression}. We also put $ \word{\rho}^i \defl \word{\rho}^{\basvec_i} $ for all $ i \in \numint{1}{n} $ and $ \word{\rho}^{ij} \defl \word{\rho}^{\basvec_i - \basvec_j} $ for all distinct $ i,j \in \numint{1}{n} $.
\end{definition}

\begin{definition}[$ \Bnsub $-extensions]\label{B:Bnsub-ext-def}
	Let $ \rootbase $ denote the standard base of $ B_n $, let $ G $ be any group with a $ B_n $-pregrading $ (\hat{U}_\alpha)_{\alpha \in B_n} $ and let $ (w_\delta)_{\delta \in \rootbase} $ be a $ \rootbase $-system of Weyl elements in $ G $. Then we define a family $ (w_\alpha)_{\alpha \in \Bnsub} $, called the \defemph*{standard $ \Bnsub $-extension of $ (w_\delta)_{\delta \in \rootbase} $}\index{Bn-extension@$ \Bnsub $-extension}, by $ w_\alpha \defl w_{\word{\rho}^{\alpha}} $ for all $ \alpha \in \Bnsub $. We will sometimes write $ w_{ij} $ for $ w_{\basvec_i - \basvec_j} $ and $ w_{i} $ for $ w_{\basvec_i} $.
\end{definition}

\begin{note}
	Since $ \refl{\alpha} = \refl{-\alpha} $ for all roots $ \alpha $, we could remove all negative roots in $ A_{n-1} $ from $ \Bnsub $ without changing the set $ \bar{S} $ of generators in \cref{B:Bnsub-def}. However, it will be more practical to have both Weyl elements $ w_{ij} $ and $ w_{ji} = w_{ij}^{-1} $ available. By contrast, we will show in \cref{B:short-weyl-center} that for all $ i \in \numint{1}{n} $, the short Weyl elements $ w_i $ and $ w_i^{-1} $ act identically on all root groups. Thus there is no need to include the short negative roots in $ \Bnsub $.
\end{note}

\begin{note}
	The main computations during which we use $ \Bnsub $-extensions are \cref{B:comm-mult-computation,B:blue:stand-signs}. The reason why $ \Bnsub $-extensions are useful in this context are \cref{B:ex-parmap-equality,B:Bnsub-conj-anygroup}.
\end{note}


\section{Construction of \texorpdfstring{$ B_n $}{B\_n}-graded Groups}

\label{sec:B-example}

\begin{secnotation}
	We denote by $ \comring $ a commutative associative ring, by $ (\module, q, v_0) $ a pointed quadratic module over $ \comring $ and by $ \map{f}{\module}{\comring}{}{} $ the linearisation of $ q $. Further, we fix an integer $ n \in \IN_{\ge 2} $, the root system $ \roots \defl B_n $ in standard representation (as in \cref{B:Bn-standard-rep}) and the standard root base $ \rootbase $ of $ B_n $. We will at some points consider the subset $ \Bnsub $ of $ B_n $ defined in \cref{B:Bnsub-def}.
\end{secnotation}

	The construction of RGD-systems of type $ B_n $ from a quadratic module over a field is well-known. For example, the required matrices can be found for $ n=2 $ in \cite[\onpage{469}]{Maldeghem-GeneralisedPolygons}. The construction in \cite{Maldeghem-GeneralisedPolygons} is, in fact, more general than what we need because it starts not from a quadratic module, but from a $ \sigma $-quadratic module where $ \sigma $ is an anti-involution of the field. Taking $ \sigma $ to be the identity map, we obtain the regular notion of quadratic modules. The construction outlined in this section is essentially an adaption of the one in \cite{Maldeghem-GeneralisedPolygons} to arbitrary $ n $ and with the field replaced by a commutative associative ring. We will change some minor specifics of the construction, like signs and the order of the matrix rows and columns, but in its essence it is the same.

\begin{convention}
	As in \cref{quadmod:ortho-grp}, we consider endomorphisms of any module to be acting from the left. Thus the composition $ \phi \circ \psi $ of two such endomorphisms is the map $ \map{}{}{}{x}{\phi(\psi(x))} $.
\end{convention}

\begin{note}[Implementation]\label{B:ex:GAP}
	An implementation of the elementary orthogonal group in the computer algebra system GAP \cite{GAP4} can be found in \cite{RGG-GAP}. This includes proofs of all computational results in this section which are not explicitly verified by hand.
\end{note}

\subsection{Construction}

Our example of a $ B_n $-graded group will be a group of orthogonal automorphisms, but not of the quadratic module $ (\module, q) $. Instead, we have to add $ n $ hyperbolic pairs to $ \module $ (in the sense of \cref{quadmod:hyperbolic-ex}).

\begin{construction}
	We put $ V_+ \defl \comring^n $, $ V_- \defl \comring^n $ and $ V \defl \module \dirsum V_+ \dirsum V_- $. We denote by $ (b_1, \ldots, b_n) $ the standard basis of $ V_+ $, by $ (b_{-1}, \ldots, b_{-n}) $ the standard basis of $ V_- $ and we will always consider $ V_+ $ and $ V_- $ to be subsets of $ V $ without specifying the natural embedding. We will usually denote elements of $ V_+, V_- $ by the letters $ v,w $, elements of $ \module $ by the letters $ m,u $ and elements of $ V $ by the letters $ x,y $. Further, we define a map
	\[ \map{}{V}{\comring}{m \dirsum (v_1, \ldots, v_n) \dirsum (v_{-1}, \ldots, v_{-n})}{q(m) + \sum_{i=1}^n v_i v_{-i}} \]
	which clearly extends $ q $, and which we denote by $ q $ as well.
\end{construction}

\begin{lemma}
	The map $ \map{q}{V}{\comring}{}{} $ is a quadratic form on $ V $.
\end{lemma}
\begin{proof}
	In fact, $ V $ is simply the direct sum of $ (\module, q) $ with the hyperbolic space of dimension $ 2n $ from \cref{quadmod:hyperbolic-ex}, so this follows from \cref{quadmod:dirsum}.
\end{proof}

\begin{remark}[Generalised matrices]\label{generalised-matrices}
	Since $ \module $ is not free, we cannot represent automorphisms of $ V $ by matrices with coefficients in $ \comring $. However, using the direct sum decomposition of $ V $, we can represent them by \enquote{generalised matrices}\index{generalised matrix} which are elements of
	\[ \begin{pmatrix}
		\Hom(\module, \module) & \Hom(V_+, \module) & \Hom(V_-, \module) \\
		\Hom(\module, V_+) & \Hom(V_+, V_-) & \Hom(V_-, V_+) \\
		\Hom(\module, V_-) & \Hom(V_+, V_-) & \Hom(V_-, V_-)
	\end{pmatrix}, \]
	where $ \Hom $ stands for the group of $ \comring $-linear homomorphisms. Further, since $ V_+ $ and $ V_- $ are free of rank $ n $, we can represent elements of $ \Hom(V_\epsilon, V_\sigma) $ by $ (n \times n) $-matrices over $ \comring $, elements of $ \Hom(V_\epsilon, \module) $ by $ (1 \times n) $-matrices with entries in $ \module $ and elements of $ \Hom(\module, V_\epsilon) $ by $ (n \times 1) $-matrices with entries in $ \Hom(\module, \comring) $. Here $ \epsilon $ and $ \sigma $ denote arbitrary signs. The elements of $ \Hom(\module, \comring) $ that we encounter will always be of the form $ f(u, \mapdot) $ for some $ u \in \module $, which is the map $ \map{}{}{}{v}{f(u,v)} $.
\end{remark}

We can now construct the desired root homomorphisms.

\begin{construction}[Long root homomorphisms]\label{B:ex-roothom-def-long}
	Let $ i,j \in \numint{1}{n} $ be distinct and let $ a \in \comring $. We denote by $ \risom{\basvec_i - \basvec_j}(a) $ the unique $ \comring $-linear endomorphism of $ V $ which is given on the direct summands of $ V $ as follows, where $ m $, $ v^+ $, $ v^- $ denote arbitrary elements of $ \module $, $ V_+ $, $ V_- $, respectively:
	\[ \risom{\basvec_i - \basvec_j}(a) \colon \quad m \mapsto m, \quad v^+ \mapsto v^+ + a v^+_j b_i, \quad v^- \mapsto v^- - a v^-_i b_{-j}. \]
	Now assume that, in addition, $ i <j $. Then we define $ \risom{\basvec_i + \basvec_j}(a) $ and $ \risom{-\basvec_i - \basvec_j}(a) $ by the following formulas:
	\begin{align*}
		\risom{\basvec_i + \basvec_j}(a) &\colon \quad m \mapsto m, \quad v^+ \mapsto v^+, \quad v^- \mapsto v^- + a v^-_j \cdot b_i - a v^-_i \cdot b_j, \\
		\risom{-\basvec_i - \basvec_j}(a) &\colon \quad m \mapsto m, \quad v^+ \mapsto v^+ + a v^+_i \cdot b_{-j} - a v^+_j \cdot b_{-i}, \quad v^- \mapsto v^- .
	\end{align*}
	Examples of the corresponding generalised matrices are given in \cref{B:ex-matrices-long}.
\end{construction}

\begin{construction}[Short root homomorphisms]\label{B:ex-roothom-def-short}
	Let $ i \in \numint{1}{n} $ and let $ u \in \module $. Then we define $ \risom{\basvec_i}(u) $ and $ \risom{-\basvec_i}(u) $ to be the unique $ \comring $-linear endomorphisms of $ V $ which are given on the direct summands of $ V $ as follows, where $ m $, $ v^+ $, $ v^- $ denote arbitrary elements of $ \module $, $ V_+ $, $ V_- $, respectively:
	\begin{align*}
		\risom{\basvec_i}(u) &\colon \enskip m \mapsto m + f(u, m) \cdot b_i, \enskip v^+ \mapsto v^+, \enskip v^- \mapsto v^- - v^-_i \cdot u - v^-_i q(u) \cdot b_i, \\
		\risom{-\basvec_i}(u) & \colon \enskip m \mapsto m - f(u, m) \cdot b_{-i}, \enskip v^+ \mapsto v^+ + v^+_i \cdot u - v^+_i q(v) \cdot b_{-i}, \enskip v^- \mapsto v^-.
	\end{align*}
	Examples of the corresponding generalised matrices are given in \cref{B:ex-matrices-short}.
\end{construction}

\premidfigure

\begin{figure}[htb]
	\centering$ \begin{aligned}
		\risom{\basvec_1 - \basvec_2}(a) &= \brackets*{\begin{array}{c|ccc|ccc}
			\id_\module &&&&&& \\
			\hline
			& 1 & a &&&& \\
			& & 1 &&&& \\
			& & & 1 \\
			\hline
			& & & & 1 && \\
			& & & & -a & 1 & \\
			& & & & & & 1
		\end{array}}, \\
		\risom{\basvec_1 + \basvec_2}(a) &= \brackets*{\begin{array}{c|ccc|ccc}
			\id_\module &&&&&& \\
			\hline
			& 1 & &  & & a & \\
			& & 1 & & -a && \\
			& & & 1 \\
			\hline
			& & & & 1 && \\
			& & & &  & 1 & \\
			& & & & & & 1
		\end{array}}, \\
		\risom{-\basvec_1 - \basvec_2}(a) &= \brackets*{\begin{array}{c|ccc|ccc}
			\id_\module &&&&&& \\
			\hline
			& 1 & &  & &  & \\
			& & 1 & &  && \\
			& & & 1 \\
			\hline
			& & -a & & 1 && \\
			& a & & &  & 1 & \\
			& & & & & & 1
		\end{array}}.
	\end{aligned} $
	\caption{The long root homomorphisms for $ B_3 $.}
	\label{B:ex-matrices-long}
\end{figure}

\begin{figure}[htb]
	\centering$ \begin{aligned}
		\risom{\basvec_1}(u) &= \brackets*{\begin{array}{c|ccc|ccc}
			\id_\module & && &-u&& \\
			\hline
			f(u, \mapdot)& 1&& &-q(u)&& \\
			&& 1& &&& \\
			& &&1 &&& \\
			\hline
			& && &1&& \\
			& && &&1& \\
			& && &&&1
		\end{array}}, \\
		\risom{-\basvec_1}(u) &= \brackets*{\begin{array}{c|ccc|ccc}
			\id_\module & u && &&& \\
			\hline
			& 1&& && \\
			& &1&& && \\
			& &&1& && \\
			\hline
			-f(u, \mapdot)& -q(u) && &1&& \\
			& && &&1& \\
			& && &&&1
		\end{array}}.
	\end{aligned} $
	\caption{The short root homomorphisms for $ B_3 $.}
	\label{B:ex-matrices-short}
\end{figure}

\postmidfigure

\begin{lemma}
	For each short root $ \beta $ and for each long root $ \alpha $, the maps $ \map{\risom{\alpha}}{(\comring, +)}{\End_\comring(V)}{}{} $ and $ \map{\risom{\beta}}{(\module, +)}{\End_\comring(V)}{}{} $ are injective homomorphisms. In particular, their images lie in $ \Aut_\comring(V) $.
\end{lemma}
\begin{proof}
	The injectivity of these maps can be deduced from their matrix representations. The homomorphism property follows from a straightforward matrix computation.
\end{proof}

Now the definition of the elementary orthogonal group, our main example of a $ B_n $-graded group, is not surprising.

\begin{definition}[Elementary orthogonal group]\label{B:ex:EOq-def}
	For each root $ \alpha $, we denote the image of $ \risom{\alpha} $ in $ \Aut_\comring(V) $ by $ \rootgr{\alpha} $, and we denote by $ \EO(q) $ the group which is generated by $ (\rootgr{\alpha})_{\alpha \in \roots} $. We call $ \EO(q) $ the \defemph*{elementary orthogonal group of $ q $}\index{elementary orthogonal group}.
\end{definition}

\begin{note}\label{B:ex-is-ortho}
	A straightforward computation shows that each root group is contained in the orthogonal group $ \Ortho(q) $ from \cref{quadmod:ortho-grp}. (Here $ q $ denotes the quadratic form on $ V $ and not the one on $ \module $.) It follows that the elementary orthogonal group $ \EO(q) $ is contained in $ \Ortho(q) $, just like the elementary group from \cref{chev:ex-SL} is contained in the special linear group. However, we will never formally use this fact, so we do not carry out the tedious computation.
\end{note}

\subsection{Weyl Elements}

We can define Weyl elements in $ \EO(q) $ which look just like the ones in Chevalley groups (\cref{chev:weyl-def}), except that the notion of \enquote{inverses of elements of $ \module $} is not obvious.

\begin{definition}[Weyl elements]\label{B:ex-weyl-const}
	For all distinct $ i,j \in \numint{1}{n} $ and all invertible $ a \in \comring $, we define
	\[ w_{ij}(a) \defl w_{\basvec_i - \basvec_j}(a) \defl \risom{\basvec_j - \basvec_i}(-a^{-1}) \circ \risom{\basvec_i - \basvec_j}(a) \circ \risom{\basvec_j - \basvec_i}(-a^{-1}) \in \EO(q) \]
	and $ w_{ij} \defl w_{\basvec_i - \basvec_j} \defl w_{ij}(1_\comring) $. For all $ i \in \numint{1}{n} $ and all $ u \in \module $ for which $ q(u) $ is invertible, we define
	\[ w_i(u) \defl w_{\basvec_i}(u) \defl \risom{-\basvec_i}\brackets[\big]{-q(u)^{-1} u} \circ \risom{\basvec_i}(u) \circ \risom{-\basvec_i}\brackets[\big]{-q(u)^{-1} u} \]
	and $ w_i \defl w_{\basvec_i} \defl w_i(v_0) $.
\end{definition}

\begin{note}[see also~\ref{B:Bnsub-def}]\label{B:ex-weyl-Bnsub-note}
	Observe that we have defined Weyl elements only for the roots in the set $ \Bnsub $ from \cref{B:Bnsub-def}. In theory, it would be enough to only consider the Weyl elements $ w_{i,i+1}(a) $ for $ i \in \numint{1}{n-1} $ and $ w_n(u) $, that is, the ones which correspond to the root base $ \rootbase $. However, by the considerations in \cref{B:Bnsub-def}, it will be more practical to have an explicit $ \alpha $-Weyl element for all $ \alpha \in \Bnsub $. At the same time, while we could easily define Weyl elements for all remaining roots in $ B_n $ as well, there is no additional benefit in doing so.
\end{note}

\begin{definition}[Standard system of Weyl elements]\label{B:ex-standard-weyl}
	The \defemph*{standard system of Weyl elements for $ \EO(q) $} is the family $ (w_\delta)_{\delta \in \rootbase} $ given by \cref{B:ex-weyl-const}.
\end{definition}

\begin{remark}[Short Weyl elements]\label{B:ex:short-weyl-matrix}
	Let $ u \in \module $ such that $ q(u) $ is invertible. The goal of this remark is to explicitly compute the generalised matrix of $ w_1(u) $. We will leave out the rows and columns which correspond to $ (b_i)_{i \in \numint{2}{n}} $ and $ (b_{-i})_{i \in \numint{2}{n}} $ because they are trivial. We begin the computation with
	\begin{align*}
		\risom{\basvec_1}(u) \risom{-\basvec_1}\brackets[\big]{-q(u)^{-1} u} \hspace{-3cm} & \\
		&= \begin{pmatrix}
			\id & 0 & -u \\
			-f(u, \mapdot) & 1 & -q(u) \\
			0 & 0 & 1
		\end{pmatrix} \begin{pmatrix}
			\id & -q(u)^{-1} u & 0 \\
			0 & 1 & 0 \\
			q(u)^{-1} f(u, \mapdot) & -q(u)^{-1} & 1
		\end{pmatrix} \\
		&= \begin{pmatrix}
			\id - q(u)^{-1} f(u, \mapdot) u & -q(u)^{-1}u + q(u)^{-1}u & -u \\
			-f(u, \mapdot) + q(u)q(u)^{-1} f(u, \mapdot) & r(u) & -q(u) \\
			q(u)^{-1} f(u, \mapdot) & -q(u)^{-1} & 1
		\end{pmatrix} \\
		&= \begin{pmatrix}
			\id - q(u)^{-1} f(u, \mapdot) u & 0 & -u \\
			0 & 0 & -q(u) \\
			q(u)^{-1} f(u, \mapdot) & -q(u)^{-1} & 1
		\end{pmatrix}
	\end{align*}
	where $ r(u) \defl -q(u)^{-1}f(u,u) + 1 + q(u)q(u)^{-1} = 0 $ because $ f(u,u) = 2q(u) $ by \thmitemcref{quadmod:basiclem}{quadmod:basiclem:fvv}. Now
	\begin{align*}
		w_1(u) &= \risom{-\basvec_1}\brackets[\big]{-q(u)^{-1} u} \risom{\basvec_1}(u) \risom{-\basvec_1}\brackets[\big]{-q(u)^{-1} u} \\
		&= \begin{pmatrix}
			\id & -q(u)^{-1} u & 0 \\
			0 & 1 & 0 \\
			q(u)^{-1} f(u, \mapdot) & -q(u)^{-1} & 1
		\end{pmatrix} \\
		& \hspace{2cm}\mathord{} \cdot \begin{pmatrix}
			\id - q(u)^{-1} f(u, \mapdot) u & 0 & -u \\
			0 & 0 & -q(u) \\
			q(u)^{-1} f(u, \mapdot) & -q(u)^{-1} & 1
		\end{pmatrix} \\
		&= \begin{pmatrix}
			\id - q(u)^{-1} f(u, \mapdot) u & 0 & u - q(u)^{-1} q(u)u \\
			0 & 0 & -q(u) \\
			\phi & -q(u)^{-1} & -q(u)^{-1} f(u,u) + q(u)^{-1} q(u) + 1
		\end{pmatrix}
	\end{align*}
	where $ \map{\phi}{\module}{\comring}{}{} $ is the map which sends $ x \in \module $ to
	\begin{align*}
		&\hspace{-1.5cm} q(u)^{-1} f\brackets[\big]{u, x - q(u)^{-1} f(u,x)u} + q(u)^{-1} f(u,x) \\
		&= q(u)^{-1} \brackets[\big]{f(u,x) - q(u)^{-1} f(u,x) f(u,u) + f(u,x)} = 0.
	\end{align*}
	We conclude that
	\[ w_1(u) = \begin{pmatrix}
		\id - q(u)^{-1} f(u, \mapdot) u & 0 & 0 \\
		0 & 0 & -q(u) \\
		0 & -q(u)^{-1} & 0
	\end{pmatrix}. \]
	Hence $ w_1(au) = w_1(u) $ for all $ a \in \comring $ with $ a^2 = 1_\comring $. In particular, $ w_1(u) = w_1(-u) $. We conclude that a short Weyl element may in general be represented by distinct Weyl triples.
\end{remark}

The following results show that the elements in \cref{B:ex-weyl-const} are indeed Weyl elements. Even more, we have specific formulas for their conjugation actions. We can see from these formulas that, in addition to additive inversion on $ \comring $ and $ \module $, we can expect a second kind of twisting in $ B_n $-graded groups: the reflection $ \refl{v_0} $ on $ \module $. Note that, unlike the unit element $ 1_\comring $ in the ring $ \comring $, the element $ v_0 $ is not uniquely determined by $ (\module, q) $, but rather an arbitrary (but fixed) element of $ \module $ with the property that $ q(v_0) = 1_\comring $.

\begin{lemma}\label{B:ex-weylformula-long}
	Let $ i,j \in \numint{1}{n} $ be distinct and let $ a \in \comring $ be invertible. Then $ w \defl w_{ij}(a) $ is an $ (\basvec_i - \basvec_j) $-Weyl element. It satisfies the following formulas for all $ b \in \comring $ and $ u \in \module $:
	\begin{lemenumerate}
		\item $ \risom{\basvec_i - \basvec_j}(b)^w = \risom{\basvec_j - \basvec_i}\brackets{-a^{-2}b} $ and $ \risom{\basvec_j - \basvec_i}(b)^w = \risom{\basvec_i - \basvec_j}(-a^2b) $.
		
		\item For all $ k \in \numint{1}{n} \setminus \compactSet{i,j} $, we have
		\begin{align*}
			 \risom{\basvec_i - \basvec_k}(b)^w &= \risom{\basvec_j - \basvec_k}(a^{-1} b), & \risom{\basvec_i + \basvec_k}(b)^w &= \risom{\basvec_j + \basvec_k}(\delmin{k \in \betint{i}{j}} a^{-1} b), \\
			 \risom{\basvec_k - \basvec_i}(b)^w &= \risom{\basvec_k - \basvec_j}(ab), & \risom{-\basvec_i - \basvec_k}(b)^w &= \risom{-\basvec_j - \basvec_k}(\delmin{k \in \betint{i}{j}} ab).
		\end{align*}
		
		\item For all $ k \in \numint{1}{n} \setminus \compactSet{i,j} $, we have
		\begin{align*}
			\risom{\basvec_j - \basvec_k}(b)^w &= \risom{\basvec_i - \basvec_k}(-a b), & \risom{\basvec_j + \basvec_k}(b)^w &= \risom{\basvec_i + \basvec_k}(-\delmin{k \in \betint{i}{j}} a b), \\
			\risom{\basvec_k - \basvec_j}(b)^w &= \risom{\basvec_k - \basvec_i}(-a^{-1}b), & \risom{-\basvec_j - \basvec_k}(b)^w &= \risom{-\basvec_i - \basvec_k}(-\delmin{k \in \betint{i}{j}} a^{-1}b).
		\end{align*}
		
		\item $ \risom{\basvec_i}(u)^w = \risom{\basvec_j}(a^{-1} u) $ and $ \risom{-\basvec_i}(u)^w = \risom{-\basvec_j}(au) $.
		
		\item $ \risom{\basvec_j}(u)^w = \risom{\basvec_i}(-au) $ and $ \risom{-\basvec_j}(u)^w = \risom{-\basvec_i}(-a^{-1}u) $.
		
		\item For any root $ \alpha $ for which no formula for the conjugation action of $ w $ on $ \rootgr{\alpha} $ is given above, this action is trivial.
	\end{lemenumerate}
\end{lemma}
\begin{proof}
	This follows from a straightforward computation. See also \cref{B:ex:GAP}.
\end{proof}

\begin{remark}\label{B:ex-wij-inv}
	Let $ i,j \in \numint{1}{n} $ be distinct and let $ a \in \comring $ be invertible. Since $ w \defl w_{ij}(a) $ is a Weyl element by \cref{B:ex-weylformula-long}, it follows from \thmitemcref{basic:weyl-general}{basic:weyl-general:minus} that
	\begin{align*}
		w_{ij}(a) &= \rismin{i}{j}(a) \rismin{j}{i}(-a^{-1}) \rismin{j}{i}(-a^{-1})^{w} \\
		&= \rismin{i}{j}(a) \rismin{j}{i}(-a^{-1}) \rismin{i}{j}(a) = w_{ji}(-a^{-1}).
	\end{align*}
	In particular, $ w_{ij} = w_{ji}^{-1} $.
\end{remark}

\begin{lemma}\label{B:ex-weylformula-short}
	Let $ i \in \numint{1}{n} $ and let $ v \in \module $ such that $ q(v) $ is invertible. Then $ w \defl w_i(v) $ is an $ \basvec_i $-Weyl element. It satisfies the following formulas for all $ b \in \comring $ and $ u \in \module $:
	\begin{lemenumerate}
		\item For all $ k \in \numint{1}{n} \setminus \compactSet{i} $, we have 
		\begin{align*}
			\risom{\basvec_i - \basvec_k}(b)^w &= \risom{-\basvec_i - \basvec_k}\brackets[\big]{\delmin{i>k} q(v)^{-1} b}, & \risom{\basvec_i + \basvec_k}(b)^w &= \risom{\basvec_k - \basvec_i}\brackets[\big]{\delmin{i>k} q(v)^{-1} b}, \\
			\risom{\basvec_k - \basvec_i}(b)^w &= \risom{\basvec_k + \basvec_i}\brackets[\big]{\delmin{i>k} q(v) b}, & \risom{-\basvec_i - \basvec_k}(b)^w &= \risom{\basvec_i - \basvec_k}\brackets[\big]{\delmin{i>k} q(v) b}.
		\end{align*}
		
		\item \label{B:ex-weylformula-short:on-self}$ \risom{\basvec_i}(u)^w = \risom{-\basvec_i}\brackets[\big]{q(v)^{-1} \refl{v}(u)} $ and $ \risom{-\basvec_i}(u)^w = \risom{\basvec_i}\brackets[\big]{q(v) \refl{v}(u)} $.
		
		\item \label{B:ex-weylformula-short:on-short}For all $ k \in \numint{1}{n} \setminus \compactSet{i} $, $ \risom{\basvec_k}(u)^w = \risom{\basvec_k}\brackets[\big]{\refl{v}(u)} $ and $ \risom{-\basvec_k}(u)^w = \risom{-\basvec_k}\brackets[\big]{\refl{v}(u)} $.
		
		\item For any root $ \alpha $ for which no formula for the conjugation action of $ w $ on $ \rootgr{\alpha} $ is given above, this action is trivial.
	\end{lemenumerate}
\end{lemma}
\begin{proof}
	This follows from a straightforward computation. See also \cref{B:ex:GAP}.
\end{proof}

\subsection{Parity Maps and Twisting Structures}

We can read off the correct definition of the twisting groups and parity maps in $ \EO(q) $ (with respect to the standard system of Weyl elements) from the above formulas. In other words, the parity maps and twisting groups are chosen precisely to ensure that \cref{B:ex-parmap-describes-conj} holds.

\begin{secnotation}\label{B:ex-twistgroups-def}
	From now on, we denote by $ (\twistgroup \times \invogroup, \module, \comring) $ the standard parameter system for $ (\module, q, v_0) $ (as in \cref{quadmod:standard-param}). Further, we define maps $ \map{\inverparsym}{B_n \times \Bnsub}{\twistgroup}{}{} $ and $ \map{\invoparsym}{B_n \times \Bnsub}{\invogroup}{}{} $ by the formulas in \cref{B:ex-parmap-def}, where $ \Bnsub $ is as in \cref{B:Bnsub-def}.
	By restricting the second component to $ \rootbase $, we obtain $ \rootbase $-parity maps which we also denote by $ \inverparsym $ and $ \invoparsym $, and which we call the \defemph*{standard parity maps of type $ B_n $}\index{parity map!standard!type Bn@type $ B_n $}.
\end{secnotation}

\premidfigure
\begin{figure}[htb]
	\centering\begin{tabular}[t]{ccc}
		\toprule
		$ \alpha $ & $ \inverpar{\alpha}{\basvec_i - \basvec_j} $ & $ \invopar{\alpha}{\basvec_i - \basvec_j} $ \\
		\midrule
		$ \pm (\basvec_i - \basvec_j) $ & $ -1_\twistgroup $ & $ 1_\invogroup $ \\
		$ \pm (\basvec_i + \basvec_j) $ & $ 1_\twistgroup $ & $ 1_\invogroup $ \\
		$ \pm (\basvec_i - \basvec_k) $ & $ 1_\twistgroup $ & $ 1_\invogroup $ \\
		$ \pm (\basvec_i + \basvec_k) $ & $ \delmin{k \in \betint{i}{j}} 1_\twistgroup $ & $ 1_\invogroup $ \\
		$ \pm (\basvec_j - \basvec_k) $ & $ -1_\twistgroup $ & $ 1_\invogroup $ \\
		$ \pm (\basvec_j + \basvec_k) $ & $ -\delmin{k \in \betint{i}{j}} 1_\twistgroup $ & $ 1_\invogroup $ \\
		$ \pm \basvec_k \pm \basvec_l $ & $ 1_\twistgroup $ & $ 1_\invogroup $ \\
		$ \pm \basvec_i $ & $ 1_\twistgroup $ & $ 1_\invogroup $ \\
		$ \pm \basvec_j $ & $ -1_\twistgroup $ & $ 1_\invogroup $ \\
		$ \pm \basvec_k $ & $ 1_\twistgroup $ & $ 1_\invogroup $ \\
		\bottomrule
	\end{tabular}\qquad\quad\begin{tabular}[t]{ccc}
		\toprule
		$ \alpha $ & $ \inverpar{\alpha}{\basvec_i} $ & $ \invopar{\alpha}{\basvec_i} $ \\
		\midrule
		$ \pm (\basvec_i - \basvec_k) $ & $ \delmin{i>k} 1_\twistgroup $ & $ 1_\invogroup $ \\
		$ \pm (\basvec_i + \basvec_k) $ & $ \delmin{i>k} 1_\twistgroup $ & $ 1_\invogroup $ \\
		$ \pm \basvec_k \pm \basvec_l $ & $ 1_\twistgroup $ & $ 1_\invogroup $ \\
		$ \pm \basvec_i $ & $ 1_\twistgroup $ & $ -1_\invogroup $ \\
		$ \pm \basvec_k $ & $ 1_\twistgroup $ & $ -1_\invogroup $ \\
		\bottomrule
	\end{tabular}
	\caption{The definition of $ \inverpar{\alpha}{\beta} $ and $ \invopar{\alpha}{\beta} $ for all $ \alpha \in B_n $ and $ \beta \in \Bnsub $, see \cref{B:ex-twistgroups-def}. In the left table, we assume that $ i,j,k,l $ are pairwise distinct. In the right table, we assume that $ i,k,l $ are pairwise distinct. For small values of $ n $, it is of course not possible to choose three or four pairwise distinct indices, in which case the corresponding rows should be ignored.}
	\label{B:ex-parmap-def}
\end{figure}
\postmidfigure

\begin{note}\label{B:ex-parmap-def-Bnsub-note}
	The values $ \inverpar{\alpha}{\beta} $ for $ \beta \in \Bnsub \setminus \rootbase $ in \cref{B:ex-parmap-def} are not relevant for the definition of the $ \rootbase $-parity map $ \inverparsym $ because it is defined on the set $ \roots \times \rootbase $. However, the assertion of the following \cref{B:ex-parmap-describes-conj} holds for all $ \beta \in \Bnsub $, not merely for $ \beta \in \rootbase $. We will use this observation in \cref{B:ex-parmap-equality} to express some values of the extended parity map
	\[ \map{\inverparsym}{\roots \times \frmon{\rootbase \union (-\rootbase)}}{\twistgroup}{}{} \]
	(which we defined in~\ref{param:parmap-def}) in terms of the values $ (\inverpar{\alpha}{\beta})_{\alpha \in B_n, \beta \in \Bnsub} $. We conclude that, while the values $ \inverpar{\alpha}{\beta} $ for $ \beta \in \Bnsub \setminus \rootbase $ are not needed to define the parity map $ \map{\inverparsym}{\roots \times \rootbase}{\twistgroup}{}{} $, they nonetheless express important properties of this map. The fact that we can directly read off these values from \cref{B:ex-parmap-def} will allow us to simplify some later computations, for example in \cref{B:comm-mult-computation}. Similar remarks apply for $ \invoparsym $ in place of $ \inverparsym $.
\end{note}

\begin{remark}\label{B:ex-parmap-restrict}
	Consider $ A_{n-1} = \Set{\basvec_i - \basvec_j \given i \ne j \in \numint{1}{n}} $ as a subset of $ B_n $ in the natural way. Then for all $ \alpha, \beta \in A_{n-1} $, the value $ \inverpar{\alpha}{\beta} $ in \cref{B:ex-parmap-def} is the same as the one in \cref{ADE:An-parmap}.
\end{remark}

\begin{lemma}\label{B:ex-parmap-describes-conj}
	Let $ \alpha \in B_n $ and $ \beta \in \Bnsub $. Let $ x \in \comring $ if $ \alpha $ is long and let $ x \in \module $ if $ \alpha $ is short. Then $ \risom{\alpha}(x)^{w_\beta} = \risom{\refl{\beta}(\alpha)}(\inverpar{\alpha}{\beta} \invopar{\alpha}{\beta} . x) $.
\end{lemma}
\begin{proof}
	This follows from \cref{B:ex-weylformula-long,B:ex-weylformula-short} by putting $ a \defl 1_\comring $ and $ v \defl v_0 $, respectively.
\end{proof}

\begin{remark}\label{B:ex-param-choice}
	Note that, even though the definition of $ \inverparsym $ and $ \invoparsym $ is clearly motivated by \cref{B:ex-weylformula-long,B:ex-weylformula-short}, it is completely independent of the choices of $ \comring $, $ \module $, $ q $ and $ v_0 $. This is evident from \cref{B:ex-parmap-def}. As a consequence, we can obtain information about $ \inverparsym $ and $ \invoparsym $ by performing computations in the group $ \EO(q) $ for any fixed choice of the parameters $ \comring $, $ \module $, $ q $ and $ v_0 $. In the sequel, we will often choose $ \comring \defl \IC $, $ \module \defl \IC^2 $, $ \map{q}{\module}{\IC}{(x,y)}{x^2 + y^2} $ and $ v_0 \defl (1,0) $. This choice of parameters has the following two crucial properties: Firstly, $ 1_\comring \ne -1_\comring $, so that the inversion map on $ (\comring, +) $ is non-trivial; and secondly, that $ \refl{v_0} $ is neither the inversion map nor the identity map on $ (\module, +) $. These properties imply that the standard parameter system $ (\twistgroup \times \invogroup, \module, \comring) $ from \cref{quadmod:standard-param} is $ (\inverparsym \times \invoparsym) $-faithful (in the sense of \cref{param:parsys-faithful}).
\end{remark}

\begin{lemma}\label{B:ex-param-thm}
	The root isomorphisms $ (\risom{\alpha})_{\alpha \in B_n} $ from \cref{B:ex-roothom-def-long,B:ex-roothom-def-short} form a parametrisation of $ \EO(q) $ by $ (\twistgroup \times \invogroup, \module, \comring) $ with respect to $ \inverparsym \times \invoparsym $ and $ (w_\delta)_{\delta \in \rootbase} $.
\end{lemma}
\begin{proof}
	This follows from \cref{B:ex-parmap-describes-conj}.
\end{proof}

\begin{lemma}\label{B:ex-weyl-extension}
	Let $ (w_\alpha)_{\alpha \in \Bnsub} $ be as in \cref{B:ex-weyl-const}. Then the following hold:
	\begin{lemenumerate}
		\item For all pairwise distinct $ i,j,k \in \numint{1}{n} $, we have
		\[ w_{ij}^{w_{jk}} = w_{ik} \midand w_{ij}^{w_{kj}} = w_{ki}. \]
		For all distinct $ i,j \in \numint{1}{n} $, we have
		\[ w_i^{w_{ij}} = w_j \midand w_i^{w_{ji}} = w_i^{-1}. \]
		
		\item $ (w_\alpha)_{\alpha \in \Bnsub} $ is the standard $ \Bnsub $-extension of $ (w_\delta)_{\delta \in \rootbase} $, the standard system of Weyl elements.
	\end{lemenumerate}
\end{lemma}
\begin{proof}
	The first assertion follows from \cref{B:ex-parmap-describes-conj} and an inspection of the values in \cref{B:ex-parmap-def}. The second assertion follows from the first one and from \cref{B:ex-wij-inv}.
\end{proof}

\subsection{Commutator Relations}

We now proceed to show that $ (\rootgr{\alpha})_{\alpha \in B_n} $ is a crystallographic $ B_n $-grading of $ \EO(q) $.

\begin{proposition}\label{B:ex-commrel}
	The group $ \EO(q) $ satisfies the following commutator relations. For all distinct $ i,j \in \numint{1}{n} $, all $ a \in \comring $ and all $ v \in \module $, we have
	\begin{align*}
		\commutator{\rismin{i}{j}(a)}{\risshpos{j}(v)} &= \risshpos{i}(av) \risplus{i}{j}\brackets[\big]{\delmin{i<j} a q(v)}, \\
		\commutator{\rismin{j}{i}(a)}{\risshneg{j}(v)} &= \risshneg{i}(-av) \risminmin{i}{j}\brackets[\big]{\delmin{i<j} aq(v)}, \\
		\commutator{\risplus{i}{j}(a)}{\risshneg{j}(v)} &= \risshpos{i}\brackets[\big]{\delmin{i<j} av} \rismin{i}{j}\brackets[\big]{\delmin{i<j} aq(v)}, \\
		\commutator{\risminmin{i}{j}(a)}{\risshpos{j}(v)} &= \risshneg{i}\brackets[\big]{\delmin{i>j} av} \rismin{j}{i}\brackets[\big]{\delmin{i<j} a q(v)}.
	\end{align*}
	For all distinct $ i,j \in \numint{1}{n} $ and all $ v,w \in \module $, we have
	\begin{align*}
		\commutator{\risshpos{i}(v)}{\risshpos{j}(w)} &= \risplus{i}{j}\brackets[\big]{\delmin{i<j} f(v,w)}, \\
		\commutator{\risshpos{i}(v)}{\risshneg{j}(w)} &= \rismin{i}{j}\brackets[\big]{f(v,w)}, \\
		\commutator{\risshneg{i}(v)}{\risshneg{j}(w)} &= \risminmin{i}{j}\brackets[\big]{\delmin{i>j} f(v,w)}.
	\end{align*}
	For all pairwise distinct $ i,j,k \in \numint{1}{n} $ and all $ a,b \in \comring $, we have
	\begin{align*}
		\commutator{\rismin{i}{j}(a)}{\rismin{j}{k}(b)} &= \rismin{i}{k}(ab), \\
		\commutator{\rismin{i}{j}(a)}{\risplus{j}{k}(b)} &= \risplus{i}{k}\brackets{\delmin{k \in \betint{i}{j}} ab}, \\
		\commutator{\rismin{i}{j}(a)}{\risminmin{k}{i}(b)} &= \risminmin{j}{k}\brackets{\delmin{k \nin \betint{i}{j}} ab}, \\
		\commutator{\risplus{i}{j}(a)}{\risminmin{k}{i}(b)} &= \rismin{j}{k}\brackets{\delmin{i \in \betint{j}{k}} ab}.
	\end{align*}
\end{proposition}
\begin{proof}
	This follows from a straightforward computation. See also \cref{B:ex:GAP}.
\end{proof}

\begin{proposition}\label{B:ex-triangular}
	Let $ \possys $ denote the standard positive system in $ B_n $. Then the group $ \EO(q) $ satisfies $ \rootgr{\possys} \intersect \rootgr{-\possys} = \compactSet{1} $.
\end{proposition}
\begin{proof}
	Recall that the generalised matrices of the root homomorphisms with respect to the decomposition $ V = \module \dirsum V_+ \dirsum V_- $ and with respect to the ordered bases $ \tup{b}{n} $ of $ V_+ $ and $ (b_{-1}, \ldots, b_{-n}) $ of $ V_- $ are given in \cref{B:ex-matrices-short,B:ex-matrices-long}. Consider instead the generalised matrices with respect to the decomposition $ V = V_+ \dirsum \module \dirsum V_- $ and with respect to the ordered bases $ \tup{b}{n} $ of $ V_+ $ and $ (b_{-n}, \ldots, b_{-1}) $ of $ V_- $. For $ n=2 $, some of these matrices are displayed in \cref{B:ex-matrices-triangular}. With respect to this decomposition, all generalised matrices of positive root isomorphisms are upper triangular while all generalised matrices of negative root isomorphisms are lower triangular. The assertion follows.
\end{proof}

\premidfigure
\begin{figure}[htb]
	\centering$ \begin{aligned}[]
		\rismin{1}{2}(a) &= \brackets*{\begin{array}{cc|c|cc}
			1 & a &&& \\
			& 1 &&& \\
			\hline
			&& \id_\module && \\
			\hline
			&&& 1 & -a \\
			&&&& 1
		\end{array}}, \\
		\risplus{1}{2}(a) &= \brackets*{\begin{array}{cc|c|cc}
			1 &&& a & \\
			& 1 &&& -a \\
			\hline
			&& \id_\module && \\
			\hline
			&&& 1 & \\
			&&&& 1
		\end{array}}, \\
		\risminmin{1}{2}(a) &= \brackets*{\begin{array}{cc|c|cc}
			1 &&&& \\
			& 1 &&& \\
			\hline
			&& \id_\module && \\
			\hline
			a &&& 1 & \\
			& -a &&& 1
		\end{array}}, \\
		\risshpos{1}(u) &= \brackets*{\begin{array}{cc|c|cc}
			1 && f(u, \mapdot) && -u \\
			& 1 &&& -q(u) \\
			\hline
			&& \id_\module && \\
			\hline
			&&& 1 & \\
			&&&& 1
		\end{array}}, \\
		\risshneg{1}(u) &= \brackets*{\begin{array}{cc|c|cc}
			1 &&&& \\
			& 1 &&& \\
			\hline
			u && \id_\module && \\
			\hline
			 &&& 1 & \\
			-q(u)&& -f(u, \mapdot)&& 1
		\end{array}}.
	\end{aligned} $
	\caption{The generalised matrices of some root isomorphisms with respect to the decomposition of $ V $ in the proof of \cref{B:ex-triangular}.}
	\label{B:ex-matrices-triangular}
\end{figure}
\postmidfigure

\begin{theorem}\label{B:ex-is-rgg}
	The family $ (\rootgr{\alpha})_{\alpha \in B_n} $ is a crystallographic $ B_n $-grading of $ \EO(q) $.
\end{theorem}
\begin{proof}
	By \cref{B:ex-weylformula-short,B:ex-weylformula-long}, there exist Weyl elements for all simple roots and by \cref{basic:weyl-ex-basis}, it follows that there exist Weyl elements for all roots. The remaining conditions are satisfied by \cref{B:ex-commrel,B:ex-triangular}.
\end{proof}

\subsection{Concluding Remarks}

We state another brief technical observation on the parity maps $ \inverparsym $ and $ \invoparsym $ which will be useful for later computations. Specifically, it will be used in \cref{B:Bnsub-conj-anygroup}. See also \cref{B:Bnsub-def}.

\begin{lemma}\label{B:ex-parmap-equality}
	Let $ \alpha \in \Bnsub $, let $ \word{\alpha} $ be the standard $ \rootbase $-expression of $ \alpha $ in the sense of \cref{B:Bnsub-ext-word} and let $ \zeta $ be an arbitrary root. Then the element $ \inverpar{\zeta}{\alpha} $ from \cref{B:ex-parmap-def} equals the value $ \inverpar{\zeta}{\word{\alpha}} $ of the extended parity map $ \map{\inverparsym}{\roots \times \frmon{\rootbase \union (-\rootbase)}}{\twistgroup}{}{} $ (from \cref{param:parmap-def}). The same assertion holds for $ \invoparsym $ in place of $ \inverparsym $.
\end{lemma}
\begin{proof}
	Observe that the assertion is independent of $ \comring $, $ \module $, $ q $ and $ v_0 $, so we have the freedom to choose these parameters as in \cref{B:ex-param-choice}. Let $ x \in \comring = \IC $ if $ \zeta $ is long and let $ x \in \module = \IC^2 $ if $ \zeta $ is short.
	Recall from \cref{B:ex-weyl-extension} that the family $ (w_\beta)_{\beta \in \Bnsub} $ of Weyl elements in $ \EO(q) $ from \cref{B:ex-weyl-const} coincides with the standard $ \Bnsub $-extension of $ (w_\delta)_{\delta \in \rootbase} $ by \cref{B:ex-weyl-extension}. This implies that $ w_\alpha = w_{\word{\alpha}} $, so $ \risom{\zeta}(x)^{w_\alpha} = \risom{\zeta}(x)^{w_{\word{\alpha}}} $. By \cref{B:ex-parmap-describes-conj}, we have
	\[ \risom{\zeta}(x)^{w_\alpha} = \risom{\zeta^{\reflbr{\alpha}}}(\inverpar{\zeta}{\alpha} \invopar{\zeta}{\alpha}.x). \]
	On the other hand, since $ (\risom{\beta})_{\beta \in B_n} $ is a parametrisation of $ \EO(q) $ by $ (\twistgroup \times \invogroup, \module, \comring) $ with respect to $ \inverparsym \times \invoparsym $ and $ (w_\delta)_{\delta \in \rootbase} $ by \cref{B:ex-param-thm} and since $ \word{\alpha} $ is a $ \rootbase $-expression of $ \alpha $, we also have
	\[ \risom{\zeta}(x)^{w_{\word{\alpha}}} = \risom{\zeta^{\reflbr{\word{\alpha}}}}(\inverpar{\zeta}{\word{\alpha}} \invopar{\zeta}{\word{\alpha}}.x) = \risom{\zeta^{\reflbr{\alpha}}}(\inverpar{\zeta}{\word{\alpha}} \invopar{\zeta}{\word{\alpha}}.x). \]
	As the parameter system $ (\twistgroup \times \invogroup, \module, \comring) $ is $ (\inverparsym \times \invoparsym) $-faithful (in the sense of \cref{param:parsys-faithful}) by \cref{B:ex-param-choice}, we infer that $ \inverpar{\zeta}{\word{\alpha}} \invopar{\zeta}{\word{\alpha}} = \inverpar{\zeta}{\word{\alpha}} \invopar{\zeta}{\word{\alpha}} $, as desired.
\end{proof}


\section{Rank-\texorpdfstring{$ 2 $}{2} Computations, Non-crystallographic Case}

\label{sec:B:rank2-noncry}

\begin{secnotation}\label{B:secnot:rank2-noncry}
	We denote by $ G $ a group which has $ B_2 $-commu\-ta\-tor relations with root groups $ (\rootgr{\alpha})_{\alpha \in B_2} $ (in the sense of \cref{rgg:group-commrel-def}) and we choose pairwise non-proportional roots $ \alpha, \beta, \gamma, \delta $ such that $ (\alpha, \beta, \gamma, \delta) $ is an interval ordering of $ \clrootintcox{\alpha}{\delta} $. Further, we assume that $ G $ is rank-2-injective.
\end{secnotation}

Our rank-2 computations are split into two parts. In this section, we derive several formulas in the higher generality of non-crystallographic commutator relations. In the following section, we will specialise to crystallographic commutator relations to simplify our results. This approach has the advantage that the distinction between long and short roots disappears, as we will explain in \cref{B:root-switch-note}. Thus every computation that can be performed in this generality saves us the trouble of performing two separate computations later on.

\begin{note}
	Throughout this section, we will often use that root group elements from adjacent root groups commute. However, it is a priori not clear (and, in the general situation of \cref{B:secnot:rank2-noncry}, not even true) that the root groups themselves are abelian. This will only be proven in \cref{B:abelian} for the higher-rank case.
\end{note}

\begin{note}\label{B:root-switch-note}
	Since our only assumption on $ (\alpha, \beta, \gamma, \delta) $ is that it is an interval ordering of a closed root interval, and since we do not assume that the $ B_2 $-commutator relations are crystallographic, everything that we prove in this section will also be true for the root quadruple $ (\delta, \gamma, \beta, \alpha) $. In particular, we do not assume that $ (\alpha, \beta, \gamma, \delta) $ is a $ B_2 $-quadruple, so we do not assume that $ \alpha $, $ \gamma $ are long and $ \beta $, $ \delta $ are short. 
\end{note}

Before we can perform any serious computations, we have to derive an analogue of \cref{basic:comm-add} for open root intervals with two elements.

\begin{lemma}\label{B:comm-add}
	The following relations hold for all $ x_\alpha, x_\alpha' \in \rootgr{\alpha} $ and all $ y_\delta, y_\delta' \in \rootgr{\delta} $:
	\begin{lemenumerate}
		\item \label{B:comm-add:delta-beta}$ \commpart{x_\alpha}{y_\delta y_\delta'}{\beta} = \commpart{x_\alpha}{y_\delta'}{\beta} \commpart{x_\alpha}{y_\delta}{\beta} $. In particular, $ \commpart{x_\alpha}{y_\delta^{-1}}{\beta} = \commpart{x_\alpha}{y_\delta}{\beta}^{-1} $.
		
		\item \label{B:comm-add:delta-beta2}$ \commpart{y_\delta y_\delta'}{x_\alpha}{\beta} = \commpart{y_\delta}{x_\alpha}{\beta} \commpart{y_\delta'}{x_\alpha}{\beta} $. In particular, $ \commpart{y_\delta^{-1}}{x_\alpha}{\beta} = \commpart{y_\delta}{x_\alpha}{\beta}^{-1} $.
		
		\item \label{B:comm-add:delta-gamma}$ \commpart{x_\alpha}{y_\delta y_\delta'}{\gamma} = \commpart{x_\alpha}{y_\delta'}{\gamma} \commutator[\big]{\commpart{x_\alpha}{y_\delta}{\beta}}{y_\delta'} \commpart{x_\alpha}{y_\delta}{\gamma} $.
		
		\item \label{B:comm-add:delta-gamma2}$ \commpart{y_\delta y_\delta'}{x_\alpha}{\gamma} = \commpart{y_\delta}{x_\alpha}{\gamma} \commutator[\big]{\commpart{y_\delta}{x_\alpha}{\beta}}{y_\delta'} \commpart{y_\delta'}{x_\alpha}{\gamma} $.
		
		\item \label{B:comm-add:alpha-beta}$ \commpart{x_\alpha x_\alpha'}{y_\delta}{\beta} = \commpart{x_\alpha}{y_\delta}{\beta} \commutator[\big]{\commpart{x_\alpha}{y_\delta}{\gamma}}{x_\alpha'} \commpart{x_\alpha'}{y_\delta}{\beta} $.
		
		\item \label{B:comm-add:alpha-beta2}$ \commpart{y_\delta}{x_\alpha x_\alpha'}{\beta} = \commpart{y_\delta}{x_\alpha'}{\beta} \commutator[\big]{\commpart{y_\delta}{x_\alpha}{\gamma}}{x_\alpha'} \commpart{y_\delta}{x_\alpha}{\beta} $.
		
		\item \label{B:comm-add:alpha-gamma}$ \commpart{x_\alpha x_\alpha'}{y_\delta}{\gamma} = \commpart{x_\alpha}{y_\delta}{\gamma} \commpart{x_\alpha'}{y_\delta}{\gamma} $. In particular, $ \commpart{x_\alpha^{-1}}{y_\delta}{\gamma} = \commpart{x_\alpha}{y_\delta}{\gamma}^{-1} $.
		
		\item \label{B:comm-add:alpha-gamma2}$ \commpart{y_\delta}{x_\alpha x_\alpha'}{\gamma} = \commpart{y_\delta}{x_\alpha'}{\gamma} \commpart{y_\delta}{x_\alpha}{\gamma} $. In particular, $ \commpart{y_\delta}{x_\alpha^{-1}}{\gamma} = \commpart{y_\delta}{x_\alpha}{\gamma}^{-1} $.
	\end{lemenumerate}
\end{lemma}
\begin{proof}
	Using \thmitemcref{group-rel}{group-rel:add}, we compute that
	\begin{align*}
		\commutator{x_\alpha}{y_\delta y_\delta'} &= \commutator{x_\alpha}{y_\delta'} \cdot \commutator{x_\alpha}{y_\delta}^{y_\delta'} = \commpart{x_\alpha}{y_\delta'}{\beta} \commpart{x_\alpha}{y_\delta'}{\gamma} \commpart{x_\alpha}{y_\delta}{\beta}^{y_\delta'} \commpart{x_\alpha}{y_\delta}{\gamma}^{y_\delta'}.
	\end{align*}
	Applying \thmitemcref{group-rel}{group-rel:conj1} and the fact that $ \rootgr{\gamma} $ commutes with $ \rootgr{\delta} $ and with $ \rootgr{\beta} $, we infer that
	\begin{align*}
		\commutator{x_\alpha}{y_\delta y_\delta'} &= \commpart{x_\alpha}{y_\delta'}{\beta} \commpart{x_\alpha}{y_\delta'}{\gamma} \commpart{x_\alpha}{y_\delta}{\beta} \commutator[\big]{\commpart{x_\alpha}{y_\delta}{\beta}}{y_\delta'} \commpart{x_\alpha}{y_\delta}{\gamma} \\
		&= \commpart{x_\alpha}{y_\delta'}{\beta} \commpart{x_\alpha}{y_\delta}{\beta} \cdot \commpart{x_\alpha}{y_\delta'}{\gamma} \commutator[\big]{\commpart{x_\alpha}{y_\delta}{\beta}}{y_\delta'} \commpart{x_\alpha}{y_\delta}{\gamma}.
	\end{align*}
	Since the product map on $ (\beta, \gamma) $ is injective by assumption, this implies that~\itemref{B:comm-add:delta-beta} and~\itemref{B:comm-add:delta-gamma} hold. Using \cref{commpart:inv-switch} and, in the second computation, \cref{basic:comm-add}, we infer that
	\begin{align*}
		\commpart{y_\delta y_\delta'}{x_\alpha}{\beta} &= \commpart{x_\alpha}{y_\delta y_\delta'}{\beta}^{-1} = \commpart{x_\alpha}{y_\delta}{\beta}^{-1} \commpart{x_\alpha}{y_\delta'}{\beta}^{-1} = \commpart{y_\delta}{x_\alpha}{\beta} \commpart{y_\delta'}{x_\alpha}{\beta} \rightand \\
		\commpart{y_\delta y_\delta'}{x_\alpha}{\gamma} &= \commpart{x_\alpha}{y_\delta y_\delta'}{\gamma}^{-1} = \commpart{x_\alpha}{y_\delta}{\gamma}^{-1} \commutator[\big]{\commpart{x_\alpha}{y_\delta}{\beta}}{y_\delta'}^{-1} \commpart{x_\alpha}{y_\delta'}{\gamma}^{-1} \\
		&=\commpart{y_\delta}{x_\alpha}{\gamma} \commutator[\big]{\commpart{y_\delta}{x_\alpha}{\beta}}{y_\delta'} \commpart{y_\delta'}{x_\alpha}{\gamma},
	\end{align*}
	thereby proving~\itemref{B:comm-add:delta-beta2} and~\itemref{B:comm-add:delta-gamma2}.
	
	By \cref{B:root-switch-note}, the results of the previous paragraph also hold for $ (\delta, \gamma, \beta, \alpha) $ in place of $ (\alpha, \beta, \gamma, \delta) $. Thus~\itemref{B:comm-add:delta-beta} implies~\itemref{B:comm-add:alpha-gamma2}, \itemref{B:comm-add:delta-beta2} implies~\itemref{B:comm-add:alpha-gamma}, \itemref{B:comm-add:delta-gamma} implies~\itemref{B:comm-add:alpha-beta2} and~\itemref{B:comm-add:delta-gamma2} implies~\itemref{B:comm-add:alpha-beta}.
\end{proof}

Using \cref{B:comm-add}, we can give an alternative proof of \cref{tw:6.4} for $ n=4 $, just like we did in \cref{A2Weyl:basecomp-cor} for $ n=3 $.

\begin{lemma}\label{B:basecomp-walpha-cox-cor-delta}
	Assume that $ (a_{-\alpha}, b_\alpha, c_{-\alpha}) $ is an $ \alpha $-Weyl triple and denote by $ w_\alpha \defl a_{-\alpha} b_\alpha c_{-\alpha} $ the corresponding Weyl element. Then the following statements hold for all $ x_\delta \in \rootgr{\delta} $:
	\begin{lemenumerate}
		\item \label{B:basecomp-walpha-cox-cor-delta:1}$ x_\delta^{w_\alpha} = \commpart{x_\delta}{b_\alpha}{\beta} $.
		
		\item $ \commutator[\big]{x_\delta}{\commpart{x_\delta}{b_\alpha}{\beta}} \commpart{x_\delta}{b_\alpha}{\gamma} \commpart[\big]{\commpart{x_\delta}{b_\alpha}{\beta}}{c_{-\alpha}}{\gamma}^{\commpart{x_\delta}{b_\alpha}{\gamma}} = 1_G $.
		
		\item $ x_\delta \commpart[\big]{\commpart{x_\delta}{b_\alpha}{\beta}}{c_{-\alpha}}{\delta} \commutator{\commpart{x_\delta}{b_\alpha}{\gamma}}{c_{-\alpha}} = 1_G $.
	\end{lemenumerate}
\end{lemma}
\begin{proof}
	Using \thmitemcref{group-rel}{group-rel:conj1}, we compute that
	\begin{align*}
		x_\delta^{w_\alpha} &= x_\delta^{b_\alpha c_{-\alpha}} = \brackets[\big]{x_\delta \commutator{x_\delta}{b_\alpha}}^{c_{-\alpha}} = x_\delta \commutator{x_\delta}{b_\alpha} \commutator[\big]{\commutator{x_\delta}{b_\alpha}}{c_{-\alpha}}.
	\end{align*}
	We want to express $ x_\delta^{w_\alpha} $ as a product of root group elements, so we use the maps from \cref{basic:commpart-def} to split each factor into its root group components:
	\begin{align*}
		x_\delta^{w_\alpha} &= x_\delta \commpart{x_\delta}{b_\alpha}{\beta} \commpart{x_\delta}{b_\alpha}{\gamma} \commutator[\big]{\commpart{x_\delta}{b_\alpha}{\beta} \commpart{x_\delta}{b_\alpha}{\gamma}}{c_{-\alpha}}.
	\end{align*}
	By \thmitemcref{group-rel}{group-rel:add},
	\begin{align*}
		\commutator[\big]{\commpart{x_\delta}{b_\alpha}{\beta} \commpart{x_\delta}{b_\alpha}{\gamma}}{c_{-\alpha}} &= \commutator[\big]{\commpart{x_\delta}{b_\alpha}{\beta}}{c_{-\alpha}}^{\commpart{x_\delta}{b_\alpha}{\gamma}} \commutator[\big]{\commpart{x_\delta}{b_\alpha}{\gamma}}{c_{-\alpha}} \\
		&= \commpart[\big]{\commpart{x_\delta}{b_\alpha}{\beta}}{c_{-\alpha}}{\gamma}^{\commpart{x_\delta}{b_\alpha}{\gamma}} \commpart[\big]{\commpart{x_\delta}{b_\alpha}{\beta}}{c_{-\alpha}}{\delta} \commutator[\big]{\commpart{x_\delta}{b_\alpha}{\gamma}}{c_{-\alpha}}.
	\end{align*}
	Further, $ x_\delta \commpart{x_\delta}{b_\alpha}{\beta} = \commpart{x_\delta}{b_\alpha}{\beta} x_\delta \commutator{x_\delta}{\commpart{x_\delta}{b_\alpha}{\beta}} $ by \thmitemcref{group-rel}{group-rel:comm}. Thus altogether, using that $ \rootgr{\gamma} $ commutes with $ \rootgr{\beta} $ and $ \rootgr{\delta} $, we have
	\[ x_\delta^{w_\alpha} = \commpart{x_\delta}{b_\alpha}{\beta} \cdot y_\gamma \cdot y_\delta \]
	where
	\begin{equation}\label{eq:B:basecomp-walpha-cox-cor-delta}
		\begin{aligned}
			y_\gamma &\defl \commutator[\big]{x_\delta}{\commpart{x_\delta}{b_\alpha}{\beta}} \commpart{x_\delta}{b_\alpha}{\gamma} \commpart[\big]{\commpart{x_\delta}{b_\alpha}{\beta}}{c_{-\alpha}}{\gamma}^{\commpart{x_\delta}{b_\alpha}{\gamma}} \in \rootgr{\gamma} \rightand \\
			y_\delta &\defl x_\delta \commpart[\big]{\commpart{x_\delta}{b_\alpha}{\beta}}{c_{-\alpha}}{\delta} \commutator[\big]{\commpart{x_\delta}{b_\alpha}{\gamma}}{c_{-\alpha}} \in \rootgr{\delta}.
		\end{aligned}
	\end{equation}
	Now the assertion follows from the injectivity of the product map on $ (\beta, \gamma, \delta) $ and the fact that $ x_\delta^{w_\alpha} $ lies in $ \rootgr{\gamma} $.
\end{proof}

Similarly to the strategy described in \cref{A2Weyl:spirit}, we can apply \cref{B:basecomp-walpha-cox-cor-delta} to the root quadruple $ (-\alpha, \delta, \gamma, \beta) $, using the $ (-\alpha) $-Weyl triples given by \thmitemcref{basic:weyl-general}{basic:weyl-general:minus}.

\begin{lemma}\label{B:basecomp-walpha-cox-cor-beta}
	Assume that $ (a_{-\alpha}, b_\alpha, c_{-\alpha}) $ is an $ \alpha $-Weyl triple and denote by $ w_\alpha \defl a_{-\alpha} b_\alpha c_{-\alpha} $ the corresponding Weyl element. Then the following statements hold for all $ x_\beta \in \rootgr{\beta} $:
	\begin{lemenumerate}
		\item $ x_\beta^{w_\alpha} = \commpart{x_\beta}{a_{-\alpha}}{\delta} = \commpart{x_\beta}{c_{-\alpha}}{\delta} $.
		
		\item $ \commutator[\big]{x_\beta}{\commpart{x_\beta}{a_{-\alpha}}{\delta}} \commpart{x_\beta}{a_{-\alpha}}{\gamma} \commpart[\big]{\commpart{x_\beta}{a_{-\alpha}}{\delta}}{b_\alpha}{\gamma}^{\commpart{x_\beta}{a_{-\alpha}}{\gamma}} = 1_G $.
		
		\item $ x_\beta \commpart[\big]{\commpart{x_\beta}{a_{-\alpha}}{\delta}}{b_\alpha}{\beta} \commutator{\commpart{x_\beta}{a_{-\alpha}}{\gamma}}{b_\alpha} = 1_G $.
	\end{lemenumerate}
\end{lemma}
\begin{proof}
	By \thmitemcref{basic:weyl-general}{basic:weyl-general:minus},
	\[ (c_{-\alpha}^{w_\alpha^{-1}}, a_{-\alpha}, b_\alpha) \midand (b_\alpha, c_{-\alpha}, a_{-\alpha}^{w_\alpha}) \]
	are $ (-\alpha) $-Weyl triples corresponding to the $ (-\alpha) $-Weyl element $ w_\alpha $. Applying \cref{B:basecomp-walpha-cox-cor-delta} to the second Weyl triple and to the root quadruple $ (-\alpha, \delta, \gamma, \beta) $, we see that $ x_\beta^{w_\alpha} = \commpart{x_\delta}{c_{-\alpha}}{\delta} $. Applying \cref{B:basecomp-walpha-cox-cor-delta} to the first Weyl triple, we obtain the remaining assertions.
\end{proof}

Using \cref{B:basecomp-walpha-cox-cor-delta}, we obtain a first formula for the action of squares of Weyl elements.

\begin{lemma}\label{B:walpha-cox-on-beta-comm}
	Let $ (a_{-\alpha}, b_\alpha, c_{-\alpha}) $ be an $ \alpha $-Weyl triple, denote by $ w_\alpha \defl a_{-\alpha} b_\alpha c_{-\alpha} $ the corresponding Weyl element and put $ w \defl w_\alpha^2 $. Then $ \commpart{x_\delta}{b_\alpha^{-1}}{\beta}^w = \commpart{x_\delta}{b_\alpha}{\beta} $ for all $ x_\delta \in \rootgr{\delta} $.
\end{lemma}
\begin{proof}
	Let $ x_\delta \in \rootgr{\delta} $. By \thmitemcref{B:basecomp-walpha-cox-cor-delta}{B:basecomp-walpha-cox-cor-delta:1}, we have
	\[ x_\delta^{w_\alpha} = \commpart{x_\delta}{b_\alpha}{\beta} \midand x_\delta^{w_\alpha^{-1}} = \commpart{x_\delta}{b_\alpha^{-1}}{\beta}, \]
	applying \thmitemcref{basic:weyl-general}{basic:weyl-general:inv} in the second case. Since $ w_\alpha^2 $ clearly maps $ x_\delta^{w_\alpha^{-1}} $ to $ x_\delta^{w_\alpha} $, the assertion follows.
\end{proof}

\begin{note}\label{B:walpha-cox-on-beta-comm-note}
	In \cref{B:walpha-cox-on-beta-comm}, we are not able to prove the stronger statement that $ \commpart{x_\delta}{x_\alpha}{\beta}^w = \commpart{x_\delta}{x_\alpha^{-1}}{\beta} $
	for all $ x_\alpha \in \rootgr{\alpha} $ and $ x_\delta \in \rootgr{\delta} $. However, this is true in the more specialised setting of crystallographic $ BC_2 $-graded groups: See \cref{BC:square-act:beta-on-beta,BC:square-act:alpha-beta-on-beta}.
\end{note}

\begin{proposition}\label{B:rootisom-beta-delta}
	For all $ b_\alpha \in \invset{\alpha} $, the map $ \map{}{\rootgr{\delta}}{\rootgr{\beta}}{x_\delta}{\commpart{x_\delta}{b_\alpha}{\beta}} $ is an isomorphism of root groups.
\end{proposition}
\begin{proof}
	Choose $ a_{-\alpha}, c_{-\alpha} $ such that $ w_\alpha \defl a_{-\alpha} b_\alpha c_{-\alpha} $ is an $ \alpha $-Weyl element. By \cref{B:basecomp-walpha-cox-cor-delta}, the map above is simply the map $ \map{}{}{}{x_\delta}{x_\delta^{w_\alpha}} $, which is an isomorphism by the definition of $ \invset{\alpha} $.
\end{proof}

We end this section with some purely notational observations which will make future referencing easier. They follow from the previous results by replacing $ (\alpha, \beta, \gamma, \delta) $ with $ (\delta, \gamma, \beta, \alpha) $.

\begin{lemma}\label{B:basecomp-wdelta-cox-alpha}
	Assume that $ (a_{-\delta}, b_\delta, c_{-\delta}) $ is a $ \delta $-Weyl triple and denote by $ w_\delta \defl a_{-\delta} b_\delta c_{-\delta} $ the corresponding Weyl element. Then the following statements hold for all $ x_\alpha \in \rootgr{\alpha} $:
	\begin{lemenumerate}
		\item $ x_\alpha^{w_\delta} = \commpart{x_\alpha}{b_\delta}{\gamma} $.
		
		\item $ \commutator[\big]{x_\alpha}{\commpart{x_\alpha}{b_\delta}{\gamma}} \commpart{x_\alpha}{b_\delta}{\beta} \commpart[\big]{\commpart{x_\alpha}{b_\delta}{\gamma}}{c_{-\delta}}{\beta}^{\commpart{x_\alpha}{b_\delta}{\beta}} = 1_G $.
		
		\item $ x_\alpha \commpart[\big]{\commpart{x_\alpha}{b_\delta}{\gamma}}{c_{-\delta}}{\alpha} \commutator{\commpart{x_\alpha}{b_\delta}{\beta}}{c_{-\delta}} = 1_G $.
	\end{lemenumerate}
\end{lemma}
\begin{proof}
	This is a reformulation of \cref{B:basecomp-walpha-cox-cor-delta}, using \cref{B:root-switch-note}.
\end{proof}

\begin{lemma}\label{B:basecomp-wdelta-cox-gamma}
	Assume that $ (a_{-\delta}, b_\delta, c_{-\delta}) $ is an $ \delta $-Weyl triple and denote by $ w_\delta \defl a_{-\delta} b_\delta c_{-\delta} $ the corresponding Weyl element. Then the following statements hold for all $ x_\gamma \in \rootgr{\gamma} $:
	\begin{lemenumerate}
		\item $ x_\gamma^{w_\delta} = \commpart{x_\alpha}{a_{-\delta}}{\alpha} = \commpart{x_\alpha}{c_{-\delta}}{\alpha} $.
		
		\item $ \commutator[\big]{x_\gamma}{\commpart{x_\gamma}{a_{-\delta}}{\alpha}} \commpart{x_\gamma}{a_{-\delta}}{\beta} \commpart[\big]{\commpart{x_\gamma}{a_{-\delta}}{\alpha}}{b_\delta}{\beta}^{\commpart{x_\gamma}{a_{-\delta}}{\beta}} = 1_G $.
		
		\item $ x_\gamma \commpart[\big]{\commpart{x_\gamma}{a_{-\delta}}{\alpha}}{b_\delta}{\gamma} \commutator{\commpart{x_\gamma}{a_{-\delta}}{\beta}}{b_\delta} = 1_G $.
	\end{lemenumerate}
\end{lemma}
\begin{proof}
	This is a reformulation of \cref{B:basecomp-walpha-cox-cor-beta}, using \cref{B:root-switch-note}.
\end{proof}

\begin{proposition}\label{B:rootisom:alpha-gamma}
	For all $ b_\delta \in \invset{\delta} $, the map $ \map{}{\rootgr{\alpha}}{\rootgr{\gamma}}{x_\alpha}{\commpart{x_\alpha}{b_\delta}{\gamma}} $ is an isomorphism of root groups.
\end{proposition}
\begin{proof}
	This is a reformulation of \cref{B:rootisom-beta-delta}, using \cref{B:root-switch-note}.
\end{proof}


\section{Rank-\texorpdfstring{$ 2 $}{2} Computations, Crystallographic Case}

\label{sec:B:rank2-cry}

\begin{secnotation}
	We denote by $ G $ a group which has crystallographic $ B_2 $-commutator relations with root groups $ (\rootgr{\alpha})_{\alpha \in B_2} $ and we choose a $ B_2 $-quadruple $ (\alpha, \beta, \gamma, \delta) $. Further, we assume that $ G $ is rank-2-injective.
\end{secnotation}

In this section, we mostly simplify the results from the previous section under the additional assumption that the commutator relations are crystallographic. We will also perform a handful of new computations. Due to the crystallographic assumption, the distinction of root lengths matters from now on. In particular, the roles of $ (\alpha, \beta, \gamma, \delta) $ and $ (\delta, \gamma, \beta, \alpha) $ are no longer interchangeable.

\begin{note}
	Since the root system $ B_2 $ is isomorphic to $ C_2 $, all results in this section are also valid for crystallographic $ C_2 $-gradings. We will explain how this works in \cref{BC:CisBC:C2}.
\end{note}

\begin{lemma}\label{B:comm-add-cry}
	The following relations hold for all $ x_\alpha, x_\alpha' \in \rootgr{\alpha} $ and all $ y_\delta, y_\delta' \in \rootgr{\delta} $:
	\begin{lemenumerate}
		\item \label{B:comm-add-cry:delta-beta}$ \commpart{x_\alpha}{y_\delta y_\delta'}{\beta} = \commpart{x_\alpha}{y_\delta'}{\beta} \commpart{x_\alpha}{y_\delta}{\beta} $. In particular, $ \commpart{x_\alpha}{y_\delta^{-1}}{\beta} = \commpart{x_\alpha}{y_\delta}{\beta}^{-1} $.
		
		\item \label{B:comm-add-cry:delta-beta2}$ \commpart{y_\delta y_\delta'}{x_\alpha}{\beta} = \commpart{y_\delta}{x_\alpha}{\beta} \commpart{y_\delta'}{x_\alpha}{\beta} $. In particular, $ \commpart{y_\delta^{-1}}{x_\alpha}{\beta} = \commpart{y_\delta}{x_\alpha}{\beta}^{-1} $.
		
		\item \label{B:comm-add-cry:delta-gamma}$ \commpart{x_\alpha}{y_\delta y_\delta'}{\gamma} = \commpart{x_\alpha}{y_\delta'}{\gamma} \commutator[\big]{\commpart{x_\alpha}{y_\delta}{\beta}}{y_\delta'} \commpart{x_\alpha}{y_\delta}{\gamma} $.
		
		\item \label{B:comm-add-cry:delta-gamma2}$ \commpart{y_\delta y_\delta'}{x_\alpha}{\gamma} = \commpart{y_\delta}{x_\alpha}{\gamma} \commutator[\big]{\commpart{y_\delta}{x_\alpha}{\beta}}{y_\delta'} \commpart{y_\delta'}{x_\alpha}{\gamma} $.
		
		\item \label{B:comm-add-cry:alpha-beta}$ \commpart{x_\alpha x_\alpha'}{y_\delta}{\beta} = \commpart{x_\alpha}{y_\delta}{\beta} \commpart{x_\alpha'}{y_\delta}{\beta} $. In particular, $ \commpart{x_\alpha^{-1}}{y_\delta}{\beta} = \commpart{x_\alpha}{y_\delta}{\beta}^{-1} $.
		
		\item \label{B:comm-add-cry:alpha-beta2}$ \commpart{y_\delta}{x_\alpha x_\alpha'}{\beta} = \commpart{y_\delta}{x_\alpha'}{\beta} \commpart{y_\delta}{x_\alpha}{\beta} $. In particular, $ \commpart{y_\delta}{x_\alpha^{-1}}{\beta} = \commpart{y_\delta}{x_\alpha}{\beta}^{-1} $.
		
		\item \label{B:comm-add-cry:alpha-gamma}$ \commpart{x_\alpha x_\alpha'}{y_\delta}{\gamma} = \commpart{x_\alpha}{y_\delta}{\gamma} \commpart{x_\alpha'}{y_\delta}{\gamma} $. In particular, $ \commpart{x_\alpha^{-1}}{y_\delta}{\gamma} = \commpart{x_\alpha}{y_\delta}{\gamma}^{-1} $.
		
		\item \label{B:comm-add-cry:alpha-gamma2}$ \commpart{y_\delta}{x_\alpha x_\alpha'}{\gamma} = \commpart{y_\delta}{x_\alpha'}{\gamma} \commpart{y_\delta}{x_\alpha}{\gamma} $. In particular, $ \commpart{y_\delta}{x_\alpha^{-1}}{\gamma} = \commpart{y_\delta}{x_\alpha}{\gamma}^{-1} $.
	\end{lemenumerate}
\end{lemma}
\begin{proof}
	These statements follow from \cref{B:comm-add} with the additional assumption that the commutator relations are crystallographic. More precisely, we only need that $ \rootgr{\alpha} $ and $ \rootgr{\gamma} $ commute and the only changes are in~\itemref{B:comm-add-cry:alpha-beta} and~\itemref{B:comm-add-cry:alpha-beta2}.
\end{proof}

Since there are two orbits of roots and we want to investigate the actions of (squares of) Weyl elements on all root groups, we have to distinguish four cases. For better clarity, we study each of these cases in a separate subsection. The case of long Weyl elements acting on long root groups will not be covered in this section because we can reduce it to the study of $ A_2 $-subsystems in \cref{sec:B-rank3}.

\subsection{The Action of Long Weyl Elements on Short Root Groups}

\label{subsec:B:rank2:long-on-short}

In this subsection, we show that squares of long Weyl elements act on the short root groups in the same $ B_2 $-subsystem by inversion (\cref{B:weylalpha-square-on-beta}). 

\begin{lemma}\label{B:basecomp-walpha-cor:delta}
	Assume that $ (a_{-\alpha}, b_\alpha, c_{-\alpha}) $ is an $ \alpha $-Weyl triple with corresponding Weyl element $ w_\alpha \defl a_{-\alpha} b_\alpha c_{-\alpha} $. Then the following hold:
	\begin{lemenumerate}
		\item \label{B:basecomp-walpha-cor:delta:act}$ x_\delta^{w_\alpha} = \commpart{x_\delta}{b_\alpha}{\beta} $.
		
		\item $ \commutator[\big]{x_\delta}{\commpart{x_\delta}{b_\alpha}{\beta}} \commpart{x_\delta}{b_\alpha}{\gamma} \commpart[\big]{\commpart{x_\delta}{b_\alpha}{\beta}}{c_{-\alpha}}{\gamma}^{\commpart{x_\delta}{b_\alpha}{\gamma}} = 1_G $.
		
		\item $ x_\delta = \commpart{\commpart{x_\delta}{b_\alpha}{\beta}}{c_{-\alpha}}{\delta}^{-1} $.
	\end{lemenumerate}
\end{lemma}
\begin{proof}
	Using that $ \rootgr{-\alpha} $ and $ \rootgr{\gamma} $ commute in a crystallographic $ B_2 $-grading, these statements follow from \cref{B:basecomp-walpha-cox-cor-delta}.
\end{proof}

\begin{lemma}\label{B:basecomp-walpha-cor:beta}
	Assume that $ (a_{-\alpha}, b_\alpha, c_{-\alpha}) $ is an $ \alpha $-Weyl triple with corresponding Weyl element $ w_\alpha \defl a_{-\alpha} b_\alpha c_{-\alpha} $. Then the following hold:
	\begin{lemenumerate}
		\item $ x_\beta^{w_\alpha} = \commpart{x_\beta}{a_{-\alpha}}{\delta} = \commpart{x_\beta}{c_{-\alpha}}{\delta} $.
		
		\item $ \commutator[\big]{x_\beta}{\commpart{x_\beta}{a_{-\alpha}}{\delta}} \commpart{x_\beta}{a_{-\alpha}}{\gamma} \commpart[\big]{\commpart{x_\beta}{a_{-\alpha}}{\delta}}{b_\alpha}{\gamma}^{\commpart{x_\beta}{a_{-\alpha}}{\gamma}} = 1_G $.
		
		\item $ x_\beta = \commpart{\commpart{x_\beta}{a_{-\alpha}}{\delta}}{b_\alpha}{\beta}^{-1} $.
	\end{lemenumerate}
\end{lemma}
\begin{proof}
	Using that $ \rootgr{\alpha} $ and $ \rootgr{\gamma} $ commute in a crystallographic $ B_2 $-grading, these statements follow from \cref{B:basecomp-walpha-cox-cor-beta}.
\end{proof}

\begin{lemma}\label{B:weylalpha-square-on-beta}
	Assume that $ w_\alpha $ is an $ \alpha $-Weyl element. Then $ w_\alpha^2 $ acts on $ \rootgr{\beta} $ and $ \rootgr{\delta} $ by inversion.
\end{lemma}
\begin{proof}
	We apply the same strategy as in the proof of \cref{A2Weyl:square-act-lemma-A2}. Choose an $ \alpha $-Weyl triple $ (a_{-\alpha}, b_\alpha, c_{-\alpha}) $ such that $ w_\alpha = a_{-\alpha} b_\alpha c_{-\alpha} $ and let $ x_\delta $ be an arbitrary element of $ \rootgr{\delta} $. Put $ w \defl w_\alpha^2 $. Using first \cref{B:basecomp-walpha-cor:delta} and then \cref{B:basecomp-walpha-cor:beta}, we have
	\[ x_\delta^{w} = \commpart{x_\delta}{b_\alpha}{\beta}^{w_\alpha} = \commpart{\commpart{x_\delta}{b_\alpha}{\beta}}{c_{-\alpha}}{\delta}. \]
	Again by \cref{B:basecomp-walpha-cor:delta}, we infer that $ x_\delta^w = x_\delta^{-1} $, so $ w $ acts on $ \rootgr{\delta} $ by inversion. Since $ w_\alpha $ is also a $ (-\alpha) $-Weyl element by \thmitemcref{basic:weyl-general}{basic:weyl-general:minus} and $ (-\alpha, \delta, \gamma, \beta) $ is a $ B_2 $-quadruple, it follows that $ w $ acts on $ \rootgr{\beta} $ by inversion as well.
\end{proof}

\begin{proposition}\label{B:short-abelian}
	Assume that $ \invset{\alpha} $ is non-empty. Then $ \rootgr{\beta} $ and $ \rootgr{\delta} $ are abelian.
\end{proposition}
\begin{proof}
	A group $ H $ is abelian if and only if the inversion map $ \map{}{}{}{h}{h^{-1}} $ is an endomorphism (and thus an automorphism) of $ H $. By \cref{B:weylalpha-square-on-beta}, the inversion maps on $ \rootgr{\beta} $ and $ \rootgr{\delta} $ are the same as conjugation by $ w_\alpha^2 $ for any $ \alpha $-Weyl element $ w_\alpha $, and conjugation by any group element is a homomorphism. The assertion follows.
\end{proof}

\begin{note}\label{B:abelian-note}
	In the rank-2 setting, we cannot prove that the long root groups are abelian as well, but this poses no problem: In the rank-3 setting, every long root is contained in an $ A_2 $-subsystem by \cref{B:rootsys:long-in-A2}, so every long root group is abelian by \cref{ADE:abelian}.
\end{note}

\subsection{The Action of Short Weyl Elements on Long Roots}

\label{subsec:B:rank2:short-on-long}

In this subsection, we obtain no conclusive results. We will revisit this case in the rank-3 setting.

\begin{lemma}\label{B:basecomp-wdelta-alpha}
	Assume that $ (a_{-\delta}, b_\delta, c_{-\delta}) $ is a $ \delta $-Weyl triple and denote by $ w_\delta \defl a_{-\delta} b_\delta c_{-\delta} $ the corresponding Weyl element. Then the following statements hold for all $ x_\alpha \in \rootgr{\alpha} $:
	\begin{lemenumerate}
		\item \label{B:basecomp-wdelta-alpha:act}$ x_\alpha^{w_\delta} = \commpart{x_\alpha}{b_\delta}{\gamma} $.
		
		\item \label{B:basecomp-wdelta-alpha:rel}$ \commpart{x_\alpha}{b_\delta}{\beta} \commpart[\big]{\commpart{x_\alpha}{b_\delta}{\gamma}}{c_{-\delta}}{\beta}^{\commpart{x_\alpha}{b_\delta}{\beta}} = 1_G $.
		
		\item $ x_\alpha \commpart[\big]{\commpart{x_\alpha}{b_\delta}{\gamma}}{c_{-\delta}}{\alpha} \commutator{\commpart{x_\alpha}{b_\delta}{\beta}}{c_{-\delta}} = 1_G $.
	\end{lemenumerate}
\end{lemma}
\begin{proof}
	This follows from \cref{B:basecomp-wdelta-cox-alpha}, using the additional information that $ \rootgr{\alpha} $ and $ \rootgr{\gamma} $ commute in a crystallographic $ B_2 $-grading.
\end{proof}

\begin{lemma}\label{B:basecomp-wdelta-gamma}
	Assume that $ (a_{-\delta}, b_\delta, c_{-\delta}) $ is an $ \delta $-Weyl triple and denote by $ w_\delta \defl a_{-\delta} b_\delta c_{-\delta} $ the corresponding Weyl element. Then the following statements hold for all $ x_\gamma \in \rootgr{\gamma} $:
	\begin{lemenumerate}
		\item \label{B:basecomp-wdelta-gamma:act}$ x_\gamma^{w_\delta} = \commpart{x_\gamma}{a_{-\delta}}{\alpha} = \commpart{x_\gamma}{c_{-\delta}}{\alpha} $.
		
		\item $ \commpart{x_\gamma}{a_{-\delta}}{\beta} \commpart[\big]{\commpart{x_\gamma}{a_{-\delta}}{\alpha}}{b_\delta}{\beta}^{\commpart{x_\gamma}{a_{-\delta}}{\beta}} = 1_G $.
		
		\item $ x_\gamma \commpart[\big]{\commpart{x_\gamma}{a_{-\delta}}{\alpha}}{b_\delta}{\gamma} \commutator{\commpart{x_\gamma}{a_{-\delta}}{\beta}}{b_\delta} = 1_G $.
	\end{lemenumerate}
\end{lemma}

\begin{proof}
	This follows from \cref{B:basecomp-wdelta-cox-alpha}, using the additional information that $ \rootgr{\alpha} $ and $ \rootgr{\gamma} $ commute in a crystallographic $ B_2 $-grading.
\end{proof}

\subsection{The Action of Short Weyl Elements on Short Root Groups}

Again, there are no conclusive results in this subsection because rank-3 assumptions are necessary.

\begin{lemma}\label{B:weylbeta-on-delta-nonabelian}
	Let $ (a_{-\beta}, b_\beta, c_{-\beta}) $ be a $ \beta $-Weyl triple and denote the corresponding Weyl element by $ w_\beta \defl a_{-\beta} b_\beta c_{-\beta} $. Then the following assertions hold for all $ x_\delta \in \rootgr{\delta} $:
	\begin{lemenumerate}
		\item $ \begin{aligned}[t]
			1_G &= \commutator{x_\delta}{c_{-\beta}} \commpart[\big]{\commutator{x_\delta}{b_\beta}}{c_{-\beta}}{-\alpha} \commutator{x_\delta}{a_{-\beta}} \commutator[\big]{\commpart[\big]{\commutator{x_\delta}{a_{-\beta}}}{b_\beta}{\delta}}{c_{-\beta}} \\
			& \qquad \mathord{} \cdot \commpart[\big]{\commpart[\big]{\commutator{x_\delta}{a_{-\beta}}}{b_\beta}{\gamma}}{c_{-\beta}}{-\alpha}.
		\end{aligned} $
		
		\item \label{B:weylbeta-on-delta-nonabelian:rel2}$ \commpart[\big]{\commutator{x_\delta}{a_{-\beta}}}{b_\beta}{\gamma} = \commutator{x_\delta}{b_\beta}^{-1} $.
		
		\item $ x_\delta^{w_\beta} = x_\delta \commpart[\big]{\commutator{x_\delta}{b_\beta}}{c_{-\beta}}{\delta} \commpart[\big]{\commutator{x_\delta}{a_{-\beta}}}{b_\beta}{\delta} \commpart[\big]{\commpart[\big]{\commutator{x_\delta}{a_{-\beta}}}{b_\beta}{\gamma}}{c_{-\beta}}{\delta} $.
	\end{lemenumerate}
\end{lemma}
\begin{proof}
	\begingroup\mywarning{page break equation}
	\allowdisplaybreaks
	As in \cref{A2Weyl:basecomp-cor,B:basecomp-walpha-cox-cor-delta}, we use \thmitemcref{group-rel}{group-rel:conj1} to perform the following computation:
	\begin{align*}
		x_\delta^{w_\beta} &= \brackets[\big]{x_\delta \commutator{x_\delta}{a_{-\beta}}}^{b_\beta c_{-\beta}} \\
		&= \brackets[\big]{x_\delta \commutator{x_\delta}{b_\beta} \commutator{x_\delta}{a_{-\beta}} \commpart[\big]{\commutator{x_\delta}{a_{-\beta}}}{b_\beta}{\delta} \commpart[\big]{\commutator{x_\delta}{a_{-\beta}}}{b_\beta}{\gamma}}^{c_{-\beta}} \\
		&= x_\delta \commutator{x_\delta}{c_{-\beta}} \commutator{x_\delta}{b_\beta} \commutator[\big]{\commutator{x_\delta}{b_\beta}}{c_{-\beta}} \commutator{x_\delta}{a_{-\beta}} \commpart[\big]{\commutator{x_\delta}{a_{-\beta}}}{b_\beta}{\delta} \\
		&\qquad \mathord{} \cdot \commutator[\big]{\commpart[\big]{\commutator{x_\delta}{a_{-\beta}}}{b_\beta}{\delta}}{c_{-\beta}} \commpart[\big]{\commutator{x_\delta}{a_{-\beta}}}{b_\beta}{\gamma} \commutator[\big]{\commpart[\big]{\commutator{x_\delta}{a_{-\beta}}}{b_\beta}{\gamma}}{c_{-\beta}} \\
		&= x_\delta \commutator{x_\delta}{c_{-\beta}} \commutator{x_\delta}{b_\beta} \commpart[\big]{\commutator{x_\delta}{b_\beta}}{c_{-\beta}}{-\alpha} \commpart[\big]{\commutator{x_\delta}{b_\beta}}{c_{-\beta}}{\delta} \commutator{x_\delta}{a_{-\beta}} \commpart[\big]{\commutator{x_\delta}{a_{-\beta}}}{b_\beta}{\delta} \\
		& \qquad \mathord{} \cdot \commutator[\big]{\commpart[\big]{\commutator{x_\delta}{a_{-\beta}}}{b_\beta}{\delta}}{c_{-\beta}} \commpart[\big]{\commutator{x_\delta}{a_{-\beta}}}{b_\beta}{\gamma} \commpart[\big]{\commpart[\big]{\commutator{x_\delta}{a_{-\beta}}}{b_\beta}{\gamma}}{c_{-\beta}}{-\alpha} \\
		& \qquad \mathord{} \cdot \commpart[\big]{\commpart[\big]{\commutator{x_\delta}{a_{-\beta}}}{b_\beta}{\gamma}}{c_{-\beta}}{\delta}
	\end{align*}
	\endgroup
	Note that each of the factors lies in $ \rootgr{-\alpha} $, $ \rootgr{\delta} $ or $ \rootgr{\gamma} $ and that each of these root groups commutes with the other two because $ G $ has crystallgrophic $ B_2 $-commutator relations. (However, we do not have any information on any of the groups $ \rootgr{-\alpha} $, $ \rootgr{\delta} $ or $ \rootgr{\gamma} $ being abelian.) Thus we can write $ x_\delta^{w_\beta} = y_{-\alpha} y_\delta y_\gamma $ where
	\begin{align*}
		y_{-\alpha} &\defl \commutator{x_\delta}{c_{-\beta}} \commpart[\big]{\commutator{x_\delta}{b_\beta}}{c_{-\beta}}{-\alpha} \commutator{x_\delta}{a_{-\beta}} \commutator[\big]{\commpart[\big]{\commutator{x_\delta}{a_{-\beta}}}{b_\beta}{\delta}}{c_{-\beta}} \\
		&\qquad \mathord{} \cdot \commpart[\big]{\commpart[\big]{\commutator{x_\delta}{a_{-\beta}}}{b_\beta}{\gamma}}{c_{-\beta}}{-\alpha}, \\
		y_\delta &\defl x_\delta \commpart[\big]{\commutator{x_\delta}{b_\beta}}{c_{-\beta}}{\delta} \commpart[\big]{\commutator{x_\delta}{a_{-\beta}}}{b_\beta}{\delta} \commpart[\big]{\commpart[\big]{\commutator{x_\delta}{a_{-\beta}}}{b_\beta}{\gamma}}{c_{-\beta}}{\delta}, \\
		y_\gamma &\defl \commutator{x_\delta}{b_\beta} \commpart[\big]{\commutator{x_\delta}{a_{-\beta}}}{b_\beta}{\gamma}.
	\end{align*}
	Now by the injectivity of the product map on $ (-\alpha, \delta, \gamma) $ and because $ x_\delta^{w_\beta} $ lies in $ \rootgr{\delta} $, we have $ y_{-\alpha} = 1_G = y_\gamma $ and $ x_\delta^{w_\beta} = y_\delta $. The assertion follows.
\end{proof}

We know from \cref{B:abelian-note} that all root groups will ultimately be proven to be abelian in the higher-rank setting. With this knowledge, we can simplify the statement of \cref{B:weylbeta-on-delta-nonabelian}.

\begin{lemma}\label{B:weylbeta-on-delta-abelian}
	Let everything be as in \cref{B:weylbeta-on-delta-nonabelian}, and assume additionally that the root groups $ \rootgr{-\alpha}, \rootgr{\delta}, \rootgr{\gamma} $ are abelian. Then the following hold for all $ x_\delta \in \rootgr{\delta} $:
	\begin{lemenumerate}
		\item \label{B:weylbeta-on-delta-abelian:1}$ \commutator{x_\delta}{c_{-\beta}} \commutator{x_\delta}{a_{-\beta}} \commutator[\big]{\commpart[\big]{\commutator{x_\delta}{a_{-\beta}}}{b_\beta}{\delta}}{c_{-\beta}} = 1_G $.
		
		\item \label{B:weylbeta-on-delta-abelian:2}$ \commpart[\big]{\commutator{x_\delta}{a_{-\beta}}}{b_\beta}{\gamma} = \commutator{x_\delta}{b_\beta}^{-1} $.
		
		\item $ x_\delta^{w_\beta} = x_\delta \commpart{\commutator{x_\delta}{a_{-\beta}}}{b_\beta}{\delta} $.
	\end{lemenumerate} 
\end{lemma}
\begin{proof}
	Define $ y_{-\alpha}, y_\delta, y_\gamma $ as in the proof of \cref{B:weylbeta-on-delta-nonabelian}. Then $ \commpart{\commutator{x_\delta}{a_{-\beta}}}{b_\beta}{\gamma} = \commutator{x_\delta}{b_\beta}^{-1} $ by \cref{B:weylbeta-on-delta-nonabelian}. Plugging this relation into the definitions of $ y_{-\alpha} $ and $ y_\delta $, we obtain that
	\begin{align*}
		y_{-\alpha} &= \commutator{x_\delta}{c_{-\beta}} \commpart[\big]{\commutator{x_\delta}{b_\beta}}{c_{-\beta}}{-\alpha} \commutator{x_\delta}{a_{-\beta}} \commutator[\big]{\commpart[\big]{\commutator{x_\delta}{a_{-\beta}}}{b_\beta}{\delta}}{c_{-\beta}} \commpart[\big]{\commutator{x_\delta}{b_\beta}^{-1}}{c_{-\beta}}{-\alpha}, \\
		y_\delta &= x_\delta \commpart[\big]{\commutator{x_\delta}{b_\beta}}{c_{-\beta}}{\delta} \commpart[\big]{\commutator{x_\delta}{a_{-\beta}}}{b_\beta}{\delta} \commpart[\big]{\commutator{x_\delta}{b_\beta}^{-1}}{c_{-\beta}}{\delta}.
	\end{align*}
	Note that
	\[ \commpart{\commutator{x_\delta}{b_\beta}^{-1}}{c_{-\beta}}{\delta} = \commpart{\commutator{x_\delta}{b_\beta}}{c_{-\beta}}{\delta}^{-1} \midand \commpart{\commutator{x_\delta}{b_\beta}^{-1}}{c_{-\beta}}{-\alpha} = \commpart{\commutator{x_\delta}{b_\beta}}{c_{-\beta}}{-\alpha}^{-1} \]
	by \thmitemcref{B:comm-add-cry}{B:comm-add-cry:alpha-beta} and~\thmitemref{B:comm-add-cry}{B:comm-add-cry:alpha-gamma}, respectively. Since $ \rootgr{-\alpha} $ and $ \rootgr{\delta} $ are abelian, it follows that
	\begin{align*}
		y_{-\alpha} = \commutator{x_\delta}{c_{-\beta}} \commutator{x_\delta}{a_{-\beta}} \commutator[\big]{\commpart[\big]{\commutator{x_\delta}{a_{-\beta}}}{b_\beta}{\delta}}{c_{-\beta}} \midand y_\delta = x_\delta \commpart[\big]{\commutator{x_\delta}{a_{-\beta}}}{b_\beta}{\delta}.
	\end{align*}
	Since $ y_{-\alpha} = 1_G $ and $ y_\delta = x_\delta^{w_\beta} $ by \cref{B:weylbeta-on-delta-nonabelian}, this finishes the proof.
\end{proof}

\begin{lemma}\label{B:weylbeta-on-delta-abelian-cor}
	Let $ (a_{-\beta}, b_\beta, c_{-\beta}) $ be a $ \beta $-Weyl triple and denote by $ w_\beta \defl a_{-\beta} b_\beta c_{-\beta} $ the corresponding Weyl element. Assume that the root groups $ \rootgr{-\alpha}, \rootgr{\delta}, \rootgr{\gamma} $ are abelian. Then the following hold:
	\begin{lemenumerate}
		\item \label{B:weylbeta-on-delta-abelian-cor:act}$ x_\delta^{w_\beta} = x_\delta \commpart{\commutator{x_\delta}{a_{-\beta}}}{b_\beta}{\delta} $.
		
		\item $ \commutator{x_\delta}{c_{-\beta}} \commutator{x_\delta}{a_{-\beta}} \commutator[\big]{\commpart[\big]{\commutator{x_\delta}{a_{-\beta}}}{b_\beta}{\delta}}{c_{-\beta}} = 1_G $.
		
		\item $ \commpart[\big]{\commutator{x_\delta}{a_{-\beta}}}{b_\beta}{\gamma} = \commutator{x_\delta}{b_\beta}^{-1} $.
	\end{lemenumerate}
\end{lemma}
\begin{proof}
	This is a consequence of \cref{B:weylbeta-on-delta-abelian}.
\end{proof}


\section{Rank-\texorpdfstring{$ 3 $}{3} Computations}

\label{sec:B-rank3}

\begin{secnotation}\label{B:rank3-notation}
	We denote by $ G $ a group which has crystallographic $ B_n $-commutator relations with root groups $ (\rootgr{\alpha})_{\alpha \in B_n} $ for some fixed integer $ n \ge 3 $. We assume that $ \invset{\alpha} $ is non-empty for all $ \alpha \in \roots $ and that $ G $ is rank-2-injective.
\end{secnotation}

In this section, we prove that $ G $ satisfies the square formula for Weyl elements (\cref{B:square-act:summary}).
Before we begin, we make a brief observation which will simplify many computations.

\begin{proposition}\label{B:abelian}
	All root groups of $ G $ are abelian.
\end{proposition}
\begin{proof}
	Every short root lies in a $ B_2 $-quadruple, so it follows from \cref{B:short-abelian} that all short root groups are abelian. Further, every long root lies in an $ A_2 $-subsystem by \cref{B:rootsys:long-in-A2}, so all short root groups are abelian by \cref{ADE:abelian}.
\end{proof}

\subsection{The Action of Long Weyl Elements on Short Root Groups}

In this case, everything is essentially done by the rank-2 computations, specifically by \cref{B:weylalpha-square-on-beta}. It remains to put everything together.

\begin{lemma}\label{B:square-act:long-on-short}
	Let $ \alpha $ be a long root, let $ \beta $ be a short root and let $ w_\alpha $ be an $ \alpha $-Weyl element. If $ \alpha $ and $ \beta $ are orthogonal, then $ w_\alpha^2 $ acts trivially on $ \rootgr{\beta} $. Otherwise $ w_\alpha^2 $ acts on $ \rootgr{\beta} $ by inversion.
\end{lemma}
\begin{proof}
	If $ \alpha $ and $ \beta $ are orthogonal, then $ \beta $ is adjacent to $ \alpha $ and $ -\alpha $ by \cref{B:rootsys:ortho} because $ \alpha $ is long. This implies that $ w_\alpha^2 $ acts trivially on $ \rootgr{\beta} $. If they are not orthogonal, then they lie in a common $ B_2 $-quadruple, so it follows from \cref{B:weylalpha-square-on-beta} that $ w_\alpha^2 $ acts on $ \rootgr{\beta} $ by inversion.
\end{proof}

\begin{proposition}\label{B:square-act:long-on-short-cartan}
	Let $ \alpha $ be a long root, let $ \beta $ be a short root, let $ w_\alpha $ be an $ \alpha $-Weyl element. Put $ w \defl w_\alpha^2 $ and set $ \epsilon \defl (-1)^{\cartanint{\beta}{\alpha}} $. Then $ x_\beta^{w} = x_\beta^\epsilon $ for all $ x_\beta \in \rootgr{\beta} $.
\end{proposition}
\begin{proof}
	We know from \cref{B:rootsys:cartan-int-parity} that $ \epsilon = 1 $ if $ \alpha, \beta $ are orthogonal and that $ \epsilon = -1 $ if they are not. Thus the assertion follows from \cref{B:square-act:long-on-short}.
\end{proof}

\subsection{The Action of Long Weyl Elements on Long Roots}

In this case, we can reduce everything to the case of gradings of type $ A_2 $ and $ A_1 \times A_1 $.

\begin{proposition}\label{B:square-act:long-on-long:cartan}
	Let $ \alpha, \gamma $ be two long roots and let $ w_\alpha $ be an $ \alpha $-Weyl element. Put $ w \defl w_\alpha^2 $ and set $ \epsilon \defl (-1)^{\cartanint{\gamma}{\alpha}} $. Then $ x_\gamma^{w} = x_\gamma^\epsilon $ for all $ x_\gamma \in \rootgr{\gamma} $.
\end{proposition}
\begin{proof}
	If $ \alpha $ and $ \gamma $ lie in a common $ A_2 $-subsystem, then the assertion holds by \cref{A2Weyl:cartan-comp}. If they do not lie in a common $ A_2 $-subsystem, then they must be orthogonal by \cref{B:rootsys:commonA2} and thus we have $ \epsilon = 1 $. In this case, $ \gamma $ is adjacent to $ \alpha $ and $ -\alpha $ by \cref{B:rootsys:ortho}, and thus $ w_\alpha $ acts trivially on $ \rootgr{\gamma} $. The assertion follows. 
\end{proof}

\begin{note}\label{B:rank2-note}
	The proof of \cref{B:square-act:long-on-long:cartan} fails if $ n=2 $ because \cref{B:rootsys:commonA2} does not hold in this case. In fact, it is not necessarily true that $ w_\alpha^2 $ acts trivially on $ \rootgr{\alpha} $ in a crystallographic $ B_2 $-graded group. We will see a counterexample in \cref{BC:ex-B-counterex-weyl}.
\end{note}

\subsection{The Action of Short Weyl Elements on Long Root Groups}

\label{subsec:B:short-on-long-rank3}

In \thmitemcref{B:basecomp-wdelta-alpha}{B:basecomp-wdelta-alpha:act}, we have already derived a simple formula for the action of short Weyl elements on long root groups. In order to arrive at a formula for the action of squares, we need to understand the behaviour of the map $ \map{}{}{}{(x_\alpha, x_\delta)}{\commpart{x_\alpha}{x_\delta}{\gamma}} $ under inverses in the second component. This is only possible in the rank-3 setting because we have to use \cref{B:square-act:long-on-long:cartan}.

\begin{lemma}\label{B:commpart-inv}
	Let $ (\alpha, \beta, \gamma, \delta) $ be a $ B_2 $-quadruple. Then $ \commpart{x_\alpha}{x_\delta^{-1}}{\gamma} = \commpart{x_\alpha}{x_\delta}{\gamma} $ for all $ x_\alpha \in \rootgr{\alpha} $ and all $ x_\delta \in \rootgr{\delta} $.
\end{lemma}
\begin{proof}
	Let $ x_\alpha \in \rootgr{\alpha} $ and let $ x_\delta \in \rootgr{\delta} $. Further, choose an arbitrary $ \gamma $-Weyl element $ w_\gamma $ and set $ w \defl w_{\gamma}^2 $. Since $ x_\delta^w = x_\delta^{-1} $ by \cref{B:square-act:long-on-short} and $ x_\alpha^w = x_\alpha $ by \cref{B:square-act:long-on-long:cartan}, it follows from \cref{basic:commpart-conj} that
	\begin{align*}
		\commpart{x_\alpha}{x_\delta}{\gamma}^w = \commpart{x_\alpha^w}{x_\delta^w}{\gamma} = \commpart{x_\alpha}{x_\delta^{-1}}{\gamma}.
	\end{align*}
	However, we also know that $ w $ acts trivially on $ \rootgr{\gamma} $ by \cref{B:square-act:long-on-long:cartan}, so $ \commpart{x_\alpha}{x_\delta}{\gamma}^w = \commpart{x_\alpha}{x_\delta}{\gamma} $. The assertion follows.
\end{proof}

\begin{note}
	Thinking in terms of the commutator relations that we will ultimately prove for $ G $ (see \cref{B:ex-commrel}), \cref{B:commpart-inv} corresponds to the fact that $ q(-v) = q(v) $ for a quadratic form $ q $.
\end{note}

\begin{lemma}\label{B:square-act:short-on-long}
	Let $ (\alpha, \beta, \gamma, \delta) $ be a $ B_2 $-quadruple and let $ w_\delta $ be a $ \delta $-Weyl element. Then $ w_\delta^2 $ acts trivially on $ \rootgr{\alpha} $ and on $ \rootgr{\gamma} $.
\end{lemma}
\begin{proof}
	Choose a $ \delta $-Weyl triple $ (a_{-\delta}, b_\delta, c_{-\delta}) $ whose corresponding Weyl element is $ w_\delta $ and let $ x_\alpha \in \rootgr{\alpha} $. By \thmitemcref{basic:weyl-general}{basic:weyl-general:inv}, $ (c_{-\delta}^{-1}, b_\delta^{-1}, a_{-\delta}^{-1}) $ is a $ \delta $-Weyl triple with corresponding Weyl element $ w_\delta^{-1} $. Thus it follows from \thmitemcref{B:basecomp-wdelta-alpha}{B:basecomp-wdelta-alpha:act} that
	\[ x_\alpha^{w_\delta} = \commpart{x_\delta}{b_\alpha}{\gamma} \midand x_\alpha^{w_\delta^{-1}} = \commpart{x_\delta}{b_\alpha^{-1}}{\gamma}. \]
	These terms are equal by \cref{B:commpart-inv}. This implies that $ w_\delta^2 $ acts trivially on $ x_\alpha $ and this on all of $ \rootgr{\alpha} $. Since $ w_\delta $ is also a $ (-\delta) $-Weyl element and $ (\gamma, \beta, \alpha, -\delta) $ is a $ B_2 $-quadruple, it follows that $ w_\delta^2 $ acts trivially on $ \rootgr{\gamma} $ as well.
\end{proof}

\begin{proposition}\label{B:square-act:short-on-long-cartan}
	Let $ \delta $ be a short, let $ \alpha $ be a long root and let $ w_\delta $ be a $ \delta $-Weyl element. Put $ \epsilon \defl (-1)^{\cartanint{\alpha}{\delta}} $ and $ w \defl w_\delta^2 $. Then $ x_\alpha^w = x_\alpha^\epsilon $ for all $ x_\alpha \in \rootgr{\alpha} $. Even more, we always have $ \epsilon = 1 $ and $ w $ acts trivially on all long root groups.
\end{proposition}
\begin{proof}
	By \cref{B:rootsys:cartan-int-longshort}, $ \cartanint{\alpha}{\delta} $ is an even number, so $ \epsilon = 1 $. If $ \alpha $ and $ \delta $ lie in a common $ B_2 $-subsystem, then $ w $ acts trivially on $ \rootgr{\alpha} $ by \cref{B:square-act:short-on-long}. Otherwise they are orthogonal, and thus $ \alpha $ is adjacent to $ \delta $ and $ -\delta $ by \cref{B:rootsys:ortho}. It follows that $ w $ acts trivially on $ \rootgr{\alpha} $ in this case, too.
\end{proof}

\subsection{The Action of Short Weyl Elements on Short Root Groups}

In this section, the rank-3 assumption is used implicitly when we refer to results from subsection~\ref{subsec:B:short-on-long-rank3}.

We already know that every long Weyl element in $ G $ is balanced by the corresponding result for $ A_2 $-graded groups (see \thmitemcref{A2Weyl:weyl}{A2Weyl:weyl:leftright}). For short Weyl elements, this is not yet clear and will only follow from a combination of \cref{B:standard-weyl-balanced,B:thm}. However, the following result is a first step into this direction.

\begin{lemma}\label{B:shortweyl:ac}
	Let $ (\alpha, \beta, \gamma, \delta) $ be a $ B_2 $-quadruple and let $ (a_{-\beta}, b_\beta, c_{-\beta}) $ be a $ \beta $-Weyl triple. Then $ \commutator{x_\delta}{a_{-\beta}} = \commutator{x_\delta}{c_{-\beta}} $ for all $ x_\delta \in \rootgr{\delta} $.
\end{lemma}
\begin{proof}
	By \thmitemcref{basic:weyl-general}{basic:weyl-general:inv}, $ (c_{-\beta}^{-1}, b_\beta^{-1}, a_{-\beta}^{-1}) $ is also a $ \beta $-Weyl triple. Thus it follows from \thmitemcref{B:weylbeta-on-delta-nonabelian}{B:weylbeta-on-delta-nonabelian:rel2} that
	\[ \commpart[\big]{\commutator{x_\delta}{a_{-\beta}}}{b_\beta}{\gamma} = \commutator{x_\delta}{b_\beta}^{-1} = \commutator{b_\beta}{x_\delta} \midand \commpart[\big]{\commutator{x_\delta}{c_{-\beta}^{-1}}}{b_\beta^{-1}}{\gamma} = \commutator{b_\beta^{-1}}{x_\delta}. \]
	Note that by \thmitemcref{B:comm-add-cry}{B:comm-add-cry:alpha-gamma} and \cref{basic:comm-add}, respectively, we have
	\begin{align*}
		\commpart[\big]{\commutator{x_\delta}{c_{-\beta}^{-1}}}{b_\beta^{-1}}{\gamma} = \commpart[\big]{\commutator{x_\delta}{c_{-\beta}}}{b_\beta^{-1}}{\gamma}^{-1} \midand  \commutator{b_\beta^{-1}}{x_\delta} = \commutator{b_\beta}{x_\delta}^{-1}.
	\end{align*}
	Further, $ \commpart{\commutator{x_\delta}{c_{-\beta}}}{b_\beta^{-1}}{\gamma} = \commpart{\commutator{x_\delta}{c_{-\beta}}}{b_\beta}{\gamma} $ by \cref{B:commpart-inv}. Altogether, we conclude that
	\[ \commpart[\big]{\commutator{x_\delta}{a_{-\beta}}}{b_\beta}{\gamma} = \commutator{b_\beta}{x_\delta} = \commutator{b_\beta^{-1}}{x_\delta}^{-1} = \commpart[\big]{\commutator{x_\delta}{c_{-\beta}^{-1}}}{b_\beta^{-1}}{\gamma}^{-1} = \commpart{\commutator{x_\delta}{c_{-\beta}}}{b_\beta}{\gamma}. \]
	Since the map $ \map{}{\rootgr{-\alpha}}{\rootgr{\gamma}}{x_{-\alpha}}{\commpart{x_{-\alpha}}{b_\beta}{\gamma}} $ is an isomorphism by \cref{B:rootisom:alpha-gamma}, we infer that $ \commutator{x_\delta}{a_{-\beta}} = \commutator{x_\delta}{c_{-\beta}} $, as desired.
\end{proof}

\begin{lemma}\label{B:square-act:beta-on-delta}
	Let $ (\alpha, \beta, \gamma, \delta) $ be a $ B_2 $-quadruple and let $ w_\beta $ be a $ \beta $-Weyl element. Then $ w_\beta^2 $ acts trivially on $ \rootgr{\delta} $.
\end{lemma}
\begin{proof}
	Choose a $ \beta $-Weyl triple $ (a_{-\beta}, b_\beta, c_{-\beta}) $ whose corresponding Weyl element is $ w_\beta $. By \thmitemcref{basic:weyl-general}{basic:weyl-general:inv}, $ (c_{-\beta}^{-1}, b_\beta^{-1}, a_{-\beta}^{-1}) $ is a $ (-\beta) $-Weyl triple with corresponding Weyl element $ w_\beta^{-1} $. Applying \cref{B:weylbeta-on-delta-abelian} to these two Weyl triples, we see that
	\[ x_\delta^{w_\beta} = x_\delta \commpart[\big]{\commutator{x_\delta}{a_{-\beta}}}{b_\beta}{\delta} \midand x_\delta^{w_\beta^{-1}} = x_\delta \commpart[\big]{\commutator{x_\delta}{c_{-\beta}^{-1}}}{b_\beta^{-1}}{\delta}. \]
	Using \cref{basic:comm-add} and \thmitemcref{B:comm-add-cry}{B:comm-add-cry:delta-beta}, we can simplify the last term as follows:
	\[ x_\delta^{w_\beta^{-1}} = x_\delta \commpart[\big]{\commutator{x_\delta}{c_{-\beta}}}{b_\beta}{\delta}. \]
	Since $ \commutator{x_\delta}{a_{-\beta}} = \commutator{x_\delta}{c_{-\beta}} $ by \cref{B:shortweyl:ac}, we conclude that $ w_\beta $ and $ w_\beta^{-1} $ act identically on $ \rootgr{\delta} $. In other words, $ w_\beta^2 $ acts trivially on $ \rootgr{\delta} $, which finishes the proof.
\end{proof}

\begin{lemma}\label{B:square-act:beta-on-beta}
	Let $ \beta $ be a short root and let $ w_\beta $ be a $ \beta $-Weyl elements. Then $ w_\beta^2 $ acts trivially on $ \rootgr{\beta} $ and $ \rootgr{-\beta} $.
\end{lemma}
\begin{proof}
	Choose roots $ \alpha, \gamma, \delta $ such that $ (\alpha, \beta, \gamma, \delta) $ is a $ B_2 $-quadruple. Further, choose any $ \alpha $-Weyl element $ w_\alpha $. By the definition of Weyl elements, we have $ \rootgr{\beta} = \rootgr{\delta}^{w_\alpha} $, so $ \rootgr{\beta} $ is contained in the group generated by $ \rootgr{\delta} \union \rootgr{\alpha} \union \rootgr{-\alpha} $. Since $ w_\beta^2 $ acts trivially on $ \rootgr{\delta} $ by \cref{B:square-act:beta-on-delta} and trivially on $ \rootgr{\alpha} $ and $ \rootgr{-\alpha} $ by \cref{B:square-act:short-on-long-cartan}, it follows that $ w_\beta^2 $ acts trivially on $ \rootgr{\beta} $ as well. Since $ w_\beta $ is also a $ (-\beta) $-Weyl element by \thmitemcref{basic:weyl-general}{basic:weyl-general:minus}, this implies that $ w_\beta^2 $ acts trivially on $ \rootgr{-\beta} $ as well.
\end{proof}

\begin{proposition}\label{B:square-act:short-on-short:cartan}
	Let $ \beta, \delta $ be two short roots and let $ w_\beta $ be a $ \beta $-Weyl element. Put $ \epsilon \defl (-1)^{\cartanint{\delta}{\beta}} $ and $ w \defl w_\beta^2 $. Then $ x_\delta^w = x_\delta^\epsilon $ for all $ x_\delta \in \rootgr{\delta} $. Even more, we always have $ \epsilon = 1 $ and $ w $ acts trivially on all short root groups.
\end{proposition}
\begin{proof}
	By \cref{B:rootsys:cartan-int-shortshort}, $ \cartanint{\delta}{\beta} $ is always an even number, so $ \epsilon = 1 $. It is clear that $ \beta $ and $ \delta $ lie in a common $ B_2 $-subsystem. Thus $ w $ acts trivially on $ \rootgr{\delta} $ by \cref{B:square-act:beta-on-delta,B:square-act:beta-on-beta}, which finishes the proof.
\end{proof}

\subsection{Summary}

We can summarise the results of the previous subsections as follows.

\begin{proposition}\label{B:square-act:summary}
	Let $ \xi, \zeta $ be any two roots in $ B_n $ and let $ w_\zeta $ be a $ \zeta $-Weyl element. Put $ w \defl w_\zeta^2 $ and $ \epsilon \defl (-1)^{\cartanint{\xi}{\zeta}} $. Then $ x_\xi^w = x_\xi^\epsilon $ for all $ x_\xi \in \rootgr{\xi} $. In other words, $ G $ satisfies the square formula for Weyl elements (see \cref{param:square-formula-rgg-def}).
\end{proposition}
\begin{proof}
	This is simply a summary of \cref{B:square-act:long-on-short-cartan,B:square-act:short-on-long-cartan,B:square-act:short-on-short:cartan,B:square-act:long-on-long:cartan}.
\end{proof}

The following result is simply a special case of \cref{B:square-act:summary}, but it deserves to be pointed out.

\begin{proposition}\label{B:short-weyl-center}
	Let $ \beta $ be a short root and let $ w_\beta $ be a $ \beta $-Weyl element. Then $ w_\beta^2 $ lies in the center of $ G $.
\end{proposition}
\begin{proof}
	By \cref{B:square-act:short-on-long-cartan,B:square-act:short-on-short:cartan}, $ w_\beta^2 $ acts trivially on all root groups. Since these generate $ G $, it follows that $ w_\beta^2 $ acts trivially on $ G $. In other words, it lies in the center of~$ G $.
\end{proof}

\subsection{Final Observations}

Before we end this section, we make a quick detour to prove \cref{B:short-weyl-same-action}. It will be needed in the proof of stabiliser-compatibility for $ B_n $-graded groups (\cref{B:stab-comp}).

\begin{lemma}\label{B:long-act-identically-rank2lem}
	Let $ (\alpha, \beta, \gamma, \delta) $ be a $ B_2 $-quadruple and choose an $ \alpha $-Weyl element $ w_\alpha $ and a $ \delta $-Weyl element $ w_\delta $. Let $ \rho $ be any short root which does not lie in $ \Set{\pm \beta, \pm \delta} $. Then the Weyl elements $ w_\delta $ and $ w_\delta^{w_\alpha} $ act identically on $ \rootgr{\rho} $.
\end{lemma}
\begin{proof}
	The only short roots which lie in a common $ B_2 $-subsystem with $ \alpha $ are those in $ \Set{\pm \beta, \pm \delta} $, so $ \rootgr{\rho} $ commutes with $ \rootgr{\alpha} $ and $ \rootgr{-\alpha} $. The assertion follows.
\end{proof}

\begin{lemma}\label{B:long-act-identically-rank3lem}
	Let $ \alpha, \beta, \gamma $ be roots which form a root base for a subsystem of $ B_n $ of type $ B_3 $ such that $ (\alpha, \beta) $ is an $ A_2 $-pair and $ (\beta, \gamma) $ is a $ B_2 $-pair. Choose corresponding Weyl elements $ w_\alpha $, $ w_\beta $, $ w_\gamma $. Then the actions of $ w_\gamma $ and $ w_\gamma^{w_\beta w_\alpha} $ on $ \rootgr{\beta} $ are identical.
\end{lemma}
\begin{proof}
	Without loss of generality, we can assume that the $ B_3 $-subsystem is in standard form (as in \cref{B:Bn-standard-rep}). Thus we assume that there exists an orthonormal basis $ (\basvec_1, \basvec_2, \basvec_3) $ of $ \gen{\alpha, \beta, \gamma}_\IR $ which satisfies
	\[ \alpha = \basvec_1 - \basvec_2 , \qquad \beta = \basvec_2 - \basvec_3 \midand \gamma = \basvec_3. \]
	Further, we put
	\begin{align*}
		w_{12} \defl w_\alpha, \quad w_{21} \defl w_{12}^{-1} , \quad w_{23} \defl w_\beta , \quad w_{32} \defl w_{23}^{-1}, \quad w_3 \defl w_\gamma.
	\end{align*}
	For an arbitrary $ x_2 \in \rootgr{\basvec_2} $, we have to show that
	\[ x_2^{w_3} = x_2^{w_{21} w_{32} w_3 w_{23} w_{12}}. \]
	Since $ x_2^{w_{21}} $ lies in $ \rootgr{\basvec_1} $, this element commutes with $ w_{32} $, so
	\[ x_2^{w_{21} w_{32} w_3 w_{23} w_{12}} = x_2^{w_{21} w_3 w_{23} w_{12}}. \]
	Similarly, $ x_2^{w_{21} w_3} $ lies in $ \rootgr{\basvec_1} $ and thus commutes with $ w_{23} $, so
	\[ x_2^{w_{21} w_3 w_{23} w_{12}} = x_2^{w_{21} w_3 w_{12}}. \]
	Finally, $ w_3 $ commutes with $ w_{12} $, so $ w_{21} w_3 w_{12} = w_{21} w_{12} w_3 = w_3 $. This finishes the proof.
\end{proof}

\begin{proposition}\label{B:short-weyl-same-action}
	Write $ B_n $ in standard representation and denote by $ \rootbase $ the standard root base of $ B_n $. Let $ (w_\delta)_{\delta \in \rootbase} $ be a $ \rootbase $-system of Weyl elements in $ G $ and let $ (w_\alpha)_{\alpha \in \Bnsub} $ be its standard $ \Bnsub $-extension (as in \cref{B:Bnsub-ext-def}). Then for all $ i,j,k \in \numint{1}{n} $ such that $ k \nin \Set{i,j} $, the actions of $ w_{\basvec_i} $ and $ w_{\basvec_j} $ on $ \rootgr{\basvec_k} $ and $ \rootgr{-\basvec_k} $ are identical.
\end{proposition}
\begin{proof}
	The assertion is trivial if $ i=j $, so we can assume that $ i,j,k $ are pairwise distinct. For all $ i \in \numint{1}{n-1} $, we put $ w_{i,i+1} \defl w_{\basvec_i - \basvec_{i+1}} $ and $ w_{i+1,i} \defl w_{i,i+1}^{-1} $. Further, we put $ w_i \defl w_{\basvec_i} $ for all $ i \in \numint{1}{n} $. We begin by showing that for all $ i \in \numint{1}{n} $ and all $ k \in \numint{1}{n-1} \setminus \compactSet{i} $, the actions of $ w_i $ and $ w_n $ on $ \rootgr{\basvec_k} $ and $ \rootgr{-\basvec_k} $ are identical. We prove this by induction on $ i $, the case $ i=n $ being trivial. Thus assume that $ i<n $. If $ k \ne i+1 $, then $ \basvec_k $ and $ -\basvec_k $ are not contained in $ \Set{\pm \basvec_i, \pm \basvec_{i+1}} $, so it follows from \cref{B:long-act-identically-rank2lem} that $ w_i = w_{i+1}^{w_{i+1,i}} $ and $ w_{i+1} $ act identically on $ \rootgr{\basvec_k} $ and $ \rootgr{-\basvec_k} $. By the induction hypothesis, it follows that $ w_i $ and $ w_n $ act identically on $ \rootgr{\basvec_k} $ and $ \rootgr{-\basvec_k} $. Now assume that $ k=i+1 $. Then $ i+2 = k+1 \le n $, so
	\[ w_i = w_{i+2}^{w_{i+2,i+1} w_{i+1,i}}. \]
	Thus it follows from \cref{B:long-act-identically-rank3lem} that $ w_i $ and $ w_{i+2} $ act identically on $ \rootgr{\basvec_{i+1}} $ and $ \rootgr{-\basvec_{i+1}} $. Again by the induction hypothesis, we infer that $ w_i $ and $ w_n $ act identically on $ \rootgr{\basvec_k} $ and $ \rootgr{-\basvec_k} $.
	
	Now let $ i,j,k \in \numint{1}{n} $ be pairwise distinct. If $ k \ne n $, then it follows from the conclusion of the previous paragraph that $ w_i $ and $ w_j $ act on $ \rootgr{\basvec_k} $ and $ \rootgr{-\basvec_k} $ in the same way as $ w_n $, so in particular, $ w_i $ and $ w_j $ act identically. Now assume that $ k=n $. Without loss of generality, assume that $ i<j $. Applying \cref{B:long-act-identically-rank2lem} inductively, we can easily see that all the Weyl elements
	\[ w_j, \quad w_j^{w_{j,j-1}}, \quad w_{j}^{w_{j,j-1} w_{j-1,j-2}}, \quad \ldots, \quad w_j^{w_{j,j-1} \cdots w_{i+1,i}} \]
	act identically on $ \rootgr{\basvec_n} $ and $ \rootgr{-\basvec_n} $. Since
	\[ w_i = w_j^{w_{j,j-1} \cdots w_{i+1,i}}, \]
	the assertion follows.
\end{proof}


\section{Standard Signs}

\label{sec:B:stsigns}

\begin{secnotation}
	We fix an integer $ n \ge 3 $ and consider the root system $ B_n $ in its standard representation (as in \cref{B:Bn-standard-rep}). We denote by $ G $ a group with a $ B_n $-pregrading $ (\rootgr{\alpha})_{\alpha \in B_n} $, by $ \comring $ a commutative associative ring and by $ (\module, q) $ a quadratic module over $ \comring $. We assume that there exists a coordinatisation $ (\risom{\alpha})_{\alpha \in B_n} $ of $ G $ by $ (\module, q) $ with standard signs (in the sense of the following \cref{B:standard-param-def}), and we fix this coordinatisation.
\end{secnotation}

\begin{definition}[Coordinatisation with standard signs]\label{B:standard-param-def}
	Let $ G $ be a group with a $ B_n $-pregrading $ (\rootgr{\alpha})_{\alpha \in B_n} $, let $ \comring $ be a commutative associate ring and let $ (\module, q) $ be a quadratic module over $ \comring $. A \defemph*{coordinatisation of $ G $ by $ (\module, q) $ with standard signs}\index{coordinatisation of a root graded group!with standard signs!type Bn@type $ B_n $} is a family $ (\risom{\alpha})_{\alpha \in B_n} $ with the following properties:
	\begin{stenumerate}
		\item For all roots $ \alpha $, the map $ \risom{\alpha} $ is an isomorphism from $ (\comring, +) $ to $ \rootgr{\alpha} $ if $ \alpha $ is long and it is an isomorphism from $ (\module, +) $ to $ \rootgr{\alpha} $ if $ \alpha $ is short.
		
		\item The same commutator relations as in \cref{B:ex-commrel} are satisfied. We will refer to them as the \defemph*{standard commutator relations}.
	\end{stenumerate}
\end{definition}

We will show in \cref{B:thm} that every $ B_n $-graded group is coordinatised by some quadratic module with standard signs.

In the remaining part of this section, we investigate some properties of groups with a standard coordinatisation. We begin by showing that all Weyl elements in such groups are balanced, which by \cref{B:thm} shows that Weyl elements in all $ B_n $-graded groups are balanced.

\begin{lemma}\label{B:short-weakly-balanced}
	Assume that $ G $ is rank-2-injective. Let $ \delta $ be a short root and let $ w_\delta $ be a $ \delta $-Weyl element in $ G $. If $ v',v,v'' \in \module $ are such that $ w_\delta = \risom{-\delta}(v') \risom{\delta}(v) \risom{-\delta}(v'') $, then $ q(v) $ is invertible and $ v' = v'' = -q(v)^{-1} v $. In particular, $ w_\delta $ is weakly balanced.
\end{lemma}
\begin{proof}
	We prove this for the short root $ \delta = \basvec_2 $. The assertion for all other short roots can be proven similarly. Let $ w_2 $ be a $ \basvec_2 $-Weyl element and let $ v', v, v'' \in \module $ be such that $ w_2 = \risshneg{2}(v') \risshpos{2}(v) \risshneg{2}(v'') $. For all $ a \in \comring $, we 
	know from \thmitemcref{B:basecomp-wdelta-alpha}{B:basecomp-wdelta-alpha:act}, \thmitemcref{B:basecomp-wdelta-gamma}{B:basecomp-wdelta-gamma:act} and the standard commutator relations that
	\begin{align*}
		\rismin{1}{2}(a)^{\hat{w}_2} &= \commpart{\rismin{1}{2}(a)}{\risshpos{2}(v)}{\basvec_1 + \basvec_2} = \risplus{1}{2}\brackets[\big]{-a q(v)} \rightand \\
		\risplus{1}{2}(a)^{\hat{w}_2} &= \commpart{\risplus{1}{2}(a)}{\risshneg{2}(v')}{\basvec_1 - \basvec_2} = \rismin{1}{2}\brackets[\big]{-a q(v')}.
	\end{align*}
	This implies that $ \rismin{1}{2}(1_\comring)^{w_2^2} = (q(v) q(v')) $. On the other hand, we know that $ \rismin{1}{2}(1_\comring)^{w_2^2} = \rismin{1}{2}(1_\comring) $ by \cref{B:square-act:short-on-long-cartan}, so it follows that $ q(v) q(v') = 1_\comring $. In particular, $ q(v) $ is invertible. Further, it follows from \thmitemcref{B:basecomp-wdelta-alpha}{B:basecomp-wdelta-alpha:rel} that
	\begin{align*}
		1_G &= \commpart{\rismin{1}{2}(1)}{\risshpos{2}(v)}{\basvec_1} \commpart[\big]{\commpart{\rismin{1}{2}(1)}{\risshpos{2}(v)}{\basvec_1 + \basvec_2}}{\risshneg{2}(v')}{\basvec_1}  \\
		&= \risshpos{1}(v) \commpart[\big]{\risplus{1}{2}\brackets[\big]{-q(v)}}{\risshneg{2}(v')}{\basvec_1} = \risshpos{1}\brackets[\big]{v + q(v)v'}.
	\end{align*}
	It follows that $ v = -q(v) v' $. Since $ q(v) $ is invertible, this implies that $ v' = -q(v)^{-1} v $. The same argument can be applied to the Weyl element $ w_2^{-1} = \risshneg{2}(-v'') \risshpos{2}(-v) \risshneg{2}(-v') $, which then yields that $ -v'' = -q(-v)^{-1} (-v) $. Thus we have $ v'' = -q(v)^{-1} v = v' $ as well.
\end{proof}

\begin{note}
	The element $ v \in \module $ in \cref{B:short-weakly-balanced} need not be uniquely determined by $ w_\delta $: See \cref{B:ex:short-weyl-matrix}.
\end{note}

\begin{proposition}\label{B:standard-weyl-balanced}
	If $ G $ is rank-2-injective, then all Weyl elements in $ G $ are balanced.
\end{proposition}
\begin{proof}
	For long Weyl elements, we already know this from the corresponding assertion for $ A_2 $-graded groups, see \thmitemcref{A2Weyl:weyl}{A2Weyl:weyl:leftright}. For short Weyl elements, this follows from \cref{B:short-weakly-balanced,basic:all-balanced}.
\end{proof}

\begin{proposition}\label{B:stsign:weyl-char}
	Assume that $ G $ is rank-2-injective. Then the following hold:
	\begin{proenumerate}
		\item \label{B:stsign:weyl-char:long}Let $ \alpha $ be a long root. Define
		\[ w_\alpha(r) \defl \risom{-\alpha}(-r^{-1}) \risom{\alpha}(r) \risom{-\alpha}(-r^{-1}) \]
		for all invertible $ r \in \comring $. Then the maps
		\[ \map{}{\ringinvset{\comring}}{\invset{\alpha}}{r}{\risom{\alpha}(r)} \midand \map{}{\ringinvset{\comring}}{\weylset{\alpha}}{r}{w_\alpha(r)} \]
		are well-defined bijections. Here $ \ringinvset{\comring} $, $ \invset{\alpha} $ and $ \weylset{\alpha} $ denote the sets of invertible elements in $ \comring $ (in the sense of \cref{ring:def-invertible}), the set of $ \alpha $-invertible elements in $ \rootgr{\alpha} $ and the set of $ \alpha $-Weyl elements, respectively.
		
		\item \label{B:stsign:weyl-char:short}Let $ \delta $ be a short root and put $ \ringinvset{\module} \defl \Set{u \in \module \given q(u) \in \ringinvset{\comring}} $. Define
		\[ w_\delta(u) \defl \risom{-\delta}\brackets[\big]{-q(u)^{-1} u} \risom{\delta}(u) \risom{-\delta}\brackets[\big]{-q(u)^{-1} u} \]
		for all $ u \in \ringinvset{\module} $. Then the maps
		\[ \map{}{\ringinvset{\module}}{\invset{\delta}}{u}{\risom{\delta}(u)} \midand \map{}{\ringinvset{\module}}{\weylset{\delta}}{u}{w_\delta(u)} \]
		are well-defined surjections, and the first one is a bijection.
		
		\item \label{B:stsign:weyl-char:formulas}The Weyl elements defined above satisfy the same conjugation formulas as in \cref{B:ex-weylformula-long,B:ex-weylformula-short}.
		
		\item \label{B:stsign:weyl-char:param}Let $ \rootbase $ be any root base of $ B_n $ and choose an element $ v_0 \in \module $ with $ q(v_0) = 1_\comring $. Define $ \twistgroup $, $ \invogroup $, $ \inverparsym $ and $ \invoparsym $ as in \cref{B:ex-twistgroups-def}. Put $ w_\alpha \defl w_\alpha(1_\comring) $ for all long roots in $ \rootbase $ and $ w_\beta \defl w_\beta(v_0) $ for the unique short root in $ \rootbase $. Then $ G $ is parametrised by $ (\twistgroup \times \invogroup, \module, \comring) $ with respect to $ \inverparsym \times \invoparsym $ and $ (w_\delta)_{\delta \in \rootbase} $.
	\end{proenumerate}
\end{proposition}
\begin{proof}
	By \thmitemcref{ADE:stsign:weyl-char}{ADE:stsign:weyl-char:bij}, we already know that the maps in~\ref{B:stsign:weyl-char:long} are bijective if they are well-defined. Using computations as in \cref{B:basecomp-walpha-cor:delta,B:basecomp-walpha-cor:beta}, one can show that elements of the form $ w_\alpha = w_\alpha(r) $ satisfy $ \rootgr{\beta}^{w_\alpha} = \rootgr{\refl{\alpha}(\beta)} $ for all short roots $ \beta $. This proves~\itemref{B:stsign:weyl-char:long}.
	
	We know from \cref{B:short-weakly-balanced} that the maps in~\itemref{B:stsign:weyl-char:short} are surjective if they are well-defined. Further, the first map is injective because $ \risom{\delta} $ is injective. Now let $ u \in \ringinvset{\module} $ and put $ w_\delta \defl w_\delta(u) $. Using the rank-2 computations (see the proofs of \cref{B:basecomp-wdelta-alpha,B:weylbeta-on-delta-nonabelian}) and the standard commutator relations, one can show that $ \rootgr{\rho}^{w_\delta} = \rootgr{\refl{\delta}(\rho)} $ for all $ \rho \in B_n \setminus \Set{\pm \delta} $. Even more, these computations show that the same conjugation formulas as in \cref{B:ex-weylformula-long,B:ex-weylformula-short} are satisfied. It remains to prove the assertion for $ \rho \in \Set{\pm \delta} $. We do this for $ \delta = \basvec_2 = \rho $ and the remaining cases can be covered in a similar way. Let $ v \in \module $ be arbitrary. By the standard commutator relations, we have
	\begin{align*}
		\risshpos{2}(v) &= \commpart{\rismin{2}{1}(1_\comring)}{\risshpos{1}(v)}{\basvec_2} = \commutator{\rismin{2}{1}(1_\comring)}{\risshpos{1}(v)} \commpart{\rismin{2}{1}(1_\comring)}{\risshpos{1}(v)}{\basvec_1 + \basvec_2}^{-1} \\
		&= \commutator{\rismin{2}{1}(1_\comring)}{\risshpos{1}(v)} \risplus{1}{2}\brackets[\big]{-q(v)}.
	\end{align*}
	Using that the conjugation formulas in \cref{B:ex-weylformula-long} are satisfied, we infer that
	\begin{align*}
		\risshpos{2}(v)^{w_2(u)} &= \commutator{\rismin{2}{1}(1_\comring)^{w_2(u)}}{\risshpos{1}(v)^{w_2(u)}} \risplus{1}{2}\brackets[\big]{-q(v)}^{w_2(u)} \\
		&= \commutator[\big]{\risminmin{1}{2}\brackets[\big]{-q(u)^{-1}}}{\risshpos{1}\brackets[\big]{\refl{u}(v)}} \rismin{1}{2}\brackets[\big]{q(u)^{-1} q(v)}
	\end{align*}
	where
	\begin{align*}
		\commutator[\big]{\risminmin{1}{2}\brackets[\big]{-q(u)^{-1}}}{\risshpos{1}\brackets[\big]{\refl{u}(v)}} &= \risshneg{2}\brackets[\big]{q(u)^{-1} \refl{u}(v)} \rismin{1}{2}\brackets[\big]{-q(u)^{-1} q\brackets[\big]{\refl{u}(v)}} \\
		&= \risshneg{2}\brackets[\big]{q(u)^{-1} \refl{u}(v)} \rismin{1}{2}\brackets[\big]{-q(u)^{-1} q\brackets{v}}.
	\end{align*}
	It follows that
	\begin{align*}
		\risshpos{2}(v)^{w_2(u)} &= \risshneg{1}\brackets[\big]{q(u)^{-1} \refl{u}(v)},
	\end{align*}
	which is precisely the desired conjugation formula in \thmitemcref{B:ex-weylformula-short}{B:ex-weylformula-short:on-self} for this choice of roots. This finishes the proof of~\itemref{B:stsign:weyl-char:short}.
	
	Assertions~\itemref{B:stsign:weyl-char:formulas} and~\itemref{B:stsign:weyl-char:param} follow from similar computations as above, or have already been proven.
\end{proof}

\begin{note}
	In the proof of \thmitemcref{B:stsign:weyl-char}{B:stsign:weyl-char:short}, it is not valid to argue that
	\begin{align*}
		\risshpos{2}(v)^{w_2(u)} &= \commpart{\rismin{2}{1}(1_\comring)}{\risshpos{1}(v)}{\basvec_2}^{w_2(u)} = \commpart{\rismin{2}{1}(1_\comring)^{w_2(u)}}{\risshpos{1}(v)^{w_2(u)}}{-\basvec_2} \\
		&= \commpart[\big]{\risminmin{1}{2}\brackets[\big]{-q(u)^{-1}}}{\risshpos{1}\brackets[\big]{\refl{u}(v)}}{-\basvec_2} = \risshneg{2}\brackets[\big]{q(u)^{-1} \refl{u}(v)}
	\end{align*}
	because we would have to use \cref{basic:commpart-conj} for the second equality. This would not be correct because $ w_2(u) $ is not yet proven to be a Weyl element (and also because $ G $ is not assumed to be rank-2-injective). Instead, we have to use the workaround outlined above.
\end{note}

\begin{remark}[compare~\ref{ADE:stsign:weyl-char:rem}]\label{B:weyl-char:rem}
	If $ G $ is not assumed to be rank-2-injective, it is no longer clear that the maps in \thmitemcref{B:stsign:weyl-char}{B:stsign:weyl-char:short} and~\thmitemref{B:stsign:weyl-char}{B:stsign:weyl-char:short} are bijective, but they are still well-defined. To prove this, we have to perform the rank-2-computations in the proofs of \cref{A2Weyl:basecomp-cor,B:basecomp-walpha-cox-cor-delta} for elements of the form $ w_\alpha(r) $ and $ w_\delta(u) $ to show that these are Weyl elements. The computation in \cref{B:basecomp-walpha-cox-cor-delta} is problematic at first glance because it involves the commutator maps from \cref{basic:commpart-def}, which are not well-defined without the assumption that $ G $ is rank-2-injective. However, we can simply interpret any term of the form, say, $ \commpart{\risom{\alpha}(r)}{\risom{\delta}(v)}{\beta} $ to be the term given by the standard commutator relations, and then we obtain the desired result.
\end{remark}

\begin{remark}
	It follows from \cref{B:stsign:weyl-char} that $ G $ satisfies $ \invset{\alpha} = \rootgr{\alpha} \setminus \compactSet{1_G} $ (the additional condition of being an RGD-system) if and only if $ \comring $ is a field and $ (\module, q) $ is anisotropic.
\end{remark}

\begin{remark}\label{universal:B}
	In a similar way as in \cref{universal:ADE}, we can construct a functor $ \hat{G} $ which assigns a $ B_n $-graded group to each pair $ (\module, q, \comring) $ consisting of a commutative associative ring $ \comring $, a $ \comring $-module $ \module $ and a quadratic form $ \map{q}{\module}{\comring}{}{} $. We call $ \hat{G}(\module, q, \comring) $ the \defemph*{universal $ B_n $-graded group over $ (\module, q, \comring) $}\index{root graded group!universal}. By \cref{B:thm}, every $ B_n $-graded group is the quotient of one such group.
\end{remark}

The following result is a special case of \thmitemcref{B:stsign:weyl-char}{B:stsign:weyl-char:short}. We prove it explicitly because it will be used in \cref{B:blue:short-weyl}.

\begin{lemma}\label{B:short-act-refl}
	Let $ \beta $ be a short root, let $ w_\beta $ be a $ \beta $-Weyl element and let $ v \in \module $ such that $ w_\beta = \risom{-\beta}\brackets[\big]{-q(v)^{-1} v} \risom{\beta}(v) \risom{-\beta}\brackets[\big]{-q(v)^{-1} v} $. Then for all short roots $ \delta $ which are orthogonal to $ \beta $, we have
	\[ \risom{\delta}(u)^{w_\beta} = \risom{\delta}(\refl{v}(u)) \]
	for all $ u \in \module $ where $ \map{\refl{v}}{\module}{\module}{u}{u - q(v)^{-1} f(v,u) v} $ denotes the reflection associated to $ v $ (as in \cref{quadmod:refl-def}).
\end{lemma}
\begin{proof}
	Let $ u \in \module $. By \cref{B:rootsys:short-ortho}, there exist unique long roots $ \alpha, \gamma $ such that $ (\alpha, \beta, \gamma, \delta) $ is a $ B_2 $-quadruple. Then by \thmitemcref{B:weylbeta-on-delta-abelian-cor}{B:weylbeta-on-delta-abelian-cor:act}, we have
	\[ \risom{\delta}(u)^{w_\beta} = \risom{\delta}(u) \commpart[\big]{\commutator[\big]{\risom{\delta}(u)}{\risom{-\beta}\brackets[\big]{-q(v)^{-1} v}}}{\risom{\beta}(v)}{\delta}. \]
	Choose $ \epsilon_1, \epsilon_2 \in \compactSet{\pm 1} $ and $ i,j \in \numint{1}{n} $ such that $ \delta = \epsilon_1 \basvec_i $ and $ \beta = \epsilon_2 \basvec_j $. Using the standard commutator relations, one can show that
	\[ \commpart[\big]{\commutator[\big]{\risom{\delta}(u)}{\risom{-\beta}\brackets[\big]{-q(v)^{-1} v}}}{\risom{\beta}(v)}{\delta} = \risom{\delta}\brackets[\big]{-q(v)^{-1} f(u,v)v} \]
	in every possible case. See \cref{fig:B:short-act-refl} for a detailed computation. We conclude that
	\[ \risom{\delta}(u)^{w_\beta} = \risom{\delta}\brackets[\big]{u - q(v)^{-1} f(u,v)v} = \risom{\delta}\brackets[\big]{\refl{v}(u)}, \]
	as desired.
\end{proof}

\begin{figure}[htb]
	\centering\begin{tabular}{cccc}
		\toprule
		$ \epsilon_1 $ & $ \epsilon_2 $ & $ x \defl \commutator[\big]{\risom{\epsilon_1 \basvec_i}(u)}{\risom{-\epsilon_2 \basvec_j}\brackets[\big]{-q(v)^{-1} v}} $ & $ \commpart[\big]{x}{\risom{\epsilon_2 \basvec_j}(v)}{\epsilon_1 \basvec_i} $ \\
		\midrule
		$ - $ & $ + $ & $ \risminmin{i}{j}\brackets[\big]{-\delmin{i>j} q(v)^{-1} f(u,v)} $ & $ \risshneg{i}\brackets[\big]{-\delmin{i>j} \delmin{i>j} q(v)^{-1} f(u,v) v} $ \\
		$ + $ & $ + $ & $ \rismin{i}{j}\brackets[\big]{-q(v)^{-1} f(u,v)} $ & $ \risshpos{i}\brackets[\big]{-q(v)^{-1} f(u,v)v} $ \\
		$ + $ & $ - $ & $ \risplus{i}{j}\brackets[\big]{-\delmin{i<j} q(v)^{-1} f(u,v)} $ & $ \risshpos{i}\brackets[\big]{- \delmin{i<j} \delmin{i<j} q(v)^{-1} f(u,v) v} $ \\
		$ - $ & $ - $ & $ \rismin{j}{i}\brackets[\big]{q(v)^{-1} f(v,u)} $ & $ \risshneg{i}\brackets[\big]{-q(v)^{-1} f(v,u)v} $ \\
		\bottomrule
	\end{tabular}
	\caption{The computation in \cref{B:short-act-refl}.}
	\label{fig:B:short-act-refl}
\end{figure}


\section{Admissible and Standard Partial Twisting Systems}

\label{sec:B:sttwist}

\begin{secnotation}\label{B:secnot:sttwist}
	We fix an integer $ n \ge 3 $ and consider the root system $ B_n $ in its standard representation (as in \cref{B:Bn-standard-rep}) with its standard root base $ \rootbase $.
	We choose $ \comring \defl \IC $, $ \module \defl \IC^2 $, $ \map{q}{\module}{\IC}{}{} $ and $ v_0 \defl (1,0) $ as in \cref{B:ex-param-choice}.
	We define the $ B_n $-graded group $ \EO(q) $ with root groups $ (\rootgr{\alpha})_{\alpha \in B_n} $ as in \cref{B:ex:EOq-def}. We denote the root isomorphisms from \cref{B:ex-roothom-def-long,B:ex-roothom-def-short} by $ (\risom{\alpha})_{\alpha \in B_n} $ and the standard system of Weyl elements from \cref{B:ex-standard-weyl} by $ (w_\delta)_{\delta \in \rootbase} $. Further, we denote by $ (\twistgroup, \inverparsym, \invogroup, \invoparsym) $ the standard partial twisting system of type $ B_n $ in the sense of the following \cref{B:standard-partwist-def}.
\end{secnotation}

In this section, we introduce the standard partial twisting system for any $ B_n $-graded group. Its definition is motivated by the twisting structure of the elementary orthogonal group in \cref{sec:B-example}. We will show that it has some desirable properties which serve as the axioms of admissible partial twisting systems. The goal of the following \cref{sec:B-param} is that any admissible partial twisting system satisfies the conditions of the parametrisation theorem. In the final sections of this chapter, we will then use the standard partial twisting system to construct a parametrisation of any $ B_n $-graded group with standard signs.

\begin{definition}[Standard partial twisting system]\label{B:standard-partwist-def}
	The \defemph*{standard partial twisting system of type $ B_n $ (with respect to $ \rootbase $)}\index{twisting system!partial!standard (type Bn)@standard (type $ B_n $)} is the tuple $ (\twistgroup, \inverparsym, \invogroup, \invoparsym) $ where $ \twistgroup \defl \invogroup \defl \compactSet{\pm 1} $ and where $ \inverparsym $, $ \invoparsym $ are the $ \rootbase $-parity maps from \cref{B:ex-twistgroups-def}. If $ G $ is a group with a $ B_n $-pregrading, then the \defemph*{standard partial twisting system for $ G $ (with respect to $ \rootbase $)} is the same tuple together with the additional information that $ \twistgroup $ acts on all root groups of $ G $ by inversion.
\end{definition}

\begin{note}\label{B:standard-partwist-note}
	Technically, the standard partial twisting system of type $ B_n $ is not a partial twisting system because a partial twisting system has to be defined with respect to a pair $ (G, (w_\delta)_{\delta \in \rootbase}) $ where $ G $ is a group with a $ B_n $-pregrading and $ (w_\delta)_{\delta \in \rootbase} $ is a $ \rootbase $-system of Weyl elements. In particular, the group $ \twistgroup $ has to be equipped with actions on the root groups of $ G $. For this reason, we distinguish between the \enquote{abstract} standard partial twisting system of type $ B_n $ and the standard partial twisting system for a specific group $ G $. Note further that it is not obvious that the standard partial twisting system for a $ B_n $-graded group satisfies all the axioms of a partial twisting system in \cref{param:partwist-def}. We will verify this throughout this section.
\end{note}

\begin{remark}\label{B:standard-param-note}
	We can also regard the groups $ \twistgroup $ and $ \invogroup $ in the standard partial twisting system $ (\twistgroup, \inverparsym, \invogroup, \invoparsym) $ as part of the standard parameter system $ (\twistgroup \times \invogroup, \module, \comring) $ from \cref{quadmod:standard-param}. Note that in the former interpretation, $ \twistgroup $ acts on the root groups and $ \invogroup $ is equipped with no action while in the latter interpretation, $ \twistgroup \times \invogroup $ acts on $ \module $ and $ \comring $.
\end{remark}

\begin{reminder}\label{B:parmap:reminder}
	Recall from \cref{B:ex-is-rgg} that $ (\rootgr{\alpha})_{\alpha \in B_n} $ is a crystallographic $ B_n $-grading of $ \EO(q) $. Thus we can apply all rank-2 and rank-3 computations from the previous sections to $ \EO(q) $. Further, we know from \cref{B:ex-param-thm} that $ \EO(q) $ is parametrised by the standard parameter system $ (\twistgroup \times \invogroup, \module, \comring) $ with respect to $ \inverparsym \times \invoparsym $ and $ (w_\delta)_{\delta \in \rootbase} $ and from \cref{B:ex-param-choice} that $ (\twistgroup \times \invogroup, \module, \comring) $ is $ (\inverparsym \times \invoparsym) $-faithful.
\end{reminder}

We now prove some basic properties of the parity maps $ \inverparsym $ and $ \invoparsym $. Since all values of $ \inverparsym $ and $ \invoparsym $ are given explicitly in \cref{B:ex-parmap-def}, all these properties could be proven with a straightforward but lengthy computation. However, we will see that they are much easier to derive by performing certain computations in the group $ \EO(q) $. 

\begin{lemma}\label{B:parmap-properties}
	$ \inverparsym $ is braid-invariant and adjacency-trivial and $ \invoparsym $ is Weyl-invariant and adjacency-trivial.
\end{lemma}
\begin{proof}
	 By \cref{braid:all}, the system $ (w_\delta)_{\delta \in \rootbase} $ satisfies the braid relations. Thus it follows from \cref{param:param-parmap-has-properties2} that $ \inverparsym \times \invoparsym $ is braid-invariant and adjacency-trivial, which by \thmitemcref{param:prod-eq-lem}{param:prod-eq-lem:eq} implies that $ \inverparsym $ and $ \invoparsym $ have the same properties. Further, an inspection of \cref{B:ex-parmap-def} shows that $ \invoparsym $ is even square-invariant. Hence $ \mu $ is Weyl-invariant, which finishes the proof.
\end{proof}

\begin{lemma}\label{B:invoparmap-semicomp}
	$ \invoparsym $ is semi-complete.
\end{lemma}
\begin{proof}
	Note that the only subgroups of $ \invogroup $ are $ \compactSet{1} $ and $ \invogroup $, both of which have a complement. Thus every parity map with values in $ \invogroup $ is semi-complete, which makes the assertion trivial.
\end{proof}

\begin{lemma}\label{B:parmap-trans-invar}
	$ \inverparsym \times \invoparsym $ is transporter-invariant and $ \inverparsym $, $ \invoparsym $ are independent.
\end{lemma}
\begin{proof}
	At first, we consider the orbit of short roots. Put $ \hat{\alpha} \defl \basvec_1 $. Then
	\begin{align*}
		\inverinvopar{\hat{\alpha}}{(\basvec_1 - \basvec_2, \basvec_1 - \basvec_2)} &= \inverinvopar{\basvec_1}{\basvec_1 - \basvec_2} \inverinvopar{\basvec_2}{\basvec_1 - \basvec_2} = (1_\twistgroup, 1_\twistgroup) (-1_\twistgroup, 1_\invogroup) \\
		&= (-1_\twistgroup, 1_\invogroup) \rightand \\
		\inverinvopar{\hat{\alpha}}{\basvec_2} &= (1_\twistgroup, -1_\invogroup).
	\end{align*}
	Since these elements generate $ \twistgroup \times \invogroup $, it follows that $ \parmoveset{(\twistgroup \times \invogroup)}{\hat{\alpha}}{\hat{\alpha}} = \twistgroup \times \invogroup $. By \cref{parmap:transport-invar-crit}, this implies that $ \inverparsym \times \invoparsym $ is transporter-invariant on the orbit of short roots.
	
	Now we consider the orbit of long roots. By an inspection of \cref{B:ex-parmap-def}, it is clear that
	\[ \parmoveset{(\twistgroup \times \invogroup)}{\alpha}{\beta} \subs \twistgroup \times \compactSet{1_\invogroup} \]
	for all long roots $ \alpha $, $ \beta $. Further, we have
	\[ \parmoveset{(\twistgroup \times \invogroup)}{\basvec_1 - \basvec_2}{\basvec_2 - \basvec_1} = \twistgroup \times \compactSet{1_\invogroup} \]
	because $ \inverinvopar{\basvec_1 - \basvec_2}{\basvec_1 - \basvec_2} = (-1_\twistgroup, 1_\invogroup) $. By criterion~\thmitemref{parmap:transport-invar-char}{parmap:transport-invar-char:bound}, we infer that $ \inverparsym \times \invoparsym $ is transporter-invariant on the orbit of long roots as well.
	
	Finally, the previous computations together with \cref{param:transporter-proj} show that $ \inverparsym $ and $ \invoparsym $ are independent in the sense of \cref{param:parmap-indep-def}.
\end{proof}

\begin{lemma}\label{B:parmap-square}
	For all roots $ \alpha \in \roots $ and $ \delta \in \rootbase $, we have $ \inverpar{\alpha}{\delta \delta} = (-1_\twistgroup)^{\cartanint{\alpha}{\delta}} $ where $ \cartanint{\alpha}{\delta} $ denotes the Cartan integer for $ (\alpha, \delta) $. In other words, $ \inverparsym $ satisfies the square formula (in the sense of \thmitemcref{parmap:prop-def}{parmap:prop-def:square-formula}).
\end{lemma}
\begin{proof}
	It suffices to verify that the conditions of \thmitemcref{param:param-parmap-has-properties2}{param:param-parmap-has-properties2:square} are satisfied. We have proven in \cref{B:square-act:summary} that $ G $ satisfies the square formula for Weyl elements. Further, by definition, $ \twistgroup = \compactSet{\pm 1} $ and $ -1_\twistgroup $ acts on both $ \comring $ and $ \module $ by inversion. Finally, by the transporter-invariance of $ \inverparsym $ and an inspection of \cref{B:ex-parmap-def}, we have $ \parmoveset{\twistgroup}{\alpha}{\alpha} = \twistgroup $ for all roots $ \alpha $. This finishes the proof.
\end{proof}

\begin{lemma}\label{B:short-weyl-square-inverpar}
	Let $ i,j \in \numint{1}{n} $ and let $ k \in \numint{1}{n} \setminus \compactSet{i,j} $. Then the words $ \word{\rho}^i $, $ \word{\rho}^j $ from \cref{B:Bnsub-ext-word} satisfy $ \inverinvopar{\epsilon \basvec_k}{\word{\rho}^i} = \inverinvopar{\epsilon \basvec_k}{\word{\rho}^j} $ for all $ \epsilon \in \compactSet{\pm 1} $.
\end{lemma}
\begin{proof}
	We know from \cref{B:short-weyl-same-action} that $ w_{\word{\rho}^i} $ and $ w_{\word{\rho}^j} $ act identically on $ \rootgr{\epsilon \basvec_k} $. Since the action of $ (\twistgroup \times \invogroup, \module, \comring) $ is $ (\inverparsym \times \invoparsym) $-faithful by \cref{B:ex-param-choice}, the assertion follows.
\end{proof}

\begin{lemma}\label{B:invo-parmap-short}
	If $ \rho, \zeta $ are short roots in $ B_n $, then $ \invopar{\rho}{\reflbr{\zeta}} = -1_{\invogroup} $.
\end{lemma}
\begin{proof}
	For the simple root $ \zeta = \basvec_n $, we can immediately read this off from \cref{B:ex-parmap-def}. Now let $ \zeta = \basvec_m $ for some $ m \in \numint{1}{n-1} $ and put $ \delta_i \defl \basvec_i - \basvec_{i+1} $ for all $ i \in \numint{1}{n-1} $. Then $ \zeta = \basvec_n^{\reflbr{\delta_{n-1} \cdots \delta_{m}}} $, so
	\begin{align*}
		\invopar{\rho}{\reflbr{\pm \zeta}} &= \invopar{\rho}{\delta_m \cdots \delta_{n-1} \basvec_n \delta_{n-1} \cdots \delta_m} = \invopar{\rho}{\delta_m \cdots \delta_{n-1}} \invopar{\rho^{\reflbr{\delta_m \cdots \delta_{n-1}}}}{\basvec_n} \invopar{\rho^{\reflbr{\delta_m \cdots \delta_{n-1} \basvec_n}}}{\delta_{n-1} \cdots \delta_m} \\
		&= 1_\invogroup (-1_\invogroup) 1_\invogroup = -1_\invogroup.
	\end{align*}
	The assertion follows.
\end{proof}

\begin{note}
	\cref{B:invo-parmap-short} also follows from \thmitemcref{B:ex-weylformula-short}{B:ex-weylformula-short:on-self} and~\thmitemref{B:ex-weylformula-short}{B:ex-weylformula-short:on-short} because the parameter system is $ (\inverparsym \times \invoparsym) $-faithful.
\end{note}

\begin{definition}[Admissible partial twisting system]\label{B:admissible-parmap-def}
	Let $ G $ be a group with a $ B_n $-pregrading $ (\rootgr{\alpha})_{\alpha \in B_n} $ and let $ (w_\delta')_{\delta \in \rootbase} $ be a $ \rootbase $-system of Weyl elements in $ G $. A \defemph*{$ B_n $-admissible partial twisting system for $ (G, (w_\delta')_{\delta \in \rootbase}) $}\index{twisting system!partial!admissible (type Bn)@admissible (for $ B_n $)} is a partial twisting system $ (\twistgroup', \inverparsym', \invogroup', \invoparsym') $ for $ (G, (w_\delta')_{\delta \in \rootbase}) $ with the following additional properties:
	\begin{stenumerate}
		\item $ \twistgroup' = \compactSet{\pm 1} $ and $ -1_{\twistgroup'} $ acts on all root groups of $ G $ by inversion.
		
		\item \label{B:admissible-parmap-def:act}$ \invogroup' = \compactSet{\pm 1} $.
		
		\item \label{B:admissible-parmap-def:adj}$ \inverparsym' $ and $ \invoparsym' $ are adjacency-trivial.
		
		\item \label{B:admissible-parmap-def:square}$ \inverparsym' $ satisfies the square formula.
		
		\item \label{B:admissible-parmap-def:same}For all $ i,j \in \numint{1}{n} $ and $ k \in \numint{1}{n} \setminus \compactSet{i,j} $, the words $ \word{\rho}^i $, $ \word{\rho}^j $ from \cref{B:Bnsub-ext-word} satisfy $ \inverpar{\epsilon \basvec_k}{\word{\rho}^i}' = \inverpar{\epsilon \basvec_k}{\word{\rho}^j}' $ for all $ \epsilon \in \compactSet{\pm 1} $.
		
		\item \label{B:admissible-parmap-def:ortho}If $ \rho, \zeta $ are orthogonal short roots in $ B_n $, then $ \invopar{\rho}{\reflbr{\zeta}}' = -1_{\invogroup'} $.
	\end{stenumerate}
	We will sometimes refer to such objects as \defemph*{admissible partial twisting systems} if the root system $ B_n $ is clear from the context.
\end{definition}

\begin{proposition}
	Let $ G $ be a group with a $ B_n $-pregrading $ (\rootgr{\alpha})_{\alpha \in B_n} $. Then for any $ \rootbase $-system $ (w_\delta)_{\delta \in \rootbase} $ of Weyl elements in $ G $, the standard partial twisting system $ (\twistgroup, \inverparsym, \invogroup, \invoparsym) $ is a $ B_n $-admissible partial twisting system for $ (G, (w_\delta)_{\delta \in \rootbase}) $. In particular, admissible partial twisting systems exist for each group with a $ B_n $-pregrading.
\end{proposition}
\begin{proof}
	We have to check that $ (\twistgroup, \inverparsym, \invogroup, \invoparsym) $ satisfies the axioms of a partial twisting system for $ (G, (w_\delta)_{\delta \in \rootbase}) $ (see \cref{param:partwist-def}) and the axioms in \cref{B:admissible-parmap-def}. By \cref{param:pargroup-inv-example}, $ \twistgroup $ is a twisting group for $ (G, (w_\delta)_{\delta \in \rootbase}) $. All the remaining axioms follow from the results in this section.
\end{proof}

\begin{note}
	We have proven in \cref{B:invo-parmap-short} that $ \invopar{\rho}{\reflbr{\zeta}} = -1_{\invogroup} $ holds for any pair of short roots. In Axiom~\thmitemref{B:admissible-parmap-def}{B:admissible-parmap-def:ortho}, however, we only require that this holds for orthogonal pairs $ (\rho, \zeta) $. In other words, we do not require that $ \invopar{\rho}{\reflbr{\rho}} = -1_{\invogroup} $. In fact, we could for any $ i \in \numint{1}{n} $ replace the homomorphism $ \risshneg{i} $ in \cref{B:ex-roothom-def-short} by $ \map{\risshneg{i}'}{}{}{v}{\risshneg{i}(\refl{v_0}(v))} $. Then the corresponding Weyl elements would have to be redefined as
	\[ w_i(u) \defl w_{\basvec_i}(u) \defl \risom{-\basvec_i}\brackets[\big]{-q(u)^{-1} \refl{v_0}(u)} \circ \risom{\basvec_i}(u) \circ \risom{-\basvec_i}\brackets[\big]{-q(u)^{-1} \refl{v_0}(u)} \]
	for all $ u \in \module $ for which $ q(u) $ is invertible. Further, the resulting parity map $ \invoparsym $ would then satisfy $ \invopar{\pm\basvec_i}{\reflbr{\basvec_i}} = 1_{\invogroup} $, but we would still have $ \invopar{\pm\basvec_i}{\reflbr{\basvec_j}} = -1_{\invogroup} $ for all $ j \in \numint{1}{n} \setminus \compactSet{i} $. This illustrates that property~\thmitemref{B:admissible-parmap-def}{B:admissible-parmap-def:ortho} is intrinsic to the group $ \EO(q) $ (that is, invariant under a twisting of the parametrisation) while the value of $ \invopar{\rho}{\reflbr{\rho}} $ can be different for distinct parametrisations of $ \EO(q) $.
\end{note}


\section{The Parametrisation}

\label{sec:B-param}

\begin{secnotation}
	We fix an integer $ n \ge 3 $ and consider the root system $ B_n $ in its standard representation (as in \cref{B:Bn-standard-rep}) with its standard root base $ \rootbase $. We denote by $ G $ a group which has crystallographic $ B_n $-commutator relations with root groups $ (\rootgr{\alpha})_{\alpha \in B_n} $ such that $ \invset{\alpha} $ is non-empty for all $ \alpha \in \roots $ and such that $ G $ is rank-2-injective. We fix a $ \rootbase $-system $ (w_\delta)_{\delta \in \rootbase} $ of Weyl elements and we denote by $ (\twistgroup, \inverparsym, \invogroup, \invoparsym) $ a $ B_n $-admissible partial twisting system for $ (G, (w_\delta)_{\delta \in \rootbase}) $.
\end{secnotation}

In this section, we show that $ G $ satisfies the conditions in the parametrisation theorem with respect to any $ B_n $-admissible partial twisting system. In particular, this holds for the standard partial twisting system, which will be the only case that we are interested in. After the work in the previous sections, all that is truly left to do is the verification of stabiliser-compatibility. For this we use the criterion from \cref{param:stabcomp-crit-ortho}.

\begin{proposition}\label{B:square-comp}
	$ G $ is square-compatible with respect to $ \inverparsym $.
\end{proposition}
\begin{proof}
	Since both $ G $ and $ \inverparsym $ satisfy the square formula (by \cref{B:square-act:summary} and Axiom~\thmitemref{B:admissible-parmap-def}{B:admissible-parmap-def:square}), this follows from \cref{param:square-formula-comp}. Here we have to use Axiom~\thmitemref{BC:admissible-parmap-def}{BC:admissible-parmap-def:act}.
\end{proof}

\begin{proposition}\label{B:stab-comp}
	$ G $ is stabiliser-compatible with respect to $ (\inverparsym, \invoparsym) $ and $ (w_\delta)_{\delta \in \rootbase} $.
\end{proposition}
\begin{proof}
	Let $ \alpha $ be any root. We know from Axiom~\thmitemref{B:admissible-parmap-def}{B:admissible-parmap-def:adj} that $ \inverparsym $ is $ \alpha $-adjacency-trivial and we want to show that $ G $ is $ \alpha $-stabiliser-compatible. If $ \alpha $ is long, then it has the property that all roots which are orthogonal to $ \alpha $ are adjacent to $ \alpha $, and thus it follows from \cref{param:adj-implies-stab} that $ G $ is $ \alpha $-stabiliser-compatible with respect to $ \inverparsym $. Thus it is also $ \alpha $-stabiliser-compatible with respect to $ (\inverparsym, \invoparsym) $ by \cref{param:stab-only-inver-rem}, as desired.
	
	Now assume that $ \alpha $ is short. Write $ \alpha = \epsilon \basvec_k $ for some $ \epsilon \in \compactSet{\pm 1} $ and $ k \in \numint{1}{n} $. As in \cref{param:stabcomp-crit-ortho}, we consider the following sets:
	\begin{align*}
		\calO &\defl \Set{\beta \in \roots \given \alpha \cdot \beta = 0} = \Set{\sigma \basvec_i \given i \in \numint{1}{n} \setminus \compactSet{k}, \sigma \in \compactSet{\pm 1}}, \\
		\calA &\defl \Set{\beta \in \calO \given \alpha \text{ is crystallographically adjacent to } \beta} = \emptyset, \\
		\bar{\calA} & \defl \calO \setminus \calA = \calO.
	\end{align*}
	For all $ \beta \in \bar{\calA} $, we have $ \invopar{\alpha}{\reflbr{\beta}} = -1_{\invogroup} $ by Axiom~\thmitemref{B:admissible-parmap-def}{B:admissible-parmap-def:ortho}. Now let $ \beta, \beta' \in \bar{\calA} $. Let $ i,j \in \numint{1}{n} $ and $ \epsilon, \epsilon' \in \compactSet{\pm 1} $ such that $ \beta = \epsilon \basvec_i $ and $ \beta' = \epsilon' \basvec_j $. Let $ \word{\delta} \defl \word{\rho}^i $ and $ \word{\delta}' \defl \word{\rho}^j $ be the words from \cref{B:Bnsub-ext-word}. By definition, they satisfy $ \reflbr{\word{\delta}} = \reflbr{\basvec_i} = \reflbr{\beta} $ and $ \reflbr{\word{\delta}'} = \reflbr{\basvec_j} = \reflbr{\beta'} $. Further, we have $ \inverpar{\alpha}{\word{\delta}} = \inverpar{\alpha}{\word{\delta}'} $ by Axiom~\thmitemref{B:admissible-parmap-def}{B:admissible-parmap-def:same} and the actions of $ w_{\word{\delta}} $ and $ w_{\word{\delta}'} $ are identical by \cref{B:short-weyl-same-action}. Thus all conditions in \cref{param:stabcomp-crit-ortho} are satisfied, and we conclude that $ G $ is $ \alpha $-stabiliser-compatible. This finishes the proof.
\end{proof}

\begin{proposition}\label{B:param-exists}
	There exist abelian groups $ (\comring, +) $ and $ (\module, +) $ (each equipped, as a set, with an action of $ \twistgroup \times \invogroup $) and a parametrisation $ (\risom{\alpha})_{\alpha \in B_n} $ of $ G $ by $ (\twistgroup \times \invogroup, \module, \comring) $ with respect to $ \inverparsym \times \invoparsym $ and $ (w_\delta)_{\delta \in \rootbase} $ such that the action of $ \twistgroup $ on $ \comring $ and $ \module $ is given by group inversion.
\end{proposition}
\begin{proof}
	This follows from the parametrisation theorem (\cref{param:thm}), whose assumptions are satisfied by \cref{B:stab-comp,B:square-comp,braid:all}.
\end{proof}


\section{Computation of the Blueprint Rewriting Rules}

\label{sec:B-bluerules}

\begin{secnotation}\label{B:blue-secnotation}
	We fix an integer $ n \ge 3 $ and consider the root system $ B_n $ in its standard representation (as in \cref{B:Bn-standard-rep}) with its standard root base $ \rootbase $.
	We denote by $ G $ a group with a crystallographic $ B_n $-grading $ (\rootgr{\alpha})_{\alpha \in B_n} $. We fix a $ \rootbase $-system of Weyl elements $ (w_\delta)_{\delta \in \rootbase} $ and denote its standard $ \Bnsub $-extension by $ (w_\beta)_{\beta \in \Bnsub} $. To simplify notation, we put $ w_{ij} \defl w_{\basvec_i - \basvec_j} $ for all distinct $ i,j \in \numint{1}{n} $ and we put $ w_i \defl w_{\basvec_i} $ for all $ i \in \numint{1}{n} $. We denote the standard partial twisting system for $ G $ by $ (\twistgroup, \inverparsym, \invogroup, \invoparsym) $, by $ (\comring, +) $, $ (\module, +) $ any groups which satisfy the assertion of \cref{B:param-exists} and by $ (\risom{\alpha})_{\alpha \in B_n} $ the corresponding parametrisation of $ G $. The groups $ \comring $ and $ \module $ are equipped with actions of $ \twistgroup \times \invogroup $, and we call the map
	\[ \map{}{\module}{\module}{v}{\modinv{v} \defl -1_\invogroup.v} \]
	the \defemph*{involution on $ \module $}. Further, we choose elements $ v_{-1}, v_0, v_1 \in \module $ such that $ w_{n} = \risshneg{n}(v_{-1}) \risshpos{n}(v_0) \risshneg{n}(v_1) $.
\end{secnotation}

In this section, we define the commutation maps of the parametrisation $ (\risom{\alpha})_{\alpha \in B_n} $ and derive their rank-2 identities. Our goal is to define the blueprint rewriting rules for $ B_n $ and to prove that they are indeed blueprint rewriting rules.

We will frequently apply the following result without reference. It says that the action of the $ \Bnsub $-extension $ (w_\beta)_{\beta \in \Bnsub} $ of $ (w_\delta)_{\delta \in \rootbase} $ on the root groups is determined by the values in \cref{B:ex-parmap-def}.

\begin{proposition}\label{B:Bnsub-conj-anygroup}
	Let $ \alpha \in B_n $ and let $ \beta \in \Bnsub $. Let $ x \in \comring $ if $ \alpha $ is long and let $ x \in\module $ if $ \alpha $ is short. Then $ \risom{\alpha}(x)^{\hat{w}_\beta} = \risom{\refl{\beta}(\alpha)}(\inverpar{\alpha}{\beta} \invopar{\alpha}{\beta}.x) $ where $ \inverpar{\alpha}{\beta} $ and $ \invopar{\alpha}{\beta} $ are the values in \cref{B:ex-parmap-def}.
\end{proposition}
\begin{proof}
	Denote by $ \word{\beta} $ the standard $ \rootbase $-expression of $ \beta $ from \cref{B:Bnsub-ext-word}. Since $ G $ is parametrised by $ (\twistgroup \times \invogroup, \module, \comring) $ with respect to $ \inverparsym \times \invoparsym $ and $ (w_\delta)_{\delta \in \rootbase} $, we have
	\[ \risom{\alpha}(x)^{\hat{w}_\beta} = \risom{\refl{\beta}(\alpha)}(\inverpar{\alpha}{\word{\beta}} \invopar{\alpha}{\word{\beta}}.x). \]
	We know from \cref{B:ex-parmap-equality} that $ \inverpar{\alpha}{\word{\beta}} \invopar{\alpha}{\word{\beta}} = \inverpar{\alpha}{\beta} \invopar{\alpha}{\beta} $, so the assertion follows.
\end{proof}

\begin{lemma}\label{B:invo-lem}
	The involution $ \map{}{\module}{\module}{v}{\modinv{v}} $ is an automorphism of $ \module $ with $ \modinv{\modinv{v}} = v $ for all $ v \in \module $.
\end{lemma}
\begin{proof}
	Since the involution is induced by the action of $ -1_\invogroup $ and $ (-1_\invogroup)^2 = 1_\invogroup $, we have $ \modinv{\modinv{v}} = v $ for all $ v \in \module $. Further, since $ \inverinvopar{\basvec_1}{\basvec_n} = (1_\twistgroup, -1_\invogroup) $, we have
	\begin{align*}
		\risshpos{1}\brackets[\big]{\modinv{v+u}} &= \risshpos{1}(v+u)^{w_{n}} = \risshpos{1}(v)^{w_n} \risshpos{1}(u)^{w_n} = \risshpos{1}(\modinv{v}) \risshpos{1}(\modinv{u}) = \risshpos{1}(\modinv{v} + \modinv{u})
	\end{align*}
	for all $ v,u \in \module $, so $ \modinv{v+u} = \modinv{v} + \modinv{u} $.
\end{proof}

\begin{lemma}
	For all $ i \in \numint{1}{n} $, we have $ w_{i} = \risshneg{i}(v_{-1}) \risshpos{i}(v_0) \risshneg{i}(v_1) $ for all $ i \in \numint{1}{n} $.
\end{lemma}
\begin{proof}
	By the choice of $ v_{-1}, v_0, v_1 $, this is clear for $ i=n $. Since $ w_{i} = w_{n}^{w_{ni}} $ and $ \inverinvopar{\pm \basvec_n}{\basvec_n - \basvec_i} = (1,1) $, the general assertion follows from \cref{B:Bnsub-conj-anygroup}.
\end{proof}

\begin{definition}[Commutation maps]\label{B:commmap-def}
	We define
	\begin{align*}
		\map{\mapdot}{\comring \times \comring}{\comring&}{(a,b)}{ab \defl a \rmult b}, \\
		\map{g}{\comring \times \module}{\module&}{(a,v)}{g(a,v)}, \\
		\map{q}{\comring \times \module}{\comring&}{(a,v)}{q(a,v)}, \\
		\map{f}{\module \times \module}{\comring&}{(u,v)}{f(u,v)}
	\end{align*}
	to be the unique maps which satisfy
	\begin{align*}
		\commutator{\rismin{1}{2}(a)}{\rismin{2}{3}(b)} &= \rismin{1}{3}(a \cdot b), \\
		\commutator{\rismin{1}{2}(a)}{\risshpos{2}(v)} &= \risshpos{2}\brackets[\big]{g(a,v)} \risplus{1}{2}\brackets[\big]{-q(a,v)}, \\
		\commutator{\risshpos{1}(u)}{\risshpos{2}(v)} &= \risplus{1}{2}\brackets[\big]{-f(u,v)}
	\end{align*}
	for all $ a,b \in \comring $ and all $ u,v \in \module $.
\end{definition}

Our goal is to show that $ \rmult $ turns $ \comring $ into a commutative associative ring, that $ g $ defines a $ \comring $-modules structure on $ \module $, that $ \map{}{\module}{\comring}{v}{q(1_\comring, v)} $ is a quadratic form with linearisation $ f $ and that $ q(r,v) = r q(1_\comring, v) $ for all $ r \in \comring $ and $ v \in \module $.

\begin{lemma}[Rank-2 identities, part 1]\label{B:param:rank2commrel}
	The maps $ \rmult $, $ g $, $ q $ and $ f $ are additive in all components except for the second component of $ q $, which satisfies
	\[ q(a, v+w) = q(a,v) + q(a,w) + f\brackets[\big]{g(a,v), w} \]
	for all $ a \in \comring $ and all $ v,w \in \module $.
\end{lemma}
\begin{proof}
	The additivity assertions follow from \cref{basic:comm-add}, as in the proof of \cref{ADE:distributivity}. Now let $ a \in \comring $ and let $ v,w \in \module $. By \thmitemcref{B:comm-add-cry}{B:comm-add-cry:delta-gamma} and the definition of the maps $ q $, $ g $ and $ f $, we have
	\begin{align*}
		\risplus{1}{2}\brackets[\big]{-q(r,v+w)} \hspace{-3cm}& \\
		&= \commpart{\rismin{1}{2}(a)}{\risshpos{2}(v) \risshpos{2}(w)}{\basvec_1 + \basvec_2} \\
		&= \commpart{\rismin{1}{2}(a)}{\risshpos{2}(w)}{\basvec_1 + \basvec_2} \commutator[\big]{\commpart{\rismin{1}{2}(a)}{\risshpos{2}(v)}{\basvec_1}}{\risshpos{2}(w)} \commpart{\rismin{1}{2}(a)}{\risshpos{2}(v)}{\basvec_1 + \basvec_2} \\
		&= \risplus{1}{2}\brackets[\big]{-q(a,w)} \commutator[\big]{\risshpos{1}\brackets[\big]{g(a,v)}}{\risshpos{2}(w)} \risplus{1}{2}\brackets[\big]{-q(a,v)} \\
		&= \risplus{1}{2}\brackets[\big]{-q(a,w) - f\brackets[\big]{g(a,v), w} - q(a,v)}.
	\end{align*}
	By the bijectivity of $ \risplus{1}{2} $, the second assertion follows.
\end{proof}

\begin{lemma}\label{B:comm-mult-computation}
	For all pairwise distinct $ i,j,k \in \numint{1}{n} $ and all $ a,b \in \comring $, the following commutator relations hold:
	\begin{align*}
		\commutator{\rismin{i}{j}(a)}{\rismin{j}{k}(b)} &= \rismin{i}{k}(ab), \\
		\commutator{\rismin{i}{j}(a)}{\risplus{j}{k}(b)} &= \risplus{i}{k}\brackets{\delmin{k \in \betint{i}{j}} ab}, \\
		\commutator{\rismin{i}{j}(a)}{\risminmin{k}{i}(b)} &= \risminmin{j}{k}\brackets{\delmin{k \nin \betint{i}{j}} ab}, \\
		\commutator{\risplus{i}{j}(a)}{\risminmin{k}{i}(b)} &= \rismin{j}{k}\brackets{\delmin{i \in \betint{j}{k}} ab}.
	\end{align*}
\end{lemma}
\begin{proof}
	Let $ a,b \in \comring $. We will frequently apply \cref{B:Bnsub-conj-anygroup} without explicitly saying so. By definition, the first assertion holds if $ (i,j,k) = (1,2,3) $. Thus it suffices to show that if $ i,j,k \in \numint{1}{n} $ are distinct and $ \sigma $ is a permutation of $ \numint{1}{n} $, then the validity of the first equation for $ (i,j,k) $ implies the validity of the first equation for $ (i^\sigma, j^\sigma, k^\sigma) $. Since the group of permutations is generated by transpositions, it suffices to prove this for transpositions $ \sigma $.
	Thus let $ i,j,k \in \numint{1}{n} $ be distinct and such that
	\begin{equation}\label{eq:ringmult-comm}\tag{$ * $}
		\commutator{\risom{\basvec_i - \basvec_j}(a)}{\risom{\basvec_j - \basvec_k}(b)} = \risom{\basvec_i - \basvec_k}(a \rmult b) \quad \text{for all } a,b \in \comring
	\end{equation}
	and let $ \sigma $ be the transposition which interchanges the distinct numbers $ p,q \in \numint{1}{n} $. If $ \Set{i,j,k} \intersect \Set{p,q} = \emptyset $, the assertion is clear. If $ \abs{\Set{i,j,k} \intersect \Set{p,q}} = 1 $, conjugating \eqref{eq:ringmult-comm} by $ w_{pq} $ (if $ p \in \Set{i,j,k} $) or by $ w_{qp} $ (if $ q \in \Set{i,j,k} $) shows that the assertion holds. Now assume that $ \abs{\Set{i,j,k} \intersect \Set{p,q}} = 2 $. If $ \Set{p,q} = \Set{i,j} $, then conjugating \eqref{eq:ringmult-comm} by $ w_{ij} $ and applying \cref{basic:comm-add} yields
	\begin{align*}
		\commutator{\risom{\basvec_i - \basvec_j}(a)}{\risom{\basvec_j - \basvec_k}(b)}^{w_{ij}} &= \commutator{\risom{\basvec_i - \basvec_j}(r)^{w_{ij}}}{\risom{\basvec_j - \basvec_k}(s)^{w_{ij}}} = \commutator{\risom{\basvec_j - \basvec_i}(-r)}{\risom{\basvec_i - \basvec_k}(-s)} \\
		&= \commutator{\rismin{j}{i}(r)}{\rismin{i}{k}(s)} \rightand \\
		\commutator{\risom{\basvec_i - \basvec_j}(a)}{\risom{\basvec_j - \basvec_k}(b)}^{w_{ij}} &= \rismin{i}{k}(r \rmult s)^{w_{ij}} = \risom{\basvec_j - \basvec_k}(r \rmult s).
	\end{align*}
	The assertion follows in this case. Similarly, conjugating \eqref{eq:ringmult-comm} by $ w_{jk} $ or by $ w_{ik} $ shows that the assertion holds for the transposition interchanging $ j $ and $ k $ or $ i $ and $ k $, respectively. This finishes the proof of the first equation.
	
	Now conjugating the first equation by $ w_k $ yields
	\[ \commutator{\rismin{i}{j}(a)}{\risplus{j}{k}(\delmin{k>j}b)} = \risplus{i}{k}\brackets{\delmin{k > i} ab} \quad \text{for all } a,b \in \comring, \]
	which (by replacing $ b $ with $ \delmin{k>j}b $) implies that
	\[ \commutator{\rismin{i}{j}(a)}{\risplus{j}{k}(b)} = \risplus{i}{k}\brackets{\delmin{k > i} \delmin{k>j} ab} \quad \text{for all } a,b \in \comring. \]
	By going through all possible cases, it is easy to observe that $ \delmin{k > i} \delmin{k>j} = \delmin{k \in \betint{i}{j}} $, which proves the second equation. In a similar manner, we can conjugate the first equation by $ w_j $ to obtain
	\[ \commutator{\risplus{i}{j}(\delmin{j>i}a)}{\risminmin{j}{k}(\delmin{j>k}b)} = \rismin{i}{k}(ab) \quad \text{for all } a,b \in \comring. \]
	This is equivalent to
	\[ \commutator{\risplus{i}{j}(a)}{\risminmin{j}{k}(b)} = \rismin{i}{k}(\delmin{j \in \betint{i}{k}}) \quad \text{for all } a,b \in \comring. \]
	Interchanging the roles of $ i $ and $ j $ yields the fourth equation. Finally, conjugating the fourth equation by $ w_j $, we see that
	\[ \commutator{\rismin{i}{j}(\delmin{j>i}a)}{\risminmin{k}{i}(b)} = \risminmin{j}{k}(\delmin{j>k} \delmin{i \in \betint{j}{k}} ab) \quad \text{for all } a,b \in \comring. \]
	This implies that
	\[ \commutator{\rismin{i}{j}(a)}{\risminmin{k}{i}(b)} = \risminmin{j}{k}(\delmin{j>i} \delmin{j>k} \delmin{i \in \betint{j}{k}} ab) \quad \text{for all } a,b \in \comring. \]
	Consulting the truth table in \cref{B:param:truth}, we see that $ \delmin{j>i} \delmin{j>k} \delmin{i \in \betint{j}{k}} = \delmin{k \nin \betint{i}{j}} $, which finishes the proof.
\end{proof}

\begin{figure}[htb]
	\centering\begin{tabular}{cccccc}
		\toprule
		& $ \delmin{j>i} $ & $ \delmin{j>k} $ & $ \delmin{i \in \betint{j}{k}} $ & $ \delmin{j>i} \delmin{j>k} \delmin{i \in \betint{j}{k}} $ & $ \delmin{k \nin \betint{i}{j}} $ \\
		\midrule
		$ i<j<k $ & $ - $ & $ + $ & $ + $ & $ - $ & $ - $ \\
		$ i<k<j $ & $ - $ & $ - $ & $ + $ & $ + $ & $ + $ \\
		$ j<i<k $ & $ + $ & $ + $ & $ - $ & $ - $ & $ - $ \\
		$ j<k<i $ & $ + $ & $ + $ & $ + $ & $ + $ & $ + $ \\
		$ k<i<j $ & $ - $ & $ - $ & $ - $ & $ - $ & $ - $ \\
		$ k<j<i $ & $ + $ & $ - $ & $ + $ & $ - $ & $ - $ \\
		\bottomrule
	\end{tabular}
	\caption{A truth table which shows that $ \delmin{j>i} \delmin{j>k} \delmin{i \in \betint{j}{k}} = \delmin{k \nin \betint{i}{j}} $ for all pairwise distinct $ i,j,k \in \numint{1}{n} $.}
	\label{B:param:truth}
\end{figure}

\begin{remark}\label{B:comm-formula-firststep}
	We could perform the same computations as in \cref{B:comm-mult-computation} for the other commutation maps. However, it will be more efficient to do so at a later point, when we have acquired more information about the maps $ f, q $ and $ g $. For the moment, we only note that for all $ i<j \in \numint{1}{n} $, conjugating the equations in \cref{B:commmap-def} by $ w_{2j} w_{1i} $ (where $ w_{kk} $ is interpreted as $ 1_G $) yields that
	\begin{align*}
		\commutator{\rismin{i}{j}(a)}{\risshpos{j}(v)} &= \risshpos{i}\brackets[\big]{g(a,v)} \risplus{i}{j}\brackets[\big]{-q(a,v)}, \\
		\commutator{\risshpos{i}(u)}{\risshpos{j}(v)} &= \risplus{i}{j}\brackets[\big]{-f(u,v)}
	\end{align*}
	for all $ a \in \comring $ and all $ u,v \in \module $. We delay the remaining computations until \cref{B:blue:stand-signs}.
\end{remark}

\begin{lemma}[Rank-2 identities, part 2]\label{B:isring}
	$ (\comring, +, \rmult) $ is a ring. If we denote its identity element by $ 1_\comring $, we have $ w_{ij} = \rismin{j}{i}(-1_\comring) \rismin{i}{j}(1_\comring) \rismin{j}{i}(-1_\comring) $ for all distinct $ i, j \in \numint{1}{n} $.
\end{lemma}
\begin{proof}
	We already know from \cref{B:param:rank2commrel} that the multiplication satisfies the distributive law, so we only have to verify the existence of an identity element. Let $ i,j \in \numint{1}{n} $ be distinct and let $ x \in \comring $ be arbitrary. Since $ n \ge 3 $, we can choose some $ k \in \numint{1}{n} \setminus \compactSet{i,j} $. By \cref{A2Weyl:weyl}, there exist unique $ a,b \in \comring $ such that $ w_{ij} = \rismin{j}{i}(a) \rismin{i}{j}(b) \rismin{j}{i}(a) $ and we have $ \rismin{j}{i}(a) = \rismin{i}{j}(b)^{w_{ij}} = \rismin{j}{i}(-b) $, so $ a=-b $. Using \thmitemcref{A2Weyl:basecomp-cor}{A2Weyl:basecomp-cor:1}, we can now compute that
	\begin{align*}
		\rismin{k}{j}(x \rmult b) &= \commutator{\rismin{k}{i}(x)}{\rismin{i}{j}(b)} = \rismin{k}{i}(x)^{w_{ij}} = \rismin{k}{j}(x), \\
		\rismin{i}{k}(b \rmult x) &= \commutator{\rismin{i}{j}(b)}{\rismin{j}{k}(x)} = \commutator{\rismin{j}{k}(x)}{\rismin{i}{j}(b)}^{-1} \\
		&= \brackets{\rismin{j}{k}(x)^{w_{ij}}}^{-1} = \rismin{i}{k}(x).
	\end{align*}
	It follows that $ \comring $ is a ring with identity element $ 1_\comring \defl b $ and that
	\[ w_{ij} = \rismin{j}{i}(-1_\comring) \rismin{i}{j}(1_\comring) \rismin{j}{i}(-1_\comring). \]
	Since the identity element of a ring is uniquely determined by the ring structure, our construction does not depend on our choice of $ i $ and $ j $. Thus we also have
	\[ w_{kl} = \rismin{l}{k}(-1_\comring) \rismin{k}{l}(1_\comring) \rismin{l}{k}(-1_\comring) \]
	for any choice of distinct indices $ k,l \in \numint{1}{n} $.
\end{proof}

\begin{note}
	The assertion of \cref{B:isring} involves only roots which lie in the canonical $ A_{n-1} $-subsystem of $ B_n $. Further, the signs of $ \inverparsym $ on this subsystem are the same as in \cref{ADE:An-parmap} by \cref{B:ex-parmap-restrict} and $ \invoparsym $ is trivial on this subsystem. Thus \cref{B:isring} is also a direct consequence of \cref{ADE:ring-identity}.
\end{note}

\begin{lemma}[Rank-2 identities, part 3]
	We have $ g(1_\comring, v) = v $ for all $ v \in \module $ and $ q(a,v_0) = a $ for all $ a \in \comring $.
\end{lemma}
\begin{proof}
	Let $ v \in \module $. By the definition of $ g $, we have
	\[ \risshpos{2}(g(1, v)) = \commpart{\rismin{1}{2}(1)}{\risshpos{2}(v)}{\basvec_2}. \]
	By an application of \cref{commpart:inv-switch} and \thmitemcref{B:basecomp-walpha-cor:delta}{B:basecomp-walpha-cor:delta:act}, we can compute that
	\begin{align*}
		\commpart{\rismin{1}{2}(1)}{\risshpos{2}(v)}{\basvec_2} &= \commpart{\risshpos{2}(v)}{\rismin{1}{2}(1)}{\basvec_2}^{-1} = \brackets[\big]{\risshpos{2}(v)^{w_{12}}}^{-1} = \risshpos{2}(v),
	\end{align*}
	so $ g(1,v) = v $. Now let $ a \in \comring $. In a similar way as above, it follows from \cref{B:Bnsub-conj-anygroup} and \thmitemcref{B:basecomp-wdelta-alpha}{B:basecomp-wdelta-alpha:act} that
	\begin{align*}
		\risplus{1}{2}\brackets[\big]{-q(a, v_0)} &= \commpart{\rismin{1}{2}(a)}{\risshpos{2}(v_0)}{\basvec_1 + \basvec_2} = \rismin{1}{2}(a)^{w_2} = \risplus{1}{2}(-a),
	\end{align*}
	so $ q(a,v_0) = a $.
\end{proof}

We now compute the blueprint rewriting rules. As explained in \cref{blue:rank-3}, we will only perform the blueprint computation in the rank-3 case: Namely, in the subgroup of $ G $ corresponding to the root subsystem spanned by
\[ \Set{\basvec_1 - \basvec_2, \basvec_2 - \basvec_3, \basvec_3}. \]

\begin{definition}[Blueprint rewriting rules]\label{B:rewriting-def}
	We define the following rewriting rules:
	\begin{align*}
		\map{\blutrans{12}}{\rootgrmin{1}{2} \times \rootgrmin{2}{3} \times \rootgrmin{1}{2}&}{\rootgrmin{2}{3} \times \rootgrmin{1}{2} \times \rootgrmin{1}{2}}{\\\brackets[\big]{\rismin{1}{2}(a), \rismin{2}{3}(b), \rismin{1}{2}(c)}&}{\brackets[\big]{\rismin{2}{3}(c), \rismin{1}{2}(-b-ca), \rismin{2}{3}(a)}}, \\
		\map{\blutrans{12}^{-1}}{\rootgrmin{2}{3} \times \rootgrmin{1}{2} \times \rootgrmin{1}{2}&}{\rootgrmin{1}{2} \times \rootgrmin{2}{3} \times \rootgrmin{1}{2}}{\\\brackets[\big]{\rismin{2}{3}(a), \rismin{1}{2}(b), \rismin{2}{3}(c)}&}{\brackets[\big]{\rismin{1}{2}(c), \rismin{2}{3}(-b-ac), \rismin{1}{2}(a)}}, \\
		\map{\phi_{13}}{\rootgrmin{1}{2} \times \rootgrshpos{3}&}{\rootgrshpos{3} \times \rootgrmin{1}{2}}{\\\brackets[\big]{\rismin{1}{2}(a), \risshpos{3}(v)}&}{\brackets[\big]{\risshpos{3}(v), \rismin{1}{2}(a)}}, \\
		\map{\phi_{13}^{-1}}{\rootgrshpos{3} \times \rootgrmin{1}{2}&}{\rootgrmin{1}{2} \times \rootgrshpos{3}}{\\\brackets[\big]{\risshpos{3}(v), \rismin{1}{2}(a)}&}{\brackets[\big]{\rismin{1}{2}(a), \risshpos{3}(v)}}
	\end{align*}
	and
	\[ \map{\blutrans{23}}{\rootgrshpos{3} \times \rootgrmin{2}{3} \times \rootgrshpos{3} \times \rootgrmin{2}{3}}{\rootgrmin{2}{3} \times \rootgrshpos{3} \times \rootgrmin{2}{3} \times \rootgrshpos{3}}{}{} \]
	which maps $ \brackets[\big]{\risshpos{3}(v), \rismin{2}{3}(a), \risshpos{3}(u), \rismin{2}{3}(b)} $ to
	\[ \brackets[\bigg]{\rismin{2}{3}(b), \risshpos{3}\brackets[\big]{\modinv{u} - \modinv{g(b, \modinv{v})}}, \rismin{2}{3}\brackets[\big]{x}, \risshpos{3}(-\modinv{v})} \]
	where $ x \defl a - q(b, -\modinv{v}) - f\brackets[\big]{u - g(b, \modinv{v}), \modinv{v}} $.
\end{definition}

\begin{lemma}
	The maps $ \blutrans{12} $, $ \blutrans{12}^{-1} $, $ \phi_{13} $, $ \phi_{13}^{-1} $ and $ \blutrans{23} $ in \cref{B:rewriting-def} are blueprint rewriting rules (with respect to $ (w_\delta)_{\delta \in \rootbase'} $). Further, $ \blutrans{12} $ and $ \blutrans{12}^{-1} $ are inverses of each other, and the same holds for $ \phi_{13} $ and $ \phi_{13}^{-1} $.
\end{lemma}
\begin{proof}
	The second statement about inverses is easy to verify. Further, we know from \cref{blue:rewriting-switch} that $ \phi_{13} $ and $ \phi_{13}^{-1} $ are blueprint rewriting rules. Observe that the rewriting rule $ \blutrans{12} $ and its inverse are the same rules that we have used for $ A_3 $. Since the restriction of $ \inverparsym $ to $ A_2 = \Set{\basvec_i - \basvec_{i+1} \given i \in \numint{1}{2}} $ yields the same parity map that we have used for $ A_3 $ (see \cref{B:ex-parmap-restrict}), the same computation as in \cref{blue:A3-rewriting-comp} shows that $ \blutrans{12} $ and $ \blutrans{12}^{-1} $ are blueprint rewriting rules. It remains to show that $ \blutrans{23} $ is a blueprint rewriting rule. For this, let $ a,b,c,d \in \comring $, let $ v,u,h,k \in \module $ and put
	\[ \word{\alpha} \defl \brackets{\basvec_3, \basvec_2 - \basvec_3, \basvec_3, \basvec_2 - \basvec_3} \midand \word{\beta} \defl \brackets{\basvec_2 - \basvec_3, \basvec_3, \basvec_2 - \basvec_3, \basvec_3}. \]
	Further, we set
	\[ x \defl \brackets[\big]{\risshpos{3}(v), \rismin{2}{3}(a), \risshpos{3}(u), \rismin{2}{3}(b)}. \]
	On the one hand, we have
	\begin{align*}
		\blumapG{\word{\alpha}}(x) &= w_3 \risshpos{3}(v) w_{23} \rismin{2}{3}(a) w_3 \risshpos{3}(u) w_{23} \rismin{2}{3}(b) \\
		&= w_3 w_{23} w_3 w_{23} \risshpos{3}(v)^{w_{23} w_3 w_{23}} \rismin{2}{3}(a)^{w_3 w_{23}} \risshpos{3}(u)^{w_{23}} \rismin{2}{3}(b) \\
		&= w_3 w_{23} w_3 w_{23} \risshpos{2}(-v)^{w_3 w_{23}} \risplus{2}{3}(-a)^{w_{23}} \risshpos{2}(-u) \rismin{2}{3}(b) \\
		&= w_3 w_{23} w_3 w_{23}\risshpos{3}(-\modinv{v}) \risplus{2}{3}(-a) \risshpos{2}(-u) \rismin{2}{3}(b).
	\end{align*}
	On the other hand, we put
	\[ x' \defl \brackets[\big]{\rismin{2}{3}(c), \risshpos{3}(h), \rismin{2}{3}(d), \risshpos{3}(k)}. \]
	Then
	\begin{align*}
		\blumapG{\word{\beta}}(x') &= w_{23} \rismin{2}{3}(c) w_3 \risshpos{3}(h) w_{23} \rismin{2}{3}(d) w_3 \risshpos{3}(k) \\
		&= w_{23} w_3 w_{23} w_3 \rismin{2}{3}(c)^{w_3 w_{23} w_3} \risshpos{3}(h)^{w_{23} w_3} \rismin{2}{3}(d)^{w_3} \risshpos{3}(k) \\
		&= w_{23} w_3 w_{23} w_3 \risplus{2}{3}(-c)^{w_{23} w_3} \risshpos{2}(-h)^{w_3} \risplus{2}{3}(-d) \risshpos{3}(k) \\
		&= w_{23} w_3 w_{23} w_3 \rismin{2}{3}(c) \risshpos{2}(-\modinv{h}) \risplus{2}{3}(-d) \risshpos{3}(k).
	\end{align*}
	In order to compare these two terms, we have to change the order of the product in the first one, using the commutator formulas in \cref{B:comm-formula-firststep} as well as \thmitemcref{group-rel}{group-rel:comm}:
	\begin{align*}
		\risshpos{3}(-\modinv{v}) \risplus{2}{3}(-a) \risshpos{2}(-u) \rismin{2}{3}(b) \hspace{-5cm}& \\
		&= \risshpos{3}(-\modinv{v}) \rismin{2}{3}(b) \risplus{2}{3}(-a) \risshpos{2}(-u) \\
		&= \rismin{2}{3}(b) \risshpos{3}(-\modinv{v}) \commutator[\big]{\risshpos{3}(-\modinv{v})}{\rismin{2}{3}(b)} \risplus{2}{3}(-a) \risshpos{2}(-u) \\
		&= \rismin{2}{3}(b) \risshpos{3}(-\modinv{v}) \risshpos{2}\brackets[\big]{-g(b, -\modinv{v})} \risplus{2}{3}\brackets[\big]{q(b, -\modinv{v})} \risplus{2}{3}(-a) \risshpos{2}(-u) \\
		&= \rismin{2}{3}(b) \risshpos{2}\brackets[\big]{g(b, \modinv{v})-u} \risplus{2}{3}\brackets[\big]{q(b, -\modinv{v}) -a} \risshpos{3}(-\modinv{v}) \\
		& \hspace{2cm} \mathord{} \cdot \commutator[\big]{\risshpos{3}(-\modinv{v})}{\risshpos{2}\brackets[\big]{g(b, \modinv{v})-u}} \\
		&= \rismin{2}{3}(b) \risshpos{2}\brackets[\big]{g(b, \modinv{v})-u} \risplus{2}{3}\brackets[\big]{q(b, -\modinv{v}) -a} \risshpos{3}(-\modinv{v}) \\
		& \hspace{2cm} \mathord{} \cdot \risplus{2}{3}\brackets[\big]{f\brackets[\big]{g(b, \modinv{v})-u, -\modinv{v}}} \\
		&= \rismin{2}{3}(b) \risshpos{2}\brackets[\big]{g(b, \modinv{v})-u} \\
		& \hspace{2cm} \mathord{} \cdot \risplus{2}{3}\brackets[\bigg]{q(b, -\modinv{v}) -a +f\brackets[\big]{u-g(b, \modinv{v}), \modinv{v}}} \risshpos{3}(-\modinv{v}).
	\end{align*}
	Since $ w_3 w_{23} w_3 w_{23} = w_{23} w_3 w_{23} w_3 $ by \cref{braid:all} and since $ G $ is rank-2-injective, we conclude that $ \blumapG{\word{\alpha}}(x) = \blumapG{\word{\beta}}(x') $ if and only if
	\begin{gather*}
		c=b, \qquad -\modinv{h} = g(b, \modinv{v}) - u, \\
		-d = q(b, -\modinv{v}) - a + f\brackets[\big]{u - g(b, \modinv{v}), \modinv{v}} \midand k=-\modinv{v}.
	\end{gather*}
	In other words, $ \blumapG{\word{\alpha}}(x) = \blumapG{\word{\beta}}(x') $ if and only if $ x' = \blutrans{23}(x) $, which shows that $ \blutrans{23} $ is a blueprint rewriting rule. This finishes the proof.
\end{proof}


\section{The Blueprint Computation}

\label{sec:B-bluecomp}

\begin{secnotation}
	\cref{B:blue-secnotation} continues to hold.
\end{secnotation}

Finally, we have everything ready to perform the blueprint computation. This computation takes place within the $ B_3 $-subsystem which is spanned by $ \rootbase' $ from \cref{B:blue-secnotation}, and there are no new identities which could be obtained by performing the computation in larger subsystems.

\begin{remark}[Blueprint computation]\label{B:blue-comp}
	A homotopy cycle of the longest word in the Weyl group of $ B_3 $ is given in \cref{fig:B:hom-cycle}. We begin with the tuple
	\[ \brackets[\big]{\risom{3}(u), \risom{2}(a), \risom{1}(b), \risom{3}(v), \risom{2}(c), \risom{3}(w), \risom{2}(d), \risom{1}(r), \risom{2}(s)} \]
	where $ \risom{1} \defl \rismin{1}{2} $, $ \risom{2} \defl \rismin{2}{3} $, $ \risom{3} \defl \risshpos{3} $ and where $ u,v,w \in \module $ and $ a,b,c,d,r,s \in \comring $ are arbitrary. Working down rows~1 to~12 in the homotopy cycle and applying the respective blueprint rewriting rules in the process, we obtain a tuple
	\[ \brackets[\big]{\risom{2}(x_1), \risom{1}(x_2), \risom{2}(x_3), \risom{3}(x_4), \risom{2}(x_5), \risom{3}(x_6), \risom{1}(x_7), \risom{2}(x_8), \risom{3}(x_9)}. \]
	Conversely, working up from row~23 to row~12, we obtain a tuple
	\[ \brackets[\big]{\risom{2}(y_1), \risom{1}(y_2), \risom{2}(y_3), \risom{3}(y_4), \risom{2}(y_5), \risom{3}(y_6), \risom{1}(y_7), \risom{2}(y_8), \risom{3}(y_9)}. \]
	Here $ x_i, y_i \in \module $ for all $ i \in \Set{4, 6, 9} $ and $ x_i, y_i \in \comring $ for all $ i \in \Set{1,2,3,5,7,8} $. Now $ x_i = y_i $ for all $ i \in \numint{1}{9} $ by \cref{blue:thm}. The precise results of this computation can be found in \cref{fig:B:blue-short,fig:B:blue-mid,fig:B:blue-long}. The intermediate steps of the computation have been performed with GAP \cite{GAP4}.
\end{remark}

\premidfigure
\begin{figure}[htb]
	\centering\begin{tabular}{rc@{\hspace{2cm}}rc}
		(1) & 321\underline{3232}12 & (13) & 21323\underline{212}3 \\
		(2) & 321232\underline{31}2 & (14) & 2\underline{13}2\underline{31}213 \\
		(3) & 3\underline{212}32132 & (15) & 23\underline{121}3213 \\
		(4) & \underline{31}2\underline{13}2132 & (16) & 2321232\underline{13} \\
		(5) & 1323\underline{121}32 & (17) & 2321\underline{2323}1 \\
		(6) & 1\underline{3232}1232 & (18) & 232\underline{13}2321 \\
		(7) & 1232\underline{31}232 & (19) & \underline{2323}12321 \\
		(8) & 12321\underline{3232} & (20) & 323\underline{212}321 \\
		(9) & 123\underline{212}323 & (21) & 32\underline{31}2\underline{13}21 \\
		(10) & 12\underline{31}2\underline{13}23 & (22) & 321323\underline{121} \\
		(11) & \underline{121}323123 & (23) & 321323212 \\
		(12) & 21\underline{2323}123 & &
	\end{tabular}
	\caption{A homotopy cycle of the longest word in $ \Weyl(B_3) $ (taken from \cite[Fig.~8, \onpage{472}]{MoufangPolygons}).}
	\label{fig:B:hom-cycle}
\end{figure}

\begin{figure}[htb]
	\centering$ \begin{gathered}
		(1) \; s=s, \qquad (2) \; -r-sd=-r-ds, \qquad (3) \; d=d, \\
		(6) \; -v + \modinv{g(s,u)} = -v + g(s, \modinv{u}), \qquad (9) \; u=u.
	\end{gathered} $
	\caption{The equations \enquote{$ (i) \; x_i = y_i $} for $ i \in \Set{1,2,3,6,9} $ in \cref{B:blue-comp}.}
	\label{fig:B:blue-short}
\end{figure}

\begin{figure}[htb]
	\centering$ \begin{gathered}
		w-g\brackets{d,\overline{v}}-\overline{g\brackets{r,\overline{u}}} = w-g\brackets{r,u}-g\brackets{ds,u}-\overline{g\brackets{d,v}}+\overline{g\brackets[\big]{d,g\brackets{s,\overline{u}}}} \\
		a-q\brackets{s,u}-f\brackets{\overline{v},u}+f\brackets[\big]{g\brackets{s,u},u} = a-q\brackets{s,-\overline{u}}-f\brackets{v,\overline{u}}+f\brackets[\big]{g\brackets{s,\overline{u}},\overline{u}}
	\end{gathered} $
	\caption{The equations $ x_4 = y_4 $ and $ x_8 = y_8 $ in \cref{B:blue-comp}, respectively.}
	\label{fig:B:blue-mid}
\end{figure}

\begin{figure}[htb]
	\centering$ \begin{aligned}
		x_5 &= c-q\brackets{d,-\overline{v}}-f\brackets{w,\overline{v}}+f\brackets{g\brackets{d,\overline{v}},\overline{v}}+ra\\
		& \hspace{.75cm} \mathord{}-\brackets[\bigg]{-b-ad-q\brackets{r,-\overline{u}}-f\brackets{\overline{w},\overline{u}}+f\brackets[\big]{\overline{g\brackets{d,\overline{v}}},\overline{u}}+f\brackets[\big]{g\brackets{r,\overline{u}},\overline{u}}}s, \\
		y_5 &= c+sb-\brackets[\big]{a-q\brackets{s,-\overline{u}}-f\brackets{v,\overline{u}}+f\brackets{g\brackets{s,\overline{u}},\overline{u}}}\brackets{-r-ds}\\
		&\hspace{.75cm}\mathord{}-q\brackets[\big]{d,-v+g\brackets{s,\overline{u}}}-f\brackets{\overline{w},v}+f\brackets{\overline{w},g\brackets{s,\overline{u}}}+f\brackets[\big]{\overline{g\brackets{r,u}},v} \\
		&\hspace{.75cm}\mathord{}-f\brackets[\big]{\overline{g\brackets{r,u}},g\brackets{s,\overline{u}}}+f\brackets[\big]{\overline{g\brackets{ds,u}},v}-f\brackets[\big]{\overline{g\brackets{ds,u}},g\brackets{s,\overline{u}}}\\
		&\hspace{0.75cm}\mathord{}+f\brackets[\big]{g\brackets{d,v},v}-f\brackets{g\brackets{d,v},g\brackets{s,\overline{u}}}-f\brackets[\big]{g\brackets{d,g\brackets{s,\overline{u}}},v} \\
		&\hspace{0.75cm}\mathord{} +f\brackets[\big]{g\brackets{d,g\brackets{s,\overline{u}}},g\brackets{s,\overline{u}}}, \\
		x_7 &= -b-ad-q\brackets{r,-\overline{u}}-f\brackets{\overline{w},\overline{u}}+f\brackets[\big]{\overline{g\brackets{d,\overline{v}}},\overline{u}}+f\brackets[\big]{g\brackets{r,\overline{u}},\overline{u}}, \\
		y_7 &= -b-q\brackets{r,u}-q\brackets{ds,u}-f\brackets{w,u}+f\brackets[\big]{g\brackets{r,u},u}+f\brackets[\big]{g\brackets{ds,u},u}\\
		&\hspace{1cm} \mathord{}-d\brackets[\big]{a-q\brackets{s,-\overline{u}}-f\brackets{v,\overline{u}}+f\brackets[\big]{g\brackets{s,\overline{u}},\overline{u}}}.
	\end{aligned} $
	\caption{The values of $ x_5, y_5, x_7, y_7 $ in \cref{B:blue-comp}.}
	\label{fig:B:blue-long}
\end{figure}
\postmidfigure

\begin{note}
	Throughout this section, we denote by $ u,v,w $ arbitrary elements of $ \module $ and by $ a $, $ b $, $ c $, $ d $, $ r $, $ s $ arbitrary elements of $ \comring $.
\end{note}

Since all identities in \cref{fig:B:blue-short,fig:B:blue-mid,fig:B:blue-long} hold for arbitrary values of the variables, they hold in particular if we replace some variables by $ 0 $, $ 1_\comring $ or $ v_0 $. This allows us to derive shorter identities. We can then use these shorter identities to simplify the original identities. We continue this process until we have proven that the set of identities in \cref{fig:B:blue-short,fig:B:blue-mid,fig:B:blue-long} is equivalent to the axioms of a quadratic module $ (\module, \map{}{}{}{v}{q(1_\comring, v)}) $ over a commutative ring $ \comring $ with linearisation $ f $ and scalar multiplication $ g $. The map $ \map{}{}{}{v}{\modinv{v}} $ will turn out to be the reflection~$ \refl{v_0} $.

\begin{lemma}[Equations 2, 6]\label{B:commutative}
	The ring $ \comring $ is commutative, and we have $ g(s, \modinv{u}) = \modinv{g(s,u)} $ for all $ s \in \comring $, $ u \in \module $.
\end{lemma}
\begin{proof}
	The first equation follows from equation~2 with $ r=0 $, the second one from equation~6 with $ v=0 $ (see \cref{fig:B:blue-short}).
\end{proof}

\cref{B:commutative} clearly covers all non-trivial identities which can be derived from \cref{fig:B:blue-short}.

\begin{lemma}[Equation 4]\label{B:eq4lem}
	We have $ g(ds, u) = g\brackets[\big]{d, g(s,u)} $ for all $ d,s \in \comring $ and $ u \in \module $.
\end{lemma}
\begin{proof}
	Using the relation $ g(s, \modinv{u}) = \modinv{g(s,u)} $ from \cref{B:commutative} and the fact that the involution on $ \module $ is of order at most~2 by \cref{B:invo-lem}, we see that both sides of equation~4 in \cref{fig:B:blue-mid} contain the term $ w - \modinv{g(d,v)} - g(r,u) $. Removing these terms from the equation, we are left with $ g(ds,u) = \modinv{g(d, g(s, \modinv{u}))} $, which by the same considerations as before yields the desired relation.
\end{proof}

Again, it is clear that no more identities can be deduced from equation~4.

\begin{lemma}[Equation 8]\label{B:eq8lem}The following hold:
	\begin{lemenumerate}
		\item \label{B:eq8lem:f}$ f(\modinv{v}, u) = f(v, \modinv{u}) $ and $ f(\modinv{v}, \modinv{u}) = f(v,u) $ for all $ u,v \in \module $.
		
		\item \label{B:eq8lem:q}$ q(s, -\modinv{u}) = q(s,u) $ for all $ s \in \comring $, $ u \in \module $.
	\end{lemenumerate}
\end{lemma}
\begin{proof}
	Putting $ a \defl s \defl 0_\comring $ in equation~8 in \cref{fig:B:blue-mid}, we see that $ f(\modinv{v}, u) = f(v, \modinv{u}) $, which implies that $ f(\modinv{v}, \modinv{u}) = f(v, \modinv{\modinv{u}}) = f(v,u) $. Using \cref{B:commutative}, we infer that
	\[ f\brackets[\big]{g(s, \modinv{u}), \modinv{u}} = f\brackets[\big]{\modinv{g(s,u)}, \modinv{u}}. \]
	Thus we can cancel  nearly all terms from equation~8, and what remains is $ q(s, -\modinv{u}) = q(s,u) $. This finishes the proof.
\end{proof}

Only equations~5 and~7 remain. We begin with equation~7 because it is shorter, except that we quickly derive the following fact from equation~5 to simplify our notation.

\begin{lemma}
	The ring $ \comring $ is associative.
\end{lemma}
\begin{proof}
	Putting all variables except for $ a,d,s $ in equation~5 in \cref{fig:B:blue-long} to zero, we obtain that $ (ad)s = a(ds) $.
\end{proof}

\begin{remark}[Equation 7, simplification]\label{B:eq7-simp}
	With the knowledge obtained so far, we can cancel many terms in equation~7 in \cref{fig:B:blue-long}, which yields the following simplified equation:
	\begin{align*}
		f\brackets[\big]{\modinv{g(d, \modinv{v})}, \modinv{u}} &= -q(ds, u) + f\brackets[\big]{g(ds, u), u} + dq(s, -\modinv{u})  \\
		&\hspace{3cm} \mathord{}+ df(v, \modinv{u})- df\brackets[\big]{g(s, \modinv{u}), \modinv{u}}.
	\end{align*}
	By \cref{B:commutative,B:eq8lem}, this is equivalent to the following equation:
	\begin{align*}
		f\brackets[\big]{g(d,v), \modinv{u}} &= -q(ds, u) + f\brackets[\big]{g(ds, u), u} + dq(s, u)  \\
		&\hspace{3cm} \mathord{}+ df(v, \modinv{u})- df\brackets[\big]{g(s, u), u}.
	\end{align*}
\end{remark}

\begin{lemma}[Equation 7]\label{B:eq7-lem}
	The following hold:
	\begin{lemenumerate}
		\item \label{B:eq7-lem:f}$ f\brackets[\big]{g(d,v), u} = d f(v,u) $ for all $ d \in \comring $ and all $ u,v \in \module $.
		
		\item $ q(d,u) = d q(1_\comring, u) $ for all $ d \in \comring $ and all $ u \in \module $.
	\end{lemenumerate}
\end{lemma}
\begin{proof}
	Putting $ s \defl 0 $ in equation~7 in \cref{B:eq7-simp}, we see that $ f\brackets[\big]{g(d,v), \modinv{u}} = d f(v, \modinv{u}) $. This proves the first assertion because the involution on $ \module $ is a bijection. Using the first assertion and the associativity of $ \comring $, we see that
	\[ f\brackets[\big]{g(ds, u), u} = (ds) f(u,u) = d \brackets[\big]{s f(u,u)} = d f\brackets[\big]{g(s,u), u}. \]
	Thus we can now simplify equation~7 to obtain $ q(ds,u) = dq(s,u) $, from which the second assertion follows by putting $ s \defl 1_\comring $.
\end{proof}

\begin{remark}
	As a variation of \thmitemcref{B:eq7-lem}{B:eq7-lem:f}, we also have $ f\brackets[\big]{\modinv{g(d,v)}, u} = df(\modinv{v}, u) $ because, by \thmitemcref{B:eq8lem}{B:eq8lem:f},
	\begin{align*}
		f\brackets[\big]{\modinv{g(d,v)}, u} &= f\brackets[\big]{g(d,v), \modinv{u}} = d f(v, \modinv{u}) = d f(\modinv{v}, u).
	\end{align*}
\end{remark}

\begin{remark}
	At first glance, the equation $ q(ds,u) = dq(s,u) $ that we have proven in \cref{B:eq7-lem} may seem stronger than the equation $ q(d,u) = d q(1_\comring, u) $. However, the first equation follows from the second one together with the associativity of $ \comring $:
	\[ q(ds,u) = (ds) q(1_\comring, u) = d \brackets[\big]{sq(1_\comring, u)} = d q(s,u). \]
	Thus we do not lose any information by only recording the second equation in \cref{B:eq7-lem}.
\end{remark}

\begin{notation}
	From now on, we put $ q(v) \defl q(1,v) $ for all $ v \in \module $.
\end{notation}

\begin{remark}[Equation 5, simplification, part 1]\label{B:eq5-simp}
	We are now left with only equation~5 in \cref{fig:B:blue-long}, the most complicated one. Using the previously derived identities, we can simplify this equation. First of all, note that the right-hand side of this equation contains the terms
	\begin{gather*}
		f\brackets[\big]{\modinv{g(ds, u)}, v} -f\brackets[\big]{g\brackets{d,g\brackets{s,\overline{u}}},v} \rightand \\
		-f\brackets[\big]{\overline{g\brackets{ds,u}},g\brackets{s,\overline{u}}} +f\brackets[\big]{g\brackets{d,g\brackets{s,\overline{u}}},g\brackets{s,\overline{u}}},
	\end{gather*}
	both of which are zero. Further, observe that the left-hand side contains
	\[ c - f(w, \modinv{v}) + bs + f\brackets[\big]{g(d, \modinv{v}), \modinv{v}} \]
	while the right-hand side contains
	\[ c - f(\modinv{w}, v) + sb + f\brackets[\big]{g(d,v), v}. \]
	Since these terms are equal, we can cancel them. What remains are the following expressions for $ x_5 $ and $ y_5 $: 
	\begin{align*}
		x_5 &= -q\brackets{d,-\overline{v}}+ra\\
		& \hspace{.75cm} \mathord{}-\brackets[\bigg]{-ad-q\brackets{r,-\overline{u}}-f\brackets{\overline{w},\overline{u}}+f\brackets[\big]{\overline{g\brackets{d,\overline{v}}},\overline{u}}+f\brackets[\big]{g\brackets{r,\overline{u}},\overline{u}}}s, \\
		y_5 &= -\brackets[\big]{a-q\brackets{s,-\overline{u}}-f\brackets{v,\overline{u}}+f\brackets{g\brackets{s,\overline{u}},\overline{u}}}\brackets{-r-ds}\\
		&\hspace{.75cm}\mathord{}-q\brackets[\big]{d,-v+g\brackets{s,\overline{u}}}+f\brackets{\overline{w},g\brackets{s,\overline{u}}}+f\brackets[\big]{\overline{g\brackets{r,u}},v} \\
		&\hspace{.75cm}\mathord{}-f\brackets[\big]{\overline{g\brackets{r,u}},g\brackets{s,\overline{u}}} -f\brackets{g\brackets{d,v},g\brackets{s,\overline{u}}} 
	\end{align*}
	Applying the distributive law and the previously established identities, we see that
	\begin{align*}
		x_5 &= - dq(v) + ra + ads +rsq(u) + sf(w,u) - ds f(v, \modinv{u}) - rs f(u,u), \\
		y_5 &=ar + ads - rs q(u) -ds^2 q(u) -rf(v, \modinv{u}) - dsf(v, \modinv{u}) + rs f(u,u) \\
		&\hspace{0.75cm} \mathord{}+ ds^2 f(u,u) -d q\brackets[\big]{g(s,\modinv{u}) - v} + f\brackets[\big]{\modinv{w}, g(s, \modinv{u})} + rf\brackets{\modinv{u}, v} - rsf(u,u) \\
		&\hspace{.75cm} \mathord{}- ds f(v, \modinv{u}).
	\end{align*}
	Clearly, the term
	\[ ra +ads - dsf(v, \modinv{u}) - rsf(u,u) \]
	on the left-hand side cancels the term
	\[ ar+ads - dsf(v, \modinv{u}) - rsf(u,u) \]
	on the right-hand side, so we are left with
	\begin{align*}
		x_5 &= - dq(v)  +rsq(u) + sf(w,u), \\
		y_5 &= - rs q(u) -ds^2 q(u) -rf(v, \modinv{u}) + rsf(u,u) + ds^2 f(u,u)  \\
		&\hspace{0.75cm} \mathord{} -d q\brackets[\big]{g(s,\modinv{u}) - v} + f\brackets[\big]{\modinv{w}, g(s, \modinv{u})} + rf\brackets{\modinv{u}, v} - ds f(v, \modinv{u}).
	\end{align*}
\end{remark}

\begin{lemma}[Equation 5, part 1]\label{B:e5lem-1}
	The following hold:
	\begin{lemenumerate}
		\item $ f(u,v) = f(v,u) $ for all $ u,v \in \module $.
		
		\item $ f\brackets[\big]{v, g(d,u)} = df(v,u) $ for all $ u,v \in \module $, $ d \in \comring $.
	\end{lemenumerate}
\end{lemma}
\begin{proof}
	Putting $ s \defl 0 $, $ d \defl 0 $, $ r \defl 1 $ in equation~5 (see \cref{B:eq5-simp}), we obtain $ 0 = -f(v, \modinv{u}) + f(\modinv{u}, v) $. Since the involution on $ \module $ is a bijection, this implies the first assertion. Now the second assertion together with \thmitemcref{B:eq7-lem}{B:eq7-lem:f} yields that
	\[ f\brackets[\big]{v, g(d,u)} = f\brackets[\big]{g(d,u), v} = df(u,v) = df(v,u), \]
	as desired.
\end{proof}

\begin{remark}[Equation 5, simplification, part 2]\label{B:eq5-simp2}
	\cref{B:e5lem-1} allows us to simplify equation~5 even further: The term $ sf(w,u) $ on the left-hand side cancels the term $ f\brackets[\big]{\modinv{w}, g(s, \modinv{u})} $ on the right-hand side. Further, we see that the term $ -rf(v, \modinv{u}) + rf(\modinv{u}, v) $ on the right-hand side equals zero. Thus we obtain
	\begin{align*}
		x_5 &= - dq(v)  +rsq(u), \\
		y_5 &= - rs q(u) -ds^2 q(u) + rsf(u,u) + ds^2 f(u,u) -d q\brackets[\big]{g(s,\modinv{u}) - v} - ds f(v, \modinv{u}).
	\end{align*}
\end{remark}

\begin{lemma}[Equation 5, part 2]\label{B:eq5lem-2}
	The following hold:
	\begin{lemenumerate}
		\item $ f(u,u) = 2q(u) $ for all $ u \in \module $.
		\item \label{B:eq5lem-2:q}$ q(-v) = q(\modinv{v}) = q(v) $ for all $ v \in \module $.
	\end{lemenumerate}
\end{lemma}
\begin{proof}
	Putting $ r \defl 1 $, $ s \defl 1 $, $ d \defl 0 $ and $ v \defl 0 $ in equation~5 in \cref{B:eq5-simp2}, we see that $ q(u) = -q(u) + f(u,u) $, which implies the first assertion. Using this identity, we obtain the following simplified version of equation~5:
	\[ -dq(v) = -ds^2 q(u) - dq\brackets[\big]{g(s, \modinv{u}) -v} - ds f(v, \modinv{u}). \]
	Putting $ u \defl 0 $ and $ d \defl 1 $, we see that $ q(v) = q(-v) $. Since we know that $ q(-\modinv{v}) = q(v) $ by \thmitemcref{B:eq8lem}{B:eq8lem:q}, the second assertion follows.
\end{proof}

\begin{remark}
	It is not immediately obvious that no further relations can be deduced from equation~5. However, using the rank-2-computation from \cref{B:param:rank2commrel}, we see that
	\begin{align*}
		dq\brackets[\big]{g(s, \modinv{u}) -v} &= dq\brackets[\big]{g(s, \modinv{u})} + dq(-v) + df\brackets[\big]{g(s, \modinv{u}), -v} \\
		&= ds^2 q(u) +dq(v) -dsf(\modinv{u}, v).
	\end{align*}
	Thus the simplified version of equation~5 in the proof of \cref{B:eq5lem-2} reduces to $ 0=0 $. It is noteworthy that, up to this point, we have never used the non-linear identity from \cref{B:param:rank2commrel} in our computations. This shows that, even if we had not proven this identity in our rank-2-computations, we would have derived it from the blueprint computations.
\end{remark}

\begin{summary}\label{B:blue:summary}
	The following hold:
	\begin{proenumerate}
		\item The ring $ \comring $ is associative and commutative.
		
		\item The map $ \map{g}{\comring \times \module}{\module}{}{} $ defines a $ \comring $-scalar multiplication in $ \module $ which turns it into a $ \comring $-module.
		
		\item The involution on $ \module $ is a $ \comring $-linear automorphism of $ \module $ with $ \modinv{\modinv{v}} = v $ and $ q(\modinv{v}) = q(v) $ for all $ v \in \module $.
		
		\item The map $ \map{f}{\module \times \module}{\module}{}{} $ is a symmetric $ \comring $-bilinear form on $ \module $ and the map $ \map{q}{\module}{\comring}{v}{q(1_\comring, v)} $ is a $ \comring $-quadratic form on $ \module $ with linearisation $ f $.
	\end{proenumerate}
\end{summary}

\begin{proposition}\label{B:blue:stand-signs}
	$ G $ satisfies the same commutator relations as the group $ \EO(q) $ in \cref{B:ex-commrel}. In other words, $ (\risom{\alpha})_{\alpha \in B_n} $ is a parametrisation of $ G $ by $ (\module, q) $ with standard signs.
\end{proposition}
\begin{proof}
	The third set of relations (concerning commutators of medium-length roots) is satisfied by \cref{B:comm-mult-computation}. Now let $ i<j \in \numint{1}{n} $, let $ r \in \ring $ and let $ v \in \module $. Throughout this computation, we will frequently apply \thmitemcref{B:eq8lem}{B:eq8lem:f} and \thmitemcref{B:eq5lem-2}{B:eq5lem-2:q}. We know from \cref{B:comm-formula-firststep} that
	\begin{align*}
		\commutator{\rismin{i}{j}(r)}{\risshpos{j}(v)} &= \risshpos{i}\brackets{rv} \risplus{i}{j}\brackets[\big]{-rq(v)}.
	\end{align*}
	Conjugating this equation by $ w_{ij} $, we infer that
	\begin{align*}
		\commutator{\rismin{j}{i}(-r)}{\risshpos{i}(-v)} &= \risshpos{j}(rv) \risplus{i}{j}\brackets[\big]{-rq(v)}.
	\end{align*}
	In other words,
	\begin{align*}
		\commutator{\rismin{j}{i}(r)}{\risshpos{i}(v)} &= \risshpos{j}(rv) \risplus{i}{j}\brackets[\big]{rq(v)}.
	\end{align*}
	This says precisely that for all distinct $ i,j \in \numint{1}{n} $, we have
	\begin{equation}\label{eq:B:long-short-commrel}
		\commutator{\rismin{i}{j}(r)}{\risshpos{j}(v)} = \risshpos{i}(rv) \risplus{i}{j}\brackets[\big]{\delmin{i<j} r q(v)},
	\end{equation}
	which is the first relation in \cref{B:ex-commrel}.
	
	Continue to assume that $ r \in \ring $, $ v \in \module $ and that $ i,j \in \numint{1}{n} $ are distinct (but not necessarily that $ i<j $). Conjugating~\eqref{eq:B:long-short-commrel} by $ w_j $, we see that
	\begin{align*}
		\commutator[\big]{\risplus{i}{j}(\delmin{j>i}r)}{\risshneg{j}\brackets[\big]{\modinv{v}}} &= \risshpos{i}\brackets[\big]{\modinv{rv}} \rismin{i}{j}\brackets[\big]{\delmin{j>i}\delmin{i<j} r q(v)}.
	\end{align*}
	In other words,
	\begin{align*}
		\commutator{\risplus{i}{j}(r)}{\risshneg{j}(v)} &= \risshpos{i}\brackets{\delmin{j>i}rv} \rismin{i}{j}\brackets[\big]{\delmin{j>i} r q(v)}.
	\end{align*}
	This is the third relation in \cref{B:ex-commrel}. Conjugating it by $ w_i $, we obtain
	\begin{align*}
		\commutator[\big]{\rismin{j}{i}(\delmin{i>j}r)}{\risshneg{j}\brackets[\big]{\modinv{v}}} &= \risshneg{i}\brackets[\big]{\delmin{j>i}\modinv{rv}} \risminmin{i}{j}\brackets[\big]{\delmin{i>j}\delmin{j>i} r q(v)}.
	\end{align*}
	This says that
	\begin{align*}
		\commutator[\big]{\rismin{j}{i}(r)}{\risshneg{j}\brackets{v}} &= \risshneg{i}\brackets{-rv} \risminmin{i}{j}\brackets[\big]{\delmin{j>i} r q(v)},
	\end{align*}
	which is the second relation in \cref{B:ex-commrel}. Now we conjugate~\eqref{eq:B:long-short-commrel} by $ w_i $. This yields
	\begin{equation*}
		\commutator{\risminmin{i}{j}(\delmin{i>j} r)}{\risshpos{j}\brackets{\modinv{v}}} = \risshneg{i}\brackets[\big]{\modinv{rv}} \rismin{j}{i}\brackets[\big]{\delmin{i>j}\delmin{i<j} r q(v)}.
	\end{equation*}
	Hence
	\begin{align*}
		\commutator{\risminmin{i}{j}(r)}{\risshpos{j}\brackets{v}} = \risshneg{i}\brackets{\delmin{i>j}rv} \rismin{j}{i}\brackets[\big]{\delmin{i<j} r q(v)},
	\end{align*}
	which is the fourth relation in \cref{B:ex-commrel}. This finishes the proof of the first set of relations in \cref{B:ex-commrel}.
	
	Now let $ i,j \in \numint{1}{n} $ be distinct and let $ u,v \in \module $. We know from \cref{B:comm-formula-firststep} that
	\begin{equation*}
		\commutator{\risshpos{i}(v)}{\risshpos{j}(u)} = \risplus{i}{j}\brackets[\big]{-f(v,u)}
	\end{equation*}
	if $ i<j $. In this case, conjugating this equation by $ w_{ij} $ yields
	\begin{align*}
		\commutator{\risshpos{j}(v)}{\risshpos{i}(-u)} = \risplus{i}{j}\brackets[\big]{-f(v,u)}.
	\end{align*}
	We conclude that
	\begin{equation*}
		\commutator{\risshpos{i}(v)}{\risshpos{j}(u)} = \risplus{i}{j}\brackets[\big]{\delmin{i<j}f(v,u)}
	\end{equation*}
	holds for all distinct $ i,j \in \numint{1}{n} $. This is the fifth relation in \cref{B:ex-commrel}. We can conjugate it by $ w_i $ to obtain
	\begin{align*}
		\commutator{\risshneg{i}(\modinv{v})}{\risshpos{j}(\modinv{u})} = \rismin{j}{i}\brackets[\big]{\delmin{i>j}\delmin{i<j}f(v,u)},
	\end{align*}
	which means that
	\begin{equation}\label{eq:B:short-short-comm}
		\commutator{\risshneg{i}(v)}{\risshpos{j}(u)} = \rismin{j}{i}\brackets[\big]{-f(v,u)}.
	\end{equation}
	Since $ \commutator{\risshneg{i}(v)}{\risshpos{j}(u)} = \commutator{\risshpos{j}(u)}{\risshneg{i}(v)}^{-1} $ and $ f $ is symmetric, this is the sixth relation in \cref{B:ex-commrel}. Finally, we can conjugate~\eqref{eq:B:short-short-comm} by $ w_j $ to obtain
	\begin{align*}
		\commutator{\risshneg{i}(\modinv{v})}{\risshneg{j}(\modinv{u})} = \risminmin{j}{i}\brackets[\big]{-\delmin{j>i}f(v,u)}.
	\end{align*}
	In other words,
	\begin{align*}
		\commutator{\risshneg{i}(v)}{\risshneg{j}(u)} = \risminmin{j}{i}\brackets[\big]{\delmin{j<i}f(v,u)},
	\end{align*}
	which is the seventh relation in \cref{B:ex-commrel}. This finishes the proof.
\end{proof}

Before we state the main result of this chapter, we make a final observation.

\begin{lemma}\label{B:blue:short-weyl}
	We have $ v_{-1} = v_1 = -v_0 $ and $ q(v_0) = 1 $ as well as $ \modinv{u} = \refl{v_0}(u) $ for all $ u \in \module $. In particular, $ (\module, q, e) $ is a pointed quadratic module.
\end{lemma}
\begin{proof}
	As in the proof of \cref{B:short-weakly-balanced}, we can easily compute that
	\begin{align*}
		\rismin{1}{2}(1)^{w_2} = \risplus{1}{2}(-q(v_0)).
	\end{align*}
	On the other hand, by \cref{B:Bnsub-conj-anygroup}, we know that
	\[ \rismin{1}{2}(1)^{w_2} = \risplus{1}{2}(-1). \]
	 It follows that $ q(v_0) = 1 $. By \cref{B:short-weakly-balanced}, this implies that
	 \[ v_{-1} = v_1 = -q(v_0)^{-1} v_0 = -v_0. \]
	 Finally, we have shown in \cref{B:short-act-refl} for every $ u \in \module $ that $ \risshpos{1}(u)^{w_2} = \risshpos{1}(\refl{v_0}(u)) $, but we also know that $ \risshpos{1}(u)^{w_2} = \risshpos{1}(\modinv{u}) $. This finishes the proof.
\end{proof}

\begin{theorem}[Coordinatisation theorem for $ B_n $]\label{B:thm}
	Let $ G $ be a group with a crystallographic $ B_n $-grading $ (\rootgr{\alpha})_{\alpha \in B_n} $ for some $ n \in \IN_{\ge 3} $. Then there exist a commutative associative ring $ \comring $ and a pointable quadratic module $ (\module, q) $ over $ \comring $ such that $ G $ is coordinatised by $ (\module, q) $ with standard signs (in the sense of \cref{B:standard-param-def}). Further, if we fix a $ \rootbase $-system of Weyl elements in $ G $, then we can choose the root isomorphisms $ (\risom{\alpha})_{\alpha \in B_n} $ so that $ w_\alpha = \risom{-\alpha}(-1_\comring) \risom{\alpha}(1_\comring) \risom{-\alpha}(-1_\comring) $ for all long simple roots $ \alpha $ and $ w_\delta = -\risom{-\delta}(-v_0) \risom{\delta}(v_0) \risom{-\delta}(-v_0) $ for the short simple root $ \delta $ where $ v_0 $ is some element of $ \module $ with $ q(v_0) = 1 $.
\end{theorem}
\begin{proof}[Summary of the proof]
	Choose a root base $ \rootbase $ of $ B_n $ and a $ \rootbase $-system $ (w_\delta)_{\delta \in \rootbase} $ of Weyl elements. Denote by $ (\twistgroup, \inverparsym, \invogroup, \invoparsym) $ the standard admissible partial twisting system for $ G $, as in \cref{B:standard-partwist-def}. Then by \cref{B:param-exists}, there exist abelian groups $ (\comring, +) $ and $ (\module, +) $ on which $ \twistgroup \times \invogroup $ acts and a parametrisation $ (\risom{\alpha})_{\alpha \in \roots} $ of $ G $ by $ (\twistgroup \times \invogroup, \module, \comring) $ with respect to $ (w_\delta)_{\delta \in \rootbase} $ and $ \inverparsym \times \invoparsym $. We can define commutation maps as in \cref{B:commmap-def}. By \cref{B:blue:summary}, these maps equip $ \comring $ with the structure of a commutative associative ring and $ \module $ with a $ \comring $-module structure and a $ \comring $-quadratic form $ \map{q}{\module}{\comring}{}{} $. By \cref{B:blue:stand-signs}, $ (\risom{\alpha})_{\alpha \in \roots} $ is a coordinatisation of $ G $ by $ (\module, q) $ with standard signs. By \cref{B:isring,B:blue:short-weyl}, the Weyl elements $ (w_\delta)_{\delta \in \rootbase} $ have the desired form.
\end{proof}

	\chapter{Jordan Modules and Related Structures}
	
	\Cref{chap:simply-laced,chap:B} each started with a section which introduced the coordinatising algebraic structures of relevance in the respective chapter. For root gradings of types $ C $ and $ BC $, the corresponding algebraic structures are the so-called Jordan modules. The theory of these objects is rather involved, and we will need several other algebraic structures to define and investigate them properly. Further, Jordan modules are a new algebraic structure which we introduce specifically to describe root gradings of types $ C $ and $ BC $. For these reasons, we dedicate an own chapter to the study of these objects.
	
	\Cref{sec:weak-quad} begins with an introduction of weakly quadratic maps, which generalise quadratic maps. We will also introduce square-modules $ (\module, \joradd) $, which are groups equipped with a weakly quadratic scalar multiplication. They serve as the fundamental language in which we describe the other algebraic structures in this chapter.
	In \cref{sec:BC:altring,sec:rinv}, we list some known properties of alternative rings and (nonassociative) rings with involution, respectively. \Cref{sec:invset} is dedicated to involutory sets, which are alternative rings $ \ring $ equipped with an involution and a certain subset of $ \ring $. In \cref{sec:pseud}, we define pseudo-quadratic modules over involutory sets. Any pseudo-quadratic module $ \module $ gives rise to a group $ \psgr(\module) $ which is our main (and in fact, the only known) example of a Jordan module.
	
	In \cref{sec:jordan-modules}, we can finally define Jordan modules. A Jordan module is a pair of groups $ (\jormod, \ring) $ equipped with a family of maps satisfying exactly the identities which result from the blueprint computation for root gradings of type $ BC $. The group $ \ring $ in a Jordan module is always equipped with the structure of an alternative ring with involution. We will show in \cref{sec:jordan-class} that if $ 2_\ring $ is invertible, then every Jordan module is of the form $ \psgr(\module) $ for some pseudo-quadratic module $ \module $. Thus Jordan modules should be regarded as a slight generalisation of the class of groups $ \psgr(\module) $ which allows us to capture certain \enquote{characteristic-2 phenomena}.
	It should be noted that we are not aware of examples of Jordan modules which are not of the form $ \psgr(\module) $, but we cannot exclude their existence.
	
	In \cref{sec:pure-alt}, we will present a notion of \enquote{purely alternative rings} which was introduced by Slater in \cite{Slater_NuclCentAlt}. We will show that such rings do not admit non-zero modules. This fact will not formally be needed in the sequel. Rather, it illustrates that the non-associativity of the base ring in the definition of pseudo-quadratic modules should be seen as a technicality: Alternative rings admit pseudo-quadratic modules only insofar that they are not purely alternative.
	
	The present chapter is independent from \cref{chap:BC}, with a few exceptions: The material on pseudo-quadratic modules will be needed in the construction in \cref{sec:BC-example}, and some results on alternative rings will be used in \cref{sec:BC:blue-comp}. Hence the the enthusiastic reader may wish to skip ahead to \cref{chap:BC} and return to the current chapter only later.
	
	\label{chap:BC-alg}
	

\section{Weakly Quadratic Maps}

\label{sec:weak-quad}

\begin{note}
	Throughout this chapter, we will often consider groups $ (M, \joradd) $ which are not assumed to be commutative but which are still written additively because we think of them as \enquote{modules} and because, in some important cases, they turn out to be abelian. We indicate this non-commutativity by the hat on the operator. The neutral element of such groups will always be denoted by $ 0 $ and the inverse of an element $ x \in M $ will be denoted by $ \jormin x $. In particular, we have $ \jormin (x \joradd y) = \jormin y \jormin x $. We will later see that in root gradings of type $ C $, only the commutative case is relevant.
\end{note}

In this section, we introduce weakly quadratic maps and square-modules. Their theory provides the basic but essential language in which we describe Jordan modules and related algebraic structures. None of the results in this section are difficult, but we are not aware of any literature in which weakly quadratic maps appear.

Weakly quadratic maps generalise quadra\-tic maps (\cref{quadmod:polarisation-def,quadmod:quadmap-def}) in two ways: Firstly, the involved groups are allowed to be nonabelian. Secondly, we do not require Axiom~\thmitemref{quadmod:quadmap-def}{quadmod:quadmap-def:scalar}.

\begin{definition}[Weakly quadratic maps]\label{BC:weakquad-def}
	Let $ (M, \joradd), (N, \joradd) $ be groups and let $ \map{q}{M}{N}{}{} $ be a map. The \defemph*{polarisation of $ q $}\index{polarisation} is the map
	\[ \map{f}{M \times M}{N}{(v,w)}{\jormin q(w) \jormin q(v) \joradd q(v\joradd w)}. \]
	We say that $ q $ is \defemph*{weakly quadratic}\index{weakly quadratic map} if its polarisation is bi-additive, in which case the map $ f $ is also called the the \defemph*{linearisation of $ q $}\index{linearisation}.
\end{definition}

\begin{remark}\label{BC:weakly-quad-zero}
	Equivalently, the polarisation of a map $ \map{q}{M}{N}{}{} $ is the unique map $ \map{f}{M \times M}{N}{}{} $ which satisfies
	\[ q(v\joradd w) = q(v) \joradd q(w) \joradd f(v,w) \]
	for all $ v,w \in M $. In particular, if $ q $ is weakly quadratic, we have
	\[ q(0) = q(0 \joradd 0) = q(0) \joradd q(0) \joradd f(0,0) = q(0) \joradd q(0), \]
	so that $ q(0) = 0 $.
\end{remark}

\begin{example}
	Let $ M,N $ be groups. Clearly, any additive map $ \map{q}{M}{N}{}{} $ is weakly quadratic because its polarisation is the zero map. If $ M $ and $ N $ are abelian, then they are $ \IZ $-modules and any $ \IZ $-quadratic map $ \map{q}{M}{N}{}{} $ is weakly quadratic. Further, if $ N $ is abelian, then the sum of two weakly quadratic maps from $ M $ to $ N $ is weakly quadratic. Thus in the abelian setting, sums of quadratic and linear maps are examples of weakly quadratic. However, not every weakly quadratic map can be written in this way.
\end{example}

Our interest in weakly quadratic maps stems from the fact that in this chapter, we will frequently see groups which are \enquote{modules} over some ring except that the scalar multiplication is only weakly quadratic in the scalar. Such objects will be called square-modules.

\begin{definition}[Square-module]\label{BC:square-mod-def}
	Let $ \ring $ be a ring. A \defemph*{(right) square-module over $ \ring $}\index{square-module} is a group $ (N, \joradd) $ together with a map $ \map{\omega}{N \times \ring}{N}{}{} $ (called the \defemph*{square-scalar multiplication}\index{square-scalar multiplication}) satisfying the following properties:
	\begin{stenumerate}
		\item $ \omega $ is additive in the first component.
		
		\item \label{BC:square-mod-def:weak-quad}For all $ v \in N $, the map $ \map{}{\ring}{N}{r}{\omega(v,r)} $ is weakly quadratic.
		
		\item $ \omega(v,1_\ring) = v $ for all $ v \in N $.
	\end{stenumerate}
\end{definition}

\begin{definition}[Properties of square-modules]
	Let $ (N, \joradd) $ be a square-module over a ring $ \ring $ with square-scalar multiplication $ \omega $.
	\begin{defenumerate}
		\item We say that $ N $ is \defemph*{properly quadratic}\index{square-module!properly quadratic} if $ N $ is abelian and the map in Axiom~\thmitemref{BC:square-mod-def}{BC:square-mod-def:weak-quad} is $ \IZ $-quadratic in the sense of \cref{quadmod:quadmap-def}.
		
		\item We say that $ N $ is \defemph*{multiplicative}\index{square-module!multiplicative} if $ \omega\brackets[\big]{\omega(v,r),s} = \omega(v,rs) $ for all $ v \in N $ and $ r,s \in \ring $.
		
		\item A subgroup $ U $ of $ N $ is called a \defemph*{square-submodule}\index{square-module!submodule} if $ \omega(u,r) \in U $ for all $ u \in U $ and $ r \in \ring $.
	\end{defenumerate}
\end{definition}

\begin{remark}
	Let $ \ring $ be a ring. A square-module $ (\module, \omega) $ over $ \ring $ is an $ \ring $-module (in the regular sense of \cref{ring:module}) if and only if it is abelian, multiplicative and $ \omega $ is additive in the second component.
\end{remark}

The most important examples of square modules in our setting come from alternative rings with nuclear involutions, their involutory sets and pseudo-quadratic modules over these objects: see \cref{jormod:invring-jordan,rinv:sym-submod,rinv:nucl-submod,rinv:tr-submod,rinv:norm-tr-submod,pseud:T-squaremod}. All interesting examples of square-modules will be multiplicative. However, these examples often arise as square-submodules of (less interesting) non-multiplicative square-modules, which justifies the higher generality in our definition of square-modules. An example of a square-module which is not known to be multiplicative but which contains interesting multiplicative square-submodules is given in~\ref{jormod:invring-jordan}.

\begin{remark}\label{squaremod:zero-scalar}
	By \cref{BC:weakly-quad-zero}, the scalar multiplication on a square-module $ N $ satisfies $ \omega(v, 0_\ring) = 0_N $ for all $ v \in N $.
\end{remark}

\begin{definition}
	Let $ (N, \omega) $, $ (N', \omega') $ be two square-modules over a ring $ \ring $ and let $ \map{f}{N}{N'}{}{} $ be a map.
	\begin{defenumerate}
		\item We say that $ f $ \defemph*{preserves the square-scalar multiplication}\index{square-scalar multiplication!-preserving} if
		\[ f\brackets[\big]{\omega(v,r)} = \omega'\brackets[\big]{f(v), r} \]
		for all $ v \in N $ and all $ r \in \ring $.
		
		\item We say that $ f $ is a \defemph*{homomorphism of square-modules from $ N $ to $ N' $}\index{square-module!homomorphism} if it is a homomorphism of additive groups which preserves the scalar multiplication. It is called an \defemph*{isomorphism of square-modules} if it is a bijective homomorphism of square-modules, and it is called an \defemph*{automorphism of square-modules} if, in addition, $ (N, \omega) = (N', \omega') $.
	\end{defenumerate}
\end{definition}

We record the following standard facts from module theory which remain valid in this generality.

\begin{lemma}\label{jormod:ker-im}
	Kernels and images of homomorphisms of square-modules are square-submodules. Further, intersections of square-submodules are square-submodules.
\end{lemma}

\begin{lemma}\label{weakquad:quot-lem}\index{square-module!quotient module}
	Let $ (\module, \omega) $ be a square-module over some ring $ \ring $ and let $ U $ be a square-submodule of $ \module $. Denote by $ \map{\pi}{\module}{\module/U}{}{} $ the canonical projection, which is a homomorphism of groups. Then there exists a unique map
	\[ \map{\tilde{\omega}}{(\module/U) \times \ring}{\module/U}{}{} \]
	such that
	\[ \tilde{\omega}\brackets[\big]{\pi(x), r} = \pi\brackets[\big]{\omega(x,r)} \]
	for all $ x \in \module $ and $ r \in \ring $. Further, $ (\module/U, \tilde{\omega}) $ is a square-module, $ \pi $ is a homomorphism of square-modules and $ (\module/U, \tilde{\omega}) $ is multiplicative if $ (\module, \omega) $ is.
\end{lemma}

On any multiplicative square-module, we have a canonical \enquote{involution}.

\begin{definition}[Involution on square-modules]\label{weakquad:inv-def}
	Let $ (\module, \omega) $ be a multiplicative square-module over some ring $ \ring $. The map
	\[ \map{\modinvmap}{\module}{\module}{x}{\omega(x, -1_\ring)} \]
	is called the \defemph*{canonical involution on $ (\module, \omega) $}.\index{involution!on square-module}
\end{definition}

\begin{lemma}\label{weakquad:inv-lem}
	Let $ (\module, \omega) $ be a multiplicative square-module over some ring $ \ring $. Then the canonical involution on $ (\module, \omega) $ is an automorphism of the square-module $ (\module, \omega) $ which is its own inverse.
\end{lemma}
\begin{proof}
	For all $ x \in \module $, we have
	\[ \modinv{\modinv{x}} = \omega\brackets[\big]{\omega(x, -1_\ring), -1_\ring} = \omega\brackets[\big]{x, (-1_\ring) (-1_\ring)} = \omega(x,1) = x, \]
	so $ \modinvmap $ is its own inverse. Further, it is a homomorphism of groups because $ \omega $ is additive in the first component. Finally, we have for all $ x \in \module $ and $ r \in \ring $ that
	\begin{align*}
		\omega(\modinv{x}, r) &= \omega\brackets[\big]{\omega(x, -1_\ring), r} = \omega(x, -r) = \omega\brackets[\big]{\omega(x, r), -1_\ring} = \modinv{\omega(x,r)}.
	\end{align*}
	Thus $ \modinvmap $ is a homomorphism of square-modules, which finishes the proof.
\end{proof}


\section{Alternative Rings}

\label{sec:BC:altring}

\begin{secnotation}
	We denote by $ \ring $ an arbitrary ring.
\end{secnotation}

We have already introduced nonassociative rings and the basic related notions in \cref{sec:ring}. In this section, we record some elementary facts about the special case of alternative rings. We will also see some other related classes of nonassociative rings. A basic reference on this subject is \cite[Chapter~III]{Schafer}, but we will also cite various other sources.

\begin{definition}[Associativity properties, {\cite[Chapter~III]{Schafer}}]\label{ring:alternative-def}
	We say that $ \ring $ is \defemph*{alternative}\index{ring!alternative} if
	\[ x(xy) = (xx)y \midand (xy)y = x(yy) \]
	for all $ x,y \in \ring $. We say that $ \ring $ is \defemph*{weakly alternative}\index{ring!weakly alternative} if
	\[ \assoc{x}{y}{z} = -\assoc{y}{x}{z} \midand \assoc{x}{y}{z} = -\assoc{x}{z}{y} \]
	for all $ x,y,z \in \ring $. We say that $ \ring $ is \defemph*{flexible}\index{ring!flexible} if
	\[ (xy)x = x(yx) \]
	for all $ x,y \in \ring $. Further, for all $ a,x,y,z \in \ring $, the following relations are called the \defemph{Moufang identities}:
	\begin{align}
		\brackets[\big]{x(ax)}y &= x\brackets[\big]{a(xy)}, \label{eq:moufangl} \tag{LM} \\
		y\brackets[\big]{x(ax)} &= \brackets[\big]{(yx)a}x, \label{eq:moufangr} \tag{RM} \\
		(xy)(ax) &= x\brackets[\big]{(ya)x}. \label{eq:moufangm} \tag{MM}
	\end{align}
\end{definition}

The terminology of \enquote{weakly alternative rings} is borrowed from \cite[\onpage{170}]{Faulkner-Barb}. It is not standard.

\begin{example}
	The standard examples of rings which are alternative but not associative are octonion algebras over fields, that is, 8-dimensional composition algebras. In fact, if $ \ring $ is an alternative simple ring, then \cite{Kleinfeld_SimpleAltRings} says that $ \ring $ is either associative or an octonion algebra over a field. We will encounter a generalisation of composition algebras over commutative associative rings in \cref{sec:F4:comp-alg}.
\end{example}

\begin{note}
	We will mostly be concerned with alternative rings in this book. The notion of weakly alternative rings matters for technical reasons: In \cref{sec:BC:blue-comp}, we will prove that any ring $ \ring $ which arises from a $ (B)C_3 $-graded group is always alternative. In the course of the proof, we will first show that $ \ring $ is weakly alternative (\cref{BC:blue:weak-alt}), which allows us to apply \cref{ring:alternative-nucleus}. Only later can we show that $ \ring $ is actually alternative (\cref{BC:blue:alternative}).
\end{note}

\begin{reminder}[Alternating maps]\label{ring:alternating-map}
	Let $ n \in \Npos $ and let $ \map{f}{\ring^n}{\ring}{}{} $ be a map which is additive in each component. Then $ f $ is \defemph*{alternating}\index{alternating map} if for all $ v_1, \ldots, v_n \in \ring $ such that at least two of the elements $ \listing{v}{n} $ are equal, we have $ f(v_1, \ldots, v_n) = 0_\ring $. Further, it is called \defemph*{weakly alternating}\index{alternating map!weakly} or \defemph*{antisymmetric}\index{antisymmetric map} if
	\[ f(v_{\sigma(1)}, \ldots, v_{\sigma(n)}) = (-1)^{\sgn(\sigma)} f(v_1, \ldots, v_n) \]
	for any permutation $ \sigma $ of $ \numint{1}{n} $ where $ \sgn(\sigma) $ denotes the signum of $ \sigma $. (It suffices to verify this condition for transpositions $ \sigma $ because they generate the permutation group.) Any alternating map is weakly alternating, and if $ 2_\ring $ is not a zero divisor in $ \ring $, then every weakly alternating map is alternating.
\end{reminder}

We have the following characterization of alternative rings. The main parts of the proof can be found in \cite[III.1]{Schafer}. The fact that alternative rings satisfy property~\thmitemref{ring:alternative-char}{ring:alternative-char:artin} is usually referred to as \defemph{Artin's Theorem}. We will see variations of this theorem in \cref{ring:alternative-assoc2,rinv:artin-inv}.

\begin{proposition}[Artin's Theorem]\label{ring:alternative-char}
	The following conditions on $ \ring $ are equivalent:
	\begin{stenumerate}
		\item $ \ring $ is alternative.
		\item $ \ring $ is alternative and flexible.
		\item The associator map is alternating.
		\item \label{ring:alternative-char:artin}Every subring of $ \ring $ which is generated by at most two elements is associative.
		\item $ \ring $ satisfies the Moufang identities \eqref{eq:moufangl}, \eqref{eq:moufangr} and \eqref{eq:moufangm}.
		\item $ \ring $ satisfies the left and right Moufang identities \eqref{eq:moufangl} and \eqref{eq:moufangr}.
	\end{stenumerate}
\end{proposition}

\begin{lemma}
	Every alternative ring is weakly alternative. If $ \ring $ is a weakly alternative ring in which $ 2_\ring $ is not a zero divisor, then $ \ring $ is alternative.
\end{lemma}
\begin{proof}
	This follows from \cref{ring:alternating-map,ring:alternative-char}.
\end{proof}

\begin{remark}\label{ring:alternative-nucleus}
	Recall the different notions of nuclei from \cref{ring:nucleus-def}. If $ \ring $ is weakly alternative, then
	\[ \nucleus(\ring) = \lnucleus(\ring) = \mnucleus(\ring) = \rnucleus(\ring). \]
	We will use this fact in \cref{BC:blue:weak-alt}, in a situation in which $ \ring $ is weakly alternative but not yet proven to be alternative.
\end{remark}

Recall from \cref{ring:def-invertible} that the notion of invertibility is slightly surprising (and non-standard) in nonassociative rings. In alternative rings, it turns out that we can use the usual notion of invertibility from associative ring theory.

\begin{lemma}\label{altring:invert-char}
	Assume that $ \ring $ is alternative. For all $ x,y \in \ring $, the following assertions are equivalent:
	\begin{stenumerate}
		\item $ x $ is an inverse of $ y $ in the sense of \cref{ring:def-invertible}.
		
		\item $ xy=1_\ring = yx $.
	\end{stenumerate}
\end{lemma}
\begin{proof}
	Assume that $ xy = 1_\ring = yx $. Using the Moufang identities, it is proven in \cite[\onpage{38}]{Schafer} that $ \assoc{x}{y}{z} = 0_\ring $ for all $ z \in \ring $. Since $ \ring $ is alternative, it follows that $ \assoc{y}{x}{z} = \assoc{z}{x}{y} = \assoc{z}{y}{x} = 0_\ring $ for all $ z \in \ring $ as well. Thus the assertion follows from \cref{ring:inv-char}.
\end{proof}

We will have occasion to use the following identities.

\begin{lemma}\label{ring:alt-sar}
	Assume that $ \ring $ is alternative. Then the following assertions hold for all $ x,x',y,z \in \ring $:
	\begin{lemenumerate}
		\item \label{ring:alt-sar:rightbump}The \defemph*{right bumping formula:}\index{bumping formula} $ \assoc{x}{y}{z} x = \assoc{x}{xy}{z} = \assoc{x}{y}{xz} $. 
		
		\item \label{ring:alt-sar:leftbump}The \defemph*{left bumping formula:} $ x \assoc{x}{y}{z} = \assoc{x}{yx}{z} = \assoc{x}{y}{zx} $.
		
		\item \label{ring:alt-sar:rightbump-lin}$ \assoc{x}{y}{z} x' + \assoc{x'}{y}{z}x = \assoc{x}{x'y}{z} + \assoc{x'}{xy}{z} = \assoc{x}{y}{x'z} + \assoc{x'}{y}{xz} $.
		
		\item \label{ring:alt-sar:leftbump-lin}$ x \assoc{x'}{y}{z} + x' \assoc{x}{y}{z} = \assoc{x}{yx'}{z} + \assoc{x'}{yx}{z} = \assoc{x}{y}{zx'} + \assoc{x'}{y}{zx} $.
	\end{lemenumerate}
\end{lemma}
\begin{proof}
	The first two identities are (2.13) and (2.14) in \cite[Lemma~2.2]{BruckKleinfeldAltDiv}, albeit with a typo in (2.14). The last two identities are linearisations of~\ref{ring:alt-sar:rightbump} and~\ref{ring:alt-sar:leftbump}. By this we mean the following: At first, observe that the formula in~\ref{ring:alt-sar:rightbump} is additive in $ y $ and $ z $ but \enquote{quadratic} in $ x $. Replacing $ x $ by $ x+x' $ in this formula, we obtain
	\begin{align*}
		\assoc{x}{y}{z}x + \assoc{x'}{y}{z}x + \assoc{x}{y}{z}x' + \assoc{x'}{y}{z}x'& \\
		&\hspace{-4cm} \mathord{}= \assoc{x}{y}{xz} + \assoc{x'}{y}{xz} + \assoc{x}{y}{x'z} + \assoc{x'}{y}{x'z}.
	\end{align*}
	All summands in which only one of the variables $ x $ and $ x' $ appears cancel by~\ref{ring:alt-sar:rightbump}, so we remain with
	\[ \assoc{x'}{y}{z}x + \assoc{x}{y}{z}x' = \assoc{x'}{y}{xz} + \assoc{x}{y}{x'z}. \]
	This is precisely one of the equations in~\itemref{ring:alt-sar:rightbump-lin}. The remaining assertions can be proven similarly.
\end{proof}

We now state some formulas concerning nuclear elements.

\begin{lemma}[\enquote{Nuclear Slipping Formula}, {\cite[3.1.6]{McCrimmonAltUnpublished}}, {\cite[21.2.1]{McCrimmon_TasteJordan}}]\label{altring:nucl-slip}
	Assume that $ \ring $ is alternative. Then for all $ x,y,z \in \ring $ and all $ n \in \nucleus(\ring) $, we have\index{nuclear slipping formula}
	\begin{align*}
		n \assoc{x}{y}{z} &= \assoc{nx}{y}{z} = \assoc{xn}{y}{z} = \assoc{x}{ny}{z} = \assoc{x}{yn}{z} \\
		&= \assoc{x}{y}{nz} = \assoc{x}{y}{zn} = \assoc{x}{y}{z}n.
	\end{align*}
	In particular, nuclear elements commute with associators.
\end{lemma}
\begin{proof}
	Let $ x,y,z \in \ring $ and let $ n \in \nucleus(\ring) $. We know from \cref{ring:nucl-slip} that in any nonassociative ring, we have
	\begin{align*}
		n \assoc{x}{y}{z} &= \assoc{nx}{y}{z}, & \assoc{xn}{y}{z} &= \assoc{x}{ny}{z}, \\
		\assoc{x}{yn}{z} &= \assoc{x}{y}{nz}, & \assoc{x}{y}{zn} &= \assoc{x}{y}{z}n.
	\end{align*}
	Further, since the associator is alternating, we also have
	\begin{align*}
		n\assoc{x}{y}{z} &= -n\assoc{y}{x}{z} = -\assoc{ny}{x}{z} = \assoc{x}{ny}{z}.
	\end{align*}
	A similar computations yields that we also have $ n\assoc{x}{y}{z} = \assoc{x}{y}{nz} $. Further,
	\begin{align*}
		\assoc{x}{y}{nz} &= \assoc{x}{yn}{z} = - \assoc{x}{z}{yn} = -\assoc{x}{z}{y}n = \assoc{x}{y}{z}n.
	\end{align*}
	Combining all statements, the assertion follows.
\end{proof}

The following lemma will (only) be used to prove that the nucleus of an alternative ring is a (necessarily multiplicative) square-module, see \cref{rinv:nucl-submod}.

\begin{lemma}[{\cite[(1.8), (1.9)]{McCrimmonAltUnpublished}}]\label{altring:nucl-comm}
	Assume that $ \ring $ is alternative. Then the following statements hold, where $ \commutator{x}{y} \defl xy - yx $ denotes the commutator for $ x,y \in \ring $:
	\begin{lemenumerate}
		\item \label{altring:nucl-comm:clos}$ \commutator[\big]{\nucleus(\ring)}{\ring} \subs \nucleus(\ring) $.
		
		\item \label{altring:nucl-comm:center}$ \commutator{x}{n} \assoc{x}{y}{z} = 0_\ring $ for all $ x,y,z \in \ring $ and all $ n \in \nucleus(\ring) $.
	\end{lemenumerate}
\end{lemma}
\begin{proof}
	Let $ x,y,z \in \ring $ and let $ n \in \nucleus(\ring) $. Then by \cref{altring:nucl-slip},
	\begin{align*}
		\assoc[\big]{\commutator{n}{x}}{y}{z} &= \assoc{nx}{y}{z} - \assoc{xn}{y}{z} = n \assoc{x}{y}{z} - n \assoc{x}{y}{z} = 0,
	\end{align*}
	which implies the first assertion. The second assertion follows from another application of \cref{altring:nucl-slip} in combination with the left bumping formula (\thmitemcref{ring:alt-sar}{ring:alt-sar:leftbump}): We have
	\begin{align*}
		(xn) \assoc{x}{y}{z} &= x \brackets[\big]{n \assoc{x}{y}{z}} = x \brackets[\big]{\assoc{x}{y}{z}n} = \brackets[\big]{x\assoc{x}{y}{z}}n \\
		&= \assoc{x}{yx}{z}n = n \assoc{x}{yx}{z} = n \brackets[\big]{x \assoc{x}{y}{z}} = (nx) \assoc{x}{y}{z},
	\end{align*}
	which implies that $ \commutator{x}{n} \assoc{x}{y}{z} = 0_\ring $.
\end{proof}

We end this section with some results on associative subrings of alternative rings.

\begin{definition}[Strongly associative set, {\cite[14.1]{GPR_AlbertRing}}]
	A subset $ A $ of $ \ring $ is called \defemph*{strongly associative}\index{strongly associative set} if $ \assoc{a}{a'}{x} = 0_\ring $ for all $ a,a' \in A $ and all $ x \in \ring $.
\end{definition}

\begin{lemma}[{\cite[{Theorem I.2}]{BruckKleinfeldAltDiv}}]\label{ring:alternative-assoc1}
	Assume that $ \ring $ is alternative and let $ A,B,C $ be subsets of $ \ring $ satisfying $ \assoc{A}{A}{\ring} = \assoc{B}{B}{\ring} = \assoc{C}{C}{\ring} = \assoc{A}{B}{C} = \compactSet{0} $. (That is, $ A,B,C $ are strongly associative sets such that $ \assoc{A}{B}{C} = \compactSet{0} $.) Then the subring of $ \ring $ which is generated by $ A \union B \union C $ is associative.
\end{lemma}

\begin{note}
	Property~\thmitemref{ring:alternative-char}{ring:alternative-char:artin} is a special case of \cref{ring:alternative-assoc1}: Assume that $ \ring $ is alternative and let $ x,y \in \ring $. Then the sets $ A \defl \Set{x} $, $ B \defl \Set{y} $ and $ C \defl \Set{0} $ satisfy the conditions of \cref{ring:alternative-assoc1}, and so the subring of $ \ring $ generated by $ \Set{x,y} $ is associative. Similarly, it follows from \cref{ring:alternative-assoc1} that the subring generated by $ \Set{x,y,z} $ where $ x,y,z \in \ring $ satisfy $ \assoc{x}{y}{z} = 0 $ is alternative.
\end{note}

\begin{lemma}[{\cite[{Theorem I.3}]{BruckKleinfeldAltDiv}}]\label{ring:alternative-assoc2}
	Assume that $ \ring $ is alternative and let $ A,B $ be subsets of $ \ring $ such that $ A $ is an associative subring of $ \ring $ and $ \assoc{A}{A}{B} = \assoc{B}{B}{\ring} = \compactSet{0} $. (That is, $ B $ is a strongly associative subset of $ \ring $ and $ A $ is an associative subring of $ \ring $ such that $ \assoc{A}{A}{B} = \compactSet{0} $.) Then the subring of $ \ring $ generated by $ A \union B $ is associative.
\end{lemma}

The following corollary of \cref{ring:alternative-assoc2} is a useful generalisation of Artin's Theorem (\cref{ring:alternative-char}).

\begin{proposition}
	Assume that $ \ring $ is alternative and let $ x,y \in \ring $. Then the subring generated by $ \Set{x,y} \union \nucleus(\ring) $ is associative.
\end{proposition}
\begin{proof}
	This follows from \cref{ring:alternative-assoc2} for $ A \defl \nucleus(\ring) $ and $ B \defl \Set{x,y} $.
\end{proof}


\section{Rings with Involution}

\label{sec:rinv}

\begin{secnotation}
	Unless otherwise specified, we denote by $ \ring $ a ring with an involution $ \map{\rinvmap}{\ring}{\ring}{}{} $ in the sense of the following \cref{rinv:invo-def}.
\end{secnotation}

In this section, we study involutions on rings. All rings with involution that we are interested in will be alternative. However, we will be faced with situations in which a ring with involution is alternative, but not yet proven to be alternative. For this reason, it is more practical to introduce involutions on arbitrary rings. A valuable, albeit not easily accessible reference on this subject is \cite[Section~9.5]{McCrimmonAltUnpublished}. Another comprehensive book on the subject in the setting of algebras over a field is \cite{BookOfInvolutions}.

\subsection{Basic Definitions and Observations}

\begin{definition}[Involution]\label{rinv:invo-def}
	An \defemph*{involution of $ \ring $}\index{involution} is a map $ \map{\rinvmap}{\ring}{\ring}{r}{\rinv{r}} $
	satisfying $ \rinv{(r+s)} = \rinv{r} + \rinv{s} $, $ \rinv{(r \rmult s)} = \rinv{s} \rmult \rinv{r} $ and $ \rinv{(\rinv{r})} = r $ for all $ r,s \in \ring $. We will refer to the second property as \defemph*{anti-compatibility with the ring multiplication}\index{anti-compatible}. Elements $ r \in \ring $ with the property that $ \rinv{r} = r $ are called \defemph*{symmetric}\index{symmetric element}, and we denote the set of all symmetric elements by $ \symring $.
\end{definition}

\begin{example}\label{rinv:id-comm}
	The identity map is an involution on $ \ring $ if and only if $ \ring $ is commutative. If this is the case, it is called the \defemph*{trivial involution (on $ \ring $)}.\index{involution!trivial}
\end{example}

\begin{example}
	The conjugation map $ \map{}{\IC}{\IC}{a+bi}{a-bi} $ is an involution on the field of complex numbers. More generally, any conic algebra has a conjugation map which is an involution in certain situations. See \cref{conic:invo-char,conic:normass-flex-conj,conic:mult-is-norm} for more details.
\end{example}

\begin{remark}\label{rinv:1-symmetric}
	The set of symmetric elements is an additive subgroup of $ (\ring, +) $, and it contains $ 1_\ring $ because
	\[ 1_\ring = \rinv{(\rinv{1_\ring})} = \rinv{(\rinv{1_\ring} \rmult 1_\ring)} = \rinv{1_\ring} \rmult \rinv{(\rinv{1_\ring})} = \rinv{1_\ring} \rmult 1_\ring = \rinv{1_\ring}. \]
	However, it is in general not a subring of $ \ring $: For all $ r,s \in \symring $, we have $ \rinv{(rs)} = \rinv{s} \rinv{r} = sr $ which, in general, is distinct from $ rs $.
\end{remark}

\begin{lemma}
	Let $ r \in \ring $ be invertible. Then $ \rinv{r} $ is invertible and $ (\rinv{r})^{-1} = \rinv{(r^{-1})} $.
\end{lemma}
\begin{proof}
	For all $ s \in \ring $, we have
	\begin{align*}
		\rinv{s} &= \rinv{\brackets[\big]{r^{-1} (rs)}} = (\rinv{s} \rinv{r}) \rinv{(r^{-1})}.
	\end{align*}
	It follows that $ (s\rinv{r}) \rinv{(r^{-1})} = s $ for all $ s \in \ring $. Similarly, one can show that the other properties in \cref{ring:def-invertible} are satisfied. Thus $ \rinv{r} $ is invertible and its inverse is $ \rinv{(r^{-1})} $.
\end{proof}

\begin{notation}
	For any invertible $ r \in \ring $, we denote the element $ (\rinv{r})^{-1} = \rinv{(r^{-1})} $ by $ \rinvmin{r} $.
\end{notation}

We can define norm and trace maps on any ring with involution. The latter in particular will be crucial in the theory of Jordan modules.

\begin{definition}[Traces and norms]\label{rinv:trace-norm-def}
	The map $ \map{\ringTr = \ringTr_\rinvmap}{\ring}{\ring}{r}{r+\rinv{r}} $ is called the \defemph*{trace (on $ \ring $)}\index{trace!on a ring with involution} and the map $ \map{\ringnorm = \ringnorm_\rinvmap}{\ring}{\ring}{r}{\rinv{r}r} $ is called the \defemph*{norm (on $ \ring $)}\index{norm!on a ring with involution}. We will also refer to ring elements which lie in the image of these maps as \defemph*{traces} and \defemph*{norms} and denote the sets of such elements by $ \ringTr(\ring) $ and $ \ringnorm(\ring) $, respectively. Further, we define a bi-additive map
	\[ \map{}{\ring \times \ring}{\ringTr(\ring)}{(a,b)}{\ringTr(a, b) \defl \ringTr(a \rinv{b})} \]
	which we also denote by $ \ringTr $ or $ \ringTr_\rinvmap $.
\end{definition}

\begin{note}
	In the definition of the norm, we make the somewhat arbitrary choice that $ \ringnorm(r) = \rinv{r}r $ and not $ \ringnorm(r) = r \rinv{r} $. Observe that the set $ \ringnorm(\ring) $ is independent of this choice because the norm of $ r $ in the first sense equals the norm of $ \rinv{r} $ in the second sense. In particular, the notions of nuclear and central involutions that we will introduce in \cref{rinv:nucl-def} do not depend on this choice. Further, we will show in \cref{rinv:cent-norm-comm} that the two definitions coincide if the involution $ \rinvmap $ is central.
\end{note}

\begin{lemma}\label{rinv:trace-norm-sym}
	The following hold:
	\begin{lemenumerate}
		\item Every trace and every norm in $ \ring $ is symmetric, and the trace is an additive map.
		
		\item \label{rinv:trace-norm-sym-2}If $ 2_\ring $ is invertible, then $ \ringnorm(\ring) \subs \ringTr(\ring) = \symring $.
	\end{lemenumerate}
\end{lemma}
\begin{proof}
	The first assertion is clear. Now assume that $ 2_\ring $ is invertible and let $ r \in \ring $ be symmetric. Then
	\[ r = \frac{r+r}{2} = \frac{r}{2} + \frac{\rinv{r}}{2} = \ringTr(r/2), \]
	which shows that $ \Fix(\rinvmap) $ is contained in $ \ringTr(\ring) $. Thus $ \ringnorm(\ring) \subs \Fix(\rinvmap) = \ringTr(\rinvmap) $, as desired.
\end{proof}

The following result shows that the mere existence of an involution on a ring makes it easier to prove that a ring is (weakly) alternative. We will use these facts in the blueprint computations for $ BC $, see \cref{BC:blue:weak-alt,BC:blue:alternative}.

\begin{lemma}\label{rinv:inv-assoc}
	The following hold:
	\begin{lemenumerate}
		\item We have $ \rinv{\assoc{x}{y}{z}} = -\assoc{\rinv{z}}{\rinv{y}}{\rinv{x}} $ for all $ x,y,z \in \ring $.
		
		\item \label{rinv:inv-assoc:weak}Let $ \ring' $ be any ring which admits an involution $ \rinvmap $ and which satisfies $ \assoc{x}{y}{z} = - \assoc{x}{z}{y} $ for all $ x,y,z \in \ring' $. Then $ \ring' $ is weakly alternative.
		
		\item \label{rinv:inv-assoc:alt}Let $ \ring' $ be any ring which admits an involution $ \rinvmap $ and which satisfies $ \assoc{x}{y}{y} = 0 $ for all $ x,y \in \ring $. Then $ \ring' $ is alternative.
	\end{lemenumerate}
\end{lemma}
\begin{proof}
	The first assertion is a simple computation:
	\begin{align*}
		\rinv{\assoc{x}{y}{z}} &= \rinv{\brackets[\big]{(xy)z - x(yz)}} = \rinv{z}(\rinv{y} \rinv{x}) - (\rinv{z} \rinv{y}) \rinv{x} = -\assoc{\rinv{z}}{\rinv{y}}{\rinv{x}}.
	\end{align*}
	The remaining two assertions follow from the first one.
\end{proof}

If $ \ring $ is not associative, we usually need additional assumptions on the involution $ \rinvmap $ to develop a satisfactory theory. In the context of Jordan modules, nuclear involutions are particularly important. The stronger notion of central involutions is not relevant for Jordan modules, but we will briefly see it again in \cref{sec:F4:conic}.

\begin{definition}[Nuclear and central involutions]\label{rinv:nucl-def}
	The involution $ \rinvmap $ is called \defemph*{nuclear}\index{involution!nuclear} if every norm lies in the nucleus of $ \ring $. It is called \defemph*{central}\index{involution!central} if every norm lies in the center of~$ \ring $.
\end{definition}

If the ring $ \ring $ is equipped with the structure of an algebra over some commutative associative ring $ \comring $, then we can also define the notion of scalar involutions in a similar way. See \cref{F4:scal-inv-def}.

\begin{remark}\label{rinv:nucl-cent-lem}
	Assume that $ \rinvmap $ is nuclear. Then for all $ r \in \ring $, the element
	\[ \rinv{(1+r)}(1+r) = 1 + \rinv{r} + r + \rinv{r} r \]
	lies in the nucleus. As the nucleus is closed under addition, this implies that the trace $ r+\rinv{r} $ lies in the nucleus. Thus if $ \rinvmap $ is nuclear, then all traces lie in the nucleus as well. Similarly, if $ \rinvmap $ is central, then all traces lie in the center. If $ 2_\ring $ is invertible, then it is clear that the converses of these statements are also true because then every norm is a trace by \thmitemcref{rinv:trace-norm-sym}{rinv:trace-norm-sym-2}. In \cref{ring:nucl-invo-char,ring:cent-invo-char}, we will see that the same remains true for alternative rings even without assumptions on the invertibility of $ 2_\ring $.
\end{remark}

\begin{remark}\label{rinv:cent-norm-comm}
	If $ \rinvmap $ is central, then it follows from $ a (a+\rinv{a}) = (a+\rinv{a})a $ that $ a \rinv{a} = \rinv{a}a = \ringnorm_\rinvmap(a) $ for all $ a \in \ring $. Linearising this identity, we obtain that $ \ringTr(a \rinv{b}) = \ringTr(\rinv{a} b) $ for all $ a,b \in \ring $.
\end{remark}

\subsection{Nuclear Involutions}

We now study nuclear involutions. A recurring theme here is that many associativity results involving a ring element $ r $ remain true if we replace some (but not necessarily all) occurrences of $ r $ by $ \rinv{r} $. The main example of this behaviour is Artin's Theorem for nuclear involutions, see \cref{rinv:artin-inv}. Another main result in this subsection is \cref{ring:nucl-invo-char}.

\begin{lemma}\label{rinv:invo-assoc}
	Assume that every trace is nuclear. Then for all $ r,s,t \in \ring $, we have $ \assoc{\rinv{r}}{s}{t} = -\assoc{r}{s}{t} $, $ \assoc{r}{\rinv{s}}{t} = -\assoc{r}{s}{t} $ and $ \assoc{r}{s}{\rinv{t}} = -\assoc{r}{s}{t} $. In particular, if $ \ring $ is alternative, then $ \assoc{r}{\rinv{r}}{s} $, $ \assoc{r}{s}{\rinv{r}} $ and $ \assoc{s}{r}{\rinv{r}} $ are zero for all $ r,s \in \ring $.
\end{lemma}
\begin{proof}
	Let $ r,s,t \in \ring $. Since $ \rinvmap $ is nuclear, we have $ \assoc{r+\rinv{r}}{s}{t} = 0 $. This implies that $ \assoc{\rinv{r}}{s}{t} = -\assoc{r}{s}{t} $. The other assertions can be proven similarly.
\end{proof}

\begin{lemma}[Bumping formulas for nuclear involutions]\label{rinv:bumping}
	Assume that $ \ring $ is alternative and that every trace is nuclear. Then for all $ x,y,z \in \ring $, the following variations of the bumping formulas\index{bumping formula} from \cref{ring:alt-sar} are satisfied:
	\[ \assoc{\rinv{x}}{y}{z} x = \assoc{\rinv{x}}{xy}{z} = \assoc{\rinv{x}}{y}{xz} \midand x \assoc{\rinv{x}}{y}{z} = \assoc{\rinv{x}}{yx}{z} = \assoc{\rinv{x}}{y}{zx} \]
\end{lemma}
\begin{proof}
	This is a consequence of \cref{rinv:invo-assoc} and the original bumping formulas. For example,
	\[ \assoc{\rinv{x}}{y}{z} x = -\assoc{x}{y}{z}x = -\assoc{x}{xy}{z} = \assoc{\rinv{x}}{xy}{z}. \]
	The other formulas can be derived in the same fashion.
\end{proof}

\begin{proposition}[Artin's Theorem for nuclear involutions]\label{rinv:artin-inv}
	Assume that $ \ring $ is alternative and that the involution on $ \ring $ is nuclear. Then for all $ r,s \in \ring $, the subring of $ \ring $ generated by $ \Set{r,s,\rinv{r}, \rinv{s}} \union \nucleus(\ring) $ is associative.\index{Artin's Theorem!for nuclear involutions}
\end{proposition}
\begin{proof}
	Set $ A \defl \Set{r, \rinv{r}} $, $ B \defl \Set{s, \rinv{s}} $ and $ C \defl \Set{0} $. Then by \cref{rinv:invo-assoc}, we have
	\[ \assoc{A}{A}{\ring} = \assoc{B}{B}{\ring} = \assoc{C}{C}{\ring} = \assoc{A}{B}{C} = \compactSet{0}. \]
	Thus it follows from \cref{ring:alternative-assoc1} that the subring $ D $ of $ \ring $ which is generated by $ \Set{r, \rinv{r}, s, \rinv{s}} $ is associative. Applying \cref{ring:alternative-assoc2}, we conclude that the subring generated by $ D \union \nucleus{\ring} $ is associative.
\end{proof}

\begin{proposition}[{\cite[Section~9.5]{McCrimmonAltUnpublished}}]\label{ring:nucl-invo-char}
	Assume that $ \ring $ is alternative. Then the following assertions on the involution $ \rinvmap $ are equivalent.
	\begin{stenumerate}
		\item Every norm lies in the nucleus of $ \ring $. That is, $ \rinvmap $ is nuclear.
		
		\item Every trace lies in the nucleus of $ \ring $.
	\end{stenumerate}
\end{proposition}
\begin{proof}
	We already know from \cref{rinv:nucl-cent-lem} that every trace lies in the nucleus if $ \rinvmap $ is nuclear. Conversely, assume that every trace is nuclear. Let $ x,y,z \in \ring $. By the linearised right bumping formula (\thmitemcref{ring:alt-sar}{ring:alt-sar:rightbump-lin}), we have
	\begin{align*}
		\assoc{x}{\rinv{x}}{y}z + \assoc{z}{\rinv{x}}{y}x = \assoc{x}{z\rinv{x}}{y} + \assoc{z}{x \rinv{x}}{y}.
	\end{align*}
	Since the associator is alternating, we can shift all associators on the right-hand side of the equation one step to the left, which yields that
	\begin{align*}
		\assoc{x \rinv{x}}{y}{z} &= \assoc{x}{\rinv{x}}{y}z + \assoc{z}{\rinv{x}}{y}x - \assoc{z \rinv{x}}{y}{x}.
	\end{align*}
	Using \cref{rinv:invo-assoc}, we obtain that
	\begin{align*}
		\assoc{x \rinv{x}}{y}{z} &= - \assoc{x}{x}{y}z + \assoc{\rinv{z}}{x}{y}x + \assoc{x \rinv{z}}{y}{x} =  \assoc{x}{y}{\rinv{z}}x - \assoc{x}{y}{x \rinv{z}}
	\end{align*}
	By the right bumping formula (\thmitemcref{ring:alt-sar}{ring:alt-sar:rightbump}), we conclude that $ \assoc{x \rinv{x}}{y}{z} = 0 $. This shows that every norm is nuclear as well.
\end{proof}

\begin{remark}\label{ring:cent-invo-char}
	The assertion of \cref{ring:nucl-invo-char} is also true if we replace \enquote{nuclear} by \enquote{central}, see \cite[Section~9.5]{McCrimmonAltUnpublished}.
\end{remark}

\begin{lemma}\label{rinv:trace-norm-sym-nucl}
	If $ \rinvmap $ is nuclear and $ 2_\ring $ is invertible, then
	\[ \ringnorm(\ring) \subs \ringTr(\ring) = \symring \subs \nucleus(\ring). \]
\end{lemma}
\begin{proof}
	We have seen in \thmitemcref{rinv:trace-norm-sym}{rinv:trace-norm-sym-2} that $ \ringnorm(\ring) \subs \ringTr(\ring) = \symring $. Since $ \rinvmap $ is nuclear, we have that $ \ringTr(\ring) \subs \nucleus(\ring) $.
\end{proof}

\begin{lemma}\label{rinv:Tr-assoc}
	Assume that $ \ring $ is alternative and that the involution $ \rinvmap $ is nuclear. Then the map $ \map{\ringTr}{\ring}{\ring}{a}{a+\rinv{a}} $ is associative in the sense that $ \ringTr((ab)c) = \ringTr(a(bc)) $ for all $ a,b,c \in \ring $.
\end{lemma}
\begin{proof}
	Let $ a,b,c \in \ring $. Since
	\begin{align*}
		\ringTr\brackets[\big]{(ab)c} &= (ab)c + \rinv{c} (\rinv{b} \rinv{a}) \midand \ringTr\brackets[\big]{a(bc)} = a(bc) + (\rinv{c} \rinv{b}) \rinv{a},
	\end{align*}
	we have to show that $ \assoc{a}{b}{c} = \assoc{\rinv{c}}{\rinv{b}}{\rinv{a}} $. This follows from \cref{rinv:invo-assoc} and the fact that the associator is alternating.
\end{proof}


\section{Involutory Sets}

\label{sec:invset}

\begin{secnotation}
	Unless otherwise specified, we denote by $ \ring $ an alternative ring with nuclear involution $ \rinvmap $.
\end{secnotation}

In this section, we study involutory sets, which are a slight refinement of alternative rings with a nuclear involution. The terminology of \enquote{involutory sets} was introduced in \cite[(11.1)]{MoufangPolygons}, but more general notions had been studied beforehand. We will discuss this in \cref{jormod:invoset-literature}.

The following map plays a major role in the context of involutory sets. We will always interpret it in the language of square-scalar multiplications from \cref{sec:weak-quad}.

\begin{example}\label{jormod:invring-jordan}
	Let $ \ring $ be an alternative ring with nuclear involution. We define a square-scalar multiplication on $ \ring $ by
	\[ \map{\omega}{\ring \times \ring}{\ring}{(r,s)}{\rinv{s}rs = (\rinv{s}r)s = \rinv{s}(rs)}, \]
	which is well-defined by \cref{rinv:artin-inv}. It is called the \defemph*{canonical square-scalar multiplication on $ \ring $}\index{square-scalar multiplication!on a ring}. With this structure, $ \ring $ is a properly quadratic square-module over itself. It is also multiplicative if $ \ring $ is associative, but it is not known to us whether the same is true for alternative rings.
\end{example}

In the following, we will frequently apply \cref{rinv:artin-inv} without comment to avoid excessive use of brackets. Further, we will always use the square-module structure on $ \ring $ given by \cref{jormod:invring-jordan}. We begin by identifying several interesting square-submodules of $ \ring $.

\begin{lemma}\label{rinv:sym-submod}
	$ \symring $ is a square-submodule of $ \ring $.
\end{lemma}
\begin{proof}
	It is clear that $ \symring $ is an additive subgroup of $ \ring $, and we further have $ \rinv{(\rinv{r}sr)} = \rinv{r} \rinv{s} r = \rinv{r} s r $ for all $ r \in \ring $, $ s \in \symring $.
\end{proof}

\begin{lemma}\label{rinv:norm-submod}
	$ \ringnorm(\ring) $ is stable under the square-scalar multiplication of $ \ring $.
\end{lemma}
\begin{proof}
	Let $ r,s \in \ring $. Then $ \rinv{r}\ringnorm(s) r = \rinv{r}(\rinv{s} s)r = (\rinv{r} \rinv{s}) (sr) = \ringnorm(sr) $ by \cref{rinv:artin-inv}, so the assertion follows.
\end{proof}

\begin{lemma}[{\cite[Section~9.5]{McCrimmonAltUnpublished}}]\label{rinv:nucl-submod}
	The nucleus of $ \ring $ is a multiplicative square-submodule of $ \ring $.
\end{lemma}
\begin{proof}
	Let $ x,y,z \in \ring $ and let $ n \in \nucleus(\ring) $. We want to show that $ \rinv{x}nx $ lies in the nucleus. Observe that $ \rinv{x}nx = \commutator{\rinv{x}}{n}x + n \rinv{x}x $ where $ n\rinv{x}x $ lies in the nucleus because the involution is nuclear. Thus
	\begin{align*}
		\assoc{\rinv{x}nx}{y}{z} &= \assoc[\big]{\commutator{\rinv{x}}{n}x}{y}{z} + \assoc{n \rinv{x}x}{y}{z} = \assoc[\big]{\commutator{\rinv{x}}{n}x}{y}{z}.
	\end{align*}
	Since $ \commutator{\rinv{x}}{n} $ lies in the nucleus by \thmitemcref{altring:nucl-comm}{altring:nucl-comm:clos}, it follows from \cref{ring:nucl-slip,rinv:invo-assoc} yields that
	\begin{align*}
		\assoc{\rinv{x}nx}{y}{z} &=\assoc[\big]{\commutator{\rinv{x}}{n}x}{y}{z} = \commutator{\rinv{x}}{n} \assoc{x}{y}{z} = -\commutator{\rinv{x}}{n} \assoc{\rinv{x}}{y}{z}.
	\end{align*}
	Using \thmitemcref{altring:nucl-comm}{altring:nucl-comm:center}, we infer that $ \assoc{\rinv{x}nx}{y}{z} = 0_\ring $. By \cref{ring:alternative-nucleus}, this shows that $ \rinv{x}nx $ lies in the nucleus. Thus the nucleus is indeed a square-submodule of $ \ring $, and it is multiplicative by Artin's Theorem for nuclear involutions (\cref{rinv:artin-inv}).
\end{proof}

\begin{lemma}\label{rinv:tr-submod}
	The map $ \map{\ringTr}{\ring}{\nucleus(\ring) \intersect \symring}{r}{r+\rinv{r}} $ is a homomorphism of square-modules. In particular, $ \ringTr(\ring) $ is a multiplicative square-submodule of $ \nucleus(\ring) \intersect \symring $.
\end{lemma}
\begin{proof}
	For all $ r,s \in \ring $, we have
	\begin{align*}
		\ringTr(\rinv{r}sr) = \rinv{r}sr + \rinv{(\rinv{r}sr)} = \rinv{r}sr + \rinv{r}\rinv{s}r = \rinv{r}(s+\rinv{s})r = \rinv{r} \ringTr(s)r
	\end{align*}
	which proves the first statement. It follows that $ \ringTr(\ring) $ is the image of a homomorphism of square-modules, so it is a square-submodule of $ \nucleus(\ring) \intersect \symring $ by \cref{jormod:ker-im}. Since the square-module $ \nucleus(\ring) $ is multiplicative, so is every square-submodule.
\end{proof}

\begin{lemma}\label{rinv:norm-tr-submod}
	The additive subgroup generated by $ \ringnorm(\ring) \union \ringTr(\ring) $ is a multiplicative square-submodule of $ \symring \intersect \nucleus(\ring) $.
\end{lemma}
\begin{proof}
	It is a square-submodule of $ \ring $ by \cref{rinv:tr-submod,rinv:norm-submod}, it is clearly contained in $ \symring $ and it is contained in $ \nucleus(\ring) $ by \cref{ring:nucl-invo-char}.
\end{proof}

With these elementary observations in mind, we can now define involutory sets.

\begin{definition}[Involutory set]\label{rinv:inv-set-def}
	A tuple $ (\ring, \ringzero, \rinvmap) $ is called a \defemph{pre-involutory set (in $ (\ring, \rinvmap) $)}\index{involutory set!pre-} if $ \ring $ is an alternative ring, $ \rinvmap $ is a nuclear involution on $ \ring $ and $ \ringzero $ is a square-submodule of $ \ring $ such that
	\[ \ringTr(\ring) \subs \ringzero \subs \symring \intersect \nucleus(\ring). \]
	It is called an \defemph{involutory set} if, in addition, $ \ringzero $ contains $ 1_\ring $. It is called \defemph*{associative} if the ring $ \ring $ is associative and it is called an \defemph*{involutory division set}\index{involutory division set} if $ \ring $ is a division ring (in the sense of \cref{ring:division-ring-def}).
\end{definition}

\begin{remark}[on \cref{rinv:inv-set-def}]
	Let $ (\ring, \ringzero, \rinvmap) $ be a pre-involutory set. If it is an involutory set, then for all $ r \in \ring $, we have
	\[ \ringnorm(r) = \rinv{r}1_\ring r \in \rinv{r} \ringzero r \subs \ringzero. \]
	Thus $ \ringzero $ contains $ \ringnorm(\ring) $. Conversely, if $ \ringzero $ contains $ \ringnorm(\ring) $, then it contains $ \ringnorm(1_\ring) = 1_\ring $.
\end{remark}

\begin{note}[Involutory sets in the literature]\label{jormod:invoset-literature}
	The terminology of \enquote{involutory sets} was introduced in \cite[(11.1)]{MoufangPolygons}, but in a less general context: Involutory sets in the sense of \cite{MoufangPolygons} are exactly the associative involutory division sets in our sense. Since \cite{MoufangPolygons} considers only RGD-systems, the division assumption on involutory sets is not surprising. Further, since every simple alternative ring is either associative or a Cayley-Dickson algebra by \cite{Kleinfeld_SimpleAltRings}, an involutory division set is either associative or comes from a Cayley-Dickson algebra. Involutory division sets of the latter form are called \defemph*{honorary involutory sets}\index{involutory set!honorary} in \cite[(38.11)]{MoufangPolygons}. Hence the non-associative case is not actually absent from \cite{MoufangPolygons}, but simply given a different name.
	
	It should be noted that involutory sets in our sense are the same thing as ample subspaces (of $ \ring $) in \cite[Section~9.5]{McCrimmonAltUnpublished}. Further, they are a special case of form parameters\index{form parameter} in the sense of \cite[5.1.C]{HahnOMeara-ClassicalGroups} and \cite[\onpage{23}]{Knus_QuadFormRing}, except that all rings in the latter references are associative.
\end{note}

The following observation illustrates that the notion of involutory sets is only a slight refinement of alternative rings with a nuclear involution and that its main purpose is to capture certain \enquote{characteristic-2 phenomena}.

\begin{remark}\label{rinv:inv-set-2inv}
	Assume that $ 2_\ring $ is invertible in $ \ring $. Then
	\[ \ringTr(\ring) \union \ringnorm(\ring) = \ringTr(\ring) = \symring = \symring \intersect \nucleus(\ring) \]
	by \cref{rinv:trace-norm-sym-nucl}. It follows that $ (\ring, \symring, \rinvmap) $ is the only pre-involutory set in $ (\ring, \rinvmap) $, and it is in fact an involutory set. Thus if we restrict our attention to rings in which $ 2_\ring $ is invertible, then an involutory set carries exactly the same data as an alternative ring with nuclear involution.
	
	In the general situation, denote by $ T_1 $ the additive group that is generated by $ \ringTr(\ring) \union \ringnorm(\ring) $. Then
	\[ \brackets[\big]{\ring, \ringTr(\ring), \rinvmap} \midand (\ring, T_1, \rinvmap) \]
	are the minimal pre-involutory set and the minimal involutory set in $ (\ring, \rinvmap) $, respectively. Further,
	\[ (\ring, \nucleus(\ring) \intersect \symring, \rinvmap) \]
	is the maximal pre-involutory set in $ (\ring, \rinvmap) $, and it is also an involutory set.
\end{remark}

We close this section with some basic facts about involutory sets.

\begin{remark}\label{rinv:r0-equals-r}
	Let $ (\ring, \ringzero, \rinvmap) $ be a pre-involutory set and assume that $ \ring = \ringzero $. Then $ \ring = \symring \intersect \nucleus(\ring) $, so $ \rinvmap $ must be the identity map and $ \ring $ must be associative. By \cref{rinv:id-comm}, the first fact implies that $ \ring $ must be commutative as well. We conclude that the case $ \ring = \ringzero $ can only arise under very specific circumstances.
\end{remark}

\begin{remark}\label{rinvset:dirsum}
	Let $ (\ring, \ringzero, \rinvmap) $ and $ (\ring', \ringzero', \rinvmap') $ be two pre-involutory sets. Then $ (\ring \dirsum \ring', \ringzero \dirsum \ringzero', \rinvmap \dirsum \rinvmap') $ is also an involutory set where $ \rinvmap \dirsum \rinvmap' $ is the nuclear involution which maps $ (r,s) $ to $ (\rinv{r}, s^{\rinvmap'}) $. Note that this involutory set is never an involutory division set unless $ \ring $ or $ \ring' $ is zero. Further, it is an involutory set if and only if both $ (\ring, \ringzero, \rinvmap) $ and $ (\ring', \ringzero', \rinvmap') $ are involutory sets.
\end{remark}

\begin{lemma}\label{rinv:sum-in-R0-lem}
	Let $ (\ring, \ringzero, \rinvmap) $ be a pre-involutory set, let $ h \in \ringzero $ and let $ r,s \in \ring $. Then $ \rinv{r} h s + \rinv{s} hr $ lies in $ \ringzero $.
\end{lemma}
\begin{proof}
	Note that, since $ h $ is nuclear, the expression $ \rinv{r} h s + \rinv{s} hr $ is well-defined. Further,
	\begin{align*}
		\rinv{r} h s + \rinv{s} hr = \rinv{(r+s)} h(r+s) - \rinv{r}hr - \rinv{s}hs
	\end{align*}
	where the element on the right-hand side lies in $ \ringzero $ because each summand does. This finishes the proof.
\end{proof}

\begin{remark}
	Even when $ h $ is an arbitrary element of $ \ring $, the expression on the right-hand side of the equation in the proof of \cref{rinv:sum-in-R0-lem} is still well-defined by \cref{rinv:artin-inv}. This observation yields that
	\[ (\rinv{r} h) s + (\rinv{s} h)r = \rinv{r} (h s) + \rinv{s} (hr) \]
	for all $ r,s,h \in \ring $, but this element does not necessarily lie in $ \ringzero $.
\end{remark}

The following argument is used several times in \cite[Chapter~11]{MoufangPolygons}, where all involutory sets are assumed to have division. It has some important consequences for pseudo-quadratic modules over these involutory sets, which we collect in \cref{pseud:divring-fvv,pseud:divring-f-unique}. In our more general settings, these arguments do not apply. We merely state them to illustrate the additional difficulties that arise in the non-division setting.

\begin{lemma}\label{pseud:divring-R0}
	Assume that $ \ring $ is a division ring and that $ \ringzero \ne \ring $. Let $ s \in \ring $ such that $ sr \in \ringzero $ for all $ r \in \ring $. Then $ s = 0_\ring $.
\end{lemma}
\begin{proof}
	Assume that $ s \ne 0 $ and choose $ t \in \ring \setminus \ringzero $. Then by putting $ r \defl s^{-1}t $, we infer that $ t = sr \in \ringzero $, which is a contradiction.
\end{proof}


\section{Pseudo-quadratic Modules}

\label{sec:pseud}

In this section, we study pseudo-quadratic forms on modules. These objects were introduced over associative division rings by Tits in \cite[Section~8.2]{Tits-LectureNotes74} during the classification of thick spherical buildings of rank at least~3. Unsurprisingly, they also play in important role in the classification of thick spherical Moufang buildings of rank~2 in \cite{MoufangPolygons}. Most arguments in this section stem from the latter reference, though in some cases, they only prove weaker results in our more general setting.

Any pseudo-quadratic module is defined with respect to some involutory set. As usual, we cannot restrict our attention to division rings in the context of root graded groups. In contrast, it is a surprising observation that we even have to allow alternative rings in this setting. (Recall from \cref{rinv:inv-set-def} that the ring $ \ring $ in an involutory set is assumed to be alternative, and also that the involution $ \rinvmap $ is nuclear.) However, it will turn out that this is more of a technical nuisance and that all \enquote{interesting} phenomena occur within the associative world. We will investigate this more deeply in \cref{sec:pure-alt}.

It is noteworthy that pseudo-quadratic modules over associative rings are a special case of quadratic modules over rings with a form parameter in the sense of \cite[Section~5.1D]{HahnOMeara-ClassicalGroups} (see also \cref{jormod:invoset-literature}). However, we are not aware of a reference which specifically covers the special case of pseudo-quadratic forms in the generality of non-division base rings. Further, we do not know a previous reference which treats this subject in the nonassociative setting.

\subsection{Definition and Examples}

Before we can define pseudo-quadratic forms, we have to introduce skew-hermitian forms. These objects relate to pseudo-quadratic maps in a similar way as symmetric bilinear forms relate to quadratic maps. In this section, we will only encounter skew-hermitian forms on $ \ring $-modules. However, we will also need skew-hermitian forms on square-modules over $ \ring $ to define Jordan modules in \cref{sec:jordan-modules}, which is why we define them in this generality.

\begin{definition}[Sesquilinear form]\label{pseud:sesqui-def}
	Let $ \ring $ be an alternative ring with involution~$ \rinvmap $ and let $ (\module, \omega) $ be a right square-module over $ \ring $. A \defemph*{sesquilinear form on $ \module $ with respect to $ \rinvmap $}\index{sesquilinear form} is a bi-additive map $ \map{f}{\module \times \module}{\nucleus(\ring)}{}{} $ which satisfies
	\[ f\brackets[\big]{\omega(u,r), \omega(v,s)} = \rinv{r} f(u,v) s \]
	for all $ r,s \in \ring $ and $ u,v \in \module $. It is called \defemph*{hermitian with respect to $ \rinvmap $}\index{hermitian form} if, in addition, it satisfies $ \rinv{f(u,v)} = f(v,u) $ for all $ u,v \in \module $. It is called \defemph*{skew-hermitian with respect to $ \rinvmap $}\index{skew-hermitian form} if it satisfies $ \rinv{f(u,v)} = -f(v,u) $. The zero map $ \map{}{\module \times \module}{\nucleus(\ring)}{}{} $ is called the \defemph*{trivial sesquilinear form}.\index{sesquilinear form!trivial}
\end{definition}

\begin{note}\label{pseud:sesqui-nuclear-note}
	The requirement in \cref{pseud:sesqui-def} that the image of $ f $ lies in the nucleus of $ \ring $ is necessary because otherwise, the formula $ f(ur, vs) = \rinv{r} f(u,v) s $ is not well-defined. Alternatively, we could define a sesquilinear form on $ \module $ to be a bi-additive map $ \map{f}{\module \times \module}{\ring}{}{} $ which satisfies the two separate formulas $ f(ur,v) = \rinv{r} f(u,v) $ and $ f(u,vs) = f(u,v) s $ for all $ u,v \in \module $ and all $ r,s \in \ring $. If $ f $ is a sesquilinear form in this sense, then
	\begin{align*}
		\brackets[\big]{\rinv{r} f(u,v)} s &= f(ur,v)s = f(ur,vs) = \rinv{r} f(u,vs) = \rinv{r} \brackets[\big]{f(u,v)s}
	\end{align*}
	for all $ u,v \in \module $ and $ r,s \in \ring $. In other words, the image of $ f $ necessarily lies in the middle nucleus of $ \ring $. For alternative rings, we know from \cref{ring:alternative-nucleus} that the middle nucleus coincides with the nucleus, so the two definitions of sesquilinear forms are equivalent.
\end{note}

\begin{definition}[Pseudo-quadratic module]\label{pseud:pseudform-def}
	Let $ (\ring, \ringzero, \rinvmap) $ be a pre-involutory set and let $ \module $ be a right module over $ \ring $. A \defemph*{pseudo-quadratic form on $ \module $ with respect to $ (\ring, \ringzero, \rinvmap) $}\index{pseudo-quadratic form} is a map $ \map{q}{\module}{\nucleus(\ring)}{}{} $ for which there exists a skew-hermitian form $ f $ on $ \module $ such that the following properties are satisfied:
	\begin{stenumerate}
		\item $ q(u+v) \equiv q(u) + q(v) + f(u,v) \pmod \ringzero $ for all $ u,v \in \module $.
		
		\item \label{pseud:pseudform-def:scalar}$ q(ur) \equiv \rinv{r} q(u) r \pmod \ringzero  $ for all $ u \in \module $, $ r \in \ring $.
	\end{stenumerate}
	A map $ f $ with these properties is called a \defemph*{skew-hermitian form associated to $ q $}\index{pseudo-quadratic form!associated skew-hermitian form}. The form $ q $ is called \defemph*{anisotropic}\index{pseudo-quadratic form!anisotropic} if, in addition, the following axiom is satisfied:
	\begin{stenumerate}
		\setcounter{enumi}{2}
		
		\item $ q(u) \equiv 0 \pmod \ringzero $ implies $ u=0_\module $ for all $ u \in \module $.
	\end{stenumerate}
	Two pseudo-quadratic forms $ q,q' $ are called \defemph*{equivalent}\index{pseudo-quadratic form!equivalence} if $ q(v) \equiv q(v') \pmod \ringzero $ for all $ v \in \module $. A pseudo-quadratic form is called \defemph*{trivial}\index{pseudo-quadratic form!trivial} if it is equivalent to the zero map.
	A \defemph*{pseudo-quadratic module (over $ (\ring, \ringzero, \rinvmap) $)}\index{pseudo-quadratic module} is a triple $ (\module, q, f) $ consisting of right module $ \module $ over $ \ring $, a pseudo-quadratic form $ q $ on $ \module $ and a skew-hermitian form $ f $ associated to $ q $. If $ \ring $ is an associative division ring, these objects will also be called \defemph*{pseudo-quadratic spaces}.\index{pseudo-quadratic space} A pseudo-quadratic module $ (\module, q, f) $ is called \defemph*{anisotropic}\index{pseudo-quadratic module!anisotropic} if $ q $ is anisotropic.
\end{definition}

\begin{remark}\label{pseud:tw-rem}
	We have already seen in \cref{rinv:r0-equals-r} that a pre-involutory set can satisfy $ \ring = \ringzero $ only under very specific assumptions. Assume that this is the case. Then all axioms in \cref{pseud:pseudform-def} are trivially satisfied, so every map $ \map{q}{\module}{\ring}{}{} $ is a pseudo-quadratic form. Further, the only anisotropic pseudo-quadratic module in this situation is the zero module.
	
	In \cite{MoufangPolygons}, all pseudo-quadratic modules of interest are anisotropic. It follows from previous remarks that if $ (\module, q, f) $ is an anisotropic pseudo-quadratic module over $ (\ring, \ringzero, \rinvmap) $, then either $ \module = \compactSet{0} $ or $ \ring \ne \ringzero $.
\end{remark}

\begin{note}
	In \cite[(11.19)]{MoufangPolygons}, it is shown that in any pseudo-quadratic module $ (\module, q, f) $ over an associative involutory division set $ (\ring, \ringzero, \rinvmap) $ with $ \ring \ne \ringzero $, the map $ f $ is uniquely determined by $ q $. (We will repeat the proof of this fact in \cref{pseud:divring-f-unique}.) Using \cref{pseud:tw-rem}, it follows that the same is true for every anisotropic pseudo-quadratic module. As a consequence, the map $ f $ is not part of the structure of a pseudo-quadratic space in the definition in \cite[(11.17)]{MoufangPolygons}.
\end{note}

\begin{example}
	For any $ \ring $-module $ \module $, the zero map from $ \module $ to $ \ring $ is a pseudo-quadratic form. Thus $ (\module, 0, 0) $ is a pseudo-quadratic module.
\end{example}

\begin{example}\label{pseud:example}
	Let $ (\ring, \ringzero, \rinvmap) $ be an  involutory set and let $ n \in \Npos $. Assume that there exists $ h \in \nucleus(\ring) $ such that $ \rinv{h} = -h $, and choose such an element. We define the following maps:
	\begin{align*}
		\map{f}{\ring^n \times \ring^n}{\ring&}{(u,v)}{2 \sum_{i=1}^n \rinv{u_i} h v_i}, \\
		\map{q}{\ring^n}{\ring&}{u}{\sum_{i=1}^n \rinv{u_i} h u_i \; \brackets[\big]{\mathord{} = \frac{1}{2} f(u,u) \text{ if } 2_\ring \text{ is invertible}}}.
	\end{align*}
	Note that if $ 2_\ring = 0_\ring $, then $ f $ is the zero map, but $ q $ is not necessarily trivial. It is clear that $ f $ is skew-hermitian with respect to $ \rinvmap $. We want to show that $ q $ is pseudo-quadratic with associated pseudo-quadratic form $ f $. Let $ u,v \in \ring^n $ and let $ r \in \ring $. It is clear that
	\begin{align*}
		q(ur) &= \sum_{i=1}^n \rinv{(u_i r)} h u_i r = \rinv{r} \brackets*{\sum_{i=1}^n \rinv{u_i} h u_i} r = \rinv{r} q(u) r.
	\end{align*}
	(Note that this equality holds in $ \ring $, not merely modulo $ \ringzero $.) Further,
	\begin{align*}
		q(u+v) &= \sum_{i=1}^n \brackets[\big]{\rinv{u_i}h u_i + \rinv{v_i} hu_i + \rinv{u_i} hv_i + \rinv{v_i} hv_i} \\
		&= q(u) + q(v) + \sum_{i=1}^n \brackets[\big]{\rinv{v_i} hu_i + \rinv{u_i} hv_i}.
	\end{align*}
	By the choice of $ h $, we have $ \rinv{u_i} h v_i = -\rinv{(\rinv{v_i} h u_i)} $ for all $ i \in \numint{1}{n} $. Put $ r_i \defl \rinv{v_i} h u_i $ for all $ i \in \numint{1}{n} $. Then we see that
	\begin{align*}
		q(u+v) &= q(u) + q(v) + \sum_{i=1}^n (r_i - \rinv{r_i}) = q(u) + q(v) + 2\sum_{i=1}^n r_i - \sum_{i=1}^n (r_i + \rinv{r_i}) \\
		&= q(u) + q(v) + f(u,v) - \sum_{i=1}^n (r_i + \rinv{r_i}).
	\end{align*}
	Since $ \ringzero $ contains all traces, it follows that $ q(u+v) \equiv q(u) + q(v) + f(u,v) \pmod{\ringzero} $.
\end{example}

\begin{example}\label{pseud:example:C}
	We can take $ (\ring, \ringzero, \rinvmap) \defl (\IC, \IR, \rinvmap) $ (where $ \rinvmap $ denotes complex conjugation), $ h \defl i $ and $ n \defl 1 $ in \cref{pseud:example}. Then for all $ a,b,c,d \in \IR $, we have
	\begin{align*}
		f(a+bi, c+di) &= 2 (a-bi) i (c+di) = 2(b+ai)(c+di) \\
		&= 2(bc-ad) + 2(bd+ac)i, \\
		q(a+bi) &= \frac{1}{2} f(a+bi, a+bi) = (a^2+b^2) i.
	\end{align*}
	That is, $ q $ is precisely the usual norm function on $ \IC $, but multiplied with $ i $.
\end{example}

The fact that the images of $ q $ and $ f $ need to lie in the nucleus of $ \ring $ illustrates that pseudo-quadratic modules do not really belong to the nonassociative world. However, we can still obtain examples of pseudo-quadratic modules over alternative rings which are not associative using direct sum constructions.

\begin{remark}[Direct sums of pseudo-quadratic modules]\label{pseud:dirsum}
	Let $ X=(\ring, \ringzero, \rinvmap) $ and $ X'=(\ring', \ringzero', \rinvmap') $ be two involutory sets and let $ (\module, q, f) $ and $ (\module', q', f') $ be pseudo-quadratic modules over $ X $ and $ X' $, respectively. We have seen in \cref{rinvset:dirsum} that $ X \dirsum X' \defl (\ring \dirsum \ring', \ringzero \dirsum \ringzero', \rinvmap \dirsum \rinvmap') $ is also an involutory set. We can equip the group $ \module \dirsum \module' $ with an $ (\ring \dirsum \ring') $-module structure by putting
	\[ (u,u') \cdot (r,r') \defl (ur, u'r') \]
	for all $ u \in \module $, $ u' \in \module' $, $ r \in \ring $ and $ r' \in \ring' $.
	It is now easy to verify that $ (\module \dirsum \module', q \dirsum q', f \dirsum f') $ is a pseudo-quadratic module over $ X \dirsum X' $. This pseudo-quadratic module is anisotropic if and only if both $ (\module, q, f) $ and $ (\module', q', f') $ are anisotropic.
\end{remark}

\begin{note}[Pseudo-quadratic modules over alternative rings]\label{pseud:alternative-note}
	Let everything be as in \cref{pseud:dirsum}, and assume in addition that $ (\module', q', f') = (\compactSet{0}, 0, 0) $. Then $ \module \dirsum \module' $ can be identified with $ \module $, so \cref{pseud:dirsum} yields that the pseudo-quadratic module $ (\module, q, f) $ over $ (\ring, \ringzero, \rinvmap) $ can also be regarded as a pseudo-quadratic module over $ (\ring, \ringzero, \rinvmap) \dirsum (\ring', \ringzero', \rinvmap') $. Here $ (\ring', \ringzero', \rinvmap') $ can be an arbitrary involutory set. Taking $ (\ring', \ringzero', \rinvmap') $ to not be associative, we conclude that every pseudo-quadratic module can be regarded as a pseudo-quadratic module over an involutory set which is not associative.
	
	While the previous example shows that pseudo-quadratic modules over rings which are not associative exist, it should be regarded as a technicality: Only the associative ideal $ \ring \dirsum \compactSet{0} $ of the alternative ring $ \ring \dirsum \ring' $ contributes to the (pseudo-quadratic) module structure of $ \module $. We would like to say that, while $ \ring \dirsum \ring' $ is alternative but not associative, it is not \enquote{purely alternative} either. In this example, the ring splits into an alternative part and an associative, but we cannot expect a similar decomposition to exist for arbitrary alternative rings. Thus it is not at all clear what a reasonable definition of \enquote{purely alternative rings} should be.
	
	Following \cite{Slater_NuclCentAlt}, we will present a reasonable abstract notion of pure alternativity for rings in \cref{sec:pure-alt}. We will show that purely alternative rings do not admit any non-zero modules, so in particular, they do not admit non-zero pseudo-quadratic modules. This fact substantiates our previous claim that \enquote{pseudo-quadratic modules do not belong to the nonassociative world}. However, we emphasise that we are not aware of a theorem which allows us to decompose an arbitrary alternative ring into an associative part and a purely alternative part. Thus there is still some interest in allowing the underlying involutory set $ (\ring, \ringzero, \rinvmap) $ to not be associative.
\end{note}

\subsection{Basic Properties}

\begin{secnotation}\label{pseud:conv}
	From now on, we denote by $ (\ring, \ringzero, \rinvmap) $ a pre-involutory set, by $ \module $ a right module over $ \ring $ and by $ q $ a pseudo-quadratic form with respect to $ (\ring, \ringzero, \rinvmap) $. Further, unless otherwise specified, the symbol \enquote{$ \equiv $} will always be used to denote congruence modulo $ \ringzero $.
\end{secnotation}

\begin{remark}\label{pseud:coset}
	Let $ \map{h}{\module}{\ringzero}{}{} $ be any map. Then $ q+h $ is also a pseudo-quadratic form on $ \module $ because all axioms are defined modulo $ \ringzero $. Further, a skew-hermitian form $ \map{f}{\module \times \module}{\ring}{}{} $ is associated to $ q+h $ if and only if it is associated to $ q $. Thus the equivalence class of $ q $ can be identified with the induced map $ \map{\tilde{q}}{\module}{\nucleus(\ring) / \ringzero}{}{} $. (In fact, this is exactly the approach that is taken in \cite[5.1D]{HahnOMeara-ClassicalGroups}.) However, it will be important in some contexts that we have chosen a fixed coset representative $ q(v) $ in $ \ring $ for every $ v \in \module $. See, for example, \cref{pseud:T-equiv-note}.
\end{remark}

\begin{remark}\label{pseud:0}
	It follows from Axiom~\thmitemref{pseud:pseudform-def}{pseud:pseudform-def:scalar} that
	\begin{align*}
		q(0_\module) = q(0_\module 0_\ring) \equiv \rinv{0_\ring} q(0_\module) 0_\ring = 0_\ring.
	\end{align*}
\end{remark}

\begin{remark}\label{pseud:modulo-mult}
	Let $ r,s \in \ring $ such that $ r \equiv s $ and let $ t \in \ring $. Since $ \rinv{t} \ringzero t \subs \ringzero $, it follows that $ \rinv{t} r t \equiv \rinv{t} st $. If $ \ringzero t \subs \ringzero $, then we also have $ rt \equiv st $, but this is not true in general. Similarly, all $ t \in \ring $ such that $ t \ringzero \subs \ringzero $ have the property that $ r \equiv s $ implies $ ts \equiv ts $. For example, these properties hold for all $ t $ in the image of $ \IZ $ in~$ \ring $.
\end{remark}

The following result is the analogue of \cref{quadmod:basiclem,quadmod:symbil-bij} for pseudo-quadratic modules.

\begin{lemma}\label{pseud:2}
	Let $ f $ be a skew-hermitian form associated to $ q $. Then the following hold:
	\begin{lemenumerate}
		\item $ q(0_\module) \equiv 0_\ring $.
		
		\item $ f(v,v) \equiv 2q(v) $ for all $ v \in \module $.
		
		\item $ \rinv{r} \equiv -r $ for all $ r \in \ring $.
		
		\item \label{pseud:2:f-sym}$ f(v,u) \equiv f(u,v) $ for all $ u,v \in \module $.
		
		\item \label{pseud:2:inv}If $ 2_\ring $ is invertible, then the map $ \map{q'}{\module}{\ring}{v}{\frac{f(v,v)}{2}} $ is a pseudo-quadratic form which is equivalent to $ q $ and which satisfies $ q'(vr) = \rinv{r} q(v) r $ (in $ \ring $, not just modulo $ \ringzero $) for all $ v \in \module $ and all $ r \in \ring $.
		
		\item \label{pseud:2:inv-triv}If $ 2_\ring $ is invertible and $ f $ is trivial, then $ q $ is trivial as well.
	\end{lemenumerate}
\end{lemma}
\begin{proof}
	First of all, note that
	\begin{align*}
		q(0_\module) = q(0_\module + 0_\module) \equiv q(0_\module) + q(0_\module) + f(0_\module, 0_\module) = 2q(0_\module).
	\end{align*}
	This implies that $ q(0_\module) \equiv 0_\ring $. In a similar way, we have for all $ v \in \module $ that
	\begin{align*}
		q(2v) \equiv \rinv{2} q(v) 2 = 4q(v) \midand q(2v) = q(v+v) \equiv 2q(v) + f(v,v).
	\end{align*}
	It follows that $ 2q(v) \equiv f(v,v) $. Since $ \ringzero $ contains all traces, we have $ r + \rinv{r} \equiv 0 $ for all $ r \in \ring $. In other words, $ \rinv{r} \equiv -r $ for all $ r \in \ring $. This implies that $ f(u,v) \equiv f(v,u) $ for all $ u,v \in \module $ because $ f $ is skew-hermitian.
	
	Now assume that $ 2_\ring $ is invertible and define $ \map{q'}{\module}{\ring}{v}{\frac{f(v,v)}{2}} $. Then by~\itemref{pseud:2:f-sym} and the definition of $ q' $,
	\begin{align*}
		q'(v+u) &= \frac{f(v+u,v+u)}{2} = \frac{f(v,v)}{2} + \frac{f(u,u)}{2} + \frac{f(v,u)}{2} + \frac{f(u,v)}{2} \\
		&= q'(v)+q'(u) + \frac{f(v,u)}{2} + \frac{f(u,v)}{2} \equiv q'(v)+q'(u) + f(v,u)
	\end{align*}
	for all $ u,v \in \module $. Further,
	\begin{align*}
		q'(vr) &= \frac{f(vr,vr)}{2} = \frac{\rinv{r} f(v,v) r}{2} = \rinv{r} q(v) r
	\end{align*}
	for all $ v \in \module $ and $ r \in \ring $. Thus~\itemref{pseud:2:inv} holds. Now assume in addition that $ f $ is trivial. Then $ q' $ is trivial by definition. Since $ q' $ is equivalent to $ q $, it follows that $ q $ is trivial as well, which finishes the proof.
\end{proof}

In analogy to the orthogonal group of a quadratic module (\cref{quadmod:ortho-grp}), we can define the unitary group of a pseudo-quadratic module.

\begin{definition}[Unitary group, {\cite[5.2A]{HahnOMeara-ClassicalGroups}}]\label{pseud:unitary-def}
	The group
	\[ \Unitary(\module) \defl \Set[\big]{\phi \in \Aut_\comring(\module) \given q\brackets[\big]{\phi(v)} \equiv q(v) \text{ and } f\brackets[\big]{\phi(u), \phi(v)} \text{ for all } u,v \in \module} \]
	is called \defemph*{unitary group of $ \module = (\module, q, f) $}\index{unitary group}. We use the convention that it acts on $ \module $ from the left side, so that the composition $ \phi \circ \psi $ of $ \phi, \psi \in \Ortho(q) $ is the map $ \map{}{}{}{x}{\phi(\psi(x))} $.
\end{definition}

\begin{remark}\label{pseud:equiv-nucl}
	In the following, we will often consider elements $ h \in \ring $ with the property that $ h \equiv q(v) $ for some $ v \in \module $. Since $ \ringzero $ and the image of $ q $ are, both by definition, contained in the nucleus of $ \ring $, all such elements are nuclear.
\end{remark}

An important result about pseudo-quadratic spaces over associative involutory set is \cite[(11.19)]{MoufangPolygons}. It is not true in the present generality. The following \cref{pseud:fvsr-lem,pseud:fvv-equiv-qv} are the parts which remain valid in our general setting. \Cref{pseud:divring-fvv,pseud:divring-f-unique} are exactly the assertion of \cite[(11.19)]{MoufangPolygons}, except that the assumption on $ \ring $ to be associative turns out to be unnecessary. We will not have occasion to apply any of these results, but we state them for the record.

\begin{lemma}\label{pseud:fvsr-lem}
	Let $ v \in \module $ and let $ r,s,h \in \ring $ such that $ h \equiv q(v) $. Then $ f(vs, vr) \equiv \rinv{s}hr + \rinv{r}hs $.
\end{lemma}
\begin{proof}
	On the one hand,
	\begin{align*}
		q\brackets[\big]{v(s+r)} &\equiv q(vs) + q(vr) + f(vs,vr) \equiv \rinv{s}q(v)s + \rinv{r}q(v)r + f(vs,vr).
	\end{align*}
	On the other hand,
	\begin{align*}
		q\brackets[\big]{v(s+r)} &\equiv \rinv{(s+r)} q(v) (s+r) \equiv \rinv{s} q(v) s + \rinv{r} q(v) r + \rinv{s} q(v) r + \rinv{r} q(v) s.
	\end{align*}
	This implies that the assertion holds for $ h = q(v) $. In general, we have $ h= q(v) +t $ for some $ t \in \ringzero $, and it remains to show that $ \rinv{s}tr + \rinv{r}ts \equiv 0 $. This is exactly the statement of \cref{rinv:sum-in-R0-lem}, so the proof is finished.
\end{proof}

\begin{lemma}\label{pseud:fvv-equiv-qv}
	Let $ v \in \module $ and let $ h \in \ring $ such that $ h \equiv q(v) $. Then
	\[ \brackets[\big]{f(v,v) - h + \rinv{h}}r \equiv 0 \quad \text{for all } r \in \ring. \]
	In particular, $ f(v,v) \equiv h - \rinv{h} \equiv q(v) - \rinv{q(v)} $.
\end{lemma}
\begin{proof}
	Putting $ s \defl 1_\ring $ in \cref{pseud:fvsr-lem} yields $ f(v,v)r - hr - \rinv{r}h \equiv 0 $. Since $ \rinv{r}h \equiv -\rinv{(\rinv{r}h)} = -\rinv{h}r $, the assertion follows.
\end{proof}

Observe that the assertion of \cref{pseud:fvv-equiv-qv} is stronger than the special case $ f(v,v) \equiv h - \rinv{h} $ by \cref{pseud:modulo-mult}. However, for division rings, we have the following result.

\begin{lemma}\label{pseud:divring-fvv}
	Assume that $ \ring $ is a division ring and that $ \ringzero \ne \ring $. Then for all $ v \in \module $ and all $ h \in \ring $ with $ h \equiv q(v) $, we have $ f(v,v) = h - \rinv{h} $ (in $ \ring $, not just modulo $ \ringzero $). In particular, $ f(v,v) = q(v) - \rinv{q(v)} $ for all $ v \in \module $.
\end{lemma}
\begin{proof}
	For any $ v \in \module $ and any $ h \in q(v) + \ringzero $, it follows from \cref{pseud:fvv-equiv-qv,pseud:divring-R0} that $ f(v,v)- h + \rinv{h} $ is zero. Since $ q(v) $ lies in $ q(v) + \ringzero $, the second assertion follows.
\end{proof}

\begin{lemma}\label{pseud:divring-f-unique}
	Assume that $ \ring $ is a division ring and that $ \ringzero \ne \ring $. Then there exists exactly one skew-hermitian form which is associated to $ q $.
\end{lemma}
\begin{proof}
	Assume that $ f,f' $ are two such forms. Then $ g \defl f-f' $ is a sesquilinear form whose image is contained in $ \ringzero $. For all $ u,v \in \module $ and $ r \in \ring $, it follows that $ g(u,v)r = g(u,vr) \in \ringzero $. By \cref{pseud:divring-R0}, this implies that $ g(u,v) = 0_\ring $ for all $ u,v \in \ring $, so $ f=f' $.
\end{proof}

Motivated by \cref{pseud:divring-fvv}, we introduce the following definition.

\begin{definition}[Standard pseudo-quadratic module]
	A pseudo-quadratic module $ (\module, q, f) $ over $ (\ring, \ringzero, \rinvmap) $ is called \defemph*{standard}\index{pseudo-quadratic module!standard} if $ f(v,v) = q(v) - \rinv{q(v)} $ for all $ v \in \module $.
\end{definition}

\begin{note}
	Being standard is a rather weak condition which is automatically satisfied in many situations:
	\begin{remenumerate}
		\item We have seen in \cref{pseud:divring-fvv} that every pseudo-quadratic module over a division ring with $ \ring \ne \ringzero $ is standard. Further, the zero module is clearly standard as well. By \cref{pseud:tw-rem}, this means that every pseudo-quadratic module that appears in \cite{MoufangPolygons} is standard.
		
		\item Assume that $ 2_\ring $ is invertible and let $ (\module, q, f) $ be a pseudo-quadratic module. It follows from \thmitemcref{pseud:2}{pseud:2:inv} that $ q $ is equivalent to a pseudo-quadratic form $ q' $ such that $ (\module, q', f) $ is a standard pseudo-quadratic module.
	\end{remenumerate}
	We will see that the standard condition is necessary to construct a Jordan module from a pseudo-quadratic module (\cref{BC:jordanmodule-pseudquad-ex}). Further, every pseudo-quadratic module which arises from a Jordan module (that is, from a root grading of type $ C $ or $ BC $) is standard (\cref{BC:jordanmodule:class-pseud}). We conclude that standardness is the right condition to extend the theory of pseudo-quadratic modules to our general setting.
\end{note}

\subsection{The Group \texorpdfstring{$ \psgr(M) $}{T(M)}}

\begin{secnotation}
	In addition to \cref{pseud:conv}, we will from now on fix a skew-hermitian form $ f $ such that $ (\module, q, f) $ is a pseudo-quadratic module over $ (\ring, \ringzero, \rinvmap) $.
\end{secnotation}

For any pseudo-quadratic module $ \module $, we can define a group $ \psgr(\module) $. If $ \module $ is standard, this group can be equipped with the structure of a Jordan module (\cref{BC:jordanmodule-pseudquad-ex}). In \cref{sec:jordan-class}, we will show that every Jordan module arises in this way if $ 2_\ring $ is invertible.

\begin{definition}\label{pseud:T-def}
	We define
	\[ \psgr(\module) \defl \psgr(\module, q, f) \defl \Set{(u,h) \in \module \times \ring \given q(u) \equiv h} \]
	and $ (u,h) \cdot (v,k) \defl \brackets[\big]{u+v, h+k + f(u,v)} $ for all $ (u,h), (v,k) \in \psgr(q) $.
\end{definition}

Observe that the multiplication on $ \psgr(\module) $ depends on the choice of the skew-hermitian form $ f $. Further, this multiplication is in general not abelian because $ f $ is, in general, not symmetric.

\begin{remark}\label{pseud:T-set-bij}
	We have a canonical bijection
	\[ \map{}{\module \times \ringzero}{\psgr(\module)}{(u,h)}{\brackets[\big]{u, q(u)+h}} \]
	of sets which, in general, is not a group homomorphism with respect to the canonical group structure on $ \module \times \ringzero $.
\end{remark}

\begin{note}\label{pseud:T-equiv-note}
	If $ q' $ is another pseudo-quadratic form which is equivalent to $ q $, then $ \psgr(\module, q, f) $ and $ \psgr(\module, q', f) $ are not only isomorphic, but in fact the same set with the same group multiplication. However, the bijection in \cref{pseud:T-set-bij} depends on the choice of $ q $.
\end{note}

\begin{remark}\label{pseud:T-q-triv}
	If $ q $ is trivial, then $ \psgr(\module) = \module \times \ringzero $ as sets and the bijection in \cref{pseud:T-set-bij} is the identity map. However, the multiplication on $ \psgr(\module) $ is not necessarily the component-wise addition on $ \module \times \ringzero $. This is only the case if $ f $ is trivial as well. Note that if $ 2_\ring $ is invertible, then the triviality of $ f $ implies the triviality of $ q $ by \thmitemcref{pseud:2}{pseud:2:inv-triv}.
\end{remark}

\begin{remark}\label{pseud:T-R0-embed}
	The map $ \map{}{\ringzero}{\psgr(\module)}{r_0}{(0_\module, r_0)} $ is an injective homomorphism of square-modules. It is well-defined because $ q(0_\module) \equiv 0_\ring $ by \cref{pseud:0}.
\end{remark}

\begin{remark}\label{pseud:T-comm}
	Let $ (u,h), (v,k) \in \psgr(\module) $. Then
	\begin{align*}
		(u,h) \cdot (v,k) &= (v,k) \cdot (u,h) \cdot \brackets[\big]{0, f(u,v) - f(v,u)}
	\end{align*}
	where $ f(u,v) - f(v,u) = f(u,v) + \rinv{f(u,v)} = \ringTr_\rinvmap(f(u,v)) $.
\end{remark}

\begin{lemma}\label{pseud:T-group}
	The multiplication on $ \psgr(\module) $ defines a group structure with identity element $ (0_\module,0_\ring) $ and inverses given by
	\[ (u,h)^{-1} = \brackets[\big]{-u, f(u,u) - h} \]
	for all $ (u,v) \in \psgr(q) $.
\end{lemma}
\begin{proof}
	For all $ (u,h), (v,k) \in \psgr(\module) $, it follows from the axioms of a pseudo-quadratic form that
	\[ q(u+v) \equiv q(u) + q(v) + f(u,v) \equiv h + k + f(u,v), \]
	so $ (u+v, h+k+f(u,v)) $ lies in $ \psgr(q) $. We infer that the multiplication on $ \psgr(\module) $ is well-defined. Further, it is associative because for any $ (u_1, h_1), (u_2, h_2), (u_3, h_3) \in \psgr(\module) $, both $ \brackets[\big]{(u_1, h_1) (u_2, h_2)} (u_3, h_3) $ and $ (u_1, h_1) \brackets[\big]{(u_2, h_2) (u_3, h_3)} $ equal
	\begin{align*}
		\brackets[\big]{u_1 + u_2 + u_3, h_1 + h_2 + h_3 + f(u_1, u_2) + f(u_1, u_3) + f(u_2, u_3)}.
	\end{align*}
	It is clear that $ (0,0) $ is an identity element in $ \psgr(q) $. Finally, for all $ (u,h) \in \psgr(q) $, we have
	\begin{align*}
		(u, h) \brackets[\big]{-u, f(u,u) - h} &= \brackets[\big]{0, h -h + f(u,u) - f(u,u)} = (0,0) \\
		&= \brackets[\big]{-u, f(u,u) - h} (u,h),
	\end{align*}
	so $ (u,h) $ is invertible with inverse $ \brackets[\big]{-u, f(u,u) - h} $. This finishes the proof.
\end{proof}

\begin{remark}\label{pseud:T-group:standard}
	Note that, if $ (\module, q, f) $ is standard, then the inverse of $ (u,h) \in \psgr(q) $ is given by $ (-u, -\rinv{h}) $.
\end{remark}

The following observation will be used to define (non-balanced) Weyl elements in root gradings of type $ BC $. See \cref{BC:ex-weyl-def}.

\begin{lemma}\label{pseud:T-weyl-inverse}
	Let $ (\module, q, f) $ be a pseudo-quadratic module over $ (\ring, \ringzero, \rinvmap) $ and let $ (u,h) \in \psgr(\module) $ such that $ h $ is invertible. Then $ (uh^{-1}, \rinvmin{h}) $ and $ (-u \rinvmin{h}, \rinvmin{h}) $ lie in $ \psgr(\module) $.
\end{lemma}
\begin{proof}
	Since $ (u,h) $ lies in $ \psgr(\module) $, there exists $ r_0 \in \ringzero $ such that $ q(u) = h+r_0 $. Thus
	\begin{align*}
		q(uh^{-1}) &\equiv \rinvmin{h} q(u) h^{-1} = \rinvmin{h} + \rinvmin{h} r_0 h^{-1} \equiv \rinvmin{h} \rightand \\
		q(-u\rinvmin{h}) &\equiv h^{-1} q(u) \rinvmin{h} = \rinvmin{h} + h^{-1} r_0 \rinvmin{h} \equiv \rinvmin{h}.
	\end{align*}
	This proves the claim.
\end{proof}

We now define several maps on $ \psgr(\module) $ which will later be seen to turn $ \psgr(\module) $ into a Jordan module over $ (\ring, \rinvmap) $. As a first step, we define a square-scalar multiplication on $ \psgr(\module) $.

\begin{lemma}\label{pseud:T-squaremod}
	Define
	\[ \map{\jorsc}{\psgr(\module) \times \ring}{\psgr(\module)}{\brackets[\big]{(v, h), r}}{(vr, \rinv{r}hr)}. \]
	Then $ (\psgr(\module), \jorsc) $ is a well-defined multiplicative square-module.
\end{lemma}
\begin{proof}
	For all $ (v,h) \in \psgr(\module) $ and all $ r \in \ring $, it follows from Axiom~\thmitemref{pseud:pseudform-def}{pseud:pseudform-def:scalar} and \cref{pseud:modulo-mult} that
	\begin{align*}
		q(vr) \equiv \rinv{r} q(v)r \equiv \rinv{r}hr,
	\end{align*}
	so that $ (vr, \rinv{r}hr) $ lies in $ \psgr(\module) $. Thus $ \jorsc $ is well-defined. It is clear that $ \jorsc(t, 1_\ring) = t $ for all $ t \in \psgr(\module) $. For all $ r,s \in \ring $ and all $ (v,h) \in \psgr(\module) $, we have $ (vr)s = v(rs) $ by Axiom~\thmitemref{ring:module}{ring:module:assoc} and $ \rinv{s} (\rinv{r} hr) s = (\rinv{s} \rinv{r}) h (rs) $ by \cref{rinv:artin-inv,pseud:equiv-nucl}. Hence $ \jorsc $ is multiplicative.
	
	To see that $ \jorsc $ is weakly quadratic in the second component, let $ (v,h) \in \psgr(\module) $ and let $ r,s \in \ring $. By the formula for inverses in $ \psgr(\module) $ from \cref{pseud:T-group}, we have
	\begin{align*}
		\jorsc\brackets[\big]{(v,h),s}^{-1} \jorsc\brackets[\big]{(v,h),r}^{-1} \hspace{-3cm}& \\
		&= (vs, \rinv{s}hs)^{-1} (vr, \rinv{r}hr)^{-1} \\
		&= \brackets[\big]{-vs, f(vs, vs) - \rinv{s}hs} \brackets[\big]{-vr, f(vr, vr) - \rinv{r}hr} \\
		&= \brackets[\big]{-vs-vr, f(vs,vs) - \rinv{s}hs + f(vr,vr) - \rinv{r}hr + f(vs,vr)}.
	\end{align*}
	Further,
	\begin{align*}
		\jorsc\brackets[\big]{(v,h), r+s} &= (vr+vs, \rinv{r}hr + \rinv{r}hs + \rinv{s}hr + \rinv{s}hs).
	\end{align*}
	Thus the polarisation of $ \jorsc $ in the second component is
	\begin{align*}
		\jorsc\brackets[\big]{(v,h),s}^{-1} \jorsc\brackets[\big]{(v,h),r}^{-1} \jorsc\brackets[\big]{(v,h), r+s} \hspace{-6cm}& \\
		&=\big(0_\module, f(vs,vs) - \rinv{s}hs + f(vr,vr) - \rinv{r}hr + f(vs,vr) \\
		&\hspace{3cm} \mathord{}+ \rinv{r}hr + \rinv{r}hs + \rinv{s}hr + \rinv{s}hs - f(vs+vr, vs+vr) \big) \\
		&=\brackets[\big]{0_\module, \rinv{r}hs + \rinv{s}hr - f(vr,vs)}.
	\end{align*}
	This expression is bi-additive in $ (r,s) $, so $ \jorsc $ is weakly quadratic in the second component.
	Finally, we have for all $ (v,h), (v',h') \in \psgr(\module) $ and $ r \in \ring $ that
	\begin{align*}
		\jorsc\brackets[\big]{(v,h) (v',h'), r} &= \jorsc\brackets[\big]{\brackets{v+v', h+h' + f(v,v')}, r} \\
		&= \brackets[\big]{vr + v'r, \rinv{r}h r + \rinv{r}h'r + \rinv{r} f(v,v')r} \\
		&= \brackets[\big]{vr + v'r, \rinv{r}h r + \rinv{r}h'r + f(vr,v'r)} \\
		&= (vr, \rinv{r}hr) (v'r, \rinv{r}h'r) = \jorsc\brackets[\big]{(v,h), r} \jorsc\brackets[\big]{(v',h'), r}.
	\end{align*}
	Hence $ \jorsc $ is additive in the first component. This finishes the proof.
\end{proof}

\begin{remark}\label{pseud:T-squaremod-standard}
	If $ (\module, q, f) $ is standard, then the polarisation of $ \jorsc $ in \cref{pseud:T-squaremod} in the second component is
	\begin{align*}
		\brackets[\big]{0_\module, \rinv{r}hs + \rinv{s}hr - f(vr,vs)} = \brackets[\big]{0_\module, \rinv{s}hr + \rinv{s}\rinv{h}s} = \brackets[\big]{0_\module, \ringTr_\rinvmap(\rinv{s}hr)}.
	\end{align*}
	Here the \enquote{polarisation in the second component} is defined as in~\ref{quadmod:polarisation-def}.
\end{remark}

\begin{observation}\label{pseud:T-invo}
	By \cref{weakquad:inv-def}, we have a canonical involution on the multiplicative square-module $ (\psgr(\module), \jorsc) $. It is precisely the map
	\[ \map{\modinvmap}{\psgr(\module)}{\psgr(\module)}{(u,h)}{(-u,h)}. \]
\end{observation}

As a next step, we define maps $ \jorprojone $ and $ \jorTrone $.

\begin{lemma}\label{pseud:T-proj}
	The projection map
	\[ \map{\jorprojone}{\psgr(\module)}{\nucleus(\ring)}{(v,h)}{h} \]
	is compatible with the square-scalar multiplication $ \jorsc $ from \cref{pseud:T-squaremod}.
\end{lemma}
\begin{proof}
	We denote the canonical square-scalar multiplication on $ \nucleus(\ring) $ by $ \omega $ (see \cref{jormod:invring-jordan,rinv:nucl-submod}). At first, observe that the image of $ \jorprojone $ is indeed contained in the nucleus of $ \ring $ by \cref{pseud:equiv-nucl}. Now for all $ t=(v,h) \in \psgr(\module) $ and for all $ r \in \ring $, we have
	\begin{align*}
		\jorprojone\brackets[\big]{\jorsc(t,r)} &= \jorprojone\brackets[\big]{(vr, \rinv{r}hr)} = \rinv{r}hr = \omega(h,r) = \omega\brackets[\big]{\jorprojone(t), r}.
	\end{align*}
	This finishes the proof.
\end{proof}

\begin{lemma}\label{pseud:T-Tr}
	The map
	\[ \map{\jorTrone}{\ring}{\psgr(\module)}{r}{\brackets[\big]{0_\module, \ringTr_\rinvmap(r)} = (0_\module, r+\rinv{r})}  \]
	is a homomorphism of square-modules.
\end{lemma}
\begin{proof}
	This map is the composition of the map $ \map{\ringTr_\rinvmap}{\ring}{\ringzero}{}{} $ from \cref{rinv:tr-submod} with the embedding $ \map{}{\ringzero}{\psgr(\module)}{}{} $ from \cref{pseud:T-R0-embed}. Since both maps are homomorphisms of square-modules, the assertion follows.
\end{proof}

We have now defined all maps that are needed to equip $ \module $ with the structure of a Jordan module. Before we turn to the definition of Jordan modules, we make a final observation.

\begin{remark}[Direct products and $ \psgr(\module) $]
	Write $ X \defl (\ring, \ringzero, \rinvmap) $, let $ X'=(\ring', \ringzero', \rinvmap') $ be another involutory sets and let $ (\module', q', f') $ be a pseudo-quadratic module over $ X' $. By \cref{pseud:dirsum}, we can equip $ \module \dirsum \module' $ with the structure of a pseudo-quadratic module over $ X \dirsum X' $. We then have
	\begin{align*}
		\psgr(\module \dirsum \module') = \Set*{\brackets[\big]{(u,u'), (h,h')} \given \begin{gathered}
			u \in \module, u' \in \module', h \in \ring, h' \in \ring', \\
			(q \dirsum q')(u,u') \equiv (h,h') \pmod{\ringzero \dirsum \ringzero'}
		\end{gathered}},
	\end{align*}
	which is clearly isomorphic to
	\begin{align*}
		\psgr(\module) \dirsum \psgr(\module') &= \Set*{\brackets[\big]{(u,h), (u', h')} \given \begin{gathered}
			u \in \module, u' \in \module', h \in \ring, h' \in \ring', \\
			q(u) \equiv h \pmod \ringzero, \\
			q'(u') \equiv h' \pmod{\ringzero'}
		\end{gathered}}.
	\end{align*}
	Under this identification, the maps $ \jorsc $, $ \jorprojone $ and $ \jorTrone $ defined on $ \module \dirsum \module' $ are precisely the respective sums of the maps $ \jorsc $, $ \jorprojone $, $ \jorTrone $ defined on $ \module $ and $ \module' $. In the (not yet introduced) terminology of Jordan modules, this says precisely that the identification between $ \psgr(\module \dirsum \module') $ and $ \psgr(\module) \dirsum \psgr(\module') $ is an isomorphism of Jordan modules.
\end{remark}


\section{Jordan Modules}

\label{sec:jordan-modules}

Finally, we can define to the protagonists of this chapter.

\subsection{Definition and Examples}

\begin{definition}[Jordan module]\label{BC:jordanmodule-def}
	Let $ \ring $ be an alternative ring with a nuclear involution $ \rinvmap $ and denote the square-scalar multiplication on $ \ring $ from \cref{jormod:invring-jordan} by $ \omega $. A \defemph*{Jordan module over $ (\ring, \rinvmap) $}\index{Jordan module} is a tuple $ \jormodtup = (\jormod, \jorsc, \jorprojone, \jorTrone, \psi) $ consisting of a multiplicative right square-module $ (\jormod, \jorsc) $ over $ \ring $, a map
	\[ \map{\jorprojone}{(\jormod, \jorsc)}{(\nucleus(\ring), \omega)}{}{} \]
	which preserves the square-scalar multiplication, a homomorphism
	\[ \map{\jorTrone}{(\ring, \omega)}{(\jormod, \jorsc)}{}{} \]
	of square-modules and a skew-hermitian form
	\[ \map{\psi}{\jormod \times \jormod}{\nucleus(\ring)}{}{} \]
	with respect to $ \rinvmap $ such that the following conditions are satisfied:
	\begin{stenumerate}
		\item \label{BC:jordanmodule-def:tr}$ \jorTrone(\rinv{a}) = \jorTrone(a) $ and $ \jorTrone\brackets[\big]{(ab)c} = \jorTrone\brackets[\big]{a(bc)} $ for all $ a,b,c \in \ring $. Thus we can omit brackets in the argument of $ \jorTrone $ without ambiguity.
		
		\item \label{BC:jordanmodule-def:proj-rinv}$ \rinv{\jorprojone(u)} = \jorprojone(u) - \psi(u,u) $ for all $ u \in \jormod $.
		
		\item \label{BC:jordanmodule-def:psi}$ \psi\brackets[\big]{\jorTrone(a), u} = 0 = \psi\brackets[\big]{u, \jorTrone(a)} $ for all $ a \in \ring $, $ u \in \jormod $.
		
		\item \label{BC:jordanmodule-def:lin}The following \enquote{linearisation properties} hold for all $ u,v \in \jormod $ and all $ a,b \in \ring $:
		\begin{align}
			\jorprojone(u \joradd v) &= \jorprojone(u) + \jorprojone(v) + \psi(u,v), \label{eq:jormod-lin:proj} \\
			\jorsc(u, a+b) &= \jorsc(u, a) \joradd \jorsc(u,b) \joradd \jorTrone\brackets[\big]{\rinv{b} \jorprojone(u) a}, \label{eq:jormod-lin:jorsc} \\
			u \joradd v &= v \joradd u \joradd \jorTrone\brackets[\big]{\psi(u,v)}. \label{eq:jormod-lin:add}
		\end{align}
		
		\item \label{BC:jordanmodule-def:square}For all $ r \in \ring $, we have $ \jorprojone(\jorTrone(r)) = r + \rinv{r} $.
		
		\item \label{BC:jordanmodule-def:inv}There exists $ v_0 \in \jormod $ such that $ \jorprojone(v_0) = 1_\ring $ and $ \psi(v_0, u) = 0_\ring = \psi(u, v_0) $ for all $ u \in \jormod $.
	\end{stenumerate}
	Further, the Jordan module $ \jormodtup $ is called \defemph*{abelian}\index{Jordan module!abelian} if $ (\jormod, \joradd) $ is an abelian group, it is called \defemph*{of type $ C $}\index{Jordan module!of type C@of type $ C $} if $ \psi $ is the zero map and it is called \defemph*{anisotropic}\index{Jordan module!anisotropic} if $ \jorprojone(u) = 0_\ring $ implies $ u=0 $ for all $ u \in \jormod $. If the Jordan module $ \jormodtup $ is fixed, we will refer to $ \jorprojone $ as the \defemph*{Jordan module projection}\index{Jordan module!projection}, to $ \jorTrone $ as the \defemph*{Jordan module trace}\index{Jordan module!trace} and to $ \psi $ as the \defemph*{Jordan module skew-hermitian form}\index{Jordan module!skew-hermitian form}. We will sometimes refer to $ \jormod $ as a Jordan module if the remaining objects in $ \jormodtup $ are clear from the context.
\end{definition}

\begin{remark}\label{BC:jordanmodule-cor}
	The requirements on the maps $ \jorprojone, \jorTrone, \psi $ which are given in the first paragraph of \cref{BC:jordanmodule-def} can be expressed by the following formulas:
	\begin{align*}
		\jorprojone\brackets[\big]{\jorsc(u, a)} &= \rinv{a} \jorprojone(u) a, \\
		\jorTrone(a+b) &= \jorTrone(a) \joradd \jorTrone(b), \\
		\jorsc\brackets[\big]{\jorTrone(a), b} &= \jorTrone\brackets{\rinv{b} a b}, \\
		\psi(u \joradd u', v \joradd v') &= \psi(u,v) + \psi(u,v') + \psi(u',v) + \psi(u',v'), \\
		\rinv{\psi(u,v)} &= -\psi(v,u), \\
		\psi\brackets[\big]{\jorsc(u,a), \jorsc(v,b)} &= \rinv{a} \psi(u,v) b
	\end{align*}
	for all $ u,v,u',v' \in \jormod $ and all $ a,b \in \ring $.
\end{remark}

\begin{note}
	We will show in \cref{chap:BC} that Jordan modules of type $ C $ are precisely those Jordan modules which coordinatise $ C_n $-graded groups for $ n \ge 3 $. General Jordan modules appear in the coordinatisation of $ BC_n $-graded groups.
\end{note}

\begin{note}
	The name \enquote{Jordan module} is motivated by the fact that every Jordan module of type $ C $ can be equipped with the structure of a quadratic unital Jordan algebra, except that these Jordan algebras are only \enquote{weakly quadratic}. We will make this observation more precise in \cref{subsec:jor-mod-alg}.
\end{note}

We now consider a sequence of three examples of Jordan modules, each being a generalisation of the previous one. Hence it suffices to verify that the last example satisfies the axioms of a Jordan module, which we do in \cref{BC:jordanmodule-pseudquad-axioms}. We will show in \cref{sec:jordan-class} that if $ 2_\ring $ is invertible, then every Jordan module is of the form in \cref{BC:jordanmodule-pseudquad-ex}. We do not know an example of a Jordan module which is not of this form.

\begin{example}\label{BC:jordanmodule-invset-ex}
	Let $ (\ring, \ringzero, \rinvmap) $ be an involutory set. Denote by $ (\jormod, \jorsc) $ the square-module $ \ringzero $ with its canonical square-scalar multiplication and by $ \map{\jorprojone}{\jormod}{\nucleus(\ring)}{}{} $ the canonical embedding of $ \jormod $ into $ \ring $. Further, we put
	\[ \map{\jorTrone \defl \ringTr_\rinvmap}{\ring}{\jormod}{r}{r+\rinv{r}} \]
	and we denote by $ \map{\psi}{\jormod \times \jormod}{\ring}{}{} $ the zero map. Then $ \jormodtup \defl (\jormod, \jorsc, \jorprojone, \jorTrone, \psi) $ is an abelian Jordan module over $ (\ring, \ringzero, \rinvmap) $ of type $ C $. Since $ \jorprojone $ is an injective homomorphism of groups, $ \jormodtup $ is anisotropic.
\end{example}

\begin{example}\label{BC:jordanmodule-module-ex}
	Let $ (\ring, \ringzero, \rinvmap) $ be an involutory set and let $ \module $ be a right $ \ring $-module. Put $ \jormod \defl \module \times \ringzero $, which is an additive abelian group, and
	\[ \map{\jorsc}{\jormod \times \ring}{\jormod}{\brackets[\big]{(v,h), r}}{(vr, \rinv{r}hr)}. \]
	Further, we denote by $ \map{\jorprojone}{\jormod}{\ring}{}{} $ the canonical projection on the second component and by $ \map{\psi}{\jormod \times \jormod}{\ring}{}{} $ the zero map. Finally, set
	\[ \map{\jorTrone}{\ring}{\jormod}{r}{(0_\module, r + \rinv{r})}. \]
	Then $ \jormodtup \defl (\jormod, \jorsc, \jorprojone, \jorTrone, \psi) $ is an abelian Jordan module over $ (\ring, \ringzero, \rinvmap) $ of type $ C $. The special case $ \module = \compactSet{0} $ is precisely \cref{BC:jordanmodule-invset-ex}. Further, $ \jormodtup $ is anisotropic if and only if $ \module = \compactSet{0} $.
\end{example}

\begin{example}\label{BC:jordanmodule-pseudquad-ex}
	Let $ (\ring, \ringzero, \rinvmap) $ be an involutory set and let $ (\module, q, f) $ be a standard pseudo-quadratic module over $ (\ring, \ringzero, \rinvmap) $. Denote by $ \jormod \defl \psgr(\module) \subs \module \times \ring $ the group from \cref{pseud:T-def} and by $ \map{\jorsc}{\jormod \times \ring}{\jormod}{}{} $ the square-scalar multiplication from \cref{pseud:T-squaremod}. Further, denote by $ \map{\jorprojone}{\jormod}{\nucleus(\ring)}{}{} $ the projection map from \cref{pseud:T-proj}. Finally, we put
	\begin{align*}
		\map{\jorTrone}{\ring}{\jormod}{r}{(0_\module, r + \rinv{r})} \midand \map{\psi}{\jormod \times \jormod}{\ring}{\brackets[\big]{(v,h), (v', h')}}{f(v,v')}.
	\end{align*}
	Then $ \jormodtup \defl (\jormod, \jorsc, \jorprojone, \jorTrone, \psi) $ is a Jordan module over $ (\ring, \ringzero, \rinvmap) $. A Jordan module which is isomorphic to one of this form is called \defemph*{of pseudo-quadratic type}\index{Jordan module!of pseudo-quadratic type}.
	
	Note that if $ q $ and $ f $ are trivial (so that the pseudo-quadratic module $ (\module, q, f) $ is essentially an $ \ring $-module with no additional structure), then it follows from \cref{pseud:T-q-triv} that this Jordan module is precisely the one from \cref{BC:jordanmodule-module-ex}. Note further that $ f $ is trivial if and only if $ \psgr(\module) $ is of type $ C $, and that the triviality of $ f $ implies the triviality of $ q $ if $ 2_\ring $ is invertible (\thmitemcref{pseud:2}{pseud:2:inv-triv}).
\end{example}

\begin{lemma}\label{BC:jordanmodule-pseudquad-axioms}
	Let $ (\ring, \ringzero, \rinvmap) $ be an involutory set and let $ (\module, q, f) $ be a standard pseudo-quadratic module over $ (\ring, \ringzero, \rinvmap) $. Then the tuple $ \jormodtup \defl (\jormod, \jorsc, \jorprojone, \jorTrone, \psi) $ from \cref{BC:jordanmodule-pseudquad-ex} is indeed a Jordan module.
\end{lemma}
\begin{proof}
	By \cref{pseud:T-squaremod}, $ (\jormod, \jorsc) $ is a multiplicative square-module. The map $ \jorprojone $ is compatible with the square-scalar multiplication by \cref{pseud:T-proj}. By \cref{pseud:T-Tr}, the map $ \jorTrone $ is a homomorphism of square-modules. Since $ f $ is a skew-hermitian form, the same holds for $ \psi $. Axiom~\thmitemref{BC:jordanmodule-def}{BC:jordanmodule-def:tr} holds by \cref{rinv:Tr-assoc}. Axiom~\thmitemref{BC:jordanmodule-def}{BC:jordanmodule-def:proj-rinv} holds because $ (\jormod, q, f) $ is standard. By the definition of $ \psi $, $ \jorTrone $ and $ \jorprojone $, Axioms~\thmitemref{BC:jordanmodule-def}{BC:jordanmodule-def:psi} and~\thmitemref{BC:jordanmodule-def}{BC:jordanmodule-def:square} are satisfied. Equation~\eqref{eq:jormod-lin:proj} holds by the definition of the multiplication on $ \psgr(\module) $. Equation~\eqref{eq:jormod-lin:jorsc} holds by \cref{pseud:T-squaremod-standard}. Equation~\eqref{eq:jormod-lin:add} holds by \cref{pseud:T-comm}. In Axiom~\thmitemref{BC:jordanmodule-def}{BC:jordanmodule-def:inv}, we can take $ v_0 \defl (0, 1_\ring) $. This finishes the proof.
\end{proof}

\begin{lemma}
	Let $ (\ring, \ringzero, \rinvmap) $ be an involutory set, let $ (\module, q, f) $ be a standard pseudo-quadratic module over $ (\ring, \ringzero, \rinvmap) $ and let $ \jormodtup \defl (\jormod, \jorsc, \jorprojone, \jorTrone, \psi) $ be the Jordan module from \cref{BC:jordanmodule-pseudquad-ex}. Then $ (\module, q, f) $ is anisotropic if and only if $ \jormodtup $ is anisotropic.
\end{lemma}
\begin{proof}
	At first, assume that $ (\module, q, f) $ is anisotropic. Let $ (u,h) \in \jormod $ such that $ \jorprojone(u,h) = 0 $, that is, such that $ h = 0_\ring $. Then $ q(u) \equiv h = 0 $ modulo $ \ringzero $. Since $ q $ is anisotropic, it follows that $ u= 0 $. Thus $ (u,h) = 0_\jormod $, so $ \jormodtup $ is anisotropic.
	
	Now assume that $ \jormodtup $ is anisotropic and let $ u \in \module $ such that $ q(u) \equiv 0 $ modulo $ \ringzero $. Then $ (u,0) $ lies in $ \jormod $ and we have $ \jorprojone(u, 0) = 0_\ring $. Since $ \jormodtup $ is anisotropic, it follows that $ (u,0) = 0_\jormod $, so $ u = 0_\module $. Hence $ (\module, q, f) $ is anisotropic.
\end{proof}

\subsection{Basic Properties}

\begin{secnotation}\label{secnot:jormod}
	From now on, we denote by $ \ring $ an alternative ring with a nuclear involution $ \rinvmap $ and by $ \calJ = (\jormod, \jorsc, \jorprojone, \jorTrone, \psi) $ a Jordan module over $ (\ring, \rinvmap) $.
\end{secnotation}

\begin{definition}[Homomorphism of Jordan modules]\label{BC:jormod:hom-def}
	Let
	\[ \jormodtup' = (\jormod', \jorsc', \jorprojone', \jorTroneprime, \psi') \]
	be another Jordan module over $ (\ring, \rinvmap) $. A \defemph*{homomorphism from $ \jormodtup $ to $ \jormodtup' $}\index{Jordan module!homomorphism} is a homomorphism $ \map{f}{(\jormod, \jorsc)}{(\jormod', \jorsc')}{}{} $ of square-modules such that $ \jorprojone = \jorprojone' \circ f $, $ \jorTroneprime = f \circ \jorTrone $ and $ \psi = \psi' \circ (f \times f) $. An \defemph*{isomorphism of Jordan modules} is a bijective homomorphism of Jordan modules and an \defemph*{automorphism of Jordan modules}\index{Jordan module!automorphism} is an isomorphism from a Jordan module to itself.
\end{definition}

\begin{definition}[Jordan submodule]
	A \defemph*{Jordan submodule of $ \jormodtup $}\index{Jordan module!submodule} is a square-submo\-dule $ U $ of $ (\jormod, \jorsc) $ which contains the image of $ \jorTrone $ and which satisfies Axiom~\thmitemref{BC:jordanmodule-def}{BC:jordanmodule-def:inv}.
\end{definition}

Clearly, every Jordan submodule of $ \jormodtup $ is a Jordan submodule itself with respect to the restrictions of the maps $ \jorsc $, $ \jorprojone $, $ \jorTrone $ and $ \psi $.

Any Jordan module has a canonical involution. It will appear as one of the twisting actions in the standard parameter system of type $ BC $ (\cref{BC:jordanmod:standard-paramsys}).

\begin{definition}[Jordan module involution]
	The canonical involution on the multiplicative square-module $ (\jormod, \jorsc) $ in the sense of \cref{weakquad:inv-def} is called the \defemph*{Jordan module involution (on $ \jormodtup $)}\index{Jordan module!involution}. It is the map
	\[ \map{\modinvmap}{\jormod}{\jormod}{u}{\modinv{u} \defl \jorsc(u, -1)}. \]
\end{definition}

\begin{lemma}\label{BC:jordanmod:inv-auto}
	The Jordan module involution is an automorphism of $ \jormod $ which fixes the image of $ \jorTrone $ element-wise.
\end{lemma}
\begin{proof}
	Let $ u,v \in \jormod $ and $ r \in \ring $. Then we have
	\begin{align*}
		\jorprojone(\modinv{u}) &= \jorprojone\brackets[\big]{\jorsc(u, -1_\ring)} = \rinv{(-1_\ring)} \jorprojone(u) (-1_\ring) = \jorprojone(u), \\
		\psi(\modinv{u}, \modinv{v}) &= \psi\brackets[\big]{\jorsc(u,-1_\ring), \jorsc(v, -1_\ring)} = \rinv{(-1_\ring)} \psi(u,v) (-1_\ring) = \psi(u,v), \\
		\modinv{\jorTrone(r)} &= \jorsc\brackets[\big]{\jorTrone(r), -1_\ring} = \jorTrone\brackets[\big]{\rinv{(-1_\ring)} r (-1_\ring)} = \jorTrone(r).
	\end{align*}
	Since $ \modinvmap $ is a homomorphism of square-modules by \cref{weakquad:inv-lem}, the assertion follows.
\end{proof}

\begin{remark}\label{BC:jordanmod:comm-square}
	For all $ u \in \jormod $, we have
	\begin{align*}
		0_\jormod &= \jorsc(u, 0_\ring) = \jorsc\brackets[\big]{u, 1_\ring + (-1_\ring)} \\
		&= \jorsc(u, 1_\ring) \joradd \jorsc(u, -1_\ring) \joradd \jorTrone\brackets[\big]{\rinv{(-1)} \jorprojone(u) 1_\ring} \\
		&= u \joradd \modinv{u} \jormin \jorTrone\brackets[\big]{\jorprojone(u)}.
	\end{align*}
	This shows that $ \jorTrone(\jorprojone(u)) = u \joradd \modinv{u} $. Together with Axiom~\thmitemref{BC:jordanmodule-def}{BC:jordanmodule-def:square}, we conclude that the diagram in \cref{fig:BC:jordanmodule-def:square} commutes.
\end{remark}

\premidfigure
\begin{figure}[htb]
	\centering$ \begin{tikzcd}[column sep=large]
		\ring \arrow[r, "a \mapsto a+\rinv{a}"] \arrow[d, "\jorTrone"'] & \ring \arrow[d, "\jorTrone"] \\
		\jormod \arrow[ur, "\jorprojone"] \arrow[r, "u \mapsto u + \modinv{u}"'] & \jormod.
	\end{tikzcd} $
	\caption{The commutative square in \cref{BC:jordanmod:comm-square}.}
	\label{fig:BC:jordanmodule-def:square}
\end{figure}
\postmidfigure

The properties of being abelian and of being of type $ C $ are closely related, but not equivalent.

\begin{lemma}\label{BC:jordanmod:prop}
	The following statements hold:
	\begin{lemenumerate}
		\item \label{BC:jordanmod:prop:abel-char}$ \calJ $ is abelian if and only if $ \mathord{\jorTrone} \circ \psi = 0 $.
		
		\item \label{BC:jordanmod:prop:C-char}$ \calJ $ is of type $ C $ if and only if $ \map{\jorprojone}{\jormod}{\ring}{}{} $ is a homomorphism of square-modules.
		
		\item \label{BC:jordanmod:prop:C-abel}If $ \calJ $ is of type $ C $, then it is abelian.
		
		\item \label{BC:jordanmod:prop:abel-sym}If $ \calJ $ is abelian, then $ \psi $ is symmetric.
		
		\item \label{BC:jordanmod:prop:psi-char}$ \psi $ is symmetric if and only if its image lies in the set $ \Set{r \in \ring \given \rinv{r} = -r} $ of skew-symmetric elements.
	\end{lemenumerate}
\end{lemma}
\begin{proof}
	Assertions~\itemref{BC:jordanmod:prop:abel-char} and~\itemref{BC:jordanmod:prop:C-char} follow from Axiom~\thmitemref{BC:jordanmodule-def}{BC:jordanmodule-def:lin}. Further,~\itemref{BC:jordanmod:prop:abel-char} implies~\itemref{BC:jordanmod:prop:C-abel}. Now assume that $ \calJ $ is abelian. Then for all $ u,v \in \jormod $, it follows from Axiom~\thmitemref{BC:jordanmodule-def}{BC:jordanmodule-def:lin} that
	\begin{align*}
		\jorprojone(v) + \jorprojone(u) + \psi(v,u) = \jorprojone(v \joradd u) = \jorprojone(u \joradd v) = \jorprojone(u) + \jorprojone(v) + \psi(u,v).
	\end{align*}
	Since $ (\ring, +) $ is abelian, we infer that $ \psi $ is symmetric. This proves~\itemref{BC:jordanmod:prop:abel-sym}.
	
	Now let $ u,v \in \jormod $. Since $ \psi $ is skew-hermitian, we have $ \rinv{\psi(u,v)} = -\psi(v,u) $. Thus the following assertions are equivalent:
	\begin{align*}
		\psi(u,v) = \psi(v,u) \Longleftrightarrow \rinv{\psi(u,v)} = \rinv{\psi(v,u)} \Longleftrightarrow \rinv{\psi(u,v)} = -\psi(u,v).
	\end{align*}
	It follows that~\itemref{BC:jordanmod:prop:psi-char} holds.
\end{proof}

\begin{remark}[The radical]\label{BC:jordanmod:radical}
	Denote by
	\[ R \defl \Rad(\psi) = \Set{u \in R \given \psi(x,u) = 0 = \psi(u,x) \text{ for all } x \in \jormod} \]
	the radical\index{radical} of the skew-hermitian form $ \psi $. This is a square-submodule of $ \jormod $. Further, Axiom~\thmitemref{BC:jordanmodule-def}{BC:jordanmodule-def:psi} says precisely that the image of $ \jorTrone $ is contained in $ R $. Similarly, Axiom~\thmitemref{BC:jordanmodule-def}{BC:jordanmodule-def:inv} says that there exists $ v_0 \in \Rad(\psi) $ such that $ \jorprojone(v_0) = 1_\ring $. Thus
	\[ \calJ' \defl \brackets[\big]{R, \restrict{\jorsc}{R \times \ring}, \restrict{\jorprojone}{R}, \jorTrone, \restrict{\psi}{R \times R}} \]
	is also a Jordan module over $ (\ring, \rinvmap) $. In fact, $ \restrict{\psi}{R \times R} $ is the zero map, so $ \calJ' $ is of type $ C $.
\end{remark}

\begin{remark}\label{BC:jordanmodule:2-inv-Trinv}
	Assume that $ 2_\ring $ is invertible and put $ v_0' \defl \jorTrone(2_\ring^{-1}) \in \jormod $. Then $ \jorprojone(v_0') = 2_\ring^{-1} + \rinv{(2_\ring^{-1})} = 1_\ring $ by Axiom~\thmitemref{BC:jordanmodule-def}{BC:jordanmodule-def:square} and $ v_0' \in \Rad(\psi) $ by Axiom~\thmitemref{BC:jordanmodule-def}{BC:jordanmodule-def:psi}. We conclude that Axiom~\thmitemref{BC:jordanmodule-def}{BC:jordanmodule-def:inv} is redundant if $ 2_\ring $ is invertible.
\end{remark}

\begin{remark}\label{BC:jordanmod:Tr-jorproj-rem}
	We can define the following maps:
	\begin{align*}
		\map{\jorproj}{\jormod \times \ring}{\ring}{(u,a)}{a \jorprojone(u)} \midand \map{\jorTr}{\ring \times \ring}{\jormod}{(a,b)}{\jorTrone(a \rinv{b})}.
	\end{align*}
	It is straightforward to verify that the following properties, which we will need later, follow from the axioms of a Jordan module:
	\begin{align*}
		\modinv{\jorsc(u,a)} &= \jorsc\brackets[\big]{\jorsc(u,a), -1} = \jorsc(u,-a) = \jorsc\brackets[\big]{\jorsc(u,-1),a} = \jorsc(\modinv{u}, a), \\
		\modinv{\jorTrone(a)} &= \jorsc\brackets[\big]{\jorTrone(a), -1} = \jorTrone\brackets[\big]{\rinv{(-1)}a (-1)} = \jorTrone(a), \\
		\psi(\modinv{u}, v) &= \rinv{(-1)} \psi(u,v) = \psi(u,v)(-1) = \psi(u, \modinv{v}), \\
		\jorproj(\modinv{u}, a) &= a \jorprojone\brackets[\big]{\jorsc(u, -1)} = a \rinv{(-1)} \jorprojone(u) (-1) = \jorproj(u,a), \\
		\jorproj\brackets[\big]{\jorTr(a,b),c} &= c \jorprojone\brackets[\big]{\jorTrone(a \rinv{b})} = c(a \rinv{b}) + c(\rinv{b}a) \\
		&= (ca) \rinv{b} + (c\rinv{b})a - \assoc{c}{a}{\rinv{b}} - \assoc{c}{\rinv{b}}{a} = (ca) \rinv{b} + (c\rinv{b})a, \\
		\jorproj\brackets[\big]{\jorsc(u,a), b} &= b \jorprojone\brackets[\big]{\jorsc(u,a)} = b \rinv{a} \jorprojone(u) a, \\
		\jorTr(a,b) &= \jorTrone(a \rinv{b}) = \jorTrone\brackets[\big]{\rinv{(a \rinv{b})}} = \jorTrone(b \rinv{a}) = \jorTr(b, a), \\
		\jorTr(ab,c) &= \jorTrone(ab \rinv{c}) = \jorTrone(c\rinv{b} \rinv{a}) = \jorTr(c \rinv{b}, a) = \jorTr(a, c \rinv{b}), \\
		\jorsc\brackets[\big]{\jorTr(a,b), c} &= \jorsc\brackets[\big]{\jorTrone(a \rinv{b}), c} = \jorTrone(\rinv{c} a \rinv{b} c) = \jorTr(\rinv{c}a, \rinv{c}b) = \jorTr(\rinv{c}b, \rinv{c}a)
	\end{align*}
	for all $ a,b,c \in \ring $ and all $ u \in \jormod $. (Here we have used \cref{rinv:artin-inv} to write the term $ b \rinv{a} \jorprojone(u) a $ without brackets.)
	
	The maps $ \jorproj $ and $ \jorTr $ are the ones which describe the commutator relations in root gradings of type $ BC $. Thus in some sense, they are more fundamental for our group-theoretic interests than the maps $ \jorprojone $ and $ \jorTrone $. However, the latter allow for a more concise treatment of Jordan modules in the algebraic setting. In particular, the fact that $ \jorprojone $ and $ \jorTrone $ are homomorphisms of square-modules (except for the minor nuisance that $ \jorprojone $ is not necessarily additive) and that we have the commutative square in \cref{BC:jordanmod:comm-square} cannot be expressed as succinctly in terms of the maps $ \jorproj $ and $ \jorTr $.
\end{remark}

Jordan modules arise from root gradings of type $ (B)C $ in the form of the following parameter system.

\begin{definition}[Standard parameter system]\label{BC:jordanmod:standard-paramsys}
	Define $ \twistgroup \defl \compactSet{\pm 1}^2 $ and $ \invogroup \defl \compactSet{\pm 1} $. Denote the elements of $ \twistgroup $ by $ (\pm 1_\twistgroup, \pm 1_\twistgroup) $ and the elements of $ \invogroup $ by $ \pm 1_\invogroup $. We declare that the first component of $ \twistgroup $ acts on $ (\ring, +) $ and $ (\jormod, \joradd) $ by inversion and that the second component of $ \twistgroup $ acts trivially on $ \ring $ and by the Jordan module involution on $ \jormod $. Further, we declare that $ \invogroup $ acts by $ \rinvmap $ on $ \ring $ and trivially on $ \psgr $. More precisely, these declarations mean that
	\begin{gather*}
		(-1_\twistgroup, 1_\twistgroup).r = -r, \quad (-1_\twistgroup, 1_\twistgroup).u = \jormin u, \quad (1_\twistgroup, -1_\twistgroup).r = r, \quad (1_\twistgroup, -1_\twistgroup).u = \modinv{u}, \\
		-1_\invogroup.r = \rinv{r}, \qquad -1_\invogroup.u = u
	\end{gather*}
	for all $ r \in \ring $ and all $ u \in \jormod $. Then the triple $ (\twistgroup \times \invogroup, \jormod, \ring) $ is called the \defemph*{standard parameter system for $ \jormodtup $}.\index{parameter system!standard!type BCn@type $ BC_n $}
\end{definition}

\begin{remark}
	The action of $ \twistgroup \times \invogroup $ on $ \psgr $ and $ \ring $ in \cref{BC:jordanmod:standard-paramsys} is induced by the respective actions of $ \twistgroup $ and $ \invogroup $. Such an induced action exists because the actions of $ \twistgroup $ and $ \invogroup $ commute.
\end{remark}

\begin{reminder}
	If the Jordan module $ \jormodtup $ is of the form $ \psgr(\module, q, f) $ for a standard pseudo-quadratic module $ (\module, q, f) $ (as in \cref{BC:jordanmodule-pseudquad-ex}), then for all $ (u,h) \in \psgr(\module, q, f) $, the inverse $ \jormin (u, h) $ equals $ (-u, -\rinv{h}) $ (see \cref{pseud:T-group,pseud:T-group:standard}). In this situation, the Jordan module involution is given by $ \modinv{(u,h)} = (-u, h) $ (see \cref{pseud:T-invo}).
\end{reminder}

\subsection{Jordan Modules and Jordan Algebras}

\label{subsec:jor-mod-alg}

We now investigate the relation of Jordan modules to Jordan algebras. \Cref{secnot:jormod} continues to hold.

\begin{definition}[Weak Jordan algebra]\label{jordanalg-def}
	A \defemph*{weakly quadratic Jordan algebra}\index{Jordan algebra!weakly quadratic} is a tuple $ (\jormod, U, 1_\jormod) $ consisting of an abelian group $ \jormod = (\jormod, +) $, an element $ 1_\jormod \in \jormod $ and a weakly quadratic map $ \map{U}{\jormod}{\Hom(\jormod, \jormod)}{x}{U_x} $ such that the following properties are satisfied for all $ x,y,z,v \in \jormod $, where
	\[ \map{\tribrace{}{}{}}{\jormod \times \jormod \times \jormod}{\jormod}{(x,y,z)}{ U_{x+z}y - U_x y - U_z y} \]
	\begingroup\mywarning{page break equation}
	\allowdisplaybreaks
	denotes the tri-linearisation of $ U $:
	\begin{align*}
		U_{U_x y} v &= U_x U_y U_x v, \\
		U_x(\tribrace{y}{x}{z}) &= \tribrace{U_x y}{z}{x}, \\
		U_1 v &= v, \\
		\tribrace{U_x y}{v}{\tribrace{x}{y}{z}} &= U_x U_y(\tribrace{x}{v}{z}) + \tribrace{x}{U_y U_x v}{z}, \\
		\tribrace{U_x y}{v}{U_z y} + U_{\tribrace{x}{y}{z}}v &= U_x U_y U_z v + U_z U_y U_x v + \tribrace{x}{U_y(\tribrace{x}{v}{z})}{z}, \\
		U_x(\tribrace{y}{v}{z}) + \tribrace{x}{\tribrace{y}{x}{z}}{v} &= \tribrace{U_x y}{z}{v} + \tribrace{\tribrace{x}{y}{v}}{z}{x}.
	\end{align*}
	\endgroup
	Further, $ (\jormod, U, 1_\jormod) $ is called a \defemph{Jordan algebra} if, in addition, the map $ U $ is (properly) $ \IZ $-quadratic in the sense of \cref{quadmod:quadmap-def}.
\end{definition}

\begin{note}
	What we call \enquote{Jordan algebra} could, in more detail, be called \enquote{unital quadratic Jordan algebra over $ \IZ $}. Our definition agrees with the one in \cite[1.3.4]{Jacobson_JordanStrucTheory}, but it is phrased in a slightly different way. Observe that the last three identities in \cref{jordanalg-def} are precisely the polarisations of the first two identities. (See, for example, \cite[p.~1072]{McCrimmon66}.) It follows from this fact that $ (\jormod, U, 1_\jormod) $ is a Jordan algebra if and only if for all commutative associative rings $ \comring $, the scalar extension $ (\jormod_\comring, U_\comring, (1_\jormod)_\comring) $ satisfies the first three identities in \cref{jordanalg-def}. This is precisely the definition that is given in \cite{Jacobson_JordanStrucTheory}.
	
	The observation in the previous paragraph relies on the fact that a quadratic map $ U $ has a unique and well-defined extension $ U_\comring $ to any scalar extension. (For a proof, see for example \cite[11.5]{GPR_AlbertRing}.) This is not clear for weakly quadratic maps. Hence there is no way to phrase the definition of weakly quadratic Jordan algebras in the language of scalar extensions. In particular, scalar extension of weakly quadratic Jordan algebras do not necessarily have the structure of a weakly quadratic Jordan algebra, and need not even exist.
\end{note}

\begin{remark}[Jordan algebras from Jordan modules]
	Assume that $ \jormodtup $ is of type $ C $ (that is, that $ \psi = 0 $). We define a map $ \map{U}{\jormod}{\Hom_\IZ(\jormod,\jormod)}{x}{U_x} $ by
	\[ \map{U_x}{\jormod}{\jormod}{v}{\phi\brackets[\big]{v, \pi(x)}} \]
	for all $ x \in \jormod $. Further, we define a tri-additive map
	\begin{align*}
		\map{\tribrace{}{}{}}{\jormod \times \jormod \times \jormod&}{\jormod}{\\(x,y,z)&}{U_{x+z}y - U_x y - U_z y = \jorTrone\brackets[\big]{\rinv{\pi(x)} \pi(y) \pi(z)}}.
	\end{align*}
	Since the map $ \map{}{\jormod \times \jormod}{\Hom(\jormod,\jormod)}{(x,z)}{\tribrace{x}{}{z}} $ is the polarisation of $ U $, it follows that $ U $ is weakly quadratic. A straightforward computation shows that $ (\jormod, U, v_0) $ is a weakly quadratic Jordan algebra where $ v_0 $ is any element as in Axiom~\thmitemref{BC:jordanmodule-def}{BC:jordanmodule-def:inv}.
\end{remark}


\section{A Near-classification of Jordan Modules}

\label{sec:jordan-class}

\begin{secnotation}\label{secnot:jordan-class}
	Unless otherwise stated, we denote by $ \ring $ an alternative ring with a nuclear involution $ \rinvmap $ and by $ \jormodtup = (\jormod, \jorsc, \jorprojone, \jorTrone, \psi) $ a Jordan module over $ (\ring, \rinvmap) $.
\end{secnotation}
 
In this section, we show that under some weak assumptions, every Jordan module $ \jormodtup $ is of pseudo-quadratic type in the sense of \cref{BC:jordanmodule-pseudquad-ex}. These assumptions will always be satisfied if $ 2_\ring $ is invertible.

Note that, a priori, the Jordan module $ \jormodtup $ is defined over the ring $ \ring $ with a nuclear involution $ \rinvmap $ whereas a pseudo-quadratic module $ \module $ is defined over an involutory set $ (\ring, \ringzero, \rinvmap) $. Thus as a first step, we should locate an appropriate subset $ \ringzero $ of $ \ring $. We know from \cref{rinv:inv-set-2inv} that this set is uniquely determined if $ 2_\ring $ is invertible, but in general it is not. It turns out that our choice of $ \ringzero $ will depend on the choice of a so-called involutory submodule of $ \jormod $.

\begin{definition}[Involutory submodule]\label{BC:jordanmod:inv-submodule-def}
	An \defemph*{involutory submodule of $ \jormodtup $}\index{involutory submodule} is a square-submodule $ \invsub $ of $ \jormod $ which satisfies
	\[ \jorTrone(\ring) \subs \invsub \subs \Rad(\psi) \]
	where $ \Rad(\psi) $ denotes the radical of $ \psi $ (see \cref{BC:jordanmod:radical}). It is called \defemph*{unital}\index{involutory submodule!unital} if there exists $ d \in \invsub $ such that $ \jorprojone(d) = 1_\ring $. It is called an \defemph*{embedding involutory submodule}\index{involutory submodule!embedding} if the restriction $ \map{\restrict{\jorprojone}{\invsub}}{\invsub}{\ring}{}{} $ is injective.
\end{definition}

\begin{remark}\label{BC:jordanmodule:invsub-abel}
	Let $ \invsub $ be an involutory submodule of $ \jormodtup $. Since $ \invsub $ is contained in $ \Rad(\psi) $, it follows from \cref{BC:jordanmod:radical,BC:jordanmod:prop} that $ \psi $ is zero on $ \invsub \times \invsub $, that $ \invsub $ is abelian and that the restriction of $ \jorprojone $ to $ \invsub $ is a homomorphism of square-modules.
\end{remark}

The following result justifies the terminology of \enquote{involutory submodules}.

\begin{lemma}\label{BC:jordanmodule:invsub-invset}
	Let $ \invsub $ be an involutory submodule of $ \jormodtup $. Then the triple $ (\ring, \jorprojone(\invsub), \rinvmap) $ is a pre-involutory set. It is an involutory set if and only if $ \invsub $ is unital.
\end{lemma}
\begin{proof}
	It follows from \cref{jormod:ker-im,BC:jordanmodule:invsub-abel} that $ \jorprojone(\invsub) $ is a square-submodule of $ \nucleus(\ring) $. Further, since $ \restrict{\psi}{\invsub \times \invsub} $ is the zero map, it follows from Axiom~\thmitemref{BC:jordanmodule-def}{BC:jordanmodule-def:proj-rinv} that $ \jorprojone(\invsub) $ is contained in $ \symring $. Finally, as $ \invsub $ contains $ \jorTrone(\ring) $, Axiom~\thmitemref{BC:jordanmodule-def}{BC:jordanmodule-def:square} implies that $ \jorprojone(\invsub) $ contains $ \ringTr(\ring) $. We conclude that $ (\ring, \jorprojone(\invsub), \rinvmap) $ is a pre-involutory set. The second assertion is trivial.
\end{proof}

\begin{example}[of an involutory submodule]\label{BC:jordanmodule:invsub-ex}
	Assume that $ \jormod $ is the Jordan module $ \psgr(\module) $ from \cref{BC:jordanmodule-pseudquad-ex} for some pseudo-quadratic module $ (\module, q, f) $ over a pre-involutory set $ (\ring, \ringzero, \rinvmap) $. Then $ \invsub \defl \compactSet{0_\jormod} \times \ringzero \subs \psgr(\module) $ is an embedding involutory submodule of $ \psgr(\module) $ (which is unital if and only if $ (\ring, \ringzero, \rinvmap) $ is an involutory set).
\end{example}

We have seen in \cref{BC:jordanmodule:invsub-ex} that the existence of a embedding involutory submodule of $ \jormod $ is a necessary condition for $ \jormod $ to be of pseudo-quadratic type. The main result of this section is precisely that this is also a sufficient condition.

Before we turn to the proof of our main result, we investigate the existence of embedding involutory submodules. It is clear that both $ \jorTrone(\ring) $ and $ \Rad(\psi) $ are involutory submodules of $ \jormodtup $. Further, since the element $ v_0 $ from Axiom~\thmitemref{BC:jordanmodule-def}{BC:jordanmodule-def:inv} lies in $ \Rad(\psi) $, the square-submodule generated by $ \jorTrone(\ring) \union \compactSet{v_0} $ is a unital involutory submodule of $ \jormodtup $. Thus the main difficulty lies finding an involutory submodule which is embedding. Under suitable assumptions on $ 2 $-torsion and related phenomena, we will show that $ \jorTrone(\ring) $ satisfies this condition. In the general situation, however, we have no result which guarantees the existence of an embedding involutory submodule in $ \jormod $.

Recall from \cref{BC:jordanmodule:invsub-abel} that the restriction of $ \jorprojone $ to $ \jorTrone(\ring) $ is a homomorphism of groups. Hence to show that $ \jorTrone(\ring) $ is embedding, it suffices to verify that the preimage of $ 0_\ring $ in $ \jorTrone(\ring) $ under $ \jorprojone $ is trivial.

\begin{lemma}
	Assume that there exists no 2-torsion on $ \jorTrone(\ring) $. (That is, if $ d \in \jorTrone(\ring) $ satisfies $ d \joradd d = 0_\jormod $, then $ d = 0_\jormod $.) Then $ \jorTrone(\ring) $ is a embedding involutory submodule of~$ \jormodtup $.
\end{lemma}
\begin{proof}
	Let $ d \in \jorTrone(\ring) $ such that $ \jorprojone(d) = 0_\ring $. We want to show that $ d=0_\jormod $. Applying \cref{BC:jordanmod:comm-square}, we see that
	\[ 0_\jormod = \jorTrone(0_\ring) = \jorTrone\brackets[\big]{\jorprojone(d)} = d \joradd \modinv{d}. \]
	Since $ d $ lies in $ \jorTrone(\ring) $, we have $ \modinv{d} = d $ by \cref{BC:jordanmod:inv-auto}. Thus $ d \joradd d = 0_\jormod $, which implies that $ d=0_\jormod $. We infer that the restriction of $ \jorprojone $ to $ \jorTrone(\ring) $ is injective.
\end{proof}

\begin{lemma}\label{BC:jordanmodule:invsub-crit-2inv}
	Assume that $ 2_\ring $ is invertible. Then $ \jorTrone(\ring) $ is a unital embedding involutory submodule of $ \jormodtup $.
\end{lemma}
\begin{proof}
	By \cref{BC:jordanmodule:2-inv-Trinv}, we have $ \jorprojone(v_0') = 1_\ring $ for $ v_0' \defl \jorTrone(2_\ring^{-1}) $, so $ \jorTrone(\ring) $ is unital. Now let $ d \in \jorTrone(\ring) $ such that $ \jorprojone(d) = 0_\ring $. We want to show that $ d = 0_\jormod $. Choose $ r \in \ring $ such that $ d = \jorTrone(r) $. It follows from Axiom~\thmitemref{BC:jordanmodule-def}{BC:jordanmodule-def:square} that
	\[ r + \rinv{r} = \jorprojone\brackets[\big]{\jorTrone(r)} = \jorprojone(d) = 0_\ring. \]
	Hence $ d' \defl \jorTrone(2_\ring^{-1} r) $ also satisfies
	\[ \jorprojone(d') = \jorprojone\brackets[\big]{\jorTrone(2_\ring^{-1} r)} = \frac{r}{2} + \rinv{\brackets*{\frac{r}{2}}} = \frac{1}{2} (r + \rinv{r}) = 0_\ring. \]
	Again by Axiom~\thmitemref{BC:jordanmodule-def}{BC:jordanmodule-def:square} and \cref{BC:jordanmod:inv-auto}, we infer that
	\begin{align*}
		0_\jormod &= \jorTrone(0_\ring) = \jorTrone\brackets[\big]{\jorprojone(d')} = d' \joradd \modinv{d'} = d' \joradd d' \\
		&= \jorTrone(2_\ring^{-1}r) \joradd \jorTrone(2_\ring^{-1}r) = \jorTrone(2_\ring^{-1}r + 2_\ring^{-1}r) = \jorTrone(r) = d.
	\end{align*}
	Thus the restriction of $ \jorprojone $ to $ \jorTrone(\ring) $ is injective, as desired.
\end{proof}

\begin{secnotation}
	From now on, we fix an involutory submodule $ \invsub $ of $ \jormodtup $. We put $ \quotmod \defl \jormod / \invsub $ and we denote by $ \map{\tilde{\jorsc}}{\quotmod \times \ring}{\quotmod}{}{} $ the multiplicative square-scalar multiplication on $ \quotmod $ from \cref{weakquad:quot-lem}. That is, we have
	\[ \quotepi\brackets[\big]{\jorsc(u,r)} = \tilde{\jorsc}\brackets[\big]{\quotepi(u), r} \]
	for all $ u \in \quotmod $, $ r \in \ring $ where $ \map{\quotepi}{\jormod}{\quotmod}{}{} $ denotes the canonical homomorphism of square-modules.
\end{secnotation}

For the first few observations, we do not have to assume that $ \invsub $ is embedding. This will only be necessary in \cref{BC:jordanmod:class-inj}.

We begin by equipping $ \quotmod $ with the structure of a pseudo-quadratic module. First of all, we have to verify that $ (\quotmod, \tilde{\jorsc}) $ is an $ \ring $-module and not merely a square-module.

\begin{lemma}
	The group $ (\quotmod, \joradd) $ (whose group structure is induced by the one on $ \jormod $) is abelian.
\end{lemma}
\begin{proof}
	Let $ u,v \in \jormod $. Then
	\[ u \joradd v = v \joradd u \joradd \jorTrone\brackets[\big]{\psi(u,v)} \]
	by Axiom~\thmitemref{BC:jordanmodule-def}{BC:jordanmodule-def:lin}. Applying $ \quotepi $, and using that $ \invsub $ contains the image of $ \jorTrone $, we infer that $ \quotepi(u) \joradd \quotepi(v) = \quotepi(v) \joradd \quotepi(u) $. Since $ \quotepi $ is surjective, it follows that $ \quotmod $ is abelian.
\end{proof}

\begin{lemma}
	$ (\quotmod, \tilde{\jorsc}) $ is an $ \ring $-module (in the regular sense of \cref{ring:module}).
\end{lemma}
\begin{proof}
	By construction, $ (\quotmod, \tilde{\jorsc}) $ is a square-module. Thus it only remains to show that $ \tilde{\jorsc} $ is additive in the second component. Let $ u \in \jormod $ and let $ r,s \in \ring $. We know from Axiom~\thmitemref{BC:jordanmodule-def}{BC:jordanmodule-def:lin} that
	\begin{align*}
		\tilde{\jorsc}(u, r+s) &= \tilde{\jorsc}(u,r) \joradd \tilde{\jorsc}(u,s) \joradd \jorTrone\brackets[\big]{\rinv{s} \jorprojone(u) r}.
	\end{align*}
	Applying $ \quotepi $ and using that $ \invsub $ contains the image of $ \jorTrone $, we infer that
	\[ \tilde{\jorsc}\brackets[\big]{\quotepi(u), r+s} = \tilde{\jorsc}\brackets[\big]{\quotepi(u), r} \joradd \tilde{\jorsc}\brackets[\big]{\quotepi(u), s}. \]
	Since $ \quotepi $ is surjective, the assertion follows.
\end{proof}

We now construct the pseudo-quadratic structure on $ \quotmod $.

\begin{lemma}
	There exists a unique map $ \map{f}{\quotmod \times \quotmod}{\nucleus(\ring)}{}{} $ such that
	\[ f\brackets[\big]{\quotepi(u), \quotepi(v)} = \psi(u,v) \]
	for all $ u,v \in \jormod $. This map is a skew-hermitian form with respect to $ \rinvmap $.
\end{lemma}
\begin{proof}
	For all $ u,v \in \module $ and all $ d,d' \in \invsub $, we have $ \psi(u+d, v+d') = \psi(u,v) $ because $ \invsub $ is contained in the radical of $ \psi $. Thus such a map $ f $ exists. Since $ \quotepi $ is surjective, $ f $ is uniquely determined. Further, it is skew-hermitian because $ \psi $ is.
\end{proof}

Recall from \cref{pseud:coset} that the equivalence class of a pseudo-quadratic form $ \map{q}{\quotmod}{\nucleus(\ring)}{}{} $ can be identified in a canonical way with the induced map $ \map{\tilde{q}}{\module}{\nucleus(\ring) / \ringzero}{}{} $. Since the group $ \psgr(N) $ depends only on the equivalence class of a pseudo-quadratic module $ N $, it thus makes sense to begin by constructing $ \tilde{q} $ and then lifting it (in an arbitrary way) to a map $ \map{q}{\quotmod}{\nucleus(\ring)}{}{} $.

\begin{lemma}
	There exists a unique map $ \map{\tilde{q}}{\quotmod}{\nucleus(\ring) / \jorprojone(\invsub)}{}{} $ such that $ \tilde{q}(\quotepi(u)) = \jorprojone(u) + \jorprojone(\invsub) $ for all $ u \in \jormod $.
\end{lemma}
\begin{proof}
	For any subgroup $ A $ of $ \jormod $, the map $ \map{\jorprojone}{\jormod}{\nucleus(\ring)}{}{} $ induces a map $ \map{}{\jormod/A}{\nucleus(\ring)/\jorprojone(A)}{}{} $. The assertion follows by taking $ A \defl \invsub $.
\end{proof}

\begin{notation}
	We choose a map $ \map{q}{\quotmod}{\nucleus(\ring)}{}{} $ such that $ q(u) \in \tilde{q}(u) $ for all $ u \in \quotmod $.
\end{notation}

\begin{lemma}\label{BC:jordanmodule:class-pseud}
	The tuple $ (\quotmod, q, f) $ is a standard pseudo-quadratic module with respect to the pre-involutory set $ (\ring, \jorprojone(\invsub), \rinvmap) $.
\end{lemma}
\begin{proof}
	For all $ u,v \in \jormod $, it follows from Axiom~\thmitemref{BC:jordanmodule-def}{BC:jordanmodule-def:lin} that
	\begin{align*}
		\tilde{q}\brackets[\big]{\quotepi(u) + \quotepi(v)} &= \tilde{q}\brackets[\big]{\quotepi(u+v)} = \jorprojone(u+v) + \jorprojone(\invsub) \\
		&= \jorprojone(u) + \jorprojone(v) + \psi(u,v) + \jorprojone(\invsub) \\
		&= \tilde{q}\brackets[\big]{\quotepi(u)} +\tilde{q}\brackets[\big]{\quotepi(v)} + f\brackets[\big]{\quotepi(u), \quotepi(v)} + \jorprojone(\invsub).
	\end{align*}
	Further, since $ \jorprojone $ is compatible with the square-scalar multiplication, we have
	\begin{align*}
		\tilde{q}\brackets[\big]{\tilde{\jorsc}\brackets[\big]{\quotepi(u), r}} &= \tilde{q}\brackets[\big]{\quotepi\brackets[\big]{\jorsc(u,r)}} = \jorprojone\brackets[\big]{\jorsc(u,r)} + \jorprojone(\invsub) \\
		&= \rinv{r} \jorprojone(u) r + \jorprojone(\invsub) \sups \rinv{r} \brackets[\big]{\jorprojone(u) + \jorprojone(\invsub)} r = \rinv{r} \tilde{q}\brackets[\big]{\quotepi(u)} r.
	\end{align*}
	for all $ u \in \jormod $ and $ r \in \ring $. It follows that
	\[ q(x+y) \equiv q(x) + q(y) + f(x,y) \midand q\brackets[\big]{\tilde{\jorsc}(x,r)} \equiv \rinv{r} q(x) r \]
	modulo $ \jorprojone(\invsub) $ for all $ x,y \in \quotmod $ and all $ r \in \ring $. We conclude that $ (\quotmod, q, f) $ is a pseudo-quadratic module.
	
	To see that $ (\quotmod, q, f) $ is standard, let $ \tilde{u} \in \quotmod $ be arbitrary and choose $ u \in \jormod $ such that $ \tilde{u} = \quotepi(u) $. Then it follows from Axiom~\thmitemref{BC:jordanmodule-def}{BC:jordanmodule-def:proj-rinv} and the definition of $ f $ that
	\begin{align*}
		\jorprojone(u) - \rinv{\jorprojone(u)} &= \jorprojone(u) - \jorprojone(u) + \psi(u,u) = f(\tilde{u}, \tilde{u}).
	\end{align*}
	Since $ q(\tilde{u}) $ lies in $ \tilde{q}(\tilde{u}) = \jorprojone(u) + \jorprojone(\invsub) $, there exists $ d \in \invsub $ such that $ q(\tilde{u}) = \jorprojone(u) + \jorprojone(d) $. Then
	\[ q(\tilde{u}) - \rinv{q(\tilde{u})} = \jorprojone(u) - \rinv{\jorprojone(u)} \]
	because $ \jorprojone(d) $ lies in $ \symring $ by \cref{BC:jordanmodule:invsub-invset}. It follows that $ q(\tilde{u}) - \rinv{q(\tilde{u})} = f(\tilde{u}, \tilde{u}) $, so we conclude that $ (\module, q, f) $ is standard.
\end{proof}

Having constructed a pseudo-quadratic structure on $ \quotmod $, we can now consider the group $ \psgr(\quotmod) = \psgr(\quotmod, q, f) $. It remains to find an isomorphism between this group and $ \jormod $.

\begin{definition}
	We define a map
	\[ \map{\alpha}{\jormod}{\psgr(\quotmod)}{u}{\brackets[\big]{\quotepi(u), \jorprojone(u)}}. \]
\end{definition}

By the definition of $ q $, the image of $ \alpha $ is indeed contained in $ \psgr(\quotmod) $, so $ \alpha $ is well-defined.

\begin{lemma}
	The map $ \alpha $ is a homomorphism of groups.
\end{lemma}
\begin{proof}
	Let $ u,v \in \jormod $ and put $ \tilde{u} \defl \quotepi(u) $, $ \tilde{v} \defl \quotepi(v) $. On the one hand,
	\begin{align*}
		\alpha(u \joradd v) &= \brackets[\big]{\quotepi(u \joradd v), \jorprojone(u \joradd v)} = \brackets[\big]{\tilde{u} \joradd \tilde{v}, \jorprojone(u) + \jorprojone(v) + \psi(u,v)}
	\end{align*}
	by Axiom~\thmitemref{BC:jordanmodule-def}{BC:jordanmodule-def:lin}. On the other hand,
	\begin{align*}
		\alpha(u) \alpha(v) &= \brackets[\big]{\tilde{u}, \jorprojone(u)} \brackets[\big]{\tilde{v}, \jorprojone(v)} = \brackets[\big]{\tilde{u} \joradd \tilde{v}, \jorprojone(u) + \jorprojone(v) + f(\tilde{u}, \tilde{v})}.
	\end{align*}
	We conclude that $ \alpha $ is a homomorphism.
\end{proof}

\begin{lemma}
	The map $ \alpha $ is a homomorphism of Jordan modules (with respect to the Jordan module structure on $ \psgr(\quotmod) $ from \cref{BC:jordanmodule-pseudquad-ex} and in the sense of \cref{BC:jormod:hom-def}).
\end{lemma}
\begin{proof}
	Let $ u,v \in \module $, let $ r \in \ring $ and put $ \tilde{u} \defl \quotepi(u) $, $ \tilde{v} \defl \quotepi(v) $. Denote the structure maps of the Jordan module $ \psgr(\quotmod) $ by $ \jorsc' $, $ \jorprojone' $, $ \jorTroneprime $ and $ \psi' $. Then we have
	\begin{align*}
		\alpha\brackets[\big]{\jorsc(u,r)} &= \brackets[\big]{\quotepi\brackets[\big]{\jorsc(u,r)}, \jorprojone\brackets[\big]{\jorsc(u,r)}} = \brackets[\big]{\tilde{\jorsc}(\tilde{u}, r), \rinv{r} \jorprojone(u) r} \\
		&= \phi'\brackets[\big]{\brackets[\big]{\tilde{u}, \jorprojone(u)}, r} = \jorsc'\brackets[\big]{\alpha(u), r}, \\
		\jorprojone'\brackets[\big]{\alpha(u)} &= \jorprojone'\brackets[\big]{\tilde{u}, \jorprojone(u)} = \jorprojone(u), \\
		\alpha\brackets[\big]{\jorTrone(r)} &= \brackets[\big]{\quotepi\brackets[\big]{\jorTrone(r)}, \jorprojone\brackets[\big]{\jorTrone(r)}} = (0_\quotmod, r+\rinv{r}) = \jorTroneprime(r), \\
		\psi'\brackets[\big]{\alpha(u), \alpha(v)} &= \psi'\brackets[\big]{\brackets[\big]{\tilde{u}, \jorprojone(u)}, \brackets[\big]{\tilde{v}, \jorprojone(v)}} = f(\tilde{u}, \tilde{v}) = \psi(u,v).
	\end{align*}
	This shows that $ \alpha $ is a homomorphism of Jordan modules.
\end{proof}

\begin{lemma}
	The map $ \alpha $ is surjective.
\end{lemma}
\begin{proof}
	Let $ (\tilde{u}, h) \in \psgr(\quotmod) $ be arbitrary and choose $ u \in \jormod $ such that $ \tilde{u} = \quotepi(u) $. By the definition of $ \psgr(\quotmod) $, we have $ h \equiv q(\tilde{u}) $ modulo $ \jorprojone(\invsub) $, which says precisely that $ h \in \tilde{q}(\tilde{u}) = \jorprojone(u) + \jorprojone(\invsub) $. Thus there exists $ d \in \invsub $ such that $ h = \jorprojone(u) + \jorprojone(d) $. Now
	\begin{align*}
		\alpha(u \joradd d) &= \brackets[\big]{\tilde{u} + \quotepi(d), \jorprojone(u + d)} = \brackets[\big]{\tilde{u}, \jorprojone(u) + \jorprojone(d) + \psi(u,d)} = (\tilde{u}, h).
	\end{align*}
	This shows that $ \alpha $ is surjective.
\end{proof}

\begin{lemma}\label{BC:jordanmod:class-inj}
	The kernel of $ \alpha $ is the kernel of $ \restrict{\jorprojone}{\invsub} $. In particular, $ \alpha $ is injective if and only if $ \invsub $ is embedding.
\end{lemma}
\begin{proof}
	Let $ u \in \jormod $. Then $ \alpha(u) = (\quotepi(u), \jorprojone(u)) $. Thus $ \alpha $ lies in the kernel of $ \alpha $ if and only if it lies in the kernels of $ \jorprojone $ and $ \quotepi $. Since the kernel of $ \quotepi $ is $ \invsub $, the assertion follows.
\end{proof}

We can now put everything together to obtain the following main result.

\begin{theorem}\label{BC:jordanmodule:class-general}
	Let $ \ring $ be an alternative ring, let $ \rinvmap $ be a nuclear involution on $ \ring $ and let $ \jormodtup = (\jormod, \jorsc, \jorprojone, \jorTrone, \psi) $ be a Jordan module over $ (\ring, \rinvmap) $. Choose any involutory submodule $ \invsub $ of $ \jormodtup $. Put $ \module \defl \jormod / \invsub $. Then $ \jorprojone $ induces a pseudo-quadratic module structure on $ \module $ with respect to $ (\ring, \jorprojone(\invsub), \rinvmap) $ and there exists a surjective homomorphism $ \map{\alpha}{\jormod}{\psgr(\module)}{}{} $ of Jordan modules whose kernel is precisely the kernel of $ \restrict{\jorprojone}{\invsub} $.
\end{theorem}

Using the criterion from \cref{BC:jordanmodule:invsub-crit-2inv}, we can simplify the main result as follows.

\begin{theorem}\label{BC:jordanmodule:class-2inv}
	Let $ \ring $ be an alternative ring such that $ 2_\ring $ is invertible, let $ \rinvmap $ be a nuclear involution on $ \ring $ and let $ \jormodtup $ be a Jordan module over $ (\ring, \rinvmap) $. Then $ \jormodtup $ is of pseudo-quadratic type. That is, there exists a pseudo-quadratic module $ \module $ over $ (\ring, \ringTr_\rinvmap(\ring), \rinvmap) $ such that $ \jormodtup $ is isomorphic to the Jordan module $ \psgr(\module) $ from \cref{BC:jordanmodule-pseudquad-ex}.
\end{theorem}

We can also give a precise classification of Jordan modules of type $ C $ if $ 2_\ring $ is invertible.

\begin{theorem}\label{BC:jordanmodule:class-2inv-C}
	Let $ \ring $ be an alternative ring such that $ 2_\ring $ is invertible, let $ \rinvmap $ be a nuclear involution on $ \ring $ and let $ \jormodtup $ be a Jordan module over $ (\ring, \rinvmap) $ of type $ C $. Put $ \ringzero \defl \ringTr_\rinvmap(\ring) = \symring $. Then there exist an $ \ring $-module $ \module $ over $ \ring $ such that $ \jormodtup $ is isomorphic to the Jordan module $ \ringzero \times \module $ which was constructed in \cref{BC:jordanmodule-module-ex}.
\end{theorem}
\begin{proof}
	This follows from \cref{BC:jordanmodule:class-2inv} and the concluding remarks in \cref{BC:jordanmodule-pseudquad-ex}.
\end{proof}

\begin{remark}\label{BC:jordanmodule:class-pseudquad}
	Let $ (\ring, \ringzero, \rinvmap) $ be an involutory set and let $ N $ be a pseudo-qua\-dra\-tic module over $ (\ring, \ringzero, \rinvmap) $. Recall from \cref{BC:jordanmodule:invsub-ex} that $ \invsub \defl \compactSet{0_\jormod} \times \ringzero $ is a unital embedding involutory submodule of $ \psgr(N) $, independent of any assumption on the invertibility of $ 2_\ring $. For this choice of $ \invsub $, the isomorphism $ \alpha $ is (up to canonical identifications) the identity map.
\end{remark}

Using \cref{BC:jordanmodule:class-pseudquad}, we can prove the following corollary of \cref{BC:jordanmodule:class-general}.

\begin{proposition}
	Let $ (\ring, \ringzero, \rinvmap) $ be an involutory set, let $ N $ be a pseudo-quadratic module over $ (\ring, \ringzero, \rinvmap) $ and assume that $ \jormodtup $ is the Jordan module $ \psgr(N) $. Let $ \ringzero' $ be a subset of $ \ringzero $ such that $ (\ring, \ringzero', \rinvmap) $ is an involutory set. Then there exists a pseudo-quadratic module $ N' $ over $ (\ring, \ringzero', \rinvmap) $ such that $ \psgr(N) $ is isomorphic to $ \psgr(N') $ as a Jordan module.
\end{proposition}
\begin{proof}
	By \cref{BC:jordanmodule:invsub-ex}, $ \invsub' \defl \compactSet{0_N} \times \ringzero' \subs \psgr(N) $ is an embedding involutory submodule of $ \psgr(\module) $. Further, we have $ \jorprojone(\invsub') = \ringzero' $. Thus it follows from \cref{BC:jordanmodule:class-general} that $ \psgr(N) $ is isomorphic as a Jordan module to $ \psgr(N') $ for some pseudo-quadratic module $ N' $ over $ (\ring, \ringzero', \rinvmap) $.
\end{proof}


\section{Purely Alternative Rings and Jordan Modules}

\label{sec:pure-alt}

Recall from \cref{pseud:alternative-note} that we can easily construct non-trivial pseudo-quadratic modules $ \module $ over rings of the form $ \ring \dirsum \ring' $ where $ \ring $ is associative and $ \ring' $ is alternative. However, we observed that only the associative ideal $ \ring \dirsum \compactSet{0} $ of $ \ring \dirsum \ring' $ is relevant for the pseudo-quadratic module structure on $ \module $. Thus we attributed the existence of a pseudo-quadratic module over $ \ring \dirsum \ring' $ to the fact that $ \ring \dirsum \ring' $ is not \enquote{purely alternative}, though we gave no precise definition of this property.

In this section, we present two closely related definitions of pure alternativity which were proposed by Slater in \cite{Slater_NuclCentAlt} for not necessarily unital rings. We will show that purely alternative rings in the stronger sense do not admit any non-zero modules and that purely alternative rings in the weaker sense do not admit pseudo-quadratic modules with non-trivial skew-hermitian form. Some more information on purely alternative rings can be found in \cite[Section~9.4]{McCrimmonAltUnpublished}.

The material in this section is largely independent from the remaining part of this book and will not be referenced in any formal argument. Instead, its main purpose is to convince the reader that pseudo-quadratic modules are objects which do not properly belong to the world of alternative rings.

To avoid confusion, we explicitly point out that Jordan modules, unlike pseudo-quadratic modules, have their place in the alternative setting: We have already seen an important class of Jordan modules over arbitrary alternative rings with nuclear involution in \cref{BC:jordanmodule-invset-ex}. This class consists precisely of the Jordan modules $ \psgr(\module) $ where $ \module = \compactSet{0} $.

\begin{definition}[Nuclear ideal]
	Let $ \ring $ be a ring. An ideal of $ \ring $ is called \defemph*{nuclear}\index{ideal (in a ring)!nuclear} if it is contained in the nucleus of $ \ring $.
\end{definition}

\begin{remark}
	The nucleus of a ring $ \ring $ is not (necessarily) a nuclear ideal because it is not (necessarily) an ideal. However, there exists a unique maximal nuclear ideal of $ \ring $: the sum of all nuclear ideals.
\end{remark}

\begin{lemma}\label{altring:annihil-ideal}
	Let $ \ring $ be an alternative ring and let $ N $ be a nuclear ideal. Then the left annihilator
	$ L \defl \Set{a \in \ring \given aN = \compactSet{0}} $
	is an ideal of $ \ring $.
\end{lemma}
\begin{proof}
	It is clear that $ L $ is closed under addition. Now let $ a \in L $ and $ x \in \ring $. We want to show that $ ax $ and $ xa $ lie in $ L $, so let $ n \in N $ be arbitrary. We clearly have $ (xa)n = x(an) = x 0_\ring = 0_\ring $ because $ n $ lies in the nucleus. Hence $ xa $ lies in $ L $. Further, since $ N $ is an ideal, we have $ xn \in N $. Thus $ (ax)n = a(xn) = 0_\ring $ because $ a \in L $. The assertion follows.
\end{proof}

\begin{definition}[Purely alternative ring, {\cite[4.1, 4.2]{Slater_NuclCentAlt}}]
	Let $ \ring $ be an alternative ring.
		We say that $ \ring $ is \defemph*{purely alternative in the weak sense}\index{ring!alternative!purely} if it has no non-zero nuclear ideal.
		We say that $ \ring $ is \defemph*{purely alternative in the strong sense}\index{ring!alternative!purely} if $ \assoc{\ring}{\ring}{\ring} = \ring $. Here $ \assoc{\ring}{\ring}{\ring} $ denotes the ideal generated by all associators.
\end{definition}

\begin{note}
	Purely alternative rings in the weak sense are simply called purely alternative rings in \cite[4.1]{Slater_NuclCentAlt}. Our notion of purely alternative rings in the strong sense is suggested in \cite[4.2]{Slater_NuclCentAlt} as another \enquote{reasonable} definition. Further, it is then shown that any purely alternative ring in the strong sense is also purely alternative in the weak sense if a mild extra condition is satisfied. This extra condition is always satisfied in our setting because rings are assumed to be unital. We will see the precise formulation of this argument in \cref{altring:pure-implication}.
\end{note}

\begin{example}[Direct sums]
	Let $ \ring $ be an associative ring and let $ \ring' $ be an alternative ring. Assume that both rings are not zero. Then $ \ring \dirsum \compactSet{0} $ is a nuclear ideal in $ \ring \dirsum \ring' $, so $ \ring \dirsum \ring' $ is not purely alternative in the weak sense. Further, the associator ideal of $ \ring \dirsum \ring' $ is contained in $ \compactSet{0} \dirsum \ring' $, so $ \ring \dirsum \ring' $ is not purely alternative in the strong sense either.
\end{example}

\begin{example}
	Let $ \ring $ be a ring which is simple in the sense that it contains no non-trivial ideal. Then $ \ring $ is purely alternative (in either sense) if and only if it is alternative and not associative.
	
	(Note that in the famous theorem from \cite{Kleinfeld_SimpleAltRings} that every simple alternative ring is either associative or a Cayley-Dickson algebra over a field, a ring is called simple if it has no non-trivial ideal \emph{and} contains at least one element which is not nilpotent.)
\end{example}

\begin{remark}
	An alternative ring is purely alternative in the strong sense if and only if it has no non-zero homomorphic image which is associative.
\end{remark}

\begin{lemma}[{\cite[Section~3]{Slater_NuclCentAlt}}]\label{altring:pure-implication-lem}
	Let $ \ring $ be an alternative ring and let $ N $ be a nuclear ideal of $ \ring $. Then $ \assoc{\ring}{\ring}{\ring} N = \compactSet{0} $.
\end{lemma}
\begin{proof}
	Let $ x,y,z \in \ring $ and let $ n \in N $. It follows from the nuclear slipping formula (\cref{altring:nucl-slip}) that $ \assoc{x}{y}{z}n = \assoc{x}{y}{zn} $. Since $ N $ is an ideal, $ zn $ lies in $ N $ and thus in the nucleus. Hence $ \assoc{x}{y}{zn} = 0_\ring $. This shows that $ \assoc{x}{y}{z} $ lies in the left annihilator of $ N $. By \cref{altring:annihil-ideal}, it follows that the whole ideal $ \assoc{\ring}{\ring}{\ring} $ which is generated by all associators is trivial.
\end{proof}

\begin{lemma}\label{altring:pure-implication}
	Let $ \ring $ be a purely alternative ring in the strong sense. Then $ \ring $ is also a purely alternative ring in the weak sense.
\end{lemma}
\begin{proof}
	It follows from \cref{altring:pure-implication-lem} that $ \ring N = \compactSet{0} $ for any nuclear ideal $ N $. Since $ \ring $ contains a unit element $ 1_\ring $, this implies that every nuclear ideal is zero. This says precisely that $ \ring $ is purely alternative in the weak sense.
\end{proof}

We can now prove the results that we announced in the introduction of this section.

\begin{proposition}
	Let $ \ring $ be a purely alternative ring in the strong sense and let $ \module $ be an $ \ring $-module. Then $ \module = \compactSet{0} $.
\end{proposition}
\begin{proof}
	We have seen in \cref{ring:module-obs} that modules over $ \ring $ are actually modules over $ \ring / \assoc{\ring}{\ring}{\ring} $. In other words, we have $ vr = 0_\module $ for all $ v \in \module $ and all $ r \in \assoc{\ring}{\ring}{\ring} $. Since $ \ring $ is purely alternative in the strong sense, this means that $ \module \ring = \compactSet{0_\module} $. At the same time, we have $ \module 1_\ring = \module $. We infer that $ \module = \compactSet{0} $.
\end{proof}

\begin{proposition}\label{pure-alt:weak-thm}
	Let $ (\ring, \ringzero, \rinvmap) $ be an involutory set in which $ \ring $ is purely alternative in the weak sense and let $ (\module, q, f) $ be a pseudo-quadratic module over $ (\ring, \ringzero, \rinvmap) $. Then $ f=0 $. In particular, $ q $ is trivial if $ 2_\ring $ is invertible.
\end{proposition}
\begin{proof}
	Denote by $ I $ the additive subgroup of $ (\ring, +) $ which is generated by the image of $ f $. It follows from the sesquilinearity of $ f $ that $ I $ is an ideal of $ \ring $. Further, $ I $ is nuclear by \cref{pseud:sesqui-def}. Since $ \ring $ is purely alternative in the weak sense, we infer that $ I = \compactSet{0} $. Hence $ f=0 $. If $ 2_\ring $ is invertible, then the triviality of $ f $ implies the triviality of $ q $ by \thmitemcref{pseud:2}{pseud:2:inv-triv}.
\end{proof}

	\chapter{Root Gradings of Types \texorpdfstring{$ C $}{C} and \texorpdfstring{$ BC $}{BC}}
	
	\label{chap:BC}
	
	In this chapter, we study $ C_n $-graded groups and $ BC_n $-graded groups for $ n \ge 2 $, though we will assume that $ n \ge 3 $ for the main results. We will see that $ C_n $-graded groups can be regarded as special cases of $ BC_n $-graded groups. This allows us to restrict our attention to $ BC_n $-gradings most of the time. Similarly to the situation in root gradings of type $ B $, every medium-length root in $ BC_n $ is contained in an $ A_2 $-subsystem (\cref{BC:rootsys:med-in-A2}), so it is clear that there exists a ring $ \ring $ which coordinatises the medium-length root groups. We will show that $ \ring $ is alternative and that it has a nuclear involution $ \rinvmap $. Further, we will parametrise the short root groups by a (not necessarily abelian) group $ (\jormod, \joradd) $. The commutator relations yield maps $ \jorsc $, $ \jorprojone $, $ \jorTrone $, $ \psi $ which equip $ \jormod $ with the structure of a Jordan module over~$ (\ring, \rinvmap) $.
	
	Some partial results on $ C_n $-gradings for $ n \ge 3 $ have already been obtained by Zhang in his PhD thesis \cite{Zhang}. See~\ref{rgg-lit:coord} for more details. For RGD-systems of type $ BC_2 $, the standard reference are Chapters~25 and~26 in \cite{MoufangPolygons}.
	
	As usual, this chapter follows the outline described in \cref{sec:param:outline}, except that we have already covered all the relevant algebraic structures in \cref{chap:BC-alg}. Hence we can begin with the investigation of the root systems $ C_n $ and $ BC_n $ in \cref{sec:BC-rootsys}. In \cref{sec:CisBC}, we will uncover several connections between root gradings of types $ B $, $ C $ and $ BC $. In particular, we will show how root gradings of type $ C $ can be regarded as root gradings of type $ BC $ with additional properties. We will also see that the rank-2 results from \cref{chap:B} can also be applied, in a modified way, to root gradings of type $ BC $. \Cref{sec:BC-example} presents the construction of a $ BC_n $-graded group from a standard pseudo-quadratic module. This does not provide a complete solution of the construction problem for $ BC_n $-graded groups, and we will summarise what precisely remains to be done.
	
	After this point, the coordinatisation of $ BC_n $-graded groups begins. In \cref{sec:BC:rank2,sec:BC:rank3}, we will perform the necessary rank-2 and rank-3 computations. We introduce the notions of standard signs and standard partial twisting systems for $ BC_n $-graded groups in \cref{sec:BC:stsigns,sec:BC:sttwist}, respectively. In \cref{sec:BC-param}, we construct the parametrising groups $ (\ring, +) $ and $ (\jormod, \joradd) $. The commutation maps, their rank-2 identities and the blueprint rewriting rules are defined, derived and computed in \cref{sec:BC-bluerules}. Finally, we perform the blueprint computations for $ BC_n $ in \cref{sec:BC:blue-comp}, and we state our final result in \cref{BC:thm}.
	

\section{Root Systems of Types \texorpdfstring{$ C $}{C} and \texorpdfstring{$ BC $}{BC}}

\label{sec:BC-rootsys}

\begin{secnotation}
	We denote by $ n $ an integer at least $ 1 $.
\end{secnotation}

In this section, we collect some basic facts about root systems of types $ C $ and $ BC $ which will be needed later on. As usual, all results in this section are straightforward to verify. We will see in \cref{BC:CisBC:CasBC} that root gradings of type $ C_n $ are special cases of $ BC_n $-graded groups. Hence we can restrict ourselves to root systems of type $ BC $ most of the time.

\begin{remark}[Standard representation of $ BC_n $]\label{BC:BCn-standard-rep}
	Let $ V $ be a Euclidean space of dimension $ n $ with orthonormal basis $ \tup{\basvec}{n} $. The \defemph*{standard representation of $ BC_n $}\index{standard representation!of BCn@of $ BC_n $} is
	\begin{align*}
		BC_n &\defl \Set{\epsilon_1 \basvec_i +\epsilon_2 \basvec_j \given i \ne j \in \numint{1}{n}, \epsilon_1, \epsilon_2 \in \compactSet{\pm 1}} \union \Set{\epsilon \basvec_i \given i \in \numint{1}{n}, \epsilon \in \compactSet{\pm 1}} \\
		&\hspace{2cm} \mathord{} \union \Set{2\epsilon \basvec_i \given i \in \numint{1}{n}, \epsilon \in \compactSet{\pm 1}}
	\end{align*}
	The long roots are exactly those which lie in the third set, the medium-length roots are those in the first set and the short roots are those in the second one. The \defemph*{standard root base} is
	\[ \rootbase \defl \Set{\basvec_i - \basvec_{i+1} \given i \in \numint{1}{n-1}} \union \compactSet{\basvec_n}, \]
	the \defemph*{standard rescaled root base} is
	\[ \rootbase \defl \Set{\basvec_i - \basvec_{i+1} \given i \in \numint{1}{n-1}} \union \compactSet{2\basvec_n} \]
	and the corresponding positive system is
	\begin{align*}
		\possys &\defl \Set{\basvec_i - \basvec_{j} \given i<j \in \numint{1}{n}} \union \Set{\basvec_i + \basvec_j \given i \ne j \in \numint{1}{n}} \\
		& \qquad \mathord{}\union \Set{\lambda\basvec_i \given i \in \numint{1}{n}, \lambda \in \compactSet{1,2}}.
	\end{align*}
	Recall \cref{param:motiv:rescaled} for a discussion why we usually prefer the standard rescaled root base over the standard root base.
\end{remark}

\begin{remark}[Standard representation of $ C_n $]\label{BC:Cn-standard-rep}
	Let $ V $ be a Euclidean space of dimension $ n $ with orthonormal basis $ \tup{\basvec}{n} $. Assume that $ n \ge 2 $. (As for $ B_n $, we exclude the case $ n=1 $ because the resulting root system $ C_1 $ would be isomorphic to $ A_1 $.) The \defemph*{standard representation of $ C_n $}\index{standard representation!of Cn@of $ C_n $} is
	\begin{align*}
		C_n &\defl \Set{\epsilon_1 \basvec_i +\epsilon_2 \basvec_j \given i \ne j \in \numint{1}{n}, \epsilon_1, \epsilon_2 \in \compactSet{\pm 1}} \union \Set{2\epsilon \basvec_i \given i \in \numint{1}{n}, \epsilon \in \compactSet{\pm 1}}.
	\end{align*}
	We can see it as a crystallographically closed subset of (the standard representation of) $ BC_n $ in a natural way. The long roots are exactly those which lie in the second set and the short roots are those in the first one. The \defemph*{standard root base} is
	\[ \rootbase \defl \Set{\basvec_i - \basvec_{i+1} \given i \in \numint{1}{n-1}} \union \compactSet{2\basvec_n} \]
	and the corresponding positive system is
	\[ \possys \defl \Set{\basvec_i - \basvec_{j} \given i<j \in \numint{1}{n}} \union \Set{\basvec_i + \basvec_j \given i \ne j \in \numint{1}{n}} \union \Set{2\basvec_i \given i \in \numint{1}{n}}. \]
\end{remark}

\begin{warning}
	Roots of the form $ \pm \basvec_i \pm \basvec_j $ are short roots when considered as roots in $ C_n $ but of medium length when considered roots in $ BC_n $. However, it will always be clear from the context in which root system we are working.
\end{warning}

\begin{observation}
	Since $ \refl{\lambda \alpha} = \refl{\alpha} $ for all $ \lambda \in \IR \setminus \compactSet{0} $, the root systems $ B_n $, $ C_n $ and $ BC_n $ have the same Weyl group.
\end{observation}

\begin{definition}[$ BC_2 $-pairs and $ BC_2 $-quadruples]\label{BC:BC2-pair-def}
	Let $ \roots $ be any root system. A \defemph*{$ BC_2 $-pair (in $ \roots $)}\index{BC2-pair@$ BC_2 $-pair} is a pair $ (\alpha, \delta) $ of roots such that $ \gen{\alpha, \delta}_\IR \intersect \roots $ is a root subsystem of $ \roots $ of type $ BC_2 $ with root base $ (\alpha, \delta) $ with $ \alpha $ being the medium-length simple root in this subsystem and $ \delta $ being the short simple root. A \defemph*{$ BC_2 $-quadruple}\index{BC2-quadruple@$ BC_2 $-quadruple} is a quadruple $ (\alpha, \beta, \gamma, \delta) $ of roots in $ \roots $ such that $ (\alpha, \delta) $ is a $ BC_2 $-pair, $ \beta = \alpha + \delta $ and $ \gamma = \alpha + 2\delta $.
\end{definition}

We can define $ C_2 $-pairs and $ C_2 $-quadruples in an analogous way, although they will rarely be needed. 

\begin{definition}[$ C_2 $-pairs and $ C_2 $-quadruples]\label{C:BC2-pair-def}
	Let $ \roots $ be any root system. A \defemph*{$ C_2 $-pair (in $ \roots $)}\index{BC2-pair@$ BC_2 $-pair} is a pair $ (\alpha, \delta) $ of roots such that $ \gen{\alpha, \delta}_\IZ \intersect \roots $ is a root subsystem of $ \roots $ of type $ C_2 $ with root base $ (\alpha, \delta) $ such that $ \delta $ is longer than $ \alpha $. A \defemph*{$ C_2 $-quadruple}\index{C2-quadruple@$ C_2 $-quadruple} is a quadruple $ (\alpha, \beta, \gamma, \delta) $ of roots in $ \roots $ such that $ (\alpha, \delta) $ is a $ C_2 $-pair, $ \beta = 2\alpha + \delta $ and $ \gamma = \alpha + \delta $.
\end{definition}

\begin{warning}
	Note that we consider $ \gen{\alpha, \delta}_\IZ \intersect \roots $ in \cref{C:BC2-pair-def}, not $ \gen{\alpha, \delta}_\IR \intersect \roots $: The latter subsystem is allowed to be of type $ BC_2 $.
\end{warning}

\begin{remark}
	Assume that $ (\alpha, \beta, \gamma, \delta) $ is a $ BC_2 $-quadruple. Then the scaled tuple $ (\alpha, 2\beta, \gamma, 2\delta) $ is a $ C_2 $-quadruple.
\end{remark}

\begin{remark}[compare \cref{B:B2-cry-crit}]\label{BC:C2-cry-crit}
	Let $ G $ be a group with a $ C_2 $-grading $ (\rootgr{\alpha})_{\alpha \in B_2} $ and let $ (\alpha, \beta, \gamma, \delta) $ be a $ C_2 $-quadruple. Then this grading is crystallographic if and only if the commutators $ \commutator{\rootgr{\epsilon\beta}}{\rootgr{\sigma\delta}} $ are trivial for all $ \epsilon,\sigma \in \compactSet{\pm 1} $.
\end{remark}

\premidfigure
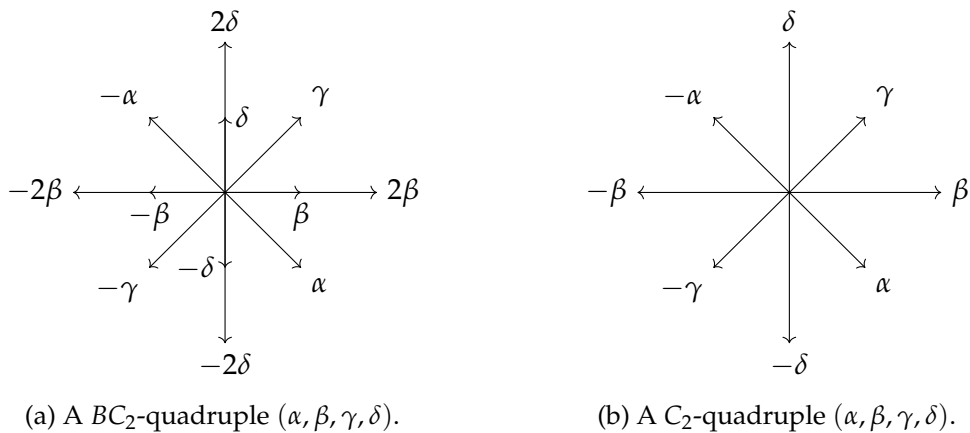
\begin{figure}[htb]
	\centering\begin{subfigure}{0.45\linewidth}
		\centering\begin{tikzpicture}
			\draw[->] (0,0) -- (1,-1);
			\draw[->] (0,0) -- (1,0);
			\draw[->] (0,0) -- (2,0);
			\draw[->] (0,0) -- (1,1);
			\draw[->] (0,0) -- (0,1);
			\draw[->] (0,0) -- (0,2);
			\draw[->] (0,0) -- (-1,1);
			\draw[->] (0,0) -- (-1,0);
			\draw[->] (0,0) -- (-2,0);
			\draw[->] (0,0) -- (-1,-1);
			\draw[->] (0,0) -- (0,-1);
			\draw[->] (0,0) -- (0,-2);
			
			\node[below right] at (1,-1){$ \alpha $};
			\node[below] at (1,0){$ \beta $};
			\node[right] at (2,0){$ 2\beta $};
			\node[above right] at (1,1){$ \gamma $};
			\node[right] at (0,1){$ \delta $};
			\node[above] at (0,2){$ 2\delta $};
			\node[above left] at (-1,1){$ -\alpha $};
			\node[below] at (-1,0){$ -\beta $};
			\node[left] at (-2,0){$ -2\beta $};
			\node[below left] at (-1,-1){$ -\gamma $};
			\node[left] at (0,-1){$ -\delta $};
			\node[below] at (0,-2){$ -2\delta $};
		\end{tikzpicture}
		\caption{A $ BC_2 $-quadruple $ (\alpha, \beta, \gamma, \delta) $.}
		\label{fig:BC2-quad}
	\end{subfigure}\hspace{0.1\linewidth}%
	\begin{subfigure}{0.45\linewidth}
		\centering\begin{tikzpicture}
			\draw[->] (0,0) -- (1,-1);
			\draw[->] (0,0) -- (2,0);
			\draw[->] (0,0) -- (1,1);
			\draw[->] (0,0) -- (0,2);
			\draw[->] (0,0) -- (-1,1);
			\draw[->] (0,0) -- (-2,0);
			\draw[->] (0,0) -- (-1,-1);
			\draw[->] (0,0) -- (0,-2);
			
			\node[below right] at (1,-1){$ \alpha $};
			\node[right] at (2,0){$ \beta $};
			\node[above right] at (1,1){$ \gamma $};
			\node[above] at (0,2){$ \delta $};
			\node[above left] at (-1,1){$ -\alpha $};
			\node[left] at (-2,0){$ -\beta $};
			\node[below left] at (-1,-1){$ -\gamma $};
			\node[below] at (0,-2){$ -\delta $};
		\end{tikzpicture}
		\caption{A $ C_2 $-quadruple $ (\alpha, \beta, \gamma, \delta) $.}
		\label{fig:C2-quad}
	\end{subfigure}
	\caption{$ C_2 $- and $ BC_2 $-quadruples.}
	\label{fig:C2BC2-quad}
\end{figure}
\postmidfigure

We now prove similar properties as in \cref{sec:B:rootsys}. Unless otherwise specified, the word \enquote{root} always refers to roots in $ BC_n $.

\begin{lemma}\label{BC:rootsys:med-in-A2}
	Assume that $ n \ge 3 $ and let $ \alpha $ be a medium-length root in $ BC_n $. Then there exists an $ A_2 $-subsystem of $ BC_n $ which contains $ \alpha $.
\end{lemma}

\begin{lemma}
	Every subsystem of $ BC_n $ of type $ A_2 $ or $ BC_2 $ is parabolic.
\end{lemma}

\begin{lemma}\label{BC:rootsys:BC2-trans}
	The Weyl group of $ BC_n $ acts transitively on the set of $ BC_2 $-pairs and on the set of $ BC_2 $-quadruples.
\end{lemma}
\begin{proof}
	Let $ (\alpha, \delta) $, $ (\alpha', \delta') $ be two $ BC_2 $-pairs. By \thmitemcref{rootsys:subsys}{rootsys:subsys:base}, there exist root bases $ \rootbase $ and $ \rootbase' $ containing $ (\alpha, \delta) $ and $ (\alpha', \delta') $, respectively. By \cref{rootsys:weyl-base-trans}, there exists an element $ u $ of the Weyl group which maps $ \rootbase $ to $ \rootbase' $. Note that $ (\alpha, \delta) $ is the unique pair of elements in $ \rootbase $ which forms a $ BC_2 $-pair, and similarly for $ (\alpha', \delta') $ and $ \rootbase' $. It follows from this characterisation that $ u $ maps $ (\alpha, \delta) $ to $ (\alpha', \delta') $. Thus the Weyl group acts transitively on the set of $ BC_2 $-pairs. This implies that it acts transitively on the set of $ BC_2 $-quadruples as well.
\end{proof}

\begin{lemma}\label{BC:rootsys:med-decomp}
	Let $ (\alpha, \beta, \gamma, \delta) $ be a $ BC_2 $-quadruple in $ BC_n $ and assume that $ n \ge 3 $. Then there exist medium-length roots $ \gamma_1, \gamma_2 $ such that the following conditions are satisfied:
	\begin{stenumerate}
		\item $ (\gamma_1, \gamma, \gamma_2) $ is an $ A_2 $-triple. In particular, $ \gamma = \gamma_1 + \gamma_2 $.
		
		\item $ (\alpha, \gamma_1) $ is an $ A_2 $-pair.
		
		\item $ (\alpha, -\gamma_2) $ is an $ A_2 $-pair.
	\end{stenumerate}
\end{lemma}
\begin{proof}
	Without loss of generality, we can use the standard representation of $ BC_n $ and assume that $ \alpha = \basvec_1 - \basvec_2 $ and $ \gamma = \basvec_1 + \basvec_2 $. Then the roots $ \gamma_1 \defl \basvec_2 - \basvec_3 $ and $ \gamma_2 \defl \basvec_1 + \basvec_3 $ have the desired properties.
\end{proof}

\begin{lemma}\label{BC:rootsys:common-BC2}
	Assume that $ n \ge 2 $. Let $ \alpha \in BC_n $ be a medium-length root and let $ \delta \in BC_n $ be short. Then $ \alpha $ and $ \delta $ lie in a common $ BC_2 $-subsystem if and only if $ \alpha \cdot \delta \ne 0 $.
\end{lemma}

Let $ v, w $ be any non-zero vectors in the Euclidean space surrounding $ BC_n $. By the definition of the Cartan numbers, we have $ \cartanint{\lambda v}{w} = \lambda \cartanint{v}{w} $ and $ \cartanint{v}{\lambda w} = \lambda^{-1} \cartanint{v}{w} $ for all $ \lambda \in \IR \setminus \compactSet{0} $. Thus as soon as we have computed the Cartan integers for all short and medium-length roots, we can easily derive the Cartan integers for all roots in $ BC_n $. Since the Cartan integer of short and medium-length roots stay the same when we consider them as roots in $ B_n $, we can deduce the following characterisation from \cref{B:rootsys:cartan-int-parity}.

\begin{proposition}\label{BC:rootsys:cartan-int-parity}
	Let $ \rho, \zeta $ be two roots in $ BC_n $. Then $ \cartanint{\rho}{\zeta} $ is an even number if and only if one of the following conditions is satisfied:
	\begin{lemenumerate}
		\item $ \rho $ is long.
		
		\item $ \zeta $ is short.
	
		\item $ \rho \in \Set{\pm \zeta} $.
		
		\item $ \rho $ and $ \zeta $ are orthogonal.
	\end{lemenumerate}
\end{proposition}

The following lemma will be used in \cref{BC:stab-comp} to prove stabiliser-compatibility for $ BC_n $-graded groups. Observe that the set $ \bar{\calA} $ is defined exactly as in \cref{param:stabcomp-crit-ortho}.

\begin{lemma}\label{BC:rootsys:ortho-adj}
	Assume that $ n \ge 2 $. Let $ \alpha $ be a medium-length root, choose roots $ \beta, \gamma, \delta $ such that $ (\alpha, \beta, \gamma, \delta) $ is a $ BC_2 $-quadruple and put
	\[ \bar{\calA} \defl \alpha^\perp \intersect \Set{\rho \in BC_n \given \rho \text{ is not crystallographically adjacent to } \alpha}. \]
	Then $ \bar{\calA} = \Set{\gamma, -\gamma} $. In particular, for any root base $ \rootbase $ of $ BC_n $, there exists exactly one $ \rootbase $-positive root in $ \bar{\calA} $, and this root has medium length.
\end{lemma}

\begin{lemma}\label{BC:rootsys:ortho-adj-2}
	Assume that $ n \ge 2 $. Let $ \alpha, \alpha' $ be two orthogonal medium-length roots which are crystallographically adjacent. Then $ \alpha $ and $ -\alpha' $ are crystallographically adjacent as well.
\end{lemma}

\begin{lemma}\label{BC:rootsys:orth-adj-trans}
	Assume that $ n \ge 2 $. Then the Weyl group of $ BC_n $ acts transitively on the set
	\[ S \defl \Set*{(\alpha, \gamma) \given \begin{gathered}
		\alpha, \gamma \text{ are of medium length, orthogonal} \\
		\text{and not crystallographically adjacent}
	\end{gathered}}. \]
\end{lemma}
\begin{proof}
	The set $ S $ consists precisely of the pairs $ (\alpha, \gamma) $ for which there exist roots $ \beta, \delta $ such that $ (\alpha, \beta, \gamma, \delta) $ is a $ BC_2 $-quadruple. Thus the assertion follows from \cref{BC:rootsys:BC2-trans}.
\end{proof}

Recall that we have defined a certain subset $ \Bnsub $ of $ B_n $ in \cref{B:Bnsub-def}. Since $ C_n $ and $ BC_n $ have the same Weyl group as $ B_n $, the same reasoning leads us to the definition of the set $ \Cnsub \subs C_n \subs BC_n $. We will mainly use this set in \cref{BC:comm-mult-computation,BC:blue:stand-signs}. In this context, it is worth recalling \cref{param:motiv:rescaled}.

\begin{definition}\label{BC:Cnsub-def}
	For all $ n \in \IN_{\ge 2} $, we define the following subset of $ C_n $ and $ BC_n $:\index{Cn@$ \Cnsub $}
	\[ \Cnsub \defl \Set{\basvec_i - \basvec_j \given i \ne j \in \numint{1}{n}} \union \Set{2\basvec_i \given i \in \numint{1}{n}} \subs C_n. \]
\end{definition}

\begin{definition}[Standard $ \rootbase $-expression]\label{BC:Cnsub-ext-word}
	Let $ \rootbase $ denote the standard root base of $ C_n $, considered as a subset of $ BC_n $. For any root $ \alpha \in \Cnsub $, we define a $ \rootbase $-expression $ \word{\rho}^{\alpha} $ of $ \alpha $ (in the sense of \cref{param:Delta-expr}) as follows:
	\begin{defenumerate}
		\item If $ \alpha \in \rootbase $, we put $ \word{\rho}^\alpha \defl (\alpha) $.
		
		\item If $ \alpha = 2\basvec_i $ for some $ i \in \numint{1}{n-1} $, we define $ \word{\rho}^{\basvec_i} \defl (\basvec_i - \basvec_{i+1}, \word{\rho}^{\basvec_{i+1}}, \basvec_{i+1}-\basvec_i) $.
		
		\item If $ \alpha = \basvec_i - \basvec_j $ for some $ i<j \in \numint{1}{n-1} $ with $ i+1<j $, we define $ \word{\rho}^{\basvec_i - \basvec_j} \defl (\basvec_j - \basvec_{j-1}, \word{\rho}^{\basvec_i - \basvec_{j-1}}, \basvec_{j-1} - \basvec_j) $.
		
		\item If $ \alpha = \basvec_j - \basvec_i $ for some $ i<j \in \numint{1}{n-1} $ with $ i+1<j $, we define $ \word{\rho}^{\basvec_j - \basvec_i} \defl (\word{\rho}^{\basvec_i - \basvec_j})^{-1} $.
	\end{defenumerate}
	The word $ \word{\rho}^{\alpha} $ will also be called the \defemph*{standard $ \rootbase $-expression of $ \alpha $}\index{standard Delta-expression@standard $ \rootbase $-expression}. We also put $ \word{\rho}^i \defl \word{\rho}^{\basvec_i} $ for all $ i \in \numint{1}{n} $ and $ \word{\rho}^{ij} \defl \word{\rho}^{\basvec_i - \basvec_j} $ for all distinct $ i,j \in \numint{1}{n} $.
\end{definition}

\begin{note}
	Let $ n \in \IN_{\ge 2} $. Denote by $ \rootbase_B $ the standard root base of $ B_n $ and by $ \rootbase_C $ the standard root base of $ C_n $. Let $ \alpha \in A_{n-1} $. Then $ \alpha $ can be regarded as a root in $ B_n $ and as a root in $ C_n $, so we have a standard $ \rootbase_B $-expression $ \word{\rho}^\alpha_B $ and a standard $ \rootbase_C $-expression $ \word{\rho}^\alpha_C $. In fact, these expression only contain roots from the canonical $ A_{n-1} $-subsystems of $ B_n $ and $ C_n $, and we have $ \word{\rho}^\alpha_B = \word{\rho}^\alpha_C $.
	
	Now let $ i \in \numint{1}{n} $. Denote by $ \word{\rho}^i_B $ the standard $ \rootbase_B $-expression of $ \basvec_i $ and by $ \word{\rho}^i_C $ the standard $ \rootbase_C $-expression of $ 2\basvec_i $. Then $ \word{\rho}^i_C $ is the word obtained from $ \word{\rho}^i_B $ by replacing each occurrence of $ \basvec_n $ by $ 2\basvec_n $.
\end{note}

\begin{definition}[$ \Cnsub $-extensions]\label{BC:Cnsub-ext-def}
	Denote by $ \rootbase $ the standard rescaled root base of $ BC_n $. Let $ G $ be any group with a $ BC_n $-pregrading $ (\hat{U}_\alpha)_{\alpha \in BC_n} $ and let $ (w_\delta)_{\delta \in \rootbase} $ be a $ \rootbase $-system of Weyl elements in $ G $. Then we define a family $ (w_\alpha)_{\alpha \in \Cnsub} $, called the \defemph*{standard $ \Cnsub $-extension of $ (w_\delta)_{\delta \in \rootbase} $}\index{Cn-extension@$ \Cnsub $-extension}, by $ w_\alpha \defl w_{\word{\rho}^{\alpha}} $ for all $ \alpha \in \Cnsub $. We will sometimes write $ w_{ij} $ for $ w_{\basvec_i - \basvec_j} $ and $ w_{i} $ for $ w_{2\basvec_i} $.
\end{definition}


\section{The Relationship between Root Gradings of Types \texorpdfstring{$ B $}{B}, \texorpdfstring{$ C $}{C} and \texorpdfstring{$ BC $}{BC}}

\label{sec:CisBC}

\begin{secnotation}
	We denote by $ n $ and integer at least~$ 2 $ and by $ G $ an arbitrary group.
\end{secnotation}

In this section, we investigate how (crystallographic) root gradings of types $ B $, $ C $ and $ BC $ can be constructed from each other. All assertions in this sections are immediate consequences of the definition of root gradings. The most important results are \cref{BC:CisBC:C2}, which says that crystallographic $ B_2 $-gradings are the same as crystallographic $ C_2 $-gradings, and \cref{BC:CisBC:CinBC}, which says that crystallographic $ BC_n $-graded groups have crystallographic $ C_n $-graded subgroups. Both statements together yield that a crystallographic $ BC_n $-graded group has many crystallographic $ B_2 $-graded subgroups to which we can apply the results of \cref{sec:B:rank2-cry}. Further, \cref{BC:CisBC:CasBC} yields that we do not have to consider crystallographic $ C_n $-gradings separately because they can be seen as a special case of crystallographic $ BC_n $-gradings.

We begin with the situation of non-crystallographic gradings. These observations are not intrinsically relevant for us because we are ultimately interested in crystallographic root gradings. Rather, they illustrate that the crystallographic condition is the key property which differentiates root gradings of types $ B $, $ C $ and $ BC $ from each other.

\begin{remark}[The non-crystallographic case for $ B $ and $ C $]\label{BC:CisBC:C-is-B-noncry}
	If we consider root gradings which are not assumed to be crystallographic, the distinction between $ B_n $ and $ C_n $ disappears. To make this precise, consider $ B_n $ and $ C_n $ as root systems in the same Euclidean space and denote by $ \map{f}{B_n}{C_n}{}{} $ the map which sends $ \pm \basvec_i $ to $ \pm 2\basvec_i $ and which fixes all other roots. Let $ (\rootgr{\alpha})_{\alpha \in B_n} $ be a $ B_n $-pregrading of $ G $. Put $ \rootgr{\beta}' \defl \rootgr{f^{-1}(\beta)} $ for all $ \beta \in C_n $. Then $ (\rootgr{\alpha})_{\alpha \in B_n} $ is a (not necessarily crystallographic) $ B_n $-grading of $ G $ if and only if $ (\rootgr{\beta}')_{\beta \in C_n} $ is a (not necessarily crystallographic) $ C_n $-grading of $ G $.
\end{remark}

\begin{remark}[The non-crystallographic case for $ BC $]\label{BC:BisBC:BC-noncry}
	Let $ (\rootgr{\alpha})_{\alpha \in BC_n} $ be a (not necessarily crystallographic) $ BC_n $-grading of $ G $. Then the group $ \rootgr{\pm 2\basvec_i} $ is a subgroup of $ \rootgr{\pm \basvec_i} $ for all $ i \in \numint{1}{n} $. Further, $ (\rootgr{\alpha})_{\alpha \in B_n} $ is a $ B_n $-gradins of $ G $ (and thus also a $ C_n $-grading of $ G $).
	
	Conversely, let $ (\rootgr{\alpha})_{\alpha \in C_n} $ be a $ C_n $-grading of $ G $. To this grading we can associate a $ BC_n $-grading $ (\rootgr{\alpha})_{\alpha \in BC_n} $ by setting $ \rootgr{\pm \basvec_i} \defl \rootgr{\pm 2\basvec_i} $. By the observations in the previous paragraph, we can construct $ BC_n $-gradings from $ B_n $-gradings in a similar way.
	
	We conclude that a $ BC_n $-grading $ (\rootgr{\alpha})_{\alpha \in BC_n} $ is the same thing as a $ B_n $-grading $ (\rootgr{\alpha})_{\alpha \in B_n} $ (or a $ C_n $-grading) together with an additional family $ (\rootgr{\beta})_{\beta \in BC_n \setminus B_n} $ of subgroups of the short root groups such that the axioms of a root grading are satisfied for $ (\rootgr{\alpha})_{\alpha \in BC_n} $. This additional requirement says precisely that the following statements are satisfied:
	\begin{stenumerate}
		\item For each $ \alpha \in B_n $, there exists an $ \alpha $-Weyl element $ w_\alpha $ with respect to $ (\rootgr{\beta})_{\beta \in B_n} $ such that $ \rootgr{\pm 2\basvec_i}^{w_\alpha} = \rootgr{\pm \refl{\alpha}(2\basvec_i)} $ for all $ i \in \numint{1}{n} $.
		
		\item For all $ i \in \numint{1}{n} $, there exists a $ \pm 2\basvec_i $-Weyl element. That is, there exists a $ (\pm \basvec_i) $-Weyl element which lies in $ \rootgr{\mp 2\basvec_i} \rootgr{\pm 2\basvec_i} \rootgr{\mp 2\basvec_i} $.
	\end{stenumerate}
\end{remark}

We now turn to the case of crystallographic root gradings.

\begin{remark}[The crystallographic case for $ B $ and $ C $]
	It is clear from \cref{B:B2-cry-crit,BC:C2-cry-crit} that the map $ \map{f}{B_n}{C_n}{}{} $ from \cref{BC:CisBC:C-is-B-noncry} no longer induces an equivalence between crystallographic $ B_n $-gradings and crystallographic $ C_n $-gradings. If $ n=2 $, then we can replace $ f $ by a different map to obtain an equivalence between crystallographic $ B_2 $-gradings and crystallographic $ C_2 $-gradings. We will see this in \cref{BC:CisBC:C2}. If $ n>2 $, however, crystallographic $ B_n $-gradings and crystallographic $ C_n $-gradings are not related.
\end{remark}

Some of the constructions in \cref{BC:CisBC:C-is-B-noncry} remain valid in the crystallographic setting.

\begin{remark}[$ C_n $-graded subgroups in $ BC_n $-graded groups]\label{BC:CisBC:CinBC}
	Let $ (\rootgr{\alpha})_{\alpha \in BC_n} $ be a crystallographic $ BC_n $-grading of $ G $. Since $ C_n $ is a closed root subsystem of $ BC_n $, it follows from \cref{rgg:subgroup} that $ (\rootgr{\alpha})_{\alpha \in C_n} $ is a crystallographic $ C_n $-grading, but of $ \gen{\rootgr{\alpha} \given \alpha \in C_n} $ and not of $ G $. However, $ (\rootgr{\alpha})_{\alpha \in B_n} $ is not necessarily a crystallographic $ B_n $-grading of $ G $: For example, it is not necessarily true that $ \rootgr{\basvec_1 - \basvec_2} $ commutes with $ \rootgr{\basvec_1 + \basvec_2} $.
\end{remark}

\begin{remark}[$ C_n $-graded groups as $ BC_n $-graded groups]\label{BC:CisBC:CasBC}
	Now let $ (\rootgr{\alpha})_{\alpha \in C_n} $ be a crystallographic $ C_n $-grading of $ G $. Putting $ \rootgr{\pm \basvec_i} \defl \rootgr{\pm 2\basvec_i} $ for all $ i \in \numint{1}{n} $, we obtain a crystallographic $ BC_n $-grading $ (\rootgr{\alpha})_{\alpha \in BC_n} $ of $ G $. Thus we can regard crystallographic $ C_n $-gradings as special cases of crystallographic $ BC_n $-gradings. For this reason, we will not specifically consider $ C_n $-gradings most of the time in this chapter.
\end{remark}

\begin{remark}[$ BC_n $-graded groups as $ C_n $-graded groups]\label{BC:CisBC:BCasC}
	Let $ (\rootgr{\alpha})_{\alpha \in BC_n} $ be a crystallographic $ BC_n $-grading of $ G $. Put $ \rootgr{2\alpha}' \defl \rootgr{\alpha} $ for all short roots $ \alpha $ in $ BC_n $ and $ \rootgr{\beta}' \defl \rootgr{\beta} $ for all medium-length roots $ \beta $ in $ BC_n $. Then $ (\rootgr{\alpha}')_{\alpha \in C_n} $ is a crystallographic $ C_n $-grading of $ G $ if and only if for all distinct $ i,j \in \numint{1}{n} $ and for all signs $ \epsilon, \sigma \in \Set{\pm 1} $, we have that $ \rootgr{\epsilon\basvec_i} $ commutes with $ \rootgr{\sigma \basvec_j} $. In other words, a crystallographic $ BC_n $-grading can be regarded as a crystallographic $ C_n $-grading if and only if pairs of orthogonal short root groups commute. A crystallographic $ BC_n $-grading is called \defemph*{proper}\index{proper BCn-grading@proper $ BC_n $-grading} if it cannot be regarded as a crystallographic $ C_n $-gradings in this way.
\end{remark}

\begin{remark}[$ B_2 $-graded groups are $ C_2 $-graded groups]\label{BC:CisBC:C2}
	Let $ (\alpha, \delta) $ be a $ B_2 $-pair in $ B_2 $ and let $ (\alpha', \delta') $ be a $ C_2 $-pair in $ C_2 $. Then the unique vector space isomorphism $ \map{g}{\gen{\alpha, \delta}_\IR}{\gen{\alpha', \delta'}_\IR}{}{} $ which maps $ \alpha $ to $ \delta' $ and $ \delta $ to $ \alpha' $ induces an isomorphism from $ B_2 $ to $ C_2 $. Now let $ G $ be a group, let $ (\rootgr{\alpha})_{\alpha \in B_2} $ be a $ B_2 $-pregrading of $ G $ and put $ \rootgr{\beta}' \defl \rootgr{g^{-1}(\beta)} $ for all $ \beta \in C_2 $. Then $ (\rootgr{\alpha})_{\alpha \in B_2} $ is a crystallographic $ B_2 $-grading if and only if $ (\rootgr{\beta}')_{\beta \in C_2} $ is a crystallographic $ C_2 $-grading. Thus crystallographic $ B_2 $-gradings and crystallographic $ C_2 $-gradings are essentially the same thing.
	
	We will usually apply the previous observation in the following form. Let $ G $ be a group with a crystallographic $ C_2 $-grading $ (\rootgr{\zeta})_{\zeta \in C_2} $ and let $ (\alpha, \beta, \gamma, \delta) $ be a $ C_2 $-quadruple. Then we can regard $ (\rootgr{\zeta})_{\zeta \in C_2} $ as a crystallographic $ B_2 $-grading of $ G $ for which $ (\delta, \gamma, \beta, \alpha) $ is a $ B_2 $-quadruple. To this grading, we can apply the results from \cref{sec:B:rank2-cry}. In other words, we can use the results from \cref{sec:B:rank2-cry} by reversing the roles of $ \alpha $ and $ \delta $ (and consequently, of $ \beta $ and $ \gamma $).
	
	Recall from \cref{BC:CisBC:CinBC} that any $ BC_n $-graded group has a canonical crystallographic $ C_n $-graded subgroup and hence many crystallographic $ C_2 $-graded subgroups. We can apply the observations of the previous paragraph to these subgroups as well.
\end{remark}


\section{Construction of \texorpdfstring{$ BC_n $}{BC\_n}-graded Groups}

\label{sec:BC-example}

\begin{secnotation}\label{secnot:BC-example}
	We denote by $ (\ring, \ringzero, \rinvmap) $ and associative involutory set, by $ (\module, q, f) $ a standard pseudo-quadratic module over $ (\ring, \ringzero, \rinvmap) $ and by $ \psgr \defl \psgr(\module) $ the group from \cref{pseud:T-def}. Further, we fix an integer $ n \in \IN_{\ge 2} $, the root system $ \roots \defl BC_n $ in standard representation (as in \cref{BC:BCn-standard-rep}) and we denote by $ \rootbase $ the standard rescaled root base of $ BC_n $.
\end{secnotation}

The goal of this section is the construction of a $ BC_n $-graded group $ \EU(q) $, called the \defemph*{elementary unitary group of $ q $}, from the data in \cref{secnot:BC-example}. This construction is similar to the one of the elementary orthogonal group in \cref{sec:B-example}. As a consequence, both sections follow essentially the same outline, and some remarks remain true almost verbatim. In \cref{BC:ex-constprob-restriction}, we will explain to which extent the construction in this section solves the existence problem for crystallographic $ BC_n $-graded groups. Observe that both $ f $ and $ q $ are allowed to be zero. In this case, $ (\module, q, f) $ is simply an $ \ring $-module with no additional structure, and the group $ \psgr $ equals $ \ringzero \times \module $ by \cref{pseud:T-q-triv}.

Just like the elementary orthogonal orthogonal group, we have implemented the group $ \EU(q) $ in GAP \cite{GAP4} to perform all the computations in this section. See \cref{BC:ex:GAP} for a few details on the implementation.

\begin{convention}
	In this section, we will consider modules over the ring $ \ring $. We use the convention that $ \ring $-scalars act on such modules from the right while endomorphisms act from the left. In particular, the composition $ \phi \circ \psi $ of two such endomorphisms $ \phi, \psi $ is the map $ \map{}{}{}{x}{\phi(\psi(x))} $.
\end{convention}

\begin{note}[Implementation]\label{BC:ex:GAP}
	An implementation of the elementary orthogonal group in the computer algebra system GAP \cite{GAP4} can be found in \cite{RGG-GAP}. This includes proofs of all computational results in this section which are not explicitly verified by hand.
\end{note}

\subsection{Construction}

\begin{construction}
	We put $ V_+ \defl \ring^n $, $ V_- \defl \ring^n $ and $ V \defl \module \dirsum V_+ \dirsum V_- $. We denote by $ (b_1, \ldots, b_n) $ the standard basis of $ V_+ $, by $ (b_{-1}, \ldots, b_{-n}) $ the standard basis of $ V_- $ and we will always consider $ V_+ $ and $ V_- $ to be subsets of $ V $ without specifying the natural embedding. We will usually denote elements of $ V_+, V_- $ by the letters $ v,w $, elements of $ \module $ by the letters $ m,u $ and elements of $ V $ by the letters $ x,y $. We can represent automorphisms of $ V $ using generalised matrices in the same way as in \cref{generalised-matrices}.
\end{construction}

\begin{construction}[Medium-length root homomorphisms]\label{BC:ex-roothom-def-med}
	Let $ i,j \in \numint{1}{n} $ be distinct and let $ a \in \comring $. We denote by $ \risom{\basvec_i - \basvec_j}(a) $ the unique $ \ring $-linear endomorphism of $ V $ which is given on the direct summands of $ V $ as follows, where $ m $, $ v^+ $, $ v^- $ denote arbitrary elements of $ \module $, $ V_+ $, $ V_- $, respectively:
	\[ \risom{\basvec_i - \basvec_j}(a) \colon \quad m \mapsto m, \quad v^+ \mapsto v^+ + b_i \cdot a v^+_j, \quad v^- \mapsto v^- -  b_j \cdot \rinv{a} v^-_i. \]
	Now assume that, in addition, $ i <j $. Then we define $ \risom{\basvec_i + \basvec_j} $ and $ \risom{-\basvec_i - \basvec_j} $ by the following formulas:
	\begin{align*}
		\risom{\basvec_i + \basvec_j}(a) &\colon \quad m \mapsto m, \quad v^+ \mapsto v^+, \quad v^- \mapsto v^- + b_i \cdot a v^-_j + b_j \cdot \rinv{a} v^-_i, \\
		\risom{-\basvec_i - \basvec_j}(a) &\colon \quad m \mapsto m, \quad v^+ \mapsto v^+ + b_{-j} \cdot a v^+_i + b_{-i} \cdot \rinv{a} v^+_j, \quad v^- \mapsto v^- .
	\end{align*}
	Examples of the corresponding generalised matrices are given in \cref{BC:ex-matrices-medium}.
\end{construction}

\begin{figure}[htb]
	\centering$ \begin{aligned}
		\risom{\basvec_1 - \basvec_2}(a) &= \brackets*{\begin{array}{c|ccc|ccc}
			\id_\module &&&&&& \\
			\hline
			& 1 & a &&&& \\
			& & 1 &&&& \\
			& & & 1 \\
			\hline
			& & & & 1 && \\
			& & & & -\rinv{a} & 1 & \\
			& & & & & & 1
		\end{array}}, \\
		\risom{\basvec_1 + \basvec_2}(a) &= \brackets*{\begin{array}{c|ccc|ccc}
			\id_\module &&&&&& \\
			\hline
			& 1 & &  & & a & \\
			& & 1 & & \rinv{a} && \\
			& & & 1 \\
			\hline
			& & & & 1 && \\
			& & & &  & 1 & \\
			& & & & & & 1
		\end{array}}, \\
		\risom{-\basvec_1 - \basvec_2}(a) &= \brackets*{\begin{array}{c|ccc|ccc}
			\id_\module &&&&&& \\
			\hline
			& 1 & &  & &  & \\
			& & 1 & &  && \\
			& & & 1 \\
			\hline
			& & \rinv{a} & & 1 && \\
			& a & & &  & 1 & \\
			& & & & & & 1
		\end{array}}
	\end{aligned} $
	\caption{The medium-length root homomorphisms for $ BC_3 $.}
	\label{BC:ex-matrices-medium}
\end{figure}

\begin{construction}[Short root homomorphisms]\label{BC:ex-roothom-def-short}
	Let $ i \in \numint{1}{n} $ and let $ (u,h) \in \psgr $. Then we define $ \risom{2\basvec_i}(u,h) $ and $ \risom{-2\basvec_i}(u,h) $ to be the unique $ \ring $-linear endomorphisms of $ V $ which are given on the direct summands of $ V $ as follows, where $ m $, $ v^+ $, $ v^- $ denote arbitrary elements of $ \module $, $ V_+ $, $ V_- $, respectively:
	\begin{align*}
		\risom{2\basvec_i}(u,h) &\colon \enskip m \mapsto m - b_i \cdot f(u,m), \enskip v^+ \mapsto v^+, \enskip v^- \mapsto v^- - u \cdot v^-_i + b_i \cdot h v^-_i, \\
		\risom{-2\basvec_i}(u,h) & \colon \enskip m \mapsto m - b_{-i} \cdot f(u, m), \enskip v^+ \mapsto v^+ + u \cdot v^+_i - b_{-i} \cdot h v^+_i, \enskip v^- \mapsto v^-.
	\end{align*}
	Examples of the corresponding generalised matrices are given in \cref{BC:ex-matrices-short}.
\end{construction}

\begin{figure}[htb]
	\centering$ \begin{aligned}
		\risom{\basvec_1}(u, h) &= \brackets*{\begin{array}{c|ccc|ccc}
			\id_\module & && &-u&& \\
			\hline
			-f(u, \mapdot)& 1&& &h&& \\
			&& 1& &&& \\
			& &&1 &&& \\
			\hline
			& && &1&& \\
			& && &&1& \\
			& && &&&1
		\end{array}}, \\
		\risom{-\basvec_1}(u, h) &= \brackets*{\begin{array}{c|ccc|ccc}
			\id_\module & u && &&& \\
			\hline
			& 1&& && \\
			& &1&& && \\
			& &&1& && \\
			\hline
			-f(u, \mapdot)& -h && &1&& \\
			& && &&1& \\
			& && &&&1
		\end{array}}
	\end{aligned} $
	\caption{The short root homomorphisms for $ BC_3 $.}
	\label{BC:ex-matrices-short}
\end{figure}

\begin{lemma}
	For each medium-length root $ \alpha $ and for each short root $ \beta $, the maps $ \map{\risom{\alpha}}{(\ring, +)}{\End_\ring(V)}{}{} $ and $ \map{\risom{\beta}}{\psgr}{\End_\ring(V)}{}{} $ are injective homomorphisms. In particular, their images lie in $ \Aut_\ring(V) $.
\end{lemma}
\begin{proof}
	The injectivity of these maps can be deduced from their matrix representations. The homomorphism property follows from a straightforward matrix computation.
\end{proof}

We can now define our main example of a $ BC_n $-graded group.

\begin{definition}[Elementary unitary group]\label{BC:ex:EUq-def}
	For each root $ \alpha $, we denote the image of $ \risom{\alpha} $ by $ \rootgr{\alpha} $, and we denote by $ \EU(q) $ the group which is generated by $ (\rootgr{\alpha})_{\alpha \in \roots} $. We call $ \EU(q) $ the \defemph*{elementary unitary group of $ q $}\index{elementary unitary group}.
\end{definition}

\begin{warning}
	The group $ \EU(q) $ depends not only on $ q $ but also on the choice of $ f $, which we suppress in our notation.
\end{warning}

\begin{note}
	The same remarks as in \cref{B:ex-is-ortho} hold for the elementary unitary group: A straightforward computation shows that this group is contained in the unitary group $ \Unitary(q) $ from \cref{pseud:unitary-def}, but we will have no need to formally use this fact.
\end{note}

\subsection{Weyl Elements}

Weyl elements for medium-length roots are defined in the same way as usual. However, the formula for short Weyl elements, which is nearly the same as in \cite[(32.9)]{MoufangPolygons}, is more intricate. In particular, observe that short Weyl elements are not necessarily balanced (in the sense of \cref{weyl:balanced}).

\begin{definition}[Weyl elements]\label{BC:ex-weyl-def}
	For all distinct $ i,j \in \numint{1}{n} $ and all invertible $ a \in \ring $, we define
	\[ w_{ij}(a) \defl w_{\basvec_i - \basvec_j}(a) \defl \risom{\basvec_j - \basvec_i}(-a^{-1}) \circ \risom{\basvec_i - \basvec_j}(a) \circ \risom{\basvec_j - \basvec_i}(-a^{-1}) \]
	and $ w_{ij} \defl w_{\basvec_i - \basvec_j} \defl w_{ij}(1_\ring) $. For all $ i \in \numint{1}{n} $ and all $ (u,h) \in \psgr $ for which $ h $ is invertible, we define
	\[ w_i(u,h) \defl w_{2\basvec_i}(u,h) \defl \risom{-2\basvec_i}\brackets[\big]{u h^{-1}, \rinvmin{h}} \circ \risom{2\basvec_i}(u,h) \circ \risom{-2\basvec_i}\brackets[\big]{-u \rinvmin{h}, \rinvmin{h}} \]
	and $ w_i \defl w_{2\basvec_i} \defl w_i(0_\module, 1_\ring) $. Note that $ w_i(u,h) $ is well-defined by \cref{pseud:T-weyl-inverse} and that $ (0_\module, 1_\ring) $ lies in $ \psgr $ because $ (\ring, \ringzero, \rinvmap) $ is an involutory set (and not merely a pre-involutory set).
\end{definition}

\begin{definition}[Standard system of Weyl elements]\label{BC:ex-standard-weyl}
	The \defemph*{standard system of Weyl elements for $ \EU(q) $} is the family $ (w_\delta)_{\delta \in \rootbase} $ given by \cref{BC:ex-weyl-def}.
\end{definition}

Similar remarks as in \cref{B:ex-weyl-Bnsub-note} apply to \cref{BC:ex-weyl-def}: We have only defined Weyl elements for the subset $ \Cnsub $ of $ BC_n $, but this will pose no problem in the sequel.

\begin{remark}[Short Weyl elements, see also~\ref{B:ex:short-weyl-matrix}]\label{BC:ex:short-weyl-matrix}
	Let $ (u,h) \in \psgr $ be such that $ h $ is invertible. The goal of this remark is to explicitly compute the generalised matrix of $ w_1(u,h) $. We will leave out the rows and columns which correspond to $ (b_i)_{i \in \numint{2}{n}} $ and $ (b_{-i})_{i \in \numint{2}{n}} $ because they are trivial. We begin the computation with
	\begin{align*}
		\risom{\basvec_1}(u,h) \risom{-\basvec_1}(-u \rinvmin{h}, \rinvmin{h}) \hspace{-3.5cm} & \\
		&= \begin{pmatrix}
			\id & 0 & -u \\
			-f(u, \mapdot) & 1 & h \\
			0 & 0 & 1
		\end{pmatrix} \begin{pmatrix}
			\id & -u\rinvmin{h} & 0 \\
			0 & 1 & 0 \\
			\rinvmin{h} f(u, \mapdot) & -\rinvmin{h} & 1
		\end{pmatrix} \\
		&= \begin{pmatrix}
			\id - u h^{-1} f(u, \mapdot) & -u\rinvmin{h} + u\rinvmin{h} & -u \\
			-f(u, \mapdot) + hh^{-1} f(u, \mapdot) & f(u,u) \rinvmin{h} + 1 - h\rinvmin{h} & h \\
			h^{-1} f(u, \mapdot) & -\rinvmin{h} & 1
		\end{pmatrix} \\
		&= \begin{pmatrix}
			\id - u h^{-1} f(u, \mapdot) & 0 & -u \\
			0 & 0 & h \\
			h^{-1} f(u, \mapdot) & -\rinvmin{h} & 1
		\end{pmatrix}.
	\end{align*}
	In the last step, we have used that $ f(u,u) = h - \rinv{h} $ because $ h \in q(u) + \ringzero $ and $ (\module, q, f) $ is standard. Now
	\begin{align*}
		w_1(u,h) &= \risom{-\basvec_1}(uh^{-1}, \rinvmin{h}) \risom{\basvec_1}(u,h) \risom{-\basvec_1}(-u \rinvmin{h}, \rinvmin{h}) \\
		&= \begin{pmatrix}
			\id & uh^{-1} & 0 \\
			0 & 1 & 0 \\
			-\rinvmin{h} f(u, \mapdot) & -\rinvmin{h} & 1
		\end{pmatrix} \begin{pmatrix}
			\id - u h^{-1} f(u, \mapdot) & 0 & -u \\
			0 & 0 & h \\
			h^{-1} f(u, \mapdot) & -\rinvmin{h} & 1
		\end{pmatrix} \\
		&= \begin{pmatrix}
			\id - u h^{-1} f(u, \mapdot) & 0 & -u + uh^{-1} h \\
			0 & 0 & h \\
			\phi & -\rinvmin{h} & - \rinvmin{h} f(u,u) - \rinvmin{h}h + 1
		\end{pmatrix}
	\end{align*}
	where $ \map{\phi}{\module}{\ring}{}{} $ is the map which sends $ x \in \module $ to
	\begin{align*}
		&\hspace{-1.5cm}-\rinvmin{h} f\brackets[\big]{u, x - uh^{-1} f(u,x)} + h^{-1} f(u,x) \\
		&= -\rinvmin{h} f(u,x) + \rinvmin{h} f(u,u) h^{-1} f(u,x) + h^{-1} f(u,x) = 0.
	\end{align*}
	We conclude that
	\[ w_1(u,h) = \begin{pmatrix}
		\id - u h^{-1} f(u, \mapdot) & 0 & 0 \\
		0 & 0 & h \\
		0 & -\rinvmin{h} & 0
	\end{pmatrix}. \]
	Hence $ w_1(ua,h) = w_1(u,h) $ where $ a $ is any element of the center of $ \ring $ with $ \rinv{a}a = 1_\ring $ and $ a \rinv{a} = 1_\ring $. (The first equation implies that $ (ua,h) $ lies in $ \psgr $ and the second one implies that $ w_1(ua,h) = w_1(u,h) $.) In particular, $ w_1(u,v) = w_1(-u,v) $ where $ (-u,v) $ is the element obtained from $ (u,v) $ by applying the Jordan module involution. Hence there exist short Weyl elements which can be represented by distinct Weyl triples.
\end{remark}

We now prove that the elements in \cref{BC:ex-weyl-def} are indeed Weyl elements. Since the resulting formulas are rather intricate, we explicitly state the formulas for the standard Weyl elements separately. Only these simplified formulas will be needed to see the twisting actions on the root groups.

\begin{lemma}\label{BC:ex-weylformula-med}
	Let $ i,j \in \numint{1}{n} $ be distinct and let $ a \in \ring $ be invertible. Then $ w \defl w_{ij}(a) $ is an $ (\basvec_i - \basvec_j) $-Weyl element. It satisfies the following formulas for all $ b \in \ring $ and all $ (u,h) \in \psgr $:
	\begin{lemenumerate}
		\item We have
		\begin{gather*}
			\rismin{i}{j}(b)^w = \rismin{j}{i}(-a^{-1} b a^{-1}), \quad \risplus{i}{j}(b)^w = \risplus{i}{j}\brackets[\big]{- \delinv{i>j}\brackets[\big]{a \delinv{i<j}(b) \rinvmin{a}}}, \\
			\rismin{j}{i}(b)^w = \rismin{i}{j}(-aba), \quad \risminmin{i}{j}(a)^w = \risminmin{i}{j}\brackets[\big]{-\delinv{i>j}\brackets[\big]{\rinv{a} \delinv{i<j}(b) a^{-1}}}.
		\end{gather*}
		
		\item For all $ k \in \numint{1}{n} \setminus \Set{i,j} $, we have
		\begin{align*}
			\rismin{i}{k}(b)^w &= \rismin{j}{k}(a^{-1}b), & \risplus{i}{k}(b)^w &= \risplus{j}{k}\brackets[\big]{\delinv{k<j}\brackets[\big]{a^{-1} \delinv{k<i}(b)}}, \\
			\rismin{k}{i}(b)^w &= \rismin{k}{j}(ba), & \risminmin{i}{k}(b)^w &= \risminmin{j}{k}\brackets[\big]{\delinv{j<k}\brackets[\big]{\rinv{a} \delinv{i<k}(b)}}.
		\end{align*}
		
		\item For all $ k \in \numint{1}{n} \setminus \Set{i,j} $, we have
		\begin{align*}
			\rismin{j}{k}(b)^w &= \rismin{i}{k}(-ab), & \risplus{j}{k}(b)^w &= \risplus{i}{k}\brackets[\big]{-\delinv{k<i}\brackets[\big]{a \delinv{k<j}(b)}}, \\
			\rismin{k}{j}(b)^w &= \rismin{k}{i}(-ba^{-1}), & \risminmin{j}{k}(b)^w &= \risminmin{i}{k}\brackets[\big]{-\delinv{i<k}\brackets[\big]{\rinvmin{a} \delinv{j<k}(b)}}.
		\end{align*}
		
		\item \label{BC:ex-weylformula-med:on-short}$ \risshpos{i}(u,h)^w = \risshpos{j}\brackets[\big]{u\rinvmin{a}, a^{-1} h \rinvmin{a}} $ and $ \risshneg{i}(u,h)^w = \risshneg{j}\brackets[\big]{ua, \rinv{a}ha} $.
		
		\item $ \risshpos{j}(u,h)^w = \risshpos{i}\brackets[\big]{-u\rinv{a}, ah\rinv{a}} $ and $ \risshneg{j}(u,h)^w = \risshneg{i}\brackets[\big]{-ua^{-1}, \rinvmin{a} ha^{-1}} $.
		
		\item For any root $ \alpha $ for which no formula for the conjugation action of $ w $ on $ \rootgr{\alpha} $ is given above, this action is trivial.
	\end{lemenumerate}
\end{lemma}
\begin{proof}
	This follows from a straightforward computation. See also \cref{BC:ex:GAP}.
\end{proof}

\begin{lemma}[Corollary of \ref{BC:ex-weylformula-med}]\label{BC:ex-weylformula-med-cor}
	Let $ i,j \in \numint{1}{n} $ be distinct. Then $ w \defl w_{ij} $ is an $ (\basvec_i - \basvec_j) $-Weyl element. It satisfies the following formulas for all $ b \in \ring $ and all $ (u,h) \in \psgr $:
	\begin{lemenumerate}
		\item We have
		\begin{align*}
			\rismin{i}{j}(b)^w &= \rismin{j}{i}(-b), & \risplus{i}{j}(b)^w &= \risplus{i}{j}\brackets[\big]{-\rinv{b}}, \\
			\rismin{j}{i}(b)^w &= \rismin{i}{j}(-b), & \risminmin{i}{j}(a)^w &= \risminmin{i}{j}\brackets[\big]{-\rinv{b}}.
		\end{align*}
		
		\item For all $ k \in \numint{1}{n} \setminus \Set{i,j} $, we have
		\begin{align*}
			\rismin{i}{k}(b)^w &= \rismin{j}{k}(b), & \risplus{i}{k}(b)^w &= \risplus{j}{k}\brackets[\big]{\delinv{k \in \betint{i}{j}}(b)}, \\
			\rismin{k}{i}(b)^w &= \rismin{k}{j}(b), & \risminmin{i}{k}(b)^w &= \risminmin{j}{k}\brackets[\big]{\delinv{k \in \betint{i}{j}}(b)}.
		\end{align*}
		
		\item For all $ k \in \numint{1}{n} \setminus \Set{i,j} $, we have
		\begin{align*}
			\rismin{j}{k}(b)^w &= \rismin{i}{k}(-b), & \risplus{j}{k}(b)^w &= \risplus{i}{k}\brackets[\big]{-\delinv{k \in \betint{i}{j}}(b)}, \\
			\rismin{k}{j}(b)^w &= \rismin{k}{i}(-b), & \risminmin{j}{k}(b)^w &= \risminmin{i}{k}\brackets[\big]{-\delinv{k \in \betint{i}{j}}(b)}.
		\end{align*}
		
		\item \label{BC:ex-weylformula-med-cor:short1}$ \risshpos{i}(u,h)^w = \risshpos{j}\brackets[\big]{u, h} $ and $ \risshneg{i}(u,h)^w = \risshneg{j}\brackets[\big]{u, h} $.
		
		\item \label{BC:ex-weylformula-med-cor:short2}$ \risshpos{j}(u,h)^w = \risshpos{i}\brackets[\big]{-u, h} $ and $ \risshneg{j}(u,h)^w = \risshneg{i}\brackets[\big]{-u, h} $.
		
		\item For any root $ \alpha $ for which no formula for the conjugation action of $ w $ on $ \rootgr{\alpha} $ is given above, this action is trivial.
	\end{lemenumerate}
\end{lemma}
\begin{proof}
	This follows from \cref{BC:ex-weylformula-med} by putting $ a \defl 1_\ring $.
\end{proof}

\begin{lemma}\label{BC:ex-weylformula-short}
	Let $ i \in \numint{1}{n} $ and let $ (v,t) \in \psgr $ be such that $ t $ is invertible. Then $ w \defl w_i(v,t) $ is an $ \basvec_i $-Weyl element. It satisfies the following formulas for all $ b \in \ring $ and all $ (u,h) \in \psgr $:
	\begin{lemenumerate}
		\item For all $ k \in \numint{1}{n} \setminus \compactSet{i} $, we have
		\begin{align*}
			\rismin{i}{k}(b)^w &= \risminmin{i}{k}\brackets[\big]{\delinv{i>k}(\rinv{b} \rinvmin{t})}, & \risplus{i}{k}(b)^w &= \rismin{k}{i}\brackets[\big]{-\delinv{i<k}(b) \rinvmin{t}}, \\
			\rismin{k}{i}(b)^w &= \risplus{k}{i}\brackets[\big]{\delinv{i<k}(bt)}, & \risminmin{i}{k}(b)^w &= \rismin{i}{k}\brackets[\big]{-\rinv{t} \delinv{i<k}(b)}.
		\end{align*}
		
		\item $ \risshpos{i}(u,h)^w = \risshneg{i}\brackets[\big]{u \rinvmin{t} + v \rinvmin{t} f(v,u) \rinvmin{t}, t^{-1} h \rinvmin{t}} $ and \\ $ \risshneg{i}(u,h)^w = \risshpos{i}\brackets[\big]{-ut - v \rinvmin{t} f(v,u) t, \rinv{t} ht} $.
		
		\item $ \risshpos{j}(u,h)^w = \risshpos{j}\brackets[\big]{u + v\rinvmin{t} f(v,u), h} $ and \\ $ \risshneg{j}(u,h)^w = \risshneg{j}\brackets[\big]{u + v \rinvmin{t} f(v,u), h} $.
		
		\item For any root $ \alpha $ for which no formula for the conjugation action of $ w $ on $ \rootgr{\alpha} $ is given above, this action is trivial.
	\end{lemenumerate}
\end{lemma}
\begin{proof}
	This follows from a straightforward computation. See also \cref{BC:ex:GAP}.
\end{proof}

\begin{lemma}[Corollary of \ref{BC:ex-weylformula-short}]\label{BC:ex-weylformula-short-cor}
	Let $ i \in \numint{1}{n} $ . Then $ w \defl w_i $ is an $ \basvec_i $-Weyl element. It satisfies the following formulas for all $ b \in \ring $ and all $ (u,h) \in \psgr $:
	\begin{lemenumerate}
		\item For all $ k \in \numint{1}{n} \setminus \compactSet{i} $, we have
		\begin{align*}
			\rismin{i}{k}(b)^w &= \risminmin{i}{k}\brackets[\big]{\delinv{i<k}(b)}, & \risplus{i}{k}(b)^w &= \rismin{k}{i}\brackets[\big]{-\delinv{i<k}(b)}, \\
			\rismin{k}{i}(b)^w &= \risplus{k}{i}\brackets[\big]{\delinv{i<k}(b)}, & \risminmin{i}{k}(b)^w &= \rismin{i}{k}\brackets[\big]{- \delinv{i<k}(b)}.
		\end{align*}
		
		\item \label{BC:ex-weylformula-short-cor:short-self}$ \risshpos{i}(u,h)^w = \risshneg{i}\brackets[\big]{u, h } $ and $ \risshneg{i}(u,h)^w = \risshpos{i}\brackets[\big]{-u, h} $.
		
		\item For any root $ \alpha $ for which no formula for the conjugation action of $ w $ on $ \rootgr{\alpha} $ is given above, this action is trivial.
	\end{lemenumerate}
\end{lemma}
\begin{proof}
	This follows from \cref{BC:ex-weylformula-short} by putting $ (v,t) \defl (0_\module, 1_\ring) $.
\end{proof}

\begin{remark}[compare \ref{B:ex-wij-inv}]\label{BC:ex-wij-inv}
	Let $ i,j \in \numint{1}{n} $ be distinct and let $ a \in \ring $ be invertible. Since $ w \defl w_{ij}(a) $ is a Weyl element by \cref{BC:ex-weylformula-med}, it follows from \thmitemcref{basic:weyl-general}{basic:weyl-general:minus} that
	\begin{align*}
		w_{ij}(a) &= \rismin{i}{j}(a) \rismin{j}{i}(-a^{-1}) \rismin{j}{i}(-a^{-1})^{w} \\
		&= \rismin{i}{j}(a) \rismin{j}{i}(-a^{-1}) \rismin{i}{j}(a) = w_{ji}(-a^{-1}).
	\end{align*}
	In particular, $ w_{ij} = w_{ji}^{-1} $.
\end{remark}

\begin{remark}\label{BC:ex-jordan-twist}
	Let $ i,j \in \numint{1}{n} $ be distinct and put $ w \defl w_i^2 $. It follows from \thmitemcref{BC:ex-weylformula-short-cor}{BC:ex-weylformula-short-cor:short-self} that $ \risshpos{i}(u,h)^w = \risshpos{i}(-u,h) $ for all $ (u,h) \in \psgr $. That is, $ w $ acts on $ \rootgrshpos{i} $ by the Jordan module involution (see \cref{pseud:T-invo}). In particular, the action of $ w $ on $ \rootgrshpos{i} $ is, in general, neither trivial nor by inversion, so it does not satisfy the square formula for Weyl elements. By \thmitemcref{BC:ex-weylformula-med-cor}{BC:ex-weylformula-med-cor:short1} and~\thmitemref{BC:ex-weylformula-med-cor}{BC:ex-weylformula-med-cor:short2}, $ w_{ij}^2 $ acts on $ \rootgrshpos{i} $ by the Jordan module involution as well.
	
	The observations in this remark also say that the Jordan module involution coincides with the short involution that we will define in~\ref{BC:def:short-invo} on the short root groups of any $ BC_n $-graded group.
\end{remark}

\subsection{Parity Maps and Twisting Structures}

As in \cref{B:ex-twistgroups-def}, the formulas above tell us how to define the twisting groups and parity maps for $ \EU(q) $. They are chosen precisely to make sure that \cref{BC:ex-parmap-describes-conj} holds.

\begin{secnotation}\label{BC:ex-twistgroups-def}
	From now on, we denote by $ (\twistgroup \times \invogroup, \psgr, \ring) $ the standard parameter system (in the sense of \cref{BC:jordanmod:standard-paramsys}) for the Jordan module $ \psgr $ which is constructed from the pseudo-quadratic module $ (\module, q, f) $ as in \cref{BC:jordanmodule-pseudquad-ex}. Further, we define maps
	\[ \map{\inverparsym}{BC_n \times \Cnsub}{\twistgroup}{}{} \midand \map{\invoparsym}{B_n \times \Cnsub}{\invogroup}{}{} \]
	by the formulas in \cref{B:ex-parmap-def}, where $ \Cnsub $ is as in \cref{BC:Cnsub-def}.
	By restricting the second components of $ \inverparsym $ and $ \invoparsym $ to $ \rootbase $, we obtain $ \rootbase $-parity maps which we also denote by $ \inverparsym $ and $ \invoparsym $, and which we call the \defemph*{standard parity maps of type $ BC_n $}\index{parity map!standard!type BCn@type $ BC_n $}.
\end{secnotation}

\premidfigure
\begin{figure}[htb]
	\centering\begin{tabular}[t]{ccc}
		\toprule
		$ \alpha $ & $ \inverpar{\alpha}{\basvec_i - \basvec_j} $ & $ \invopar{\alpha}{\basvec_i - \basvec_j} $ \\
		\midrule
		$ \pm (\basvec_i - \basvec_j) $ & $ (-1_\twistgroup, 1_\twistgroup) $ & $ 1_\invogroup $ \\
		$ \pm (\basvec_i + \basvec_j) $ & $ (-1_\twistgroup, 1_\twistgroup) $ & $ -1_\invogroup $ \\
		$ \pm (\basvec_i - \basvec_k) $ & $ (1_\twistgroup, 1_\twistgroup) $ & $ 1_\invogroup $ \\
		$ \pm (\basvec_i + \basvec_k) $ & $ ( 1_\twistgroup, 1_\twistgroup) $ & $ \delmin{k \in \betint{i}{j}}1_\invogroup $ \\
		$ \pm (\basvec_j - \basvec_k) $ & $ (-1_\twistgroup, 1_\twistgroup) $ & $ 1_\invogroup $ \\
		$ \pm (\basvec_j + \basvec_k) $ & $ (- 1_\twistgroup, 1_\twistgroup) $ & $ \delmin{k \in \betint{i}{j}}1_\invogroup $ \\
		$ \pm \basvec_k \pm \basvec_l $ & $ (1_\twistgroup, 1_\twistgroup) $ & $ 1_\invogroup $ \\
		$ \pm (2)\basvec_i $ & $ (1_\twistgroup, 1_\twistgroup) $ & $ 1_\invogroup $ \\
		$ \pm (2)\basvec_j $ & $ (1_\twistgroup, -1_\twistgroup) $ & $ 1_\invogroup $ \\
		$ \pm (2)\basvec_k $ & $ (1_\twistgroup, 1_\twistgroup) $ & $ 1_\invogroup $ \\
		\bottomrule
	\end{tabular}\qquad\begin{tabular}[t]{ccc}
		\toprule
		$ \alpha $ & $ \inverpar{\alpha}{\basvec_i} $ & $ \invopar{\alpha}{2\basvec_i} $ \\
		\midrule
		$ \basvec_i - \basvec_k $ & $ (1_\twistgroup, 1_\twistgroup) $ & $ \delmin{i<k} 1_\invogroup $ \\
		$ \basvec_k - \basvec_i $ & $ (1_\twistgroup, 1_\twistgroup) $ & $ \delmin{i<k} 1_\invogroup $ \\
		$ \pm (\basvec_i + \basvec_k) $ & $ (-1_\twistgroup, 1_\twistgroup) $ & $ \delmin{i<k} 1_\invogroup $ \\
		$ \pm \basvec_k \pm \basvec_l $ & $ (1_\twistgroup, 1_\twistgroup) $ & $ 1_\invogroup $ \\
		$ (2)\basvec_i $ & $ (1_\twistgroup, 1_\twistgroup) $ & $ 1_\invogroup $ \\
		$ -(2)\basvec_i $ & $ (1_\twistgroup, -1_\twistgroup) $ & $ 1_\invogroup $ \\
		$ \pm (2)\basvec_k $ & $ (1_\twistgroup, 1_\twistgroup) $ & $ 1_\invogroup $ \\
		\bottomrule
	\end{tabular}
	\caption{The definition of $ \inverpar{\alpha}{\beta} $ and $ \invopar{\alpha}{\beta} $ for all $ \alpha \in BC_n $ and $ \beta \in \Cnsub $, see \cref{BC:ex-twistgroups-def}. In the left table, we assume that $ i,j,k,l $ are pairwise distinct. In the right table, we assume that $ i,k,l $ are pairwise distinct. For small values of $ n $, it is of course not possible to choose three or four pairwise distinct indices, in which case the corresponding rows should be ignored.}
	\label{BC:ex-parmap-def}
\end{figure}
\postmidfigure

\begin{note}
	Similar remarks as in \cref{B:ex-parmap-def-Bnsub-note} apply in this situation: The values $ \inverpar{\alpha}{\beta} $ for $ \alpha \in BC_n $ and $ \beta \in \Cnsub \setminus \rootbase $ are not need in the definition of the parity map $ \map{\inverparsym}{BC_n \times \rootbase}{\twistgroup}{}{} $, but they will appear in \cref{BC:ex-parmap-equality}.
\end{note}

\begin{remark}[compare \ref{B:ex-parmap-restrict}]\label{BC:ex-parmap-restrict}
	Consider the subset
	\[ A_{n-1} = \Set{\basvec_i - \basvec_j \given i \ne j \in \numint{1}{n}} \]
	of $ B_n $. Then for all $ \alpha, \beta \in A_{n-1} $, the first component of $ \inverpar{\alpha}{\beta} $ in \cref{B:ex-parmap-def} is the same as the value in \cref{ADE:An-parmap} while the second component is trivial.
\end{remark}

\begin{lemma}\label{BC:ex-parmap-describes-conj}
	Let $ \alpha \in BC_n $ and $ \beta \in \Cnsub $. Let $ x \in \ring $ if $ \alpha $ is medium-length and let $ x \in \psgr $ if $ \alpha $ is short. Then $ \risom{\alpha}(x)^{w_\beta} = \risom{\refl{\beta}(\alpha)}(\inverpar{\alpha}{\beta} \invopar{\alpha}{\beta} . x) $.
\end{lemma}
\begin{proof}
	This follows from \cref{BC:ex-weylformula-med-cor,BC:ex-weylformula-short-cor}.
\end{proof}

\begin{remark}\label{BC:ex-param-choice}
	As in \cref{B:ex-param-choice}, we can perform computations in the group $ \EU(q) $ for some fixed choice of the parameters in \cref{secnot:BC-example} to obtain information about the parity maps $ \inverparsym $ and $ \invoparsym $. In the sequel, we will often choose $ (\ring, \ringzero, \rinvmap) \defl (\IC, \IR, \rinvmap) $ where $ \rinvmap $ is complex conjugation and $ \module \defl \IC $ with $ q $ and $ f $ defined as in \cref{pseud:example:C}. Then the inversion map on $ \ring $ and the ring involution are non-trivial and distinct. Further, the inversion map on $ \psgr $ and the Jordan module involution on $ \psgr $ are non-trivial and distinct. These observations imply that the parameter system $ (\twistgroup \times \invogroup, \psgr, \ring) $ is $ (\inverparsym \times \invoparsym) $-faithful (in the sense of \cref{param:parsys-faithful}).
\end{remark}

\begin{lemma}\label{BC:ex-param-thm}
	The root isomorphisms $ (\risom{\alpha})_{\alpha \in BC_n} $ from \cref{B:ex-roothom-def-long,B:ex-roothom-def-short} form a parametrisation of $ \EU(q) $ by $ (\twistgroup \times \invogroup, \psgr, \ring) $ with respect to $ \inverparsym \times \invoparsym $ and $ (w_\delta)_{\delta \in \rootbase} $. 
\end{lemma}
\begin{proof}
	This follows from \cref{BC:ex-parmap-describes-conj}.
\end{proof}

\begin{lemma}\label{BC:ex-weyl-extension}
	Let $ (w_\alpha)_{\alpha \in \Cnsub} $ be as in \cref{B:ex-weyl-const}. Then the following hold:
	\begin{lemenumerate}
		\item $ w_{ij}^{w_{jk}} = w_{ik} $ and $ w_{ij}^{w_{kj}} = w_{ki} $ for all pairwise distinct $ i,j,k \in \numint{1}{n} $, and $ w_i^{w_{ij}} = w_j $ and $ w_i^{w_{ji}} = w_i^{-1} $ for all distinct $ i,j \in \numint{1}{n} $.
		
		\item $ (w_\alpha)_{\alpha \in \Cnsub} $ is the standard $ \Cnsub $-extension of $ (w_\delta)_{\delta \in \rootbase} $.
	\end{lemenumerate}
\end{lemma}
\begin{proof}
	The first assertion follows from \cref{BC:ex-parmap-describes-conj} and an inspection of \cref{BC:ex-parmap-def}. The second assertion follows from the first one and \cref{BC:ex-wij-inv}.
\end{proof}

\subsection{Commutator Relations}

We now proceed to show that $ (\rootgr{\alpha})_{\alpha \in BC_n} $ is a crystallographic $ BC_n $-grading of $ \EU(q) $.

\begin{proposition}\label{BC:ex-commrel}
	The group $ \EU(q) $ satisfies the following commutator relations. For all pairwise distinct $ i,j,k \in \numint{1}{n} $ and all $ a,b \in \ring $, we have
	\begin{align*}
		\commutator{\rismin{i}{j}(a)}{\rismin{j}{k}(b)} &= \rismin{i}{k}(ab), \\
		\commutator{\rismin{i}{j}(a)}{\risplus{j}{k}(b)} &= \risplus{i}{k}\brackets[\big]{\delinv{k<i}\brackets{a \delinv{k<j}(b)}}, \\
		\commutator{\rismin{i}{j}(a)}{\risplus{i}{j}(b)} &= \risshpos{i}\brackets[\big]{0_\module, \ringTr\brackets{a, \delinv{i>j}(b)}}, \\
		\commutator{\rismin{i}{j}(a)}{\risminmin{i}{j}(b)} &= \risshneg{j}\brackets[\big]{0_\module, \ringTr\brackets{\rinv{a}, \delinv{j<i}(b)}}, \\
		\commutator{\rismin{i}{j}(a)}{\risminmin{i}{k}(b)} &= \risminmin{j}{k}\brackets[\big]{-\delinv{j<k}\brackets{\rinv{a} \delinv{i<k}(b)}}, \\
		\commutator{\risplus{i}{j}(a)}{\risminmin{j}{k}(b)} &= \rismin{i}{k}\brackets[\big]{\delinv{j<i}(a) \delinv{j<k}(b)}
	\end{align*}
	where $ \ringTr(a,b) \defl a \rinv{b} + b \rinv{a} $. For all distinct $ i,j \in \numint{1}{n} $ and all $ (u,h), (v,k) \in \psgr(q) $, we have
	\begin{align*}
		\commutator{\risshpos{i}(u,h)}{\risshpos{j}(v,k)} &= \begin{cases}
			\risplus{i}{j}\brackets[\big]{f(u,v)} & \text{if } i<j, \\
			\risplus{i}{j}\brackets[\big]{-f(v,u)} & \text{if } i>j,
		\end{cases} \\
		\commutator{\risshpos{i}(u,h)}{\risshneg{j}(v,k)} &= \rismin{i}{j}\brackets[\big]{-f(u,v)}, \\
		\commutator{\risshneg{i}(u,h)}{\risshneg{j}(v,k)} &= \begin{cases}
			\risminmin{i}{j}\brackets[\big]{f(v,u)} & \text{if } i<j, \\
			\risminmin{i}{j}\brackets[\big]{-f(u,v)} & \text{if } i>j.
		\end{cases}
	\end{align*}
	For all distinct $ i,j \in \numint{1}{n} $, all $ a \in \ring $ and all $ (u,h) \in \psgr(q) $, we have
	\begin{align*}
		\commutator{\risshpos{j}(u,h)}{\rismin{i}{j}(a)} &= \risshpos{i}\brackets[\big]{-u \rinv{a}, a h \rinv{a}} \risplus{i}{j}\brackets[\big]{-\delinv{i>j}(ah)}, \\
		\commutator{\risshneg{i}(u,h)}{\rismin{i}{j}(a)} &= \risshneg{j}\brackets[\big]{ua, \rinv{a} h a} \risminmin{i}{j}\brackets[\big]{\delinv{i>j}(\rinv{a} h)}, \\
		\commutator{\risshneg{i}(u,h)}{\risplus{i}{j}(a)} &= \risshpos{j}\brackets[\big]{-u \delinv{i>j}(a), \delinv{i<j}(a) h \delinv{i>j}(a)} \rismin{j}{i}\brackets[\big]{\delinv{i<j}(a) h}, \\
		\commutator{\risshpos{i}(u,h)}{\risminmin{i}{j}(a)} &= \risshneg{j}\brackets[\big]{-u \delinv{i<j}(a), \delinv{i>j}(a) h \delinv{i<j}(a)} \rismin{i}{j}\brackets[\big]{\rinv{h} \delinv{i<j}(a)}.
	\end{align*}
\end{proposition}
\begin{proof}
	This follows from a straightforward but lengthy computation. See also \cref{BC:ex:GAP}.
\end{proof}

\begin{proposition}\label{BC:ex-triangular}
	Let $ \possys $ denote the standard positive system in $ BC_n $. Then the group $ \EU(q) $ satisfies $ \rootgr{\possys} \intersect \rootgr{-\possys} = \compactSet{1} $.
\end{proposition}
\begin{proof}
	We consider the same rearrangement of the decomposition of $ V $ as in \cref{B:ex-triangular}. With respect to this rearrangement, the group $ \rootgr{\possys} $ consists of upper triangular matrices while $ \rootgr{-\possys} $ consists of lower triangular matrices. The assertion follows.
\end{proof}

\begin{theorem}\label{BC:ex-is-rgg}
	The family $ (\rootgr{\alpha})_{\alpha \in BC_n} $ is a crystallographic $ BC_n $-grading of $ \EU(q) $.
\end{theorem}
\begin{proof}
	By \cref{BC:ex-weylformula-short,BC:ex-weylformula-med}, there exist Weyl elements for all roots in the standard rescaled root base of $ BC_n $. By \cref{basic:weyl-ex-basis:BC}, it follows that there exist Weyl elements for all roots. The remaining conditions are satisfied by \cref{BC:ex-commrel,BC:ex-triangular}.
\end{proof}

\subsection{Concluding Remarks}

\begin{remark}\label{BC:ex:C}
	It follows from the commutator relations in \cref{BC:ex-commrel} that $ \EU(q) $ can be regarded as a crystallographic $ C_n $-graded group (in the sense of \cref{BC:CisBC:BCasC}) if and only if $ f=0 $. If $ 2_\ring $ is invertible, this property implies that $ q $ is trivial by \thmitemcref{pseud:2}{pseud:2:inv-triv}. In the case that both $ f $ and $ q $ are zero, our construction produces a crystallographic $ C_n $-graded group whose long root groups are parametrised by $ \ringzero \times \module $ where $ \module $ is an $ \ring $-module with no additional structure.
\end{remark}

\begin{note}\label{BC:ex-B-counterex-weyl}
	Assume that $ \module \ne \compactSet{0} $, $ q = 0 $, $ f=0 $ and $ n=2 $, so that $ \EU(q) $ is a crystallographic $ C_2 $-graded group by \cref{BC:ex:C}. Then $ \psgr = \module \times \ringzero $ by \cref{pseud:T-q-triv}. We have seen in \cref{BC:ex-jordan-twist} that $ w \defl w_1^2 $ and $ v \defl w_{21}^2 $ act on $ \rootgrshpos{1} $ by the Jordan module involution on $ \psgr $. Since $ \module \ne \compactSet{0} $, this involution is neither the identity nor the inversion map.
	
	Now put $ \alpha \defl \basvec_1 $ and $ \delta \defl \basvec_2 - \basvec_1 $ and regard $ \EU(q) $ as a crystallographic $ B_2 $-graded group, as in \cref{BC:CisBC:C2}. Then $ (\alpha, \delta) $ is a $ B_2 $-pair in the $ B_2 $-grading of $ \EU(q) $. The observation in the previous paragraph implies that $ \EU(q) $ is a crystallographic $ B_2 $-graded group in which the squares of $ \alpha $-Weyl elements and $ \delta $-Weyl do not act trivially on $ \rootgr{\alpha} $. This is the counterexample that was promised in \cref{B:rank2-note}.
\end{note}

The following result is the analogue of \cref{B:ex-parmap-equality} for the root system $ BC $. It will be used in the proof of \cref{BC:Cnsub-conj-anygroup}

\begin{lemma}\label{BC:ex-parmap-equality}
	Let $ \alpha \in \Cnsub $, let $ \word{\alpha} $ be the standard $ \rootbase $-expression of $ \alpha $ in the sense of \cref{BC:Cnsub-def} and let $ \zeta $ be an arbitrary root. Then the element $ \inverpar{\zeta}{\alpha} $ from \cref{BC:ex-parmap-def} equals the value $ \inverpar{\zeta}{\word{\alpha}} $ of the extended parity map $ \map{\inverparsym}{\roots \times \frmon{\rootbase \union (-\rootbase)}}{\twistgroup}{}{} $ (from \cref{param:parmap-def}). The same assertion holds for $ \invoparsym $ in place of $ \inverparsym $.
\end{lemma}
\begin{proof}
	This can be proven in exactly the same way as \cref{B:ex-parmap-equality}, using the choice of parameters in \cref{BC:ex-param-choice}.
\end{proof}

\begin{note}[The existence problem]\label{BC:ex-constprob-restriction}
	For a complete solution of the existence problem for $ BC_n $-graded groups, we would have to construct a $ BC_n $-graded group from an arbitrary Jordan module $ \psgr $ over an alternative ring $ \ring $ with nuclear involution $ \rinvmap $. This section covers exactly the special case that $ \ring $ is associative and that $ \psgr $ is of pseudo-quadratic type in the sense of \cref{BC:jordanmodule-pseudquad-ex}.
	
	The second restriction is relatively minor: We know from \cref{BC:jordanmodule:class-2inv} that every Jordan module is of pseudo-quadratic type if $ 2 $ is invertible in the base ring. Further, we are not aware of an example of a Jordan module which is not of pseudo-quadratic type.
	
	The restriction on $ \ring $ to be associative is more intricate. Recall from \cref{pure-alt:weak-thm} that a purely alternative ring $ \ring $ does not admit a pseudo-quadratic module $ (\module, q, f) $ with $ f \ne 0 $. Thus it seems reasonable to exclude alternative rings which are not associative from the discussion if $ f \ne 0 $. By \cref{BC:ex:C}, the assumption that $ f \ne 0 $ says precisely that $ (\rootgr{\alpha})_{\alpha \in BC_n} $ is a \emph{proper} crystallographic $ BC_n $-grading in the sense of \cref{BC:CisBC:BCasC}. Equivalently, it says that the Jordan module $ \psgr(\module) $ is not of type $ C $.
	
	We conclude that the construction in this section is close to a complete solution of the existence problem for proper crystallographic $ BC_n $-gradings. However, the assumption that $ \ring $ is associative is an important restriction in the setting of crystallographic $ C_n $-gradings. In \cite[Section~4.2]{Zhang}, Zhang gives a construction of crystallographic $ C_n $-graded groups for alternative rings $ \ring $ with nuclear involution in which $ 2_\ring $ is invertible. This construction builds on \cite[21.12]{LoosNeherBook} and the standard Jordan matrix algebras in \cite{McCrimmon66}. The Jordan module which coordinatises the long root groups in this example is always $ \ringTr_\rinvmap(\ring) = \symring $. In other words, this construction covers precisely the Jordan modules in \cref{BC:jordanmodule-invset-ex} for which $ 2_\ring $ is invertible. We are confident that a generalisation of Zhang's strategy can be used to solve the existence problem for arbitrary Jordan modules of type~$ C $. We plan to address this problem in future work.
\end{note}


\section{Rank-\texorpdfstring{$ 2 $}{2} Computations}

\label{sec:BC:rank2}

\begin{secnotation}
	We denote by $ G $ a group which has crystallographic $ BC_2 $-commutator relations with root groups $ (\rootgr{\alpha})_{\alpha \in BC_2} $ (in the sense of \cref{rgg:group-commrel-def}) and we choose a $ BC_2 $-quadruple $ (\alpha, \beta, \gamma, \delta) $. Further, we assume that $ G $ is rank-2-injective.
\end{secnotation}

This section is structured similarly to \cref{sec:B:rank2-cry}. Our goal is to understand the action of squares of Weyl elements on all root groups. In most cases, we will see that this action obeys the square formula. However, we have already seen in the example of the elementary unitary group that this is not always the case (see \cref{BC:ex-jordan-twist}).

Recall from \cref{BC:BisBC:BC-noncry} that $ (\rootgr{\alpha})_{\alpha \in B_2} $ is a (non-crystallographic) $ B_2 $-grading of $ G $. Thus we can freely apply the results from \cref{sec:B:rank2-noncry}. By \cref{B2-is-C2}, we can also apply all results from \cref{sec:B:rank2-cry} if we reverse the roles of $ \alpha $ and $ 2\delta $ (and, consequently, of $ 2\beta $ and $ \gamma $).

\begin{reminder}\label{B2-is-C2}
	Recall from \cref{BC:CisBC:CinBC} that $ (\rootgr{\alpha})_{\alpha \in C_2} $ is a crystallographic $ C_2 $-grading of the subgroup of $ G $ that it generates. Thus it is also a crystallographic $ B_2 $-grading by \cref{BC:CisBC:C2}. In this $ B_2 $-grading, $ (2\delta, \gamma, 2\beta, \alpha) $ is a $ B_2 $-quadruple.
\end{reminder}

In addition to \cref{B:comm-add}, we will frequently use the following computation rule.

\begin{lemma}\label{BC:comm-add}
	The following relations hold for all $ x_\alpha \in \rootgr{\alpha} $ and $ y_{2\delta}, y_{2\delta}' \in \rootgr{2\delta} $:
	\begin{lemenumerate}
		\item \label{BC:comm-add:delta-gamma}$ \commpart{x_\alpha}{y_{2\delta} y_{2\delta}'}{\gamma} = \commpart{x_\alpha}{y_{2\delta}'}{\gamma} \commpart{x_\alpha}{y_{2\delta}}{\gamma} $.
		
		\item \label{BC:comm-add:delta-gamma2}$ \commpart{y_{2\delta} y_{2\delta}'}{x_\alpha}{\gamma} = \commpart{y_{2\delta}}{x_\alpha}{\gamma} \commpart{y_{2\delta}'}{x_\alpha}{\gamma} $.
	\end{lemenumerate}
\end{lemma}
\begin{proof}
	By \cref{B:comm-add}, we have the following relations:
	\begin{align*}
		\commpart{x_\alpha}{y_{2\delta} y_{2\delta}'}{\gamma} &= \commpart{x_\alpha}{y_{2\delta}'}{\gamma} \commutator[\big]{\commpart{x_\alpha}{y_{2\delta}}{\beta}}{y_{2\delta}'} \commpart{x_\alpha}{y_{2\delta}}{\gamma}, \\
		\commpart{y_{2\delta} y_{2\delta}'}{x_\alpha}{\gamma} &= \commpart{y_{2\delta}}{x_\alpha}{\gamma} \commutator[\big]{\commpart{y_{2\delta}}{x_\alpha}{\beta}}{y_{2\delta}'} \commpart{y_{2\delta}'}{x_\alpha}{\gamma}.
	\end{align*}
	Since $ y_{2\delta}' $ commutes with $ \rootgr{\beta} $, the assertions follow.
\end{proof}

\begin{note}\label{BC:square-act:length-restrict-note}
	We will only investigate the actions of squares of $ \zeta $-Weyl elements on root groups $ \rootgr{\rho} $ when $ \rho $ is not long and $ \zeta $ is not short. We do not have to consider long roots $ \rho $ because the short root groups are contained in the long root groups. Thus if we know the action of $ w_\zeta^2 $ on all short and medium-length root groups, then we know its action on all root groups. Further, we do not consider short roots $ \zeta $ by the observations in \cref{param:motiv:rescaled}.
\end{note}

\subsection{The Action of Long Weyl Elements on Medium-length Root Groups}

All results in this subsection follow from the corresponding ones in \cref{subsec:B:rank2:long-on-short} using \cref{B2-is-C2}. Thus we obtain that squares of long Weyl elements act on the medium-length root groups in the same $ BC_2 $-subsystem by inversion.

\begin{lemma}\label{C:basecomp-wdelta-alpha-cor}
	Assume that $ (a_{-2\delta}, b_{2\delta}, c_{-2\delta}) $ is a $ 2\delta $-Weyl triple and denote by $ w_{2\delta} \defl a_{-2\delta} b_{2\delta} c_{-2\delta} $ the corresponding Weyl element. Then the following statements hold for all $ x_\alpha \in \rootgr{\alpha} $:
	\begin{lemenumerate}
		\item \label{C:basecomp-wdelta-alpha-cor:conjformula}$ x_\alpha^{w_{2\delta}} = \commpart{x_\alpha}{b_{2\delta}}{\gamma} $.
		
		\item \label{C:basecomp-wdelta-alpha-cor:1G}$ \commutator[\big]{x_\alpha}{\commpart{x_\alpha}{b_{2\delta}}{\gamma}} \commpart{x_\alpha}{b_{2\delta}}{\beta} \commpart[\big]{\commpart{x_\alpha}{b_{2\delta}}{\gamma}}{c_{-2\delta}}{\beta}^{\commpart{x_\alpha}{b_{2\delta}}{\beta}} = 1_G $.
		
		\item $ x_\alpha = \commpart[\big]{\commpart{x_\alpha}{b_{2\delta}}{\gamma}}{c_{-2\delta}}{\alpha}^{-1} $.
	\end{lemenumerate}
\end{lemma}
\begin{proof}
	Since $ \rootgr{\beta} $ and $ \rootgr{2\delta} $ commute in a crystallographic $ B_2 $-grading, these statements follow from \cref{B:basecomp-wdelta-cox-alpha}. Alternatively, this follows from \cref{B:basecomp-walpha-cor:delta} by the considerations in \cref{B2-is-C2}.
\end{proof}

\begin{lemma}\label{C:basecomp-wdelta-gamma-cor}
	Assume that $ (a_{-2\delta}, b_{2\delta}, c_{-2\delta}) $ is a $ 2\delta $-Weyl triple and denote by $ w_{2\delta} \defl a_{-2\delta} b_{2\delta} c_{-2\delta} $ the corresponding Weyl element. Then the following statements hold for all $ x_\gamma \in \rootgr{\gamma} $:
	\begin{lemenumerate}
		\item \label{C:basecomp-wdelta-gamma-cor:conjformula}$ x_\gamma^{w_{2\delta}} = \commpart{x_\gamma}{a_{-2\delta}}{\alpha} = \commpart{x_\gamma}{c_{-2\delta}}{\alpha} $.
		
		\item $ \commutator[\big]{x_\gamma}{\commpart{x_\gamma}{a_{-2\delta}}{\alpha}} \commpart{x_\gamma}{a_{-2\delta}}{\beta} \commpart[\big]{\commpart{x_\gamma}{a_{-2\delta}}{\alpha}}{b_{2\delta}}{\beta}^{\commpart{x_\gamma}{a_{-2\delta}}{\beta}} = 1_G $.
		
		\item $ x_\gamma = \commpart[\big]{\commpart{x_\gamma}{a_{-2\delta}}{\alpha}}{b_{2\delta}}{\gamma}^{-1} $.
	\end{lemenumerate}
\end{lemma}
\begin{proof}
	Since $ \rootgr{\beta} $ and $ \rootgr{2\delta} $ commute in a crystallographic $ B_2 $-grading, these statements follow from \cref{B:basecomp-wdelta-cox-gamma}. Alternatively, this follows from \cref{B:basecomp-walpha-cor:beta} by the considerations in \cref{B2-is-C2}.
\end{proof}

\begin{lemma}\label{BC:weyldelta-square-on-alpha}
	Let $ w_{2\delta} $ be a $ 2\delta $-Weyl element. Then $ w_{2\delta}^2 $ acts on $ \rootgr{\alpha} $ and $ \rootgr{\gamma} $ by inversion.
\end{lemma}
\begin{proof}
	Using \cref{B2-is-C2}, this follows from \cref{B:weylalpha-square-on-beta}. Alternatively, we can mimic the proof of \cref{B:weylalpha-square-on-beta} with references to \cref{B:basecomp-walpha-cor:delta,B:basecomp-walpha-cor:beta} replaced by references to \cref{C:basecomp-wdelta-alpha-cor,C:basecomp-wdelta-gamma-cor}.
\end{proof}

\begin{proposition}\label{C:short-abelian}
	Assume that $ \invset{2\delta} $ is non-empty. Then $ \rootgr{\alpha} $ and $ \rootgr{\gamma} $ are abelian.
\end{proposition}
\begin{proof}
	This follows from \cref{B:short-abelian} using \cref{B2-is-C2}, or alternatively by mimicing the proof of \cref{B:short-abelian} with references to \cref{B:weylalpha-square-on-beta} replaced by references to \cref{BC:weyldelta-square-on-alpha}.
\end{proof}

\begin{note}
	In the case of rank at least~$ 3 $, the result of \cref{C:short-abelian} is also a consequence of \cref{A2Weyl:rootisom}.
\end{note}

\subsection{The Action of Long Weyl Elements on Short Root Groups}

\begin{lemma}\label{BC:square-act:beta-on-delta}
	Let $ w_{2\beta} $ be a $ 2\beta $-Weyl element. Then $ w_{2\beta}^2 $ acts trivially on $ \rootgr{\delta} $.
\end{lemma}
\begin{proof}
	This is trivial because $ \delta $ is crystallographically adjacent to $ 2\beta $ and $ -2\beta $.
\end{proof}

We have seen in \cref{BC:ex-jordan-twist} that in the elementary unitary group, the square of a certain \enquote{standard} $ 2\beta $-Weyl element acts on $ \rootgr{\beta} $ by the Jordan module involution. The following formula provides a description of this involution in the general situation. Note that it holds not only for some specific Weyl element, but in fact for all $ 2\beta $-Weyl elements.

\begin{lemma}\label{BC:square-act:beta-on-beta}
	Let $ w_{2\beta} $ be a $ 2\beta $-Weyl element and put $ w \defl w_{2\beta}^2 $. Then $ \commpart{x_\delta}{x_\alpha}{\beta}^w = \commpart{x_\delta}{x_\alpha^{-1}}{\beta} $ for all $ x_\alpha \in \rootgr{\alpha} $ and $ x_\delta \in \rootgr{\delta} $. In particular, the action of $ w $ on $ \rootgr{\beta} $ does not depend on the choice of $ w_{2\beta} $ if $ \invset{\alpha} $ is non-empty.
\end{lemma}
\begin{proof}
	We know that $ w $ acts trivially on $ \rootgr{\delta} $ (by \cref{BC:square-act:beta-on-delta}) and by inversion on $ \rootgr{\alpha} $ (by \cref{BC:weyldelta-square-on-alpha}). Thus by an application of \cref{basic:commpart-conj}, we see that the action of $ w $ on $ \rootgr{\beta} $ satisfies the first assertion.
	
	Now assume, in addition, that $ \invset{\alpha} $ is non-empty. Then every $ x_\beta \in \rootgr{\beta} $ can be written as $ x_\beta = \commpart{x_\delta}{x_\alpha}{\beta} $ for some $ x_\delta \in \rootgr{\delta} $ and $ x_\alpha \in \invset{\alpha} $ by \cref{B:rootisom-beta-delta}. It follows that the action of $ w $ on $ \rootgr{\beta} $ is completely determined by the formula in the first assertion. In particular, it does not depend on the choice of $ w_{2\beta} $.
\end{proof}

\subsection{The Action of Medium-length Weyl Elements on Short Root Groups}

We begin with a repetition of some of the formulas in \cref{subsec:B:rank2:short-on-long}.

\begin{lemma}\label{C:basecomp-walpha-cor-delta}
	Assume that $ (a_{-\alpha}, b_\alpha, c_{-\alpha}) $ is an $ \alpha $-Weyl triple and denote by $ w_\alpha \defl a_{-\alpha} b_\alpha c_{-\alpha} $ the corresponding Weyl element. Then the following statements hold for all $ x_\delta \in \rootgr{\delta} $:
	\begin{lemenumerate}
		\item \label{C:basecomp-walpha-cor-delta:1}$ x_\delta^{w_\alpha} = \commpart{x_\delta}{b_\alpha}{\beta} $.
		
		\item $ \commpart{x_\delta}{b_\alpha}{\gamma}^{-1} = \commpart[\big]{\commpart{x_\delta}{b_\alpha}{\beta}}{c_{-\alpha}}{\gamma}^{\commpart{x_\delta}{b_\alpha}{\gamma}} $.
		
		\item $ x_\delta \commpart[\big]{\commpart{x_\delta}{b_\alpha}{\beta}}{c_{-\alpha}}{\delta} \commutator{\commpart{x_\delta}{b_\alpha}{\gamma}}{c_{-\alpha}} = 1_G $.
	\end{lemenumerate}
\end{lemma}
\begin{proof}
	This follows from \cref{B:basecomp-walpha-cox-cor-delta} with the additional information that $ \rootgr{\beta} $ and $ \rootgr{\delta} $ commute. Alternatively, this follows from \cref{B:basecomp-wdelta-alpha} using the considerations in \cref{B2-is-C2}.
\end{proof}

\begin{lemma}\label{C:basecomp-walpha-cor-beta}
	Assume that $ (a_{-\alpha}, b_\alpha, c_{-\alpha}) $ is an $ \alpha $-Weyl triple and denote by $ w_\alpha \defl a_{-\alpha} b_\alpha c_{-\alpha} $ the corresponding Weyl element. Then the following statements hold for all $ x_\beta \in \rootgr{\beta} $:
	\begin{lemenumerate}
		\item \label{C:basecomp-walpha-cor-beta:1}$ x_\beta^{w_\alpha} = \commpart{x_\beta}{a_{-\alpha}}{\delta} = \commpart{x_\beta}{c_{-\alpha}}{\delta} $.
		
		\item $ \commutator[\big]{x_\beta}{\commpart{x_\beta}{a_{-\alpha}}{\delta}} \commpart{x_\beta}{a_{-\alpha}}{\gamma} \commpart[\big]{\commpart{x_\beta}{a_{-\alpha}}{\delta}}{b_\alpha}{\gamma}^{\commpart{x_\beta}{a_{-\alpha}}{\gamma}} = 1_G $.
		
		\item $ x_\beta \commpart[\big]{\commpart{x_\beta}{a_{-\alpha}}{\delta}}{b_\alpha}{\beta} \commutator{\commpart{x_\beta}{a_{-\alpha}}{\gamma}}{b_\alpha} = 1_G $.
	\end{lemenumerate}
\end{lemma}
\begin{proof}
	This follows from \cref{B:basecomp-walpha-cox-cor-beta} with the additional information that $ \rootgr{\beta} $ and $ \rootgr{\delta} $ commute. Alternatively, this follows from \cref{B:basecomp-wdelta-gamma} using the considerations in \cref{B2-is-C2}.
\end{proof}

We can prove the same formula as in \cref{BC:square-act:beta-on-beta} for the actions of $ w_\alpha^2 $ on $ \rootgr{\beta} $, but with a caveat: A priori, it is not clear that the element $ x_\alpha \in \rootgr{\alpha} $ in this formula can be chosen arbitrarily. This will only follow from an application of \cref{BC:square-act:beta-on-beta}.

\begin{lemma}\label{BC:square-act:alpha-beta-on-beta}
	Let $ w_\alpha $ be an $ \alpha $-Weyl element and let $ w_{2\beta} $ be a $ 2\beta $-Weyl element. Then the actions of $ w_\alpha^2 $ and $ w_{2\beta}^2 $ on $ \rootgr{\beta} $ are identical. More precisely, we have $ \commpart{x_\delta}{x_\alpha}{\beta}^w = \commpart{x_\delta}{x_\alpha^{-1}}{\beta} $ for all $ x_\alpha \in \rootgr{\alpha} $, $ x_\delta \in \rootgr{\delta} $ and $ w \in \Set{w_\alpha^2, w_{2\beta}^2} $.
\end{lemma}
\begin{proof}
	Choose $ a_{-\alpha}, c_{-\alpha} \in \rootgr{-\alpha} $ and $ b_\alpha \in \rootgr{\alpha} $ such that $ w_\alpha = a_{-\alpha} b_\alpha c_{-\alpha} $. Then by \cref{BC:square-act:beta-on-beta,B:walpha-cox-on-beta-comm}, both $ w_\alpha^2 $ and $ w_{2\beta}^2 $ map $ \commpart{x_\delta}{b_\alpha^{-1}}{\beta} $ to $ \commpart{x_\delta}{b_\alpha}{\beta} $ for all $ x_\delta \in \rootgr{\delta} $. Since $ b_\alpha^{-1} $ lies in $ \invset{\alpha} $ by \thmitemcref{basic:weyl-general}{basic:weyl-general:inv}, we know from \cref{B:rootisom-beta-delta} that any $ x_\beta \in \rootgr{\beta} $ can be written as $ x_\beta = \commpart{x_\delta}{b_\alpha^{-1}}{\beta} $ for some $ x_\delta \in \rootgr{\delta} $. Thus the actions of $ w_\alpha^2 $ and $ w_{2\beta}^2 $ are completely determined by this property. In particular, these actions are the same, so it follows from \cref{BC:square-act:beta-on-beta} that they map $ \commpart{x_\delta}{x_\alpha}{x_\beta} $ to $ \commpart{x_\delta}{x_\alpha^{-1}}{\beta} $ for all $ x_\alpha \in \rootgr{\alpha} $ and all $ x_\delta \in \rootgr{\delta} $. This finishes the proof.
\end{proof}

\begin{lemma}\label{BC:alpha-beta-on-beta:4-power}
	Let $ w_\alpha $ be an $ \alpha $-Weyl element and let $ w_\beta $ be a $ \beta $-Weyl element. Then $ w_\alpha^4 $ and $ w_\beta^4 $ act trivially on $ \rootgr{\beta} $.
\end{lemma}
\begin{proof}
	This is a consequence of \cref{BC:square-act:alpha-beta-on-beta}.
\end{proof}

Using the results of this subsection, we can show that all long root groups are abelian if sufficiently many medium-length Weyl elements exist.

\begin{lemma}\label{BC:long-inv-gen}
	Let $ b_\alpha \in \invset{\alpha} $. Then $ \rootgr{2\beta} $ is generated by $ \rootgr{\gamma} \union \rootgr{2\delta} \union \compactSet{b_\alpha} $.
\end{lemma}
\begin{proof}
	Let $ x_{2\beta} \in \rootgr{2\beta} $ and let $ a_{-\alpha}, c_{-\alpha} \in \rootgr{-\alpha} $ such that $ w_\alpha \defl a_{-\alpha} b_\alpha c_{-\alpha} $ is an $ \alpha $-Weyl element. Put $ x_{2\delta} \defl x_{2\beta}^{w_\alpha} \in \rootgr{2\delta} $. Then by \thmitemcref{C:basecomp-walpha-cor-delta}{C:basecomp-walpha-cor-delta:1},
	\[ x_{2\beta} = x_{2\delta}^{w_\alpha^{-1}} = \commpart{x_{2\delta}}{b_\alpha^{-1}}{\beta}. \]
	Thus $ \commutator{x_{2\delta}}{b_\alpha^{-1}} = x_{2\beta} \commpart{x_{2\delta}}{b_\alpha^{-1}}{\gamma} $. It follows that
	\[ x_{2\beta} = \commutator{x_{2\delta}}{b_\alpha^{-1}} \commpart{x_{2\delta}}{b_\alpha^{-1}}{\gamma}^{-1} \in \gen{b_\alpha, \rootgr{\gamma}, \rootgr{2\delta}}. \]
	This finishes the proof.
\end{proof}

\begin{lemma}\label{BC:long-abelian-lem}
	If $ \invset{\alpha} $ is not empty, then $ \rootgr{2\beta} $ is abelian.
\end{lemma}
\begin{proof}
	Choose an arbitrary element $ b_\alpha \in \invset{\alpha} $. By the crystallographic commutator relations, $ \rootgr{2\beta} $ commutes with $ \rootgr{\gamma} \union \rootgr{2\delta} \union \compactSet{b_\alpha} $. Hence it follows from \cref{BC:long-inv-gen} that $ \rootgr{2\beta} $ commutes with itself.
\end{proof}


\section{Rank-\texorpdfstring{$ 3 $}{3} Computations}

\label{sec:BC:rank3}

\begin{secnotation}\label{BC:rank3-notation}
	We denote by $ G $ a group which has crystallographic $ BC_n $-commutator relations with root groups $ (\rootgr{\alpha})_{\alpha \in BC_n} $ for some fixed integer $ n \ge 3 $. We assume that $ \invset{\alpha} $ is non-empty for all roots $ \alpha $ and that $ G $ is rank-2-injective.
\end{secnotation}

The main result of this section is \cref{BC:square-act:summary}, which provides a formula for the action of squares of Weyl elements. In contrast to root gradings of types $ A $, $ B $, $ D $ and $ E $, this formula does not obey the square formula for Weyl elements in all cases. Instead, we will see an additional involution on the short root groups, which we call the \defemph*{short involution}. The short involution will later turn out to be precisely the Jordan module involution on the Jordan module which coordinatises the short root groups.

The main part of the proof of \cref{BC:square-act:summary} is already done by the rank-2 computations in the previous section. We will need the rank-3 assumption only to determine the action of medium-length Weyl elements on medium-length root groups.

Before we begin, we make a quick observation.

\begin{proposition}\label{C:abelian}
	All long and medium-length root groups are abelian.
\end{proposition}
\begin{proof}
	For long root groups, this follows from \cref{BC:long-abelian-lem}. Every medium-length root is contained in a subsystem of type $ A_2 $, so the second assertion holds by \cref{ADE:abelian}.
\end{proof}

\subsection{The Short Involution}

\begin{note}
	In the following, we will often say \enquote{let $ \delta $ be a short root} and then make statements about the long root $ 2\delta $. Since every long root can be written this way, this simply means that we study an arbitrary long root and call it $ 2\delta $.
\end{note}

\begin{definition}[Short involution]\label{BC:def:short-invo}
	Let $ \beta $ be a any short root. The \defemph*{(short) involution on $ \rootgr{\beta} $}\index{short involution} is the map $ \map{}{\rootgr{\beta}}{\rootgr{\beta}}{x_\beta}{\shortinvo{x_\beta} \defl x_\beta^w} $ where $ w \defl w_{2\beta}^2 $ for an arbitrary $ 2\beta $-Weyl element $ w_{2\beta} $. If the root $ \beta $ is not specified or if it is clear from the context, we will sometimes call this map the \defemph*{short involution}.
\end{definition}

\begin{note}
	It is of course misleading to talk about \emph{the} short involution since there is not only one such map but one for each short root group. Whenever we say \enquote{the short involution has some property}, we mean \enquote{for every short root $ \beta $, the involution on $ \rootgr{\beta} $ has this property}.
\end{note}

\begin{note}\label{BC:short-invo-C-note}
	Recall from \cref{BC:CisBC:CasBC} that to any group $ G' $ with a $ C_n $-grading $ (\rootgr{\alpha}')_{\alpha \in C_n} $ we can associate a canonical $ BC_n $-grading $ (\rootgr{\alpha}')_{\alpha \in BC_n} $. This allows us to define the short involution on $ C_n $-graded groups as well. That is, for every long root $ \beta $ of $ C_n $ we have an endomorphism of $ \rootgr{\beta} = \rootgr{\beta/2} $ which is induced by the conjugation with the square of an arbitrary $ \beta $-Weyl element, and we call this endomorphism the short involution on $ \rootgr{\beta} $. Thus confusingly, the short involution on $ C_n $-graded groups is defined on the long root groups. We accept this because, over the course of this chapter, we will nearly always work with the more general case of $ BC_n $-graded groups and rarely refer to the special case of $ C_n $-graded groups.
\end{note}

\begin{lemma}\label{BC:short-invo-prop}
	Let $ \beta $ be a short root. Then the short involution (on $ \rootgr{\beta} $) has the following properties:
	\begin{lemenumerate}
		\item It is well-defined. That is, it does not depend on the choice of $ w_{2\beta} $.
		
		\item \label{BC:short-invo-prop:order}It is a group automorphism of $ \rootgr{\beta} $ of order $ 1 $ or $ 2 $. In particular, $ \shortinvo{(x_\beta^{-1})} = (\shortinvo{x_\beta})^{-1} $ for all $ x_\beta \in \rootgr{\beta} $.
		
		\item \label{BC:short-invo-prop:comm-formula}If $ \alpha, \gamma, \delta $ are roots such that $ (\alpha, \beta, \gamma, \delta) $ is a $ BC_2 $-quadruple, then $ \shortinvo{\commpart{x_\delta}{x_\alpha}{\beta}} = \commpart{x_\delta}{x_\alpha^{-1}}{\beta} $ for all $ x_\alpha \in \rootgr{\alpha} $ and $ x_\delta \in \rootgr{\delta} $.
		
		\item \label{BC:short-invo-prop:weyl-comp}Let $ \rho $ be any root and let $ w_\rho $ be a $ \rho $-Weyl element. Then $ (\shortinvo{x_\beta})^{w_\rho} = \shortinvo{(x_\beta^{w_\rho})} $ for all $ x_\beta \in \rootgr{\beta} $.
	\end{lemenumerate}
\end{lemma}
\begin{proof}
	The first and the third assertion were proven in \cref{BC:square-act:beta-on-beta}. Since the short involution is simply conjugation by some group element, it is a group automorphism. Further, by the formula in~\itemref{BC:short-invo-prop:comm-formula} and by \cref{B:rootisom-beta-delta}, it is clear that $ \shortinvo{(\shortinvo{x_\beta})} = x_\beta $ for all $ x_\beta \in \rootgr{\beta} $, so the short involution is of order $ 1 $ or $ 2 $.
	
	The last assertion is, again, essentially a consequence of the formula in~\itemref{BC:short-invo-prop:comm-formula} and \cref{B:rootisom-beta-delta}. Let $ x_\beta \in \rootgr{\beta} $ be arbitrary, choose roots $ \alpha, \gamma, \delta $ such that $ (\alpha, \beta, \gamma, \delta) $ is a $ BC_2 $-quadruple and choose $ x_\alpha \in \rootgr{\alpha} $, $ x_\delta \in \rootgr{\delta} $ such that $ x_\beta = \commutator{x_\delta}{x_\alpha} $. Then by \cref{basic:commpart-conj}, we have
	\begin{align*}
		(\shortinvo{x_\beta})^{w_\rho} &= (\shortinvo{\commpart{x_\delta}{x_\alpha}{\beta}})^{w_\rho} = \commpart{x_\delta}{x_\alpha^{-1}}{\beta}^{w_\rho} = \commpart{x_\delta^{w_\rho}}{(x_\alpha^{-1})^{w_\rho}}{\refl{\rho}(\beta)} \\
		&= \commpart{x_\delta^{w_\rho}}{(x_\alpha^{w_\rho})^{-1}}{\refl{\rho}(\beta)} = \shortinvo{\commpart{x_\delta^{w_\rho}}{x_\alpha^{w_\rho}}{\refl{\rho}(\beta)}} = \shortinvo{(\commpart{x_\delta}{x_\alpha}{\beta}^{w_\rho})} = \shortinvo{(x_\beta^{w_\rho})}.
	\end{align*}
	This finishes the proof.
\end{proof}

\subsection{The Action of Long Weyl Elements on Medium-length Root Groups}

\begin{lemma}\label{BC:square-act:long-on-med}
	Let $ \delta $ be a short root, let $ \alpha $ be a medium-length root and let $ w_{2\delta} $ be a $ 2\delta $-Weyl element. If $ \alpha $ and $ \delta $ lie in a common $ BC_2 $-subsystem, then $ w_{2\delta}^2 $ acts on $ \rootgr{\alpha} $ by inversion. Otherwise $ w_{2\delta}^2 $ acts trivially on $ \rootgr{\alpha} $.
\end{lemma}
\begin{proof}
	The first assertion follows from \cref{BC:weyldelta-square-on-alpha}. In the second case, $ \alpha $ is adjacent to $ 2\delta $ and $ -2\delta $, so $ w_{2\delta}^2 $ acts trivially on $ \rootgr{\alpha} $.
\end{proof}

\begin{proposition}\label{BC:square-act:long-on-med:cartan}
	Let $ \delta $ be a short root, let $ \alpha $ be a medium-length root and let $ w_{2\delta} $ be a $ 2\delta $-Weyl element. Put $ w \defl w_{2\delta}^2 $ and $ \epsilon \defl (-1)^{\cartanint{\alpha}{2\delta}} $. Then $ x_\alpha^w = x_\alpha^\epsilon $ for all $ x_\alpha \in \rootgr{\alpha} $.
\end{proposition}
\begin{proof}
	It follows from \cref{BC:rootsys:cartan-int-parity} that $ \epsilon = 1 $ if and only if $ \alpha \cdot \delta = 0 $. By \cref{BC:rootsys:common-BC2}, this means that $ \epsilon=-1 $ if and only if $ \alpha $ and $ \delta $ lie in a common $ BC_2 $-subsystem. Thus the assertion is simply a reformulation of \cref{BC:square-act:long-on-med}.
\end{proof}

\subsection{The Action of Long Weyl Elements on Short Root Groups}

\begin{proposition}\label{BC:square-act:long-on-short:cartan}
	Let $ \beta, \delta $ be two short roots, let $ w_{2\delta} $ be a $ 2\delta $-Weyl element and put $ w \defl w_{2\delta}^2 $. Then the following hold:
	\begin{proenumerate}
		\item If $ \beta \in \Set{\pm \delta} $, then $ w $ acts on $ \rootgr{\beta} $ through the short involution. That is, we have $ x_\beta^w = \shortinvo{x_\beta} $ for all $ x_\beta \in \rootgr{\beta} $.
		
		\item If $ \beta \nin \Set{\pm \delta} $, then $ w $ acts trivially on $ \rootgr{\beta} $. In particular, we have $ x_\beta^w = x_\beta^\epsilon $ for all $ x_\beta \in \rootgr{\beta} $ where $ \epsilon \defl (-1)^{\cartanint{\beta}{\delta}} = (-1)^0 = 1 $.
	\end{proenumerate}
\end{proposition}
\begin{proof}
	The first statement is in fact the definition of the short involution, see \cref{BC:def:short-invo}. The second statement is trivial because $ \beta $ is crystallographically adjacent to $ -\delta $ and $ \delta $.
\end{proof}

\subsection{The Action of Medium-length Weyl Elements on the Short Root Groups}

\begin{proposition}\label{BC:square-act:med-on-short}
	Let $ \alpha $ be a medium-length root, let $ \delta $ be a short root, let $ w_\alpha $ be an $ \alpha $-Weyl element and put $ w \defl w_\alpha^2 $. Then the following hold:
	\begin{proenumerate}
		\item \label{BC:square-act:med-on-short:BC2}If $ \alpha $ and $ \delta $ lie in a common $ BC_2 $-subsystem, then $ w $ acts on $ \rootgr{\delta} $ through the short involution. That is, we have $ x_\delta^w = \shortinvo{x_\delta} $ for all $ x_\delta \in \rootgr{\delta} $.
		
		\item \label{BC:square-act:med-on-short:ortho}If $ \alpha $ and $ \delta $ do not lie in a common $ BC_2 $-subsystem, then $ w $ acts trivially on $ \rootgr{\delta} $. In particular, we have $ x_\delta^w = x_\delta^\epsilon $ for all $ x_\delta \in \rootgr{\delta} $ where $ \epsilon \defl (-1)^{\cartanint{\delta}{\alpha}}=1 $.
	\end{proenumerate}
\end{proposition}
\begin{proof}
	At first, assume that $ \alpha $ and $ \delta $ lie in a common $ BC_2 $-subsystem. Then there exist roots $ \beta, \gamma $ such that either $ (\alpha, \delta, \gamma, \beta) $ or $ (\alpha, \beta, \gamma, \delta) $ is a $ BC_2 $-quadruple. In the first case, assertion~\itemref{BC:square-act:med-on-short:BC2} follows directly from \cref{BC:square-act:alpha-beta-on-beta}. In the second case, it is also a consequence of \cref{BC:square-act:alpha-beta-on-beta}, but we have to use in addition that $ w_\alpha $ is also a $ (-\alpha) $-Weyl element (which holds by \thmitemcref{basic:weyl-general}{basic:weyl-general:minus}) and that $ (-\alpha, \delta, \gamma, \beta) $ is a $ BC_2 $-quadruple. This finishes the proof of~\itemref{BC:square-act:med-on-short:BC2}.
	
	Now assume that $ \alpha $ and $ \delta $ do not lie in a common $ BC_2 $-subsystem. Then they must be orthogonal by \cref{BC:rootsys:common-BC2}, so $ \cartanint{\delta}{\alpha} = 0 $ and $ \delta $ is adjacent to $ \alpha $ and $ -\alpha $. This proves the second assertion.
\end{proof}

\subsection{The Action of Medium-length Weyl Elements on Medium-length Root Groups}

In contrast to the previous subsections, we now have to use that $ n \ge 3 $.

\begin{proposition}\label{BC:square-act:med-on-med}
	Let $ \alpha, \gamma $ be two medium-length roots and let $ w_\alpha $ be an $ \alpha $-Weyl element. Put $ w \defl w_\alpha^2 $ and $ \epsilon \defl (-1)^{\cartanint{\gamma}{\alpha}} $. Then $ x_\gamma^w = x_\gamma^\epsilon $ for all $ x_\gamma \in \rootgr{\gamma} $.
\end{proposition}
\begin{proof}
	Denote by $ \roots' $ the root subsystem of $ BC_n $ which is spanned by $ \alpha $ and $ \gamma $. Since $ \alpha, \gamma $ are both of medium length, $ \roots' $ is of type $ A_1 $, $ A_1 \times A_1 $, $ A_2 $ or $ BC_2 $. If it is of type $ A_1 \times A_1 $, then $ \alpha, \gamma $ are orthogonal (so $ \cartanint{\gamma}{\alpha} = 0 $ and $ \epsilon = 1 $) and $ w $ acts trivially on $ \rootgr{\gamma} $. If $ \roots' $ is of type $ A_1 $ or $ A_2 $, then it lies in a subsystem of type $ A_2 $ and the assertion holds by \cref{A2Weyl:cartan-comp}.
	
	Now assume that $ \roots' $ is of type $ BC_2 $. Then there exist short roots $ \beta $, $ \delta $ such that $ (\alpha, \beta, \gamma, \delta) $ is a $ BC_2 $-quadruple. Choose roots $ \gamma_1 $ and $ \gamma_2 $ as in \cref{BC:rootsys:med-decomp} and let $ x_\gamma \in \rootgr{\gamma} $ be arbitrary. Since $ (\gamma_1, \gamma, \gamma_2) $ is an $ A_2 $-triple, it follows from \cref{A2Weyl:rootisom} that there exist $ x_1 \in \rootgr{\gamma_1} $ and $ x_2 \in \rootgr{\gamma_2} $ such that $ x_\gamma = \commutator{x_1}{x_2} $. Note that by the choice of $ \gamma_1, \gamma_2 $ and by \cref{A2Weyl:square-act-lemma-A2}, $ w $ acts on $ \rootgr{\gamma_1} $ and $ \rootgr{\gamma_2} $ by inversion. Thus $ x_\gamma^w = \commutator{x_1^w}{x_2^w} = \commutator{x_1^{-1}}{x_2^{-1}} $, which equals $ \commutator{x_1}{x_2} $ by \cref{basic:comm-add}. In other words, $ w $ acts trivially on $ \rootgr{\gamma} $. Since $ \alpha $ and $ \gamma $ are orthogonal, we further have $ \epsilon = 1 $, so this finishes the proof.
\end{proof}

\subsection{Summary}

We can summarise the results of this section as follows. By \cref{BC:square-act:length-restrict-note}, all interesting cases are covered.

\begin{proposition}\label{BC:square-act:summary}
	Let $ \rho, \zeta $ be any two roots in $ BC_n $ such that $ \rho $ is either short or medium-length and such that $ \zeta $ is either medium-length or long. Let $ w_\zeta $ be a $ \zeta $-Weyl element and put $ w \defl w_\zeta^2 $ and $ \epsilon \defl (-1)^{\cartanint{\rho}{\zeta}} $. Then $ x_\rho^w = x_\rho^\epsilon $ for all $ x_\rho \in \rootgr{\rho} $ unless one of the following conditions is satisfied:
	\begin{stenumerate}
		\item $ \rho $ is short, $ \zeta $ is long and $ \zeta \in \compactSet{\pm 2 \rho} $.
		
		\item $ \rho $ is short, $ \zeta $ is medium-length and $ \rho \cdot \zeta \ne 0 $.
	\end{stenumerate}
	In each of these exceptional cases, we have $ \epsilon = 1 $ and $ w $ acts on $ \rootgr{\rho} $ through the short involution from \cref{BC:def:short-invo}.
\end{proposition}
\begin{proof}
	By the combined statements of \cref{BC:square-act:long-on-med:cartan,BC:square-act:long-on-short:cartan,BC:square-act:med-on-short,BC:square-act:med-on-med}, it only remains to prove the final assertion on $ \epsilon $. In the first case, we have
	\[ \cartanint{\rho}{\zeta} = 2 \frac{\rho \cdot \zeta}{\zeta \cdot \zeta} = 2 \cdot \frac{\pm 2}{4} = \pm 1 \]
	and in the second case,
	\[ \cartanint{\rho}{\zeta} = 2 \frac{\pm 1}{2} = \pm 1. \]
	Thus $ \epsilon = 1 $ in both cases.
\end{proof}


\section{Standard Signs}

\label{sec:BC:stsigns}

\begin{secnotation}
	We fix an integer $ n \ge 3 $ and consider the root system $ BC_n $ in its standard representation (as in \cref{BC:BCn-standard-rep}). We denote by $ G $ a group with a $ BC_n $-pregrading $ (\rootgr{\alpha})_{\alpha \in B_n} $, by $ \ring $ an alternative ring with a nuclear involution $ \rinvmap $ and by $ \jormodtup = (\jormod, \jorsc, \jorprojone, \jorTrone, \psi) $ a Jordan module over $ (\ring, \rinvmap) $. We assume that there exists a coordinatisation $ (\risom{\alpha})_{\alpha \in BC_n} $ of $ G $ by $ \jormodtup $ with standard signs (in the sense of the following \cref{BC:standard-param-def}), and we fix this coordinatisation.
\end{secnotation}

\begin{definition}[Standard signs]\label{BC:standard-param-def}
	Let $ n \in \IN_{\ge 2} $ and let $ G $ be a group with a $ BC_n $-pregrading $ (\rootgr{\alpha})_{\alpha \in BC_n} $. Let $ \ring $ be an alternative ring with nuclear involution $ \rinvmap $ and let $ \jormodtup = (\jormod, \jorsc, \jorprojone, \jorTrone, \psi) $ be a Jordan module over $ (\ring, \rinvmap) $. A \defemph*{coordinatisation of $ G $ by $ \jormodtup $ with standard signs}\index{coordinatisation of a root graded group!with standard signs!type BCn@type $ BC_n $} is a family $ (\risom{\alpha})_{\alpha \in BC_n} $ with the following properties:
	\begin{stenumerate}
		\item For all roots $ \alpha $, the map $ \risom{\alpha} $ is an isomorphism from $ (\ring, +) $ to $ \rootgr{\alpha} $ if $ \alpha $ is medium-length and it is an isomorphism from $ (\jormod, \joradd) $ to $ \rootgr{\alpha} $ if $ \alpha $ is short.
		
		\item The following commutator relations, called the \defemph*{standard commutator relations}, are satisfied.
		For all pairwise distinct $ i,j,k \in \numint{1}{n} $ and all $ a,b \in \ring $, we have
		\begin{align*}
			\commutator{\rismin{i}{j}(a)}{\rismin{j}{k}(b)} &= \rismin{i}{k}(ab), \\
			\commutator{\rismin{i}{j}(a)}{\risplus{j}{k}(b)} &= \risplus{i}{k}\brackets[\big]{\delinv{k<i}\brackets{a \delinv{k<j}(b)}}, \\
			\commutator{\rismin{i}{j}(a)}{\risminmin{k}{i}(b)} &= \risminmin{j}{k}\brackets[\big]{-\delinv{j<k}\brackets{\rinv{a} \delinv{i<k}(b)}}, \\
			\commutator{\risplus{i}{j}(a)}{\risminmin{k}{i}(b)} &= \rismin{j}{k}\brackets[\big]{\delinv{i<j}(a) \delinv{i<k}(b)}.
		\end{align*}
		For all distinct $ i,j \in \numint{1}{n} $ and all $ a,b \in \ring $, we have
		\begin{align*}
			\commutator{\rismin{i}{j}(a)}{\risplus{i}{j}(b)} &= \risshpos{i}\brackets[\big]{\jorTrone(a \delinv{i<j}(b))}, \\
			\commutator{\rismin{i}{j}(a)}{\risminmin{i}{j}(b)} &= \risshneg{j}\brackets[\big]{\jorTrone\brackets[\big]{\rinv{a} \delinv{i<j}(b)}}.
		\end{align*}
		For all distinct $ i,j \in \numint{1}{n} $, all $ v \in \module $ and all $ a \in \ring $, we have
		\begin{align*}
			\commutator{\risshpos{j}(v)}{\rismin{i}{j}(a)} &= \risshpos{i}\brackets[\big]{\jorsc(v, - \rinv{a})} \risplus{i}{j}\brackets[\big]{- \delinv{i>j}\brackets[\big]{a \jorprojone(v)}}, \\
			\commutator{\risshpos{j}(v)}{\risminmin{i}{j}(a)} &= \risshneg{i}\brackets[\big]{\jorsc(v, -\delinv{i>j}(a))} \rismin{j}{i}\brackets[\big]{\rinv{\jorprojone(v)} \delinv{i>j}(a)}, \\
			\commutator{\risshneg{i}(v)}{\rismin{i}{j}(a)} &= \risshneg{j}\brackets[\big]{\jorsc(v,a)} \risminmin{i}{j}\brackets[\big]{-\delinv{i>j}\brackets[\big]{\rinv{a} \jorprojone(v)}}, \\
			\commutator{\risshneg{i}(v)}{\risplus{i}{j}(a)} &= \risshpos{j}\brackets[\big]{\jorsc(v, -\delinv{i>j}(a))} \rismin{j}{i}\brackets[\big]{\delinv{i<j}(a) \jorprojone(u)}.
		\end{align*}
		For all distinct $ i,j \in \numint{1}{n} $ and all $ u,v \in \module $, we have
		\begin{align*}
			\commutator{\risshpos{i}(u)}{\risshpos{j}(v)} &= \begin{cases}
				\risplus{i}{j}\brackets[\big]{\psi(u,v)} & \text{if } i<j, \\
				\risplus{i}{j}\brackets[\big]{-\psi(v,u)} & \text{if } i>j,
			\end{cases} \\
			\commutator{\risshpos{i}(u)}{\risshneg{j}(v)} &= \rismin{i}{j}\brackets[\big]{-\psi(u,v)} \\
			\commutator{\risshneg{i}(u)}{\risshneg{j}(v)} &= \begin{cases}
				\risminmin{i}{j}\brackets[\big]{\psi(v,u)} & \text{if } i<j, \\
				\risminmin{i}{j}\brackets[\big]{-\psi(u,v)} & \text{if } i>j.
			\end{cases}
		\end{align*}
	\end{stenumerate}
\end{definition}

\begin{note}\label{BC:stsign:jor-language}
	The commutator relations in \cref{BC:standard-param-def} are exactly the same that we have seen in \cref{BC:ex-commrel}, but expressed in the language of Jordan modules. In other words, if we take $ \jormodtup $ to be a Jordan module of pseudo-quadratic type as in \cref{BC:jordanmodule-module-ex}, then the commutator relations in \cref{BC:standard-param-def,BC:ex-commrel} are the same. In a similar way, the conjugation formulas in root gradings of type $ BC $ are the same as the ones in in \cref{BC:ex-weylformula-med,BC:ex-weylformula-short}, but expressed in the language of Jordan modules.  For example, the formula $ \risshpos{i}(u,h)^w = \risshpos{j}\brackets[\big]{u\rinvmin{a}, a^{-1} h \rinvmin{a}} $ for all $ (u,h) \in \psgr $ in \thmitemcref{BC:ex-weylformula-med}{BC:ex-weylformula-med:on-short} should be interpreted as $ \risshpos{i}(u)^w = \risshpos{j}\brackets[\big]{\jorsc(u, \rinvmin{a})} $ for all $ u \in \jormod $.
\end{note}

We will show in \cref{BC:thm} that every $ BC_n $-graded group for $ n \ge 3 $ is coordinatised by a Jordan module with standard signs.

As in \cref{B:short-weakly-balanced}, we can show that all long Weyl elements in $ BC_n $-graded groups with a standard coordinatisation are balanced.

\begin{lemma}\label{BC:long-weakly-balanced}
	Assume that $ G $ is rank-2-injective. Let $ \delta $ be a short root and let $ w_{2\delta} $ be a $2 \delta $-Weyl element in $ G $. Choose $ v',v,v'' \in \jormod $ such that $ w_{2\delta} = \risom{-\delta}(v') \risom{\delta}(v) \risom{-\delta}(v'') $. Then $ \jorprojone(v) $ is invertible with inverse $ \jorprojone(v') $ and $ v' = v'' = \jorsc(\modinv{v}, \rinvmin{\jorprojone(v)}) $. In particular, $ w_{2\delta} $ is weakly balanced.
\end{lemma}
\begin{proof}
	We prove this only for $ \delta = \basvec_2 $, the remaining cases being similar. Let $ w_2 $ be a $ 2\basvec_2 $-Weyl element and let $ v', v, v'' \in \jormod $ be such that $ w_2 = \risshneg{2}(v') \risshpos{2}(v) \risshneg{2}(v'') $. For all $ a \in \ring $, we know from \namecrefs{C:basecomp-wdelta-alpha-cor}~\thmitemref{C:basecomp-wdelta-alpha-cor}{C:basecomp-wdelta-alpha-cor:conjformula} and~\thmitemref{C:basecomp-wdelta-gamma-cor}{C:basecomp-wdelta-alpha-cor:conjformula} and the standard commutator relations that
	\begin{align*}
		\rismin{1}{2}(a)^{w_2} &= \commpart{\rismin{1}{2}(a)}{\risshpos{2}(v)}{\basvec_1 + \basvec_2} = \risplus{1}{2}\brackets[\big]{a \jorprojone(v)}, \\
		\risplus{1}{2}(a)^{w_2} &= \commpart{\risplus{1}{2}(a)}{\risshneg{2}(v')}{\basvec_1 - \basvec_2} = \rismin{1}{2}\brackets[\big]{-a \jorprojone(v')}.
	\end{align*}
	Together with \cref{BC:square-act:summary}, this implies that
	\begin{align*}
		\rismin{1}{2}(-1_\ring) &= \rismin{1}{2}(1_\ring)^{w_2^2} = \rismin{1}{2}\brackets[\big]{- \jorprojone(v) \jorprojone(v')}, \\
		\risplus{1}{2}(-1_\ring) &= \risplus{1}{2}(1_\ring)^{w_2^2} = \risplus{1}{2}\brackets[\big]{-\jorprojone(v') \jorprojone(v)}.
	\end{align*}
	We infer that $ \jorprojone(v) $ is invertible with inverse $ \jorprojone(v') $.
	
	Further, we know from \thmitemcref{C:basecomp-wdelta-alpha-cor}{C:basecomp-wdelta-alpha-cor:1G} that
	\begin{align*}
		1_G &= \commutator[\big]{\rismin{1}{2}(a)}{\commpart{\rismin{1}{2}(a)}{\risshpos{2}(v)}{\basvec_1 + \basvec_2}} \commpart{\rismin{1}{2}(a)}{\risshpos{2}(v)}{\basvec_1} \\
		& \hspace{1cm} \mathord{} \cdot \commpart[\big]{\commpart{\rismin{1}{2}(a)}{\risshpos{2}(v)}{\basvec_1 + \basvec_2}}{\risshneg{2}(v'')}{\basvec_1}^{\commpart{\rismin{1}{2}(a)}{\risshpos{2}(v)}{\basvec_1}}
	\end{align*}
	for all $ a \in \ring $ where
	\begin{align*}
		\commutator[\big]{\rismin{1}{2}(a)}{\commpart{\rismin{1}{2}(a)}{\risshpos{2}(v)}{\basvec_1 + \basvec_2}} &= \commutator[\big]{\rismin{1}{2}(a)}{\risplus{1}{2}\brackets[\big]{a \jorprojone(v)}} \\
		&= \risshpos{1}\brackets[\big]{\jorTrone(a \rinv{\jorprojone(v)} \rinv{a})}, \\
		\commpart{\rismin{1}{2}(a)}{\risshpos{2}(v)}{\basvec_1} &= \risshpos{1}\brackets[\big]{\jormin \jorsc(v, -\rinv{a})}, \\
		\commpart[\big]{\commpart{\rismin{1}{2}(a)}{\risshpos{2}(v)}{\basvec_1 + \basvec_2}}{\risshneg{2}(v'')}{\basvec_1} &= \commpart[\big]{\risplus{1}{2}\brackets[\big]{a \jorprojone(v)}}{\risshneg{2}(v'')}{\basvec_1} \\
		&= \risshpos{1}\brackets[\big]{\jormin \jorsc\brackets[\big]{v'', -\rinv{\jorprojone(v)} \rinv{a}}}.
	\end{align*}
	Observe that
	\begin{align*}
		\risshpos{1}\brackets[\big]{\jormin \jorsc\brackets[\big]{v'', -\rinv{\jorprojone(v)} \rinv{a}}} &= \risshpos{1}\brackets[\big]{\jorsc\brackets[\big]{v'', - \rinv{\brackets[\big]{a \jorprojone(v)}}}}^{-1} \\
		&= \commpart{\risshpos{2}(v'')}{\rismin{1}{2}\brackets[\big]{a \jorprojone(v)}}{\basvec_1},
	\end{align*}
	so this element lies in $ \rootgr{2\basvec_1} $ because $ \risshpos{2}(v'') $ lies in $ \rootgr{2\basvec_2} $. In particular, it commutes with $ \commpart{\rismin{1}{2}(a)}{\risshpos{2}(v)}{\basvec_1} $. Thus it follows from the previous computations that
	\begin{align*}
		0_\jormod &= \jorTrone\brackets[\big]{a \rinv{\jorprojone(v)} \rinv{a}} \jormin \jorsc(v, -\rinv{a}) \jormin \jorsc\brackets[\big]{v'', -\rinv{\jorprojone(v)} \rinv{a}}
	\end{align*}
	for all $ a \in \ring $. Putting $ a \defl 1_\ring $ in this equation, we infer that
	\begin{align*}
		0_\jormod &= v \joradd \modinv{v} \jormin \modinv{v} \jormin \modinv{\jorsc\brackets[\big]{v'', \rinv{\jorprojone(v)}}} = v \jormin \modinv{\jorsc\brackets[\big]{v'', \rinv{\jorprojone(v)}}}.
	\end{align*}
	In other words,
	\begin{align*}
		v = \jorsc\brackets[\big]{\modinv{v''}, \rinv{\jorprojone(v)}}.
	\end{align*}
	Applying the map $ \jorsc(\modinv{\mapdot}, \rinvmin{\jorprojone(v)}) $ to both sides of the equation, we obtain
	\begin{align*}
		\jorsc\brackets[\big]{v, \rinvmin{\jorprojone(v)}} &= v''.
	\end{align*}
	Since the $ 2\basvec_2 $-Weyl element $ w_2 = \risshneg{2}(v') \risshpos{2}(v) \risshneg{2}(v'') $ is arbitrary, the same conclusion holds for the Weyl element $ w_2 = \risshneg{2}(\jormin v'') \risshpos{2}(\jormin v) \risshneg{2}(\jormin v') $, which yields
	\begin{align*}
		\jorsc\brackets[\big]{\jormin v, \rinvmin{\jorprojone(v)}} &= \jormin v'.
	\end{align*}
	This implies that
	\begin{align*}
		v' = \jorsc\brackets[\big]{v, \rinvmin{\jorprojone(v)}} = v''
	\end{align*}
	because $ \jorsc $ is additive in the first component. This finishes the proof for the root $ 2\delta = 2\basvec_2 $. The remaining cases can be covered in a similar way.
\end{proof}

We also have the usual identification between Weyl-invertible group elements and \enquote{invertible} elements of the coordinatising algebraic structure.

\begin{proposition}\label{BC:stsign:weyl-char}
	Assume that $ G $ is rank-2-injective. Then the following hold:
	\begin{proenumerate}
		\item \label{BC:stsign:weyl-char:long}Let $ \alpha $ be a medium-length root. Define
		\[ w_\alpha(r) \defl \risom{-\alpha}(-r^{-1}) \risom{\alpha}(r) \risom{-\alpha}(-r^{-1}) \]
		for all invertible $ r \in \ring $. Then the maps
		\[ \map{}{\ringinvset{\ring}}{\invset{\alpha}}{r}{\risom{\alpha}(r)} \midand \map{}{\ringinvset{\ring}}{\weylset{\alpha}}{r}{w_\alpha(r)} \]
		are well-defined bijections. Here $ \ringinvset{\ring} $, $ \invset{\alpha} $ and $ \weylset{\alpha} $ denote the sets of invertible elements in $ \ring $ (in the sense of \cref{ring:def-invertible}), the set of $ \alpha $-invertible elements in $ \rootgr{\alpha} $ and the set of $ \alpha $-Weyl elements, respectively.
		
		\item \label{BC:stsign:weyl-char:short}Let $ \delta $ be a short root and put $ \ringinvset{\jormod} \defl \Set{u \in \jormod \given \jorprojone(u) \in \ringinvset{\ring}} $. Define
		\[ w_\delta(u) \defl \risom{-\delta}\brackets[\big]{\jorsc\brackets[\big]{\modinv{u}, \rinvmin{\jorprojone(u)}}} \risom{\delta}(u) \risom{-\delta}\brackets[\big]{\jorsc\brackets[\big]{\modinv{u}, \rinvmin{\jorprojone(u)}}} \]
		for all $ u \in \ringinvset{\jormod} $. Then the maps
		\[ \map{}{\ringinvset{\jormod}}{\invset{\delta}}{u}{\risom{\delta}(u)} \midand \map{}{\ringinvset{\jormod}}{\weylset{\delta}}{u}{w_\delta(u)} \]
		are well-defined surjections, and the first one is a bijection.
		
		\item \label{BC:stsign:weyl-char:formulas}The Weyl elements defined above satisfy the same conjugation formulas as in \cref{BC:ex-weylformula-med,BC:ex-weylformula-short}, but expressed in the language of Jordan modules. (See \cref{BC:stsign:jor-language}.)
		
		\item \label{BC:stsign:weyl-char:param}Let $ \rootbase $ be the standard rescaled root base of $ BC_n $ and choose an element $ v_0 \in \jormod $ with $ \jorprojone(v_0) = 1_\ring $. Define $ \twistgroup $, $ \invogroup $, $ \inverparsym $ and $ \invoparsym $ as in \cref{BC:ex-twistgroups-def}. Put $ w_\alpha \defl w_\alpha(1_\ring) $ for all medium-length roots in $ \rootbase $ and $ w_\beta \defl w_\beta(v_0) $ for the unique long root in $ \rootbase $. Then $ G $ is parametrised by $ (\twistgroup \times \invogroup, \jormod, \ring) $ with respect to $ \inverparsym \times \invoparsym $ and $ (w_\delta)_{\delta \in \rootbase} $.
	\end{proenumerate}
\end{proposition}
\begin{proof}
	This can be proven in a similar way as \cref{B:stsign:weyl-char}. For the proof of~\itemref{BC:stsign:weyl-char:short}, we need \cref{BC:long-weakly-balanced}.
\end{proof}

\begin{remark}
	It follows from \cref{BC:stsign:weyl-char} that $ G $ satisfies $ \invset{\alpha} = \rootgr{\alpha} \setminus \compactSet{1_G} $ (the additional condition of being an RGD-system) if and only if $ \ring $ is a division ring and $ \jormodtup $ is anisotropic.
\end{remark}

\begin{remark}\label{BC:weyl-char:rem}
	The same arguments as in \cref{ADE:stsign:weyl-char:rem,B:weyl-char:rem} apply to \cref{BC:stsign:weyl-char}: If $ G $ is not rank-2-injective, then it is still true that the maps in \thmitemcref{BC:stsign:weyl-char}{BC:stsign:weyl-char:long} and~\thmitemref{BC:stsign:weyl-char}{BC:stsign:weyl-char:short} are well-defined, but it is no longer clear that they are surjective or bijective.
\end{remark}

\begin{remark}\label{universal:BC}
	As in \cref{universal:ADE,universal:B}, we can use a presentation to define a group $ \hat{G}(\jormodtup, \ring, \rinvmap) $ for every Jordan module $ \jormodtup $ over an alternative ring $ \ring $ with involution $ \rinvmap $\index{root graded group!universal}, and this construction is functorial in $ (\jormodtup, \ring, \rinvmap) $. However, for those pairs $ (\jormodtup, \ring, \rinvmap) $ for which the existence problem is not solved (see \cref{BC:ex-constprob-restriction}), it is not clear that $ \hat{G}(\jormodtup, \ring, \rinvmap) $ is a $ BC_n $-graded group.
\end{remark}


\section{Admissible and Standard Partial Twisting Systems}

\label{sec:BC:sttwist}

\begin{secnotation}\label{BC:secnot:sttwist}
	We fix an integer $ n \ge 3 $ and consider the root system $ BC_n $ in its standard representation with its standard rescaled root base $ \rootbase $ (as in \cref{BC:BCn-standard-rep}). We choose $ (\ring, \ringzero, \rinvmap) \defl (\IC, \IR, \rinvmap) $, $ \module \defl \IC $, $ \map{q}{\module}{\IC}{}{} $ and $ \map{f}{\module \times \module}{\ring}{}{} $ as in \cref{BC:ex-param-choice}.
	We define the $ BC_n $-graded group $ \EU(q) $ with root groups $ (\rootgr{\alpha})_{\alpha \in BC_n} $ as in \cref{BC:ex:EUq-def}. We denote the root isomorphisms from \cref{BC:ex-roothom-def-med,BC:ex-roothom-def-short} by $ (\risom{\alpha})_{\alpha \in B_n} $ and the standard system of Weyl elements from \cref{B:ex-standard-weyl} by $ (w_\delta)_{\delta \in \rootbase} $. Further, we denote by $ (\twistgroup, \inverparsym, \invogroup, \invoparsym) $ the standard partial twisting system of type $ BC_n $ in the sense of the following \cref{BC:standard-partwist-def}.
\end{secnotation}

This section is structured similarly to \cref{sec:B:sttwist}: Motivated by the elementary unitary group from \cref{sec:BC-example}, we introduce the standard partial twisting system of type $ BC $. We will show that it is admissible in the sense that it satisfies a number of desirable properties. In \cref{sec:BC-param}, we will see that every admissible partial twisting system satisfies the condition of the parametrisation theorem.

\begin{definition}[Standard partial twisting system]\label{BC:standard-partwist-def}
	The \defemph*{standard partial twisting system of type $ BC_n $ (with respect to $ \rootbase $)}\index{twisting system!partial!standard (type BCn)@standard (type $ BC_n $)} is the tuple $ (\twistgroup, \inverparsym, \invogroup, \invoparsym) $ where $ \twistgroup \defl \compactSet{\pm 1}^2 $,  $ \invogroup \defl \compactSet{\pm 1} $ and where $ \inverparsym $, $ \invoparsym $ are the $ \rootbase $-parity maps from \cref{BC:ex-twistgroups-def}. If $ G $ is a group with a crystallographic $ BC_n $-grading, then the \defemph*{standard partial twisting system for $ G $ (with respect to $ \rootbase $)} is the same tuple together with the additional information that $ \twistgroup $ acts on all root groups of $ G $ as follows: The first component of $ \twistgroup $ acts on all root groups by inversion while the second component acts trivially on all medium-length root groups and by the short involution (from \cref{BC:def:short-invo}) on the short root groups.
\end{definition}

Similar comments as in \cref{B:standard-partwist-note,B:standard-param-note} apply to \cref{BC:standard-partwist-def}.

\begin{remark}\label{BC:parmap:twistgroup}
	Let $ G $ be a group with a crystallographic $ BC_n $-grading $ (\varrootgr{\alpha})_{\alpha \in BC_n} $ and let $ (\twistgroup, \inverparsym, \invogroup, \invoparsym) $ be the standard partial twisting system for $ G $. We have to verify that the action of $ \twistgroup $ on the root groups of $ G $ is well-defined. For the action on the medium-length root groups, this is evident. For the action on the short root groups, we have to verify that the short involution is of order at most 2 and that it commutes with the group inversion. Both assertions hold by \thmitemcref{BC:short-invo-prop}{BC:short-invo-prop:order}. Further, the twisting action of $ \twistgroup $ commutes with conjugation by Weyl elements by \thmitemcref{BC:short-invo-prop}{BC:short-invo-prop:weyl-comp}. We conclude that $ \twistgroup $ satisfies all the axioms of a twisting group for $ (G, (w_\delta)_{\delta \in \rootbase}) $ (see \cref{param:twist-grp-def}).
\end{remark}

\begin{note}
	The assumption in \cref{BC:standard-partwist-def} that $ G $ has a crystallographic $ BC_n $-grading is necessary to ensure that the short involution exists. In the $ B_n $-situation in \cref{B:standard-partwist-def}, no such assumption was necessary.
\end{note}

\begin{reminder}[compare~\ref{B:parmap:reminder}]\label{BC:parmap:reminder}
	Recall from \cref{BC:ex-is-rgg} that $ (\rootgr{\alpha})_{\alpha \in BC_n} $ is a crystallographic $ BC_n $-grading of $ \EU(q) $. Thus we can apply all rank-2 and rank-3 computations from the previous sections to $ \EU(q) $. Further, we know from \cref{BC:ex-param-thm} that $ \EU(q) $ is parametrised by the standard parameter system $ (\twistgroup \times \invogroup, \psgr(\module), \ring) $
	(from \cref{BC:jordanmod:standard-paramsys}) with respect to $ \inverparsym \times \invoparsym $ and $ (w_\delta)_{\delta \in \rootbase} $ and from \cref{B:ex-param-choice} that $ (\twistgroup \times \invogroup, \psgr(\module), \ring) $ is $ (\inverparsym \times \invoparsym) $-faithful.
\end{reminder}

We now verify some basic properties of the standard partial twisting system of type $ BC $ by performing computations in the group $ \EU(q) $.

\begin{lemma}\label{BC:parmap-properties}
	$ \inverparsym $ is braid-invariant and adjacency-trivial and $ \invoparsym $ is Weyl-invariant and adjacency-trivial.
\end{lemma}
\begin{proof}
	Using that $ (w_\delta)_{\delta \in \rootbase} $ satisfies the braid relations by \cref{braid:all}, this can be proven in the same way as \cref{B:parmap-properties}.
\end{proof}

\begin{lemma}
	$ \invoparsym $ is semi-complete.
\end{lemma}
\begin{proof}
	This follows from the fact that the only subgroups of $ \invogroup $ are $ \compactSet{1} $ and $ \invogroup $, just like in \cref{B:invoparmap-semicomp}.
\end{proof}

\begin{lemma}\label{BC:parmap-trans-invar}
	$ \inverparsym \times \invoparsym $ is transporter-invariant and $ \inverparsym $, $ \invoparsym $ are independent.
\end{lemma}
\begin{proof}
	At first, we consider the orbit of short roots. Put $ \hat{\alpha} \defl \basvec_1 $. We can read off from \cref{BC:ex-parmap-def} that $ \parmoveset{(\twistgroup \times \invogroup)}{\hat{\alpha}}{\beta} $ is contained in $ \compactSet{1_\twistgroup} \times \compactSet{\pm 1_\twistgroup} \times \compactSet{1_\invogroup} $ for all short roots $ \beta $. Further, we have $ \parmoveset{(\twistgroup \times \invogroup)}{\hat{\alpha}}{\hat{\alpha}} = \compactSet{1_\twistgroup} \times \compactSet{\pm 1_\twistgroup} \times \compactSet{1_\invogroup} $ because
	\begin{align*}
		\inverinvopar{\hat{\alpha}}{(\basvec_1 - \basvec_2, \basvec_1 - \basvec_2)} &= \inverinvopar{\basvec_1}{\basvec_1 - \basvec_2} \inverinvopar{\basvec_2}{\basvec_1 - \basvec_2} = (1_\twistgroup, 1_\twistgroup, 1_\invogroup) (1_\twistgroup, -1_\twistgroup, 1_\invogroup) \\
		&= (1_\twistgroup, -1_\twistgroup, 1_\invogroup).
	\end{align*}
	Thus it follows from criterion~\thmitemref{parmap:transport-invar-char}{parmap:transport-invar-char:weak-bound} (with $ \hat{\beta} \defl \hat{\alpha} $) that $ \inverparsym \times \invoparsym $ is transporter-invariant on the orbit of short roots. The same computations apply for $ \hat{\alpha} \defl 2\basvec_1 $, so it is transporter-invariant on the orbit of long roots as well.
	
	Now we consider the orbit of medium-length roots. Put $ \hat{\alpha} \defl \basvec_1 + \basvec_2 $. Again, we can read off from \cref{BC:ex-parmap-def} that $ \parmoveset{(\twistgroup \times \invogroup)}{\hat{\alpha}}{\beta} $ is contained in $ \compactSet{\pm 1_\twistgroup} \times \compactSet{1_\twistgroup} \times \compactSet{\pm 1_\invogroup} $ for all medium-length roots $ \beta $. Further,
	\begin{align*}
		\inverinvopar{\hat{\alpha}}{\basvec_1 - \basvec_2} &= (-1_\twistgroup, 1_\twistgroup, -1_\invogroup) \rightand \\
		\inverinvopar{\hat{\alpha}}{(\basvec_2 - \basvec_3, \basvec_2 - \basvec_3)} &= \inverinvopar{\basvec_1 + \basvec_2}{\basvec_2 - \basvec_3} \inverinvopar{\basvec_1 + \basvec_3}{\basvec_2 - \basvec_3} \\
		&= (1_\twistgroup, 1_\twistgroup, 1_\invogroup) (-1_\twistgroup, 1_\twistgroup, 1_\invogroup) = (-1_\twistgroup, 1_\twistgroup, 1_\invogroup).
	\end{align*}
	Since these two elements generate $ \compactSet{\pm 1_\twistgroup} \times \compactSet{1_\twistgroup} \times \compactSet{\pm 1_\invogroup} $, it follows that
	\[ \parmoveset{(\twistgroup \times \invogroup)}{\hat{\alpha}}{\hat{\alpha}} = \compactSet{\pm 1_\twistgroup} \times \compactSet{1_\twistgroup} \times \compactSet{\pm 1_\invogroup}. \]
	Hence again by criterion~\thmitemref{parmap:transport-invar-char}{parmap:transport-invar-char:weak-bound} (with $ \hat{\beta} \defl \hat{\alpha} $), we infer that $ \inverparsym \times \invoparsym $ is transporter-invariant on the orbit of medium-length roots.
	
	Finally, the previous computations together with \cref{param:transporter-proj} show that $ \inverparsym $ and $ \invoparsym $ are independent in the sense of \cref{param:parmap-indep-def}.
\end{proof}

\begin{lemma}\label{BC:parmap-square}
	Let $ \rho, \zeta $ be any two roots in $ BC_n $ such that $ \rho $ is either short or medium-length and such that $ \zeta $ is either medium-length or long. Then the following hold:
	\begin{lemenumerate}
		\item If $ \rho $ is short, $ \zeta $ is long and $ \zeta \in \compactSet{\pm 2\rho} $, then $ \inverpar{\rho}{\zeta \zeta} = (1,-1) $.
		
		\item If $ \rho $ is short, $ \zeta $ is medium-length and $ \rho \cdot \zeta \ne 0 $, then $ \inverpar{\rho}{\zeta \zeta} = (1,-1) $.
		
		\item In all cases not covered by the previous assertions, we have $ \inverpar{\rho}{\zeta \zeta} = \brackets[\big]{(-1)^{\cartanint{\rho}{\zeta}}, 1} $.
	\end{lemenumerate}
\end{lemma}
\begin{proof}
	This follows from \cref{BC:square-act:summary} because $ (\twistgroup \times \invogroup, \psgr(\module), \ring) $ is $ (\inverparsym \times \invoparsym) $-faithful. In the first two parts, we also use that the Jordan module involution on $ \psgr(\module) $ (through which the second component of $ \twistgroup $ acts on $ \psgr(\module) $) coincides with the short involution on $ \EU(q) $ by \cref{BC:ex-jordan-twist}.
\end{proof}

\begin{lemma}\label{BC:parmap:invo-ortho}
	Let $ \alpha, \gamma $ be two medium-length roots which are orthogonal. Then $ \invopar{\alpha}{\reflbr{\gamma}} = 1_\invogroup $ if $ \alpha $ and $ \gamma $ are crystallographically adjacent and $ \invopar{\alpha}{\reflbr{\gamma}} = -1_\invogroup $ if they are not.
\end{lemma}
\begin{proof}
	At first, assume that $ \alpha $ and $ \gamma $ are crystallographically adjacent. By \cref{rootsys:any-in-rootbase}, there exist $ \lambda \in \IR_{>0} $ such that $ \lambda \gamma $ is a root and a $ \rootbase $-expression $ \word{\delta} $ of $ \lambda \gamma $. Since $ \gamma $ is of medium length, we have $ \lambda = 1 $. By \cref{BC:rootsys:ortho-adj-2}, $ \alpha $ is crystallographically adjacent to $ -\gamma $ as well, so the $ \gamma $-Weyl element $ w_{\word{\delta}} $ in $ \EU(q) $ centralises $ \rootgr{\alpha} $. On the other hand, we have
	\begin{align*}
		\risom{\alpha}(r)^{w_{\word{\delta}}} = \risom{\alpha}(\inverpar{\alpha}{\word{\delta}} \invopar{\alpha}{\word{\delta}}.r)
	\end{align*}
	for all $ r \in \ring $ by \cref{BC:ex-param-thm}. Since the parameter system $ (\twistgroup \times \invogroup, \psgr(\module), \ring) $ is $ (\inverparsym \times \invoparsym) $-faithful by \cref{BC:parmap:reminder}, it follows that $ \inverpar{\alpha}{\word{\delta}} = (1_\twistgroup, 1_\twistgroup) $ and $ \invopar{\alpha}{\word{\delta}} = 1_\invogroup $. This implies that $ \invopar{\alpha}{\reflbr{\gamma}} = 1_\invogroup $ because $ \word{\delta} $ is a $ \rootbase $-expression of $ \gamma $.
	
	Now consider the case that $ \alpha $ and $ \gamma $ are not crystallographically adjacent. By \cref{BC:rootsys:orth-adj-trans}, there exists an element $ u $ of $ \Weyl(BC_n) $ such that $ (\alpha, \gamma)^u = (\basvec_1 + \basvec_2, \basvec_1 - \basvec_2) $. Let $ \word{\rho} $ be any word over $ \rootbase $ such that $ u = \reflbr{\word{\rho}} $. Then it follows from \cref{parmap:stab-conj} that
	\begin{align*}
		\invopar{\alpha}{\reflbr{\gamma}} &= \invopar{\alpha^{\reflbr{\word{\rho}}}}{\reflbr{\word{\rho}}^{-1} \reflbr{\gamma} \reflbr{\word{\rho}}} = \invopar{\alpha^u}{\reflbr{\gamma^u}} = \invopar{\basvec_1 + \basvec_2}{\reflbr{\basvec_1 - \basvec_2}} = \invopar{\basvec_1 + \basvec_2}{\basvec_1 - \basvec_2}.
	\end{align*}
	We can read off from \cref{BC:ex-parmap-def} that $ \invopar{\basvec_1 + \basvec_2}{\basvec_1 - \basvec_2} = -1_\invogroup $. Thus the claimed assertion follows.
\end{proof}

\begin{definition}[Admissible partial twisting system]\label{BC:admissible-parmap-def}
	Let $ G $ be a group with a $ BC_n $-grading $ (\rootgr{\alpha})_{\alpha \in B_n} $ and let $ (w_\delta')_{\delta \in \rootbase} $ be a $ \rootbase $-system of Weyl elements in $ G $. A \defemph*{$ BC_n $-admissible partial twisting system for $ (G, (w_\delta')_{\delta \in \rootbase}) $}\index{twisting system!partial!admissible (type BCn)@admissible (for $ BC_n $)} is a partial twisting system $ (\twistgroup', \inverparsym', \invogroup', \invoparsym') $ for $ (G, (w_\delta')_{\delta \in \rootbase}) $ with the following additional properties:
	\begin{stenumerate}
		\item $ \twistgroup' = \compactSet{\pm 1}^2 $ and $ \invogroup' = \compactSet{\pm 1} $. 
		
		\item \label{BC:admissible-parmap-def:act}The element $ (-1_{\twistgroup'}, 1_{\twistgroup'}) $ acts on all root groups of $ G $ by inversion whereas $ (1_{\twistgroup'}, -1_{\twistgroup'}) $ acts trivially on the medium-length root groups and by the short involution on the short and long root groups.
		
		\item \label{BC:admissible-parmap-def:adj}$ \inverparsym' $ and $ \invoparsym' $ are adjacency-trivial.
		
		\item \label{BC:admissible-parmap-def:square}$ \inverparsym' $ satisfies the formulas in \cref{BC:parmap-square}.

		\item \label{BC:admissible-parmap-def:ortho}If $ \alpha, \gamma $ are orthogonal medium-length roots, then $ \invopar{\alpha}{\gamma} $ is given by the formula in \cref{BC:parmap:invo-ortho}.
	\end{stenumerate}
	We will sometimes refer to such objects as \defemph*{admissible partial twisting systems} if the root system $ BC_n $ is clear from the context.
\end{definition}

\begin{proposition}
	Let $ G $ be a group with a $ BC_n $-grading $ (\rootgr{\alpha})_{\alpha \in B_n} $. Then for any $ \rootbase $-system $ (w_\delta)_{\delta \in \rootbase} $ of Weyl elements in $ G $, the standard partial twisting system $ (\twistgroup, \inverparsym, \invogroup, \invoparsym) $ is a $ BC_n $-admissible partial twisting system for $ (G, (w_\delta)_{\delta \in \rootbase}) $. In particular, admissible partial twisting systems exist for each group with a $ BC_n $-grading.
\end{proposition}
\begin{proof}
	We have to check that $ (\twistgroup, \inverparsym, \invogroup, \invoparsym) $ satisfies the axioms of a partial twisting system for $ (G, (w_\delta)_{\delta \in \rootbase}) $ (see \cref{param:partwist-def}) and the axioms in \cref{B:admissible-parmap-def}. These properties follow from the results in this section.
\end{proof}


\section{The Parametrisation}

\label{sec:BC-param}

\begin{secnotation}
	We fix an integer $ n \ge 3 $ and consider the root system $ BC_n $ in its standard representation with its standard rescaled root base $ \rootbase $ (as in \cref{BC:BCn-standard-rep}). We denote by $ G $ a group which has crystallographic $ BC_n $-commutator relations with root groups $ (\rootgr{\alpha})_{\alpha \in B_n} $ such that $ \invset{\alpha} $ is non-empty for all $ \alpha \in \roots $ and such that $ G $ is rank-2-injective. We fix a $ \rootbase $-system $ (w_\delta)_{\delta \in \rootbase} $ of Weyl elements and we denote by $ (\twistgroup, \inverparsym, \invogroup, \invoparsym) $ a $ BC_n $-admissible partial twisting system for $ (G, (w_\delta)_{\delta \in \rootbase}) $.
\end{secnotation}

The goal of this section is to show that $ G $ satisfies the conditions in the parametrisation theorem with respect to the admissible partial twisting system $ (\twistgroup, \inverparsym, \invogroup, \invoparsym) $. As in \cref{sec:B-param}, the main effort lies in the verification of stabiliser compatibility, for which we will use the criterion from \cref{param:stabcomp-crit-ortho}.

\begin{proposition}\label{BC:square-comp}
	$ G $ is square-compatible with respect to $ \inverparsym $.
\end{proposition}
\begin{proof}
	This follows from \cref{BC:square-act:summary,BC:parmap-square}. Here we use Axiom~\thmitemref{BC:admissible-parmap-def}{BC:admissible-parmap-def:act}.
\end{proof}

\begin{proposition}\label{BC:stab-comp}
	$ G $ is stabiliser-compatible with respect to $ (\inverparsym, \invoparsym) $.
\end{proposition}
\begin{proof}
	Let $ \alpha $ be any root. We know from Axiom~\thmitemref{BC:admissible-parmap-def}{BC:admissible-parmap-def:adj} that $ \inverparsym $ is $ \alpha $-adjacency-trivial and we want to show that $ G $ is $ \alpha $-stabiliser-compatible. At first, assume that $ \alpha $ is short. Then it follows from \cref{param:adj-implies-stab,param:adj-implies-stab:length-rem} that $ G $ is $ \alpha $-stabiliser compatible. Since any long root group is contained in a short root group, it follows that $ G $ is $ \alpha $-stabiliser-compatible for any long root as well.
	
	Now assume that $ \alpha $ is medium-length. Then there exist distinct $ i,j \in \numint{1}{n} $ and signs $ \epsilon_i, \epsilon_j \in \compactSet{\pm 1} $ such that $ \alpha = \epsilon_i \basvec_i + \epsilon_j \basvec_j $. As in \cref{param:stabcomp-crit-ortho}, we consider the following sets:
	\begin{align*}
		\calO \defl{}& \Set{\beta \in \roots \given \alpha \cdot \beta = 0}, \\
		\calA \defl{} & \Set{\beta \in \calO \given \alpha \text{ is crystallographically adjacent to } \beta \text{ and } -\beta} \\
		\bar{\calA}  \defl{}& \calO \setminus \calA.
	\end{align*}
	By \cref{BC:rootsys:ortho-adj}, we have
	\begin{align*}
		\bar{\calA} &= \Set{\epsilon_i\basvec_i - \epsilon_j \basvec_j, \epsilon_j \basvec_j - \epsilon_i \basvec_i} \midand		\calA = \roots \intersect \gen{\basvec_k \given k \in \numint{1}{n} \setminus \compactSet{i,j}}.
	\end{align*}
	Recall that $ \invopar{\alpha}{\reflbr{\beta}} = 1_\invogroup $ for all $ \beta \in \calA $ and $ \invopar{\alpha}{\reflbr{\beta}} = -1_\invogroup $ for all $ \beta \in \bar{\calA} $ by Axiom~\thmitemref{BC:admissible-parmap-def}{BC:admissible-parmap-def:ortho}. Further, observe that $ \bar{\calA} $ contains exactly one $ \rootbase $-positive root. This implies that the condition in \cref{param:stabcomp-crit-ortho} for all $ \rootbase $-positive roots $ \beta, \beta' \in \bar{\calA} $ is trivially satisfied (because $ \beta = \beta' $). Hence all conditions in \cref{param:stabcomp-crit-ortho} hold. We conclude that $ G $ is $ \alpha $-stabiliser-compatible with respect to $ (\inverparsym, \invoparsym) $ and $ (w_\delta)_{\delta \in \rootbase} $.
\end{proof}

\begin{proposition}\label{BC:param-exists}
	There exist an abelian group $ (\ring, +) $, a group $ (\jormod, \joradd) $ (both equipped, as sets, with an action of $ \twistgroup \times \invogroup $) and a family $ (\risom{\alpha})_{\alpha \in BC_n} $ of maps such that the following conditions are satisfied:
	\begin{stenumerate}
		\item \label{BC:param-exists:1}$ (\risom{\alpha})_{\alpha \in BC_n} $ is a parametrisation of $ G $ by $ (\twistgroup \times \invogroup, \jormod, \ring) $ with respect to $ \inverparsym \times \invoparsym $ and $ (w_\delta)_{\delta \in \rootbase} $.
		
		\item The action of $ (-1_\twistgroup, 1_\twistgroup) $ on $ \ring $ and on $ \jormod $ is given by group inversion, the action of $ (1_\twistgroup, -1_\twistgroup) $ on $ \ring $ is trivial and the action of $ (1_\twistgroup, -1_\twistgroup) $ on $ \jormod $ is given by the formula
		\begin{align*}
			\risom{\beta}\brackets[\big]{(1_\twistgroup, -1_\twistgroup).v} = \shortinvo{\risom{\beta}(v)}
		\end{align*}
		for all $ v \in \jormod $ and all short roots $ \beta $ where $ \shortinvo{\risom{\beta}(v)} $ denotes the short involution (from \cref{BC:def:short-invo}).
	\end{stenumerate}
\end{proposition}
\begin{proof}
	This follows from the parametrisation theorem (\cref{param:thm}), whose assumptions are satisfied by \cref{BC:stab-comp,BC:square-comp}.
\end{proof}


\section{Computation of the Blueprint Rewriting Rules}

\label{sec:BC-bluerules}

\begin{secnotation}\label{BC:blue-secnotation}
	We fix an integer $ n \ge 3 $ and consider the root system $ BC_n $ in its standard representation with its standard rescaled root base $ \rootbase $ (as in \cref{BC:BCn-standard-rep}).
	We denote by $ G $ a group with a crystallographic $ BC_n $-grading $ (\rootgr{\alpha})_{\alpha \in B_n} $. We fix a $ \rootbase $-system of Weyl elements $ (w_\delta)_{\delta \in \rootbase} $ and denote its standard $ \Cnsub $-extension by $ (w_\beta)_{\beta \in \Bnsub} $. To simplify notation, we put $ w_{ij} \defl w_{\basvec_i - \basvec_j} $ for all distinct $ i,j \in \numint{1}{n} $ and we put $ w_i \defl w_{2\basvec_i} $ for all $ i \in \numint{1}{n} $. We denote the standard partial twisting system for $ G $ (from \cref{BC:standard-partwist-def}) by $ (\twistgroup, \inverparsym, \invogroup, \invoparsym) $, by $ (\ring, +) $, $ (\jormod, \joradd) $ any groups which satisfy the assertion of \cref{BC:param-exists} and by $ (\risom{\alpha})_{\alpha \in BC_n} $ the corresponding parametrisation of $ G $. The groups $ \ring $ and $ \jormod $ are equipped with actions of $ \twistgroup \times \invogroup $, and we call the maps
	\[ \map{\rinvmap}{\ring}{\ring}{r}{\rinv{r} \defl -1_\invogroup.r} \midand \map{\modinvmap}{\jormod}{\jormod}{v}{\modinv{v} \defl (1_\twistgroup, -1_\twistgroup).v} \]
	the \defemph*{involutions on $ \ring $ and $ \jormod $}, respectively. Further, we choose elements $ v_{-1}, v_0, v_1 \in \jormod $ such that $ w_{n} = \risshneg{n}(v_{-1}) \risshpos{n}(v_0) \risshneg{n}(v_1) $.
\end{secnotation}

We will frequently use the following result without reference. It says that the action of the $ \Cnsub $-extension $ (w_\beta)_{\beta \in \Cnsub} $ of $ (w_\delta)_{\delta \in \rootbase} $ on the root groups is determined by the values in \cref{BC:ex-parmap-def}.

\begin{proposition}\label{BC:Cnsub-conj-anygroup}
	Let $ \alpha \in BC_n $ and let $ \beta \in \Cnsub $. Let $ x \in \ring $ if $ \alpha $ is of medium length and let $ x \in\jormod $ if $ \alpha $ is short. Then $ \risom{\alpha}(x)^{\hat{w}_\beta} = \risom{\refl{\beta}(\alpha)}(\inverpar{\alpha}{\beta} \invopar{\alpha}{\beta}.x) $ where $ \inverpar{\alpha}{\beta} $ and $ \invopar{\alpha}{\beta} $ are the values in \cref{BC:ex-parmap-def}.
\end{proposition}
\begin{proof}
	Denote by $ \word{\beta} $ the standard $ \rootbase $-expression of $ \beta $ from \cref{BC:Cnsub-ext-word}. Since $ G $ is parametrised by $ (\twistgroup \times \invogroup, \jormod, \ring) $ with respect to $ \inverparsym \times \invoparsym $ and $ (w_\delta)_{\delta \in \rootbase} $, we have
	\[ \risom{\alpha}(x)^{\hat{w}_\beta} = \risom{\refl{\beta}(\alpha)}(\inverpar{\alpha}{\word{\beta}} \invopar{\alpha}{\word{\beta}}.x). \]
	We know from \cref{BC:ex-parmap-equality} that $ \inverpar{\alpha}{\word{\beta}} \invopar{\alpha}{\word{\beta}} = \inverpar{\alpha}{\beta} \invopar{\alpha}{\beta} $, so the assertion follows.
\end{proof}

\begin{lemma}\label{BC:invo-lem}
	The involutions on $ \ring $ and $ \jormod $ are group automorphisms with $ \rinv{(\rinv{r})} = r $ for all $ r \in \ring $ and $ \modinv{\modinv{v}} = v $ for all $ v \in \jormod $.
\end{lemma}
\begin{proof}
	Since the involutions are induced by the actions of $ -1_\invogroup $ and $ (1_\twistgroup, -1_\twistgroup) $, respectively, and since both these elements have order $ 2 $, we have $ \rinv{(\rinv{r})} = r $ for all $ r \in \ring $ and $ \modinv{\modinv{v}} = v $ for all $ v \in \jormod $. Further, since $ \inverinvopar{\basvec_2}{\basvec_1 - \basvec_2} = (1, -1, 1) $, we have
	\begin{align*}
		\risshpos{1}\brackets{\modinv{u \joradd v}} &= \risshpos{2}(u \joradd v)^{w_{12}} = \risshpos{2}(u)^{w_{12}} \risshpos{2}(v)^{w_{12}} = \risshpos{1}(\modinv{u}) \risshpos{1}(\modinv{v}) = \risshpos{1}(\modinv{u} \joradd \modinv{v})
	\end{align*}
	for all $ u,v \in \jormod $, which proves that the involution on $ \jormod $ is an automorphism. Similarly,
	\begin{align*}
		\risplus{2}{1}\brackets[\big]{\rinv{(r+s)}} &= \rismin{2}{1}(r+s)^{w_1} = \rismin{2}{1}(r)^{w_1} \rismin{2}{1}(s)^{w_1} \\
		&= \risplus{2}{1}(\rinv{r}) \risplus{2}{1}(\rinv{s}) = \risplus{2}{1}(\rinv{r}+\rinv{s})
	\end{align*}
	for all $ r,s \in \ring $.
\end{proof}

\begin{lemma}\label{BC:blue:short-weyl-decomp}
	For all $ i \in \numint{1}{n} $, we have $ w_{i} = \risshneg{i}(v_{-1}) \risshpos{i}(v_0) \risshneg{i}(v_1) $.
\end{lemma}
\begin{proof}
	By the choice of $ v_{-1}, v_0, v_1 $, this is clear for $ i=n $. Since
	\[ w_{i} = w_{n}^{w_{ni}} \midand \inverinvopar{\pm\basvec_n}{\basvec_n - \basvec_i} = (1,1,1) \]
	the general assertion follows from \cref{BC:Cnsub-conj-anygroup}.
\end{proof}

\begin{definition}[Commutator maps]\label{BC:commmap-def}
	We define
	\begin{align*}
		\map{\mapdot}{\ring \times \ring}{\ring&}{(a,b)}{ab \defl a \rmult b}, \\
		\map{\jorproj}{\jormod \times \ring}{\ring&}{(v,a)}{\jorproj(v,a)}, \\
		\map{\jorsc}{\jormod \times \ring}{\jormod&}{(v,a)}{\jorsc(v,a)}, \\
		\map{\jorTr}{\ring \times \ring}{\jormod&}{(a,b)}{\jorTr(a,b)}, \\
		\map{\psi}{\jormod \times \jormod}{\ring&}{(u,v)}{\psi(u,v)}
	\end{align*}
	to be the unique maps which satisfy the commutator relations
	\begin{align*}
		\commutator{\rismin{1}{2}(a)}{\rismin{2}{3}(b)} &= \rismin{1}{3}(ab), \\
		\commutator{\risshpos{2}(v)}{\rismin{1}{2}(a)} &= \risplus{1}{2}\brackets[\big]{-\jorproj(v,a)} \risshpos{1}\brackets[\big]{\jorsc(v, -\rinv{a})}, \\
		\commutator{\rismin{1}{2}(a)}{\risplus{1}{2}(b)} &= \risshpos{1}\brackets[\big]{\jorTr(a,b)}, \\
		\commutator{\risshpos{1}(u)}{\risshpos{2}(v)} &= \risplus{1}{2}\brackets[\big]{\psi(u,v)}
	\end{align*}
	for all $ a,b \in \ring $ and all $ u,v \in \jormod $.
\end{definition}

\begin{goal}
	In \cref{BC:blue:pi-one-def,BC:blue:Tr-one-def}, we will use $ \jorproj $ and $ \jorTr $ to define maps $ \jorprojone $ and $ \jorTrone $ by specialising one of the components to the element $ 1_\ring $ (which is not yet defined). See also \cref{BC:jordanmod:Tr-jorproj-rem}. Our goal is to show that $ \rmult $ equips $ \ring $ with an alternative ring structure with respect to which $ \rinvmap $ is a nuclear involution and that the maps $ \jorprojone $, $ \jorsc $, $ \jorTrone $, $ \psi $ equip $ \jormod $ with a Jordan module structure over $ (\ring, \rinvmap) $. The involution $ \modinvmap $ will turn out to be the Jordan module involution.
\end{goal}

\begin{lemma}[Rank-2 computations, part 1]\label{BC:param:rank2commrel}
	The following statements hold:
	\begin{lemenumerate}
		\item \label{BC:param:rank2commrel:add}The maps $ \mapdot $, $ \jorproj $, $ \jorsc $, $ \jorTr $ and $ \psi $ are additive in all components except (possibly) for the first component of $ \jorproj $ and the second component of $ \jorsc $.
		
		\item \label{BC:param:rank2commrel:cent}For all $ u $ in the image of $ \jorTr $, we have $ \psi(u,v) = 0 = \psi(v,u) $ and $ u \joradd v = v \joradd u $ for all $ v \in \jormod $.
		
		\item \label{BC:param:rank2commrel:proj}$ \jorproj(u \joradd v,a) = \jorproj(u,a) + \jorproj(v,a) + \psi\brackets[\big]{\jorsc(u, \rinv{a}), v} $ for all $ u,v \in \jormod $, $ a \in \ring $.
		
		\item \label{BC:param:rank2commrel:sc}$ \jorsc(u,a+b) = \jorsc(u, a) \joradd \jorsc(u,b) \joradd \jorTr\brackets[\big]{\rinv{a}, \jorproj(u, \rinv{b})} $ for all $ u \in \jormod $, $ a,b \in \ring $.
	\end{lemenumerate}
\end{lemma}
\begin{proof}
	The maps $ \mapdot $, $ \jorTr $ and $ \psi $ are bi-additive by \cref{basic:comm-add,basic:comm-add:abel}. Further, $ \jorproj $ is additive in the second component by \thmitemcref{B:comm-add}{B:comm-add:alpha-gamma2} and $ \jorsc $ is additive in the first component by \thmitemcref{B:comm-add}{B:comm-add:delta-beta2}. This proves~\itemref{BC:param:rank2commrel:add}. 
	
	Now let $ u $ be in the image of $ \jorTr $. Then $ \risshpos{1}(u) $ lies in $ \commutator{\rootgrmin{1}{2}}{\rootgrplus{1}{2}} $ and thus in $ \rootgrlongpos{1} $. Since $ \rootgrlongpos{1} $ commutes with $ \rootgrshpos{1} $ and with $ \rootgrshpos{2} $, this implies that $ u \joradd v = v \joradd u $ and $ \psi(u,v) = 0 $ for all $ v \in \jormod $. Further, since
	\[ \risshpos{2}(u) = \risshpos{1}(u)^{w_{12}} \in \rootgrlongpos{1}^{w_{12}} = \rootgrlongpos{2} \]
	we also have $ \psi(v,u) = 0 $ for all $ v \in \jormod $.
	
	Now let $ u,v \in \jormod $ and let $ a \in \ring $. Then by \thmitemcref{B:comm-add}{B:comm-add:delta-gamma2}, we have
	\begin{align*}
		\risplus{1}{2}\brackets[\big]{\jorproj(u \joradd v, a)} &= \commpart{\risshpos{2}(u \joradd v)}{\rismin{1}{2}(-a)}{\basvec_1 + \basvec_2} \\
		&= \commpart{\risshpos{2}(u)}{\rismin{1}{2}(-a)}{\basvec_1 + \basvec_2} \commutator[\big]{\commpart{\risshpos{2}(u)}{\rismin{1}{2}(-a)}{\basvec_1}}{\risshpos{2}(v)} \\
		& \hspace{1cm}\mathord{} \cdot \commpart{\risshpos{2}(v)}{\rismin{1}{2}(-a)}{\basvec_1 + \basvec_2} \\
		&= \risplus{1}{2}\brackets[\big]{\jorproj(u, a)} \commutator[\big]{\risshpos{1}\brackets[\big]{\jorsc(u, \rinv{a})}}{\risshpos{2}(v)} \risplus{1}{2}\brackets[\big]{\jorproj(v, a)} \\
		&= \risplus{1}{2}\brackets[\big]{\jorproj(u,a) + \jorproj(v,a) + \psi\brackets[\big]{\jorsc(u, \rinv{a}), v}}.
	\end{align*}
	It follows that~\itemref{BC:param:rank2commrel:proj} holds. Finally, let $ u \in \jormod $ and let $ a,b \in \ring $. Then by \thmitemcref{B:comm-add}{B:comm-add:alpha-beta2},
	\begin{align*}
		\risshpos{1}\brackets[\big]{\jorsc(u, a+b)} &= \commpart[\big]{\risshpos{2}(u)}{\rismin{1}{2}\brackets[\big]{-\rinv{(a+b)}}}{\basvec_1} \\
		&= \commpart[\big]{\risshpos{2}(u)}{\rismin{1}{2}\brackets{-\rinv{b}} \rismin{1}{2}(-\rinv{a})}{\basvec_1} \\
		&= \commpart[\big]{\risshpos{2}(u)}{\rismin{1}{2}(-\rinv{a})}{\basvec_1} \commutator[\big]{\commpart{\risshpos{2}(u)}{\rismin{1}{2}(-\rinv{b})}{\basvec_1 + \basvec_2}}{\rismin{1}{2}(-\rinv{a})} \\
		& \hspace{1cm} \mathord{} \cdot \commpart{\risshpos{2}(u)}{\rismin{1}{2}(-\rinv{b})}{\basvec_1} 
	\end{align*}
	Here
	\begin{align*}
		\commpart[\big]{\risshpos{2}(u)}{\rismin{1}{2}(-\rinv{a})}{\basvec_1} &= \risshpos{1}\brackets[\big]{\jorsc(u,a)}, \\
		\commutator[\big]{\commpart{\risshpos{2}(u)}{\rismin{1}{2}(-\rinv{b})}{\basvec_1 + \basvec_2}}{\rismin{1}{2}(-\rinv{a})} &= \commutator[\big]{\risplus{1}{2}\brackets[\big]{\jorproj(u, \rinv{b})}}{\rismin{1}{2}(-\rinv{a})} \\
		&= \commutator[\big]{\rismin{1}{2}(-\rinv{a})}{\risplus{1}{2}\brackets[\big]{\jorproj(u, \rinv{b})}}^{-1} \\
		&= \risshpos{1}\brackets[\big]{\jorTr\brackets[\big]{\rinv{a}, \jorproj(u, \rinv{b})}}, \\
		\commpart{\risshpos{2}(u)}{\rismin{1}{2}(-\rinv{b})}{\basvec_1} &= \risshpos{1}\brackets[\big]{\jorsc(u,b)}.
	\end{align*}
	Since the image of $ \jorTr $ centralises $ \jormod $, we conclude that
	\begin{align*}
		\risshpos{1}\brackets[\big]{\jorsc(u, a+b)} &= \risshpos{1}\brackets[\big]{\jorsc(u,a) \joradd \jorsc(u,b) \joradd \jorTr\brackets[\big]{\rinv{a}, \jorproj(u, \rinv{b})}}.
	\end{align*}
	This finishes the proof.
\end{proof}

The following computation will not be needed until the end of the following section, but there is no reason to delay it either.

\begin{lemma}\label{BC:comm-mult-computation}
	For all pairwise distinct $ i,j,k \in \numint{1}{n} $ and all $ a,b \in \ring $, the following commutator relations hold:
	\begin{align*}
		\commutator{\rismin{i}{j}(a)}{\rismin{j}{k}(b)} &= \rismin{i}{k}(ab), \\
		\commutator{\rismin{i}{j}(a)}{\risplus{j}{k}(b)} &= \risplus{i}{k}\brackets[\big]{\delinv{k<i}\brackets{a \delinv{k<j}(b)}}, \\
		\commutator{\rismin{i}{j}(a)}{\risminmin{k}{i}(b)} &= \risminmin{j}{k}\brackets[\big]{-\delinv{j<k}\brackets{\rinv{a} \delinv{i<k}(b)}}, \\
		\commutator{\risplus{i}{j}(a)}{\risminmin{k}{i}(b)} &= \rismin{j}{k}\brackets[\big]{\delinv{i<j}(a) \delinv{i<k}(b)}.
	\end{align*}
\end{lemma}
\begin{proof}
	We can conjugate the equations in \cref{BC:commmap-def} by the same Weyl elements as in \cref{B:comm-mult-computation}, but the result is slightly different because the values of the involved parity maps are different.
\end{proof}

\begin{remark}\label{BC:comm-formula-firststep}
	We could perform the same computations as in \cref{BC:comm-mult-computation} for the other commutator maps. However, it will be more efficient to do so at a later point, when we have acquired more information about the maps $ \jorproj $, $ \jorsc $, $ \jorTr $ and $ \psi $. For the moment, we only note that for all $ i<j \in \numint{1}{n} $, conjugating the equations in \cref{BC:commmap-def} by $ w_{2j} w_{1i} $ (where $ w_{kk} $ is interpreted as $ 1_G $) yields that
	\begin{align*}
		\commutator{\risshpos{j}(v)}{\rismin{i}{j}(a)} &= \risplus{i}{j}\brackets[\big]{-\jorproj(v,a)} \risshpos{i}\brackets[\big]{\jorsc(v, -\rinv{a})}, \\
		\commutator{\rismin{i}{j}(a)}{\risplus{i}{j}(b)} &= \risshpos{i}\brackets[\big]{\jorTr(a,b)}, \\
		\commutator{\risshpos{i}(u)}{\risshpos{j}(v)} &= \risplus{i}{j}\brackets[\big]{\psi(u,v)}
	\end{align*}
	for all $ a \in \ring $ and all $ u,v \in \jormod $.
\end{remark}

\begin{lemma}[Rank-2 computations, part 2]\label{BC:isring}
	$ (\ring, +, \rmult) $ is a ring. If we denote its identity element by $ 1_\ring $, we have $ w_{ij} = \rismin{j}{i}(-1_\comring) \rismin{i}{j}(1_\comring) \rismin{j}{i}(-1_\comring) $ for all distinct $ i, j \in \numint{1}{n} $.
\end{lemma}
\begin{proof}
	We have already proven in \thmitemcref{BC:param:rank2commrel}{BC:param:rank2commrel:add} that the multiplication satisfies the distributive law. It remains to verify the existence of an identity element. This can be done in the same way as for $ B_n $ in the proof of \cref{B:isring}.
\end{proof}

\begin{lemma}[Rank-2 computations, part 3]\label{BC:blue:rank2:inv}
	We have
	\[ \jorsc(v, \rinv{1_\ring}) = v \midand \jorsc(v, -\rinv{1_\ring}) = \modinv{v} \]
	for all $ v \in \jormod $ and $ \jorproj(v_0, a) = a $ for all $ a \in \ring $.
\end{lemma}
\begin{proof}
	Let $ v \in \jormod $. Since $ w_{12} = \rismin{2}{1}(-1_\ring) \rismin{1}{2}(1_\ring) \rismin{2}{1}(-1_\ring) $ by \cref{BC:isring}, it follows from \thmitemcref{C:basecomp-walpha-cor-delta}{C:basecomp-walpha-cor-delta:1} and \thmitemcref{C:basecomp-walpha-cor-beta}{C:basecomp-walpha-cor-beta:1} that
	\begin{align*}
		\risshpos{1}(\modinv{v}) &= \risshpos{2}(v)^{w_{12}} = \commpart{\risshpos{2}(v)}{\rismin{1}{2}(1_\ring)}{\basvec_1} = \risshpos{1}\brackets[\big]{\jorsc(v, -\rinv{1_\ring})}, \\
		\risshpos{1}(v) &= \risshpos{2}(v)^{w_{21}} = \commpart{\risshpos{2}(v)}{\rismin{1}{2}(-1_\ring)}{\basvec_1} = \risshpos{1}\brackets[\big]{\jorsc(v, \rinv{1_\ring})}.
	\end{align*}
	This proves the first assertion. Now let $ a \in \ring $ and recall from \cref{BC:blue:short-weyl-decomp} that $ w_{2} = \risshneg{2}(v_{-1}) \risshpos{i}(v_0) \risshneg{2}(v_1) $. Thus by \thmitemcref{C:basecomp-wdelta-alpha-cor}{C:basecomp-wdelta-alpha-cor:conjformula},
	\begin{align*}
		\risplus{1}{2}(a) &= \rismin{1}{2}(a)^{w_2} = \commpart{\rismin{1}{2}(a)}{\risshpos{2}(v_0)}{\basvec_1 + \basvec_2} = \commpart{\risshpos{2}(v_0)}{\rismin{1}{2}(a)}{\basvec_1 + \basvec_2}^{-1} \\
		&= \risplus{1}{2}\brackets[\big]{\jorproj(v_0, a)}.
	\end{align*}
	This proves the second assertion.
\end{proof}

We will show in \cref{BC:blue-inv} that $ \rinv{1_\ring} = 1_\ring $, thereby simplifying the formulas in \cref{BC:blue:rank2:inv}.

We can now compute the rewriting rules. As was the case for $ B_n $, we will only consider the $ BC_3 $-subsystem spanned by $ \Set{\basvec_1 - \basvec_2, \basvec_2 - \basvec_2, \basvec_3} $. See, however, \cref{BC:A3-assoc-rem}.

\begin{definition}[Blueprint rewriting rules]\label{BC:rewriting-def}
	We define the following rewriting rules:
	\begin{align*}
		\map{\blutrans{12}}{\rootgrmin{1}{2} \times \rootgrmin{2}{3} \times \rootgrmin{1}{2}&}{\rootgrmin{2}{3} \times \rootgrmin{1}{2} \times \rootgrmin{1}{2}}{\\\brackets[\big]{\rismin{1}{2}(a), \rismin{2}{3}(b), \rismin{1}{2}(c)}&}{\brackets[\big]{\rismin{2}{3}(c), \rismin{1}{2}(-b-ca), \rismin{2}{3}(a)}}, \\
		\map{\blutrans{12}^{-1}}{\rootgrmin{2}{3} \times \rootgrmin{1}{2} \times \rootgrmin{1}{2}&}{\rootgrmin{1}{2} \times \rootgrmin{2}{3} \times \rootgrmin{1}{2}}{\\\brackets[\big]{\rismin{2}{3}(a), \rismin{1}{2}(b), \rismin{2}{3}(c)}&}{\brackets[\big]{\rismin{1}{2}(c), \rismin{2}{3}(-b-ac), \rismin{1}{2}(a)}}, \\
		\map{\phi_{13}}{\rootgrmin{1}{2} \times \rootgrlongpos{3}&}{\rootgrlongpos{3} \times \rootgrmin{1}{2}}{\\\brackets[\big]{\rismin{1}{2}(a), \rislongpos{3}(v)}&}{\brackets[\big]{\rislongpos{3}(v), \rismin{1}{2}(a)}}, \\
		\map{\phi_{13}^{-1}}{\rootgrlongpos{3} \times \rootgrmin{1}{2}&}{\rootgrmin{1}{2} \times \rootgrlongpos{3}}{\\\brackets[\big]{\rislongpos{3}(v), \rismin{1}{2}(a)}&}{\brackets[\big]{\rismin{1}{2}(a), \rislongpos{3}(v)}}
	\end{align*}
	and
	\[ \map{\blutrans{23}}{\rootgrlongpos{3} \times \rootgrmin{2}{3} \times \rootgrlongpos{3} \times \rootgrmin{2}{3}}{\rootgrmin{2}{3} \times \rootgrlongpos{3} \times \rootgrmin{2}{3} \times \rootgrlongpos{3}}{}{} \]
	which maps $ \brackets[\big]{\risshpos{3}(v), \rismin{2}{3}(a), \risshpos{3}(u), \rismin{2}{3}(b)} $ to
	\[ \brackets[\bigg]{\rinv{b}, \modinv{\jorsc(\modinv{v}, -\rinv{b})} \joradd u \joradd \modinv{\jorTr(b, \rinv{a})}, -\jorproj(\modinv{v}, b) - \rinv{a} - \psi\brackets[\big]{\jorsc(\modinv{v}, -\rinv{b}) \joradd \modinv{u}, \modinv{v}}, \modinv{v}}. \]
\end{definition}

\begin{lemma}
	The maps $ \blutrans{12} $, $ \blutrans{12}^{-1} $, $ \phi_{13} $, $ \phi_{13}^{-1} $ and $ \blutrans{23} $ in \cref{BC:rewriting-def} are blueprint rewriting rules (with respect to $ (w_\delta)_{\delta \in \rootbase'} $). Further, $ \blutrans{12} $ and $ \blutrans{12}^{-1} $ are inverses of each other, and the same holds for $ \phi_{13} $ and $ \phi_{13}^{-1} $.
\end{lemma}
\begin{proof}
	The second statement about inverses is straightforward to verify. Further, we know from \cref{blue:rewriting-switch} that $ \phi_{13} $ and $ \phi_{13}^{-1} $ are blueprint rewriting rules. Observe that the rewriting rule $ \blutrans{12} $ and its inverse are the same rules that we have used for $ A_3 $. Since the restriction of $ \inverparsym $ to $ A_2 = \Set{\basvec_i - \basvec_{i+1} \given i \in \numint{1}{2}} $ yields the same parity map that we have used for $ A_3 $ (see \cref{BC:ex-parmap-restrict}), the same computation as in \cref{blue:A3-rewriting-comp} shows that $ \blutrans{12} $ and $ \blutrans{12}^{-1} $ are blueprint rewriting rules. It remains to show that $ \blutrans{23} $ is a blueprint rewriting rule. For this, let $ a,b,c,d \in \ring $, let $ v,u,h,k \in \module $ and put
	\[ \word{\alpha} \defl \brackets{\basvec_3, \basvec_2 - \basvec_3, \basvec_3, \basvec_2 - \basvec_3} \midand \word{\beta} \defl \brackets{\basvec_2 - \basvec_3, \basvec_3, \basvec_2 - \basvec_3, \basvec_3}. \]
	Further, we set $ x \defl \brackets[\big]{\risshpos{3}(v), \rismin{2}{3}(a), \risshpos{3}(u), \rismin{2}{3}(b)} $. On the one hand, we have
	\begin{align*}
		\blumapG{\word{\alpha}}(x) &= w_3 \risshpos{3}(v) w_{23} \rismin{2}{3}(a) w_3 \risshpos{3}(u) w_{23} \rismin{2}{3}(b) \\
		&= w_3 w_{23} w_3 w_{23} \risshpos{3}(v)^{w_{23} w_3 w_{23}} \rismin{2}{3}(a)^{w_3 w_{23}} \risshpos{3}(u)^{w_{23}} \rismin{2}{3}(b) \\
		&= w_3 w_{23} w_3 w_{23} \risshpos{2}(\modinv{v})^{w_3 w_{23}} \risplus{2}{3}(a)^{w_{23}} \risshpos{2}(\modinv{u}) \rismin{2}{3}(b) \\
		&= w_3 w_{23} w_3 w_{23} \risshpos{3}(\modinv{v}) \risplus{2}{3}(-\rinv{a}) \risshpos{2}(\modinv{u}) \rismin{2}{3}(b).
	\end{align*}
	On the other hand, we put $ x' \defl \brackets[\big]{\rismin{2}{3}(c), \risshpos{3}(h), \rismin{2}{3}(d), \risshpos{3}(k)} $. Then
	\begin{align*}
		\blumapG{\word{\beta}}(x') &= w_{23} \rismin{2}{3}(c) w_3 \risshpos{3}(h) w_{23} \rismin{2}{3}(d) w_3 \risshpos{3}(k) \\
		&= w_{23} w_3 w_{23} w_3 \rismin{2}{3}(c)^{w_3 w_{23} w_3} \risshpos{3}(h)^{w_{23} w_3} \rismin{2}{3}(d)^{w_3} \risshpos{3}(k) \\
		&= w_{23} w_3 w_{23} w_3 \risplus{2}{3}(c)^{w_{23} w_3} \risshpos{2}(\modinv{h})^{w_3} \risplus{2}{3}(d) \risshpos{3}(k) \\
		&= w_{23} w_3 w_{23} w_3 \rismin{2}{3}(\rinv{c}) \risshpos{2}(\modinv{h}) \risplus{2}{3}(d) \risshpos{3}(k).
	\end{align*}
	In order to compare these two terms, we have to change the order of the product in the first one, using the commutator formulas in \cref{BC:comm-formula-firststep} as well as \thmitemcref{group-rel}{group-rel:comm}:
	\begin{align*}
		\risshpos{3}(\modinv{v}) \risplus{2}{3}(-\rinv{a}) \risshpos{2}(\modinv{u}) \rismin{2}{3}(b) \hspace{-5cm} & \\
		&= \risshpos{3}(\modinv{v}) \risplus{2}{3}(-\rinv{a}) \rismin{2}{3}(b) \risshpos{2}(\modinv{u}) \\
		&= \risshpos{3}(\modinv{v}) \rismin{2}{3}(b) \risplus{2}{3}(-\rinv{a}) \commutator{\risplus{2}{3}(-\rinv{a})}{\rismin{2}{3}(b)} \risshpos{2}(\modinv{u}) \\
		&= \risshpos{3}(\modinv{v}) \rismin{2}{3}(b) \risplus{2}{3}(-\rinv{a})  \commutator{\rismin{2}{3}(b)}{\risplus{2}{3}(\rinv{a})} \risshpos{2}(\modinv{u}) \\
		&= \risshpos{3}(\modinv{v}) \rismin{2}{3}(b) \risplus{2}{3}(-\rinv{a}) \risshpos{2}\brackets[\big]{\modinv{u} \joradd \jorTr(b, \rinv{a})} \\
		&= \rismin{2}{3}(b) \risshpos{3}(\modinv{v}) \commutator{\risshpos{3}(\modinv{v})}{\rismin{2}{3}(b)} \risshpos{2}\brackets[\big]{\modinv{u} \joradd \jorTr(b, \rinv{a})} \risplus{2}{3}(-\rinv{a}) \\
		&= \rismin{2}{3}(b) \risshpos{3}(\modinv{v}) \risshpos{2}\brackets[\big]{\jorsc(\modinv{v}, b)} \risplus{2}{3}\brackets[\big]{-\jorproj(\modinv{v}, b)} \risshpos{2}\brackets[\big]{\modinv{u} \joradd \jorTr(b, \rinv{a})} \risplus{2}{3}(-\rinv{a}) \\
		&= \rismin{2}{3}(b)  \risshpos{2}\brackets[\big]{\jorsc(\modinv{v}, b) \joradd \modinv{u} \joradd \jorTr(b, \rinv{a})} \risplus{2}{3}(-\jorproj(\modinv{v}, b)-\rinv{a}) \\
		& \hspace{1cm} \mathord{} \cdot \commutator[\big]{\risshpos{3}(\modinv{v})}{\risshpos{2}\brackets[\big]{\jorsc(\modinv{v}, b) \joradd \modinv{u} \joradd \jorTr(b, \rinv{a})}}  \risshpos{3}(\modinv{v}) \\
		&= \rismin{2}{3}(b)  \risshpos{2}\brackets[\big]{\jorsc(\modinv{v}, b) \joradd \modinv{u} \joradd \jorTr(b, \rinv{a})} \\
		& \hspace{1cm} \mathord{} \cdot \risplus{2}{3}\brackets[\big]{-\jorproj(\modinv{v}, b)-\rinv{a} - \psi\brackets[\big]{\jorsc(\modinv{v}, b) \joradd \modinv{u} \joradd \jorTr(b, \rinv{a}), \modinv{v}}} \risshpos{3}(\modinv{v})
	\end{align*}
	Note that
	\[ \psi\brackets[\big]{\jorsc(\modinv{v}, b) \joradd \modinv{u} \joradd \jorTr(b, \rinv{a}), \modinv{v}} = \psi\brackets[\big]{\jorsc(\modinv{v}, b) \joradd \modinv{u}, \modinv{v}} \]
	by \thmitemcref{BC:param:rank2commrel}{BC:param:rank2commrel:cent}. Recall that $ w_3 w_{23} w_3 w_{23} = w_{23} w_3 w_{23} w_3 $ by \cref{braid:all} and that $ G $ is rank-2-injective. We conclude that $ \blumapG{\word{\alpha}}(x) = \blumapG{\word{\beta}}(x') $ if the following conditions are satisfied:
	\begin{gather*}
		c = \rinv{b}, \qquad h = \modinv{\jorsc(\modinv{v}, b)} \joradd u \joradd \modinv{\jorTr(b, \rinv{a})}, \\
		d = \modinv{v}, \qquad k = -\jorproj(\modinv{v}, b)-\rinv{a} - \psi\brackets[\big]{\jorsc(\modinv{v}, b) \joradd \modinv{u}, \modinv{v}}.
	\end{gather*}
	These conditions are equivalent to the equation $ x' = \blutrans{23}(x) $, which shows that $ \blutrans{23} $ is a blueprint rewriting rule. This finishes the proof.
\end{proof}


\section{Blueprint Computations}

\label{sec:BC:blue-comp}

\Cref{BC:blue-secnotation} continues to hold.

\begin{remark}[Blueprint computation]\label{BC:blue-comp}
	We now have everything ready to perform the blueprint computation for $ BC_3 $. Since $ B_3 $ and $ BC_3 $ have the same Weyl group, we can use the same homotopy cycle as in \cref{fig:B:hom-cycle} on page~\pageref{fig:B:hom-cycle}. We begin with the tuple
	\[ \brackets[\big]{\risom{3}(u), \risom{2}(a), \risom{1}(b), \risom{3}(v), \risom{2}(c), \risom{3}(w), \risom{2}(d), \risom{1}(r), \risom{2}(s)} \]
	where $ \risom{1} \defl \rismin{1}{2} $, $ \risom{2} \defl \rismin{2}{3} $, $ \risom{3} \defl \risshpos{3} $ and where $ u,v,w \in \jormod $ and $ a,b,c,d,r,s \in \ring $ are arbitrary. Working down rows~1 to~12 in the homotopy cycle and applying the respective blueprint rewriting rules in the process, we obtain a tuple
	\[ \brackets[\big]{\risom{2}(x_1), \risom{1}(x_2), \risom{2}(x_3), \risom{3}(x_4), \risom{2}(x_5), \risom{3}(x_6), \risom{1}(x_7), \risom{2}(x_8), \risom{3}(x_9)}. \]
	Conversely, working up from row~23 to row~12, we obtain a tuple
	\[ \brackets[\big]{\risom{2}(y_1), \risom{1}(y_2), \risom{2}(y_3), \risom{3}(y_4), \risom{2}(y_5), \risom{3}(y_6), \risom{1}(y_7), \risom{2}(y_8), \risom{3}(y_9)} \]
	where $ x_i, y_i \in \jormod $ for all $ i \in \Set{4, 6, 9} $ and $ x_i, y_i \in \ring $ for all $ i \in \Set{1,2,3,5,7,8} $. Now $ x_i = y_i $ for all $ i \in \numint{1}{9} $ by \cref{blue:thm}. The results of this computation can be found in \cref{fig:BC:blue-short,fig:BC:blue-mid,fig:BC:blue-long}. The intermediate steps of the above computation have been performed with GAP \cite{GAP4}.
\end{remark}

\premidfigure

\begin{figure}[htb]
	\centering$ \begin{gathered}
		(1) \; \rinv{s} = \rinv{s}, \quad (2) \; -\rinv{r} - \rinv{s} \rinv{d} = -\rinv{r} - \rinv{(ds)}, \quad (3) \; \rinv{d} = \rinv{d}, \quad (9) \; u=u
	\end{gathered} $
	\caption{The equations \enquote{$ (i) \; x_i = y_i $} for $ i \in \Set{1,2,3,9} $ in \cref{BC:blue-comp}.}
	\label{fig:BC:blue-short}
\end{figure}

\begin{figure}[htb]
	\centering$ \begin{aligned}
		\modinv{\jorsc\brackets{u,-\rinv{s}}}\joradd \modinv{v}\joradd \modinv{\jorTr\brackets{s,\rinv{a}}} &= \jorsc\brackets{\modinv{u},-\rinv{s}}\joradd \modinv{v}\joradd \jorTr\brackets{s,\rinv{a}}, \\
		-\jorproj\brackets{u,s}-\rinv{a}-\psi\brackets[\big]{\jorsc\brackets{u,-\rinv{s}}\joradd v,u} &= -\jorproj\brackets{\modinv{u},s}-\rinv{a}-\psi\brackets[\big]{\jorsc\brackets{\modinv{u},-\rinv{s}}\joradd \modinv{v},\modinv{u}}
	\end{aligned} $
	\caption{The equations $ x_6 = y_6 $ and $ x_8 = y_8 $ in \cref{BC:blue-comp}, respectively.}
	\label{fig:BC:blue-mid}
\end{figure}

\begin{figure}[htb]
	\centering$ \begin{aligned}
		x_4 &= \modinv{\jorsc\brackets{\modinv{u},-\rinv{r}}}\joradd \modinv{\jorsc\brackets{\modinv{v},-\rinv{d}}}\joradd w\joradd \modinv{\jorTr\brackets{d,\rinv{c}}}\joradd \modinv{\jorTr\brackets{r,-\rinv{b}-\rinv{\brackets{a\rinv{d}}}}} \\
		y_4 &= \modinv{\jorsc\brackets[\big]{\jorsc\brackets{\modinv{u},-\rinv{s}}\joradd \modinv{v}\joradd \jorTr\brackets{s,\rinv{a}},-\rinv{d}}}\joradd \modinv{\jorsc\brackets[\big]{u,\rinv{r}+\rinv{\brackets{ds}}}}\joradd w  \\
		&\hspace{1cm}\mathord{} \joradd \modinv{\jorTr\brackets{-r-ds,\rinv{b}}} \joradd \modinv{\jorTr\brackets[\big]{d,\rinv{c}+\rinv{\brackets{sb}}-\rinv{\sqbrackets{A\brackets{-\rinv{r}-\rinv{\brackets{ds}}}}}}} \\
		&\hspace{1cm}\text{where } A \defl -\jorproj\brackets{\modinv{u},s}-\rinv{a}-\psi\brackets[\big]{\jorsc\brackets{\modinv{u},-\rinv{s}}\joradd \modinv{v},\modinv{u}}, \\
		x_5 &= -\jorproj\brackets{\modinv{v},d}-\rinv{c}-\psi\brackets[\big]{\jorsc\brackets{\modinv{v},-\rinv{d}}\joradd \modinv{w},\modinv{v}}+ra\\
		&\hspace{1cm}\mathord{}-\brackets[\big]{-\jorproj\brackets{\modinv{u},r}+\rinv{b}+\rinv{\brackets{a\rinv{d}}}-\psi\brackets{B,\modinv{u}}}\rinv{s} \\
		& \hspace{1cm} \text{where } B \defl \jorsc\brackets{\modinv{u},-\rinv{r}}\joradd \jorsc\brackets{\modinv{v},-\rinv{d}}\joradd \modinv{w}\joradd \jorTr\brackets{d,\rinv{c}}, \\
		y_5 &= -\jorproj\brackets[\big]{\jorsc\brackets{\modinv{u},-\rinv{s}}\joradd \modinv{v}\joradd \jorTr\brackets{s,\rinv{a}},d}-\rinv{c}-\rinv{\brackets{sb}}\\
		&\hspace{1cm}\mathord{}+\rinv{\sqbrackets[\big]{\brackets[\big]{-\jorproj\brackets{\modinv{u},s}-\rinv{a}-\psi\brackets[\big]{\jorsc\brackets{\modinv{u},-\rinv{s}}\joradd \modinv{v},\modinv{u}}}\brackets[\big]{-\rinv{r}-\rinv{\brackets{ds}}}}} \\
		&\hspace{1cm}\mathord{}-\psi\brackets[\big]{C \joradd \modinv{w}\joradd \jorTr\brackets{-r-ds,\rinv{b}},\jorsc\brackets{\modinv{u},-\rinv{s}}\joradd \modinv{v}\joradd \jorTr\brackets{s,\rinv{a}}}, \\
		&\hspace{1cm} \text{where } C \defl \jorsc\brackets[\big]{\jorsc\brackets{\modinv{u},-\rinv{s}}\joradd \modinv{v}\joradd \jorTr\brackets{s,\rinv{a}},-\rinv{d}} \joradd \jorsc\brackets[\big]{u,\rinv{r}+\rinv{\brackets{ds}}}, \\
		x_7 &= -\jorproj\brackets{\modinv{u},r}+\rinv{b}+\rinv{\brackets{a\rinv{d}}}-\psi\brackets[\big]{\jorsc\brackets{\modinv{u},-\rinv{r}}\joradd \jorsc\brackets{\modinv{v},-\rinv{d}}\joradd \modinv{w}\joradd \jorTr\brackets{d,\rinv{c}},\modinv{u}}, \\
		y_7 &= \jorproj\brackets{u,-r-ds}+\rinv{b}+\psi\brackets[\big]{\jorsc\brackets{u,\rinv{r}+\rinv{\brackets{ds}}}\joradd \modinv{w},u}\\
		&\hspace{1cm}\mathord{}-d\brackets[\big]{-\jorproj\brackets{\modinv{u},s}-\rinv{a}-\psi\brackets{\jorsc\brackets{\modinv{u},-\rinv{s}}\joradd \modinv{v},\modinv{u}}}
	\end{aligned} $
	\caption{The values of $ x_4, y_4, x_5, y_5, x_7, y_7 $ in \cref{BC:blue-comp}.}
	\label{fig:BC:blue-long}
\end{figure}

\postmidfigure

\begin{note}
	In the following computations, we have to be careful to remember that the group $ \jormod $ is not abelian. For example, it is not possible to \enquote{subtract $ \modinv{v} $ from equation~6} to obtain
	\[ \modinv{\jorsc\brackets{u,-\rinv{s}}}\joradd \modinv{\jorTr\brackets{s,\rinv{a}}} = \jorsc\brackets{\modinv{u},-\rinv{s}}\joradd  \jorTr\brackets{s,\rinv{a}} \]
	because $ \modinv{v} $ lies \enquote{in the middle} of equation~6. (However, it is still true that the identity above is a consequence of equation~6 because we can simply put $ v \defl 0_\jormod $.) Similarly, we will avoid refering to inverses in the group $ (\jormod, \joradd) $ whenever possible.
\end{note}

As in \cref{sec:B-bluecomp}, the goal of this section is to show that equations~1 to~9 are equivalent to the axioms of a Jordan module. Clearly, equations~1, 3 and 9 are trivial.

\begin{remark}[Associativity for $ n \ge 4 $]\label{BC:A3-assoc-rem}
	If $ n \ge 4 $, then $ BC_n $ contains a parabolic subsystem of type $ A_3 $ whose root groups are parametrised by the ring $ \ring $. In this case, the blueprint computations for $ A_3 $ yield that $ \ring $ must be associative. However, this is the only identity that can be obtained from the assumption that $ n \ge 4 $. Note that there is no need to perform a (long) blueprint computation in a rank-4 subsystem: It suffices to perform two separate (short) blueprint computations in the subsystems $ A_3 $ and $ BC_3 $. Of course, the blueprint computation for $ A_3 $ has already been done in \cref{sec:A3-blue}.
\end{remark}

\begin{lemma}\label{BC:blue-inv}
	The map $ \rinvmap $ is an involution of the ring $ \ring $. In particular, $ \rinv{1_\ring} = 1_\ring $.
\end{lemma}
\begin{proof}
	The first assertion follows from equation~2 in \cref{fig:BC:blue-short}. The second assertion is a consequence of the first one by \cref{rinv:1-symmetric}.
\end{proof}

\cref{BC:blue-inv} clearly covers all non-trivial identities which can be derived from \cref{fig:BC:blue-short}.

\begin{lemma}[Equation~6]\label{BC:blue:eq6lem}
	We have $ \modinv{\jorsc(u,s)} = \jorsc(\modinv{u}, s) $ and $ \modinv{\jorTr(a,b)} = \jorTr(a,b) $ for all $ a,b,s \in \ring $ and all $ u \in \jormod $.
\end{lemma}
\begin{proof}
	The first assertion follows from equation~6 in \cref{fig:BC:blue-mid} by putting $ v \defl 0 $ and $ a \defl 0 $. The second assertion follows from the same equation by putting $ u \defl 0 $ and $ v \defl 0 $.
\end{proof}

\begin{lemma}[Equation~8]\label{BC:blue:eq8lem}
	We have $ \psi(\modinv{u}, \modinv{v}) = \psi(u,v) $ and $ \jorproj(\modinv{u}, s) = \jorproj(u, s) $ for all $ u,v \in \jormod $ and all $ s \in \ring $.
\end{lemma}
\begin{proof}
	The first assertion follows from equation~8 in \cref{fig:BC:blue-mid} by putting $ s \defl 0 $. Using \cref{BC:blue:eq6lem}, we can now simplify equation~8 to obtain
	\[ -\jorproj\brackets{u,s}-\rinv{a}-\psi\brackets[\big]{\jorsc\brackets{u,-\rinv{s}}\joradd v,u} = -\jorproj\brackets{\modinv{u},s}-\rinv{a}-\psi\brackets[\big]{\jorsc\brackets{u,-\rinv{s}}\joradd v,u}. \]
	This is equivalent to the second assertion.
\end{proof}

Again, it is clear that \cref{BC:blue:eq6lem,BC:blue:eq8lem} cover everything that can be deduced from \cref{fig:BC:blue-mid}.

\begin{lemma}[Equation~7]\label{BC:blue:eq7lem}
	The following hold for all $ u,v \in \jormod $ and all $ a,b \in \ring $:
	\begin{lemenumerate}
		\item \label{BC:blue:eq7lem:psi-min}$ \psi(\modinv{u}, v) = \psi(u,\modinv{v}) = -\psi(u,v) $.
		
		\item \label{BC:blue:eq7lem:psi-jor}$ \psi\brackets[\big]{\jorsc(u, a), v} = \rinv{a} \psi(u,v) $.
		
		\item \label{BC:blue:eq7lem:pi}$ \jorproj(u,a) = a \jorproj(u, 1_\ring) $.
		
		\item \label{BC:blue:eq7lem:nucl-mix}$ (ab) \jorproj(u, 1_\ring) + (ab) \psi(u,u) = a\brackets[\big]{b\jorproj(u,1_\ring)}+ a\brackets[\big]{b \psi(u,u)} $.
	\end{lemenumerate}
\end{lemma}
\begin{proof}
	Putting all variables except for $ u $ and $ w $ in equation~7 to zero, we see that $ -\psi(\modinv{w}, \modinv{u}) = \psi(\modinv{w}, u) $. Together with \cref{BC:blue:eq8lem} and the fact that $ \modinv{\mapdot} $ is of order at most $ 2 $, we infer that~\itemref{BC:blue:eq7lem:psi-min} holds. Putting all variables except for $ u $, $ v $, and $ d $ to zero, we obtain
	\begin{align*}
		-\psi\brackets[\big]{\jorsc(\modinv{v}, -\rinv{d}), \modinv{u}} = d\psi(\modinv{v}, \modinv{u}),
	\end{align*}
	which proves \itemref{BC:blue:eq7lem:psi-jor}.
	
	Observe that the term $ \rinv{b} + d \rinv{a} $ appears on both sides of equation~7, that $ \psi(\modinv{w}, \modinv{u}) = -\psi(\modinv{w}, u) $ by~\itemref{BC:blue:eq7lem:psi-min}, that $ \psi\brackets[\big]{\jorTr(d, \rinv{c}), \modinv{u}} = 0_\ring $ by~\thmitemcref{BC:param:rank2commrel}{BC:param:rank2commrel:cent} and that $ \jorproj(\modinv{u},r) = \jorproj(u,r) $ and $ \jorproj(\modinv{u},s) = \jorproj(u,s) $ by \cref{BC:blue:eq8lem}. It follows that equation~7 is equivalent to $ x_7' = y_7' $ where
	\begin{align*}
		x_7' &=-\jorproj\brackets{u,r}-\psi\brackets[\big]{\jorsc\brackets{\modinv{u},-\rinv{r}}\joradd \jorsc\brackets{\modinv{v},-\rinv{d}},\modinv{u}}, \\
		y_7' &= \jorproj\brackets{u,-r-ds}+\psi\brackets[\big]{\jorsc\brackets{u,\rinv{r}+\rinv{\brackets{ds}}},u}\\
		& \hspace{2cm} \mathord{}-d\brackets[\big]{-\jorproj\brackets{u,s}-\psi\brackets{\jorsc\brackets{\modinv{u},-\rinv{s}}\joradd \modinv{v},\modinv{u}}}.
	\end{align*}
	Using~\itemref{BC:blue:eq7lem:psi-min}, the bi-additivity of $ \psi $ and the additivity of $ \jorproj $ in the second component as well as \cref{BC:blue:eq6lem}, we see that
	\begin{align*}
		x_7' &= -\jorproj\brackets{u,r}-\psi\brackets[\big]{\jorsc\brackets{\modinv{u},-\rinv{r}},\modinv{u}} - \psi\brackets[\big]{\jorsc\brackets{\modinv{v},-\rinv{d}}, \modinv{u}} \\
		&= -\jorproj\brackets{u,r} +r \psi(u, u) + d\psi(v,u) \rightand \\
		y_7' &= -\jorproj\brackets{u,r} - \jorproj(u,ds) + (r+ds)\psi(u,u) + d \jorproj(u, s) - d\brackets[\big]{s \psi(u,u)} + d\psi(v,u).
	\end{align*}
	Cancelling the term $ -\jorproj(u,r)+r\psi(u,u) + d\psi(v,u) $ on both sides, it follows that equation~7 is equivalent to
	\begin{align*}
		0_\ring &= -\jorproj\brackets{u,ds} + (ds) \psi(u,u) + d\jorproj(u,s) - d\brackets[\big]{s \psi(u,u)}.
	\end{align*}
	Putting $ s \defl 1_\ring $, we infer that~\itemref{BC:blue:eq7lem:pi} is true. With this knowledge, equation~7 is equivalent to
	\[ d\brackets[\big]{s\jorproj(u,1_\ring)} + (ds) \psi(u,u) = (ds) \jorproj(u, 1_\ring) + d\brackets[\big]{s \psi(u,u)}, \]
	which is exactly~\itemref{BC:blue:eq7lem:nucl-mix}.
\end{proof}

We have seen in \thmitemcref{BC:blue:eq7lem}{BC:blue:eq7lem:pi} that the map $ \jorproj $ is uniquely determined by the map $ \jorproj(\mapdot, 1_\ring) $. This motivates the following definition.

\begin{definition}\label{BC:blue:pi-one-def}
	We define a map $ \map{\jorprojone}{\jormod}{\ring}{u}{\jorproj(u, 1_\ring)} $.
\end{definition}

\begin{note}
	We will later see that a stronger version of \thmitemcref{BC:blue:eq7lem}{BC:blue:eq7lem:nucl-mix} holds: In fact, we have $ (ab) \jorproj(u, 1) = a\brackets[\big]{b \jorproj(u, 1)} $ and $ (ab) \psi(u,u) = a \brackets[\big]{b \psi(u,u)} $. We will also see that $ \ring $ is alternative, so by \cref{ring:alternative-nucleus}, these equations imply that $ \jorproj(u,1) $ and $ \psi(u,u) $ lie in the nucleus of $ \ring $. Even more, we will see that $ \psi(u,v) $ lies in the nucleus for all $ u,v \in \jormod $.
\end{note}

Once more, the proof of \cref{BC:blue:eq7lem} shows that no further identities can be obtained from equation~7. Now the only remaining equations are~4 and~5. Even though equation~4 is shorter, we will continue with equation~5 since it is an equation in $ \ring $ and not in $ \jormod $. This makes it easier to evaluate because $ (\ring, +) $ is abelian.

\begin{lemma}[Equation~5, part 1]\label{BC:blue:eq5lem}
	The following identities hold for all $ u,v \in \jormod $ and all $ a,s,d \in \ring $:
	\begin{lemenumerate}
		\item \label{BC:blue:eq5lem:psi-skew}$ \rinv{\psi(u,v)} = -\psi(v,u) $.
		
		\item \label{BC:blue:eq5lem:psi-skew-cor}
		$ \psi\brackets[\big]{u, \jorsc(v, a)} = \psi(u,v)a $.
		
		\item \label{BC:blue:eq5lem:pi-Tr}$ \jorproj\brackets[\big]{\jorTr(s,a), d} = (ds)\rinv{a} + (da) \rinv{s} $.
	\end{lemenumerate}
\end{lemma}
\begin{proof}
	We consider equation~5 in \cref{fig:BC:blue-long}. Putting all variables except for $ u $, $ v $ and $ r $ to zero, we obtain
	\[ 0 = \rinv{\brackets[\big]{\psi(\modinv{v}, \modinv{u}) \rinv{r}}} - \psi\brackets[\big]{\jorsc(u, \rinv{r}), \modinv{v}}. \]
	With $ r \defl 1_\ring $, \cref{BC:blue:rank2:inv} and \thmitemcref{BC:blue:eq7lem}{BC:blue:eq7lem:psi-min}, this yields~\itemref{BC:blue:eq5lem:psi-skew}. Now~\itemref{BC:blue:eq5lem:psi-skew-cor} is a consequence of~\itemref{BC:blue:eq5lem:psi-skew} and the previous result~\thmitemref{BC:blue:eq7lem}{BC:blue:eq7lem:psi-jor} because
	\begin{align*}
		\psi\brackets[\big]{u, \jorsc(v, a)} &= -\rinv{\psi\brackets[\big]{\jorsc(v, a), u}} = -\rinv{\brackets[\big]{\rinv{a} \psi(v,u)}} = \psi(u,v) a.
	\end{align*}
	Replacing all variables in equation~5 except for $ a $, $ d $ and $ s $ by zero, we see that
	\begin{align*}
		-\rinv{(a \rinv{d})} \rinv{s} &= -\pi\brackets[\big]{\jorTr(s, \rinv{a}), d} + \rinv{\brackets[\big]{\rinv{a} \rinv{(ds)}}} - \psi\brackets[\big]{\jorsc\brackets[\big]{\jorTr(s, \rinv{a}), -\rinv{d}}, \jorTr(s, \rinv{a})}.
	\end{align*}
	Recall that the last summand on the right-hand side is trivial by \thmitemcref{BC:param:rank2commrel}{BC:param:rank2commrel:cent}. Thus we obtain
	\[ \jorproj\brackets[\big]{\jorTr(s,\rinv{a}), d} = (ds)a + (d\rinv{a}) \rinv{s}. \]
	This proves~\itemref{BC:blue:eq5lem:pi-Tr}.
\end{proof}

\begin{remark}\label{BC:blue:psi-mnucl}
	Let $ u,v \in \jormod $ and $ a,b \in \ring $. By the same computation as in \cref{pseud:sesqui-nuclear-note}, it follows from \thmitemcref{BC:blue:eq7lem}{BC:blue:eq7lem:psi-jor} and \thmitemcref{BC:blue:eq5lem}{BC:blue:eq5lem:psi-skew-cor} that
	\[ \brackets[\big]{a \psi(u,v)}b = a \brackets[\big]{\psi(u,v)b}. \]
	In other words, the image of $ \psi $ lies in the middle nucleus of $ \ring $.
\end{remark}

\begin{lemma}[Equation~5, part 2]\label{BC:blue:eq5lem2}
	The following hold for all $ u \in \jormod $:
	\begin{lemenumerate}
		\item \label{BC:blue:eq5lem2:pi-inv}$ \rinv{\jorprojone(u)} = \jorproj(u) - \psi(u,u) $.
		
		\item \label{BC:blue:eq5lem2:pi-mnucl}$ \jorprojone(u) $ lies in the middle nucleus of $ \ring $.
	\end{lemenumerate}
\end{lemma}
\begin{proof}
	We consider equation~5 in \cref{fig:BC:blue-long} and put all variables except for $ u $, $ r $ and $ s $ to zero. This yields
	\begin{align*}
		-\brackets[\big]{-\jorproj\brackets{\overline{u},r}-\psi\brackets[\big]{\jorsc\brackets{\overline{u},-\rinv{r}},\overline{u}}}\rinv{s} &= \rinv{\brackets[\bigg]{\brackets[\big]{-\jorproj\brackets{\overline{u},s}-\psi\brackets[\big]{\jorsc\brackets{\overline{u},-\rinv{s}},\overline{u}}}\brackets{-\rinv{r}}}} \\
		&\hspace{1cm} \mathord{}-\psi\brackets[\big]{\jorsc\brackets{u,\rinv{r}},\jorsc\brackets{\overline{u},-\rinv{s}}}.
	\end{align*}
	We can simplify this equation to obtain
	\begin{align*}
		\jorproj(u,r)\rinv{s} -r \psi(u,u) \rinv{s} &= r \rinv{\jorproj(u,s)} + r \psi(u,u)\rinv{s} - r \psi(u,u) \rinv{s},
	\end{align*}
	or in other words,
	\[ \jorproj(u,r) \rinv{s} - r \psi(u,u) \rinv{s} = r \rinv{\jorproj(u,s)}. \]
	For $ r \defl s \defl 1_\ring $, this is exactly assertion~\itemref{BC:blue:eq5lem2:pi-inv}. Using~\itemref{BC:blue:eq5lem2:pi-inv} and \thmitemcref{BC:blue:eq7lem}{BC:blue:eq7lem:pi}, the equation that we have just established simplifies as follows:
	\begin{align*}
		\brackets[\big]{r \jorprojone(u)} \rinv{s} - r \psi(u,u) \rinv{s} &= r \rinv{\brackets[\big]{s \jorprojone(u)}} = r \brackets[\big]{\jorprojone(u) \rinv{s}} - r \psi(u,u) \rinv{s}.
	\end{align*}
	This implies that $ \brackets[\big]{r \jorprojone(u)} \rinv{s} = r \brackets[\big]{\jorprojone(u) \rinv{s}} $. In other words, $ \jorprojone(u) $ lies in the middle nucleus, which finishes the proof of~\itemref{BC:blue:eq5lem2:pi-mnucl}.
\end{proof}

\begin{remark}\label{BC:blue:pi-inv}
	\thmitemcref{BC:blue:eq7lem}{BC:blue:eq7lem:pi} together with \thmitemcref{BC:blue:eq5lem2}{BC:blue:eq5lem2:pi-inv} yields that
	\begin{align*}
		\rinv{\jorproj(u,a)} &= \rinv{\brackets[\big]{a \jorprojone(u)}} = \jorprojone(u)\rinv{a} - \psi(u,u)\rinv{a}
	\end{align*}
	for all $ u \in \jormod $ and $ a \in \ring $.
\end{remark}

Recall the definition of nuclei from \cref{ring:nucleus-def}.

\begin{lemma}\label{BC:blue:weak-alt}
	The ring $ \ring $ is weakly alternative and satisfies
	\[ \nucleus(\ring) = \lnucleus(\ring) = \mnucleus(\ring) = \rnucleus(\ring). \]
	Further, the images of $ \jorprojone $ and $ \psi $ lie in $ \nucleus(\ring) $ and the involution $ \rinvmap $ is nuclear.
\end{lemma}
\begin{proof}
	For any $ a \in \ring $, we have $ \jorprojone(\jorTr(1_\ring, a)) = a + \rinv{a} $ by \thmitemcref{BC:blue:eq5lem}{BC:blue:eq5lem:pi-Tr}. In particular, $ a+\rinv{a} $ lies in the image of $ \jorprojone $, so it lies in the middle nucleus by \thmitemcref{BC:blue:eq5lem2}{BC:blue:eq5lem2:pi-mnucl}. By the same computation as in  \cref{rinv:invo-assoc}, this implies that
	\begin{equation}\label{eq:BC:blue:weak-alt-inv}
		\assoc{d}{\rinv{a}}{s} = -\assoc{d}{a}{s} \quad \text{for all } a,d,s \in \ring.
	\end{equation}
	Now let $ a,d,s \in \ring $ be arbitrary. On the one hand, \thmitemcref{BC:blue:eq5lem}{BC:blue:eq5lem:pi-Tr} yields that
	\[ \jorproj\brackets[\big]{\jorTr(s,a), d} = (ds) \rinv{a} (da) \rinv{s}. \]
	On the other hand, it follows from \thmitemcref{BC:blue:eq7lem}{BC:blue:eq7lem:pi} that
	\begin{align*}
		\jorproj\brackets[\big]{\jorTr(s,a), d} &= d \jorproj\brackets[\big]{\jorTr(s,a), 1_\ring} = d \brackets[\big]{s \rinv{a} + a \rinv{s}}.
	\end{align*}
	Together, these two identities imply that $ \assoc{d}{s}{\rinv{a}} = -\assoc{d}{a}{\rinv{s}} $. Invoking~\eqref{eq:BC:blue:weak-alt-inv}, we infer that $ -\assoc{d}{\rinv{s}}{\rinv{a}} = \assoc{d}{\rinv{a}}{\rinv{s}} $. This says precisely that the conditions of \thmitemcref{rinv:inv-assoc}{rinv:inv-assoc:weak} are satisfied, so $ \ring $ is weakly alternative.
	
	By \cref{ring:alternative-nucleus}, it follows from the weak alternativity of $ \ring $ that nucleus, left nucleus, middle nucleus and right nucleus coincide. Since we have already seen in \cref{BC:blue:psi-mnucl}, \thmitemcref{BC:blue:eq5lem2}{BC:blue:eq5lem2:pi-mnucl} and~\eqref{eq:BC:blue:weak-alt-inv} that the images of $ \jorprojone $ and $ \psi $ and all traces of $ \ring $ lie in the middle nucleus, the remaining assertions follow.
\end{proof}

\begin{lemma}[Equation 5, part 3]\label{BC:blue:eq5lem3}
	We have $ \jorproj\brackets[\big]{\jorsc(u, s), d} = (d\rinv{s})\jorprojone(u) s $ for all $ u \in \jormod $ and $ d,s \in \ring $.
\end{lemma}
\begin{proof}
	We put all variables in equation~5 except for $ u $, $ d $ and $ s $ to zero. This yields
	\begin{align*}
		0 &= -\jorproj\brackets[\big]{\jorsc\brackets{\overline{u},-\rinv{s}},d}+\rinv{\brackets[\bigg]{\brackets[\big]{-\jorproj\brackets{\overline{u},s}-\psi\brackets{\jorsc\brackets{\overline{u},-\rinv{s}},\overline{u}}}\brackets[\big]{-\rinv{\brackets{ds}}}}} \\
		& \hspace{1cm} \mathord{}-\psi\brackets[\bigg]{\jorsc\brackets[\big]{\jorsc\brackets{\overline{u},-\rinv{s}},-\rinv{d}}\joradd \jorsc\brackets[\big]{u,\rinv{\brackets{ds}}},\jorsc\brackets{\overline{u},-\rinv{s}}}.
	\end{align*}
	Using \cref{BC:blue:eq6lem,BC:blue:eq8lem,BC:blue:eq7lem} and also~\thmitemcref{BC:blue:eq5lem}{BC:blue:eq5lem:psi-skew}, we observe that
	\begin{align*}
		\rinv{\brackets[\bigg]{\brackets[\big]{-\jorproj\brackets{\overline{u},s}-\psi\brackets{\jorsc\brackets{\overline{u},-\rinv{s}},\overline{u}}}\brackets[\big]{-\rinv{\brackets{ds}}}}} \hspace{-6cm} & \\
		&= -(ds) \brackets[\bigg]{-\rinv{\jorproj(u,s)} - \rinv{\psi\brackets[\big]{\jorsc(u, -\rinv{s}), u}}} \\
		&= (ds) \brackets[\bigg]{\rinv{\jorproj(u,s)} - \psi\brackets[\big]{u, \jorsc(u, -\rinv{s})}} \\
		&= (ds) \brackets[\big]{\rinv{\jorproj(u,s)} + \psi\brackets[\big]{u, u} \rinv{s}} = (ds) \rinv{\jorproj(u,s)} + (ds) \psi(u,u) \rinv{s}
	\end{align*}
	and
	\begin{align*}
		\psi\brackets[\bigg]{\jorsc\brackets[\big]{\jorsc\brackets{\overline{u},-\rinv{s}},-\rinv{d}}\joradd \jorsc\brackets[\big]{u,\rinv{\brackets{ds}}},\jorsc\brackets{\overline{u},-\rinv{s}}} \hspace{-7cm}& \\
		&= \psi\brackets[\bigg]{\jorsc\brackets[\big]{\jorsc\brackets{\overline{u},-\rinv{s}},-\rinv{d}},\jorsc\brackets{\overline{u},-\rinv{s}}} + \psi\brackets[\bigg]{\jorsc\brackets[\big]{u,\rinv{\brackets{ds}}},\jorsc\brackets{\overline{u},-\rinv{s}}} \\
		&= d \psi\brackets[\big]{\jorsc(\modinv{u}, -\rinv{s}), \modinv{u}} \rinv{s} - (ds) \psi(u, \modinv{u}) \rinv{s} \\
		&= -\brackets[\big]{ds \psi(u,u)} \rinv{s} + (ds) \psi(u,u) \rinv{s}.
	\end{align*}
	Therefore,
	\begin{align*}
		0 &= -\jorproj\brackets[\big]{\jorsc(u, -\rinv{s}), d} + (ds) \rinv{\jorproj(u,s)} + (ds) \psi(u,u) \rinv{s} \\
		& \hspace{3cm} \mathord{}+ \brackets[\big]{d s \psi(u,u)} \rinv{s} - (ds) \psi(u,u) \rinv{s} \\
		&= -\jorproj\brackets[\big]{\jorsc(u, -\rinv{s}), d} + (ds) \rinv{\jorproj(u,s)} + \brackets[\big]{d s \psi(u,u)} \rinv{s}.
	\end{align*}
	Since
	\begin{align*}
		(ds) \rinv{\jorproj(u,s)} &= (ds) \jorprojone(u) \rinv{s} - (ds) \psi(u,u) \rinv{s} \rightand \\
		(ds) \psi(u,u) \rinv{s} &= \brackets[\big]{d s \psi(u,u)} \rinv{s}
	\end{align*}
	by \cref{BC:blue:pi-inv} and \cref{BC:blue:weak-alt}, respectively, we infer that
	\[ \jorproj\brackets[\big]{\jorsc(u, -\rinv{s}), d} = (ds) \jorprojone(u) \rinv{s}. \]
	Replacing $ s $ by $ -\rinv{s} $, the assertion follows.
\end{proof}

\begin{proposition}\label{BC:blue:alternative}
	The ring $ \ring $ is alternative.
\end{proposition}
\begin{proof}
	Let $ s,d \in \ring $ and recall that the fixed element $ v_0 \in \jormod $ from \cref{BC:blue-secnotation} satisfies $ \jorprojone(v_0) = 1_\ring $ by \cref{BC:blue:rank2:inv}. Thus on the one hand, \cref{BC:blue:eq5lem3} yields that
	\begin{align*}
		\jorproj\brackets[\big]{\jorsc(v_0, s), d} &= (d \rinv{s}) \jorprojone(v_0)s = (d\rinv{s}) s.
	\end{align*}
	On the other hand, and together with \thmitemcref{BC:blue:eq7lem}{BC:blue:eq7lem:pi}, \cref{BC:blue:eq5lem3} implies that
	\begin{align*}
		\jorproj\brackets[\big]{\jorsc(v_0, s), d} &= d \brackets[\big]{\jorproj\brackets[\big]{\jorsc(v_0, s), 1_\ring}} = d \brackets[\big]{\rinv{s} \jorprojone(v_0) s} = d(\rinv{s}s).
	\end{align*}
	Thus $ \assoc{d}{\rinv{s}}{s} = 0 $. Since the involution $ \rinvmap $ is nuclear by \cref{BC:blue:weak-alt}, it follows that $ \assoc{d}{s}{s} = 0 $ for all $ s,d \in \ring $. By an application of \thmitemcref{rinv:inv-assoc}{rinv:inv-assoc:alt}, we conclude that $ \ring $ is alternative.
\end{proof}

\begin{remark}\label{BC:blue:eq5-finished}
	It is not clear that no further identities can be deduced from equation~5 in \cref{fig:BC:blue-long}. We verify this by simplifying this equation with our recently obtained knowledge.
	For the simplification of $ y_5 $, we begin with the first summand. Repeatedly applying \thmitemcref{BC:param:rank2commrel}{BC:param:rank2commrel:proj} as well as the identities from this section and \thmitemcref{BC:param:rank2commrel}{BC:param:rank2commrel:cent}, we see that
	\begin{align*}
		\jorproj\brackets[\big]{\jorsc\brackets{\modinv{u},-\rinv{s}}\joradd \modinv{v}\joradd \jorTr\brackets{s,\rinv{a}},d} \hspace{-4.5cm} & \\
		&= \jorproj\brackets[\big]{\jorsc\brackets{\modinv{u},-\rinv{s}}\joradd \modinv{v},d} + \jorproj\brackets[\big]{\jorTr\brackets{s,\rinv{a}}, d} + \psi\brackets[\big]{\jorsc\brackets{\jorsc\brackets{\modinv{u},-\rinv{s}}\joradd \modinv{v}, \rinv{d}}, \jorTr\brackets{s,\rinv{a}}} \\
		&= \jorproj\brackets[\big]{\jorsc\brackets{\modinv{u},-\rinv{s}}\joradd \modinv{v},d} + \jorproj\brackets[\big]{\jorTr\brackets{s,\rinv{a}}, d} \\
		&= \jorproj\brackets[\big]{\jorsc\brackets{\modinv{u},-\rinv{s}},d} + \jorproj\brackets[\big]{\modinv{v},d} + \psi\brackets[\big]{\jorsc\brackets[\big]{\jorsc\brackets{\modinv{u},-\rinv{s}}, \rinv{d}}, \modinv{v}} + \jorproj\brackets[\big]{\jorTr\brackets{s,\rinv{a}}, d} \\
		&= (ds) \jorprojone(u) \rinv{s} + d \jorprojone(v) - ds \psi(u,v) + (ds)a + (d\rinv{a})\rinv{s}.
	\end{align*}
	Further,
	\begin{align*}
		\rinv{\sqbrackets[\big]{\brackets[\big]{-\jorproj\brackets{\modinv{u},s}-\rinv{a}-\psi\brackets[\big]{\jorsc\brackets{\modinv{u},-\rinv{s}}\joradd \modinv{v},\modinv{u}}}\brackets[\big]{-\rinv{r}-\rinv{\brackets{ds}}}}} \hspace{-8cm} & \\
		&= \brackets{r+ds} \brackets[\big]{\rinv{\jorproj\brackets{u,s}}+a+\rinv{\psi\brackets[\big]{\jorsc\brackets{u,-\rinv{s}}\joradd v,u}}} \\
		&= (r+ds) \brackets[\big]{\jorprojone(u) \rinv{s} - \psi(u,u) \rinv{s} + a - \psi\brackets[\big]{u, \jorsc(u, -\rinv{s}) \joradd v}} \\
		&= (r+ds) \brackets[\big]{\jorprojone(u) \rinv{s} - \psi(u,u) \rinv{s} + a + \psi(u,u) \rinv{s} - \psi(u,v)} \\
		&= (r+ds) \brackets[\big]{\jorprojone(u) \rinv{s} + a - \psi(u,v)}.
	\end{align*}
	Finally,
	\begin{align*}
		-\psi\brackets[\big]{C \joradd \modinv{w}\joradd \jorTr\brackets{-r-ds,\rinv{b}},\jorsc\brackets{\modinv{u},-\rinv{s}}\joradd \modinv{v}\joradd \jorTr\brackets{s,\rinv{a}}} \hspace{-8cm} & \\
		&= -\psi\brackets[\big]{C \joradd \modinv{w},\jorsc\brackets{\modinv{u},-\rinv{s}}\joradd \modinv{v}} \\
		&= -\psi(C, u) \rinv{s} + \psi(C,v) + \psi(w,u) \rinv{s} - \psi(w,v)
	\end{align*}
	where, for any $ x \in \jormod $,
	\begin{align*}
		\psi(C,x) &= \psi\brackets[\bigg]{\jorsc\brackets[\big]{\jorsc\brackets{\modinv{u},-\rinv{s}}\joradd \modinv{v}\joradd \jorTr\brackets{s,\rinv{a}},-\rinv{d}} \joradd \jorsc\brackets[\big]{u,\rinv{r}+\rinv{\brackets{ds}}}, x} \\
		&= -d\psi\brackets[\big]{\jorsc\brackets{\modinv{u},-\rinv{s}}\joradd \modinv{v}\joradd \jorTr\brackets{s,\rinv{a}}, x} + (r+ds) \psi(u,x) \\
		&= -ds \psi(u,x) + d\psi(v,x) + (r+ds) \psi(u,x).
	\end{align*}
	Altogether, we conclude that
	\begin{align*}
		y_5 &= -(ds) \jorprojone(u) \rinv{s} - d \jorprojone(v) + ds \psi(u,v) - (ds)a - (d\rinv{a})\rinv{s} - \rinv{c} - \rinv{b} \rinv{s} \\
		& \hspace{1cm} \mathord{} + r \brackets[\big]{\jorprojone(u) \rinv{s} + a - \psi(u,v)} + (ds) \brackets[\big]{\jorprojone(u) \rinv{s} + a - \psi(u,v)} \\
		&\hspace{1cm} \mathord{} +(ds) \psi(u,u) \rinv{s} - d\psi(v,u) \rinv{s} - (r+ds) \psi(u,u) \rinv{s} \\
		&\hspace{1cm} \mathord{}-ds \psi(u,v) + d\psi(v,v) + (r+ds) \psi(u,v) + \psi(w,u) \rinv{s} - \psi(w,v).
	\end{align*}
	Note that all the summands $ (ds) \jorprojone(u) \rinv{s} $, $ (ds)a $, $ r\psi(u,v) $ and $ (ds) \psi(u,u) \rinv{s} $ appear once with a positive and once with a negative sign in the expression of $ y_5 $, and $ (ds) \psi(u,v) $ even appears twice with a positive and twice with a negative sign. Thus we have
	\begin{align*}
		y_5 &=  - d \jorprojone(v) - (d\rinv{a})\rinv{s} - \rinv{c} - \rinv{b} \rinv{s} + r \brackets[\big]{\jorprojone(u) \rinv{s} + a } \\
		& \hspace{1cm} \mathord{} - d\psi(v,u) \rinv{s} - r \psi(u,u) \rinv{s} + d\psi(v,v) + \psi(w,u) \rinv{s} - \psi(w,v).
	\end{align*}
	For the simplification of $ x_5 $, we first observe that
	\begin{align*}
		\psi(B,\modinv{u}) &= \psi\brackets[\big]{\jorsc\brackets{\modinv{u},-\rinv{r}}\joradd \jorsc\brackets{\modinv{v},-\rinv{d}}\joradd \modinv{w}\joradd \jorTr\brackets{d,\rinv{c}}, \modinv{u}} \\
		&= \psi\brackets[\big]{\jorsc\brackets{u,-\rinv{r}}\joradd \jorsc\brackets{v,-\rinv{d}}\joradd w, u} \\
		&= -r \psi(u,u) -d \psi(v,u) + \psi(w,u).
	\end{align*}
	Therefore,
	\begin{align*}
		x_5 &= -d \jorprojone(v) - \rinv{c} + d \psi(v,v) - \psi(w,v) + ra + r \jorprojone(u) \rinv{s} - \rinv{b} \rinv{s} - (d \rinv{a}) \rinv{s} \\
		& \hspace{1cm} \mathord{}-r \psi(u,u) \rinv{s} -d \psi(v,u) \rinv{s} + \psi(w,u) \rinv{s}.
	\end{align*}
	It is now easy to verify that the equation $ x_5 = y_5 $ is equivalent to $ 0=0 $. Thus there exist no other identities which can be deduced from equation~5.
\end{remark}

\begin{lemma}[Equation 4, part 1]\label{BC:blue:eq4lem}
	The following hold for all $ u \in \jormod $ and $ a,b,c \in \ring $:
	\begin{lemenumerate}
		\item \label{BC:blue:eq4lem:Tr-shift}$ \jorTr(ab,c) = \jorTr(a, c \rinv{b}) $.
		
		\item \label{BC:blue:eq4lem:sym}$ \jorTr(a,b) = \jorTr(b,a) $.
		
		\item \label{BC:blue:eq4lem:phi-min-inv}$ \jorsc(u, -a) = \modinv{\jorsc(u,a)} $.
		
		\item \label{BC:blue:eq4lem:phi-Tr}$ \jorsc\brackets[\big]{\jorTr(a, b), c} = \jorTr(\rinv{c}b, \rinv{c}a) $.
		
		\item \label{BC:blue:eq4lem:comm}$ u \joradd v = v \joradd u \joradd \jorTr\brackets[\big]{1_\ring, \psi(u,v)} $.
	\end{lemenumerate}
\end{lemma}
\begin{proof}
	We consider equation~4 in \cref{fig:BC:blue-long}. Putting all variables except for $ b $, $ d $ and $ s $ to zero, we obtain
	\[ 0_\jormod = \modinv{\jorTr\brackets{-ds,\rinv{b}}}\joradd \modinv{\jorTr\brackets{d,\rinv{\brackets{sb}}}}. \]
	In other words, $ \jorTr(ds, \rinv{b}) = \jorTr(d, \rinv{b} \rinv{s}) $. This proves~\itemref{BC:blue:eq4lem:Tr-shift}.
	
	Replacing all variables in equation~4 except for $ a $, $ d $ and $ r $ by zero, we see that
	\[ \modinv{\jorTr\brackets{r,-\rinv{\brackets{a\rinv{d}}}}} = \modinv{\jorTr\brackets{d,-\rinv{\brackets{\brackets{-\rinv{a}}\brackets{-\rinv{r}}}}}}. \]
	This is equivalent to the identity $ \jorTr(r, d \rinv{a}) = \jorTr(d, ra) $, which for $ a \defl 1_\ring $ implies that~\itemref{BC:blue:eq4lem:sym} holds.
	
	Now we substitute each variable in equation~4 except for $ u $ and $ r $ by zero. In this way, we obtain
	\[ \modinv{\jorsc\brackets{\modinv{u},-\rinv{r}}} = \modinv{\jorsc\brackets{u,\rinv{r}}}. \]
	Since
	\[ \modinv{\jorsc\brackets{\modinv{u},-\rinv{r}}} = \modinv{\modinv{\jorsc(u, -\rinv{r})}} \]
	by \cref{BC:blue:eq6lem}, assertion~\itemref{BC:blue:eq4lem:phi-min-inv} follows.
	
	Putting all variables in equation~4 except for $ a $, $ d $ and $ s $ to zero, we obtain
	\begin{align*}
		0_\jormod &= \modinv{\jorsc\brackets[\big]{\jorTr\brackets{s,\rinv{a}},-\rinv{d}}}\joradd \modinv{\jorTr\brackets[\big]{d,-\rinv{\brackets[\big]{\brackets{-\rinv{a}}\brackets{-\rinv{\brackets{ds}}}}}}}.
	\end{align*}
	Note that, by~\itemref{BC:blue:eq4lem:phi-min-inv} and \cref{BC:blue:eq6lem}, we have
	\begin{align*}
		\jorsc\brackets[\big]{\jorTr\brackets{s,\rinv{a}},-\rinv{d}} &= \modinv{\jorsc\brackets[\big]{\jorTr\brackets{s,\rinv{a}},\rinv{d}}} = \jorsc\brackets[\big]{\modinv{\jorTr(s, \rinv{a})}, \rinv{d}} = \jorsc\brackets[\big]{\jorTr(s, \rinv{a}), \rinv{d}}.
	\end{align*}
	Hence
	\begin{align*}
		\jorsc\brackets[\big]{\jorTr(s, \rinv{a}), \rinv{d}} &= \jorTr\brackets[\big]{d, (ds)a},
	\end{align*}
	which together with~\itemref{BC:blue:eq4lem:Tr-shift} shows that~\itemref{BC:blue:eq4lem:phi-Tr} holds.
	
	Next we put all variables in equation~4 except for $ u $, $ v $, $ d $ and $ r $ to zero. This yields
	\begin{align*}
		\modinv{\jorsc\brackets{\modinv{u},-\rinv{r}}}\joradd \modinv{\jorsc\brackets{\modinv{v},-\rinv{d}}} &= \modinv{\jorsc\brackets{\modinv{v},-\rinv{d}}}\joradd \modinv{\jorsc\brackets{u,\rinv{r}}}\joradd \modinv{\jorTr\brackets[\big]{d,-\rinv{\brackets{\brackets[\big]{-\psi\brackets{\modinv{v},\modinv{u}}}\brackets{-\rinv{r}}}}}}.
	\end{align*}
	Using \cref{BC:blue:eq6lem,BC:blue:eq8lem}, we can remove most occurrences of $ \modinv{\mapdot} $ in this equation to obtain
	\[ \jorsc\brackets{u,-\rinv{r}} \joradd \jorsc\brackets{v,-\rinv{d}} = \jorsc\brackets{v,-\rinv{d}} \joradd \modinv{\jorsc\brackets{u,\rinv{r}}} \joradd \jorTr\brackets{d,-r \rinv{\psi(v,u)}}. \]
	Putting $ d \defl r \defl -1_\ring $ and using~\itemref{BC:blue:eq4lem:phi-min-inv} and \cref{BC:blue:rank2:inv}, we infer that
	\[ u \joradd v = v \joradd u \joradd \jorTr\brackets[\big]{1_\ring, \psi(u,v)}, \]
	which is exactly the assertion of~\itemref{BC:blue:eq4lem:comm}.
\end{proof}

Motivated by \thmitemcref{BC:blue:eq4lem}{BC:blue:eq4lem:Tr-shift} and~\thmitemref{BC:blue:eq4lem}{BC:blue:eq4lem:sym}, we make the following definition, which is similar to \cref{BC:blue:pi-one-def}.

\begin{definition}\label{BC:blue:Tr-one-def}
	We define a map $ \map{\jorTrone}{\ring}{\jormod}{a}{\jorTr(a,1_\ring) = \jorTr(1_\ring, a)} $.
\end{definition}

\begin{remark}
	Let $ a,b,c\in \ring $. Then \cref{BC:blue:eq4lem} yields the following properties of the map $ \jorTrone $:
	\begin{align*}
		\jorsc\brackets[\big]{\jorTrone(a), c} &= \jorTr(\rinv{c}, \rinv{c} a) = \jorTr(1_\ring, \rinv{c} a c) = \jorTrone(\rinv{c}ac), \\
		\jorTrone(a) &= \jorTr(1_\ring a, 1_\ring) = \jorTr(1_\ring, 1_\ring \rinv{a}) = \jorTrone(\rinv{a}), \\
		\jorTrone\brackets[\big]{(ab)c} &= \jorTr\brackets[\big]{(ab)c, 1_\ring} = \jorTr(ab, \rinv{c}) = \jorTr(a, \rinv{c} \rinv{b}) = \jorTr\brackets[\big]{1_\ring, (\rinv{c} \rinv{b}) \rinv{a}} \\
		&= \jorTrone\brackets[\big]{(\rinv{c} \rinv{b}) \rinv{a}} = \jorTrone\brackets[\big]{a (bc)}.
	\end{align*}
	Further, \thmitemcref{BC:blue:eq5lem}{BC:blue:eq5lem:pi-Tr} implies that
	\begin{align*}
		\jorprojone\brackets[\big]{\jorTrone(a)} &= \jorproj\brackets[\big]{\jorTr(a, 1_\ring), 1_\ring} = a + \rinv{a}.
	\end{align*}
	This last property is in some sense dual to the following \thmitemcref{BC:blue:eq4lem2}{BC:blue:eq4lem2:Tr-pi}.
\end{remark}

\begin{lemma}[Equation 4, part 2]\label{BC:blue:eq4lem2}
	The following identities hold for all $ u \in \jormod $ and all $ s,d \in \ring $:
	\begin{lemenumerate}
		\item \label{BC:blue:eq4lem2:Tr-pi}$ \jorTrone\brackets[\big]{\jorprojone(u)} = u + \modinv{u} = \modinv{u} + u $.
		
		\item \label{BC:blue:eq4lem2:phi-mult}$ \jorsc\brackets[\big]{\jorsc(u, s), d} = \jorsc(u, sd) $.
	\end{lemenumerate}
\end{lemma}
\begin{proof}
	In this proof, we will only consider equation~4 in \cref{fig:BC:blue-long} with all variables except for $ u $, $ d $ and $ s $ replaced by zero. This produces the following identity:
	\begin{align*}
		0 &= \modinv{\jorsc\brackets[\big]{\jorsc\brackets{\modinv{u},-\rinv{s}},-\rinv{d}}}\joradd \modinv{\jorsc\brackets[\big]{u,\rinv{\brackets{ds}}}} \\
		&\hspace{1cm} \mathord{}\joradd \modinv{\jorTr\brackets[\bigg]{d,-\rinv{\sqbrackets[\big]{\brackets[\big]{-\jorproj\brackets{\modinv{u},s}-\psi\brackets[\big]{\jorsc\brackets{\modinv{u},-\rinv{s}},\modinv{u}}}\brackets[\big]{-\rinv{\brackets{ds}}}}}}}
	\end{align*}
	Using the previously computed identities, we can simplify this equation to obtain
	\begin{align*}
		0 &= \jorsc\brackets[\big]{\jorsc(u, -\rinv{s}), -\rinv{d}} \joradd \modinv{\jorsc\brackets[\big]{u, \rinv{s} \rinv{d}}} \joradd \jorTr\brackets[\bigg]{d, (ds) \brackets[\big]{- \rinv{\jorproj(u, s)} + \rinv{\brackets[\big]{s \psi(u,u)}}}}
	\end{align*}
	where
	\begin{align*}
		\jorTr\brackets[\bigg]{d, (ds) \brackets[\big]{- \rinv{\jorproj(u, s)} + \rinv{\brackets[\big]{s \psi(u,u)}}}} \hspace{-4cm} & \\
		&= \jorTr\brackets[\bigg]{d, (ds) \brackets[\big]{ - \jorprojone(u) \rinv{s} + \psi(u,u) \rinv{s} - \psi(u,u) \rinv{s}}} \\
		&= \jorTr\brackets[\big]{1_\ring, -ds\jorprojone(u) \rinv{s} \rinv{d}} = \jorTrone(-ds\jorprojone(u) \rinv{s} \rinv{d}).
	\end{align*}
	Replacing $ s $ and $ d $ by $ -\rinv{s} $ and $ -\rinv{d} $, respectively, we conclude that
	\begin{equation}\label{eq:BC:blue:eq4}
		0 = \jorsc\brackets[\big]{\jorsc(u,s), d} \joradd \modinv{\jorsc(u, sd)} \joradd \jorTrone\brackets[\big]{-\rinv{d} \rinv{s} \jorprojone(u) sd}.
	\end{equation}
	Putting $ s \defl d \defl 1_\ring $, we obtain $ \jorTrone(\jorprojone(u)) = u \joradd \modinv{u} $. Since $ \jorprojone(u) = \jorprojone(\modinv{u}) $ by \cref{BC:blue:eq8lem}, it follows that $ \jorTrone(\jorprojone(u)) = \modinv{u} \joradd u $ holds as well. This finishes the proof of~\itemref{BC:blue:eq4lem2:Tr-pi}.
	
	Using~\itemref{BC:blue:eq4lem2:Tr-pi} and \cref{BC:blue:eq5lem3}, we see that
	\begin{align*}
		\jorTrone\brackets[\big]{-\rinv{d} \rinv{s} \jorprojone(u) sd} &= \jormin \jorTrone\brackets[\big]{\jorprojone\brackets[\big]{\jorsc(u, sd)}} = \jormin \brackets[\big]{\jorsc(u,sd) \joradd \modinv{\jorsc(u,sd)}} \\
		&= \brackets[\big]{\jormin \modinv{\jorsc(u,sd)}} \joradd \brackets[\big]{\jormin \jorsc(u,sd)}.
	\end{align*}
	With this information, \eqref{eq:BC:blue:eq4} becomes
	\begin{align*}
		0 = \jorsc\brackets[\big]{\jorsc(u,s), d} \joradd \brackets[\big]{\jormin \jorsc(u,sd)},
	\end{align*}
	which proves~\itemref{BC:blue:eq4lem2:phi-mult}.
\end{proof}

\begin{remark}
	Similarly as in \cref{BC:blue:eq5-finished}, we will now show that no further identities can be deduced from equation~4 in \cref{fig:BC:blue-long}. The right-hand side of this equation can be written as
	\begin{align*}
		y_4 &= \jorsc(u, \rinv{s} \rinv{d}) \joradd \jorsc(v, -\rinv{d}) \joradd \jorTrone(d sa \rinv{d}) \joradd \modinv{\jorsc(u, \rinv{r} + \rinv{s} \rinv{d})} \\
		& \hspace{1cm} \mathord{}\joradd w \joradd \jorTrone(-rb - dsb) \joradd \jorTrone(dc) \joradd \jorTrone(dsb) \\
		& \hspace{1cm} \mathord{}\joradd\jorTrone\brackets[\bigg]{d, (r+ds) \sqbrackets[\big]{- \rinv{\jorproj(u,s)} - a - \rinv{\psi\brackets[\big]{\jorsc(u, -\rinv{s}) \joradd v, u}}}}.
	\end{align*}
	By \thmitemcref{BC:param:rank2commrel}{BC:param:rank2commrel:sc},
	\begin{align*}
		\modinv{\jorsc(u, \rinv{r} + \rinv{s} \rinv{d})} &= \modinv{\jorsc(u, \rinv{r})} \joradd \modinv{\jorsc(u, \rinv{s} \rinv{d})} \joradd \jorTrone\brackets[\big]{ds \jorprojone(u) \rinv{r}}.
	\end{align*}
	Further,
	\begin{align*}
		\jorTrone\brackets[\bigg]{d, (r+ds) \sqbrackets[\big]{- \rinv{\jorproj(u,s)} - a - \rinv{\psi\brackets[\big]{\jorsc(u, -\rinv{s}) \joradd v, u}}}} \hspace{-9cm} & \\
		&= \jorTrone\brackets[\bigg]{d, (r+ds) \sqbrackets[\big]{- \rinv{\jorprojone(u)}\rinv{s} - a + \rinv{\brackets[\big]{s \psi(u,u)}} - \rinv{\psi(v,u)}}} \\
		&= \jorTrone\brackets[\bigg]{d, (r+ds) \sqbrackets[\big]{- \rinv{\jorprojone(u)}\rinv{s} - a - \psi(u,u) \rinv{s} + \psi(u,v)}} \\
		&= \jorTrone\brackets[\big]{- r \rinv{\jorprojone(u)} \rinv{s} \rinv{d} -ra\rinv{d} - r \psi(u,u) \rinv{s} \rinv{d} + r\psi(u,v) \rinv{d}} \\
		&\hspace{1cm} \mathord{} \joradd \jorTrone\brackets[\big]{-ds \rinv{\jorprojone(u)} \rinv{s} \rinv{d} -dsa \rinv{d} - ds \psi(u,u) \rinv{s} \rinv{d} + ds \psi(u,v) \rinv{d}}.
	\end{align*}
	We conclude that
	\begin{align*}
		y_4 &= \jorsc(u, \rinv{s} \rinv{d}) \joradd \jorsc(v, -\rinv{d}) \joradd \jorTrone(d sa \rinv{d}) \joradd \modinv{\jorsc(u, \rinv{r})} \joradd \modinv{\jorsc(u, \rinv{s} \rinv{d})} \\
		& \hspace{1cm} \mathord{}  \joradd \jorTrone\brackets[\big]{ds \jorprojone(u) \rinv{r}} \joradd w \joradd \jorTrone(-rb - dsb) \joradd \jorTrone(dc) \joradd \jorTrone(dsb) \\
		&\hspace{1cm} \mathord{} \jorTrone\brackets[\big]{- r \rinv{\jorprojone(u)} \rinv{s} \rinv{d} -ra\rinv{d} - r \psi(u,u) \rinv{s} \rinv{d} + r\psi(u,v) \rinv{d}} \\
		&\hspace{1cm} \mathord{} \joradd \jorTrone\brackets[\big]{-ds \rinv{\jorprojone(u)} \rinv{s} \rinv{d} -dsa \rinv{d} - ds \psi(u,u) \rinv{s} \rinv{d} + ds \psi(u,v) \rinv{d}}.
	\end{align*}
	Using that the image of $ \jorTrone $ lies in the center of $ \jormod $ by \thmitemcref{BC:param:rank2commrel}{BC:param:rank2commrel:cent}, that the summands $ \jorTrone(dsa\rinv{d}) $ and $ \jorTrone(dsb) $ appear once with a positive and once with a negative sign, and that
	\[ \jorTrone\brackets[\big]{- r \rinv{\jorprojone(u)} \rinv{s} \rinv{d}} = \jorTrone\brackets[\big]{- \rinv{\brackets[\big]{r \rinv{\jorprojone(u)} \rinv{s} \rinv{d}}}} = -\jorTrone\brackets[\big]{ds \jorprojone(u) \rinv{r}}, \]
	we can simplify this expression of $ y_4 $ as follows:
	\begin{align*}
		y_4 &= \jorsc(u, \rinv{s} \rinv{d}) \joradd \jorsc(v, -\rinv{d}) \joradd \modinv{\jorsc(u, \rinv{r})} \joradd \modinv{\jorsc(u, \rinv{s} \rinv{d})} \joradd w \joradd \jorTrone(-rb) \\
		&\hspace{1cm} \mathord{}   \joradd \jorTrone(dc) \joradd \jorTrone(-ra\rinv{d}) \joradd \jorTrone\brackets[\big]{- r \psi(u,u) \rinv{s} \rinv{d}} \joradd \jorTrone\brackets[\big]{r \psi(u,v) \rinv{d}} \\
		&\hspace{1cm} \mathord{}  \joradd \jorTrone\brackets[\big]{-ds\jorprojone(u) \rinv{s} \rinv{d}} \joradd \jorTrone\brackets[\big]{-ds \psi(u,u) \rinv{s} \rinv{d}} \joradd \jorTrone\brackets[\big]{ds \psi(u,v) \rinv{d}}.
	\end{align*}
	The left-hand side of equation~4 is easier to simplify:
	\begin{align*}
		x_4 &= \jorsc(u, -\rinv{r}) \joradd \jorsc(v, -\rinv{d}) \joradd w \joradd \jorTrone(dc) \joradd \jorTrone(-rb) \joradd \jorTrone(-ra \rinv{d}).
	\end{align*}
	Subtracting $ w \joradd \jorTrone(dc) \joradd \jorTrone(-rb) \joradd \jorTrone(-ra \rinv{d}) $ from the right side from both $ x_4 $ and $ y_4 $, we see that equation~4 is equivalent to $ x_4' = y_4' $ where
	\begin{align*}
		x_4' &= \jorsc(u, -\rinv{r}) \joradd \jorsc(v, -\rinv{d}), \\
		y_4' &= \jorsc(u, \rinv{s} \rinv{d}) \joradd \jorsc(v, -\rinv{d}) \joradd \modinv{\jorsc(u, \rinv{r})} \joradd \modinv{\jorsc(u, \rinv{s} \rinv{d})} \joradd \jorTrone\brackets[\big]{- r \psi(u,u) \rinv{s} \rinv{d}} \\
		&\hspace{1cm} \mathord{}   \joradd \jorTrone\brackets[\big]{r \psi(u,v) \rinv{d}} \joradd \jorTrone\brackets[\big]{-ds\jorprojone(u) \rinv{s} \rinv{d}} \\
		&\hspace{1cm} \mathord{}  \joradd \jorTrone\brackets[\big]{-ds \psi(u,u) \rinv{s} \rinv{d}} \joradd \jorTrone\brackets[\big]{ds \psi(u,v) \rinv{d}}.
	\end{align*}
	Observe that by \cref{BC:blue:eq5lem3} and \thmitemcref{BC:blue:eq4lem2}{BC:blue:eq4lem2:Tr-pi},
	\begin{align*}
		\jorTrone\brackets[\big]{-ds\jorprojone(u) \rinv{s} \rinv{d}} &= \jorTrone\brackets[\big]{-\jorprojone\brackets[\big]{\jorsc(u, \rinv{s} \rinv{d})}} = \brackets[\big]{\jormin \modinv{\jorsc(u, \rinv{s} \rinv{d})}} \joradd \brackets[\big]{\jormin \jorsc(u, \rinv{s} \rinv{d})}.
	\end{align*}
	Thus
	\begin{align*}
		y_4' &= \jorsc(u, \rinv{s} \rinv{d}) \joradd \jorsc(v, -\rinv{d}) \joradd \modinv{\jorsc(u, \rinv{r})} \joradd \jorsc(\jormin u, \rinv{s} \rinv{d}) \joradd \jorTrone\brackets[\big]{- r \psi(u,u) \rinv{s} \rinv{d}}  \\
		&\hspace{1cm} \mathord{} \joradd \jorTrone\brackets[\big]{r \psi(u,v) \rinv{d}} \joradd \jorTrone\brackets[\big]{-ds \psi(u,u) \rinv{s} \rinv{d}} \joradd \jorTrone\brackets[\big]{ds \psi(u,v) \rinv{d}}.
	\end{align*}
	Further, by \thmitemcref{BC:blue:eq4lem}{BC:blue:eq4lem:comm},
	\begin{align*}
		\jorsc(u, -\rinv{r}) \joradd \jorsc(v, -\rinv{d}) &= \jorsc(v, -\rinv{d}) \joradd \jorsc(u, -\rinv{r}) \joradd \jorTrone \brackets[\big]{\psi\brackets[\big]{\jorsc(u, -\rinv{r}), \jorsc(v, -\rinv{d})}} \\
		&= \jorsc(v, -\rinv{d}) \joradd \jorsc(u, -\rinv{r}) \joradd \jorTrone \brackets[\big]{r\psi\brackets{u, v} \rinv{d}}, \\
		\jorsc(u, \rinv{s} \rinv{d}) \joradd \jorsc(v, -\rinv{d}) &= \jorsc(v, -\rinv{d}) \joradd \jorsc(u, \rinv{s} \rinv{d}) \joradd \jorTrone\brackets[\big]{-ds \psi(u,v) \rinv{d}}.
	\end{align*}
	Therefore, subtracting $ \jorsc(v, -\rinv{d}) \joradd \jorTrone(r \psi(u,v) \rinv{d}) $ from the left side of the equation $ x_4' = y_4' $, we infer that equation~4 is equivalent to $ x_4'' = y_4'' $ where $ x_4'' = \jorsc(u, -\rinv{r}) $ and
	\begin{align*}
		y_4' &= \jorsc(u, \rinv{s} \rinv{d}) \joradd \modinv{\jorsc(u, \rinv{r})} \joradd \jorsc(\jormin u, \rinv{s} \rinv{d}) \joradd \jorTrone\brackets[\big]{- r \psi(u,u) \rinv{s} \rinv{d}}  \\
		&\hspace{1cm} \mathord{} \joradd \jorTrone\brackets[\big]{-ds \psi(u,u) \rinv{s} \rinv{d}}.
	\end{align*}
	Again by \thmitemcref{BC:blue:eq4lem}{BC:blue:eq4lem:comm}, we have
	\begin{align*}
		\jorsc(u, \rinv{s} \rinv{d}) \joradd \modinv{\jorsc(u, \rinv{r})} &= \modinv{\jorsc(u, \rinv{r})} \joradd \jorsc(u, \rinv{s} \rinv{d}) \joradd \jorTrone\brackets[\big]{ds\psi(u, \modinv{u}) \rinv{r}}
	\end{align*}
	where
	\begin{align*}
		\jorTrone\brackets[\big]{ds\psi(u, \modinv{u}) \rinv{r}} &= \jorTrone\brackets[\big]{- ds \psi(u,u) \rinv{r}} = \jorTrone\brackets[\big]{-r \rinv{\psi(u,u)} \rinv{s} \rinv{d}} \\
		&= \jorTrone\brackets[\big]{r \psi(u,u) \rinv{s} \rinv{d}}.
	\end{align*}
	Therefore, $ y_4'' = \modinv{\jorsc(u, \rinv{r})} $. Since $ \modinv{\jorsc(u, \rinv{r})} = \jorsc(u, -\rinv{r}) $ by \thmitemcref{BC:blue:eq4lem}{BC:blue:eq4lem:phi-min-inv}, it follows that equation~4 is equivalent to $ 0_\jormod=0_\jormod $. Thus our computations are finished.
\end{remark}

\begin{summary}\label{BC:blue:summary}
	It follows from the collection of identities in this section and the previous one that $ \ring $ is an alternative ring, $ \rinvmap $ is a nuclear involution on $ \ring $, $ (\module, \jorsc, \jorprojone, \jorTrone, \psi) $ is a Jordan module over $ (\ring, \ringzero, \rinvmap) $ and $ \modinvmap $ is the Jordan module involution on $ \module $. If $ n \ge 4 $, then we know in addition that $ \ring $ is associative (see \cref{BC:A3-assoc-rem}). Further, the computations in \cref{BC:jordanmodule-cor} show that every Jordan module satisfies all the identities that we have collected in this section and the previous one. In other words, Jordan modules are exactly the algebraic structures which are characterised by the results of the blueprint computations. It remains to compute the commutator relations which are not covered by \cref{BC:comm-mult-computation,BC:comm-formula-firststep}, which is a straightforward task.
\end{summary}

\begin{proposition}\label{BC:blue:stand-signs}
	$ (\risom{\alpha})_{\alpha \in BC_n} $ is a coordinatisation of $ G $ by $ (\module, \jorsc, \jorprojone, \jorTrone, \psi) $ with standard signs.
\end{proposition}
\begin{proof}
	We already know from \cref{BC:comm-mult-computation,BC:comm-formula-firststep} that the first set of the standard commutator relations in \cref{BC:standard-param-def} holds and that the first relation in each of the other sets holds if $ i<j $. We begin with the second set of relations. We already know that
	\[ \commutator{\rismin{i}{j}(a)}{\risplus{i}{j}(b)} = \risshpos{i}\brackets[\big]{\jorTrone(a \rinv{b})} \]
	holds for $ i<j $. Conjugating this identity by $ w_{ij} $, we obtain that
	\[ \commutator{\rismin{j}{i}(-a)}{\risplus{i}{j}(-\rinv{b})} = \risshpos{j}\brackets[\big]{\jorTrone(a\rinv{b})}. \]
	It follows that the first relation in the second set is true. Conjugating this identity by $ w_i $, we infer that
	\[ \commutator{\risminmin{i}{j}(\delinv{i<j}(a))}{\rismin{j}{i}(-\delinv{i<j}(b))} = \risshneg{i}\brackets[\big]{\jorTrone(a \delinv{i<j}(b))}, \]
	which proves the second relation.
	
	Now we turn to the third set of equations. We already know that
	\[ \commutator{\risshpos{j}(v)}{\rismin{i}{j}(a)} = \risshpos{i}\brackets[\big]{\jorsc(v, -\rinv{a})} \risplus{i}{j}\brackets[\big]{-\jorproj(v,a)}. \]
	It follows from this equation by conjugation with $ w_{ji} $ that
	\[ \commutator{\risshpos{i}(v)}{\rismin{j}{i}(-a)} = \risshpos{j}\brackets[\big]{\jorsc(v, \rinv{a})} \risplus{i}{j}\brackets[\big]{\rinv{\brackets[\big]{a \jorprojone(v)}}}. \]
	We infer that the first equation holds. Conjugating it by $ w_i $, we see that
	\[ \commutator{\risshpos{j}(v)}{\risminmin{i}{j}\brackets[\big]{\delinv{i<j}(a)}} = \risshneg{i}\brackets[\big]{\jorsc(v, -\rinv{a})} \rismin{j}{i}\brackets[\big]{\rinv{\jorprojone(v)} \rinv{a}}. \]
	The second equation follows. Conjugating the first equation by $ w_j $, we compute that
	\[ \commutator[\big]{\risshneg{j}(v)}{\risplus{i}{j}\brackets[\big]{\delinv{i>j}(a)}} = \risshpos{i}\brackets[\big]{\jorsc(v, -\rinv{a})} \rismin{i}{j}\brackets[\big]{a \jorprojone(v)}. \]
	The fourth equation follows from this one by interchanging the roles of $ i $ and $ j $. Conjugating the fourth equation by $ w_j $, we obtain
	\[ \commutator[\big]{\risshneg{i}(v)}{\rismin{i}{j}\brackets[\big]{-\delinv{i>j}(a)}} = \risshneg{j}\brackets[\big]{\jorsc(v, -\delinv{i>j}(a))} \risminmin{i}{j}\brackets[\big]{\delinv{i>j}\brackets[\big]{\delinv{i<j}(a) \jorprojone(u)}}. \]
	This proves the third equation.
	
	Finally, we turn to the last set of equations. Once again, we already know that
	\[ \commutator{\risshpos{i}(u)}{\risshpos{j}(v)} = \risplus{i}{j}\brackets[\big]{\psi(u,v)}. \]
	Conjugating this identity by $ w_{ij} $, we obtain that
	\[ \commutator{\risshpos{i}(u)}{\risshpos{j}(\modinv{v})} = \risplus{i}{j}\brackets[\big]{-\rinv{\psi(u,v)}} = \risplus{i}{j}\brackets[\big]{\psi(v,u)}. \]
	Thus the first relation holds. Conjugating this relation by $ w_j $, it follows that
	\[ \commutator{\risshpos{i}(u)}{\risshneg{j}(v)} = \begin{cases}
		\rismin{i}{j}\brackets[\big]{-\psi(u,v)} & \text{if } i<j \\
		\rismin{i}{j}\brackets[\big]{\rinv{\psi(v,u)}} = \rismin{i}{j}\brackets[\big]{-\psi(u,v)} & \text{if } i>j,
	\end{cases} \]
	proving the second identity. Finally, the last identity follows from the second one by conjugation with $ w_i $.
\end{proof}

\begin{lemma}\label{BC:v-formula}
	We have $ v_{-1} = v_1 = \modinv{v_0} $.
\end{lemma}
\begin{proof}
	By \cref{BC:blue:stand-signs}, we can apply \cref{BC:long-weakly-balanced}. This yields that
	\[ v_{-1} = v_1 = \jorsc\brackets[\big]{\modinv{v_0}, \rinvmin{\jorprojone(v_0)}}. \]
	Since $ \jorprojone(v_0) = 1_\ring $ by \cref{BC:blue:rank2:inv}, the assertion follows.
\end{proof}

We can now state the main result of this chapter. Recall the near-classification results for Jordan modules in \cref{BC:jordanmodule:class-general,BC:jordanmodule:class-2inv,BC:jordanmodule:class-2inv-C}.

\begin{theorem}[Coordinatisation theorem for $ BC_n $]\label{BC:thm}
	Let $ n \in \IN_{\ge 3} $ and let $ G $ be a group with a crystallographic $ BC_n $-grading $ (\rootgr{\alpha})_{\alpha \in BC_n} $. Then there exist an alternative ring $ \ring $ with nuclear involution $ \rinvmap $ and a Jordan module $ \jormodtup = (\jormod, \jorsc, \jorprojone, \jorTrone, \psi) $ over $ (\ring, \rinvmap) $ such that $ G $ is coordinatised by $ \jormodtup $ with standard signs. The ring $ \ring $ must be associative if $ n \ge 4 $ and the Jordan module $ \jormodtup $ must be of type $ C $ if $ (\rootgr{\alpha})_{\alpha \in BC_n} $ arises from a crystallographic $ C_n $-grading as in \cref{BC:CisBC:CasBC}. Further, if we fix a $ \rootbase $-system of Weyl elements in $ G $, then we can choose the root isomorphisms $ (\risom{\alpha})_{\alpha \in BC_n} $ so that $ w_\delta = \risom{-\delta}(-1_\ring) \risom{\delta}(1_\ring) \risom{-\delta}(-1_\ring) $ for all long simple roots $ \delta $ and $ w_\delta = \risom{-\delta}(\modinv{v_0}) \risom{\delta}(v_0) \risom{-\delta}(\modinv{v_0}) $ for the short simple root $ \delta $ where $ v_0 $ is some element of $ \jormod $ with $ \jorprojone(v_0) = 1 $ and $ \compinvmap $ denotes the Jordan module involution.
\end{theorem}
\begin{proof}[Summary of the proof]
	Choose the standard rescaled root base $ \rootbase $ of $ BC_n $ and a $ \rootbase $-system $ (w_\delta)_{\delta \in \rootbase} $ of Weyl elements. Denote by $ (\twistgroup, \inverparsym, \invogroup, \invoparsym) $ the standard admissible partial twisting system for $ G $, as in \cref{BC:standard-partwist-def}. Then by \cref{BC:param-exists}, there exist groups $ (\ring, +) $ and $ (\jormod, \joradd) $ on which $ \twistgroup \times \invogroup $ acts and a parametrisation $ (\risom{\alpha})_{\alpha \in BC_n} $ of $ G $ by $ (\twistgroup \times \invogroup, \jormod, \ring) $ with respect to $ (w_\delta)_{\delta \in \rootbase} $ and $ \inverparsym \times \invoparsym $. We can define commutation maps as in \cref{BC:commmap-def}. By \cref{BC:blue:summary}, these maps equip $ \ring $ with the structure of an alternative ring with nuclear involution (which must be associative if $ n \ge 4 $) and the group $ \jormod $ with the structure of a Jordan module $ \jormodtup $. By \cref{BC:blue:stand-signs}, $ (\risom{\alpha})_{\alpha \in BC_n} $ is a coordinatisation of $ G $ by $ \jormodtup $ with standard signs. By \cref{BC:isring,BC:v-formula}, the Weyl elements $ (w_\delta)_{\delta \in \rootbase} $ have the desired form. If $ (\rootgr{\alpha})_{\alpha \in BC_n} $ arises from a crystallographic $ C_n $-grading, then pairs of orthogonal short root groups commute (see \cref{BC:CisBC:BCasC}). This says precisely that $ \psi = 0 $, or in other words, that $ \jormodtup $ is of type $ C $.
\end{proof}

	\chapter{Root Gradings of Type \texorpdfstring{$ F $}{F}}
	
	\label{chap:F}
	
	In this final chapter, we investigate $ F_4 $-graded groups. Our goal is to show that each such group is coordinatised by a multiplicative conic alternative algebra $ \compalg $ over a commutative associative ring $ \comring $. These objects should be thought of as suitable generalisations of composition algebras. If the base ring $ \comring $ is a field of characteristic not~2 and the norm function on $ \compalg $ is anisotropic, then $ \compalg $ is a composition division algebra in the classical sense (\cref{conic:rgd-alg}). These assumptions are always satisfied for multiplicative conic alternative algebras which coordinatise RGD-systems of type $ F_4 $.
	
	While an explicit standard representation for the root system $ F_4 $ is available, it is much less practical to work with than the standard representations for root systems of types $ A $, $ B $, $ C $, $ BC $ and $ D $. This makes it much more difficult to concisely write down all commutator relations in $ F_4 $-graded groups. For this reason, we will only specify the commutator relations which can be seen in some (standard) root base. All other pairs of non-proportional roots adhere to commutator relations of the same form, but with twisted signs.
	
	Recall that the main work in the previous chapters consisted of rank-2 and rank-3 computations. Since every proper root subsystem of $ F_4 $ is of type $ A $, $ B $ or $ C $ (or a product of such root systems), most of the necessary computations have already been performed in the previous chapters. In fact, we can directly apply the coordinatisation theorems for $ B_3 $-graded groups and $ C_3 $-graded groups to certain subgroups $ G_B $ and $ G_C $ of an arbitrary $ F_4 $-graded group $ G $. All that remains to do is to investigate the \enquote{overlap} of $ G_B $ and $ G_C $, by which we mean the root groups which are contained in both groups. By comparing the commutator relations of $ G_B $ and $ G_C $ on these root groups, we obtain additional identities which turn the coordinatising structure of the short root groups into a multiplicative conic alternative algebra.
	
	The outline of this chapter is as follows. In \cref{sec:F4:conic,sec:F4:comp-alg}, we give a brief survey of conic algebras, composition algebras and related structures, following the conventions of \cite{GPR_AlbertRing}. We describe how to realise the root system $ F_4 $ as a folding of $ E_6 $ in \cref{sec:F4:rootsys}. In \cref{sec:F4:const}, we construct an $ F_4 $-graded group for each multiplicative conic algebra which is associative, commutative and faithful. We will also define the notion of coordinatisations with standard signs and the standard parity maps in this section. In \cref{sec:F4:comp}, we perform a few computations which are specific to $ F_4 $-gradings. In \cref{sec:F4:param}, we introduce the standard partial twisting system for $ F_4 $-graded groups and show that every $ F_4 $-graded group has a parametrisation by a pair of abelian groups $ (\compalg, \comring) $. In \cref{sec:F4:coord}, we use the coordinatisation theorems for $ B_3 $ and $ C_3 $ to equip $ \comring $ with the structure of an associative commutative ring and $ \compalg $ with the structure of a multiplicative conic alternative algebra over $ \comring $. We state this main result in \cref{F4:thm}.
	

\section{Conic Algebras and Related Structures}

\label{sec:F4:conic}

\begin{secnotation}
	We denote by $ \comring $ an arbitrary commutative associative ring. Unless otherwise specified, all modules and algebras are defined over~$ \comring $.
\end{secnotation}

\subsection{Algebras and Scalar Involutions}

Before we can define conic algebras, we briefly recall the notion of algebras. We will also introduce the concept of algebras with scalar involutions, which are closely related to conic algebras. For a precise formulation of this statement, see \cref{conic:conic-scalar}.

\begin{definition}[Algebra]\label{F4:compalg:alg-def}
	A \defemph*{$ \comring $-algebra}\index{algebra} is a pair consisting of a ring $ \compalg $ and a homomorphism $ \map{}{\comring}{\zentrum(\compalg)}{}{} $, called the \defemph*{structural homomorphism of $ \compalg $}\index{structural homomorphism} and usually denoted by $ \str_\compalg $. We will often simply call $ \compalg $ a $ \comring $-algebra, leaving the structural homomorphism implicit. The elements of $ \compalg $ which lie in the image of $ \str_\compalg $ are called \defemph*{scalars}\index{scalar (of an algebra)}. Further, an algebra is called \defemph*{faithful}\index{algebra!faithful} if its structural homomorphism is injective.
\end{definition}

Recall from \cref{ring:nucleus-def} that elements in the center of a ring are not only assumed to commute with all elements of the ring, but also to lie in the nucleus.

\begin{notation}[Algebras as modules]\label{F4:compalg:alg-as-module}
	Given a $ \comring $-algebra $ (\compalg, \str_\compalg) $, we can define a right $ \comring $-module structure on $ \compalg $ by
	\[ \map{\scmult}{\compalg \times \comring}{\compalg}{(a, \lambda)}{\str_\compalg(\lambda) a = a \str_\compalg(\lambda)}. \]
	Since $ \comring $ is commutative, we can also consider $ \compalg $ as a left $ \comring $-module with the same scalar multiplication, and we well use both notations depending on the circumstances. We will always use the symbol $ \scmult $ to denote the scalar multiplication on a $ \comring $-algebra.
	
	In fact, we could define a $ \comring $-algebra to be a ring $ \compalg $ together with a $ \comring $-module structure $ \scmult $ on $ \compalg $ such that the ring multiplication of $ \compalg $ is $ \comring $-bilinear. Then the map $ \map{\str_\compalg}{\comring}{\compalg}{\lambda}{1_\compalg \scmult \lambda} $ is a ring homomorphism, and its image lies in the center of $ \compalg $ because the ring multiplication is $ \comring $-bilinear.
\end{notation}

\begin{notation}
	Let $ \compalg $ be a $ \comring $-algebra, let $ \lambda \in \comring $ and let $ a \in \compalg $. We will often write $ \lambda a $ for the element $ \lambda \scmult a = \str_\compalg(\lambda)a $ to simplify notation. In particular, we will often write $ \lambda 1_\compalg $ for $ \str_\compalg(\lambda) $. However, we will never write $ \lambda $ for $ \str_\compalg(\lambda) $ because the map $ \str_\compalg $ is not necessarily injective.
\end{notation}

The notion of faithfulness for algebras agrees with the one for modules.

\begin{lemma}\label{F4:compalg:alg-faith-mod}
	Let $ \compalg $ be a $ \comring $-algebra. Then $ \compalg $ is faithful as a $ \comring $-algebra (that is, $ \str_\compalg $ is injective) if and only if it is faithful as a $ \comring $-module (that is, if $ \lambda \in \comring $ satisfies $ \lambda \compalg = \compactSet{0} $, then $ \lambda = 0 $).
\end{lemma}
\begin{proof}
	Assume that $ \compalg $ is faithful as a module and let $ \lambda \in \Ker(\str_\compalg) $. Then we have $ \lambda \scmult a = \str_\compalg(\lambda)a = 0_\compalg \cdot a = 0_\compalg $ for all $ a \in \compalg $. Thus $ \lambda =0 $ because $ \compalg $ is faithful as a module. This shows that $ \compalg $ is faithful as an algebra.
	
	Now assume that $ \compalg $ is faithful as an algebra and let $ \lambda \in \comring $ such that $ \lambda \scmult a = 0_\compalg $ for all $ a \in \compalg $. In particular, we then have $ \str_\compalg(\lambda) = \str_\compalg(\lambda) \cdot 1_\compalg = \lambda \scmult 1_\compalg = 0_\compalg $. Since $ \compalg $ is faithful as an algebra, it follows that $ \lambda = 0 $, so $ \compalg $ is faithful as a module as well.
\end{proof}

The property of a ring involution (in the sense of \cref{rinv:invo-def}) to be compatible with the algebra structure can be phrased in two equivalent ways.

\begin{lemma}\label{F4:compalg:alg-invo-char}
	Let $ \compalg $ be a $ \comring $-algebra and let $ \rinvmap $ be an involution of the ring $ \compalg $ (in the sense of \cref{rinv:invo-def}). Then the following assertions are equivalent:
	\begin{stenumerate}
		\item The map $ \rinvmap $ is linear with respect to the $ \comring $-module structure of $ \compalg $.
		
		\item The image of the structural homomorphism $ \map{\str_\compalg}{\comring}{\compalg}{}{} $ is contained in $ \Fix(\rinvmap) $.
	\end{stenumerate}
\end{lemma}
\begin{proof}
	The involution $ \rinvmap $ is $ \comring $-linear if and only if $ \rinv{(a \scmult \lambda)} = \rinv{a} \scmult \lambda $ for all $ \lambda \in \comring $ and $ a \in \compalg $. Observe that
	\begin{align*}
		\rinv{(a \scmult \lambda)} &= \rinv{\brackets[\big]{a \str_\compalg(\lambda)}} = \rinv{\str_\compalg(\lambda)} \rinv{a} \rightand \\
		\rinv{a} \scmult \lambda &= \rinv{a} \str_\compalg(\lambda) = \str_\compalg(\lambda) \rinv{a}
	\end{align*}
	for all $ \lambda \in \comring $ and $ a \in \compalg $. Thus $ \rinvmap $ is clearly linear if the image of $ \str_\compalg $ is contained in $ \Fix(\rinvmap) $. Conversely, by putting $ a \defl 1_\compalg $ in the equations above, we see that $ \rinv{\str_\compalg(\lambda)} = \str_\compalg(\lambda) $ for all $ \lambda \in \comring $ if $ \rinvmap $ is linear.
\end{proof}

With the notion of algebra involutions at hand, we can now define scalar involutions.

\begin{definition}[Scalar involution]\label{F4:scal-inv-def}
	Let $ \compalg $ be a $ \comring $-algebra. An \defemph*{(algebra) involution of $ \compalg $}\index{involution!on an algebra} is an involution $ \rinvmap $ of the ring $ \compalg $ which satisfies the equivalent properties of \cref{F4:compalg:alg-invo-char}. It is called \defemph*{scalar involution}\index{involution!scalar} if all norms (in the sense of \cref{rinv:trace-norm-def}) are scalars. That is, $ \rinv{a}a $ lies in $ \str_\compalg(\comring) $ for all $ a \in \compalg $.
\end{definition}

Since all scalars in a $ \comring $-algebra $ \compalg $ lie in the center, every scalar involution of $ \compalg $ is central (and thus nuclear) in the sense of \cref{rinv:nucl-def}. Conversely, every central involution can be regarded as a scalar involution, though we may have to change the ground ring.

\begin{lemma}[{\cite[Section~9.5]{McCrimmonAltUnpublished}}]
	Let $ \ring $ be a ring with a central involution $ \rinvmap $ and put $ \zentruminv(\ring) \defl \zentrum(\ring) \intersect \symring $. Then $ \ring $ is a faithful $ \zentruminv(\ring) $-algebra with scalar involution $ \rinvmap $.
\end{lemma}
\begin{proof}
	Since $ \zentruminv(\ring) $ is a subring of $ \ring $, the natural inclusion $ \map{i}{\zentruminv(\ring)}{\ring}{}{} $ provides a faithful $ \zentruminv(\ring) $-algebra structure on $ \ring $. For all $ \lambda \in \zentruminv(\ring) $ and all $ a \in \ring $, we have
	\[ \rinv{(a \lambda)} = \rinv{\lambda} \rinv{a} = \lambda \rinv{a} = \rinv{a} \lambda, \]
	so $ \rinvmap $ is an involution with respect to this algebra structure on $ \ring $. Further, all norms of $ \rinvmap $ lie in the center of $ \ring $ because $ \rinvmap $ is central and they satisfy $ \rinv{(\rinv{a} a)} = \rinv{a}a $ for all $ a \in \ring $. We conclude that all norms lie in $ \zentruminv(\ring) $, which shows that $ \rinvmap $ is a scalar involution
\end{proof}

\begin{remark}
	Let $ \compalg $ be a $ \comring $-algebra with an algebra involution $ \rinvmap $. Recall from \cref{rinv:nucl-cent-lem,ring:nucl-invo-char,ring:cent-invo-char} that all traces in $ \compalg $ are nuclear (or central) if $ \rinvmap $ is nuclear (or central), and that the converse is true if $ \compalg $ is alternative or if $ 2_\compalg $ is invertible. By the same argument as in \cref{rinv:nucl-cent-lem}, all traces in $ \compalg $ are scalar if $ \rinvmap $ is scalar and the converse is true if $ 2_\compalg $ is invertible. If we know that $ \compalg $ is alternative, the assumption that $ 2_\compalg $ is invertible can be relaxed slightly: It suffices that some trace in $ \compalg $ is invertible, see \cite[Section~9.5]{McCrimmonAltUnpublished}.
\end{remark}

\subsection{Conic Algebras}

\begin{convention}
	In this section, we will denote the linearisation of any $ \comring $-quadratic form $ \map{q}{\module}{\comring}{}{} $ by $ q $ as well, so that $ q(x,y) \defl q(x+y) - q(x) - q(y) $ for all $ x,y \in \module $.
\end{convention}

We can now turn to the protagonists of this section. In this subsection, we will closely follow the exposition in \cite{GPR_AlbertRing}, which in turn follows \cite{McCrimmons_ScalarInv}.

\begin{definition}[Conic algebra, {\cite[16.1]{GPR_AlbertRing}}]
	A \defemph*{conic algebra over $ \comring $}\index{conic algebra} is a pair $ (\compalg, \compnorm) $ consisting of a $ \comring $-algebra $ \compalg $ and a $ \comring $-quadratic form $ \map{\compnorm}{\compalg}{\comring}{}{} $ (in the sense of \cref{quadmod:quadform-def}) such that $ \compnorm(1_\compalg) = 1_\comring $ and such that the following \defemph{Cayley-Hamilton equation} is satisfied for all $ a \in \compalg $:
	\begin{align*}
		a^2 - \compnorm(a, 1_\compalg) a + \compnorm(a) 1_\compalg = 0_\compalg.
	\end{align*}
	The map $ \compnorm $ is called the \defemph*{norm of $ \compalg $}.\index{norm!on a conic algebra}
\end{definition}

If $ (\compalg, \compnorm) $ is a conic algebra, then $ (\compalg, \compnorm, 1_\compalg) $ is a pointed quadratic module. Thus by \cref{quadmod:pointed-tr-conj}, we have corresponding conjugation and trace maps. 

\begin{definition}[Conjugation and trace, {\cite[16.1]{GPR_AlbertRing}}]\label{conic:conj-tr-def}
	Let $ (\compalg, \compnorm) $ be a conic algebra over $ \comring $. The map
	\[ \map{\compTr}{\compalg}{\comring}{a}{\compnorm(a,1_\compalg)} \]
	is called the \defemph*{trace of $ (\compalg, \compnorm) $}\index{trace!on a conic algebra} and the map
	\[ \map{\compinvmap = -\refl{1_\compalg}}{\compalg}{\compalg}{a}{\compTr(a)1_\compalg - x} \]
	is called the \defemph*{conjugation of $ (\compalg, \compnorm) $}.\index{conjugation!on a conic algebra}
\end{definition}

The fact that conic algebras are pointable quadratic modules also has the following consequence.

\begin{lemma}\label{conic:2-faithful}
	Assume that $ 2_\comring $ is not a zero divisor and let $ (\compalg, \compnorm) $ be a conic algebra over $ \comring $. Then $ \compalg $ is a faithful $ \comring $-algebra.
\end{lemma}
\begin{proof}
	This follows from \cref{quadmod:2-faithful}.
\end{proof}

\begin{note}
	In \cite[Sections~2.1,~2.2]{McCrimmonAltUnpublished} and \cite[\onpage{86}]{McCrimmons_ScalarInv}, an algebra $ \compalg $ is called \defemph*{of degree 2}\index{degree 2 algebra} if there exist a quadratic form $ \map{\compnorm}{\compalg}{\comring}{}{} $ and a linear form $ \map{\compTr}{\compalg}{\comring}{}{} $ such that $ \compnorm(1_\compalg) = 1_\comring $, $ \compTr(1_\compalg) = 2_\comring $ and
	\[ a^2 - \compTr(a) a + \compnorm(a) 1_\compalg = 0 \]
	for all $ a \in \compalg $. Thus conic algebras are precisely those algebras of degree 2 in which $ \compTr(a) = \compnorm(a, 1_\compalg) $ for all $ a \in \compalg $. If $ \compalg $ is an algebra of degree 2 with norm $ \compnorm $ and trace $ \compTr $, then $ \compnorm(a, 1_\compalg) 1_\compalg = \compTr(a) 1_\compalg $ for all $ a \in \compalg $ by \cite[(0.8)]{McCrimmons_ScalarInv}. Hence for faithful algebras, the two notions coincide. In particular, they coincide if $ 2_\comring $ is not a zero divisor because \cref{conic:2-faithful} is true for algebras of degree 2 as well. Even in the non-faithful setting, many results in this section remain true for algebras of degree 2, but we will have no need to use this fact.
\end{note}

Here are two examples of conic algebras. The first one motivates the name of the Cayley-Hamilton equation and the second one will be useful in the construction of $ F_4 $-graded groups in \cref{sec:F4:const}.

\begin{example}
	Let $ \compalg $ be the $ \comring $-algebra of $ (2 \times 2) $-matrices over $ \comring $ and put $ \compnorm(a) \defl \det(a) $ for all $ a \in \compalg $. Then $ (\compalg, \compnorm) $ is a conic algebra over $ \comring $. Its trace is the usual trace map
	\[ \map{\compTr}{\compalg}{\comring}{a}{a_{11} + a_{22}} \]
	and its conjugation is
	\[ \map{\compinvmap}{\compalg}{\compalg}{\begin{pmatrix}
		a_{11} & a_{12} \\
		a_{21} & a_{22}
	\end{pmatrix}}{\begin{pmatrix}
		a_{22} & -a_{12} \\
		-a_{21} & a_{11}
	\end{pmatrix}}. \]
	The \defemph{Cayley-Hamilton theorem} says that any $ a \in \compalg $ is a root of its characteristic polynomial, and this is precisely the Cayley-Hamilton equation in this situation.
\end{example}

\begin{example}[{\cite[16.2~(d)]{GPR_AlbertRing}}]\label{conic:dirprod-example}
	Put $ \compalg \defl \comring \times \comring $ and $ \compnorm(a,b) \defl ab $ for all $ (a,b) \in \compalg $. Then $ (\compalg, \compnorm) $ is a conic algebra over $ \comring $ with trace $ \compTr(a,b) = a+b $ and conjugation $ \compinv{(a,b)} = (b,a) $ for all $ (a,b) \in \compalg $. Further, the linearisation of $ \compnorm $ is
	\[ \map{}{\compalg \times \compalg}{\comring}{\brackets[\big]{(a,b), (c,d)}}{ad+cb}. \]
	Observe that, forgetting the multiplication on $ \compalg $, the quadratic module $ (\compalg, \compnorm) $ is the hyperbolic plane (that is, the hyperbolic space of dimension~2) from \cref{quadmod:hyperbolic-ex}.
\end{example}

Conic algebras arise from $ F_4 $-graded groups in the form of the following parameter system.

\begin{definition}[Standard parameter system]\label{conic:standard-param}
	Define $ \twistgroup \defl \invogroup \defl \compactSet{\pm 1} $. We declare that $ \twistgroup $ acts on $ \compalg $ and $ \comring $ by inversion and that $ \invogroup $ acts trivially on $ \comring $ and by conjugation on $ \compalg $ (that is, by the switching involution). More precisely, this means that
	\begin{align*}
		-1_\twistgroup.\lambda = -\lambda, \quad -1_\twistgroup.a = -a, \quad -1_\invogroup.\lambda = \lambda, \quad -1_\invogroup.a = \compinv{a}
	\end{align*}
	for all $ \lambda \in \comring $ and all $ a \in \compalg $. Then the triple $ (\twistgroup \times \invogroup, \compalg, \comring) $ is called the \defemph*{standard parameter system for $ (\compalg, \compnorm) $}.\index{parameter system!standard!type F4@type $ F_4 $}
\end{definition}

We now investigate some basic properties of conic algebras.
The following identities are straightforward to verify.

\begin{lemma}[{\cite[16.5]{GPR_AlbertRing}}]\label{conic:basic-id}
	Let $ (\compalg, \compnorm) $ be a conic algebra over $ \comring $. Then the following hold:
	\begin{lemenumerate}
		\item $ \compnorm(1_\compalg) = 1_\comring $ and $ \compTr(1_\compalg) = 2_\comring $.
		
		\item $ \compinv{1_\compalg} = 1_\compalg $ and $ \compinv{\compinv{a}} = a $ for all $ a \in \compalg $.
		
		\item \label{conic:basic-id:scalar}$ \compinv{a}a = \compnorm(a)1_\compalg = a \compinv{a} $ and $ a + \compinv{a} = \compTr(a)1_\compalg $ for all $ a \in \compalg $.
		
		\item \label{conic:basic-id:scalar-lin}$ a \compinv{b} + b \compinv{a} = \compnorm(a,b)1_\compalg $ for all $ a,b \in \compalg $.
		
		\item \label{conic:basic-id:invar}$ \compnorm(\compinv{a}) = \compnorm(a) $, $ \compnorm(\compinv{a}, \compinv{b}) = \compnorm(a,b) $ and $ \compTr(\compinv{a}) = \compTr(a) $ for all $ a,b \in \compalg $.
		
		\item The conjugation and trace maps are $ \comring $-linear.
	\end{lemenumerate}
\end{lemma}

In order to show that the conjugation on a conic algebra is an (algebra) involution, it only remains to show that it is anti-compatible with the multiplication. This is, however, not true in general. Instead, we have the following characterisation.

\begin{lemma}[{\cite[16.10]{GPR_AlbertRing}}]\label{conic:invo-char}
	Let $ (\compalg, \compnorm) $ be a conic algebra over $ \comring $. Then the conjugation of $ \compalg $ is an algebra involution if and only if $ \compTr(ab) 1_\compalg = \compnorm(a, \compinv{b}) 1_\compalg $ (that is, $ \compnorm(ab,1_\compalg) 1_\compalg = \compnorm(a, \compinv{b}) $) for all $ a,b \in \compalg $.
\end{lemma}

\begin{remark}[Scalar involutions and conic algebras, {\cite[1.1]{McCrimmons_ScalarInv}}]\label{conic:conic-scalar}
	Let $ \compalg $ be a $ \comring $-algebra. For the scope of this remark, a \defemph*{linear involution of $ \comring $} is defined to be the same thing as an algebra involution except that it is not assumed to be anti-compatible with the ring multiplication. The notion of scalar linear involutions is defined in the same way as for algebra involutions. It remains true in this generality that all traces with respect to a scalar involution are scalars and that, by the same proof as in \cref{rinv:cent-norm-comm}, we have $ a \rinv{a} = \rinv{a}a $ for all $ a \in \compalg $ if $ \rinvmap $ is a scalar linear involution.

	Now assume that $ (\compalg, \compnorm) $ is a conic algebra over $ \comring $. Then \thmitemcref{conic:basic-id}{conic:basic-id:scalar} says precisely that the linear involution $ \compinvmap $ is scalar, though of course it need not be a ring involution. In particular, it is also central and nuclear.
	
	Conversely, let $ \compalg $ be a $ \comring $-algebra with a scalar linear involution $ \rinvmap $. Denote the image of $ \comring $ in $ \compalg $ under the structural homomorphism by $ \tilde{\comring} $. (If $ \compalg $ is a faithful $ \comring $-algebra, then we can identify $ \comring $ with $ \tilde{\comring} $.) We define
	\[ \map{\compTr}{\compalg}{\tilde{\comring}}{a}{\ringTr_\rinvmap(a) = a+\rinv{a}} \midand \map{\compnorm}{\compalg}{\tilde{\comring}}{a}{\ringnorm_\rinvmap(a) = a\rinv{a} = \rinv{a}a}. \]
	Then for all $ a \in \compalg $, we have
	\begin{align*}
		a^2 - \compTr(a)a + \compnorm(a) 1_\compalg &= a^2 -(a + \rinv{a})a + \rinv{a} a = 0_\compalg \rightand \\
		\compnorm(a, 1_\compalg) &= (a+1_\compalg)\rinv{(a+1_\compalg)} - a \rinv{a} - 1_\compalg = a + \rinv{a}. 
	\end{align*}
	Hence $ (\compalg, \compnorm) $ is a faithful conic $ \tilde{\comring} $-algebra with trace $ \compTr $ and conjugation $ \rinvmap $. In particular, every faithful $ \comring $-algebra with a scalar linear involution is also a faithful conic algebra over $ \comring $. Without faithfulness assumptions, however, it is not at all clear that $ \compnorm $ and $ \compTr $ can be lifted from $ \tilde{\comring} $ to $ \comring $.
	
	We conclude from the previous two paragraphs that faithful $ \comring $-algebras with a scalar linear involution are the same thing as faithful conic algebras. 
\end{remark}

The following identities, which are well-known in the theory of composition algebras over a field, turn out to be equivalent.

\begin{lemma}[{\cite[16.12]{GPR_AlbertRing}}]\label{conic:norm-ass-char}
	Let $ (\compalg, \compnorm) $ be a conic algebra over $ \comring $. Then the following assertions are equivalent:
	\begin{lemenumerate}
		\item \label{conic:norm-ass-char:1}$ \compnorm(a,ba) = \compnorm(ba, a) = \compTr(b) \compnorm(a) $ for all $ a,b \in \compalg $.
		
		\item $ \compnorm(a,ab) = \compnorm(ab, a) = \compTr(b) \compnorm(a) $ for all $ a,b \in \compalg $.
		
		\item $ \compnorm(ab,c) = \compnorm(a,c \compinv{b}) $ for all $ a,b,c \in \compalg $.
		
		\item \label{conic:norm-ass-char:4}$ \compnorm(ab,c) = \compnorm(b, \compinv{a}c) $ for all $ a,b,c \in \compalg $.
		
		\item \label{conic:norm-ass-char:5}$ \compTr(ab) = \compnorm(a, \compinv{b}) $ and $ \compTr\brackets[\big]{(ab)c} = \compTr\brackets[\big]{a(bc)} $ for all $ a,b,c \in \compalg $.
	\end{lemenumerate}
\end{lemma}

\begin{definition}[Norm-associativity, {\cite[16.12]{GPR_AlbertRing}}]
	A conic algebra is called \defemph*{norm-associative} if it satisfies the equivalent conditions in \cref{conic:norm-ass-char}.
\end{definition}

By \cref{conic:invo-char} and \thmitemcref{conic:norm-ass-char}{conic:norm-ass-char:5}, the conjugation in norm-associative conic algebras is an algebra involution. It turns out that these algebras are automatically flexible (in the sense of \cref{ring:alternative-def}) as well.

\begin{lemma}[{\cite[16.14]{GPR_AlbertRing}}]\label{conic:normass-flex-conj}
	Let $ (\compalg, \compnorm) $ be a norm-associative conic algebra over $ \comring $. Then $ \compalg $ is flexible and the conjugation of $ \compalg $ is an algebra involution.
\end{lemma}

\begin{note}
	The converse of \cref{conic:normass-flex-conj} is nearly true: If $ (\compalg, \compnorm) $ is a flexible conic algebra whose conjugation is an involution, and if in addition $ \compalg $ is projective as a $ \comring $-module, then $ (\compalg, \compnorm) $ is norm-associative (\cite[16.14]{GPR_AlbertRing}). Further, the projectiveness assumption can be dropped if the algebra structure on $ \compalg $ is faithful (\cite[16.15]{GPR_AlbertRing}).
\end{note}

In general, the norm of a conic algebra $ (\compalg, \compnorm) $ is not uniquely determined by the algebra $ \compalg $, as is illustrated by the family of examples in \cite[Exercise~17.8]{GPR_AlbertRing}. However, we have the following result.

\begin{lemma}[{\cite[16.16]{GPR_AlbertRing}}]\label{conic:norm-unique}
	Let $ \compalg $ be a $ \comring $-algebra and let $ \map{\compnorm, \compnorm'}{\compalg}{\comring}{}{} $ be two quadratic forms which turn $ \compalg $ into a conic algebra. If $ \compalg $ is projective as a $ \comring $-module, then $ \compnorm = \compnorm' $.
\end{lemma}

The following property is eponymous in the theory of composition algebras.

\begin{definition}[Multiplicative conic algebra]
	A conic algebra $ (\compalg, \compnorm) $ over $ \comring $ is called \defemph*{multiplicative}\index{conic algebra!multiplicative} if its norm permits composition in the sense that $ \compnorm(ab) = \compnorm(a) \compnorm(b) $ for all $ a,b \in \compalg $.
\end{definition}

\begin{lemma}[{\cite[17.2]{GPR_AlbertRing}}]\label{conic:mult-is-norm}
	Multiplicative conic algebras are norm-associative. In particular, they are flexible and their conjugation is a scalar algebra involution.
\end{lemma}
\begin{proof}
	This holds by \cite[17.2]{GPR_AlbertRing}, except for the assertion that the conjugation is scalar, which holds by \cref{conic:conic-scalar}.
\end{proof}

In the following, we list some properties of conic algebras which are, in addition, alternative.

\begin{lemma}[Kirmse's identities, {\cite[17.4]{GPR_AlbertRing}}]\index{Kirmse's identities}
	Let $ (\compalg, \compnorm) $ be a conic alternative algebra. Then we have
	\[ a(\compinv{a}b) = \compnorm(a)b = (b \compinv{a})a \]
	for all $ a,b \in \compalg $.
\end{lemma}

\begin{note}[{\cite[17.4]{GPR_AlbertRing}}]\label{conic:alt-proj-is-mult}
	Let $ (\compalg, \compnorm) $ be a conic alternative algebra. One can show that
	\[ \compnorm(aba) = \compnorm(a)^2 \compnorm(b) \]
	for all $ a,b \in \compalg $. In this sense, conic alternative algebras are \enquote{close to being multiplicative}. Yet, \cite[Exercise~17.8]{GPR_AlbertRing} provides examples of conic alternative algebras which are not multiplicative. If $ \compalg $ is projective as a $ \comring $-module, however, $ (\compalg, \compnorm) $ is indeed multiplicative by \cite[17.6]{GPR_AlbertRing}.
\end{note}

\begin{lemma}[{\cite[17.5]{GPR_AlbertRing}}]\label{conic:norm-inv}
	Let $ (\compalg, \compnorm) $ be a conic alternative algebra over $ \comring $ and let $ a \in \compalg $. Then $ a $ is invertible (in $ \compalg $) if and only if $ \compnorm(a) $ is invertible (in $ \comring $). In this case, $ a^{-1} = \compnorm(a)^{-1} \compinv{a} $ and $ \compnorm(a^{-1}) = \compnorm(a)^{-1} $.
\end{lemma}

\subsection{\texorpdfstring{$ F_4 $}{F\_4}-data}

We now define an algebraic structure that we call \enquote{$ F_4 $-datum}. In \cref{sec:F4:coord}, we will see that $ F_4 $-data are precisely the kind of structures that coordinatise $ F_4 $-graded groups. The goal of this section is to show, on a purely algebraic level, that an $ F_4 $-datum contains the same information as a multiplicative conic alternative algebra.

\begin{definition}[$ F_4 $-datum]\label{F4:compalg:F4-datum-def}
	An \defemph*{$ F_4 $-datum}\index{F4-datum@$ F_4 $-datum} is a tuple
	\[ \brackets[\big]{\comring, \compalg, \mathord{\scmult}, \compnorm, \scp{}{}, \compinvmap, \str_\compalg, \jorsc} \]
	with the following properties:
	\begin{stenumerate}
		\item \label{F4:compalg:F4-datum-def:comring}$ \comring $ is a commutative assocative ring.
		
		\item \label{F4:compalg:F4-datum-def:compalg}$ \compalg $ is an alternative ring and $ \scmult $ is a right $ \comring $-module structure on $ \compalg $ (which, a priori, is not assumed to be compatible with the ring structure).
		
		\item \label{F4:compalg:F4-datum-def:quad}$ \map{\compnorm}{\compalg}{\comring}{}{} $ is a quadratic form with $ \compnorm(1_\compalg) = 1_\comring $, and $ \map{\scp{}{}}{\compalg \times \compalg}{\comring}{}{} $ is the linearisation of $ \compnorm $.
		
		\item \label{F4:compalg:F4-datum-def:inv}$ \map{\compinvmap}{\compalg}{\compalg}{}{} $ is a nuclear involution of the ring $ \compalg $, and we have $ \compinvmap = -\refl{1_\compalg} $ where $ \refl{1_\compalg} $ denotes the reflection corresponding to $ 1_\compalg $ in the quadratic module $ (\compalg, \compnorm) $ (in the sense of \cref{quadmod:refl-def}).
		
		\item \label{F4:compalg:F4-datum-def:str}$ \map{\str_\compalg}{\comring}{\nucleus(\compalg)}{}{} $ is a map which satisfies $ a \scmult \lambda = a \ast \str_\compalg(\lambda) $ for all $ a \in \compalg $ and $ \lambda \in \comring $ (where $ \ast $ denotes the ring multiplication on $ \compalg $).
		
		\item \label{F4:compalg:F4-datum-def:jorsc}$ \jorsc $ is the map $ \map{}{\comring \times \compalg}{\comring}{(\lambda, a)}{\lambda \compnorm(a)} $.
		
		\item \label{F4:compalg:F4-datum-def:jormod}$ (\comring, \jorsc, \str_\compalg, \scp{}{1_\compalg}, 0) $ is a Jordan module of type $ C $ over $ (\compalg, \compinvmap) $.
	\end{stenumerate}
\end{definition}

\begin{note}
	The definition of $ F_4 $-data contains a lot of redundancy. For example, the maps $ \scp{}{} $, $ \compinvmap $ and $ \jorsc $ are uniquely determined by the remaining part of an $ F_4 $-datum. For our application of $ F_4 $-data in \cref{sec:F4:coord}, however, it will be practical to keep this redundancy as part of the definition.
\end{note}

\begin{proposition}\label{F4:compalg:comp-is-F4-datum}
	Let $ (\compalg, \compnorm) $ be a multiplicative conic alternative algebra over $ \comring $ and define
	\[ \map{\jorsc}{\comring \times \compalg}{\comring}{(\lambda, a)}{\lambda \compnorm(a)}. \]
	Then $ (\comring, \compalg, \mathord{\scmult}, \compnorm, \scp{}{}, \compinvmap, \str_\compalg, \jorsc) $ is an $ F_4 $-datum where $ \scmult $ denotes the $ \comring $-module structure on $ \compalg $ (from \cref{F4:compalg:alg-as-module}), $ \scp{}{} $ denotes the linearisation of $ \compnorm $, $ \compinvmap $ is the conjugation on $ \compalg $, $ \str_\compalg $ is the structural homomorphism of the $ \comring $-algebra $ \compalg $ and the $ \comring $-module structure on $ \compalg $ is taken to be the one induced by the $ \comring $-algebra structure.
\end{proposition}
\begin{proof}
	Axioms~\thmitemref{F4:compalg:F4-datum-def}{F4:compalg:F4-datum-def:comring}, \thmitemref{F4:compalg:F4-datum-def}{F4:compalg:F4-datum-def:compalg}, \thmitemref{F4:compalg:F4-datum-def}{F4:compalg:F4-datum-def:quad} are trivial. The conjugation map $ \compinvmap $ is a nuclear involution by \cref{conic:mult-is-norm} and it equals $ -\refl{1_\compalg} $ by its definition. Hence Axiom~\thmitemref{F4:compalg:F4-datum-def}{F4:compalg:F4-datum-def:inv} is satisfied. Further, Axioms~\thmitemref{F4:compalg:F4-datum-def}{F4:compalg:F4-datum-def:str} and~\thmitemref{F4:compalg:F4-datum-def}{F4:compalg:F4-datum-def:jorsc} hold by the definitions of $ \scmult $ and of $ \jorsc $, respectively. Denote by $ \compTr \defl \scp{}{1_\compalg} $ the trace of the conic algebra $ (\compalg, \compnorm) $. It remains to show that $ (\comring, \jorsc, \str_\compalg, \compTr, 0) $ satisfies the axioms of a Jordan module of type $ C $ over $ (\compalg, \compinvmap) $ from \cref{BC:jordanmodule-def}. That is, $ \jorprojone \defl \str_{\compalg} $ plays the role of the Jordan module projection and $ \jorTrone \defl \compTr $ plays the role of the Jordan module trace.
	
	At first, we verify that the maps $ \jorsc $, $ \jorprojone $ and $ \jorTrone $ have the desired properties in the first paragraph of \cref{BC:jordanmodule-def}. We will always denote the square-scalar multiplication on $ \ring $ from \cref{jormod:invring-jordan} by $ \omega $. It is clear that $ \jorsc $ is additive in the first component and weakly quadratic in the second component, so $ (\comring, \jorsc) $ is a square-module. Since $ \comring $ is associative and $ \compnorm $ admits composition, it is a multiplicative square-module.
	
	Since $ \compalg $ is a $ \comring $-algebra, the image of $ \str_\compalg = \jorprojone $ lies in the center and thus in the nucleus of $ \compalg $. Further, $ \str_\compalg $ preserves the Jordan scalar multiplication: For all $ a \in \compalg $ and $ \lambda \in \comring $, we have
	\begin{align*}
		\str_\compalg\brackets[\big]{\jorsc(\lambda, a)} &= \str_\compalg\brackets[\big]{\lambda \compnorm(a)} = \str_\compalg(\lambda) \str_\compalg\brackets[\big]{\compnorm(a)} = \str_\compalg(\lambda) \compinv{a} a \\
		&= \compinv{a} \str_\compalg(\lambda) a = \compinv{a} \jorprojone(\lambda) a = \omega\brackets[\big]{\jorprojone(\lambda), a}
	\end{align*}
	by \thmitemcref{conic:basic-id}{conic:basic-id:scalar} and because $ \comring $ is commutative associative.
	
	It is clear that $ \jorTrone $ is additive. Now let $ a,b \in \compalg $. Then
	\begin{align*}
		\jorTrone\brackets[\big]{\omega(a, b)} &= \compTr\brackets[\big]{\compinv{b} a b} = \compnorm(\compinv{b}ab, 1_\compalg),
	\end{align*}
	which by \thmitemcref{conic:norm-ass-char}{conic:norm-ass-char:4} equals $ \compnorm(ab, b) $. By \thmitemcref{conic:norm-ass-char}{conic:norm-ass-char:1}, this term equals $ \compTr(a) \compnorm(b) $. It follows that
	\begin{align*}
		\jorTrone\brackets[\big]{\omega(a, b)} &= \compTr(a) \compnorm(b) = \jorsc\brackets[\big]{\jorTrone(a), b},
	\end{align*}
	so $ \jorTrone $ is a homomorphism of square-modules. We conclude that the maps $ \jorsc $, $ \jorprojone $ and $ \jorTrone $ have all the desired properties in the first paragraph of \cref{BC:jordanmodule-def}.
	
	We now turn to the remaining axioms in \cref{BC:jordanmodule-def}. Since $ \compinvmap $ is an algebra involution, the image of $ \str_\compalg $ lies in $ \Fix(\compinvmap) $ (by \cref{F4:compalg:alg-invo-char}), which implies that Axiom~\thmitemref{BC:jordanmodule-def}{BC:jordanmodule-def:proj-rinv} holds. Axiom~\thmitemref{BC:jordanmodule-def}{BC:jordanmodule-def:tr} holds by \thmitemcref{conic:basic-id}{conic:basic-id:invar} and~\thmitemcref{conic:norm-ass-char}{conic:norm-ass-char:5}, the latter being applicable by \cref{conic:mult-is-norm}. Axiom~\thmitemref{BC:jordanmodule-def}{BC:jordanmodule-def:psi} is trivial because the Jordan module skew-hermitian form is zero.
	
	Now consider Axiom~\thmitemref{BC:jordanmodule-def}{BC:jordanmodule-def:lin}. Equation~\eqref{eq:jormod-lin:proj} holds because $ \jorprojone = \str_{\compalg} $ is additive and equation~\eqref{eq:jormod-lin:add} holds because $ (\comring, +) $ is abelian. Further, for all $ \lambda \in \comring $ and $ a,b \in \compalg $, we have
	\begin{align*}
		\jorsc(\lambda, a+b) &= \lambda \compnorm(a+b) = \lambda \compnorm(a) + \lambda \compnorm(b) + \lambda \compnorm(a,b)
	\end{align*}
	where
	\begin{align*}
		\lambda \compnorm(a,b) &= \lambda \compnorm(b,a) = \lambda \compnorm(\compinv{b}, \compinv{a}) = \compnorm\brackets[\big]{\compinv{b} \str_\compalg(\lambda)a, 1_\compalg} = \jorTrone\brackets[\big]{\compinv{b} \jorprojone(\lambda)a}
	\end{align*}
	by \thmitemcref{conic:basic-id}{conic:basic-id:invar}. This shows that equation~\eqref{eq:jormod-lin:jorsc} holds as well, and so Axiom~\thmitemref{BC:jordanmodule-def}{BC:jordanmodule-def:lin} is satisfied.
	
	By \thmitemcref{conic:basic-id}{conic:basic-id:scalar}, we have
	\begin{align*}
		\jorprojone\brackets[\big]{\jorTrone(a)} &= \str_\compalg\brackets[\big]{\compTr(a)} = \compTr(a) 1_\compalg = a + \rinv{a}
	\end{align*}
	for all $ a \in \compalg $, so Axiom~\thmitemref{BC:jordanmodule-def}{BC:jordanmodule-def:square} holds. Finally, Axiom~\thmitemref{BC:jordanmodule-def}{BC:jordanmodule-def:inv} is satisfied for $ v_0 \defl 1_\comring $ because $ \jorprojone(1_\comring) = \str_\compalg(1_\comring) = 1_\compalg $.
\end{proof}

\begin{proposition}\label{F4:compalg:F4-datum-is-comp}
	Let $ (\comring, \compalg, \mathord{\scmult}, \compnorm, \scp{}{}, \compinvmap, \str_\compalg, \jorsc) $ be an $ F_4 $-datum. Then the following hold:
	\begin{lemenumerate}
		\item \label{F4:compalg:F4-datum-is-comp:alg}The ring $ \compalg $ together with the $ \comring $-module structure $ \scmult $ is a $ \comring $-algebra (in the sense of \cref{F4:compalg:alg-as-module}), and its structural homomorphism is exactly $ \str_\compalg $.
		
		\item \label{F4:compalg:F4-datum-is-comp:comp}$ (\compalg, \compnorm) $ is a multiplicative conic alternative algebra over $ \comring $ with conjugation $ \compinvmap $.
	\end{lemenumerate}
\end{proposition}
\begin{proof}
	Putting $ a \defl 1_\compalg $ in Axiom~\thmitemref{F4:compalg:F4-datum-def}{F4:compalg:F4-datum-def:str}, we see that
	\begin{equation}\label{eq:F4-datum}
		\str_\compalg(\lambda) = 1_\compalg * \str_{\compalg}(\lambda) = 1_\compalg \scmult \lambda
	\end{equation}
	for all $ \lambda \in \comring $. In particular, $ \str_\compalg(1_\comring) = 1_\compalg $ and $ \str_{\compalg} $ is additive.
	Further, it follows from Axiom~\thmitemref{F4:compalg:F4-datum-def}{F4:compalg:F4-datum-def:str} that
	\begin{align*}
		\str_\compalg(\lambda \mu) &= 1_\compalg * \str_\compalg(\lambda \mu) = 1_\compalg \scmult \lambda \mu = (1_\compalg \scmult \lambda) \scmult \mu = \str_\compalg(\lambda) \scmult \mu
	\end{align*}
	for all $ \lambda, \mu \in \comring $. By another application of Axiom~\thmitemref{F4:compalg:F4-datum-def}{F4:compalg:F4-datum-def:str}, we infer that
	\begin{align*}
		\str_\compalg(\lambda \mu) &= \str_\compalg(\lambda) * \str_\compalg(\mu)
	\end{align*}
	for all $ \lambda, \mu \in \comring $, so $ \str_\compalg $ is a homomorphism of rings. Further, since $ \str_\compalg $ is the Jordan module projection in an abelian Jordan module (by Axiom~\thmitemref{F4:compalg:F4-datum-def}{F4:compalg:F4-datum-def:jormod}), we know from Axiom~\thmitemref{BC:jordanmodule-def}{BC:jordanmodule-def:proj-rinv} that $ \compinv{\str_\compalg(\lambda)} = \str_\compalg(\lambda) $ for all $ \lambda \in \comring $. Moreover, $ \compinvmap $ is $ \comring $-linear because it equals $ -\refl{1_\compalg} $ by Axiom~\thmitemref{F4:compalg:F4-datum-def}{F4:compalg:F4-datum-def:inv}. In other words, $ \compinv{a \scmult \lambda} = \compinv{a} \scmult \lambda $ for all $ a \in \compalg $ and $ \lambda \in \comring $. We infer that
	\begin{align*}
		\str_\compalg(\lambda) * \compinv{a} &= \compinv{\str_\compalg(\lambda)} * \compinv{a} = \compinv{a * \str_\compalg(\lambda)} = \compinv{a \scmult \lambda} = \compinv{a} \scmult \lambda = \compinv{a} * \str_\compalg(\lambda)
	\end{align*}
	for all $ a \in \compalg $ and $ \lambda \in \comring $. This says precisely that the image of $ \str_\compalg $ is contained in the center of $ \compalg $. Thus $ (\compalg, \str_\compalg) $ is a $ \comring $-algebra, and by \eqref{eq:F4-datum}, its $ \comring $-scalar multiplication is exactly $ \scmult $. This proves~\itemref{F4:compalg:F4-datum-is-comp:alg}.
	
	We now prove that $ (\compalg, \compnorm) $ is a conic algebra. We already know from Axiom~\thmitemref{F4:compalg:F4-datum-def}{F4:compalg:F4-datum-def:quad}  that $ \compnorm $ is a quadratic form with $ \compnorm(1_\compalg) = 1_\comring $, so it only remains to verify the Cayley-Hamilton equation. As usual, we denote the square-scalar multiplication on $ \compalg $ by $ \omega $. Using that the Jordan module projection $ \str_\compalg $ preserves the square-scalar multiplications by Axiom~\thmitemref{F4:compalg:F4-datum-def}{F4:compalg:F4-datum-def:jormod}, we find that
	\begin{align*}
		\str_\compalg\brackets[\big]{\compnorm(a)} &= \str_\compalg\brackets[\big]{\jorsc(1_\comring, a)} = \omega\brackets[\big]{\str_\compalg(1_\comring), a} = \compinv{a} \str_\compalg(1_\comring) a = \compinv{a}a
	\end{align*}
	for all $ a \in \compalg $. Together with Axioms~\thmitemref{BC:jordanmodule-def}{BC:jordanmodule-def:square} and~\thmitemref{F4:compalg:F4-datum-def}{F4:compalg:F4-datum-def:str} and the fact that the image of $ \str_\compalg $ is central, this implies that
	\begin{align*}
		a \scmult \scp{a}{1_\compalg} &= a \str_\compalg\brackets{\scp{a}{1_\compalg}} = \str_\compalg\brackets{\scp{a}{1_\compalg}} a = (a+\compinv{a})a \\
		&= a^2 + \compinv{a}a = a^2 + \str_\compalg\brackets[\big]{\compnorm(a)} = a^2 + 1_\compalg \scmult \compnorm(a)
	\end{align*}
	for all $ a \in \compalg $. This says precisely that the Cayley-Hamilton equation is satisfied, so $ (\compalg, \compnorm) $ is a conic algebra. By Axiom~\thmitemref{F4:compalg:F4-datum-def}{F4:compalg:F4-datum-def:compalg}, it is also alternative.
	
	It remains to show that $ (\compalg, \compnorm) $ is multiplicative. By the definition of $ \jorsc $ (see Axiom~\thmitemref{F4:compalg:F4-datum-def}{F4:compalg:F4-datum-def:jorsc}) and the fact that it is the scalar multiplication of the multiplicative square-module $ \comring $ by Axiom~\thmitemref{F4:compalg:F4-datum-def}{F4:compalg:F4-datum-def:jormod}, we have
	\begin{align*}
		\compnorm(ab) &= \jorsc(1_\comring, ab) = \jorsc\brackets[\big]{\jorsc(1_\comring, a), b} = \jorsc\brackets[\big]{\compnorm(a), b} = \compnorm(a) \compnorm(b)
	\end{align*}
	for all $ a,b \in \compalg $. Thus $ \compnorm $ permits composition, which finishes the proof.
\end{proof}


\section{Composition Algebras and Pre-composition Algebras}

\label{sec:F4:comp-alg}

In the previous section, we have seen that multiplicative conic alternative algebras -- that is, the objects that will turn out to coordinatise $ F_4 $-gradings -- have many desirable properties. In this section, we address the question which additional properties a conic algebra should satisfy in order to deserve the name \enquote{composition algebra}. There exists no standard terminology for these objects over commutative associative rings. In this section, we present the conventions of \cite[Chapter~IV]{GPR_AlbertRing}.

The content of this section plays no role for the proof of our main results, but it is valuable to put the notion of conic algebras into a more well-known context. Further, we will see in \cref{F4:stsign:weyl-char} that RGD-systems of type $ F_4 $ correspond precisely to pre-composition algebras over fields, which are the same as regular composition algebras if the field is not of characteristic 2.

\subsection{Regularity Conditions on Quadratic Forms}

Before can state the definition of composition algebras, we have to introduce the notions of regular and non-singular quadratic forms. These properties replace the usual requirement of non-degeneracy in the classical theory, which is inadequate over rings because it is not stable under base change. Again, there is no standard terminology. In the following, we present the conventions of \cite[11.9, 11.11]{GPR_AlbertRing}.

\begin{definition}[Radical]
	Let $ \module $ be a $ \comring $-module and let $ \map{f}{\module \times \module}{\comring}{}{} $ be a symmetric $ \comring $-bilinear form. The \defemph*{radical of $ f $}\index{radical} is
	\[ \Rad(f) \defl \Set{v \in \module \given f(u,v) = 0 \text{ for all } u \in \module}. \]
	If $ \map{q}{\module}{\comring}{}{} $ is a $ \comring $-quadratic form with linearisation $ f $, then the \defemph*{radical of $ q $} is
	\[ \Rad(q) \defl \Set{v \in \Rad(f) \given q(v) = 0}. \]
\end{definition}

\begin{definition}[Dual module]
	Let $ \module $ be a $ \comring $-module. The \defemph*{dual of $ \module $}\index{dual module} is $ \module^* \defl \Hom_\comring(\module, \comring) $. If $ \map{f}{\module \times \module}{\comring}{}{} $ be a symmetric $ \comring $-bilinear form, we put
	\[ \map{\tilde{f}}{\module}{\module^*}{v}{\brackets[\big]{\map{}{}{}{u}{f(v,u)}}}. \]
\end{definition}

\begin{definition}[Regularity conditions on quadratic forms]
	Let $ \module $ be a $ \comring $-module, let $ \map{q}{\module}{\comring}{}{} $ be a $ \comring $-quadratic form and denote by $ \map{f}{\module \times \module}{\comring}{}{} $ its linearisation.
	\begin{defenumerate}
		\item $ q $ is called \defemph*{non-degenerate}\index{quadratic form!non-degenerate} if $ \Rad(q) $ is zero.
		
		\item $ q $ is called \defemph*{weakly regular}\index{quadratic form!weakly regular} if $ \Rad(f) $ is zero.
		
		\item $ q $ is called \defemph*{non-singular}\index{quadratic form!non-singular} if $ \module $ is projective and for all $ \comring $-algebras $ K $ which are fields, the scalar extension $ q_K $ is non-degenerate.
		
		\item $ q $ is called \defemph*{regular}\index{quadratic form!regular} if $ \module $ is finitely generated projective and the map $ \tilde{f} $ is an isomorphism.
	\end{defenumerate}
\end{definition}

\begin{remark}[Base change]
	Let $ \module $ be a $ \comring $-module and let $ \map{q}{\module}{\comring}{}{} $ be a $ \comring $-quadratic form. For any $ \comring $-algebra $ R $, the base change $ (\module^*)_R $ of the dual module can be identified with the dual $ (\module_R)^* $ of the base change in a canonical way if $ \module $ is finitely generated projective (see \cite[9.15]{GPR_AlbertRing}). From this it follows that the notion of regularity is invariant under base change. The same holds for non-singularity. However, non-degeneracy is not invariant under base change, and this defect is precisely the motivation to introduce non-singularity.
\end{remark}

Some important relations between the regularity conditions can be found in \cref{fig:reg-props}. They are either straightforward to prove or can be found in \cite[11.9, 11.11]{GPR_AlbertRing}.

\begin{figure}[htb]
	\centering$ \begin{tikzcd}[arrows=Rightarrow, column sep=6.5cm, row sep=5em]
		\text{regular} \arrow[d] \arrow[r, bend right=2] & \text{weakly regular} \arrow[l, bend right=10, "\comring \text{ field, } \dim \module < \infty"'] \arrow[d, bend right] \\
		\text{non-singular} \arrow[r, "\comring \text{ field}"'] & \text{non-degenerate} \arrow[ul, pos=0.65, bend left=2, "\addtolength{\jot}{-2mm}\begin{gathered}
			\scriptstyle\comring \text{ field,} \charac(\comring) \ne 2{,} \\
			\scriptstyle\dim \module < \infty
		\end{gathered}"] \arrow[u, bend right, "\scriptstyle\addtolength{\jot}{-2mm}\begin{gathered}
			\scriptstyle\comring \text{ field,} \\
			\scriptstyle\charac(\comring) \ne 2
		\end{gathered}"']
	\end{tikzcd} $
	\caption{The relationship between the regularity conditions for a quadratic form $ \map{q}{\module}{\comring}{}{} $.}
	\label{fig:reg-props}
\end{figure}

\subsection{Composition Algebras}

With the notions of non-singularity and regularity at hand, we can now define composition algebras over commutative associative rings.

\begin{definition}[Composition algebra, {\cite[19.5]{GPR_AlbertRing}}]\label{conic:compalg-def}
	A \defemph*{composition algebra over $ \comring $}\index{composition algebra} is a $ \comring $-algebra $ \compalg $ which satisfies the following conditions.
	\begin{defenumerate}
		\item \label{conic:compalg-def:proj}$ \compalg $ is projective as a $ \comring $-module.
		
		\item \label{conic:compalg-def:rank}The rank function $ \map{}{\operatorname{spec}(\comring)}{\Nzero \union \compactSet{\infty}}{\mathfrak{p}}{\operatorname{rank}(\compalg_{\mathfrak{p}})} $ is locally constant with respect to the Zariski topology on $ \operatorname{spec}(\comring) $. Here $ \operatorname{spec}(\comring) $ denotes the set of prime ideals of $ \comring $ and $ \compalg_{\mathfrak{p}} $ denotes the localisation of $ \compalg $ at $ \mathfrak{p} $.
		
		\item \label{conic:compalg-def:norm}There exists a non-singular quadratic form $ \map{\compnorm}{\compalg}{\comring}{}{} $ that permits composition in the sense that $ \compnorm(1_\compalg) = 1_\comring $ and $ \compnorm(ab) = \compnorm(a) \compnorm(b) $ for all $ a,b \in \compalg $.
	\end{defenumerate}
	Further, it is called a \defemph*{regular composition algebra}\index{composition algebra!regular} if $ \compnorm $ can be chosen to be a regular quadratic form.
\end{definition}

\begin{remark}[Composition algebras over fields]
	Assume that $ \comring $ is a field. Then Axiom~\thmitemref{conic:compalg-def}{conic:compalg-def:rank} is trivially satisfied because $ \comring $ has only one prime ideal, the zero ideal. Further, every vector space over $ \comring $ is projective. Thus~\thmitemref{conic:compalg-def}{conic:compalg-def:norm} is the only non-trivial axiom over fields. However, even over fields, there is no agreement in the literature on the right definition of composition algebras because the regularity condition on the norm may be chosen differently. See, however, \cref{conic:compalg-class-rem}.
\end{remark}

\begin{note}
	Some authors do not require the existence of an identity element in the definition of composition algebras. In the context of root graded groups, however, this convention is not practical.
\end{note}

Composition algebras fit into the framework of conic algebras as follows.

\begin{proposition}[{\cite[19.12]{GPR_AlbertRing}}]\label{conic:comp-char}
	Let $ \compalg $ be a $ \comring $-algebra. Then the following assertions are equivalent:
	\begin{stenumerate}
		\item $ \compalg $ is a composition algebra (respectively, a regular composition algebra).
		
		\item \label{conic:comp-char:conic}$ \compalg $ is finitely generated projective as a $ \comring $-module and there exists a non-singular (respectively, regular) quadratic form $ \map{\compnorm}{\compalg}{\comring}{}{} $ such that $ (\compalg, \compnorm) $ is a conic alternative algebra.
	\end{stenumerate}
	In this situation, the norm $ \compnorm $ is uniquely determined by the $ \comring $-algebra $ \compalg $. Further, the conic algebra $ (\compalg, \compnorm) $ is norm-associative and its conjugation is a scalar algebra involution.
\end{proposition}
\begin{proof}
	The equivalence of the two assertions is proven in \cite[19.12]{GPR_AlbertRing}. The uniqueness of $ \compnorm $ holds by \cref{conic:norm-unique}. It follows from \cref{conic:mult-is-norm} that the conic algebra is norm-associative and that its conjugation is a scalar algebra involution.
\end{proof}

Observe that the norm in condition~\thmitemref{conic:comp-char}{conic:comp-char:conic} is necessarily multiplicative by \cref{conic:alt-proj-is-mult}.

\begin{lemma}[{\cite[19.11]{GPR_AlbertRing}}]\label{conic:char2-reg}
	Assume that $ \comring $ is a field of characteristic not~2. Then every composition algebra over $ \comring $ is regular.
\end{lemma}

\begin{note}
	Axiom~\thmitemref{conic:compalg-def}{conic:compalg-def:rank} is invoked in \cite[19.11]{GPR_AlbertRing} to show that a composition algebra is finitely generated as a $ \comring $-module. Conversely, any finitely generated projective $ \comring $-module satisfies Axiom~\thmitemref{conic:compalg-def}{conic:compalg-def:rank}. Thus we could replace Axiom~\thmitemref{conic:compalg-def}{conic:compalg-def:rank} by \enquote{$ \compalg $ is finitely generated as a $ \comring $-module} without changing the meaning of the word \enquote{composition algebra}. However, recall that in the classical setting over fields, the possibility of non-finitely generated composition algebras is not ruled out from the beginning, and it is instead a theorem that they must be finitely generated as modules. Axiom~\thmitemref{conic:compalg-def}{conic:compalg-def:rank} makes it possible to carry this result over to the more general situation.
\end{note}

\subsection{Pre-composition Algebras}

In order to discuss composition algebras in the classical situation (that is, over fields), the notion of pre-composition algebras will be useful.

\begin{definition}[Pre-composition algebra, {\cite[19.1]{GPR_AlbertRing}}]
	A \defemph*{pre-composition algebra over $ \comring $}\index{composition algebra!pre-} is a $ \comring $-algebra $ \compalg $ which satisfies the following conditions.
	\begin{defenumerate}
		\item \label{conic:pre-compalg-def:proj}$ \compalg $ is projective as a $ \comring $-module.
		
		\item \label{conic:pre-compalg-def:norm}There exists a non-degenerate quadratic form $ \map{\compnorm}{\compalg}{\comring}{}{} $ that permits composition in the sense that $ \compnorm(1_\compalg) = 1_\comring $ and $ \compnorm(ab) = \compnorm(a) \compnorm(b) $ for all $ a,b \in \compnorm $.
	\end{defenumerate}
\end{definition}

We have the following analogue of \cref{conic:comp-char} for pre-composition algebras.

\begin{proposition}[{\cite[17.6, 19.3]{GPR_AlbertRing}}]\label{conic:pre-comp-char}
	Let $ \compalg $ be a $ \comring $-algebra. Then the following assertions are equivalent:
	\begin{stenumerate}
		\item $ \compalg $ is a pre-composition algebra.
		
		\item \label{conic:pre-comp-char:conic}$ \compalg $ is projective as a $ \comring $-module and there exists a non-degenerate quadratic form $ \map{\compnorm}{\compalg}{\comring}{}{} $ such that $ (\compalg, \compnorm) $ is a conic alternative algebra.
	\end{stenumerate}
	In this situation, the norm $ \compnorm $ is uniquely determined by the $ \comring $-algebra $ \compalg $. Further, the conic algebra $ (\compalg, \compnorm) $ is norm-associative and its conjugation is a scalar algebra involution.
\end{proposition}

Again, the norm in condition~\thmitemref{conic:pre-comp-char}{conic:pre-comp-char:conic} is necessarily multiplicative by \cref{conic:alt-proj-is-mult}.

\begin{remark}
	A composition algebra over $ \comring $ is a pre-composition algebra if and only if its norm is non-degenerate. By \cref{fig:reg-props}, this is always the case if the composition algebra is regular of if $ \comring $ is a field. However, as \cite[19.6]{GPR_AlbertRing} illustrates, it is not true in general.
\end{remark}

Unlike composition algebras, pre-composition algebras are not stable under base change, which makes them significantly less interesting.

\begin{example}[{\cite[19.4]{GPR_AlbertRing}}]\label{conic:char2-ext-ex}
	Let $ K $ be a purely inseparable field extension of a field $ \comring $ of characteristic~2 such that $ K^2 \subs \comring $. Then $ K $ is a pre-composition algebra with anisotropic norm $ \map{}{}{}{x}{x^2} $. This norm is non-degenerate, but not weakly regular because its linearisation is the zero map. Further, there exists a scalar extension of $ \comring $ over which $ K $ is no longer a pre-composition algebra. In particular, $ K $ is not a composition algebra.
\end{example}

\begin{lemma}\label{conic:char2-pre-reg}
	Assume that $ \comring $ is a field of characteristic not~2. Then every pre-composition algebra over $ \comring $ is a regular composition algebra.
\end{lemma}
\begin{proof}
	At first, assume only that $ \comring $ is a field, but make no assumption on its characteristic. In \cite[1.2.1]{SpringerVeldkamp-Octonions}, Springer-Veldkamp define a composition algebra over $ \comring $ to be a pre-composition algebra over $ \comring $ whose norm is weakly regular.\footnote{More precisely, they define a composition algebra to be an algebra which satisfies Axiom~\thmitemref{conic:compalg-def}{conic:compalg-def:norm} for a \enquote{nondegenerate} norm, where by \enquote{nondegenerate} they mean what we call \enquote{weakly regular} (see \cite[\onpage{2}]{SpringerVeldkamp-Octonions}).} If $ \comring $ is not of characteristic~2, then this definition is equivalent to the notion of pre-composition algebras.
	
	Now assume, in addition, that $ \comring $ is not of characteristic~2. Then by the conclusion of the previous paragraph, it follows from \cite[1.6.2]{SpringerVeldkamp-Octonions} that any pre-composition algebra $ \compalg $ over $ \comring $ is finite-dimensional. By \cref{fig:reg-props}, this implies that $ \compalg $ is a regular composition algebra.
\end{proof}

\begin{note}\label{conic:compalg-class-rem}
	Putting \cref{conic:char2-reg,conic:char2-pre-reg} together, we see that the notions of composition algebras, regular composition algebra and pre-composition algebras coincide over fields of characteristic not~2. As \cref{conic:char2-ext-ex} illustrates, this ceases to be true over general fields. See also \cite[1.2.2]{SpringerVeldkamp-Octonions} for more remarks on this topic.
\end{note}

The following result will be useful in the discussion of conic algebras which coordinatise RGD-systems of type $ F_4 $.

\begin{lemma}\label{conic:rgd-alg}
	Assume that $ \comring $ is a field and let $ (\compalg, \compnorm) $ be a conic alternative algebra which is anisotropic as a quadratic space in the sense of \cref{quadmod:quadform-def} (but which is not assumed to be multiplicative). Then $ \compalg $ is a pre-composition algebra in which every non-zero element is invertible. If $ \comring $ is not of characteristic~2, then $ \compalg $ is a regular composition algebra in which every non-zero element is invertible.
\end{lemma}
\begin{proof}
	Since $ \compnorm $ is anisotropic, its radical is zero, which says that $ \compnorm $ is non-degen\-er\-ate. Hence $ \compalg $ is a pre-composition algebra by \cref{conic:pre-comp-char}. If $ \comring $ is not of characteristic~2, then it follows from \cref{conic:char2-pre-reg} that $ \compalg $ is a regular composition algebra. Further, $ \compalg $ is a division algebra by \cref{conic:norm-inv}.
\end{proof}


\section{The Root System \texorpdfstring{$ F_4 $}{F\_4}}

\label{sec:F4:rootsys}

As usual, we begin our discussion of the root system $ F_4 $ with a standard representation, even though it is too cumbersome to be of practical use.

\begin{remark}[Standard representation of $ F_4 $]\label{F4:F-standard-rep}
	Let $ V $ be a Euclidean space of dimension $ n $ with orthonormal basis $ (b_1, \ldots, b_4) $. The \defemph*{standard representation of $ F_4 $}\index{standard representation!of F4@of $ F_4 $} is
	\begin{align*}
		F_4 &\defl \Set{\epsilon_1 b_i +\epsilon_2 b_j \given i \ne j \in \numint{1}{n}, \epsilon_1, \epsilon_2 \in \compactSet{\pm 1}} \union \Set{\epsilon b_i \given i \in \numint{1}{n}, \epsilon \in \compactSet{\pm 1}} \\
		&\hspace{2cm} \mathord{} \union \Set*{\frac{1}{2}\brackets[\big]{\epsilon_1 b_1 + \epsilon_2 b_2 + \epsilon_3 b_3 + \epsilon_4 b_4} \given \listing{\epsilon}{4} \in \compactSet{\pm 1}}.
	\end{align*}
	Its \defemph*{standard root base} consists of
	\begin{align*}
		f_1 \defl b_2 - b_3, \quad f_2 \defl b_3 - b_4, \quad f_3 \defl b_4, \quad f_4 \defl \frac{1}{2} \brackets{b_1 - b_2 - b_3 - b_4},
	\end{align*}
	so that $ f_1 $, $ f_2 $ are long while $ f_3 $, $ f_4 $ are short. In total, $ F_4 $ has 48 roots.
\end{remark}

\begin{figure}[htb]
	\centering\begin{tikzpicture}
		\draw (1,0) -- (2, 0) (3,0) -- (4,0);
		\draw[coxdoublearrow] (2, 0) -- (3,0);
		\foreach \x in {1, ..., 4}{
			\node[coxcircle] at (\x, 0){};
			\node at (\x, -0.4){$ f_\x $};
		}
	\end{tikzpicture}
	\caption{The Dynkin Diagram of $ F_4 $.}
	\label{fig:F4-dynkin}
\end{figure}
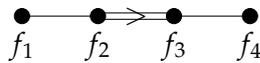

We will always consider root bases of $ F_4 $ to be in the following order, which we call standard.

\begin{definition}[Standard order]
	Let $ \rootbase = (f_1, f_2, f_3, f_4) $ be an ordered root base of $ F_4 $. We say that $ \rootbase $ is in \defemph*{standard order}\index{root base!standard order (for $ F_4 $)} if $ f_1, \ldots, f_4 $ are as in \cref{fig:F4-dynkin}. That is, $ (f_1, f_2) $ and $ (f_3, f_4) $ are $ A_2 $-pairs and $ (f_2, f_3) $ is a $ B_2 $-pair (which, in particular, means that $ f_2 $ is longer than $ f_3 $).
\end{definition}

\begin{definition}[Canonical root subsystems]\label{F4:rootsys:canon-sub}
	Let $ \rootbase = (f_1, f_2, f_3, f_4) $ be a root base of $ F_4 $ in standard order. The \defemph*{canonical subsystem of type $ B_3 $ (with respect to $ \rootbase $)} is the root subsystem generated by $ \Set{f_1, f_2, f_3} $. The \defemph*{canonical subsystem of type $ C_3 $ (with respect to $ \rootbase $)} is the root subsystem generated by $ \Set{f_2, f_3, f_4} $.
\end{definition}

\begin{lemma}\label{F4:rootsys:inBC}
	Every root in $ F_4 $ is contained in a parabolic subsystem of type $ B_3 $, in a parabolic subsystem of type $ C_3 $ and in a parabolic subsystem of type $ A_2 $.
\end{lemma}
\begin{proof}
	By an inspection of \cref{fig:F4-dynkin}, we see that the root $ f_2 $ has the desired properties. Since the Weyl group acts transitively on the set of long roots by \cref{rootsys:cry-orb-length}, it follows that the assertion holds for all long roots. By the same argument, it holds for all short roots because it holds for $ f_3 $.
\end{proof}

In \cref{sec:F4:const}, we will construct examples of $ F_4 $-graded groups as foldings of $ E_6 $-graded groups. In this setting, we will describe the root system $ F_4 $ as follows.

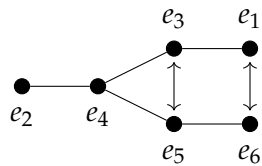
\begin{figure}[htb]
	\centering\begin{tikzpicture}
		\node[coxcircle] (e2) at (0, 0){};
		\node[coxcircle] (e4) at (1, 0){};
		\node[coxcircle] (e3) at (2, 0.5){};
		\node[coxcircle] (e1) at (3, 0.5){};
		\node[coxcircle] (e5) at (2, -0.5){};
		\node[coxcircle] (e6) at (3, -0.5){};
		\draw (e2) -- (e4) -- (e3) -- (e1);
		\draw (e4) -- (e5) -- (e6);
		\draw[<->, shorten > =2pt, shorten < =2pt] (e3) -- (e5);
		\draw[<->, shorten > =2pt, shorten < =2pt] (e1) -- (e6);
		\node at ($ (e2) + (0, -0.4) $){$ e_2 $};
		\node at ($ (e4) + (0, -0.4) $){$ e_4 $};
		\node at ($ (e5) + (0, -0.4) $){$ e_5 $};
		\node at ($ (e6) + (0, -0.4) $){$ e_6 $};
		\node at ($ (e3) + (0, 0.4) $){$ e_3 $};
		\node at ($ (e1) + (0, 0.4) $){$ e_1 $};
	\end{tikzpicture}
	\caption{The non-trivial diagram automorphism $ \rho $ of $ E_6 $.}
	\label{fig:F4:E6-auto}
\end{figure}

\begin{construction}[$ F_4 $ as a folding of $ E_6 $]\label{F4:F4-fold-const}
	Let $ \rootbaseE = (e_1, \ldots, e_6) $ be an ordered root base of $ E_6 $ which is indexed\footnote{This seemingly unnatural ordering of $ \rootbaseE $ agrees with the one in \cite{HumphreysLieAlg}. We use it because it is the one used by the GAP computer algebra system.} as in \cref{fig:F4:E6-auto}, and denote by $ \rho $ the unique non-trivial automorphism of this diagram (also as in \cref{fig:F4:E6-auto}). Denote by $ (V, \cdot) $ the six-dimensional Euclidean space surrounding $ E_6 $. Since all roots of $ E_6 $ have the same length, $ \rho $ is an isometry, so $ \rho = \tau $ in the notation of~\ref{secnot:rootsys-fold}. Denote by $ \fixspace $ the fixed space of $ \tau $. Then
	\begin{align*}
		\fixspace = \gen{e_2, e_4, e_3 + e_5, e_1 + e_6}_\IR \midand \fixspace^\perp = \gen{e_3 - e_5, e_1 - e_6}_\IR.
	\end{align*}
	Denote by $ \map{\pi}{V}{\fixspace}{}{} $ the orthogonal projection on $ \fixspace $ (whose kernel is $ \fixspace^\perp $). Then it follows from \cref{fold:rootsys} that $ \roots' \defl \pi(E_6) $ is a root system of rank~4 in $ \fixspace $ and that $ \rootbase' \defl \pi(\rootbaseE) $ is a root base of $ \roots' $. An inspection shows that $ \roots' $ is of type $ F_4 $ and that the vectors
	\begin{gather*}
		f_1 \defl e_2, \quad f_2 \defl e_4, \quad f_3 \defl \pi(e_3) = \pi(e_5) = \frac{1}{2}(e_3 + e_5), \\
		f_4 \defl \pi(e_1) = \pi(e_6) = \frac{1}{2}(e_1 + e_6)
	\end{gather*}
	are in standard order. See, for example, \cite[13.3.3]{Carter-Chev}.
\end{construction}

\begin{remark}\label{F4:fold-const-rem}
	The long roots in $ F_4 $ are precisely those whose preimage in $ E_6 $ under $ \pi $ contains exactly one root (namely, the root itself) while the short roots in $ F_4 $ are those whose preimage contains exactly two roots (whose common image in $ F_4 $ is their average by \cref{fold:proj-average}). In particular, every long root in $ F_4 $ is also a root in $ E_6 $. Further, the two roots in the preimage of a short root are orthogonal (and thus adjacent). This is clear from \cref{fig:F4:E6-auto} for simple roots and it follows for arbitrary roots by the transitivity properties of the Weyl group.
\end{remark}


\section{Construction of \texorpdfstring{$ F_4 $}{F\_4}-graded Groups via Foldings}

\label{sec:F4:const}

\begin{secnotation}
	We consider the root system $ F_4 $ as a folding of $ E_6 $ as in \cref{F4:F4-fold-const}. In particular, we denote by $ (V, \cdot) $ the Euclidean space surrounding $ E_6 $, by $ \rho = \tau $ the non-trivial diagram automorphism of $ E_6 $, by $ \fixspace $ the fixed space of $ \rho $, by $ \map{\pi}{V}{\fixspace}{}{} $ the orthogonal projection onto $ \fixspace $ and by $ \rootbaseE = (e_1, \ldots, e_6) $ and $ \rootbaseF = (f_1, \ldots, f_4) = \pi(\rootbaseE) $ the corresponding root bases.
	Further, we fix a family $ \chevstr = (\chevstr_{\alpha, \beta})_{\alpha, \beta \in E_6} $ of Chevalley structure constants of type $ E_6 $ and by $ (G, (\rootgr{\alpha})_{\alpha \in E_6}) $ a Chevalley group of type $ E_6 $ over the commutative associative ring $ \compalg $ with respect to $ \chevstr $. We denote the corresponding root isomorphisms by~$ (\risom{\alpha}^E)_{\alpha \in E_6} $.
\end{secnotation}

\begin{reminder}[see \cref{conic:dirprod-example}]
	The commutative associative ring $ \compalg \times \compalg $ is a multiplicative conic algebra over $ \compalg $ with norm $ \compnorm(a,b) \defl ab $, trace $ \compTr(a,b) = a+b $ and conjugation $ \compinv{(a,b)} = (b,a) $ for all $ (a,b) \in \compalg $.
\end{reminder}

We already know from \cref{rgg-fold:thm} that the $ E_6 $-graded group $ G $ has a crystallographic $ F_4 $-grading $ (\varrootgr{\alpha'})_{\alpha' \in F_4} $ given by
\[ \varrootgr{\alpha'} \defl \gen{\rootgr{\alpha} \given \alpha \in E_6, \pi(\alpha) = \alpha'} \]
for all $ \alpha' \in F_4 $. In this section, we want to show that this $ F_4 $-grading is coordinatised by the multiplicative conic algebra $ \compalg \times \compalg $ over $ \compalg $. In the first subsection, we describe the naive approach to this question. Unfortunately, this strategy turns out to be unsuitable because certain signs do not fit together. In the second subsection, we construct a refined coordinatisation of $ (\varrootgr{\alpha})_{\alpha \in F_4} $ by twisting certain root isomorphism of the naive coordinatisation. Since the choice of root isomorphisms which have to be twisted depends on the family $ \chevstr $, we have to fix $ \chevstr $ from this point on. For computational practicality, we choose $ \chevstr $ to be the family of structure constants that is used in the computer algebra system GAP \cite{GAP4}.

In the final subsection, we assume that $ \compalg $ is equipped with the additional structure of a faithful multiplicative conic algebra over a commutative associative ring $ \comring $. We then construct a crystallographic $ F_4 $-graded subgroup of $ G $ which is coordinatised by the conic algebra $ \compalg $. This solves the existence problem for $ F_4 $-graded groups for multiplicative conic alternative algebras which are, in addition, commutative, associative and faithful.

\begin{note}[Implementation]\label{F4:ex:GAP}
	An implementation of the the group that we construct in this section in the computer algebra system GAP \cite{GAP4} can be found in \cite{RGG-GAP}. This includes proofs of all computational results in this section which are not explicitly verified by hand.
\end{note}

\subsection{The Naive Coordinatisation}

At first, we need to introduce a total order on $ F_4 $. There is no abstract reason to prefer any one order over another, but we need to fix some ordering in the definition of the (short) root homomorphisms (see \cref{F4:ex:naiv-hom,F4:ex:refined-coord}).

\begin{reminder}[Lexicographic order]
	Let $ n \in \Npos $. The \defemph*{strict lexicographic order on $ \IR^n $}\index{lexicographic order}, denoted by $ \lexls $, is the strict total order defined as follows: For all $ a = \tup{a}{n}, b = \tup{b}{n} \in \IR^n $, we have $ a \lexls b $ if and only if there exists $ i \in \numint{1}{n} $ such that $ a_j = b_j $ for all $ j \in \numint{1}{i-1} $ and $ a_i < b_i $. Further, the \defemph*{lexicographic order on $ \IR^n $} is the total order $ \lexle $ which is defined as follows: $ a \lexle b $ holds if and only if $ a = b $ or $ a \lexls b $.
\end{reminder}

\begin{definition}[Total order on $ F_4 $]\label{F4:ex:total-order}
	We define a total order $ \Forder $ on the positive system $ \possys_F $ of $ F_4 $ corresponding to $ \rootbaseF $ by
	\begin{align*}
		1000 &\Forder 0100 \Forder 0010 \Forder 0001 \Forder 1100 \Forder 0110 \Forder 0011 \Forder 1110 \Forder 0120 \\
		&\Forder 0111 \Forder 1120 \Forder 1111 \Forder 0121 \Forder 1220 \Forder 1121 \Forder 0122 \Forder 1221 \\
		&\Forder 1122 \Forder 1231 \Forder 1222 \Forder 1232 \Forder 1242 \Forder 1342 \Forder 2342
	\end{align*}
	where a sequence $ a_1 a_2 a_3 a_4 $ represents $ \sum_{i=1}^4 a_i f_i $.
	We extend $ \Forder $ to a total order on $ F_4 $ as follows:
	\begin{stenumerate}
		\item For all $ \alpha \in \possys_F $, $ \beta \in -\possys_F $, we declare that $ \alpha \Forder \beta $.
		
		\item For all $ \alpha, \beta \in -\possys_F $, we declare that $ \alpha \Forder \beta $ if and only if $ -\alpha \Forder -\beta $.
	\end{stenumerate}
\end{definition}

\begin{note}
	The ordering of the positive roots in \cref{F4:ex:total-order} is a height ordering in the sense of \cref{rootorder:prop-def}, and it is also the one used in \cite[\onpage{96}]{VavilovPlotkin}.
\end{note}

\begin{definition}[Naive coordinatisation]\label{F4:ex:naiv-hom}
	The \defemph*{naive coordinatisation of $ (\varrootgr{\alpha'})_{\alpha' \in F_4} $} consists of the isomorphisms
	\[ \map{\risomnaiv{\alpha'}}{(\comring, +)}{\varrootgr{\alpha'} = \rootgr{\alpha'}}{a}{\risom{\alpha'}(a)} \]
	for all long roots $ \alpha' \in F_4 $ and
	\[ \map{\risomnaiv{\beta'}}{(\comring \times \comring, +)}{\varrootgr{\beta'} = \rootgr{\beta_1} \rootgr{\beta_2}}{(a,b)}{\risom{\beta_1}(a) \risom{\beta_2}(b)} \]
	for all short roots $ \beta' \in F_4 $ with corresponding preimage $ \pi^{-1}(\beta') \intersect E_6 = \Set{\beta_1, \beta_2} $ in $ E_6 $ where $ \beta_1 \Forder \beta_2 $.
\end{definition}

By \cref{F4:fold-const-rem}, the maps in \cref{F4:ex:naiv-hom} are well-defined and isomorphisms. Note that we need some total order on $ F_4 $ to decide without ambiguity whether $ (a,b) $ should be mapped to $ \risom{\beta_1}(a) \risom{\beta_2}(b) $ or $ \risom{\beta_2}(a) \risom{\beta_1}(b) $.

\begin{definition}[Naive Weyl elements]
	For all roots $ \alpha \in F_4 $, we define the \defemph*{naive standard $ \alpha $-Weyl element} to be
	\[ \whomnaiv{\alpha} \defl \risomnaiv{-\alpha}(-1_{\alpha}) \risomnaiv{\alpha}(1_{\alpha}) \risomnaiv{-\alpha}(-1_{\alpha}) \]
	where $ 1_{\alpha} = 1_\compalg $ if $ \alpha $ is long and $ 1_{\alpha} = (1_\compalg, 1_\compalg) $ if $ \alpha $ is short.
\end{definition}

\begin{remark}
	Let $ \alpha, \beta, \delta \in F_4 $ such that $ \alpha $ is long and $ \beta $ is short and let $ a,b \in \compalg $. Put $ w \defl \whomnaiv{\delta} $. A straightforward computation (see \cref{F4:ex:GAP}) shows that in all possible cases, we have
	\begin{align*}
		\risomnaiv{\alpha}(a)^w &\in \Set{\risomnaiv{\refl{\delta}(\alpha)}(\epsilon a) \given \epsilon \in \compactSet{\pm 1}} \rightand \\
		\risomnaiv{\beta}(a,b)^w &\in \Set{\risomnaiv{\refl{\delta}(\beta)}(\epsilon a, \sigma b), \risomnaiv{\refl{\delta}(\beta)}(\epsilon b, \sigma a) \given \epsilon, \sigma \in \compactSet{\pm 1}}.
	\end{align*}
	This suggests that the twisting group for the coordinatisation of $ (\varrootgr{\alpha'})_{\alpha' \in F_4} $ should be $ \IZ_2^3 $ where the actions of the three components on $ \compalg \times \compalg $ should be given by $ \map{}{}{}{(a,b)}{(-a,-b)} $, $ \map{}{}{}{(a,b)}{(b,a)} $ and $ \map{}{}{}{(a,b)}{(-a,b)} $ (or $ \map{}{}{}{(a,b)}{(a,-b)} $). However, this is a fallacy: By twisting the coordinatisation $ (\risomnaiv{\alpha'})_{\alpha' \in F_4} $ appropriately, we can arrange it so that only the inversion action $ \map{}{}{}{(a,b)}{(-a,-b)} $ and the involution $ \map{}{}{}{(a,b)}{(b,a)} $ are required. Note that the latter is precisely the conjugation of the conic algebra $ \compalg \times \compalg $.
\end{remark}

\subsection{A Refined Approach}

\begin{secnotation}
	From now on, we assume that $ \chevstr $ is the family of structure constants of the simple Lie algebra over $ \IQ $ of type $ E_6 $ which is constructed in GAP \cite{GAP4} with the command \texttt{SimpleLieAlgebra(``E'', 6, Rationals)}.
\end{secnotation}

\begin{definition}[Refined coordinatisation]\label{F4:ex:refined-coord}
	Put
	\begin{align*}
		\calT' &\defl \Set*{\begin{gathered}
			f_2 + f_3, f_3 + f_4, f_1 + f_2 + f_3, f_2 + 2f_3 + f_4, f_1 + f_2 + 2f_3 + f_4, \\
			f_1 + 2f_2 + 2f_3 + f_4, f_1 + 2f_2 + 3f_3 + 2f_4
		\end{gathered}} \subs F_4, \\
		\calT &\defl \calT' \union (-\calT').
	\end{align*}
	The \defemph*{refined coordinatisation of $ (G, (\varrootgr{\alpha'})_{\alpha' \in F_4}) $} consists of the isomorphisms
	\[ \map{\risom{\alpha} \defl \risomnaiv{\alpha'}}{(\comring, +)}{\varrootgr{\alpha'} = \rootgr{\alpha'}}{a}{\risom{\alpha'}(a)} \]
	for all long roots $ \alpha' \in F_4 $ and
	\[ \map{\risom{\beta'}}{(\comring \times \comring, +)}{\varrootgr{\beta'} = \rootgr{\beta_1} \rootgr{\beta_2}}{(a,b)}{\begin{cases}
		\risom{\beta_1}(a) \risom{\beta_2}(b) & \text{if } \beta' \nin \calT, \\
		\risom{\beta_1}(a) \risom{\beta_2}(-b) & \text{if } \beta' \in \calT
	\end{cases}} \]
	for all short roots $ \beta' \in F_4 $ with corresponding preimage $ \pi^{-1}(\beta') \intersect E_6 = \Set{\beta_1, \beta_2} $ in $ E_6 $ where $ \beta_1 \Forder \beta_2 $.
\end{definition}

\begin{note}
	The set $ \calT' $ in \cref{F4:ex:refined-coord} is chosen precisely so that the following \cref{F4:ex:parmap-exist} holds. The specific form of $ \calT' $ has been determined by trial and error.
\end{note}

\begin{definition}[Refined Weyl elements]
	For all roots $ \alpha \in F_4 $, we define the \defemph*{refined standard $ \alpha $-Weyl element} to be
	\[ \whom{\alpha} \defl \risom{-\alpha}(-1_{\alpha}) \risom{\alpha}(1_{\alpha}) \risom{-\alpha}(-1_{\alpha}) \]
	where $ 1_{\alpha} = 1_\compalg $ if $ \alpha $ is long and $ 1_{\alpha} = (1_\compalg, 1_\compalg) $ if $ \alpha $ is short.
\end{definition}

\begin{proposition}\label{F4:ex:parmap-exist}
	Denote by $ (\twistgroup \times \invogroup, \compalg \times \compalg, \compalg) $ the standard parameter system for the conic algebra $ \compalg \times \compalg $ over $ \compalg $ in the sense of \cref{conic:standard-param}.
	Then there exist unique $ \rootbaseF $-parity maps $ \map{\inverparsym}{F_4 \times \rootbaseF}{\twistgroup}{}{} $ and $ \map{\invoparsym}{F_4 \times \rootbaseF}{\invogroup}{}{} $ such that $ (\risom{\alpha})_{\alpha \in F_4} $ is a parametrisation of $ (G, (\varrootgr{\alpha})_{\alpha \in F_4}) $ by $ (\twistgroup \times \invogroup, \compalg \times \compalg, \compalg) $ with respect to $ (w_\delta)_{\delta \in \rootbaseF} $ and $ \inverparsym \times \invoparsym $ and such that $ \invopar{\alpha}{\delta} = 1_\invogroup $ for all long roots $ \alpha \in F_4 $ and all $ \delta \in \rootbaseF $. Further, these maps satisfy $ \inverparbr{\alpha}{\delta} = \inverparbr{-\alpha}{\delta} $ and $ \invoparbr{\alpha}{\delta} = \invoparbr{-\alpha}{\delta} $ for all $ \alpha \in F_4 $ and $ \delta \in \rootbaseF $.
\end{proposition}
\begin{proof}
	For the existence, we have to verify that
	\begin{align*}
		\risom{\alpha}(a)^{w_\delta} &\in \Set{\risom{\refl{\delta}(\alpha)}(\epsilon a) \given \epsilon \in \compactSet{\pm 1}} \rightand \\
		\risom{\beta}(a,b)^{w_\delta} &\in \Set{\risom{\refl{\delta}(\beta)}(\epsilon a, \epsilon b), \risom{\refl{\delta}(\beta)}(\epsilon b, \epsilon a) \given \epsilon \in \compactSet{\pm 1}}.
	\end{align*}
	for all short roots $ \beta \in F_4 $, for all long roots $ \alpha \in F_4 $, for all $ \delta \in \rootbaseF $ and all $ a,b \in \compalg $. This is a straightforward computation (see \cref{F4:ex:GAP}). The uniqueness of $ \inverparsym $ and $ \invoparsym $ follows from the fact that $ \twistgroup $ and $ \invogroup $ act faithfully on $ \compalg \times \compalg $, that $ \twistgroup $ acts faithfully on $ \compalg $ and from the additional requirement that $ \invopar{\alpha}{\delta} = 1_\invogroup $ for all long roots $ \alpha \in F_4 $ and all $ \delta \in \rootbaseF $. The property that $ \inverparbr{\alpha}{\delta} = \inverparbr{-\alpha}{\delta} $ and $ \invoparbr{\alpha}{\delta} = \invoparbr{-\alpha}{\delta} $ for all $ \alpha \in F_4 $ and $ \delta \in \rootbaseF $ is also a mere computation.
\end{proof}

The precise values of $ \inverparsym $ and $ \invoparsym $ can be found in \cref{fig:F4:parmap-table}.

\premidfigure
\begin{figure}[tb]
	\centering\begin{tabular}{ccccc}
		\toprule
		$ \alpha $ & $ \totalparbr{\alpha}{f_1} $ & $ \totalparbr{\alpha}{f_2} $ & $ \totalparbr{\alpha}{f_3} $ & $ \totalparbr{\alpha}{f_4} $ \\
		\midrule
		$ 1000 $ & $ (-1, 1) $ & $ (-1, 1) $ & $ (1, 1) $ & $ (1, 1) $ \\
		$ 0100 $ & $ (1, 1) $ & $ (-1, 1) $ & $ (-1, 1) $ & $ (1, 1) $ \\
		$ 0010 $ & $ (1, 1) $ & $ (1, 1) $ & $ (-1, 1) $ & $ (-1, 1) $ \\
		$ 0001 $ & $ (1, 1) $ & $ (1, 1) $ & $ (1, 1) $ & $ (-1, 1) $ \\
		$ 1100 $ & $ (-1, 1) $ & $ (1, 1) $ & $ (-1, 1) $ & $ (1, 1) $ \\
		\midrule
		$ 0110 $ & $ (1, 1) $ & $ (-1, 1) $ & $ (-1, -1) $ & $ (-1, 1) $ \\
		$ 0011 $ & $ (1, 1) $ & $ (1, 1) $ & $ (-1, 1) $ & $ (1, 1) $ \\
		$ 1110 $ & $ (-1, 1) $ & $ (1, 1) $ & $ (-1, -1) $ & $ (-1, 1) $ \\
		$ 0120 $ & $ (1, 1) $ & $ (1, 1) $ & $ (-1, 1) $ & $ (-1, 1) $ \\
		$ 0111 $ & $ (1, 1) $ & $ (-1, 1) $ & $ (1, 1) $ & $ (1, 1) $ \\
		\midrule
		$ 1120 $ & $ (-1, 1) $ & $ (-1, 1) $ & $ (-1, 1) $ & $ (-1, 1) $ \\
		$ 1111 $ & $ (-1, 1) $ & $ (1, 1) $ & $ (1, 1) $ & $ (1, 1) $ \\
		$ 0121 $ & $ (1, 1) $ & $ (1, 1) $ & $ (-1, 1) $ & $ (-1, -1) $ \\
		$ 1220 $ & $ (1, 1) $ & $ (1, 1) $ & $ (1, 1) $ & $ (-1, 1) $ \\
		$ 1121 $ & $ (-1, 1) $ & $ (-1, 1) $ & $ (-1, 1) $ & $ (-1, -1) $ \\
		\midrule
		$ 0122 $ & $ (1, 1) $ & $ (1, 1) $ & $ (1, 1) $ & $ (-1, 1) $ \\
		$ 1221 $ & $ (1, 1) $ & $ (1, 1) $ & $ (-1, 1) $ & $ (-1, -1) $ \\
		$ 1122 $ & $ (-1, 1) $ & $ (-1, 1) $ & $ (1, 1) $ & $ (-1, 1) $ \\
		$ 1231 $ & $ (1, 1) $ & $ (1, 1) $ & $ (1, 1) $ & $ (1, 1) $ \\
		$ 1222 $ & $ (1, 1) $ & $ (1, 1) $ & $ (-1, 1) $ & $ (-1, 1) $ \\
		\midrule
		$ 1232 $ & $ (1, 1) $ & $ (1, 1) $ & $ (-1, -1) $ & $ (-1, 1) $ \\
		$ 1242 $ & $ (1, 1) $ & $ (-1, 1) $ & $ (-1, 1) $ & $ (1, 1) $ \\
		$ 1342 $ & $ (-1, 1) $ & $ (1, 1) $ & $ (1, 1) $ & $ (1, 1) $ \\
		$ 2342 $ & $ (1, 1) $ & $ (1, 1) $ & $ (1, 1) $ & $ (1, 1) $ \\
		\bottomrule
	\end{tabular}
	\caption{The values of $ \inverparsym $ and $ \invoparsym $ in \cref{F4:ex:parmap-exist}, where $ \totalparsym \defl \inverparsym \times \invoparsym $. In the first column, a sequence $ a_1 a_2 a_3 a_4 $ represents $ \sum_{i=1}^4 a_i f_i $. Further, $ \totalparbr{\alpha}{\delta} = \totalparbr{-\alpha}{\delta} $ for all $ \alpha \in F_4 $ and $ \delta \in \rootbaseF $.}
	\label{fig:F4:parmap-table}
\end{figure}
\postmidfigure

\begin{definition}[Standard parity maps]\label{F4:ex:parmap-def}
	Let $ \roots' $ be any root system of type $ F_4 $ and let $ \rootbase' = (f_1', f_2', f_3', f_4') $ be any root base of $ \roots' $ in standard order. Denote by $ \map{\phi}{\roots'}{F_4}{}{} $ the unique isomorphism of root systems with $ \phi(f_i') = f_i $ for all $ i \in \numint{1}{4} $. Then the \defemph*{standard $ \rootbase' $-parity maps for $ \roots' $}\index{parity map!standard!type F4@type $ F_4 $} are the $ \rootbase' $-parity maps $ \map{\inverparsym', \invoparsym'}{\roots' \times \rootbase'}{\compactSet{\pm 1}}{}{} $ defined by
	\[ \inverparsym'\brackets{\alpha, \delta} \defl \inverparsym\brackets[\big]{\phi(\alpha), \phi(\delta)} \midand \invoparsym'\brackets{\alpha, \delta} \defl \invoparsym\brackets[\big]{\phi(\alpha), \phi(\delta)} \]
	for all $ \alpha \in \roots' $ and $ \delta \in \rootbase' $ where $ \inverparsym $, $ \invoparsym $ are as in \cref{F4:ex:parmap-exist}.
\end{definition}

Alternatively, the standard parity maps with respect to any root base in standard order can explicitly be defined to be the values given by \cref{fig:F4:parmap-table}.

We can now compute the commutator relations in $ G $ with respect to $ (\risom{\alpha})_{\alpha \in F_4} $.

\begin{definition}[Coordinatisation with standard signs]\label{F4:ex:coord:stsigns}
	Denote by $ \rootbase' = (f_1', f_2', f_3', f_4') $ any root base of $ F_4 $ in standard order, let $ H $ be a group with an $ F_4 $-pregrading $ (\rootgr{\alpha}')_{\alpha \in F_4} $ and let $ (w_\delta)_{\delta \in \rootbase} $ be a $ \rootbase $-system of Weyl elements. Let $ (\varcompalg, \compnorm) $ be a multiplicative conic alternative algebra with conjugation $ \compinvmap $ and let $ (\twistgroup \times \invogroup, \varcompalg, \comring) $ denote the standard parameter system for $ (\varcompalg, \compnorm) $. A \defemph*{coordinatisation of $ H $ by $ (\varcompalg, \compnorm) $ with standard signs with respect to $ (w_\delta)_{\delta \in \rootbase} $}\index{coordinatisation of a root graded group!with standard signs!type F4@type $ F_4 $} is a family $ (\risom{\alpha}')_{\alpha \in F_4} $ of maps such that the following conditions are satisfied:
	\begin{stenumerate}
		\item \label{F4:ex:coord:stsigns:param}$ (\risom{\alpha}')_{\alpha \in F_4} $ is a parametrisation of $ (H, (\rootgr{\alpha}')_{\alpha \in F_4}) $ by $ (\twistgroup \times \invogroup, \varcompalg, \comring) $ with respect to $ (w_\delta)_{\delta \in \rootbaseF} $ and $ \inverparsym \times \invoparsym $. Here $ \inverparsym $ and $ \invoparsym $ denote the standard $ \rootbase' $-parity maps of~$ F_4 $ (from \cref{F4:ex:parmap-def}).
		
		\item \label{F4:ex:coord:stsigns:commrel}The following commutator relations hold for all $ \lambda, \mu \in \comring $ and all $ c,d \in \varcompalg $:
		\begin{align*}
			\commutator{\risom{f_1}(\lambda)}{\risom{f_2}(\mu)} &= \risom{f_1 + f_2}(-\lambda \mu), \\
			\commutator{\risom{f_2}(\lambda)}{\risom{f_3}(c)} &= \risom{f_2 + f_3}(-c\str_\compalg(\lambda)) \risom{f_2 + 2f_3}\brackets[\big]{-\lambda\compnorm(c)}, \\
			\commutator{\risom{f_2 + f_3}(c)}{\risom{f_3}(d)} &= \risom{f_2 + 2f_3}\brackets[\big]{\compnorm(c,d)}, \\
			\commutator{\risom{f_4}(c)}{\risom{f_3}(d)} &= \risom{f_3 + f_4}(cd).
		\end{align*}
		In the third relation, $ \map{\compnorm}{\varcompalg \times \varcompalg}{\comring}{}{} $ denotes the linearisation of $ \compnorm $. These relations are called the \defemph*{standard commutator relations}.
	\end{stenumerate}
\end{definition}

\begin{note}\label{F4:ex:coord:stsigns:note}
	For all root systems $ \roots $ except for $ F_4 $, we defined standard coordinatisations by prescribing formulas for the commutator $ \commutator{\rootgr{\alpha}}{\rootgr{\beta}} $ for all pairs $ (\alpha, \beta) $ of non-proportional roots. The number of such pairs is $ \abs{\roots} \cdot (\abs{\roots}-2) $ (if $ \roots $ is reduced), which is a large number even if $ \roots $ is \enquote{small}. Until now, we could always bypass this problem by partitioning the set of root pairs into a manageable number of subsets and giving one commutator formula which holds for all pairs in a given subset. For example, we needed 13 formulas to describe the commutator relations in $ B_n $ (\cref{B:ex-commrel}), and only one formula to describe simply-laced root gradings (\cref{ADE:param-def}). In these examples, we did not have to require that a fixed family of Weyl elements acts on the root groups in a certain way because this is a consequence of the commutator formulas (see \namecrefs{ADE:stsign:weyl-char}~\thmitemref{ADE:stsign:weyl-char}{ADE:stsign:weyl-char:param}, \thmitemref{B:stsign:weyl-char}{B:stsign:weyl-char:param} and~\thmitemref{BC:stsign:weyl-char}{BC:stsign:weyl-char:param}).
	
	For $ F_4 $, we are not aware of a way to partition the set of root pairs in a similar manner. Instead, we only describe the commutator relations on a small number of pairs and require in addition that the coordinatisation is compatible with a certain $ \rootbase' $-system of Weyl elements. For any pair $ (\alpha, \beta) $ of non-proportional roots in $ F_4 $, there exists an element $ u $ of the Weyl group such that $ (\alpha^u, \beta^u) $ is pair with a prescribed (standard) commutator relation. Conjugating the standard commutator relation by an appropriate sequence of Weyl elements, we can thus compute the commutator relation for an arbitrary pair $ (\alpha, \beta) $. 
\end{note}

\begin{remark}[Weyl elements in standard coordinatisations]\label{F4:stsign:weyl-char}
	Let $ (\varcompalg, \compnorm) $ be a multiplicative conic alternative algebra over an associative commutative ring $ \comring $. Let $ (H, (\rootgrV{\gamma})_{\gamma \in F_4}) $ be an $ F_4 $-graded group and let $ (\risom{\gamma}^H)_{\gamma \in F_4} $ be a coordinatisation of $ H $ by $ (\varcompalg, \compnorm) $ with standard signs (and with respect to some family of Weyl elements). Recall from \cref{conic:norm-inv} that the set $ \ringinvset{\varcompalg} $ of invertible elements in $ \varcompalg $ is precisely $ \Set{b \in \varcompalg \given \compnorm(b) \in \ringinvset{\comring}} $. By a similar computation as in \cref{B:stsign:weyl-char,BC:stsign:weyl-char}, the maps
	\[ \map{\risom{\alpha}}{\ringinvset{\comring}}{\invsetV{\alpha}}{}{} \midand \map{\risom{\beta}}{\ringinvset{\varcompalg}}{\invsetV{\beta}}{}{} \]
	are well-defined bijections for all long roots $ \alpha $ and all short roots $ \beta $. It follows that if $ H $ is an RGD-system, then $ \comring $ is a field and $ \compnorm $ is anisotropic. By \cref{conic:rgd-alg}, this implies that $ (\compalg, \compnorm) $ is a pre-composition division algebra, and it is even a regular composition algebra if $ \comring $ is not of characteristic~2.
\end{remark}

\begin{proposition}\label{F4:ex:refined-coord-thm}
	The $ F_4 $-graded group $ (G, (\varrootgr{\alpha})_{\alpha \in F_4}) $ is coordinatised by the conic algebra $ \compalg \times \compalg $ over $ \compalg $ with standard signs and with respect to $ (w_\delta)_{\delta \in \rootbase} $.
\end{proposition}
\begin{proof}
	We already know from \cref{F4:ex:parmap-exist} that Axiom~\thmitemref{F4:ex:coord:stsigns}{F4:ex:coord:stsigns:param} is satisfied. It is a straightforward computation (see \cref{F4:ex:GAP}) to show that the desired commutator relations hold.
\end{proof}

\begin{remark}\label{F4:parmap-values}
	Let $ \inverparsym $ and $ \invoparsym $ denote the standard $ \rootbase $-parity maps for $ F_4 $. Only the following specific values of these parity maps will be needed in \cref{sec:F4:param,sec:F4:coord}:
	\begin{align*}
		\inverinvopar{f_2}{f_1} &= (1,1), & \inverinvopar{f_2}{f_1 f_1} &= (-1, 1), \\
		\inverinvopar{f_2}{f_1 f_2} &= (1,1), & \inverinvopar{f_2}{f_3} &= (-1, 1), \\
		\inverinvopar{f_3}{f_4} &= (-1,1), & \inverinvopar{f_3}{f_4 f_4} &= (-1, 1) \\
		\inverinvopar{f_3}{f_2} &= (1,1), & \inverinvopar{f_3}{f_4 f_3} &= (1,1), \\
		\inverinvopar{f_3}{(-f_2, f_3, f_2)} &= (-1, -1).
	\end{align*}
	They can be read off from \cref{fig:F4:parmap-table}, or follow from a short computation.
\end{remark}

The following result will be generalised in \cref{F4:parmap:ortho} and then used to apply \cref{param:stabcomp-crit-ortho}.

\begin{lemma}\label{F4:parmap-ortho-base}
	Let $ \alpha \in F_4 $ and let $ \delta \in \rootbase $ such that $ \alpha \cdot \delta = 0 $. If $ \alpha $ is crystallographically adjacent to $ \delta $, then $ \invopar{\alpha}{\delta} = 1_\invogroup $. Otherwise, $ \invopar{\alpha}{\delta} = -1_\invogroup $.
\end{lemma}
\begin{proof}
	If $ \alpha $ is crystallographically adjacent to $ \delta $, then it is also crystallographically adjacent to $ -\delta $, so $ w_\delta $ acts trivially on $ \rootgr{\alpha} $. Hence $ \invopar{\alpha}{\delta} = 1_\invogroup $ in this case. The second assertion follows from an inspection of \cref{fig:F4:parmap-table}.
\end{proof}

\subsection{A More General Construction}

\begin{secnotation}
	From now on, we assume that $ \comring $ is a commutative associative ring and that $ \compalg $ is equipped with the structure of a multiplicative conic algebra over $ \comring $ which is faithful, commutative and associative. We denote the norm on $ \compalg $ by $ \compnorm $ and its conjugation by $ \compinvmap $. Further, we denote by $ \tilde{\comring} $ the image of $ \comring $ in $ \compalg $ and by $ \tilde{\compalg} $ the subalgebra $ \Set{(a,\compinv{a}) \given a \in \compalg} $ of $ \compalg \times \compalg $.
\end{secnotation}

In \cref{F4:ex:refined-coord-thm}, we have solved the existence problem for the specific kind of conic algebra from \cref{conic:dirprod-example}. The goal of this subsection is to construct a subgroup of $ G $ which is coordinatised by the conic algebra $ \tilde{\compalg} $ over $ \tilde{\comring} $. Since $ \tilde{\compalg} $ can be identified with $ \compalg $ and $ \tilde{\comring} $ can be identified with $ \comring $, this solves the existence problem for multiplicative conic algebras which are faithful, commutative and associative. Observe that the faithfulness assumption is always satisfied by \cref{conic:2-faithful} if $ 2_\comring $ is not a zero divisor in $ \comring $.

\begin{definition}
	For all long roots $ \alpha $ in $ F_4 $, we denote by $ \varvarrootgr{\alpha} $ the image of $ \tilde{\comring} $ under $ \risom{\alpha} $ in $ \varrootgr{\alpha} $. For all short roots $ \beta $, we denote by $ \varvarrootgr{\beta} $ the image of $ \tilde{\compalg} $ under $ \risom{\beta} $ in $ \varrootgr{\beta} $. Further, we denote by $ \varvarrisom{\alpha} $ and $ \varvarrisom{\beta} $ the induced isomorphisms
	\[ \map{\varvarrisom{\alpha} \defl \risom{\alpha} \circ \str_\compalg}{\comring}{\varvarrootgr{\alpha}}{}{} \midand \map{\varvarrisom{\beta}}{\compalg}{\varvarrootgr{\beta}}{a}{\risom{\beta}(a, \compinv{a})} \]
	and by $ \bar{G} $ the subgroup of $ G $ which is generated by $ (\varvarrootgr{\gamma})_{\gamma \in F_4} $.
\end{definition}

\begin{proposition}
	The pair $ (\bar{G}, (\varvarrootgr{\gamma})_{\gamma \in F_4}) $ is an $ F_4 $-graded group which is coordinatised by the conic algebra $ (\compalg, \compnorm) $ over $ \comring $ with standard signs and with respect to $ (w_\delta)_{\delta \in \rootbase} $.
\end{proposition}
\begin{proof}
	Observe that the Weyl elements $ (w_\delta)_{\delta \in \rootbase} $ do indeed lie in $ \bar{G} $. Further, $ \tilde{\comring} $ and $ \tilde{\compalg} $ are invariant under the action of the twisting group $ \twistgroup \times \invogroup $ in the standard parameter system $ (\twistgroup \times \invogroup, \compalg \times \compalg, \compalg) $ for $ \compalg \times \compalg $ and $ \compalg $. It follows that
	\begin{equation}\label{eq:F4:ex:compalg-param}
		\varvarrootgr{\gamma}^{w_\delta} = \varvarrootgr{\refl{\delta}(\gamma)}
	\end{equation}
	for all $ \gamma \in F_4 $ and all $ \delta \in \rootbase $. Further, $ \bar{G} $ is generated by $ (\varvarrootgr{\gamma})_{\gamma \in F_4} $ by definition, and $ \bar{G} $ satisfies Axiom~\thmitemref{rgg-def}{rgg-def:nondeg} because each root group of $ \bar{G} $ is contained in a root group of $ G $.
	
	We now verify the commutator relations. For all $ a,b,c \in \compalg $, we know from \cref{F4:ex:refined-coord-thm} that the commutator relation
	\[ \commutator{\risom{f_2}(a)}{\risom{f_3}(b,c)} = \risom{f_2 + f_3}\brackets[\big]{-(b,c) \cdot (a,a)} \risom{f_2 + 2f_3}(-abc) \]
	holds. This implies that
	\begin{align*}
		\commutator{\varvarrisom{f_2}(\lambda)}{\varvarrisom{f_3}(b)} &= \commutator{\risom{f_2}(\lambda 1_\compalg)}{\risom{f_3}(b, \compinv{b})} = \risom{f_2 + f_3}(-\lambda b, -\lambda \compinv{b}) \risom{f_2 + 2f_3}(-\lambda b \compinv{b}) \\
		&= \varvarrisom{f_2 + f_3}\brackets[\big]{-b \str_{\compalg}(\lambda)} \varvarrisom{f_2 + 2f_3}\brackets[\big]{-\lambda \compnorm(b)}
	\end{align*}
	for all $ \lambda \in \comring $ and $ b \in \compalg $. Here we used that $ b \compinv{b} = \compnorm(b)1_\compalg $ by \thmitemcref{conic:basic-id}{conic:basic-id:scalar}. Thus the desired commutator relation for $ \commutator{\varvarrootgr{f_2}}{\varvarrootgr{f_3}} $ is satisfied. Further, we know that
	\[ \commutator{\risom{f_2+f_3}(a,b)}{\risom{f_3}(c,d)} = \risom{f_2 + 2f_3}\brackets[\big]{\scp[\big]{(a,b)}{(c,d)}} \]
	for all $ a,b,c,d \in \compalg $ where $ \scp{}{} $ denotes the linearisation of the norm on $ \compalg \times \compalg $. By \cref{conic:dirprod-example}, $ \scp[\big]{(a,b)}{(c,d)} = ad+bc $. We infer that
	\begin{align*}
		\commutator[\big]{\varvarrisom{f_2+f_3}(a)}{\varvarrisom{f_3}(c)} &= \commutator{\risom{f_2 + f_3}(a, \compinv{a})}{\risom{f_3}(c, \compinv{c})} = \risom{f_2 + 2f_3}(a \compinv{c} + \compinv{a}c) \\
		&= \risom{f_2 + 2f_3}\brackets[\big]{\compnorm(a,c)1_\compalg} = \varvarrisom{f_2 + 2f_3}\brackets[\big]{\compnorm(a,c)}
	\end{align*}
	for all $ a,b \in \compalg $. Here we used that $ a \compinv{c} + \compinv{a}c = \compnorm(a,c)1_\compalg $ by \thmitemcref{conic:basic-id}{conic:basic-id:scalar-lin} Hence the desired commutator relation for $ \commutator{\varvarrootgr{f_2 + f_3}}{\varvarrootgr{f_3}} $ is satisfied as well. The commutator relations for $ \commutator{\varvarrootgr{f_1}}{\varvarrootgr{f_2}} $ and $ \commutator{\varvarrootgr{f_4}}{\varvarrootgr{f_3}} $ can be verified in a similar, straightforward manner. Together with~\eqref{eq:F4:ex:compalg-param}, it follows that $ \bar{G} $ has crystallographic $ F_4 $-commutator relations with root groups $ (\varvarrootgr{\gamma})_{\gamma \in F_4} $.
	
	Finally, we have to show that the coordinatisation $ (\varvarrisom{\gamma})_{\gamma \in F_4} $ is compatible with the Weyl elements $ (w_\delta)_{\delta \in \rootbase} $ and the standard $ \rootbase $-parity maps $ \inverparsym $, $ \invoparsym $. That is, we have to show that
	\[ \varvarrisom{\alpha}(\lambda)^{w_\delta} = \varvarrisom{\refl{\delta}(\alpha)}\brackets[\big]{\inverpar{\alpha}{\delta} \invopar{\alpha}{\delta}.\lambda} \midand  \varvarrisom{\beta}(a)^{w_\delta} = \varvarrisom{\refl{\beta}(\alpha)}\brackets[\big]{\inverpar{\alpha}{\delta} \invopar{\alpha}{\delta}.a} \]
	for all long roots $ \alpha $, all short roots $ \beta $, all $ \delta \in \rootbase $, all $ \lambda \in \comring $ and all $ a \in \compalg $. Since a similar statement is true for the coordinatisation $ (\risom{\gamma})_{\gamma \in F_4} $, we only have to verify that the restriction of the standard parameter system $ (\twistgroup \times \invogroup, \compalg \times \compalg, \compalg) $ to $ (\tilde{\compalg}, \tilde{\comring}) $ is the standard parameter system for $ \tilde{\compalg} $. Evidently, the inversion action on $ \compalg \times \compalg $ and $ \compalg $ restricts to the inversion action on $ \tilde{\compalg} $ and $ \tilde{\comring} $. Further, the switching involution $ \map{}{}{}{(a,b)}{(b,a)} $ restricts to the conjugation $ \map{}{}{}{(a, \compinv{a})}{(\compinv{a}, \compinv{\compinv{a}}) = (\compinv{a}, a)} $ on $ \tilde{\compalg} $. This finishes the proof.
\end{proof}


\section{Computations in \texorpdfstring{$ F_4 $}{F\_4}-graded Groups}

\label{sec:F4:comp}

\begin{secnotation}
	We fix a root base $ \rootbase $ of $ F_4 $.
\end{secnotation}

Most rank-2 and rank-3 computations that are necessary to parametrise $ F_4 $-graded groups have already been performed in the \cref{chap:B,chap:BC}. In this section, we collect a small number of additional results that are required to apply the parametrisation theorem for $ F_4 $-graded groups.

We begin with a verification of one of the conditions in \cref{param:stabcomp-crit-ortho}. We have already shown that it holds for $ \roots \in \Set{B_3, C_3} $. By a reduction to rank-3 subsystems, we can lift the result to $ F_4 $.

\newcommand{\preworddelta}[1]{\prescript{}{#1}{\word{\delta}}}

\begin{lemma}\label{F4:word-same-act}
	Let $ \roots \in \Set{B_3, C_3, F_4} $, let $ \rootbase $ be a root base of $ \roots $ and let $ \alpha \in \roots $. As in \cref{param:stabcomp-crit-ortho}, we define sets
	\begin{align*}
		\calO &\defl \Set{\beta \in \roots \given \alpha \cdot \beta = 0}, \\
		\calA &\defl \Set{\beta \in \calO \given \alpha \text{ is (crystallographically) adjacent to } \beta \text{ and } -\beta}, \\
		\bar{\calA} &\defl \calO \setminus \calA.
	\end{align*}
	Then for all $ \rootbase $-positive roots $ \beta_1, \beta_2 \in \bar{\calO} $, there exist words $ \preworddelta{1} $, $ \preworddelta{2} $ over $ \rootbase \union (-\rootbase) $ such that $ \reflbr{\preworddelta{1}} = \reflbr{\beta_1} $, $ \reflbr{\preworddelta{2}} = \reflbr{\beta_2} $ and such that for any group $ G $ which has crystallographic $ \roots $-commutator relations with root groups $ (\rootgr{\rho})_{\rho \in F_4} $ and for any $ \rootbase $-system $ (w_\delta)_{\delta \in \rootbase} $ of Weyl elements in $ G $, the actions of $ w_{\preworddelta{1}} $ and $ w_{\preworddelta{2}} $ on $ \rootgr{\alpha} $ are identical.
\end{lemma}
\begin{proof}
	If $ \alpha $ is long, then $ \bar{\calA} $ is empty, so we can assume that $ \alpha $ is short. If $ \roots = C_3 $, then $ \bar{\calA} $ contains exactly one element $ \rootbase $-positive root (see \cref{BC:rootsys:ortho-adj}), so the desired statement is trivially satisfied (because we can choose $ \preworddelta{1} = \preworddelta{2} $). For $ \roots = B_3 $, the desired statement holds by \cref{B:short-weyl-same-action}.

	We proceed to show that the claim holds for $ \roots = F_4 $ as well. Throughout the proof, we denote by $ G $ a group which has crystallographic $ F_4 $-commutator relations with root groups $ (\rootgr{\rho})_{\rho \in F_4} $ and by $ (w_\delta)_{\delta \in \rootbase} $ a $ \rootbase $-system of Weyl elements in $ G $. We will make sure that all constructions in this proof are independent of the choices of $ G $, $ (\rootgr{\rho})_{\rho \in F_4} $ and $ (w_\delta)_{\delta \in \rootbase} $. Denote by $ \roots' $ a parabolic subsystem of $ F_4 $ of rank~3 which contains $ \alpha, \beta_1, \beta_2 $. Note that $ \roots' $ must be of type $ A_2 \times A_1 $, $ B_3 $ or $ C_3 $. In the first case, we can simply take $ \preworddelta{1} $ and $ \preworddelta{2} $ to be $ \rootbase $-expressions of $ \beta_1 $ and $ \beta_2 $, respectively, and then both $ w_{\preworddelta{1}} $ and $ w_{\preworddelta{2}} $ act trivially on $ \rootgr{\alpha} $ because orthogonal roots in $ A_2 \times A_1 $ are adjacent. Thus from now on, we can assume that $ \roots' $ is of type $ B_3 $ or $ C_3 $.
	
	Choose an arbitrary root base $ \rootbase' $ of $ \roots' $. By \cref{rootsys:subsys-bas-conj}, there exists $ u \in \Weyl(\roots) $ such that $ \rootbase' $ is a subset of $ \rootbase^u $. By \cref{rootsys:simple-gen-weyl}, we can find a word $ \word{\zeta} = \tup{\zeta}{k} $ over $ \rootbase \union (-\rootbase) $ such that $ u = \reflbr{\word{\zeta}} $. Choose roots $ \delta_1, \delta_2, \delta_3 \in \rootbase $ such that $ \rootbase' = \Set{\delta_1', \delta_2', \delta_3'} $ where $ \delta_i' \defl \delta_i^u $ for all $ i \in \Set{1,2,3} $. Further, define $ w'_{\delta_i'} \defl w_{\delta_i}^{w_{\word{\zeta}}} $ for all $ i \in \Set{1,2,3} $, so that $ (w'_{\delta'})_{\delta' \in \rootbase'} $ is a $ \rootbase' $-system of Weyl elements in the group generated by $ (\rootgr{\rho})_{\rho \in \roots'} $. Since the assertion is known to be true for $ \roots \in \Set{B_3, C_3} $, there exist words $ \preworddelta{1}', \preworddelta{2}' $ over $ \rootbase' \union (-\rootbase') $ such that $ \reflbr{\preworddelta{1}'} = \reflbr{\beta_1} $, $ \reflbr{\preworddelta{2}'} = \reflbr{\beta_2} $ and such that the actions of $ w'_{\preworddelta{1}'} $ and $ w'_{\preworddelta{2}'} $ on $ \rootgr{\alpha} $ are identical. For each $ i \in \Set{1,2} $, we can write $ \preworddelta{i}' $ as
	\[ \brackets[\big]{\epsilon_{i,1}\delta'_{n(i, 1)}, \ldots, \epsilon_{i,m_i}\delta'_{n(i, m_i)}} \]
	for integers $ m_i \in \Nzero $, $ n(i,1), \ldots, n(i, m_i) \in \Set{1,2,3} $ and $ \epsilon_{i,1}, \ldots, \epsilon_{i,m_i} \in \Set{\pm 1} $ because $ \preworddelta{i}' $ is a word over $ \rootbase' \union (-\rootbase') $. Now we can define
	\[ \preworddelta{i} \defl \brackets[\big]{\word{\zeta}^{-1}, \epsilon_{i,1}\delta_{n(i, 1)}, \ldots, \epsilon_{i,m_i}\delta_{n(i, m_i)}, \word{\zeta}} \quad \text{for all } i \in \Set{1,2}. \]
	In the following, for any word $ \word{\rho} $ over $ \rootbase \union (-\rootbase) $ we will write $ w(\word{\rho}) $ instead of $ w_{\word{\rho}} $ for better legibility, and similarly for $ w'(\word{\rho}') $. Then for all $ i \in \Set{1,2} $, we have
	\begin{align*}
		w(\preworddelta{i}) &= \brackets*{\prod_{j=1}^m w(\delta_{n(i,j)})^{\epsilon_{i,j}}}^{w(\word{\zeta})} = \prod_{j=1}^m w(\delta_{n(i,j)})^{\epsilon_{i,j}w(\word{\zeta})} = \prod_{j=1}^m w'(\delta'_{n(i,j)})^{\epsilon_{i,j}} \\
		&= w'(\preworddelta{i}')
	\end{align*}
	and, similarly,
	\begin{align*}
		\reflbr{\preworddelta{i}} &= \brackets*{\prod_{j=1}^{m_i} \reflbr{\delta_{n(i, j)}}}^{\reflbr{\word{\zeta}}} = \prod_{j=1}^{m_i} \reflbr{\delta_{n(i, j)}}^{\reflbr{\word{\zeta}}} = \prod_{j=1}^{m_i} \sigma\brackets[\big]{\delta_{n(i, j)}^{\reflbr{\word{\zeta}}}} = \prod_{j=1}^{m_i} \reflbr{\delta_{n(i, j)}'} \\
		&= \reflbr{\preworddelta{i}'} = \reflbr{\beta_i}.
	\end{align*}
	We conclude that the words $ \preworddelta{1} $ and $ \preworddelta{2} $ have the desired properties.
\end{proof}

\begin{remark}\label{F4:word-same-act:parmap}
	Denote by $ \inverparsym $ and $ \invoparsym $ the standard $ \rootbase $-parity maps for $ F_4 $. By the uniqueness of these maps, the words $ \preworddelta{1} $ and $ \preworddelta{2} $ in \cref{F4:word-same-act} automatically satisfy $ \inverparbr{\alpha}{\preworddelta{1}} = \inverparbr{\alpha}{\preworddelta{2}} $ and $ \invoparbr{\alpha}{\preworddelta{1}} = \invoparbr{\alpha}{\preworddelta{2}} $.
\end{remark}

\begin{lemma}\label{F4:short-invo-triv}
	Let $ G $ be a group with a crystallographic $ F_4 $-grading $ (\rootgr{\alpha})_{\alpha \in \roots} $ and let $ \roots_C $ be a parabolic subsystem of $ F_4 $ of type $ C_3 $. Let $ H $ be the corresponding $ \roots_C $-graded subgroup of $ G $. Then the short involution on $ H $ (which, by \cref{BC:short-invo-C-note}, is defined on the long root groups of $ H $) is trivial.
\end{lemma}
\begin{proof}
	Let $ \alpha $ be a long root in $ \roots $ and let $ w_\alpha $ be any $ \alpha $-Weyl element. By an application of \cref{F4:rootsys:inBC}, we find a parabolic $ A_2 $-subsystem $ \roots_A $ of $ F_4 $ which contains $ \alpha $. Hence it follows from \cref{A2Weyl:square-act-lemma-self} that $ w_\alpha^2 $ acts trivially on $ \rootgr{\alpha} $. This action is precisely the short involution on $ \rootgr{\alpha} $ (by \cref{BC:def:short-invo}), which finishes the proof.
\end{proof}

\begin{lemma}\label{F4:square-act}
	Let $ G $ be a group with a crystallographic $ F_4 $-grading $ (\rootgr{\alpha})_{\alpha \in \roots} $. Then $ G $ satisfies the square formula for Weyl elements.
\end{lemma}
\begin{proof}
	Let $ \alpha $, $ \beta $ be roots and let $ \roots $ be a parabolic rank-3 subsystem of $ F_4 $ which contains $ \alpha $ and $ \beta $. Then $ \roots $ is of type $ B_3 $, $ C_3 $ or $ A_2 \times A_1 $. In each case, the square formula is satisfied by \cref{A2Weyl:cartan-comp,B:square-act:summary,BC:square-act:summary}. In the case where $ \roots $ is of type $ C_3 $, we have to use in addition that the short involution in $ G $ is trivial, which holds by \cref{F4:short-invo-triv}.
\end{proof}


\section{The Parametrisation}

\label{sec:F4:param}

In this section, we define the standard partial twisting system for $ F_4 $, show that it satisfies all desired compatibility conditions, and use it to parametrise arbitrary $ F_4 $-graded groups.

\begin{secnotation}
	We choose a root base $ \rootbase $ of $ F_4 $ in standard order and we denote by $ (\twistgroup, \inverparsym, \invogroup, \invoparsym) $ the standard partial twisting system of type $ F_4 $ with respect to $ \rootbase $ in the sense of the following \cref{F4:st-pts-def}.
\end{secnotation}

\begin{definition}[Standard partial twisting system]\label{F4:st-pts-def}
	The \defemph*{standard partial twisting system of type $ F_4 $ (with respect to $ \rootbase $)}\index{twisting system!partial!standard (type F4)@standard (type $ F_4 $)} is the tuple $ (\twistgroup, \inverparsym, \invogroup, \invoparsym) $ where $ \twistgroup \defl \invogroup \defl \compactSet{\pm 1} $ and where $ \inverparsym $, $ \invoparsym $ are the standard $ \rootbase $-parity maps for $ F_4 $. If $ H $ is a group with an $ F_4 $-pregrading, then the \defemph*{standard partial twisting system for $ G $} is the same tuple together with the additional information that $ \twistgroup $ acts on all root groups of $ H $ by inversion.
\end{definition}

We begin with the purely combinatorial properties of the parity maps.

\begin{lemma}\label{F4:parmap-props}
	$ \inverparsym $ is braid-invariant, $ \invoparsym $ is Weyl-invariant and both are adjacency-trivial. Further, $ \inverparsym $ satisfies the square formula.
\end{lemma}
\begin{proof}
	Braid-invariance and adjacency-triviality follow from \cref{param:param-parmap-has-properties2}. The square formula follows from \cref{param:param-parmap-has-properties2} together with \cref{F4:square-act}. The square-invariance of $ \invoparsym $ can be verified computationally, see \cref{F4:ex:GAP}.
\end{proof}

\begin{lemma}\label{F4:parmap-transp}
	$ \inverparsym \times \invoparsym $ is transporter-invariant and $ \inverparsym $, $ \invoparsym $ are independent.
\end{lemma}
\begin{proof}
	We begin with the orbit of long roots. Put $ \hat{\alpha} \defl f_2 $. For all long roots $ \beta $, we know from \cref{F4:ex:parmap-exist} that $ \parmoveset{(\twistgroup \times \invogroup)}{\hat{\alpha}}{\beta} $ is contained in $ \twistgroup \times \compactSet{1_\invogroup} $. Further,
	\[ \inverinvopar{f_2}{f_1 f_1} = (-1, 1) \]
	by \cref{F4:parmap-values}. Hence $ \parmoveset{(\twistgroup \times \invogroup)}{\hat{\alpha}}{\hat{\alpha}} = \twistgroup \times \compactSet{1_\invogroup} $. Thus it follows from criterion~\thmitemref{parmap:transport-invar-char}{parmap:transport-invar-char:weak-bound} (with $ \hat{\beta} \defl \hat{\alpha} $) that $ \inverparsym \times \invoparsym $ is transporter-invariant on the orbit of long roots.
	
	For the orbit of short roots, we consider $ \hat{\beta} \defl f_3 $. Again by \cref{F4:parmap-values}, we have
	\[ \inverinvopar{f_3}{f_4 f_4} = (-1, 1) \midand \inverinvopar{f_3}{(-f_2, f_3, f_2)} = (-1, -1). \]
	Since $ \reflbr{f_4 f_4} $ and $ \reflbr{-f_2, f_3, f_2} $ stabilise $ f_3 $, it follows that $ \parmoveset{(\twistgroup \times \invogroup)}{\hat{\alpha}}{\hat{\alpha}} = \twistgroup \times \invogroup $. Hence again by criterion~\thmitemref{parmap:transport-invar-char}{parmap:transport-invar-char:weak-bound} (with $ \hat{\beta} \defl \hat{\alpha} $), we infer that $ \inverparsym \times \invoparsym $ is transporter-invariant on the orbit of short roots as well.
	
	Finally, the previous computations together with \cref{param:transporter-proj} show that $ \inverparsym $ and $ \invoparsym $ are independent.
\end{proof}

\begin{lemma}\label{F4:parmap:semi-comp}
	$ \invoparsym $ is semi-complete.
\end{lemma}
\begin{proof}
	This follows from the fact that the only subgroups of $ \invogroup $ are $ \compactSet{1} $ and $ \invogroup $, just like in \cref{B:invoparmap-semicomp}.
\end{proof}

We conclude that the standard partial twisting system is indeed a partial twisting system.

\begin{lemma}
	Let $ (G, (\rootgr{\alpha})_{\alpha \in F_4}) $ be a crystallographic $ F_4 $-graded group and choose any $ \rootbase $-system $ (w_\delta)_{\delta \in \rootbase} $ of Weyl elements in $ G $. Then the standard partial twisting system $ (\twistgroup, \inverparsym, \invogroup, \invoparsym) $ for $ G $ (in the sense of \cref{F4:st-pts-def}) is a partial twisting system for $ (G, (w_\delta)_{\delta \in \rootbase}) $ (in the sense of \cref{param:partwist-def}).
\end{lemma}
\begin{proof}
	By \cref{param:pargroup-inv-example}, $ \twistgroup $ is a twisting group for $ (G, (\rootgr{\alpha})_{\alpha \in F_4}) $. The remaining properties are satisfied by \cref{F4:parmap-props,F4:parmap-transp,F4:parmap:semi-comp}.
\end{proof}

It remains to verify the compatibility conditions. We begin with a generalisation of \cref{F4:parmap-ortho-base}.

\begin{lemma}\label{F4:parmap:ortho}
	Let $ \alpha, \beta \in F_4 $ be orthogonal. Then $ \invopar{\alpha}{\reflbr{\beta}} = 1_\invogroup $ if $ \alpha $ and $ \beta $ are crystallographically adjacent and $ \invopar{\alpha}{\reflbr{\beta}} = -1_\invogroup $ otherwise.
\end{lemma}
\begin{proof}
	By \cref{F4:parmap-ortho-base}, the assertion is true if $ \beta $ lies in $ \rootbase $. In general, we know from \cref{rootsys:indiv-in-rootbase} that there exists $ u \in \Weyl(F_4) $ such that $ \delta \defl \beta^u $ lies in $ \rootbase $. Choose a word $ \word{\rho} $ over $ \rootbase $ such that $ \reflbr{\word{\rho}} = u $. Then it follows from \cref{parmap:stab-conj} that
	\begin{align*}
		\inverparbr{\alpha}{\refl{\beta}} &= \inverparbr{\alpha^{\reflbr{\word{\rho}}}}{\refl{\word{\rho}}^{-1} \refl{\beta} \refl{\word{\rho}}}.
	\end{align*}
	Here $ \refl{\word{\rho}}^{-1} \refl{\beta} \refl{\word{\rho}} = \reflbr{\beta^{\reflbr{\word{\rho}}}} = \reflbr{\delta} $. Note that $ \alpha $ and $ \beta $ are orthogonal if and only if $ \alpha^u $ and $ \delta $ are orthogonal, and the same holds for crystallographic adjacency. Hence the general assertion follows from the special case in \cref{F4:parmap-ortho-base}.
\end{proof}

\begin{proposition}\label{F4:square-comp}
	Let $ (G, (\rootgr{\alpha})_{\alpha \in F_4}) $ be a crystallographic $ F_4 $-graded group and let $ (w_\delta)_{\delta \in \rootbase} $ be a $ \rootbase $-system of Weyl elements in $ G $. Then $ G $ is square-compatible with respect to $ \inverparsym $ and $ (w_\delta)_{\delta \in \rootbase} $.
\end{proposition}
\begin{proof}
	We know from \cref{F4:square-act} that $ G $ satisfies the square formula for Weyl elements and from \cref{F4:parmap-props} that $ \inverparsym $ satisfies the square formula. Hence the assertion follows from \cref{param:square-formula-comp}.
\end{proof}

\begin{proposition}\label{F4:stab-comp}
	Let $ (G, (\rootgr{\alpha})_{\alpha \in F_4}) $ be a crystallographic $ F_4 $-graded group and let $ (w_\delta)_{\delta \in \rootbase} $ be a $ \rootbase $-system of Weyl elements in $ G $. Then $ G $ is stabiliser-compatible with respect to $ (\inverparsym, \invoparsym) $ and $ (w_\delta)_{\delta \in \rootbase} $.
\end{proposition}
\begin{proof}
	By \cref{F4:word-same-act,F4:word-same-act:parmap,F4:parmap:ortho}, the conditions of \cref{param:stabcomp-crit-ortho} are satisfied. The assertion follows.
\end{proof}

Finally, we can apply the parametrisation theorem.

\begin{proposition}\label{F4:param-exists}
	Let $ (G, (\rootgr{\alpha})_{\alpha \in F_4}) $ be a crystallographic $ F_4 $-graded group, let $ (w_\delta)_{\delta \in \rootbase} $ be a $ \rootbase $-system of Weyl elements in $ G $ and let $ (\twistgroup, \inverparsym, \invogroup, \invoparsym) $ be the standard partial twisting system for $ G $. Then there exist abelian groups $ (\comring, +) $ and $ (\compalg, +) $ (both equipped, as sets, with an action of $ \twistgroup \times \invogroup $) such that $ G $ is parametrised by $ (\twistgroup \times \invogroup, \compalg, \comring) $ with respect to $ \inverparsym \times \invoparsym $ and $ (w_\delta)_{\delta \in \rootbase} $ and such that the action of $ \twistgroup $ on $ \comring $ and $ \compalg $ is given by group inversion.
\end{proposition}
\begin{proof}
	This follows from the parametrisation theorem (\cref{param:thm}), whose conditions are satisfied by \cref{F4:square-comp,F4:stab-comp}.
\end{proof}


\section{The Coordinatisation}

\label{sec:F4:coord}

\newcommand{\lin}{f}

\begin{secnotation}\label{secnot:F4-coord}
	We denote by $ \rootbase = (f_1, f_2, f_3, f_4) $ a root base of $ F_4 $ in standard order and by $ (V, \cdot) $ the Euclidean space which is generated by $ \rootbase $. We consider a group $ G $ group with a crystallographic $ F_4 $-grading $ (\rootgr{\alpha})_{\alpha \in F_4} $ and a $ \rootbase $-system of Weyl elements $ (w_\delta)_{\delta \in \rootbase} $. We denote by $ (\twistgroup, \inverparsym, \invogroup, \invoparsym) $ the standard partial twisting system for $ G $ (from \cref{F4:st-pts-def}), by $ (\comring, +) $, $ (\compalg, +) $ abelian groups as in \cref{F4:param-exists} (which are equipped, as sets, with actions of $ \twistgroup \times \invogroup $) and by $ (\risom{\alpha})_{\alpha \in F_4} $ the corresponding parametrisation of $ G $. The action of $ -1_\invogroup $ on $ \compalg $ is denoted by
	\[ \map{\compinvmap}{\compalg}{\compalg}{a}{\compinv{a} \defl -1_\invogroup.a}. \]
	We choose elements $ 1_\comring \in \comring $ and $ 1_\compalg \in \compalg $ such that
	\[ w_{f_2} \in \rootgr{-f_2} \risom{f_2}(1_\comring) \rootgr{-f_2} \midand w_{f_3} \in \rootgr{-f_3} \risom{f_3}(1_\compalg) \rootgr{-f_3}. \]
	Put $ \twistgroup' \defl \compactSet{\pm 1}^2 $. We identify $ \twistgroup $ with the subgroup $ \Set{(\pm 1, 1)} $ of $ \twistgroup' $ and we declare that the second component of $ \twistgroup' $ acts trivially on all root groups. Further, we denote by $ G_B $ and $ G_C $ the root graded subgroups of $ G $ which correspond to the canonical root subsystems of types $ B_3 $ and $ C_3 $, respectively (as in \cref{F4:rootsys:canon-sub}).
\end{secnotation}

It follows from \thmitemcref{A2Weyl:weyl}{A2Weyl:weyl:weylel-gives-triple} that $ 1_\comring $ and $ 1_\compalg $ are uniquely determined by $ w_{f_2} $ and $ w_{f_3} $.

\begin{remark}[on the proofs of \cref{F4:B-commrel,F4:C-commrel}]\label{F4:param-different-actions-rem}
	In the proof of \cref{F4:B-commrel}, we will apply the parametrisation theorem for $ B_3 $-graded groups (\cref{B:param-exists}) to the subgroup $ G_B $ of $ G $. A priori, this yields a new group $ \comring' $ which parametrises the long root groups and a new group $ \compalg' $ which parametrises the short root groups. Further, the groups $ \comring' $ and $ \compalg' $ are, as sets, equipped with an action of $ \twistgroup \times \invogroup $. Afterwards, \cref{B:thm} provides the commutator relations in $ G_B $ with respect to $ \comring' $ and $ \compalg' $.
	
	By \cref{param:thm-strengthen-rem}, we can choose $ \comring' $ to be equal (and not just isomorphic) to $ \comring $ and $ \compalg' $ to be equal to $ \compalg' $. However, it is a priori not clear that the actions of $ \twistgroup \times \invogroup $ on $ \comring' $ and $ \compalg' $ agree with the respective actions on $ \comring $ and $ \compalg $ after this identification. In other words, the notation $ (a,b).x $ is (possibly) ambiguous for $ a \in \twistgroup $, $ b \in \invogroup $ and $ x \in \comring \union \compalg $. For the action of $ \twistgroup $, no ambiguity arises because we know that $ -1_\twistgroup $ acts by group inversion in both cases. Further, the action of $ \invogroup $ on $ \comring $ will never be relevant (see \cref{param:semi-comp-note}). However, it will turn out that the action of $ \invogroup $ on $ \compalg $ does not, in general, agree with the action of $ \invogroup $ on $ \compalg' $. This inconvenience is not relevant in the proof of \cref{F4:B-commrel}, so we delay a more detailed investigation of this phenomenon until \cref{F4:B-rinv-note}.
	
	A similar remark holds for the proof of \cref{F4:C-commrel}, where we apply the parametrisation theorem for $ C_3 $-graded groups (\cref{BC:param-exists}). One notable difference is that the parametrisation theorem for $ C_3 $-graded groups provides actions of the large group $ \twistgroup' \times \invogroup $ on the parametrising structure, not of the group $ \twistgroup \times \invogroup $. However, we can quickly focus our interest on the action of this subgroup. Further, we will have to investigate the relationship between the two actions of $ \invogroup $ on $ \compalg $ during the proof of \cref{F4:C-commrel} (and not, as in \cref{F4:B-commrel}, as an afterthought). It will turn out that the two actions of $ \invogroup $ on $ \compalg $ are, in fact, identical.
\end{remark}

\begin{lemma}\label{F4:B-commrel}
	There exist a commutative associative ring structure on $ \comring $ whose identity element is $ 1_\comring $, a (right) $ \comring $-module structure $ \scmult $ on $ \compalg $ and a $ \comring $-quadratic form $ \map{\compnorm}{\compalg}{\comring}{}{} $ with $ \compnorm(1_\compalg) = 1_\comring $ such that the following commutator relations hold for all $ \lambda, \mu \in \comring $ and all $ c,d \in \compalg $, where $ \scp{}{} $ denotes the linearisation of $ \compnorm $:
	\begin{align*}
		\commutator{\risom{f_1}(\lambda)}{\risom{f_2}(\mu)} &= \risom{f_1 + f_2}(-\lambda \mu), \\
		\commutator{\risom{f_2}(\lambda)}{\risom{f_3}(c)} &= \risom{f_2 + f_3}(-c \scmult \lambda) \risom{f_2 + 2f_3}\brackets[\big]{-\lambda\compnorm(c)}, \\
		\commutator{\risom{f_2 + f_3}(c)}{\risom{f_3}(d)} &= \risom{f_2 + 2f_3}\brackets[\big]{\scp{c}{d}}.
	\end{align*}
\end{lemma}
\begin{proof}
	Denote by $ \rootbaseB \defl (b_1, b_2, b_3) \defl (f_1, f_2, f_3) $ the ordered root base of the canonical $ B_3 $-subsystem of $ F_4 $ and by $ \map{\inverparsymB}{B_3 \times \rootbaseB}{\twistgroup}{}{} $ and $ \map{\invoparsymB}{B_3 \times \rootbaseB}{\invogroup}{}{} $ the $ \rootbaseB $-parity maps which we defined in~\ref{B:ex-twistgroups-def}. This means that, if we choose an orthonormal basis $ (e_1, e_2, e_3, e_4) $ of $ V $ such that $ b_1 = e_1 - e_2 $, $ b_2 = e_2 - e_3 $ and $ b_3 = e_3 $, then the maps $ \inverparsymB $ and $ \invoparsymB $ are given by the exact same formulas as in \cref{B:ex-parmap-def} on page~\pageref{B:ex-parmap-def}. Observe that $ (\twistgroup, \inverparsymB, \invogroup, \invoparsymB) $ is the standard partial twisting system for $ G_B $ with respect to $ \rootbaseB $ (from \cref{B:standard-partwist-def}). 
	
	We know from \cref{B:param-exists} that there exist groups $ (\comring', +) $ and $ (\compalg', +) $ (each equipped with an action of $ \twistgroup \times \invogroup $), an isomorphism $ \map{\risomB{\alpha}}{(\comring', +)}{\rootgr{\alpha}}{}{} $ for each long root $ \alpha $ in $ B_3 $ and an isomorphism $ \map{\risomB{\beta}}{(\compalg', +)}{\rootgr{\beta}}{}{} $ for each short root $ \beta $ in $ B_3 $ such that $ G_B $ is parametrised by $ (\twistgroup \times \invogroup, \compalg', \comring') $ with respect to $ \inverparsymB \times \invoparsymB $ and $ (w_\delta)_{\delta \in \rootbaseB} $ and such that the action of $ \twistgroup $ on $ \compalg' $ and $ \comring' $ is given by group inversion. Denote the corresponding root isomorphisms by $ (\risomB{\alpha})_{\alpha \in B_3} $. By \cref{param:thm-strengthen-rem}, we can achieve that $ \comring' = \comring $, $ \compalg' = \compalg $, $ \risomB{f_2} = \risom{f_2} $ and $ \risomB{f_3} = \risom{f_3} $. (See, however, the warning in \cref{F4:param-different-actions-rem}.)
	
	We now turn to the relationship between the two families of root isomorphisms.
	We know from \cref{F4:parmap-values} that
	\begin{gather*}
		\inverinvopar{f_2}{f_1} = (1,1), \qquad \inverinvopar{f_2}{f_1 f_2} = (1,1), \qquad \inverinvopar{f_2}{f_3} = (-1, 1) \\
		\text{and} \quad \inverinvopar{f_3}{f_2} = (1,1).
	\end{gather*}
	Using the table in \cref{B:ex-parmap-def}, we can also compute that
	\begin{align*}
		\inverinvoparB{f_2}{f_1} &= \inverinvoparB{\basvec_2 - \basvec_3}{\basvec_1 - \basvec_2} = (-1, 1), \\
		\inverinvoparB{f_2}{f_1 f_2} &= \inverinvoparB{\basvec_2 - \basvec_3}{(\basvec_1 - \basvec_2, \basvec_2 - \basvec_3)} \\
		&= \inverinvoparB{\basvec_2 - \basvec_3}{\basvec_1 - \basvec_2} \inverinvoparB{\basvec_1 - \basvec_3}{\basvec_2 - \basvec_3} \\
		&= (-1,1) (-1,1) = (1,1), \\
		\inverinvoparB{f_2}{f_3} &= \inverinvoparB{\basvec_2 - \basvec_3}{\basvec_3} = (-1,1), \\
		\inverinvoparB{f_3}{f_2} &= \inverinvoparB{\basvec_3}{\basvec_2 - \basvec_3} = (-1, 1).
	\end{align*}
	Since
	\[ f_2^{\reflbr{f_1}} = f_1 + f_2, \quad f_2^{\reflbr{f_1 f_2}} = f_1, \quad f_2^{\reflbr{f_3}} = f_2 + 2f_3 \enskip \text{and} \enskip f_3^{\reflbr{f_2}} = f_2 + f_3, \]
	it follows that
	\begin{align*}
		\risom{f_1+f_2}(a) &= \risom{f_2}(a)^{w_{f_1}} = \risomB{f_2}(a)^{w_{f_1}} = \risomB{f_1 + f_2}(-a), \\
		\risom{f_1}(a) &= \risom{f_2}(a)^{w_{f_1} w_{f_2}} = \risomB{f_2}(a)^{w_{f_1} w_{f_2}} = \risomB{f_1}(a), \\
		\risom{f_2 + 2f_3}(u) &= \risom{f_2}(-u)^{w_{f_3}} = \risomB{f_2}(-u)^{w_{f_3}} = \risomB{f_2 + 2f_3}(u), \\
		\risom{f_2 + f_3}(u) &= \risom{f_3}(u)^{w_{f_2}} = \risomB{f_3}(u)^{w_{f_2}} = \risomB{f_2 + f_3}(-u).
	\end{align*}
	
	We can now investigate the commutator relations. From \cref{B:thm}, we know that there exist a commutative associative ring structure on $ \comring $ whose identity element is $ 1_\comring $, a (left) $ \comring $-module structure $ \scmult' $ on $ \compalg $ and a $ \comring $-quadratic form $ \map{\compnorm}{\compalg}{\comring}{}{} $ with $ \compnorm(1_\compalg) = 1_\comring $ such that $ G_B $ is coordinatised by $ (\compalg, \comring) $ with standard signs. In particular, this means that the following commutator relations hold for all $ \lambda, \mu \in \comring $ and all $ u,v \in \compalg $, where $ \scp{}{} $ denotes the linearisation of $ \compnorm $:
	\begin{align*}
		\commutator{\risomB{f_1}(\lambda)}{\risomB{f_2}(\mu)} &= \risomB{f_1 + f_2}(\lambda \mu), \\
		\commutator{\risomB{f_2}(\lambda)}{\risomB{f_3}(c)} &= \risomB{f_2 + f_3}(\lambda \scmult' c) \risomB{f_2 + 2f_3}\brackets[\big]{-\lambda\compnorm(c)}, \\
		\commutator{\risomB{f_2 + f_3}(c)}{\risomB{f_3}(d)} &= \risomB{f_2 + 2f_3}\brackets[\big]{-\scp{c}{d}}.
	\end{align*}
	Using the formulas above for the root isomorphisms, and replacing the left module structure $ \scmult' $ by the corresponding right module structure defined by $ c \scmult \lambda \defl \lambda \scmult' c $, we conclude that the desired commutator relations hold. This finishes the proof.
\end{proof}

\begin{remark}\label{F4:B-rinv-note}
	Recall from \cref{F4:param-different-actions-rem} that we have two actions of $ \invogroup $ on the module $ \compalg $ from \cref{F4:B-commrel}. The goal of this remark is two compare these two actions. First of all, we have to make our notation precise: For all $ c \in \compalg $, we define $ \compinv{c} \defl -1_\invogroup.c $ where $ -1_\invogroup.u $ denotes the action from \cref{F4:param-exists} (as in \cref{secnot:F4-coord}) and we define $ \compinvB{c} \defl -1_\invogroup.c $ where $ -1_\invogroup.c $ denotes the action from \cref{B:param-exists}.
	
	Recall from \cref{F4:parmap-values} that $ \inverinvopar{f_3}{(-f_2, f_3, f_2)} = (-1_\twistgroup,-1_\invogroup) $. Further, we can compute that
	\begin{align*}
		\inverinvoparB{f_3}{(-f_2, f_3, f_2)} &= \inverinvoparB{\basvec_3}{(\basvec_3- \basvec_2, \basvec_3, \basvec_2 - \basvec_3)} \\
		&= \inverinvoparB{\basvec_3}{\basvec_3 - \basvec_2} \inverinvoparB{\basvec_2}{\basvec_3} \inverinvoparB{\basvec_2}{\basvec_2 - \basvec_3} \\
		&= (1,1) (1, -1) (1,1) = (1,-1).
	\end{align*}
	For all $ c \in \compalg $, it follows that
	\begin{align*}
		\risom{f_3}\brackets{-\compinv{c}} &= \risom{f_3}(c)^{w_{f_2}^{-1} w_{f_3}} = \risomB{f_3}(c)^{w_{f_2}^{-1} w_{f_3}} = \risomB{f_3}(\compinvB{c}) = \risom{f_3}(\compinvB{c}),
	\end{align*}
	so $ \compinvB{c} = -\compinv{c} $.
	
	Since $ \compinvmapB $ is precisely the action of $ \invogroup $ on $ \compalg $ that we have studied in \cref{sec:B-bluecomp}, it follows from \cref{B:blue:short-weyl} that it is the reflection corresponding to $ 1_\compalg $ in the sense of \cref{quadmod:refl-def}. That it, $ \compinvmapB = \refl{1_\compalg} $. Thus it follows from the conclusion of the previous paragraph that $ \compinvmap = -\refl{1_\compalg} $.
\end{remark}

\begin{lemma}\label{F4:C-commrel}
	There exist an alternative ring structure on $ \compalg $ whose identity element is $ 1_\compalg $ and maps $ \jorsc $, $ \str_\compalg $, $ \jorTrone $ such that $ (\comring, \jorsc, \str_\compalg, \jorTrone, 0) $ is a Jordan module of type $ C $ over $ \compalg $, the map $ \compinvmap $ from \cref{secnot:F4-coord} is a nuclear involution of $ \compalg $, $ \str_\compalg(1_\comring) = 1_\compalg $ and the following commutator relations hold for all $ \lambda \in \comring $ and all $ c,d \in \compalg $:
	\begin{align*}
		\commutator{\risom{f_4}(c)}{\risom{f_3}(d)} &= \risom{f_3 + f_4}(cd), \\
		\commutator{\risom{f_2}(\lambda)}{\risom{f_3}(c)} &= \risom{f_2 + 2f_3}\brackets[\big]{-\jorsc(\lambda, -\compinv{c})} \risom{f_2 + f_3}\brackets[\big]{-c\str_\compalg(\lambda)}, \\
		\commutator{\risom{f_2 + f_3}(c)}{\risom{f_3}(c)} &= \risom{f_2 + 2f_3}\brackets[\big]{\jorTrone(c \compinv{d})}.
	\end{align*}
\end{lemma}
\begin{proof}
	Denote by $ \rootbaseC \defl (c_1, c_2, c_3) \defl (f_4, f_3, f_2) $ the ordered root base of the canonical $ C_3 $-subsystem of $ F_4 $ and by $ \map{\inverparsymC}{C_3 \times \rootbaseC}{\twistgroup'}{}{} $ and $ \map{\invoparsymC}{C_3 \times \rootbaseC}{\invogroup}{}{} $ the $ \rootbaseC $-parity maps which we defined in~\ref{BC:ex-twistgroups-def}. This means that, if we choose an orthonormal basis $ (e_1, e_2, e_3, e_4) $ of $ V $ such that $ c_1 = e_1 - e_2 $, $ c_2 = e_2 - e_3 $ and $ c_3 = 2e_3 $, then the maps $ \inverparsymB $ and $ \invoparsymB $ are given by the exact same formulas as in \cref{BC:ex-parmap-def} on page~\pageref{BC:ex-parmap-def}. Further, we denote by $ \map{\varinverparsymC}{C_3 \times \rootbaseC}{\twistgroup \le \twistgroup'}{}{} $ the projection of $ \inverparsymC $ to the first component. Observe that $ (\twistgroup', \inverparsymC, \invogroup, \invoparsymC) $ is the standard partial twisting system for $ G_C $ with respect to $ \rootbaseC $ (from \cref{BC:standard-partwist-def}). (Here we use that the short involution on $ G_C $ is trivial by \cref{F4:short-invo-triv}, so that the action of $ \twistgroup' $ on the root groups has the desired form.)
	
	We know from \cref{BC:param-exists} that there exist an abelian group $ (\comring', +) $ and a group $ (\compalg', +) $ (each group being equipped with an action of $ \twistgroup' \times \invogroup $), an isomorphism $ \map{\risomC{\alpha}}{(\comring', +)}{\rootgr{\alpha}}{}{} $ for each long root $ \alpha $ in $ C_3 $ and an isomorphism $ \map{\risomC{\beta}}{(\compalg', +)}{\rootgr{\beta}}{}{} $ for each short root $ \beta $ in $ C_3 $ such that $ G_C $ is parametrised by $ (\twistgroup' \times \invogroup, \compalg', \comring') $ with respect to $ \varinverparsymC \times \invoparsymC $ and $ (w_\delta)_{\delta \in \rootbaseC} $ and such that the action of $ (-1_\twistgroup, 1_\twistgroup) $ on $ \compalg' $ and $ \comring' $ is given by group inversion, the action of $ (1_\twistgroup, -1_\twistgroup) $ on $ \compalg' $ is trivial and the action of $ (1_\twistgroup, -1_\twistgroup) $ is given by the short involution. Since the short involution on $ G_C $ is trivial, we see that $ (1_\twistgroup, -1_\twistgroup) $ acts trivially on both groups $ \comring' $ and $ \compalg' $. Thus the second component of $ \twistgroup' $ is irrelevant. In a more technical language, this means that $ G_C $ is parametrised by $ (\compalg', \comring') $ with respect to $ \inverparsymC \times \invoparsymC $ (and not $ \varinverparsymC \times \invoparsymC $) and $ (w_\delta)_{\delta \in \rootbaseC} $. By \cref{param:thm-strengthen-rem}, we can achieve that $ \comring' = \comring $, $ \compalg' = \compalg $, $ \risomB{f_2} = \risom{f_2} $ and $ \risomB{f_3} = \risom{f_3} $ (which shows that, in particular, $ \compalg' $ must be abelian).
	
	As in the proof of \cref{F4:B-commrel}, \cref{F4:param-different-actions-rem} applies in our situation. We begin by proving that the two actions of $ \invogroup $ on $ \compalg $ are, in fact, equal. By the same arguments as in \cref{F4:B-rinv-note}, it suffices to show that
	\[ \inverinvopar{f_3}{(-f_2, f_3, f_2)} = \inverinvoparC{f_3}{(-f_2, f_3, f_2)}. \]
	We already know from \cref{F4:parmap-values} that
	\[ \inverinvopar{f_3}{(-f_2, f_3, f_2)} = (-1_\twistgroup,-1_\invogroup). \]
	Further, we can compute that
	\begin{align*}
		\inverinvoparC{f_3}{(-f_2, f_3, f_2)} &= \inverinvoparC{\basvec_2 - \basvec_3}{(-2\basvec_3, \basvec_2 - \basvec_3, 2\basvec_3)} \\
		&= \inverinvoparC{\basvec_2 - \basvec_3}{-2\basvec_3} \inverinvoparC{\basvec_2 + \basvec_3}{\basvec_2 - \basvec_3} \inverinvoparC{\basvec_2 + \basvec_3}{2\basvec_3} \\
		&= \inverinvoparC{\basvec_2 + \basvec_3}{2\basvec_3}^{-1} \inverinvoparC{\basvec_2 + \basvec_3}{\basvec_2 - \basvec_3} \inverinvoparC{\basvec_2 + \basvec_3}{2\basvec_3} \\
		&= (-1,1) (-1,-1) (-1,1) = (-1_\twistgroup,-1_\invogroup).
	\end{align*}
	This shows that the two actions of $ \invogroup $ on $ \compalg $. In the following, we will denote them by $ \compinvmap $, as in \cref{secnot:F4-coord}.
	
	We can now turn to the relationship between the two families of root isomorphisms. We know from \cref{F4:parmap-values} that
	\begin{align*}
		\inverinvopar{f_2}{f_3} &= (-1, 1), & \inverinvopar{f_3}{f_2} &= (1,1), \\
		\inverinvopar{f_3}{f_4} &= (-1,1), & \inverinvopar{f_3}{f_4 f_3} &= (1,1).
	\end{align*}
	Further, we can compute that
	\begin{align*}
		\inverinvoparC{f_2}{f_3} &= \inverinvoparC{2\basvec_3}{\basvec_2 - \basvec_3} = (1,1), \\
		\inverinvoparC{f_3}{f_2} &= \inverinvoparC{\basvec_2 - \basvec_3}{2\basvec_3} = (1,1), \\
		\inverinvoparC{f_3}{f_4} &= \inverinvoparC{\basvec_2 - \basvec_3}{\basvec_1 - \basvec_2} = (-1, 1), \\
		\inverinvoparC{f_3}{f_4 f_3} &= \inverinvoparC{\basvec_2 - \basvec_3}{(\basvec_1 - \basvec_2, \basvec_2 - \basvec_3)} \\
		&= \inverinvoparC{\basvec_2 - \basvec_3}{\basvec_1 - \basvec_2} \inverinvoparC{\basvec_1 - \basvec_3}{\basvec_2 - \basvec_3} \\
		&= (-1, 1) (-1, 1) = (1,1)
	\end{align*}
	and
	\begin{align*}
		f_2^{\reflbr{f_3}} = f_2 + 2f_3, \quad f_3^{\reflbr{f_2}} = f_2 + f_3, \quad f_3^{\reflbr{f_4}} = f_3 + f_4 \enskip \text{and} \enskip f_3^{\reflbr{f_4 f_3}} = f_4.
	\end{align*}
	As in the proof of \cref{F4:B-commrel}, it follows that
	\begin{align*}
		\risom{f_2 + 2f_3}(\lambda) &= \risomC{f_2 + 2f_3}(-\lambda), & \risom{f_2 + f_3}(c) &= \risomC{f_2 + f_3}(c) \\
		\risom{f_3 + f_4}(c) &= \risomC{f_3 + f_4}(c), & \risom{f_4}(u) &= \risomC{f_4}(c)
	\end{align*}
	for all $ \lambda \in \comring $ and all $ c \in \compalg $.
	
	Finally, we can investigate the commutator relations. By \cref{BC:thm}, there exist an alternative ring structure on $ \compalg $ whose identity element is $ 1_\compalg $, a nuclear involution on $ \compalg $ (which coincides with the map $ \compinvmap $ defined in \cref{secnot:F4-coord}) and maps $ \jorsc $, $ \str_\compalg $, $ \jorTrone $, $ \psi $ such that $ \jormodtup = (\comring, \jorsc, \str_\compalg, \jorTrone, \psi) $ is a Jordan module over $ \compalg $ (where $ \psi = 0 $ because $ G_C $ is $ C_3 $-graded), $ \str_\compalg(1_\comring) = 1_\compalg $ and $ G_C $ is coordinatised (as a $ BC_3 $-graded group) by $ \jormodtup $ with standard signs. In particular, this means that the following standard commutator relations from \cref{BC:standard-param-def} hold for all $ \lambda \in \comring $ and all $ c,d \in \compalg $:
	\begin{align*}
		\commutator{\risomC{f_4}(c)}{\risomC{f_3}(d)} &= \risomC{f_3 + f_4}(cd), \\
		\commutator{\risomC{f_2}(\lambda)}{\risomC{f_3}(c)} &= \risomC{f_2 + 2f_3}\brackets[\big]{\jorsc(\lambda, -\compinv{c})} \risomC{f_2 + f_3}\brackets[\big]{-c\str_\compalg(\lambda)}, \\
		\commutator{\risomC{f_3}(c)}{\risomC{f_2 + f_3}(d)} &= \risomC{f_2 + 2f_3}\brackets[\big]{\jorTrone(c \compinv{d})}.
	\end{align*}
	Using the formulas above for the root isomorphisms, the assertion follows.
\end{proof}

We have the following identities which relate the maps from \cref{F4:B-commrel} to the maps from \cref{F4:B-commrel}.

\begin{lemma}\label{F4:coord:overlap}
	The following hold:
	\begin{lemenumerate}
		\item \label{F4:coord:overlap:inv}$ \compinv{c} = -\refl{1_\compalg}(c) $ for all $ c \in \compalg $ where $ \refl{1_\compalg} $ denotes the reflection in the quadratic module $ (\compalg, \compnorm) $ corresponding to $ 1_\compalg $ (in the sense of \cref{quadmod:refl-def}).
		
		\item \label{F4:coord:overlap:sc}$ c \scmult \lambda = c \str_\compalg(\lambda) $ for all $ c \in \compalg $ and all $ \lambda \in \comring $.
		
		\item \label{F4:coord:overlap:phi}$ \jorsc(\lambda, c) = \lambda \compnorm(c) $ for all $ \lambda \in \comring $, $ c \in \compalg $.
		
		\item \label{F4:coord:overlap:scp}$ \scp{c}{d} = \jorTrone(c \compinv{d}) $ for all $ c,d \in \compalg $. In particular, $ \scp{c}{1_\compalg} = \jorTrone(c) $ for all $ c \in \compalg $.
	\end{lemenumerate}
\end{lemma}
\begin{proof}
	Assertion~\itemref{F4:coord:overlap:inv} holds by \cref{F4:B-rinv-note}. For the remaining assertions, let $ \lambda \in \comring $ and $ c,d \in \compalg $ be arbitrary. The rank-2 root subsystem generated by $ \Set{f_2, f_3} $ is contained in both the canonical $ B_3 $-subsystem and the canonical $ C_3 $-subsystem, so it follows from \cref{F4:C-commrel,F4:B-commrel} that we have the following \enquote{overlapping} commutator relations:
	\begin{align*}
		\risom{f_2 + f_3}(-c \scmult \lambda) \risom{f_2 + 2f_3}\brackets[\big]{-\lambda\compnorm(c)} &= \commutator{\risom{f_2}(\lambda)}{\risom{f_3}(c)} \\
		&= \risom{f_2 + 2f_3}\brackets[\big]{-\jorsc(\lambda, -\compinv{c})} \risom{f_2 + f_3}\brackets[\big]{-c\str_\compalg(\lambda)}, \\
		\risom{f_2 + 2f_3}\brackets[\big]{\scp{c}{d}} &= \commutator{\risom{f_2 + f_3}(c)}{\risom{f_3}(d)} = \risom{f_2 + 2f_3}\brackets[\big]{\jorTrone(c \compinv{d})}.
	\end{align*}
	We conclude that
	\begin{align*}
		c \scmult \lambda = c \str_\compalg(\lambda), \quad \lambda \compnorm(c) = \jorsc(\lambda, -\compinv{c}), \quad \scp{c}{d} = \jorTrone(c \compinv{d}).
	\end{align*}
	This proves~\itemref{F4:coord:overlap:sc} and~\itemref{F4:coord:overlap:scp}. For~\itemref{F4:coord:overlap:phi}, observe that $ \compnorm(c) = \compnorm(\refl{1_\compalg}(c)) $ because $ \refl{1_\compalg} $ lies in the orthogonal group $ \Ortho(\compnorm) $ by \cref{quadmod:refl-O}. Since $ \compinv{c} = -\refl{1_\compalg}(c) $ by assertion~\itemref{F4:coord:overlap:inv}, we infer that
	\[ \jorsc(\lambda, -\compinv{c}) = \lambda \compnorm(c) = \lambda \compnorm\brackets[\big]{\refl{1_\compalg}(c)} = \lambda \compnorm(-\compinv{c}). \]
	As the map $ \map{}{}{}{c}{-\compinv{c}} $ is bijective, it follows that~\itemref{F4:coord:overlap:phi} holds.
\end{proof}

\begin{summary}\label{F4:summary}
	Altogether, we have assembled the following algebraic structures which coordinatise $ G $ and describe its commutator relations:
	\begin{remenumerate}
		\item An associative commutative ring $ \comring $ with identity element $ 1_\comring $.
		
		\item An alternative ring $ \compalg $ with identity element $ 1_\compalg $ and with a $ \comring $-module structure $ \scmult $.
		
		\item A $ \comring $-quadratic form $ \map{\compnorm}{(\compalg, \scmult)}{\comring}{}{} $ with $ \compnorm(1_\compalg) = 1_\comring $ and its linearisation $ \map{\scp{}{}}{\compalg \times \compalg}{\comring}{}{} $.
		
		\item A nuclear involution $ \compinvmap $ on $ \compalg $.
		
		\item Maps $ \jorsc $, $ \str_{\compalg} $ and $ \jorTrone $ such that $ (\comring, \jorsc, \str_\compalg, \jorTrone, 0) $ is a Jordan module of type~$ C $.
	\end{remenumerate}
	Further, we know that the identities in \cref{F4:coord:overlap} hold. We conclude that
	\[ \calF \defl \brackets[\big]{\comring, \compalg, \scmult, \compnorm, \scp{}{}, \str_\compalg, \jorsc} \]
	is an $ F_4 $-datum. By \cref{F4:compalg:F4-datum-is-comp}, it follows that the ring $ \compalg $ is a $ \comring $-algebra with structural homomorphism $ \str_\compalg $ and associated scalar multiplication $ \scmult $ and that $ (\compalg, \compnorm) $ is a multiplicative conic alternative algebra over $ \comring $ with conjugation $ \compinvmap $. By \cref{F4:compalg:comp-is-F4-datum}, the multiplicative conic alternative algebra structure on $ \compalg $ determines the remaining objects in $ \calF $.
\end{summary}

\begin{theorem}[Coordinatisation theorem for $ F_4 $]\label{F4:thm}
	Let $ G $ be a group with a crystallographic $ F_4 $-grading $ (\rootgr{\alpha})_{\alpha \in F_4} $ and let $ (w_\delta)_{\delta \in \rootbase} $ be a $ \rootbase $-system of Weyl elements in $ G $. Then there exist a commutative associative ring $ \comring $ and a multiplicative conic alternative algebra $ (\compalg, \compnorm) $ over $ \comring $ such that $ (G, (\rootgr{\alpha})_{\alpha \in F_4}) $ is coordinatised $ (\compalg, \compnorm) $ with standard signs and with respect to $ (w_\delta)_{\delta \in \rootbase} $ (in the sense of \cref{F4:ex:coord:stsigns}).
\end{theorem}
\begin{proof}
	This is a consequence of \cref{F4:summary}.
\end{proof}

	\newpage
	\bookmarksetup{startatroot}
	\printbibliography
	\addcontentsline{toc}{chapter}{\bibname}
	
	\newpage
	\RaggedRight
	
	\printindex
\end{document}